\documentclass[oneside]{memo-l}

 \usepackage{amsmath}       % I think this gives me some symbols
 \usepackage{amsthm}        % Does theorem stuff
 \usepackage{amssymb}       % more symbols and fonts
 \usepackage{empheq}        % Some more extensible arrows, like \xmapsto
  \usepackage{bbm}        
 \usepackage{mathrsfs}      % Sheafy font \mathscr{}
\usepackage{rotating}		% provides "sidewaysfigure"
\usepackage{array}
\usepackage{tikz}
\usetikzlibrary{trees}
\usetikzlibrary[shapes]
\usetikzlibrary[arrows]
\usetikzlibrary{patterns}
\usetikzlibrary{fadings}
\usetikzlibrary{backgrounds}
\usetikzlibrary{decorations.pathreplacing}
\usetikzlibrary{decorations.pathmorphing}
\usetikzlibrary{positioning}

\usepackage[colorlinks, linkcolor=black, citecolor=black,   %blue,
            bookmarksnumbered=true]{hyperref}

\newcounter{Step}

\theoremstyle{plain} %%% Plain Theorem Styles.
\newtheorem{theorem}{Theorem}[chapter]
\newtheorem{lemma}[theorem]{Lemma}
\newtheorem{corollary}[theorem]{Corollary}          
\newtheorem{proposition}[theorem]{Proposition}              
\newtheorem*{openproblem}{Open Problem}     
\newtheorem{step}[Step]{Step}       
\newtheorem*{addendum*}{Addendum}   %%% We will remove this...

\theoremstyle{definition} %%%% Definition-like Commands  
\newtheorem{definition}[theorem]{Definition}
\newtheorem{example}[theorem]{Example}

\newtheorem{hypothesis*}[theorem]{Hypothesis}

\theoremstyle{remark}  %%%% Remark-like Commands
\newtheorem{remark}[theorem]{Remark}
\newtheorem{warning}[theorem]{Warning}

\numberwithin{section}{chapter}
\numberwithin{equation}{chapter}
\numberwithin{figure}{chapter}
\numberwithin{table}{chapter}

\definecolor{CSPcolor}{rgb}{0.0,0.5,0.75}	% Textcolor for CSP

			\tikzset{stringnode/.style={rectangle, inner sep=5pt, draw = black!50, fill=black!10}}

			\tikzset{chambering/.style={blue!40, thick}}
			\tikzset{hyperplane/.style={orange}}
            
			 \newcommand{\bbar}[1]{\setbox0=\hbox{$#1$}\dimen0=.2\ht0 \kern\dimen0 \overline{\kern-\dimen0 #1}}
			  
			   \DeclareMathOperator{\codim}{codim}
			 \DeclareMathOperator{\coker}{coker}
			  \DeclareMathOperator*{\colim}{colim}
			 \DeclareMathOperator{\End}{End}
			 \DeclareMathOperator{\fun}{Fun}

			 \let\hom\relax % kills the old hom, which is lowercase
			 \DeclareMathOperator{\hom}{Hom}
			 \DeclareMathOperator{\Hom}{\mathcal{H}\kern-.125em\mathpzc{om}}
			 
			 \DeclareMathOperator{\id}{id}

			 \DeclareMathOperator{\Proj}{\mathcal{P}\kern-.125em\mathpzc{roj}}
			 \DeclareMathOperator{\pro}{pro\mathrm{-}}

			 \renewcommand{\setminus}{\smallsetminus}

			 \DeclareMathOperator{\st}{st}
			 \DeclareMathOperator{\qst}{qst}

			    \DeclareMathOperator{\triplearrow} {{\; \tikz{ \foreach \y in {0, 0.1, 0.2} \draw [->] (0, \y) -- +(0.5, 0);} \; }} %Makes a stack of 3 right arrows, with appropriate space on either side. 

			 \newcommand{\udot}{\ensuremath{{\lower .183333em \hbox{\LARGE \kern -.05em$\cdot$}}}}

			  \DeclareMathOperator{\alg}	{{\sf Alg}}

			  \DeclareMathOperator{\bicat}	{{\sf Bicat}}
			    \DeclareMathOperator{\catbicat}	{{\sf bicat}}
			    \DeclareMathOperator{\bord}{{\sf Bord}}
			    \DeclareMathOperator{\BORD}{{\sf \mathbb{B}ord}}
			
			 \DeclareMathOperator{\cat}	{\sf Cat}
			  \DeclareMathOperator{\Cat}	{\sf cat}
			 
			     \DeclareMathOperator{\cob}{{\sf Cob}}
			     \DeclareMathOperator{\COB}{{\sf	 \mathbb{C}ob}}

			  \DeclareMathOperator{\comp}{{\sf Comp}}

			 \newcommand{\Fun}{\text{\sf{Fun}}} 
			  \DeclareMathOperator{\frob}	{{\sf Frob}}

			\DeclareMathOperator{\Gray}	{\sf Gray}

			\DeclareMathOperator{\gset}	{\sf gSet}

			  \DeclareMathOperator{\man}	{\sf Man}
			 \renewcommand\mod{\textrm{\sf Mod}}

			   \DeclareMathOperator{\prin}	{{\sf Prin}}

			   \DeclareMathOperator{\res}{{res}}

			 \DeclareMathOperator{\set}	{\sf Set}

			  \DeclareMathOperator{\symbicat}	{{\sf SymBicat}}

			 \DeclareMathOperator{\Top}	{\sf Top}

			 \DeclareMathOperator{\vect}	{\sf Vect}

			\newcommand{\haloman}{\tikz[baseline=(A.base)]{\node [inner sep=0pt] (A) {$\man$}; \draw (A.110) ++(-0.4mm,0.4mm) arc [start angle=-20,	end angle=230,
					x radius=.2cm, y radius=1mm, rotate=30]; }}
			\newcommand{\halomand}{{\tikz[baseline=(A.base)]{\node [inner sep=0pt] (A) {$\man$}; \draw (A.110) ++(-0.4mm,0.4mm) arc [start angle=-20,	end angle=230,
							x radius=.2cm, y radius=1mm, rotate=30]; }^{\raisebox{-4pt}{$\scriptstyle d$}}} }

			 \newcommand{\cA}{\mathcal{A}}
			 \newcommand{\cB}{\mathcal{B}}
			 \newcommand{\cC}{\mathcal{C}}
			 \newcommand{\cD}{\mathcal{D}}
			 \newcommand{\cE}{\mathcal{E}}
			 \newcommand{\cF}{\mathcal{F}}
			 
			 \newcommand{\cH}{\mathcal{H}}

			 \newcommand{\cK}{\mathcal{K}}
			 
			 \newcommand{\cM}{\mathcal{M}}
			 
			 \newcommand{\cO}{\mathcal{O}}
			\newcommand{\cP}{\mathcal{P}}
			
			\newcommand{\cR}{\mathcal{R}}
			\newcommand{\cS}{\mathcal{S}}
			 
			 \newcommand{\cU}{\mathcal{U}} 
			 \newcommand{\cV}{\mathcal{V}}
			 \newcommand{\cW}{\mathcal{W}}

			 \newcommand{\cZ}{\mathcal{Z}}

			 \newcommand{\B}{\mathbb{B}}
			
			\newcommand{\D}{\mathbb{D}}
			 
			 \newcommand{\F}{\mathbb{F}}
			\newcommand{\G}{\mathbb{G}}
			 \newcommand{\N}{\mathbb{N}}
			 
			 \renewcommand{\P}{\mathbb{P}}
			 
			 \newcommand{\R}{\mathbb{R}}
			 \let\oldS\S 						%%% Sets "\oldS" to be the section symbole (the default of "\S").
			 \renewcommand{\S}{\mathbb{S}}		%%% Renames \S to be blackboard

			  \newcommand{\Z}{\mathbb{Z}}

			  \newcommand{\sA}{\mathsf{A}}
			 \newcommand{\sB}{\mathsf{B}}
			 \newcommand{\sC}{\mathsf{C}}
			 \newcommand{\sD}{\mathsf{D}}
			 
			 \newcommand{\sF}{\mathsf{F}}

			 \newcommand{\sM}{\mathsf{M}}

			\newcommand{\sS}{\mathsf{S}}
			 
			 \newcommand{\sU}{\mathsf{U}} 
			 
			 \newcommand{\sW}{\mathsf{W}}

			 \newcommand{\m}{\mathfrak{m}}
			  
			\definecolor{MyBlue}{rgb}{0.0,0.5,0.75}
\tikzset{l/.style={font=\fontsize{8}{8}\selectfont}}
\definecolor{salmon}{rgb}{1,0.47,0.425}

\newcommand{\tensor}{\otimes}
\newcommand{\maps}{\colon}

\newcommand{\rcirc}{ {\stackrel{\mbox{\circle{2}}}{r}} } % for app-defnSymMonBicat.tex

\makeatletter
\def\l@figure{\@tocline{0}{3pt plus2pt}{0pt}{1.8pc}{}}
\makeatother

\makeatletter
\def\l@table{\@tocline{0}{3pt plus2pt}{0pt}{1.8pc}{}}
\makeatother

\makeindex

\begin{document}

\frontmatter

\title[The Classification of 2D Extended TFTs]{The Classification of Two-Dimensional  Extended Topological Field Theories} 
\author[C. J. Schommer-Pries]{Christopher John Schommer-Pries}

%%%%% It is not clear which Address is appropriate...
\address{ 
Max Planck Institute for Mathematics \\
Vivatsgasse 7 \\ 
53111 Bonn \\
Germany
%Department of Mathematics,
% Massachusetts Institute of Technology,
% Building 2, Room 179,
% 77 Massachusetts Avenue,
% Cambridge, MA 02139-4307
%University of California, Berkeley \\
%70 Evans Hall \#3840 \\
% Berkeley, CA 94720-3840 USA
}  % Address where the research was carried out.  
%\curraddr{Department of Mathematics \\
%	Harvard University\\
%	1 Oxford St.\\
%	Cambridge, MA 02138} % Current Address
\email{schommerpries.chris.math@gmail.com} 
\urladdr{https://sites.google.com/site/chrisschommerpriesmath/} 
%\dedicatory{To Natalie Schommer-Pries, for her loving support and constant encouragement. Without her, this would not have been possible.}  
%\date{} % AMS will fill in if necessary
%\translator{...} 
%\thanks{...} % This is a footnote for acknowledgement of support. 

\begin{abstract}
	We provide a complete generators and relations presentation of the 2-dimensional extended unoriented and oriented bordism bicategories as symmetric monoidal bicategories. Thereby we classify these types of 2-dimensional extended topological field theories with arbitrary target bicategory. As an immediate corollary we obtain a concrete classification when the target is the symmetric monoidal bicategory of  algebras, bimodules, and intertwiners over a fixed commutative ground ring. In the oriented case, such an extended topological field theory is equivalent to specifying a separable symmetric Frobenius algebra. 

The body of this work is divided into three chapters. The first chapter gives a detailed treatment of a variant of higher Morse theory. We modify the techniques used in the proof of Cerf theory and the classification of small codimension singularities to obtain a bicategorical decomposition theorem for surfaces. Moreover these techniques produce a finite list of local relations which are sufficient to pass between any two decompositions. We deliberately avoid the use of the classification of surfaces, and consequently our techniques are readily adaptable to  higher dimensions.

The second chapter gives an extensive treatment of the theory of symmetric monoidal bicategories. We prove several foundational results, such as a theorem which provides a simple list of criteria for determining when a morphism of symmetric monoidal bicategories is an equivalence. We also introduce several stricter variants on the notion of symmetric monoidal bicategory. We give a very general treatment of the notion of presentation by generators and relations, which applies not just to the theory of symmetric monoidal bicategories but to many related theories. Finally we provide a host of strictification and cohernece results for symmetric monoidal bicategories. 

In the final chapter we study extended topological field theories. We give a precise treatment of the extended bordism bicategory equipped with additional structure (such as framings or orientations). We apply the results of the previous two chapters to obtain a simple presentation of both the oriented and unoriented bordism bicategories. We describe the general method to obtain such classifications for other choices of structure, and we examine the consequences of our classification when the target is the bicategory of algebras, bimodules, and maps, over a fix commutative ground ring. 
\end{abstract}

\maketitle

\cleardoublepage
\thispagestyle{empty}
\vspace*{13.5pc}
\begin{center}
	To Natalie Schommer-Pries, \\[4pt] 
	for her loving support and constant encouragement. \\[2pt] 
	Without her, this would not have been possible.
 % Dedication text (use \\[2pt] for line break if necessary)
\end{center}
\cleardoublepage

%    Change page number to 7 if a dedication is present.
%\setcounter{page}{7} % was "\setcounter{page}{4}"
\setcounter{page}{3} % I changed to single sided to remove extra blank pages. This is the new setting. 

\setcounter{tocdepth}{3}
\tableofcontents
\listoffigures
\listoftables

%% Preface/ Introduction Include unnumbered chapters (preface, acknowledgments, etc.) here.

\chapter*{Introduction}

%\section*{Introduction}

Recently there have been many exciting developments in the fields of topology, quantum field theory, and higher categories. This work aims to solve a particular problem in the intersection of these three fields, the classification of extended 2-dimensional topological field theories, and has the potential for many future applications and developments. 

A {\em topological quantum field theory} (TQFT), 
\index{topological quantum field theory (TQFT)} as axiomatized by Atiyah \cite{Atiyah88} and Segal \cite{Segal04}, is a functor between two particular categories, a geometric category and an algebraic category. 
The geometric category is a category of {\em bordisms}. There is one such category for each non-negative integer, $d$, leading to a notion of $d$-dimensional TQFT. The $d$-dimensional bordism category has objects which are closed $(d-1)$-dimensional manifolds and morphisms which are $d$-dimensional bordisms between these, i.e., a compact $d$-dimensional manifold whose boundary is divided into two components, the source and the target. These bordisms are taken up to diffeomorphism, relative to the boundary, and composition is given by gluing. 

The algebraic category is usually the category of vector spaces over a fixed ground field. Both of these categories have extra structure making them into symmetric monoidal categories and a topological quantum field theory is required to be a symmetric monoidal functor. The symmetric monoidal structure on the bordism category is given by the disjoint union of manifolds, and the symmetric monoidal structure on vector spaces is given by tensor product. Thus a TQFT, $Z$, comes equipped with canonical isomorphisms $Z(M_1 \sqcup M_2) \cong Z(M_1) \otimes Z(M_2)$. 

There are many reasons to study topological quantum field theories, but one reason is that they exhibit a beautiful relationship between algebra and geometry or topology. It is not surprising that a TQFT provides some relations between the two, since it is a functor between an algebraic and a geometric category, but what is surprising is how seemingly unrelated, yet often well known algebraic structures emerge via TQFTs. Probably the best known example of this is the following theorem.

\begin{theorem}(Folklore) \label{FolkTheorem}
	The 2-dimensional oriented bordism category is the free symmetric monoidal category with a single commutative Frobenius algebra object. In particular the category of oriented 2-dimensional TQFTs is equivalent to the category of commutative Frobenius $k$-algebras.  
\end{theorem}

The history of this theorem is somewhat intricate. The first  statement that 2D TQFTs are commutative Frobenius algebras appears in Dijkgraaf's thesis  \cite{Dijkgraaf89}, but is also listed as a folk theorem by Voronov \cite{Voronov94}. This was surely known to Atiyah and Segal and many others, but mathematical proofs didn't appear until the work of  L. Abrams \cite{Abrams96}, S. Swain \cite{Sawin95}, F. Quinn \cite{Quinn95}, B.~Dubrovin~\cite{Dubrovin96}, and J.~Kock~\cite{Kock04}

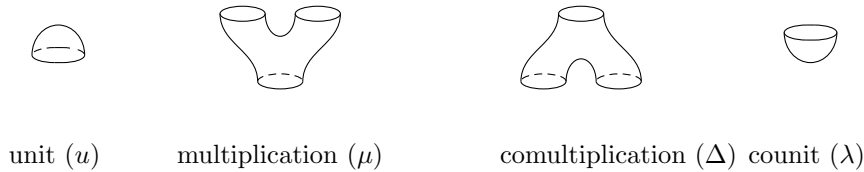
\begin{figure}[ht]
\begin{center}
 \begin{tikzpicture}
 \node at (0, 1) {unit ($u$)};
 \node at (3, 1) {multiplication ($\mu$)};
\node at (7.5, 1) {comultiplication ($\Delta$)};
\node at (10, 1) {counit ($\lambda$)};
 % Draw cap
\begin{scope}[ xshift = 0cm, yshift = 0.25cm]
\draw (.1, 2) -- (0,2);
\draw (0,2) arc (270: 180: 0.3cm and 0.1cm);
\draw [densely dashed] (0, 2.2) arc (90: 180: 0.3cm and 0.1cm); 
\draw [densely dashed] (.1, 2.2) -- (0, 2.2);
\draw (.1,2) arc (-90: 0: 0.3cm and 0.1cm);
\draw [dashed] (0.1, 2.2) arc (90:0: 0.3cm and 0.1cm);
\draw (-.3, 2.1) arc (180: 0: 0.35cm and 0.4cm);
\end{scope}
\begin{scope}[ xshift = 2.5cm, yshift = 0cm]
% CoPants
\draw (0,2.9) ellipse (0.3cm and 0.1cm);
\draw (1, 2.9) ellipse (0.3cm and 0.1cm);
\draw (-0.3, 2.9) to [out = -90, in = 90] (0.2, 2) (0.8, 2) to [out = 90, in = -90] (1.3, 2.9);
\draw (0.2, 2) arc (180: 360: 0.3cm and 0.1cm);
\draw[densely dashed]  (0.2, 2) arc (180: 0: 0.3cm and 0.1cm);
\draw (.3, 2.9) arc (180: 360: 0.2cm and 0.3cm);
\end{scope}
\begin{scope}[xshift = 4.5cm]
% Pants
\draw (2.5,2.9) ellipse (0.3cm and 0.1cm);
\draw (1.7,2) arc (-180: 0: 0.3cm and 0.1cm);
\draw [densely dashed] (2.3, 2) arc (0:180: 0.3cm and 0.1cm);
\draw (1.7, 2) to [out = 90, in = -90] (2.2, 2.9) (2.8, 2.9) to [out = -90, in = 90] (3.3, 2);
\draw (2.7,2) arc (180: 360: 0.3cm and 0.1cm);
\draw [densely dashed] (3.3, 2) arc (0: 180: 0.3cm and 0.1cm); 
\draw (2.3, 2) arc (180: 0: 0.2cm and 0.3cm);
\end{scope}
%Draw cup
\begin{scope}[ xshift = 10cm, yshift = -0.25cm]
\draw (.1, 3) arc (90: -90: 0.3cm and 0.1cm) -- (0, 2.8);
\draw (0, 3) arc (90: 270: 0.3cm and 0.1cm) ;
\draw (0, 3) -- (.1, 3);
\draw (-.3, 2.9) arc (180: 360: 0.35cm and 0.4cm);
\end{scope}
\end{tikzpicture}
\caption{2D TQFTs as Commutative Frobenius Algebras}
\label{2DTQFT=CommFrobAlgFig}
\end{center}
\end{figure}

This theorem is really two theorems packaged into one. First it is a theorem which identifies the bordism category in terms of generators and relations. In particular the 2-dimensional bordism category is equivalent to the abstract symmetric monoidal category generated by a single object, $S^1$, and the morphisms depicted in Figure \ref{2DTQFT=CommFrobAlgFig},
 together with the obvious relations coming from Morse handle slides and birth/death cancelation of Morse critical points.  The second half of this theorem then identifies these generator and relations with familiar algebraic concepts. Thus a 2-dimensional oriented TQFT is completely determined by the vector space $Z(S^1)$ and the linear maps obtained by applying the TQFT to each of the generators. These give linear maps:
\begin{align*}
	u:  Z(\emptyset) = k &\to Z(S^1) \\
	\mu:  Z(S^1) \otimes Z(S^1) &\to Z(S^1) \\
	\Delta:  Z(S^1) &\to Z(S^1) \otimes Z(S^1) \\
	\lambda:  Z(S^1) &\to Z(\emptyset) = k.
\end{align*}
The relations, which are derived from purely topological facts about the bordism category, ensure that these maps give $V = Z(S^1)$ the structure of a commutative Frobenius algebra.  

An {\em extended} topological field theory is a higher categorical version of a TQFT. The symmetric monoidal  bordism category is replaced by a symmetric monoidal bordism $n$-category and an extended topological field theory is a symmetric monoidal $n$-functor from this to a, usually algebraic, target symmetric monoidal $n$-category. 

The simplest example of this is the 2-dimensional bordism bicategory. The objects of this bicategory are (compact) 0-dimensional manifolds, the 1-morphisms are 1-dimensional bordisms between these and the 2-morphisms are 2-dimensional bordisms between the bordisms, taken up to diffeomorphism rel. boundary. 

In this manuscript we provide a generators and relations presentation of the 2-dimensional bordism bicategory as a symmetric monoidal bicategory. This classifies 2-dimensional topological field theories with any target symmetric monoidal bicategory in terms of a small amount of very explicit data. More important than this theorem itself are the methods of proof and the supporting theorems and lemmas.  
Part of these are foundational and all are designed with an eye towards generalizations. This book has been divided into three parts. We have made an effort separate out those aspects which are purely algebraic/categorical from those which are purely geometric/topological and this is reflected in the three chapters of the sequel. 

Chapter~\ref{ChapPlanarDecomp} is written entirely in the language of differential topology and makes no mention whatsoever of symmetric monoidal bicategories or the bordism bicategory, although its main results will be applied in this setting.  It concerns the careful development of \emph{higher Morse theory} which can be used to give a simple two-dimensional description of manifolds. This chapter may be of independent interest.

Chapter~\ref{SymMonBicatChapt} is purely algebraic/categorical and gives an extensive foundational treatment of symmetric monoidal bicategories. After developing the basic theory of symmetric monoidal bicategories we study a very general approach to the theory of presentations of symmetric monoidal bicategories. This approach, which in the literature is called the theory of \emph{generalized computads}, applies not just to the symmetric monoidal bicategories but to almost any similar categorical notion. We introduce two new such notions,  \emph{unbiased semistrict symmetric monoidal 2-categories} and \emph{quasistrict symmetric monoidal 2-categories} which are stricter notions of symmetric monoidal bicategory. Finally we prove several versions of coherence results for symmetric monoidal bicategories, which are essential tools in our classification of topological field theories.

%We fill in these gaps and prove two general theorems about symmetric monoidal bicategories. One of these provides a precise notion of generators and relations for symmetric monoidal bicategories. Given certain generators and relations data, we prove the existence of an abstract symmetric monoidal bicategory which is presented by that data. This symmetric monoidal bicategory satisfies a universal property making it easy to understand symmetric monoidal functors out of it. In particular this gives a precise meaning to our main theorems. Our generators and relations presentation gives an abstract symmetric monoidal bicategory  which, by the universal property, comes equipped with a symmetric monoidal functor to the bordism bicategory. The main theorem proves this is an equivalence of symmetric monoidal bicategories. To aid this proof, we prove another algebraic theorem in Chapter~\ref{SymMonBicatChapt} which provides three relatively easy to check criteria which exactly characterize those symmetric monoidal functors which are symmetric monoidal equivalences.  

In Chapter~\ref{ClassificationChapter} we unite the categorical and the geometric results. We define the symmetric monoidal bordism bicategory and several variations of it. We then prove the main theorems in the unoriented and oriented cases.  Finally, we discuss some of the consequence of these theorems. In particular we completely analyze 2D topological field theories with values in the bicategory of algebras, bimodules and intertwiners (over a fixed commutative ground ring). We identify these in terms of classically studied structures.

\section*{Topological Field Theories}

In the early 1980s there were several significant advances relating modern theoretical physics and modern mathematics. Many were spearheaded by E. Witten and his work on topological aspects of quantum field theory. It was in part  an attempt to better understand the mathematical aspects of some of these developments that led Atiyah and Segal to develop a rigorous mathematical formalism in which to study quantum field theories.

%Modern theoretical physics has had an increasingly influential effect on mathematics over the past three decades. In particular the emphasis of E. Witten's work on topological aspects of quantum field theory has lead to many significant developments in both fields, many of which continue to be influential to this day. It was in part  an attempt to better understand the mathematical aspects of some of these developments that led Atiyah and Segal to develop a rigorous mathematical formalism in which to study quantum field theories. 

	In the late 1980s, Graeme Segal wrote and circulated a manuscript, now published \cite{Segal04}, in which he gave a mathematical definition of conformal field theory. He conceptualized it as a functor from a certain category, roughly a category whose objects are 1-manifolds and morphisms are conformal bordisms, to the category of Hilbert spaces and linear operators. At approximately the same time, Atiyah, in his 1988 paper \cite{Atiyah88}, introduced the definition of a {\em topological} quantum field theory (TQFT). As we mentioned earlier, this is now understood to be equivalent to defining a TQFT as a symmetric monoidal functor $\bord_d \to \vect$. 
	
	In addition to providing a connection between algebra and geometry, TQFTs have at least two other significant uses. Firstly, they serve as toy models for more complicated non-topological field theories. For example Segal's axiomatic approach to conformal field theory is very closely related to  2-dimensional topological field theories and they have many structures and general features in common. By studying the simpler topological theories, one might hope to gleam some structure or insight applicable to more general quantum field theories. 
	
	In the Folklore Theorem~\ref{FolkTheorem} we saw an example of how we can potentially obtain complete information about the entire category of topological field theories. This is particularly important for understanding the third principal application of topological quantum field theories: manifold invariants. Indeed, one of the original motivations for axiomatizing TQFTs was an attempt to capture the mathematical structure of the 3-manifold invariants coming out of quantum Chern-Simons theory. 
	
	A symmetric monoidal functor must send the monoidal unit of the source category to the monoidal unit of the target category (at least up to canonical isomorphism). The monoidal unit of the bordism category is the empty manifold $\emptyset$, and  thus if $Z$ is a TQFT, $Z(\emptyset) \cong k$, the monoidal unit for the category of $k$-vector spaces. Moreover, any closed $d$-manifold, $W$, can be viewed as a bordism from $\emptyset$ to itself and thus $Z(W)$ is an endomorphism of $Z(\emptyset) \cong k$, i.e., an element of $k$, and this element only depends on the diffeomorphism class of $W$. 

Unfortunately these manifold invariants appear to be most interesting in low dimensions, viz. less than or equal to three. A direct check, via examples, shows that TQFTs can distinguish 0-, 1-, and 2-dimensional manifolds. Specifically for any pair of non-isomorphic closed $k$-manifolds $P$ and $Q$ ($k= 0,1$ or $2$) there exists a $k$-dimensional TQFT $Z$ such $Z(P) \neq Z(Q)$. Hence $Z$ distinguishes $P$ and $Q$.  

Work of Freedman-Kitaev-Nayak-Slingerland-Walker-Wang \cite{FKNSWW05} on universal manifold pairings shows that 4-dimensional unitary\footnote{A {\em unitary} TQFT is a particular kind of oriented TQFT, taking values in complex vector spaces. A TQFT is unitary if it intertwines orientation reversal of manifolds and bordisms with complex conjugation, and if it satisfies a certain positivity condition.} TQFTs cannot distinguish smoothly $s$-cobordant 4-manifolds. More precisely for any 4-dimensional unitary TQFT $Z$ and any closed smoothly $s$-cobordant 4-manifolds $P$ and $Q$, it is necessarily the case that $Z(P) = Z(Q)$. This applies, in particular, whenever $P$ and $Q$ are homotopy equivalent and simply connected. So for example if $P$ and $Q$ are the same simply connected topological manifold and only differ by their differentiable structure, then $Z(P) = Z(Q)$ for {\em any} unitary TQFT.\footnote{The ``topological field theories'' arising from Donaldson and Seiberg-Witten invariants do not contradict this for two reasons. First they are not honest TQFTs in the sense of Atiyah as they are not defined for all bordisms, and moreover even if they were, they would lead to non-unitary TQFTs.}

Work of Kreck and Teichner \cite{KT08}, extends this, showing that, in contrast to the 4-dimensional case, unitary 5-dimensional TQFTs can distinguish closed simply connected 5-manifolds. However in dimension six or greater they exhibit a pair (in fact families of such pairs) of simply connected manifolds $P$ and $Q$ with distinct homotopy types such that $Z(P) = Z(Q)$ for any unitary TQFT $Z$. Thus unitary TQFTs of these dimensions cannot distinguish the homotopy type of manifolds.\footnote{By a judicious choice of examples, unitary TQFTs in these dimensions can detect the cohomology of the manifold as an abelian group. The cohomology rings of the Kreck-Teichner examples differ in their cup product structure.} This leaves dimension three, where we have the following open problem. Work of Calegari, Freedman, and Walker \cite{CFW08}, provides evidence that this might be the case, but compare also with Funar \cite{Funar11}. 
\begin{openproblem}
Can 3-dimensional TQFTs distinguish closed 3-manifolds?
\end{openproblem}

One important concept in the study of quantum field theories is the notion of {\em locality}. Most quantum field theories coming from physics are defined in terms of local quantities (functions, sections of bundles, etc.) and hence have a local property not expressed in the Atiyah-Segal axioms. One method of incorporating locality into the Atiyah-Segal axioms, arising out of the work of D. Freed, F. Quinn, K. Walker, J. Baez, and J. Dolan, is to redefine a topological quantum field theory in the language of higher categories. 

A category consists of objects and morphisms between the objects. A higher category consists of objects, morphisms between the objects, morphisms between the morphisms (called 2-morphisms), morphisms between the 2-morphisms (called 3-morphisms), and so on. There are various compositions of these morphisms, and various coherence data associated with these various compositions. Making this precise has become an industry in its own right, and there are many competing definitions of (weak) $n$-category. However, these are all equivalent for $n=2$ and in particular are equivalent to the notion of {\em bicategory} introduced by B\'enebou in 1967 \cite{Benabou67}.  

One prototypical example of an $n$-category is the $d$-dimensional bordism $n$-category. The objects consist of closed $(d-n)$-manifolds, the 1-morphisms consist of $(d- n + 1)$-dimensional bordisms between these, the 2-morphisms consist of $(d-n +2)$-dimensional bordisms between the bordisms, and so on, up to dimension $n$. Thus an {\em extended} topological field theory is a functor of $n$-categories from this bordism $n$-category to some (usually algebraic) target $n$-category.  More precisely, an extended topological field theory should be a {\em symmetric monoidal} functor between symmetric monoidal $n$-categories. Making this precise is, needless to say, extremely laborious. 

One reason to be interested in the study of extended topological field theories is that, again, they serve as a toy model for {\em non-topological} extended quantum field theory. Such extended non-topological QFT has arisen, for example,  in  the Stolz-Teichner program on elliptic cohomology \cite{ST04, ST08}. Their proposal, over simplifying a bit, is that the space of 2-dimensional extended supersymmetric Euclidean field theories should be the classifying space of the generalized cohomology theory {\em topological modular forms} (which is the universal elliptic cohomology theory). Their program has had great success in lower dimensions, where they prove that 0-dimensional supersymmetric Euclidean field theories yield de Rham cohomology \cite{HKST09} and 1-dimensional supersymmetric Euclidean field theories produce K-theory \cite{HST09}. These results have so far depended heavily on being able to identify generators and relations for these more exotic bordism categories. Thus it is desirable to study extended topological field theories in such a way that either generalizes to or gives insight about the extended non-topological setting.  

Another reason to be interested in extended topological field theories  that they are inherently more computable. If we are given an extended topological field theory and we want to evaluate it on, say,  a $d$-dimensional manifold, we may start by chopping the manifold into smaller, more elementary pieces. As the category number increases, these elementary pieces become simpler and fewer in number. In principle, the topological field theory becomes easier to compute.

Finally, an important appearance of extended TQFTs is in the construction of examples of TQFTs. Many known 3-dimensional TQFTs arise out of a method known as the Reshetikhin-Turaev construction. This construction takes as input a certain algebraic gadget called a {\em modular tensor category}. This is a particular kind of braided monoidal category equipped with additional structure. Given a modular tensor category, $M$, the Reshetikhin-Turaev construction provides a recipe for constructing a 3-dimensional TQFT. However the details of this construction actually produce an {\em extended} TQFT. In addition to assigning data to 3-manifolds and surfaces, data is assigned to 1-manifolds as well. Specifically the circle is assigned the category $M$. Surfaces with boundary are then assigned linear functors, for example the second bordism of Figure \ref{2DTQFT=CommFrobAlgFig} is assigned the functor,
\begin{equation*}
	\otimes: M \times M \to M.
\end{equation*}
This means that most examples of 3D topological quantum field theories are actually extended at least to circles. For example, this is true of quantum Chern-Simons theory.  This begs the following open question, although see \cite{FHLT09} for some recent progress.

\begin{openproblem}
Which 3D topological field theories can be extended to points? Specifically, can quantum Chern-Simons theory be extended to points?
\end{openproblem}

\section*{The Hypotheses of Baez and Dolan}

Following work of J. Baez and J. Dolan \cite{BD95} and more recently M. Hopkins and J. Lurie \cite{Lurie09}, the fully-local bordism $n$-category is expected to play a particularly important role in higher category theory. We will make a few comments on this work, but we refer the reader to the original work of Baez and Dolan \cite{BD95} and to the expository survey by Lurie \cite{Lurie09} for a more detailed explanation.

 There are several important conjectures in higher category theory that are explained in \cite{BD95}. 
One of the first of these is well-known in lower dimensions, and concerns the categorification progression: set, category, bicategory, tricategory, \dots. These objects can be equipped with additiional structure. For example a set can be given the structure of a unital  multiplication (so that the set becomes a monoid) or of a unital commutative multiplication (so that the set becomes a commutative monoid). Equivalently, a category with a single objects $*$ consists of the following data:
\begin{itemize}
\item a set $S = \hom(*, *)$,
\item an associative multiplication given by composition:
\begin{equation*}
S \times S = \hom(*, *) \times \hom(*, *) \to \hom(* , *) = S
\end{equation*}
\item a unit for the multiplication. 
\end{itemize}
This is precisely the data of a monoid. Thus a monoid is the ``same thing'' as a category with one object. Similarly, the Eckmann-Hilton argument shows that a commutative monoid is the same thing as a bicategory with one object and one 1-morphism. At this point the structure stabilizes: A tricategory with one object, one 1-morphism, and one 2-morphism is again the same as a commutative monoid. 

Similarly a monoidal category is the same as a single object bicategory. Work of Gordan-Powers-Street \cite{GPS95} shows that a tricategory with one object and one 1-morphism gives data which is essentially equivalent to a braided monoidal category. Following Baez and Dolan, we call an $(n+k)$-category with only one object, one 1-morphism, \dots, and one $k$-morphism a {\em $k$-tuply monoidal $n$-category}. Based on the above well known patterns, and the relationship between higher category theory and homotopy theory (explained in \cite{BD95}), Baez and Dolan conjectured that the notion of $k$-tuply monoidal $n$-category stabilizes when $k = n+2$. This can be organized into a table describing the various sorts of $n$-categories one might encounter, see Table~\ref{TableStabHypoth}. This has been investigated in low dimensions by Cheng and Gurski \cite{MR2342826, MR2839900}.

\begin{table}[ht]
\begin{center}
\begin{tabular}{|c|c| c| c|} \hline
 & $n=0$ & $n= 1$ & $n=2$ \\ \hline
 $k= 0$ & sets & categories & bicategories \\ \hline
 $k=1$ & monoids & monoidal & monoidal  \\
 &  & categories & bicategories \\ \hline
 $k=2$ & commutative  & braided &  braided \\
 & monoids & monoidal & monoidal  \\
 &  & categories & bicategories \\ \hline
  $k = 3$ & `' & symmetric & sylleptic \\
  & & monoidal & monoidal  \\
   &  & categories & bicategories \\ \hline
$k=4$ & `' & `' & symmetric \\
&&& monoidal \\
&&& bicategories \\ \hline
$k=5$ & `' &`'& `' \\ 
&&& \\ &&& \\ \hline
\end{tabular}
\end{center}
\caption{$k$-tuply monoidal $n$-categories}
\label{TableStabHypoth}
\end{table}%

There are many variations of the bordism $n$-category. For example one may equip the bordism $n$-category with various structures, such as orientations, spin structures, or framings. One may also consider bordisms embedded into Euclidian space. Based off of low dimensional examples, Baez and Dolan made the following conjecture.  

\begin{hypothesis*}[Baez-Dolan Tangle Hypothesis] The $n$-category of framed $n$-tangles in $n+k$ dimensions is the free $k$-tuply monoidal $n$-category with duality on one object. 
\end{hypothesis*}

To be a precise conjecture the notions of $n$-category and the notions of duality need to be made precise. For this reason Baez and Dolan call it a hypothesis. In low dimensions, some of these results can be directly checked \cite{BD95, BL03, Lurie09}. The results seem to indicate that equipping the bordism with various structures yields various notions of duality. Combined with the stablization hypothesis of Table~\ref{TableStabHypoth}, this produces several more or less testable conjectures, at least for sufficiently low $n$. In particular we have:

\begin{hypothesis*}[Baez-Dolan Cobordism Hypothesis] The bordism $n$-category is the free symmetric monoidal $n$-category with duality on one object. 
\end{hypothesis*}

Recently there has been substantial progress towards these results due to work of Hopkins and Lurie. In \cite{Lurie09}, Lurie outlines a sophisticated program which reformulates and proves the Baez-Dolan Cobordism Hypothesis. Both the proof sketch and the reformulation of the Cobordism Hypothesis use the language of $(\infty, n)$-categories in an essential way. In particular they exploit the relationship between $\infty$-categories and homotopy theory to great advantage. Their solution is consequently cast in the same language. 

Other previous work on the cobordism hypothesis includes the work of Baez and Langford \cite{BL03}. Building on the work of Carter-Reiger-Saito \cite{CRS97}, Baez and Langford describe generators and relations for a strictified version of the bicategory of unoriented 0-, 1-, and 2-manifolds embedded in $\R^4$ (taken up to isotopy). They also describe a particular kind of duality, which they call {\em unframed duality} and show that this 2-category is, in a particular sense, the free braided monoidal 2-category with (unframed) duality generated by a single object. This is similar in spirit to what is proven here. 

\section*{The Main Theorem}

In this monograph we study 2-dimensional extended topological field theories. We prove the following theorem, and as a consequence classify local 2-dimensional topological field theories with arbitrary target. We also prove an analogous theorem for the oriented bordism bicategory.

\begin{figure}[ht]
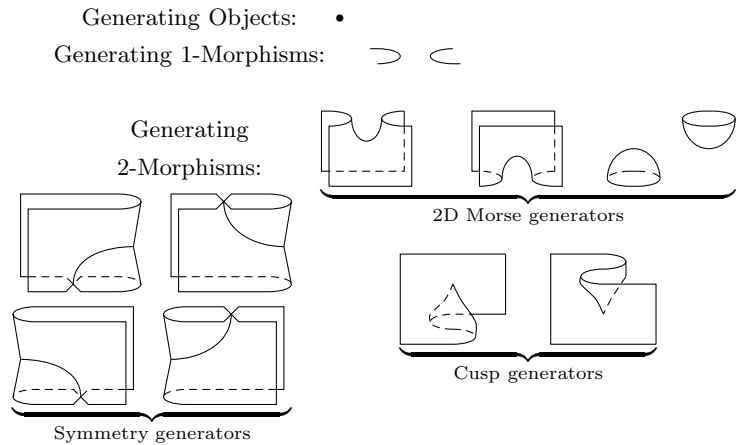

\begin{center}



\caption{Generators}
\label{MainTheormGenFig}
\end{center}
\end{figure}

\begin{theorem}[Classification of Unoriented Topological Field Theories] \label{MainTheoremIntro} 
The 2-dimensional  unoriented bordism bicategory has the following presentation as a symmetric monoidal bicategory: It has the generators depicted in Figure \ref{MainTheormGenFig} subject to the relations depicted in Figure \ref{MainTheormRelFig}.
\end{theorem}

\begin{figure}[ht]
\begin{center}
%
};

\begin{scope}

%%% Cusp Inversion AI
\begin{scope}[xshift = -2cm, yshift = -2.5cm]
\draw (0.5, 1.2) -- (-0.2,1.2) -- (-0.2,-0.5) -- (0.2, -0.5);
\draw (0.5, 1.2) arc (90: -90: 0.3cm and 0.1cm) arc (90: 270: 0.3cm and 0.1cm) -- (0.8, 0.8);
\draw [densely dashed] (0.2, -0.5) -- (0.5,-0.5) arc (90: -90: 0.3cm and 0.1cm) arc (90: 180: 0.3cm and 0.1cm);
\draw (0.2, -0.8) arc (180: 270: 0.3cm and 0.1cm) -- (0.8, -0.9);
\draw (0.8, 0.8) -- (1.2, 0.8) -- (1.2, -0.9) -- (0.8, -0.9);
\node (C) at (0.5, 0.4) {}; 
\node (D) at (0.5, -0.1) {};
\draw  plot [smooth] coordinates{(0.8,1.1) (0.8,0.9)}; 
\draw [densely dashed] plot [smooth] coordinates{ (0.8,0.9) (0.6, 0.6) (C.center)}; 
\draw  plot [smooth] coordinates{(0.2,0.9) (0.2,0.8) (0.4, 0.6) (C.center)}; 
\draw  plot [smooth] coordinates{(0.2,-0.8) (0.2,-0.6) (0.4, -0.4) (D.center)}; 
\draw [densely dashed] plot [smooth] coordinates{(0.8,-0.6) (0.8,-0.5) (0.6, -0.3) (D.center)}; 
\end{scope}

%%% Cusp Inversion AII
\begin{scope}[xshift = 1cm, yshift = -2.5cm]
\draw (0.5, 1.2) -- (-0.2,1.2) -- (-0.2,-0.5) -- (0.2, -0.5); %(0.5, -0.5) ;
\draw (0.5, 1.2) arc (90: -90: 0.3cm and 0.1cm) arc (90: 270: 0.3cm and 0.1cm) -- (0.8, 0.8);
\draw [densely dashed] (0.2, -0.5) -- (0.5,-0.5) arc (90: -90: 0.3cm and 0.1cm) arc (90: 180: 0.3cm and 0.1cm);
\draw (0.2, -0.8) arc (180: 270: 0.3cm and 0.1cm) -- (0.8, -0.9);
\draw (0.8, 0.8) -- (1.2, 0.8) -- (1.2, -0.9) -- (0.8, -0.9);
\node (C) at (0.5, 0.4) {};
\draw (0.2, 0.9) -- (0.2, -0.8) (0.8, 1.1) -- (0.8, -0.6);
\end{scope}

\node at (0, -2.5) {$=$};
\node at (0, -5) {$=$};

%%% Cusp Inversion BI
\begin{scope}[xshift = -2cm, yshift = -4.2cm]
\draw (1.2, 0.4) -- (-0.2,0.4) -- (-0.2,-0.5) -- (0.2, -0.5);
\draw [densely dashed] (0.2, -0.5) -- (0.5,-0.5) arc (90: -90: 0.3cm and 0.1cm) arc (90: 180: 0.3cm and 0.1cm);
\draw (0.2, -0.8) arc (180: 270: 0.3cm and 0.1cm) -- (0.8, -0.9);
\draw  (1.2, 0.4) -- (1.2, -0.9) -- (0.8, -0.9);
\draw (1.2, -0.9) -- (1.2, -1.6) -- (-0.2, -1.6) -- (-0.2, -0.5); 
\node (D) at (0.5, -0.1) {};
\draw  plot [smooth] coordinates{(0.2,-0.8) (0.2,-0.6) (0.4, -0.4) (D.center)}; 
\draw [densely dashed] plot [smooth] coordinates{(0.8,-0.6) (0.8,-0.5) (0.6, -0.3) (D.center)}; 
\node (C) at (0.5, -1.3) {}; 
\draw  plot [smooth] coordinates{(0.8,-0.6) (0.8,-0.7)}; 
\draw [densely dashed] plot [smooth] coordinates{ (0.8,-0.7) (0.6, -1) (C.center)}; 
\draw  plot [smooth] coordinates{(0.2,-0.8) (0.2,-0.9) (0.4, 0.-1.1) (C.center)}; 
\end{scope}

%%% Cusp Inversion BII
\begin{scope}[xshift = 1cm, yshift = -5cm]
\draw  (0.8, 0.8) to [out = 180, in = 0] (0.2, 1.2) -- (-0.2,1.2) -- (-0.2,-0.5) -- (0.2, -0.5) to [out = 0, in = 180] (0.8, -0.9);
%\draw (0.5, 1.2) arc (90: -90: 0.3cm and 0.1cm) arc (90: 270: 0.3cm and 0.1cm) -- (0.8, 0.8);
%\draw [densely dashed] (0.2, -0.5) -- (0.5,-0.5) arc (90: -90: 0.3cm and 0.1cm) arc (90: 180: 0.3cm and 0.1cm);
%\draw (0.2, -0.8) arc (180: 270: 0.3cm and 0.1cm) -- (0.8, -0.9);
\draw (0.8, 0.8) -- (1.2, 0.8) -- (1.2, -0.9) -- (0.8, -0.9);

\end{scope}

\end{scope}

%%%%%%%%%
%Symmetry Relation new
\begin{scope}[xshift = -7cm, yshift = -12.5cm, xscale = 1]
	\draw (.1,-0.9) -- (0, -0.9) --(0,1.3) 
		-- (0.6, 1.3) -- (0.7, 1.2) -- (0.8, 1.3)
		-- (1.2, 1.3) -- (1.3, 1.2) -- (1.4, 1.3)
		-- (2.1,1.3) arc (90: -90:  0.3cm and 0.1cm) 
		-- (1.4, 1.1) -- (1.3, 1.2) -- (1.2, 1.1)
		 -- (0.8, 1.1) -- (0.7, 1.2) -- (0.6, 1.1)
		-- (0.1, 1.1) -- (0.1,-1.1)
		-- (2.1, -1.1) arc (-90: 0:  0.3cm and 0.1cm)  
		-- (2.3, -0.5) -- (2.4, 0.1) 
		(0.1, 0.2) -- (0, 0.2)
	  (0.8, 1.1) -- (0.7, 1.2) -- (0.6, 1.1) -- (0.1, 1.1) -- (0.1,0)
	-- (0.6, 0) -- (0.7, 0.1) -- (0.8, 0)
	 -- (2.1, 0) 	 arc (-90: 0:  0.3cm and 0.1cm)  to (2.3, 0.6) -- (2.4, 1.2); 
\draw [densely dashed] (0.1,0.2) -- (0.6, 0.2) -- (0.7, 0.1) -- (0.8, 0.2) 
	-- (2.1, 0.2) arc (90: 0:  0.3cm and 0.1cm);
\draw (1.3, 1.2) to [out = -90, in = 180] (2.3, 0.6);
\draw (0.7, 1.2) -- (0.7, 0.1) to [out = -90, in = 180] (2.3, -0.5);
\draw [densely dashed] (0.1,-0.9) -- (2.1, -0.9) arc (90: 0:  0.3cm and 0.1cm);
\node at (3, 0) {$=$};	
% identity
\begin{scope}[xshift = 3.5cm]
\draw (.1,-0.9) -- (0, -0.9) --(0,1.3) 
	-- (0.6, 1.3) -- (0.7, 1.2) -- (0.8, 1.3)
	-- (1.2, 1.3) -- (1.3, 1.2) -- (1.4, 1.3)
	-- (2.1,1.3) arc (90: -90:  0.3cm and 0.1cm) 
	-- (1.4, 1.1) -- (1.3, 1.2) -- (1.2, 1.1)
	 -- (0.8, 1.1) -- (0.7, 1.2) -- (0.6, 1.1)
	-- (0.1, 1.1) -- (0.1,-1.1)
	 -- (2.1, -1.1) arc (-90: 0:  0.3cm and 0.1cm)  
	-- (2.4, 1.2)
	(0.7, 1.2) arc (180: 360: 0.3cm and 0.7cm);
\draw [densely dashed] (.1, -0.9) 
	-- (2.1, -0.9) arc (90: 0: 0.3cm and 0.1cm);
\end{scope}
\end{scope}

%Symmetry Relation new2
\begin{scope}[xshift = -8cm, yshift = -12.5cm, xscale = -1]
	\draw (.1,-0.9) -- (0, -0.9) --(0,1.3) 
		-- (0.6, 1.3) -- (0.7, 1.2) -- (0.8, 1.3)
		-- (1.2, 1.3) -- (1.3, 1.2) -- (1.4, 1.3)
		-- (2.1,1.3) arc (90: -90:  0.3cm and 0.1cm) 
		-- (1.4, 1.1) -- (1.3, 1.2) -- (1.2, 1.1)
		 -- (0.8, 1.1) -- (0.7, 1.2) -- (0.6, 1.1)
		-- (0.1, 1.1) -- (0.1,-1.1)
		-- (2.1, -1.1) arc (-90: 0:  0.3cm and 0.1cm)  
		-- (2.3, -0.5) -- (2.4, 0.1) 
		(0.1, 0.2) -- (0, 0.2)
	  (0.8, 1.1) -- (0.7, 1.2) -- (0.6, 1.1) -- (0.1, 1.1) -- (0.1,0)
	-- (0.6, 0) -- (0.7, 0.1) -- (0.8, 0)
	 -- (2.1, 0) 	 arc (-90: 0:  0.3cm and 0.1cm)  to (2.3, 0.6) -- (2.4, 1.2); 
\draw [densely dashed] (0.1,0.2) -- (0.6, 0.2) -- (0.7, 0.1) -- (0.8, 0.2) 
	-- (2.1, 0.2) arc (90: 0:  0.3cm and 0.1cm);
\draw (1.3, 1.2) to [out = -90, in = 180] (2.3, 0.6);
\draw (0.7, 1.2) -- (0.7, 0.1) to [out = -90, in = 180] (2.3, -0.5);
\draw [densely dashed] (0.1,-0.9) -- (2.1, -0.9) arc (90: 0:  0.3cm and 0.1cm);
\node at (3, 0) {$=$};	
% identity
\begin{scope}[xshift = 3.5cm]
\draw (.1,-0.9) -- (0, -0.9) --(0,1.3) 
	-- (0.6, 1.3) -- (0.7, 1.2) -- (0.8, 1.3)
	-- (1.2, 1.3) -- (1.3, 1.2) -- (1.4, 1.3)
	-- (2.1,1.3) arc (90: -90:  0.3cm and 0.1cm) 
	-- (1.4, 1.1) -- (1.3, 1.2) -- (1.2, 1.1)
	 -- (0.8, 1.1) -- (0.7, 1.2) -- (0.6, 1.1)
	-- (0.1, 1.1) -- (0.1,-1.1)
	 -- (2.1, -1.1) arc (-90: 0:  0.3cm and 0.1cm)  
	-- (2.4, 1.2)
	(0.7, 1.2) arc (180: 360: 0.3cm and 0.7cm);
\draw [densely dashed] (.1, -0.9) 
	-- (2.1, -0.9) arc (90: 0: 0.3cm and 0.1cm);
\end{scope}
\end{scope}
%%%%%%%%%

%Symmetry Relation B
\begin{scope}[xshift = -8cm, yshift = -3.5cm]
\draw (.1,.2) -- (0, 0.2)  --(0,1.3)-- (0.6, 1.3) -- (0.7, 1.2) -- (0.8, 1.3)
	 -- (1.3,1.3)  arc (90: -90:  0.3cm and 0.1cm) 
	 -- (0.8, 1.1) -- (0.7, 1.2) -- (0.6, 1.1) -- (0.1, 1.1) -- (0.1,0)
	 -- (1.3, 0) 	 arc (-90: 0:  0.3cm and 0.1cm)  to (1.5, 0.6) -- (1.6, 1.2); 
\draw [densely dashed] (0.1,0.2) 
	-- (1.3, 0.2) arc (90: 0:  0.3cm and 0.1cm);
\draw (0.7, 1.2) to [out = -90, in = 180] (1.5, 0.6);
	\draw [densely dashed] (0.1,-0.9)  -- (0.6, -0.9) -- (0.7, -1) -- (0.8, -0.9)
	-- (1.3, -0.9) arc (90: 0:  0.3cm and 0.1cm);
\draw  (.1,-0.9) -- (0, -0.9) --(0,0.2)
	(0.1, 0) -- (0.1,-1.1)	-- (0.6, -1.1) -- (.7, -1) -- (0.8, -1.1)
	 -- (1.3, -1.1)  arc (-90: 0:  0.3cm and 0.1cm) to (1.5, -0.5) -- (1.6, 0.1)
	 (0.7, -1) to [out = 90, in = 180] (1.5, -0.5); 	
\node at (2, 0) {$=$};	
% identity
\begin{scope}[xshift = 2.5cm]
\draw (.1,-0.9) -- (0, -0.9) --(0,1.3) 
	-- (0.6, 1.3) -- (0.7, 1.2) -- (0.8, 1.3)
	-- (1.3,1.3) arc (90: -90:  0.3cm and 0.1cm) 
	 -- (0.8, 1.1) -- (0.7, 1.2) -- (0.6, 1.1)
	-- (0.1, 1.1) -- (0.1,-1.1)
	-- (0.6, -1.1) -- (.7, -1) -- (0.8, -1.1)
	 -- (1.3, -1.1) arc (-90: 0:  0.3cm and 0.1cm)  -- (1.6, 1.2)
	 (0.7, -1) -- (0.7, 1.2);
\draw [densely dashed] (.1, -0.9) 
	-- (0.6, -0.9) -- (0.7, -1) -- (0.8, -0.9)
	-- (1.3, -0.9) arc (90: 0: 0.3cm and 0.1cm);
\end{scope}
\end{scope}

%Symmetry Relation B2
\begin{scope}[xshift = -9cm, yshift = -3.5cm, xscale = -1]
\draw (.1,.2) -- (0, 0.2)  --(0,1.3)-- (0.6, 1.3) -- (0.7, 1.2) -- (0.8, 1.3)
	 -- (1.3,1.3)  arc (90: -90:  0.3cm and 0.1cm) 
	 -- (0.8, 1.1) -- (0.7, 1.2) -- (0.6, 1.1) -- (0.1, 1.1) -- (0.1,0)
	 -- (1.3, 0) 	 arc (-90: 0:  0.3cm and 0.1cm)  to (1.5, 0.6) -- (1.6, 1.2); 
\draw [densely dashed] (0.1,0.2) 
	-- (1.3, 0.2) arc (90: 0:  0.3cm and 0.1cm);
\draw (0.7, 1.2) to [out = -90, in = 180] (1.5, 0.6);
	\draw [densely dashed] (0.1,-0.9)  -- (0.6, -0.9) -- (0.7, -1) -- (0.8, -0.9)
	-- (1.3, -0.9) arc (90: 0:  0.3cm and 0.1cm);
\draw  (.1,-0.9) -- (0, -0.9) --(0,0.2)
	(0.1, 0) -- (0.1,-1.1)	-- (0.6, -1.1) -- (.7, -1) -- (0.8, -1.1)
	 -- (1.3, -1.1)  arc (-90: 0:  0.3cm and 0.1cm) to (1.5, -0.5) -- (1.6, 0.1)
	 (0.7, -1) to [out = 90, in = 180] (1.5, -0.5); 	
\node at (2, 0) {$=$};	
% identity
\begin{scope}[xshift = 2.5cm]
\draw (.1,-0.9) -- (0, -0.9) --(0,1.3) 
	-- (0.6, 1.3) -- (0.7, 1.2) -- (0.8, 1.3)
	-- (1.3,1.3) arc (90: -90:  0.3cm and 0.1cm) 
	 -- (0.8, 1.1) -- (0.7, 1.2) -- (0.6, 1.1)
	-- (0.1, 1.1) -- (0.1,-1.1)
	-- (0.6, -1.1) -- (.7, -1) -- (0.8, -1.1)
	 -- (1.3, -1.1) arc (-90: 0:  0.3cm and 0.1cm)  -- (1.6, 1.2)
	 (0.7, -1) -- (0.7, 1.2);
\draw [densely dashed] (.1, -0.9) 
	-- (0.6, -0.9) -- (0.7, -1) -- (0.8, -0.9)
	-- (1.3, -0.9) arc (90: 0: 0.3cm and 0.1cm);
\end{scope}
\end{scope}

% Symmetry Relation A
\begin{scope}[xshift = -8cm, yshift = -6.25cm]
\draw (.1,.2) -- (0, 0.2) --(0,1.3) -- (1.3,1.3)  arc (90: -90:  0.3cm and 0.1cm) 
	 -- (0.1, 1.1) -- (0.1,0) -- (0.6, 0) -- (.7, 0.1) -- (0.8, 0) -- (1.3, 0) 
	 arc (-90: 0:  0.3cm and 0.1cm) 	 to  (1.5, 0.6) -- (1.6, 1.2); 
\draw [densely dashed] (0.1,0.2) 
	-- (0.6, 0.2) -- (0.7, 0.1) -- (0.8, 0.2)
	-- (1.3, 0.2) arc (90: 0:  0.3cm and 0.1cm);
\draw (0.7, 0.1) to [out = 90, in = 180] (1.5, 0.6);
\draw [densely dashed] (0.1,-0.9) -- (1.3, -0.9) arc (90: 0:  0.3cm and 0.1cm);
\draw (.1,-0.9) -- (0, -0.9) 	--(0,0.2)
	(0.1, 0) -- (0.1,-1.1) --(1.3, -1.1)  arc (-90: 0:  0.3cm and 0.1cm) 
	 to (1.5, -0.5)  -- (1.6, 0.1) 
	 (.7, 0.1) to [out = -90, in = 180] (1.5, -0.5);
\node at (2, 0) {$=$};	
% identity
\begin{scope}[xshift = 2.5cm]
\draw (.1,-0.9) -- (0, -0.9) --(0,1.3) -- (.3,1.3) arc (90: -90:  0.3cm and 0.1cm) -- (0.1, 1.1) 
	-- (0.1,-1.1) -- (.3, -1.1) arc (-90: 0:  0.3cm and 0.1cm)  -- (.6, 1.2);
\draw [densely dashed] (.1, -0.9) -- (.3, -0.9) arc (90: 0: 0.3cm and 0.1cm);
\end{scope}
\end{scope}

% Symmetry Relation A2
\begin{scope}[xshift = -9cm, yshift = -6.25cm, xscale = -1]
\draw (.1,.2) -- (0, 0.2) --(0,1.3) -- (1.3,1.3)  arc (90: -90:  0.3cm and 0.1cm) 
	 -- (0.1, 1.1) -- (0.1,0) -- (0.6, 0) -- (.7, 0.1) -- (0.8, 0) -- (1.3, 0) 
	 arc (-90: 0:  0.3cm and 0.1cm) 	 to  (1.5, 0.6) -- (1.6, 1.2); 
\draw [densely dashed] (0.1,0.2) 
	-- (0.6, 0.2) -- (0.7, 0.1) -- (0.8, 0.2)
	-- (1.3, 0.2) arc (90: 0:  0.3cm and 0.1cm);
\draw (0.7, 0.1) to [out = 90, in = 180] (1.5, 0.6);
\draw [densely dashed] (0.1,-0.9) -- (1.3, -0.9) arc (90: 0:  0.3cm and 0.1cm);
\draw (.1,-0.9) -- (0, -0.9) 	--(0,0.2)
	(0.1, 0) -- (0.1,-1.1) --(1.3, -1.1)  arc (-90: 0:  0.3cm and 0.1cm) 
	 to (1.5, -0.5)  -- (1.6, 0.1) 
	 (.7, 0.1) to [out = -90, in = 180] (1.5, -0.5);
\node at (2, 0) {$=$};	
% identity
\begin{scope}[xshift = 2.5cm]
\draw (.1,-0.9) -- (0, -0.9) --(0,1.3) -- (.3,1.3) arc (90: -90:  0.3cm and 0.1cm) -- (0.1, 1.1) 
	-- (0.1,-1.1) -- (.3, -1.1) arc (-90: 0:  0.3cm and 0.1cm)  -- (.6, 1.2);
\draw [densely dashed] (.1, -0.9) -- (.3, -0.9) arc (90: 0: 0.3cm and 0.1cm);
\end{scope}
\end{scope}

%%%%%%%%%
% Gluing-Cup Relation A
\begin{scope}[xshift = 2cm, yshift = -10cm, xscale = 1]
\draw (0,1) -- (0.4, 1) -- (0.5, 1.1) -- (0.6, 1) -- (1, 1)
	arc (-90: 90: 0.3cm and 0.1cm)
	-- (0.6, 1.2) -- (0.5, 1.1) -- (0.4, 1.2) 
	-- (0, 1.2) arc (90: 270: 0.3cm and 0.1cm);
\draw (1.3, 1.1) -- (1.3, 0.1) arc (0: -90: 0.3cm and 0.1cm)
	-- (0,0) arc (270: 180:  0.3cm and 0.1cm)
	-- (-0.2, 0.6) -- (-0.3, 1.1)
	(0.5, 1.1) to [out = -90, in = 0] (-0.2, 0.6)
	(-0.3, 0.1) arc (-180: 0: 0.8cm and 0.5cm);
\draw [densely dashed] (-0.3, 0.1) arc (180: 90: 0.3cm and 0.1cm) -- (1,0.2) arc (90: 0: 0.3cm and 0.1cm);
\end{scope}
\node at (1 ,-9.5) {$=$};
% Gluing-Cup Relation A2
\begin{scope}[xshift = 0cm, yshift = -10cm, xscale = -1]
\draw (0,1) -- (0.4, 1) -- (0.5, 1.1) -- (0.6, 1) -- (1, 1)
	arc (-90: 90: 0.3cm and 0.1cm)
	-- (0.6, 1.2) -- (0.5, 1.1) -- (0.4, 1.2) 
	-- (0, 1.2) arc (90: 270: 0.3cm and 0.1cm);
\draw (1.3, 1.1) -- (1.3, 0.1) arc (0: -90: 0.3cm and 0.1cm)
	-- (0,0) arc (270: 180:  0.3cm and 0.1cm)
	-- (-0.2, 0.6) -- (-0.3, 1.1)
	(0.5, 1.1) to [out = -90, in = 0] (-0.2, 0.6)
	(-0.3, 0.1) arc (-180: 0: 0.8cm and 0.5cm);
\draw [densely dashed] (-0.3, 0.1) arc (180: 90: 0.3cm and 0.1cm) -- (1,0.2) arc (90: 0: 0.3cm and 0.1cm);
\end{scope}

% Gluing-Saddle Relation A
\begin{scope}[xshift = 0.5cm, yshift = -8.5cm, xscale = 1]
\draw (0,0) -- (0,2) -- (0.4, 2) -- (0.5, 2.1) -- (0.6, 2) -- (2.4, 2)
	-- (2.4, 0) -- (2,0) arc (270: 180: 0.3cm and 0.1cm)
	-- (1.7, 1.1) arc (180: 270: 0.3cm and 0.1cm) -- (2.4, 1)
	(1.7, 1.1) arc (0: 180: 0.2cm and 0.3cm)
	-- (1.2, 0.6) -- (1.3, 0.1) arc (0: -90: 0.3cm and 0.1cm) -- (0,0);
\draw (0,0.2) -- (-0.1, 0.2) -- (-0.1, 2.2) 
	-- (0.4, 2.2) -- (0.5, 2.1) -- (0.6, 2.2)
	-- (2.5, 2.2) -- (2.5, 0.2) -- (2.4, 0.2);
\draw (-0.1, 1.2) -- (0, 1.2)
 	(0,1) -- (0.4, 1) -- (0.5, 1.1) -- (0.6, 1) -- (1,1) arc (-90: 0: 0.3cm and 0.1cm)
	(2.4, 1.2) -- (2.5, 1.2);
\draw [densely dashed] (0,0.2) -- (1, 0.2) arc (90: 0: 0.3cm and 0.1cm)
		(2.4, 0.2) -- (2, 0.2) arc (90: 180: 0.3cm and 0.1cm)
		(2.4, 1.2) -- (2, 1.2) arc (90: 180: 0.3cm and 0.1cm)
		(0,1.2) -- (0.4, 1.2) -- (0.5, 1.1) -- (0.6, 1.2) -- (1, 1.2) arc (90: 0: 0.3cm and 0.1cm);
\draw (0.5, 2.1) -- (0.5, 1.1) to [out = -90, in = 180] (1.2, 0.6);
\end{scope}
% Gluing-Saddle Relation A2
\begin{scope}[xshift = -1cm, yshift = -8.5cm, xscale = -1]
\draw (0,0) -- (0,2) -- (0.4, 2) -- (0.5, 2.1) -- (0.6, 2) -- (2.4, 2)
	-- (2.4, 0) -- (2,0) arc (270: 180: 0.3cm and 0.1cm)
	-- (1.7, 1.1) arc (180: 270: 0.3cm and 0.1cm) -- (2.4, 1)
	(1.7, 1.1) arc (0: 180: 0.2cm and 0.3cm)
	-- (1.2, 0.6) -- (1.3, 0.1) arc (0: -90: 0.3cm and 0.1cm) -- (0,0);
\draw (0,0.2) -- (-0.1, 0.2) -- (-0.1, 2.2) 
	-- (0.4, 2.2) -- (0.5, 2.1) -- (0.6, 2.2)
	-- (2.5, 2.2) -- (2.5, 0.2) -- (2.4, 0.2);
\draw (-0.1, 1.2) -- (0, 1.2)
 	(0,1) -- (0.4, 1) -- (0.5, 1.1) -- (0.6, 1) -- (1,1) arc (-90: 0: 0.3cm and 0.1cm)
	(2.4, 1.2) -- (2.5, 1.2);
\draw [densely dashed] (0,0.2) -- (1, 0.2) arc (90: 0: 0.3cm and 0.1cm)
		(2.4, 0.2) -- (2, 0.2) arc (90: 180: 0.3cm and 0.1cm)
		(2.4, 1.2) -- (2, 1.2) arc (90: 180: 0.3cm and 0.1cm)
		(0,1.2) -- (0.4, 1.2) -- (0.5, 1.1) -- (0.6, 1.2) -- (1, 1.2) arc (90: 0: 0.3cm and 0.1cm);
\draw (0.5, 2.1) -- (0.5, 1.1) to [out = -90, in = 180] (1.2, 0.6);
\end{scope}
\node at (-0.25, -7.5) {$=$};
\end{tikzpicture}
% \end{center}
\caption{Relations}
\label{MainTheormRelFig}	
\end{center}
\end{figure}

In dimension 2, the Folklore Theorem~\ref{FolkTheorem} implies that 2-dimensional topological field theories do indeed distinguish 2-manifolds. This theorem says that these field theories are the same as commutative Frobenius algebras, which are well understood and plentiful. By writing down a few examples of Frobenius algebras it is not difficult to distinguish all surfaces.  In order to answer the above open  questions, or to have applications to {\em non-topological} field theories, it would highly desirable to have proofs of Theorem~\ref{FolkTheorem} and Theorem~\ref{MainTheoremIntro} which generalize to higher dimensions. Unfortunately the proof of Theorem~\ref{FolkTheorem} which is usually given has essentially no hope of generalizing beyond 2-dimensions. 

All published proofs of  Theorem~\ref{FolkTheorem}  proceed in the following way. First one uses Morse theory to deduce that any surface bordism can be decomposed into a composition of (disjoint unions of) the elementary bordisms given in Figure \ref{2DTQFT=CommFrobAlgFig}. Hence those elementary bordisms generate the bordism category. The difficulty is then proving that the obvious relations are in fact sufficient, i.e., that two surfaces built from the elementary bordisms are diffeomorphic if and only if they can be related by a finite number of these relations.

The usual proof of this, as suggested by L. Abrams \cite{Abrams96} and fully explored by J. Kock \cite{Kock04}, is to appeal to the classification of surfaces. The classification of oriented surfaces says that two closed oriented surfaces are diffeomorphic if and only if they have the same genus. There is an analogous statement for surfaces with boundary. This allows one to introduce a {\em normal form} for each surface, and one merely needs to prove that the relations are sufficient to reduce to the normal form. This can be done in an ad hoc fashion. 

The problem with this method of proof is that it relies essentially on the classification of surfaces, and hence has no hope of generalizing to higher dimensions. While it is possible to give a proof of our main theorems  along these lines, such a proof would have little value. 

\section*{Higher Morse Theory}

%Chapter~\ref{ChapPlanarDecomp} of this manuscript concerns the development of a body of techniques which can be used to obtain complete combinatorial descriptions of manifolds. These techniques are a form of Higher Morse Theory and are geometric in nature, and do not rely on any of the category theoretical foundations developed later. As such this section may be of independent interest.  The results of this chapter do not rely on the classification of surfaces, and similar results can be adapted to any dimension. 

Fortunately, there is a second proof of the Folk Theorem~\ref{FolkTheorem}, which is the one that is generalized here. As before it begins by looking at Morse functions for surfaces to derive generators. A Morse function (which is a generic map to $\mathbb{R}$) gives a linear, and hence 1-categorical way to decompose a surface. Analyzing the singularities of Morse functions gives 1-categorical generators. Each Morse function gives a decomposition into elementary generators, and conversely such a decomposition corresponds to a Morse function. Thus determining a sufficient set of relations to pass from one decomposition to another is tantamount to determining a sufficient set of relations relating any two Morse functions. 

The new insight is to not use the classification of surfaces to get sufficient relations, but to use Cerf theory \cite{Cerf70}. The part of Cerf theory that is relevant is essentially a {\em family} version of Morse theory (a.k.a. parametrized Morse theory). In particular, a path of functions is a 1-parameter family of functions. Cerf theory then studies generic {\em paths} of functions. A generic path of functions consists of Morse functions for all but a finite number of isolated times. Analyzing these critical times more closely, we see that any two Morse functions can be related by a finite number of handle slides and birth/death moves. This proves the sufficiency of the relations. This proof approach was suggested by S. Sawin \cite{Sawin95}, but full details appear in the work of Moore and Segal, in the appendix of \cite{MS06}.

The relevant theorems of Cerf theory can be proven using classical jet transversality methods developed by Thom, Boardman, Mather, and others. These techniques can be adapted and generalized to obtain a 2-categorical generators and relations decomposition theorem as well. This is the subject of Chapter~\ref{ChapPlanarDecomp}. More specifically, given manifolds $X$ and $Y$, one may consider the {\em jet bundles} for maps from $X$ to $Y$. 
\begin{center}
\begin{tikzpicture}[thick]
	\node (LT) at (0,1.5) 	{$J^k(X, Y)$ };
	\node (LB) at (0,0) 	{$X$};
	%\node (RT) at (2,1.5) 	{$$};
	%\node (RB) at (2,0)	{$$};
	\draw [->] (LT) --  node [left] {$$} (LB);
	%\draw [->>] (LT) -- node [above] {$$} (RT);
	%\draw [right hook->] (RT) -- node [right] {$$} (RB);
	\draw [->, bend left, dashed] (LB) to node [left] {$j^kf$} (LT);
\end{tikzpicture}
\end{center}
 A given smooth map, $f: X \to Y$,  it induces a section, $j^kf$, of this bundle. If we stratify the jet bundle in some way, then Thom's jet transversality theorem ensures that the for generic maps the jet sections, $j^kf$, are transverse to these strata. These strata effectively classify the singularity type of the map $f$. If the strata are chosen judiciously, then one may derive normal coordinates around any given singularity. A more careful analysis allows one to piece together this local data to determine a global decomposition of $X$.

For Morse theory and Cerf theory one studies the jet bundles $J^k(X, \R)$ and $J^k( X \times I, \R \times I)$, respectively. The latter is important because a path of functions can be regarded as a particular kind of map $X \times I \to \R \times I$. The Thom-Boardman stratification of the jet bundles is essentially  sufficient to derive the basic results of Cerf theory, in this case.  

In the 2-categorical setting, we want to have a 2-dimension, 2-categorical decomposition of our surface, $\Sigma$. For this reason we study the jet bundles $J^k( \Sigma, \R^2)$. The usual classification of singularities, using the Thom-Boardman stratification, is not quite sufficient to obtain a 2-categorical decomposition. Instead we use the projection $\R^2 \to \R$ to obtain a {\em finer stratification} of the jet bundles, and hence a finer classification of singularities. This new stratification and the normal coordinates for each singularity type are discussed in Chapter~\ref{ChapPlanarDecomp} and yield a complete list of generators for the bordism bicategory. 

To obtain relations, we study paths of maps to $\R^2$. More precisely we study maps $f:\Sigma \times I \to \R^2 \times I$. We use the projections $\R^2 \times I \to \R \times I \to I$ to stratify the jet bundles $J^r( \Sigma \times I, \R^2 \times I)$ and thereby classify the generic singularities of maps $f$. Any two generic maps $\Sigma \to \R^2$ may be connected by a generic map $f: \Sigma \times I \to \R^2 \times I$, and the analysis of the singularities of $f$ give us a finite list of local moves which allow us to pass from one 2-categorical decomposition to another. This gives us the sufficiency of the relations. These results are summarized in the main planar decomposition theorems of Chapter~\ref{ChapPlanarDecomp}, Theorems \ref{PlanarDecompMovesThm} and \ref{PlanarDecompositionTheorem}.

%The results of Chapter~\ref{ChapPlanarDecomp}, while intended for bicategorical applications, are independent of the results on symmetric monoidal bicategories discussed in Chapter~\ref{SymMonBicatChapt} (indeed they come first in the organization of this manuscript). In Chapter~\ref{ClassificationChapter} we discuss results which unite our algebraic results on symmetric monoidal bicategories and the topological results of Chapter~\ref{ChapPlanarDecomp}.  First we introduce the symmetric monoidal bicategory of bordisms. 

\section*{Symmetric Monoidal Bicategories}

Chapter~\ref{SymMonBicatChapt} concerns the development of the theory of symmetric monoidal bicategories. The history of symmetric monodial bicategories is notoriously complicated and involves treatments and retreatments by several mathematicians \cite{KV94, KV94-2,BN96, DS97, GPS95, MR1626844,  McCrudden00, MR2717302, Stay:2013aa}. We begin by reviewing the fully weak theory of symmetric monoidal bicategories. 
 
We prove a key a theorem which characterizes in simple terms precisely those symmetric monoidal functors (also called symmetric monoidal {\em homomorphisms}) which are symmetric monoidal equivalences. We call this ``Whitehead's theorem for symmetric monoidal bicategories'' (Theorem~\ref{WhiteheadforSymMonBicats}) because it precisely mirrors Whitehead's theorem in topology. 

 Whitehead's Theorem states that a map between reasonable topological  spaces is a homotopy equivalence if and only if it induces an isomorphism of all homotopy groups and all base points. There is a well known relationship between $n$-categories in which all morphisms are invertible ($n$-groupoids) and homotopy $n$-types, i.e., reasonable spaces in which all homotopy groups above $n$ vanish. The theories of these two objects are essentially the same, a statement which has become known as the {\em homotopy hypothesis}. For small $n$,  this can be made into a precise statement, which can be verified.  One consequence of this is that a functor between ordinary 1-categorical groupoids is an equivalence if and only if it is an isomorphism on $\pi_0$ and at $\pi_1(x)$ for every object $x$. A similar statement holds for 2-groupoids, which are bicategories in which all 1-morphisms and 2-morphisms are invertible. 
 
 There is an analog of this which is valid in the non-groupoid setting. It is no longer enough to simply check homotopy groups, since we must deal with non-invertible morphisms. The analog of Whitehead's theorem for ordinary categories states that a functor between categories is an equivalence if it is essentially surjective and fully-faithful. These conditions replace the notion of isomorphism of all homotopy groups. There is a well known bicategorical analog as well. A homomorphism of bicategories is an equivalence if and only if it is essentially surjective on objects, essentially full on 1-morphisms, and fully-faithful on 2-morphisms. Just as groupoids correspond to certain spaces, symmetric monoidal groupoids correspond to certain stable, or $E_\infty$-spaces. There is an analog of Whitehead's theorem valid in this setting as well. 
 
 What we prove is that a symmetric monoidal homomorphism is a symmetric monoidal equivalence if and only if it is essentially surjective on objects, essentially full on 1-morphisms, and fully-faithful on 2-morphisms. In other words a symmetric monoidal homomorphism is a symmetric monoidal  equivalence if and only if it is an equivalence of underlying (non-monoidal) bicategories.

In addition to the theory of fully weak symmetric monoidal bicategories, we also discuss several stricter versions of this concept. We introduce two which are new, \emph{unbiased semistrict symmetric monoidal 2-categories} and \emph{quasistrict symmetric monoidal 2-categories}. In the former roughly half the coherence data of a symmetric monoidal bicategory is trivial, in the former nearly all of it is trivial. We also construct an important graphical calculus for the former notion. 

Next we turn to the theory of presentations of symmetric monoidal bicategories. We do this in the extremely general framework of \emph{generalized computads}, which applies not only to symmetric monoidal bicategories but to any theory governed by a finitary monad on globular sets. In particular we also get a theory of presentations for the stricter theories of unbiased semistrict symmetric monoidal 2-categories and quasistrict symmetric monoidal 2-categories. Each presentation gives rise to a very special sort of symmetric monoidal bicategory, which we call \emph{computadic}. Every symmetric monoidal bicategory is equivalent to a computadic one.  

We prove a {\em cofibrancy theorem} for these computadic symmetric monoidal bicategories. This states that functors out of such symmetric monoidal bicategories are easy to understand and classify. They are completely determined up to equivalence by what they do on the generating morphisms of the presentation. Computadic symmetric monoidal bicategories are analogous to CW-complexes in topology. Mapping out of a CW-complex is relatively easy as one only needs to specify where each cell goes. The same is true for computadic symmetric monoidal bicategories. Our classification theorem for topological field theories works by providing an explicit simple computadic symmetric monoidal bicategory equivalent to the bordism bicategory. 

Finally we prove several coherence theorems for symmetric monoidal bicategories. We prove that every symmetric monoidal bicategory is equivalent to a quasistrict symmetric monodial bicateory. Moreover, given a presentation of a symmetric monoidal bicategory we also get an induced quasistrict symmetric monoidal bicategory generated by the same presentation. We show that the natural comparison map between these two is an equivalence of symmetric monoidal bicategories. We also obtain strictification results for functors between symmetric monoidal bicategories.

\section*{Extended Topological Field Theories}

Previous constructions of objects similar to the bordism bicategory have been made in \cite{KL01} and \cite{Morton07}, but neither of these is quite adequate.  First, they essentially ignore the {\em symmetric monoidal} structure on the bordism bicategory, which is essential to define extended topological field theories. This is perhaps understandable, since they did not have available the results of Chapter~\ref{SymMonBicatChapt}. More importantly, neither  \cite{KL01} nor \cite{Morton07}, discuss the issue of gluing bordisms in sufficient detail to adapt their definition to a bicategory of bordisms with structure. 
%The main issue is that when one glues two smooth manifolds along a common boundary, the resulting space, $W \cup_Y W'$, which is a topological manifold, does not have a unique smooth structure extending those on $W$ and $W'$. While any two choices of such smooth structure result in diffeomorphic manifolds, this diffeomorphism is not {\em canonical}. In the 1-categorical bordism category, this issue is easily swept under the rug. But in the higher categorical setting it must be dealt with.   Kerler-Lyubashenko \cite{KL01} and Morton  \cite{Morton07}, narrowly avoid this issue. 

After introducing a naive version of the bordism bicategory, we give an improved version using manifolds which are equipped with a germ of a $d$-dimensional manifold containing them. These can be used for one method of solving the gluing problem and also for defining the bordism bicategory with structure. The types of structures that can be used are very general, and include orientations, spin structures, and $G$-principal bundles. 

Next we prove the main generators and relations classification theorem in the unoriented case, and show how to adapt our techniques  to obtain generators and relations for the bordism bicategory equipped with structures. We use the oriented case as an illustrative example. Finally we discuss some general consequences and applications of these classification results.

 Topological quantum field theories typically have linear or algebraic target categories, and in the extended context a good example of such a target bicategory is given by the bicategory $Alg^2$. The objects are algebras, the 1-morphisms from $A$ to $B$ are $A$-$B$-bimodules with composition given by tensor product of bimodules, and the 2-morphisms are given by bimodule maps. We use the generators and relations results obtained earlier in this chapter to completely classify extended topological field theories with this target, in both the oriented and unoriented settings. We identify the results with familiar classical notions in algebra. 

Finally, we would be remiss if we did not discuss future applications of this work and current related research. As we have emphasized throughout, the importance of this work is not in the 2-dimensional classification result itself, but in the developed techniques and supporting lemmas. While the 2-dimensional case is interesting, the higher dimensional cases are even more so. There are a few variations that are easy to imagine, one may study $d$-dimensional bicategorical theories in which $d > 2$. The Reshetikhin-Turaev construction produces such a theory. The techniques developed here are currently being adapted to this case \cite{BDSV1, BDSV2, BDSV3}. One may also consider $d$-dimensional theories which are $d$-categorical, hence entirely local and extended all the way to points. Again dimension 3 appears to be very interesting \cite{Douglas:2013aa, Douglas:2014aa}. 
 
As for related results, the most significant is the recent exciting work of M. Hopkins and J. Lurie proving the Baez-Dolan cobordism hypothesis in all dimensions. The Baez-Dolan cobordism hypothesis \cite{BD95} is a conjectural statement that says that the fully extended bordism $n$-category has a particularly simple algebraic characterization. Namely, it is the free symmetric monoidal $n$-category with duals on a single object. More precisely, each notion of structure: unoriented, oriented, framed, etc. will correspond to a notion of duality, and the Baez-Dolan cobordism hypothesis states that the bordism $n$-category is the free symmetric monoidal bicategory with that flavor of duality generated by a single object. 

The Hopkins-Lurie work classifies a large class of fully extended topological field theories with structure, in all dimensions. Full details of this work have not been made public, but an expository article was released by Lurie in Jan. 2009 \cite{Lurie09}. Hopkins and Lurie study a more elaborate version of the extended bordism $n$-category. Rather then taking the top dimensional bordisms up to diffeomorphism, they  incorporate the homotopy type of the diffeomorphism groups through the language of $(\infty, n)$-categories. These are a particular kind of $\infty$-category.  Their proof, which is quite different from ours, involves a delicate induction which plays several variations of the bordism category off each other (in particular the framed bordism category and the unoriented bordism category). Moreover, their proof makes use of sophisticated homotopy theoretic techniques, rather then the more elementary approach outlined here, and the statement of their results is expressed in homotopy theoretic language. 

These two approaches to the classification theorem are well matched. While the Hopkins-Lurie approach yields stronger and more conceptual results, our results are more concrete and in many cases are better adapted for constructing examples. Moreover, the techniques of the Hopkins-Lurie approach seem to require  that the extended topological field theories are {\em fully extended} in the sense that they go all the way down to points. There is no such restriction for the techniques we develop. These results compliment each other, and we hope both will be important in the exciting future of extended quantum field theory.

Finally, this work is also closely related to the work of J. Baez and L. Langford \cite{BL03}, who provide a precise generators and relations description of the braided monoidal bicategory of bordisms embedded in $\R^4$. They show that this is the free braided monoidal bicategory on a single {\em unframed} object, although their set up of the bordism bicategory is somewhat different than ours.  We prove an analogous result in the symmetric monoidal case.

%\begin{table}[htdp]
%\caption{Proofs of the Sufficiency of Relations.}
%\begin{center}
%\begin{tabular}{|c|c|} \hline
%Abrams & ``Completeness of the relations follows easily by inspection'' \\ \hline
%Sawin & ``To check that [the decomposition] does not depend on\\
%&  the Morse function, recall that by Cerf theory [cite Cerf] that \\
%& any change in Morse function has the effect of a sequence of \\
%& the following moves: [illustrated 3 moves].'' \\ \hline
%Voronov & States ``Folklore'' Theorem, no mention of relations or sufficiency proof. \\ \hline
%\end{tabular}
%\end{center}
%\label{default}
%\end{table}%

%\begin{addendum*}
%After this manuscript was written it was pointed out to me that there may be some overlap between the results and techniques of Chapter~\ref{ChapPlanarDecomp} and the results and techniques of Carter and Saito \cite{CS98} and Carter, Rieger, and Saito \cite{CRS97}.
%\end{addendum*}

\section*{Acknowledgements}

Above all, I wish to thank my Ph.D. advisor, Peter Teichner, for his support and guidance. 
Without his help this would not have been possible. %I felt that he has always provided an atmosphere of encouragement in which I was free to explore my mathematical flights of fancy. 
I would especially like to thank Peter for suggesting, after hearing an earlier attempt at the main theorems of this work, that I learn Cerf theory. %That crucial piece of advice set in motion a snowball which has grown into this dissertation. 

I would also like to extend my gratitude to Stephan Stolz. It was a conversation with him which initially sparked my serious interest in the 2-dimensional bordism bicategory. Thinking that 2-dimensional extended topological theories were fairly easy,  I was certain that the value associated to the circle should be the center of the algebra associated to the point. Stephan pointed out to me that, on the contrary,  it should be the universal trace space. It turns out we were both right, as Lemma \ref{LemmaCircle=Center} now shows. I would also like to thank him for our numerous conversations since that time. His insights have always been extremely helpful. 

I would  like to give special acknowledgment to Chris Douglas. Our conversations have had a tremendous impact on how I think about these results. %In particular they have helped shaped my view that the importance of this work rests in the potential future applications and developments. 
I would also like to thank Bruce Bartlett, Dan Berwick-Evans, Nick Gurski, Andr\'e Henriques, Mike Hopkins, Niles Johnson, Aaron Lauda, Jacob Lurie, Ang\'elica Osorno, Hendryk Pfeiffer, Piotr Pstragowski,
%Piotr Pstr\k{a}gowski,  
Urs Schreiber,  Noah Snyder,  AJ Tolland, Jamie Vicary,  and Alan Weinstein for their useful conversations and suggestions at various stages of this work. Finally I would like to extend my gratitude to the anonymous referee who so carefully read this manuscript. He or She offered numerous important suggestions for improvements, of both a mathematical and presentational nature. We would also like to extend our gratitude to Michael Stay for allowing us the use of his diagrams in Appendix~\ref{app:defnsymbicat}.

 % Introduction

\mainmatter  %% Main Chapters & Appendicies

%\chapter{Planar Decompositions of 2-Manifolds} \label{ChapPlanarDecomp}
\chapter{Higher Morse Theory} \label{ChapPlanarDecomp}

In this chapter we prove the main geometric theorems used in the classification of 2-dimensional topological field theories. We formulate these results completely in the classical language of differential topology, and in particular make no mention of 2-dimensional bordism bicategory itself. In particular the results of this chapter are independent of the specifics of the construction of the 2-dimensional bordism bicategory, a task which is taken up in Chapter \ref{ClassificationChapter}. Jet transversality, Cerf theory, and singularity theory have had a relatively long and fruitful mathematical life, and the results of this chapter are heavily influenced by analogous results in these fields.  The techniques used here are essentially classical, building on well known and heavily developed existing results. 

With this in mind, we begin with a review of the classical literature. In Sections \ref{SectJetBundles} and \ref{SectJetTransversality} we review the classical notion of jet bundles and R. Thom's jet transversality theorems, including the Parametric Transversality Lemma \ref{ParametricTransversalityLemma}. These results are completely standard and may be found, for example,  in the text \cite{GG73}. 

The main results of this chapter closely parallel standard results in Cerf theory, in particular the Cerf theory of 1-manifolds. We exploit this analogy, and in Sections \ref{SectRudMorseTheory} and \ref{SectRudCerfTheory} we provide a revisionist review of Morse theory and Cerf theory. This allows us to explain the main results, as well as introduce the key methods of proving them. For simplicity, we restrict our attention to the case of Morse and Cerf theory of 1-manifolds, which is the most analogous to the results obtained later in this chapter.  While these Morse and Cerf theory results are completely classical, we don't know of a reference which presents this material from precisely this viewpoint while simultaneously providing details at the necessary level. However the discussion of singularity theory given in the text \cite{GG73} goes most of the way there. The proof techniques are essentially taken from there. 
 
At a very informal level, Morse theory is the study of generic functions from a manifold $M$ to $\R$. Cerf theory is then the study of generic {\em families} of functions from $M$ to $\R$. One such family is given by a path of functions. One of the main results of Cerf theory is that a generic path of functions is in fact a path of Morse functions for all but a finite number of isolated critical times. Moreover at these critical times, there exist standard normal coordinates, in analogy with the results of the Morse Lemma. One may then use these results to study precisely how to change any given Morse function on $M$ into any other Morse function. One of the main techniques is to regard a path of functions as a single map $M \times I \to \R \times I$. In this way, generic paths of functions are related to generic maps into $\R \times I$ or $\R^2$, and hence to the singularity theory of such maps. 

As  we will explain more precisely in the sections below, a map from $M$ to $\R$ gives a 1-dimensional/1-categorical decomposition of $M$. The goal of this chapter is to provide a categorification of this idea. We want a 2-dimensional/2-categorical decomposition of $M$, and so we are led to study maps from $M$ to $\R^2$, and generic paths of such maps. If our goal were merely to study generic maps to $\R^2$ and generic paths of such maps, then we could appeal to well established results in singularity theory. Indeed the classification of singularities in these dimensions is a well understood, by now classical subject. However, we intend to use these results to obtain new 2-categorical information. It is not enough merely to study the singularities of maps $M \to \R^2$, but we must simultaneously study the singularities of the composite map $M \to \R^2 \to \R$, where the last map is a standard projection. The singularities and ensuing decomposition theorems can then, in Chapter \ref{ClassificationChapter}, be matched with the bicategorical framework. 

The Cerf theory results that we are mimicking proceed by stratifying the jet space by a suitable collection of submanifolds, for example by the Thom-Boardman stratification, \cite[\oldS 5]{GG73}. Then one proves a theorem analogous to the Morse Lemma, providing normal coordinates for maps whose jet sections, $j^k f$, are transverse to this stratification, see Section \ref{SectRudCerfTheory}. Because we are simultaneously considering the singularities of maps $M \to \R^2$ and the projection  $M \to \R^2 \to \R$, we are lead to a refinement of the usual Thom-Boardman stratification, and consequently to a finer classification of the possible singularities. In Sections \ref{SectStratofPlaarDecomp}, \ref{SectStratofJetSpace2D},  \ref{SectNormalCoord2D},
\ref{SectStratofJetSpace3D}, and \ref{SectNormalCoord3D} we introduce the precise stratification of the jet spaces that we will use, and derive the corresponding normal coordinates.  This results in a finer classification of singularities than the standard classification. Although obtained independently, this approach is similar to that of Carter-Saito \cite{CS98} and Carter-Rieger-Saito \cite{CRS97}. Carter, Rieger, and Saito used these techniques to obtain a complete set of `movie moves' for three-dimensional diagrams of knotted surfaces in $\R^4$. It is highly likely that, combined with the results of Chapters \ref{SymMonBicatChapt} and \ref{ClassificationChapter}, one could use the results of Carter, Rieger, and Saito to classify extended field theories of two-dimensional bordisms embedded in $\R^4$. We will not pursue this here.

In Section \ref{SectGeomofSing}, we examine more closely the geometric structure of each of these singularities. Borrowing terminology from Cerf theory, we introduce the concept of a {\em graphic}. The graphic of a generic map is the image of the singular loci in the target, $\R^2$ or $\R^2 \times I$. These images consist of a collection of immersed surfaces, arcs, and points.  In Section \ref{SectMultijetConsiderations}, we use multi-jet techniques to introduce additional strata. A map whose mulit-jet section is additionally tranverse to these strata, which holds for a dense collection of maps, has the property that these immersed surface, arcs, and points are in general position. We use this as  the basis for defining an abstract graphic as a suitably labeled collection of surfaces, arcs, and points immersed in $\R^2$ or $\R^2 \times I$, and in general position. 

Finally, in Section \ref{SectPlanarDecompThm} we prove the main theorems of this section. We introduce the notion of planar diagrams and 3D diagrams. A planar diagram is roughly an abstract graphic, equipped with some additional data (``sheet data''). A generic map from a surface $\Sigma$ to $\R^2$, with a few additional choices, induces a planar diagram. The additional data of a planar diagram is essentially combinatorial in nature, nevertheless given a planar diagram one can in fact recover a surface with a generic map to $\R^2$ inducing the original planar diagram. Moreover this surface is unique up to isomorphism. A similar analysis holds for generic maps $\Sigma \times I \to \R^2 \times I$, and this allows us to introduce an equivalence relation on planar diagrams. Two planar diagrams are equivalent if they can be related by a finite number of ``local relations'', and in fact surfaces up to isomorphism are in bijection with equivalence classes of planar diagrams. These results are encapsulated in Theorem~\ref{PlanarDecompositionTheorem}.

\section{Jet Bundles and Jet Transversality}

\index{jet bundle}
Jets and jet bundles were introduced in 1951 by Ehresmann \cite{Ehresmann51-1, Ehresmann51-2, Ehresmann51-3} and over the nearly 60 years since their inception have become an increasingly important tool in differential geometry and the study of PDEs.  In singularity theory they have a particularly prominent role, and, when combined with Thom's jet transversality theorem, have led to a classification of singularities of small codimension. %\notetoself{Add references for importance in singularity theory and otehr areas of mathematics.} 

Roughly speaking, two functions $f,g: X \to Y$ have the same $r$-jet at a point $x$ if they agree \textquotedblleft up to order $r$\textquotedblright. There are several equivalent ways to make this precise. A common method is to choose local coordinates around $x$ and $f(x) = y = g(x)$. Then $f$ and $g$ agree up to order $r$ (at $x$) if all of their partial derivatives up to order $r$ agree. For example, $f$ and $g$ agree up to $1^\textrm{st}$ order at $x$ if $f(x) = g(x)$ and the matrix of single partial derivatives of $f$ and $g$ agree at $x$. This is precisely the condition
\begin{equation*}
	df_x = dg_x.
\end{equation*}
This condition is coordinate independent, and similarly the condition that $f$ and $g$ agree up to order $r$ is also coordinate independent, which is not immediate from the above description. 

There is an alternative approach to jets and jet bundles which is algebraic. This approach has the advantage that it is manifestly functorial and coordinate free. 
%In the following treatment of jets and jet bundles we will take a different approach. Jets and jet bundles can be described algebraically in a manifestly coordinate-free manner. This is particularly useful for defining jets and jet bundles for stratifolds. 
However, in the local coordinate approach it is straightforward to show that the jets form a smooth manifold and a bundle over $X \times Y$. For a particular chart this is clear: the jets consist of the (trivial) bundle of possible values of the partial derivatives. Then one must show how these glue together when a change of coordinates is preformed. This is more difficult to show in the algebraic approach. Both approaches have their respective advantages and disadvantages and we will employ both points of view in what follows. The interested reader should consult \cite{GG73, KMS93}.

Continuing the example of 1-jets, we have that the one-jets of maps from $X$ to $Y$ form the bundle
\begin{equation*}
	J^1(X, Y) = p_1^*(T^*X) \otimes p_2^*(TY) \to X \times Y,
\end{equation*}
where $p_i$ is projection onto the $i^\textrm{th}$ factor. In this case the jets form a vector bundle over the space $X \times Y$, but in general this will not be the case. Via the projection $X \times Y \to X$, we can view the jet bundle as a bundle over $X$. Then any smooth function $f:X \to Y$  gives us a smooth section,
\begin{equation*}
	df : X \to J^1(X,Y)
\end{equation*}
whose value at $x$ is precisely the differential $df_x$. This pattern persists. A map $f: X \to Y$ gives a section, $j^r f$, of the bundle of $r$-jets over $X$, $J^r(X, Y)$.

Classical transversality has several formulations, but one version concerns the density of certain maps between smooth manifolds. Given smooth manifolds, $X$ and $Y$, and a smooth submanifold $W \subseteq Y$, the transversality theorems assert that the maps $f: X \to Y$ which are transverse to $W$ are dense in the space of all maps, and that for such a map, $f^{-1}(W)$ is a smooth submanifold of $X$ of the expected codimension.

Thom  transversality (a.k.a. jet transversality) is a \index{jet transversality} 
\index{Thom transversality theorem} statement about the sections, $j^r f$, induced by maps $f: X \to Y$. Given a submanifold $W \subseteq J^r(X,Y)$, the jet transversality theorem says that for a dense subset of maps $f: X \to Y$, the sections $j^rf : X \to J^r(X,Y)$ are transverse to $W$.

% show that the stratifolds jets are stratifolds if the originals are regular

%Recall,
%\begin{definition}
%An open cover of a topological space $X$ is {\em locally finite} if every point of the space has a neighborhood which intersects only finitely many sets in the cover. A {\em refinement} of a cover of a space $X$ is a new cover of the same space such that every set in the new cover is a subset of some set in the old cover. A  topological space is {\em paracompact} if every open cover admits an open locally finite refinement.
%\end{definition}

\subsection{Jet Bundles} \label{SectJetBundles}

Let $X$ and $Y$ be smooth manifolds. Given a point $x \in X$, we can consider the algebra $C^\infty_x(X)$ of germs of functions at $x$. This is a local ring with unique maximal ideal $\m_x$ consisting of the germs of functions which vanish at $x$. A smooth map $X \to Y$ induces a map of local algebras $C^\infty(Y) \to C^\infty(X)$, and consequently for each integer $k \geq 0$ we have an induced map of rings 
\begin{equation*}
f^*_x: C^\infty_y / \m^{k+1}_y  \rightarrow C^\infty_x / \m^{k+1}_x.
\end{equation*}
Moreover, $f$ need not be globally defined. If $U \subseteq X$ is a neighborhood of $x$, and $V \subseteq Y$ is a neighborhood of $y$, then any smooth map $f: U \to V$, mapping $x$ to $y$ induces a map of rings, $f^*_x$, as above. The map $f^*_x$ only depends on the germ of $f$ near $x$. This construction is also functorial. If $f$ is a germ of a map sending $x\in X$ to $y \in Y$ and $g$ is a germ of a map sending $y \in Y$ to $z \in Z$, then the composition
\begin{equation*}
	C^\infty_z / \m^{k+1}_z   \stackrel{g^*_x}{\longrightarrow}  C^\infty_y / \m^{k+1}_y  \stackrel{f^*_x}{\longrightarrow} C^\infty_x / \m^{k+1}_x
\end{equation*}
coincides with the homomorphism $(g \circ f)^*_x$. 

Fix an integer $r \geq 0$. We impose an equivalence relation on triples $(x, U, f)$, where $x \in X$, $U \subseteq X$ is a neighborhood of $x$ and $f: U \to Y$ is a smooth map. We set $(x , U, f) \sim_r (x', U', f')$ if $x = x'$, $f(x) = f'(x)$ (whose common value we denote $y$), and $ f^*_x = f'^*_x \in \hom( C^\infty_y / \m^{r+1}_y , C^\infty_x / \m^{r+1}_x )$. If this is the case, we say that $f$ agrees with $f'$ at $x$ up to $r^\text{th}$ order. An equivalence class is \index{jet bundle} called a {\em jet} at $x$, and the set of equivalence classes is denoted $J^r(X, Y)$. For the equivalence relation we could just as well have replaced the pair $(U, f)$ with the germ of $f$ at $x$. However by not using germs we see that there are well defined quotient maps,
\begin{equation*}
	U \times C^\infty(U, Y) \to J^r(X,Y).
\end{equation*}
These quotient maps are compatible in the sense that whenever $V \subseteq U \subseteq X$,
 the following diagram commutes:
\begin{center}
\begin{tikzpicture}[thick]
	\node (LT) at (2,2) 	{$V \times C^\infty(U, Y)$};
	\node (LB) at (0,1) 	{$V \times C^\infty(V, Y)$};
	\node (RT) at (4,1) 	{$U \times C^\infty(U, Y)$};
	\node (RB) at (2,0)	{$J^r(X,Y)$};
	\draw [->] (LT) --  node [left] {$$} (LB);
	\draw [->] (LT) -- node [above] {$$} (RT);
	\draw [->] (RT) -- node [below right] {$q$} (RB);
	\draw [->] (LB) -- node [below left] {$q$} (RB);
\end{tikzpicture}
\end{center}

More generally for a smooth manifold, $S$, we may consider $S$-families of smooth maps $U \to Y$, by which we mean smooth maps $S \times U \to Y$. We have a sequence of maps of sets,
\begin{equation*}
	S \times U \times C^\infty(S \times U, Y) \stackrel{\text{ev}}{\to} U \times C^\infty(U, Y) \stackrel{q}{\to} J^r(X, Y),
\end{equation*}
which we may use to define a topology on the spaces $J^r(X, Y)$. By adjunction this gives a map of sets,
\begin{equation*}
	C^\infty(S \times U, Y) \to Map(S \times U, J^r(X, Y))
\end{equation*}
where $Map(S \times U, J^r(X, Y))$ denotes the set theoretic maps. Denote the image by $W_{S, U}$. We give $J^r(X, Y)$ the finest topology in which each of the maps in $W_{S, U}$ is continuous for each manifold $S$ and for each open $U \subseteq X$.

%-----------
%The continuous functions $C^0(U, Y)$ become a topological space when endowed with the compactly generated compact open topology and $C^\infty(U, Y)$ inherits the subspace topology. Thus $U \times C^\infty(U, Y)$ is a topological space, where we again use the compactly generated product. 
%We give $J^r(X,Y)$ the finest topology such that all the above quotient maps are continuous. 

In particular, this implies the there is a well defined, continuous map 
\begin{equation*}
J^k(X, Y) \to X,	
\end{equation*}
which we call the projection to $X$. Moreover, taking the case $U=X$, if $f: X \to Y$ is a smooth function, we get a continuous map,
\begin{equation*}
	X = X \times \{f\} \hookrightarrow X \times C^{\infty}(X,Y) \to J^k(X, Y).
\end{equation*}
We denote this map by $j^kf: X \to J^k(X,Y)$. By construction it is a section of the projection map.

It is a classical result  \cite{GG73, KMS93} 
%\notetoself{Actually in \cite{GG73} they prove that the quotient maps are continuous when we use the compact open topology and with some compactness conditions, but using some facts from the compactly generated topology the above follows (get reference). Add this remark as a footnote when you double check it.}
 that the jet space, as defined above, agrees as a topological space with the jet bundle constructed using local coordinates. Thus the jet space is in fact a finite dimensional smooth manifold and a fiber bundle over the product $X \times Y$. In fact, associated to the pair $X$ and $Y$ is a tower of smooth fiber bundles
\begin{equation*}
	X \leftarrow X \times Y = J^0(X,Y) \leftarrow J^1(X,Y) \leftarrow J^2(X,Y) \leftarrow J^3(X,Y) \leftarrow \cdots
\end{equation*}
 which we call the {\em jet bundles}. 
 In the local coordinates picture these successive maps can be viewed as the process of forgetting the highest order partial derivatives. In the algebraic picture it is the process of taking the quotient by a larger ideal. With this smooth structure, $j^k f$ becomes a smooth section of the projection map.

If we fix a point $x \in X$, we can consider the fiber ${}_xJ^r(X, Y)$, i.e those jets whose projection to $X$ is the point $x$.  Similarly, we can consider $J^r(X, Y)_y$ consisting of all jets whose projection to $Y$ has value $y$.  We can also consider the fiber over $(x, y) \in X \times Y$, denoted ${}_xJ^r(X,Y)_y = {}_xJ^r(X, Y) \cap J^r(X, Y)_y$. For a fixed $y\in Y$, the spaces $J^r(X, Y)_y$ form bundles over $X$ and the spaces ${}_xJ^r(X, Y)$ form bundles over $Y$, for fixed $x$.

Let us return once again to the case of 1-jets. We argued from the coordinate theoretic approach that the 1-jets can be identified with the vector bundle 
\begin{equation*}
	J^1(X, Y) = p_1^*(T^*X) \otimes p_2^*(TY) \to X \times Y.
\end{equation*}
This can also be seen from the algebraic approach as well. The inclusion of the (germs of) constant functions into $C^\infty_x$ allows us to canonically split the exact sequence,
\begin{equation*}
	\m_x/ \m_x^2 \to C^\infty_x / \m^2_x \to \R = C^\infty_x / \m_x.
\end{equation*}
If we have a vector space $V$, then $\R \oplus V$ becomes an algebra via the multiplication $(t_0, v_0) \cdot (t_1, v_1) = (t_0 t_1, t_0 v_1 + t_1 v_0)$ and $C^\infty_x / \m^2_x $ is isomorphic to such an algebra. In fact, the  above splitting gives us a canonical algebra isomorphism $C^\infty_x / \m^2_x \cong \R \oplus (\m_x / \m_x^2) \cong \R \oplus T^*_xX$. An easy calculation shows that algebra homomorphisms $\R \oplus V \to \R \oplus W$ are in bijection with linear maps $ V \to W$. Thus we find that the fibers ${}_xJ^1(X, Y)_y$, calculated from the algebraic perspective, are the algebra homomorphisms
\begin{equation*}
	\R \oplus T^*_yY \to \R \oplus T^*_xX,
\end{equation*}
and so consist of precisely the space $T^*_xX \otimes T_yY$, as expected.

One of the advantages of introducing the jet bundles is that they provide a means of combining both geometric and algebraic understanding. There are two common ways of defining the tangent bundle to a space: either as equivalence classes of certain curves $\gamma: \R \to Y$ or as the derivations of the algebra of germs $C_y^\infty$. Both of these perspectives are manifest from the point of view of jet bundles. The above discussion of 1-jets implies that we have the following isomorphism
\begin{equation*}
	{}_0J^1(\R, Y)_y \cong T^*_0\R \otimes T_yY \cong T_yY.
\end{equation*}
Moreover the left hand side can be viewed in two ways, corresponding to the two ways of defining the tangent bundle. On the one hand it is defined as an equivalence class of maps  $\gamma: \R \to Y$ (i.e., equivalence classes of curves); it is precisely those curves sending $0$ to $y$,  where two are considered equivalent if they agree at $y$ up to $1^\text{st}$ order. But as we have seen, this fiber is also the algebra homomorphisms,
\begin{equation*}
	C^\infty_y / \m^2_y \to C^\infty_0 / \m_0^2 \cong \R [ x] / (x^2).
\end{equation*}
Such a homomorphism induces the composite homomorphism $ C^\infty_y \to \R [ x] / (x^2)$. Let us consider these more closely. Such a homomorphism is given by the formula,
\begin{equation*}
	h: f \mapsto f(y) + \xi_h(f) \cdot x.
\end{equation*}
The fact that this is an algebra homomorphism implies that
\begin{align*}
	 (fg)(y) + \xi_h(fg) \cdot x &= h(fg)\\
	&= h(f) h(g) \\
	&= (f(y) + \xi_h(f) \cdot x) \cdot ( g(y) + \xi_h(g) \cdot x) \\
	&= f(y) g(y) + ( \xi_h(f) g(y) + f(y) \xi_h(g) ) \cdot x.
\end{align*}
Thus we find that $h$ is equivalent to the algebra derivation $\xi_h$. If $f \in \m_y$, (so that $f(y) = 0$), then we find that
\begin{equation*}
	h(f^2) = f^2(y) + 2 f(y) \xi_h(f) \cdot x = 0.
\end{equation*}
Thus any algebra homomorphism $h$ factors as 
\begin{equation*}
	h: C^\infty_y \to C^\infty_y / \m^2_y \to\R [ x] / (x^2).
\end{equation*}
This identifies the tangent bundle,\index{tangent bundle} i.e., the fiber ${}_0J^1(\R, Y)_y$, with the algebra derivations of $C^\infty_y$, precisely the second definition of the tangent space.  A similar discussion shows that the cotangent bundle is $J^1(X, \R)_0 \to X$. 

The higher order jet spaces ${}_0J^r(\R, Y)$ and $J^r( X, \R)_0$ should be thought of as higher order analogs of the tangent and cotangent bundles. The space $\R$ is a topological ring and so, by the functoriality of the jet bundles, the fibers $ {}_xJ^r( X, \R)_0$ become rings, in fact algebras over $\R$. These algebras coincide with the algebras $C_x^\infty / \m^{r+1}_x$, as can be readily checked, see \cite{KMS93} as well. 

There are several closely related bundles, such as $J^k(X^n, \R^n)_0$ and ${}_0J^k(\R^n, Y^n)$ but these are better understood using the frame bundle of order $r$. It is well known that one can associate a $GL_n$-principle bundle to the tangent bundle (or to the cotangent bundle) of a manifold $Y$ of dimension $n$. This is given by taking the bundle of $n$-frames whose fibers at a point $y$ consist of the sets of $n$-linearly independent vectors in $T_yY$. This fiber consists exactly of the linear isomorphisms
\begin{equation*}
	\R^n \cong T_0 \R^n \to T_yY.
\end{equation*}
Equivalently these are the jets in ${}_0J^1(\R^n, Y)_y$ which are represented by a {\em local diffeomorphism}. 

We define the {\em $k^\text{th}$-order frame bundle} \index{frame bundle, of order $k$} $P^k Y$ as the sub-bundle  of ${}_0J^k(\R^n, Y^n)$ whose fiber consists of those $k$-jets represented by {\em local diffeomorphisms}. This is a $G^k_n$-principal bundle over $Y$, where $G^k_n \subseteq {}_0J^k(\R^n, \R^n)_0$ is the {\em jet group} \index{jet group} of $k^\text{th}$-order jets represented by local diffeomorphisms. Let us recall some facts about the jet groups. We have the isomorphism $GL_m = G_m^1$, familiar from the frame bundle construction. In general the projection ${}_0 J^k( \R^n, \R^n)_0 \to {}_0 J^1(\R^n, \R^n)_0$ gives us a (split) exact sequence of Lie groups,
\begin{equation*}
	1 \to B_m^k \to G_m^k \to GL_m \to 1
\end{equation*}
The normal subgroup $B_m^k$ is connected, simply connected, and nilpotent. %$G^k_m \cong B_m^k \rtimes GL_m$. 
A more detailed discussion can be found in \cite{KMS93}.

The vector space structure on $\R^n$ induces a vector space structure on $L^k_{m,n} = {}_0 J^k (\R^m, \R^n)_0$ and $G_n^k$ naturally \textquotedblleft acts\textquotedblright on the right on this space. This action is generally {\em not} a linear action, but we can still form the associated bundle.  Given a manifold $Y$, we can take the $k^\text{th}$-order frame bundle of $Y$, which is a (left) $G^k_n$-principal bundle and form the bundle associated to the \textquotedblleft representation\textquotedblright  $L^k_{m,n}$. This is precisely the bundle ${}_0J^k(\R^m, Y)$. 

A similar discussion applies to coframes. The {\em $k^\text{th}$-order coframe bundle} \index{coframe bundle, of order $k$} $P^{k *}X$ of $X^m$ consists of those jets in $J^k(X^m, \R^m)_0$ which are represented by local diffeomorphisms. This is a right  $G^k_m$-principal bundle. $L_{m,n}^k$ has a natural left $G_m^k$-action and the corresponding associated bundle is $J^k( X, \R^n)_0$. Moreover the left action of $G^k_m$ and the right action of $G^k_n$ on $L_{m, n}^k$ are compatible. We can view this as a single right action by $G^k_m \times G^k_n$ in the usual way. The bundle $P^{k}X \times P^k(Y) \to X \times Y$ is a principal $G^k_m \times G_n^k$-bundle and hence we can form the bundle associated to $L_{m,n}^k$. This is the bundle $J^k(X^m, Y^n)$ over $X \times Y$, whose typical fiber $L_{m,n}^k$ has dimension,
\begin{equation*}
	n \left[ \left( \begin{array}{c}
	m + k \\ m
\end{array} \right)  - 1 \right].
\end{equation*}
In fact $L_{m,n}^k$ can be canonically identified with the space of $n$-tuples of polynomials in $m$-variables of degree $\leq k$, whose constant coefficient vanishes. Choosing sufficiently small coordinate patches $U$ and $V$ in both $X$ and $Y$, we can trivialize the jet bundle over $U \times V$. It can be identified with the trivial bundle $U \times V \times L^k_{m,n}$. These trivializations are then glued together by the action of the jet group.

In this work we will also need to consider manifolds with boundary and corners. 
Jets and jet bundles for manifolds with corners are slightly more delicate, but are also well established \cite{Michor80}. We will only need to consider the case when the source manifold, $X$, has corners but the target manifold, $Y$, does not. Moreover, we will only consider those manifolds with corners, $X^m$, which are equipped with a germ in a neighborhood of $X \subseteq \tilde X$. Here $\tilde X$ is an $m$-manifold without boundary, and $X$ is a submanifold with corners of $\tilde X$. In this case we may define the jet bundle $J^k(X,Y)$ as the pullback to $X$ of the bundle $J^k(\tilde X, Y)$ over $\tilde X$. 

Multi-jet bundles \index{jet bundle} \index{multi-jet bundle} will also play and important role, and are defined similarly. Let $X$ be a smooth manifold. Define $X^{(s)} \subseteq X^s$ to be $X^{(s)} = \{ (x_1, \dots, x_s) \in X^s \; | \; x_i \neq x_j \; \textrm{ for } 1 \leq i < j \leq s\}$. Let ${\bf k} = (k_1, \dots, k_s)$ be a multi-index of natural numbers and let $J^{\bf k}(X, Y) = J^{k_1}(X, Y) \times \cdots \times J^{k_s}(X, Y)$.  Let $p: J^k(X, Y) \to X$ be the projection map and define $p^{\bf{k}}: J^{\bf k}(X, Y) \to X^s$ to be the obvious projection map. Then the $s$-fold ${\bf k}$-jet bundle is defined to be $J^{({\bf k})}(X, Y) =  (p^{\bf k})^{-1}(X^{(s)})$. It is a smooth bundle over the smooth manifold $X^{(s)}$.  Now let $f: X \to Y$ be a smooth map. Then we define $j^{\bf k} f: X^{(s)} \to J^{({\bf k})}(X, Y)$ to be the section,
\begin{equation*}
	j^{\bf k} f (x_1, \dots, x_s) = ( j^{k_1} f( x_1), \dots, j^{k_s}f( x_s) ). 
\end{equation*}

%Finally, we will need one last variation on jet bundles: relative jet bundles. Fix an interval, $I$. We will want to consider those maps $f: X^m \times I \to Y^n \times I$ which preserve the projection to $I$ (the {\em relative maps}),  that is $f(x,t) = (f_t(x), t)$. These form a space $C^\infty_\text{rel}( X \times I, Y \times I)$. There is a corresponding relative jet space, 
%\begin{equation*}
%	J^k_\text{rel}(X \times I, Y \times I),
%\end{equation*}
%which is defined exactly as before, but only using those equivalence classes of functions which come from relative maps. This is a fiber bundle over $X \times I \times Y$, with fiber $L^k_{m,n}$. A relative map $f$ gives rise to a section $j^kf: X \times I \to J^k_\text{rel}(X \times I, Y \times I)$. 

\subsection{Jet Transversality} \label{SectJetTransversality}

The jet transversality theorems presented below are central for the results of this chapter. The essence of these theorems is that given any submanifold $W \subseteq J^k(X, Y)$, for a generic map $f: X \to Y$, the section $j^k f: X \to J^k(X,Y)$ is transverse to $W$. Moreover, if we have a countable collection of submanifolds $W_i \subseteq J^k(X, Y)$, we may require $j^k f$ to be transverse to each $W_i$ separately,  which will also hold for generic maps.  These simple observations can have dramatic consequences. By judiciously choosing the $W_i$, we may prove that generic maps can be completely described locally. We can show that  certain  singularities occur and that others do not. In essence we will be able to achieve a local classification of maps between certain spaces. Most of the following definitions and theorems can be found the classic texts \cite{GG73, Michor80}.

\begin{definition} \index{transversality}
	Let ${ X}$ and ${ Y}$ be smooth manifolds and $f: {X} \rightarrow {Y}$ a smooth map. Let ${W}$ be a submanifold of ${Y}$ and $x$ a point in ${X}$. Then $f$ {\em  intersects ${W}$ transversally at $x$} if either
	\begin{enumerate}
		\item $f(x) \notin {W}$, or
		\item $f(x) \in  {W}$ and $T_{f(x)}  {Y} = T_{f(x)}  {W} + (df)_x (T_x  {X})$.
	\end{enumerate}
	If $A \subseteq  {X}$, then $f$ {\em intersects ${W}$ transversally on $A$} (written: $f \pitchfork W$ on $A$) if $f$ intersects ${W}$ transversally at $x$ for all $x \in A$, and finally $f$ {\em  intersects ${W}$ transversally} if $f$ intersects ${W}$ transversally on ${X}$, in which case we denote this by $f \pitchfork {W}$.
\end{definition}

One of the key features of transversality is the following classical result.
\begin{theorem}
Let $X$ and $Y$ be smooth manifolds, and $W \subseteq Y$ a smooth submanifold. Let $f: X \to Y$ be a smooth map and assume that $f \pitchfork W$. Then $f^{-1}(W)$ is a submanifold of $X$ and $\codim f^{-1}(W) = \codim(W)$.
\end{theorem}

In order to speak about generic mappings, we must be able to talk about dense subsets of mappings, and in order to talk about dense subsets we must place a topology on the space of mappings. The standard topology on $C^\infty(X, Y)$ that we will use is the Whitney $C^\infty$-topology. \index{Whitney $C^\infty$-topology} Let $S$ be a subset of $J^k(X,Y)$. We define the following set:
\begin{align*}
 M(S) :=  &\{ f \in C^\infty({\sf X}, {\sf Y}) \; | \;  j^kf (X) \subseteq S \}
\end{align*}
Note that $\cap M(S_\alpha) = M (\cap S_\alpha)$, for an arbitrary collection of sets $S_\alpha$. The collection of $M(U)$ where $U$ ranges over all the open subsets of $J^k(X, Y)$ forms a basis for a topology on $C^\infty(X, Y) $, which we call the Whitney $C^k$-topology. The topology generated by the union of these topologies for all $k$ is the Whitney $C^\infty$-topology. 

A subset $R$ of a topological space $X$ is called {\em residual} \index{residual set} if it is the intersection of a countable collection of dense open sets. If $X$ is a Baire space, \index{Baire space} then every residual set is dense. With the Whitney $C^\infty$-topology, $C^\infty(X, Y)$ is a Baire space, and moreover the map,
\begin{equation*}
	j^k: C^\infty(X, Y) \to C^\infty(X, J^k(X, Y))
\end{equation*}
 is a continuous map.

\begin{theorem} \label{OpenonclosureTheorem}
Let $X$ and $Y$ be smooth manifolds with $W$ a submanifold of $Y$. Let $K \subseteq W$ be an open subset whose closure in $Y$ is contained in $W$, and let $T_K = \{ f \in C^\infty(X, Y) \; | \; f \pitchfork W \text{ on } \overline K\}$. Then $T_K$ is an open subset of $C^\infty(X, Y)$ (in the Whitney $C^\infty$-topology).
\end{theorem}

The essence of transversality is the following parametric transversality lemma. It says that if a $B$-family of maps from $X$ to $Y$ is transverse to $W$, then for generic values of the parameter, $b \in B$, the individual maps $X \to Y$ are transverse to $W$. From this all manner of other transversality theorems follow. 

\begin{lemma}[Parametric Transversality Lemma] \label{ParametricTransversalityLemma} \index{parametric transversality lemma}
Let $X$, $B$, and $J$ be smooth manifolds, with $W \subseteq J$ a smooth submanifold. Let $j: B \to C^{\infty}(X, J)$ be a map of sets and assume that $\Phi: X \times B \to J$, defined by $\Phi(x,b) = j(b)(x)$, is smooth and $\Phi \pitchfork W$. Then the set,
%\begin{equation*}
$	\{ b \in B \; | \; j(b) \pitchfork W \}$
%\end{equation*}
is dense in $B$.
\end{lemma}

This classical lemma essentially follows from Brown's Theorem (a corollary of Sard's Theorem), and details of the proof may be found in \cite[II \oldS 4 Lemma 4.6]{GG73}. This is the key lemma which makes the proof of the following theorem possible.  

\begin{theorem}[Thom Transversality Theorem] \label{ThomTransversalityTheorem} 
\index{Thom transversality theorem}
\index{transversality}
Let $X$ and $Y$ be smooth manifolds and $\{W_i\}$ a countable collection of submanifolds of $J^k(X, Y)$. Let 
\begin{equation*}
T_W = \{ f \in C^\infty(X,Y) \; |\; j^k f \pitchfork W_i, \; \forall i \}.
\end{equation*}
Then $T_W$ is a residual subset of $C^\infty(X,Y)$ in the Whitney $C^\infty$-topology. 
\end{theorem}

\begin{proof}[Proof Outline]
First we choose a countable covering $\{W_{i,r}\}$ of each of the submanifolds $W_i$, such that each of the $W_{i,r}$ is open in $W_i$, its closure is contained in $W_i$, and for which there exist suitably nice coordinates for $X$ and $Y$ around the image of $\overline W_{i,r}$ under the projection,
\begin{equation*}
J^k(X, Y) \to X \times Y.
\end{equation*}
 Then we consider the sets,
 \begin{equation*}
T_{W_{i,r}} = \{ f \in C^\infty(X,Y) \; |\; j^k f \pitchfork W_i \; \text {on } \overline W_{i,r} \}.
\end{equation*}
 The set $T_W$ is the intersection of the sets $T_{W_{i,r}}$, hence it is enough to prove that each is open and dense. 
 
 If $T_{i,r} = \{ g \in C^\infty(X, J^k(X,Y)) \; | \; g \pitchfork W_i \text{ on } \overline W_{i,r} \}$, then $T_{i,r}$ is open by Theorem~\ref{OpenonclosureTheorem}. Since $j^k: C^\infty(X, Y) \to C^\infty(X, J^k(X, Y))$ is continuous we have that $T_{W_{i,r}} = (j^k)^{-1}(T_{i,r})$ is also open. 
 
Density is more difficult to prove. Given any map $f: X \to Y$, we let $B = \R^n \times L^k_{m,n}$, the space of all polynomial mappings $\R^m \to \R^n$ of degree $k$ (not necessarily preserving the origin). Using our judicious choice of local coordinates in $X$ and $Y$, we construct a deformation $g_b$ of $f$, parametrized by $b \in B$. We define $\Phi: X \times B \to J^k(X, Y)$ by
\begin{equation*}
	\Phi(x, b) = j^k g_b(x).
\end{equation*}
Unfortunately, we cannot immediately apply Lemma \ref{ParametricTransversalityLemma}, because the map $\Phi$ may not be transverse to $W$. However, if we shrink $B$ to a suitable neighborhood of the origin, then $\Phi$ becomes a local diffeomorphism and hence satisfies any transversality condition.
\end{proof}

The reason for expounding on the proof of the Thom Transversality Theorem is that it is now clear that essentially identical theorems hold in much more generality. For example, the jet transversality theorem in the case where $X$ has corners readily follows from the theorem for its neighborhood, $\tilde X$. A relative version also holds. 

\begin{theorem}[Thom Transversality Theorem (with corners)] \label{ThomTransversalityTheoremCorners}
	\index{Thom transversality theorem}
Let $X$ be a smooth manifold with corners,  $X \subseteq \tilde X$ a smooth embedding into  a smooth manifold without corners, $Y$ a smooth manifold, and $\{W_i\}$ a countable collection of submanifolds (possibly with corners) of $J^k(X, Y)$. Let 
\begin{equation*}
T_W = \{ f \in C^\infty(X,Y) \; |\; j^k f \pitchfork W_i, \; \forall i \}.
\end{equation*}
Then $T_W$ is a residual subset of $C^\infty(X,Y)$ in the Whitney $C^\infty$-topology. 
\end{theorem}

\begin{theorem}[Relative Thom Transversality Theorem] \index{Thom transversality theorem} \label{ThomTransversalityTheoremRelativeCorners}
\index{Thom transversality theorem}
\index{transversality}
%Fix $0 < \epsilon < \delta < 1$. 
Let $X$ be a smooth manifold, possibly with corners, equipped with a (germ of) embedding $X \subseteq \tilde X$ into a smooth manifold without corners. Let $Y$ be a smooth manifold and let and $\{W_i\}$ a countable collection of submanifolds (possibly with corners) of $J^k(X, Y)$. 
%Let $f_0, f_1: X \to Y$ be smooth maps such that there exists a map $F: X \times I \to Y \times I$ which satisfies:
%\begin{enumerate}
%%\item $F(x, t) = (f_t, t)$, i.e., $F$ preserves the projection to $I$.
%\item There exist neighborhoods $N_0$ of $0 \in I$, and $N_1$ of $1 \in I$, such that %
%%%%%%%$\overline M_0 \subseteq N_0$, $\overline M_1 \subseteq N_1$ and 
%$F$ satisfies
%%%%%%$N_0 \subseteq I$ of $ [0, \epsilon]$ and $N_1 \subseteq I$ of $ [\delta, 1]$ such that
%  $f_t(x) = f_0(x)$ for $ t \in N_0$ and $f_t(x) = f_1(x)$ for $ t \in N_1$. 
%\item $F \pitchfork W_i$ on a neighborhood of $X \times \{0, 1\}$. 
%\end{enumerate}
%Note that if there exists such an $F$, then property (2) holds for any $F$ satisfying property (1). Let $ C^\infty_\text{rel}(X \times I,Y\times I)$ denote the subspace of those maps $F$ satisfying property (1), above. Let
%\begin{equation*}
%T_W = \{ F \in C^\infty_\text{rel}(X \times I,Y\times I) \; |\;  j^k f \pitchfork W_i, \; \forall i \}.
%\end{equation*}
%Then $T_W$ is a residual subset of $C^\infty_\text{rel}(X \times I,Y \times I)$ in the Whitney $C^\infty$-topology. 
Let $X_0 \subseteq X$ be a closed subspace and let $f_0: X \to Y$ be a smooth map such that there exists and open neighborhood $N \subseteq X$ of the subspace $X_0$ so that for each $i$ we have $j^k f_0 \pitchfork W_i$ on $N$. 

Let $ C^\infty_\text{rel}(X,Y)$ denote the subspace of those maps $f$ such that there exists an open neighborhood $U \subseteq X$ of the subspace $X_0$, depending on $f$, so that $f|_U = f_0|_U$. Let
\begin{equation*}
T_W = \{ f \in C^\infty_\text{rel}(X,Y) \; |\;  j^k f \pitchfork W_i, \; \forall i \}.
\end{equation*}
Then $T_W$ is a residual subset of $C^\infty_\text{rel}(X ,Y )$ in the Whitney $C^\infty$-topology.
\end{theorem}

There is one final version of jet transversality we will need, which also generalizes to both the case of manifolds with corners and the relative case.
\begin{theorem}[Multi-Jet Transversality Theorem] \label{MultiJetTransversalityThm} 
\index{Thom transversality theorem} 
\index{multi-jet transversality theorem}
\index{transversality}
 Let ${\bf k} = (k_1, \dots, k_s)$ be a multi-index of natural numbers. Let $X$ and $Y$ be smooth manifolds with $\{ W_i \}$ a countable collection of submanifolds of $J^{({\bf k})}(X, Y)$. Let 
\begin{equation*}
T_W = \{ f \in C^\infty(X,Y) \; |\; j^{\bf k}f \pitchfork W_i, \; \forall i\}.
\end{equation*}
Then $T_W$ is a residual subset of $C^\infty(X, Y)$.
\end{theorem}

\section{Stratifications of the Jet Space}

In this section we will begin to use the results of the previous section to classify generic maps between certain manifolds $X$ and $Y$. Roughly what we will do is to choose a stratification of the jet space
\begin{equation*}
	J^k(X,Y) = \cup S_i,
\end{equation*}
then, by the Thom Transversalty Theorem~\ref{ThomTransversalityTheorem}, we are assured that a dense collection of maps will be simultaneously transverse to each of the strata $S_i$. We call such maps {\em generic}. By choosing the strata appropriately, we will be able to derive {normal coordinates} for these generic mappings, which will ultimately lead to the main result of this chapter, a planar decomposition theorem for surfaces. The techniques we will employ are essentially imported from Cerf theory, and so it is natural that we include a brief review of those aspects we will later generalize.

\subsection{Rudimentary Morse Theory} \label{SectRudMorseTheory}

Morse theory 
\index{Morse theory|(}
is a cornerstone of differential topology, and it would be difficult to improve upon the many wonderful expositions which cover it. Nevertheless, it will be useful to highlight what the previous sections on jet bundles and transversality have to offer. 

In Morse theory one studies maps from a manifold $X$ to $\R$, and so naturally we may consider the jet space $J^1(X, \R)$. We can stratify $J^1(X, \R)$ into two strata $S_0$ and $S_1$ based on the {\em corank}
\index{corank}
%\index{$S_0$, $S_1$} 
of the differential $df$. Thus $S_0$ corresponds to those jets where $df$ has rank one, and $S_1$ to those for which $df$ vanishes. $S_0$ is open and dense in $J^1(X, \R)$ and $S_1$ has codimension equal to the dimension of $X$. 

By Theorem~\ref{ThomTransversalityTheorem}, we know that for a dense subset of maps $f: X \to \R$, $j^1f$ is transverse to both $S_0$ and $S_1$. This implies that for such generic $f$, the differential $df$ vanishes at only isolated points in $X$ (the critical points).  We can obtain more refined information by stratifying the higher jet bundles. 

Let $x$ be a point of $X$ and $f: X \to Y$. Then $j^kf: X \to J^k(X, Y)$, and we can  consider its differential at $x$, i.e., the map
\begin{equation*}
	(d j^kf)_x : T_xX \to T_{j^kf(x)}J^k(X, Y).  
\end{equation*}
This map is determined by the $(k+1)$-jet of $f$ in the following sense, (see \cite{GG73}). Let $\psi = j^{k+1}f(x)$ be the $(k+1)$-jet of $f$ at $x$. Let $g$ be any other (germ of a) map such that $j^{k+1}g(x) = \psi$, then
\begin{equation*}
	(d j^k f)_x = (d j^k g)_x,
\end{equation*}
which may be checked in local coordinates. 
Hence we may reconstruct the kernel of this map simply by knowing the $(k+1)$-jet.  

Returning to the study of maps from $X$ to $\R$, let us look at the jet space $J^2(X, \R)$. The inverse images of $S_0$ and $S_1$ under the bundle projection,
\begin{equation*}
	J^2(X,Y) \to J^1(X,Y)
\end{equation*}
induce a preliminary stratification of $J^2(X,Y)$ (we will call these strata $S^{(2)}_0$ and $S^{(2)}_1$, respectively). $S^{(2)}_0$ is a perfectly fine stratum, however we may refine $S^{(2)}_1$ to obtain more useful information. 
%\index{$S^{(2)}_r$}

Recall that the fibers ${}_x J^1(X, \R)_y$ may be identified with $\hom(T_xX, T_y \R)$. Under this identification, $S_1 \cap {}_x J^1(X, \R)_y = \{ 0 \} \subseteq \hom(T_xX, T_y \R)$. Thus, the normal bundle $\nu S_1$ of $S_1$ in $J^1(X, \R)$ is canonically isomorphic to $\hom(TX, T\R)$, where by $TX$ and $T\R$ we mean their corresponding pullbacks to $S_1$. Thus by our previous discussion, we get a map of fiber bundles over $S_1$,
\begin{align*}
	S^{(2)}_1 \to \hom(TX, TJ^1(X,Y)) \to  &\hom( TX, \nu S_1) \\
	 &\cong  \hom( TX, \hom(TX, T\R)) \\
	 &\cong \hom (TX \otimes TX, T\R)
\end{align*}
where all the relevant vector bundles denote their corresponding pullbacks to $S_1$. The composition of these maps associates to $\sigma \in S^{(2)}_1$ the Hessian of $f$ at $x$, where $f$ is any map germ such that $j^2f(x) = \sigma$ (and this doesn't depend on the representative $f$). Thus the above composition factors as a map of fiber bundles over $S_1$,
\begin{equation*}
	S^{(2)}_1 \to \hom( TX \odot TX, T\R)
\end{equation*}
where $\odot$ denotes the symmetric tensor product. 
\index{$\odot$, symmetric tensor product}
An easy calculation shows that this map is a surjective submersion. We obtain a finer stratification of $S^{(2)}_1$ by stratifying $ \hom( TX \odot TX, T\R)$ and pulling it back by the above map.

Let $V$ and $W$ be vector spaces, and let $L^r(V, W)$ 
\index{$L^r(V, W)$}
denote the subset of linear maps $V \to W$ which ``drop rank by $r$''. That is, if $R$ denotes the maximal possible rank for any map $V \to W$, then 
 $L^r(V, W)$ consists of those maps which have rank $R -r$. This construction is clearly natural, and hence applies to vector bundles as well.  A straightforward computation shows that the pullback (of fiber bundles over $S_1$),
 \begin{center}
\begin{tikzpicture}[thick]
	\node (LT) at (0,1.5) 	{$\hom(TX \odot TX, T\R)_r$ };
	\node (LB) at (0,0) 	{$\hom(TX \odot TX, T\R)$};
	\node (RT) at (6,1.5) 	{$L^r(TX, \hom(TX, TR))$};
	\node (RB) at (6,0)	{$\hom(TX, \hom(TX, T\R))$};
	\node at (.5,1) {$\ulcorner$};
	\draw [->] (LT) --  node [left] {$$} (LB);
	\draw [->] (LT) -- node [above] {$$} (RT);
	\draw [->] (RT) -- node [right] {$$} (RB);
	\draw [->] (LB) -- node [below] {$$} (RB);
\end{tikzpicture}
\end{center}
is again a smooth fiber bundle. This gives us a stratification of $\hom(TX \odot TX, T\R)$ by the submanifolds $\hom(TX \odot TX, T\R)_r$, and hence a corresponding stratification of $S^{(2)}_1$ 
%\index{$S_{1, r}$}
by submanifolds $S_{1, r}$. The codimension of $S_{1,0}$ (which we call the {\em non-degenerate stratum}) is the dimension of $X$. The other strata $S_{1,r}$ (which we call {\em degenerate}) have higher codimension.  If we wish, we may further divide $S_{1,0}$ into components according to the signature of the non-degenerate Hessian.

Now Theorem~\ref{ThomTransversalityTheorem} implies that for generic $f$, not only are the only critical points isolated points in $X$, but that each of these critical points is non-degenerate.\footnote{ Due to the relative dimension of $X$ and the codimension  of the degenerate strata, transversality of $j^2f$ and $S_{1,r}$ implies that $j^2f(X) \cap S_{1,r} = \emptyset$ for all $r >0$.} Such a function is called a Morse function. The classical Morse Lemma 
\index{Morse lemma}
states that for a Morse function, $f$, in a neighborhood of each critical point there exist normal coordinates where $f$ takes the following form:
\begin{equation} \label{MorseLemmaCoordEqn}
	f(x_1, x_2, \dots, x_m) = c + x_1^2 + x_2^2 + \cdots + x_j^2 - x_{j+1}^2 - \dots - x_m^2,
\end{equation}
where $c = f(0)$ is the value of the critical point, aptly called the {\em critical value}. \index{critical value}

We will want to use this local description of $f$ to decompose $X$ into well understood elementary pieces. One problem that we would run into at this point is that it might be possible for multiple critical points to have the same critical value. We can fix this by using the multi-jet transversality theorem. The multi-jet space has a submersion,
\begin{equation*}
	\pi: J^{(1,1)}(X, \R) \to \R \times \R,
\end{equation*}
so that we may consider the submanifold,
\begin{equation*}
	S := (S_1 \times S_1) \cap J^{(1,1)}(X, \R) \cap \pi^{-1}( \Delta \R).
\end{equation*}
where $\Delta \R \subseteq\R \times \R$ is the diagonal. 
This submanifold has codimension $2 \cdot \dim X + 1$. Thus the multi-jet transversality theorem ensures that for generic $f$, $j^{(1,1)} f (X^{(2)}) \cap S = \emptyset$. Thus any two critical points of a generic map must have distinct critical values.

\begin{figure}[ht]
\begin{center}
\begin{tikzpicture}[thick]
	\draw (0,0) -- (0,1) arc (180: 0: 0.5cm) -- (1,0);
	\draw (2, 1.5) ellipse (0.5 cm and 1cm);
	\draw (1.5, 4.5) -- (1.5, 3.5) arc (180: 360: 0.5cm) arc (180: 0: 0.5cm) -- (3.5, 0);
	\draw [thin, blue] (-0.5, 0.75) -- +(4.5,0) +(8.4, 0) -- +(8.6, 0);
	\draw [thin, blue] (-0.5, 1.75) -- +(4.5,0) +(8.4, 0) -- +(8.6, 0);
	\draw [thin, blue] (-0.5, 2.75) -- +(4.5,0) +(8.4, 0) -- +(8.6, 0);
	\draw [thin, blue] (-0.5, 3.5) -- +(4.5,0) +(8.4, 0) -- +(8.6, 0);

	\draw[->] (5, 2) -- node [above] {Morse Function} +(2, 0);
	
	\begin{scope}[xshift=-1cm]
		\draw (9, 4.5) -- +(0, -4.5);

		\node (A) [circle, fill=red,inner sep=1.5pt] at (9, 0.5) {};
		\node [circle, fill=red,inner sep=1.5pt] at (9, 1.5) {};
		\node [circle, fill=red,inner sep=1.5pt] at (9, 2.5) {};
		\node [circle, fill=red,inner sep=1.5pt] at (9, 3) {};
		\node [circle, fill=red,inner sep=1.5pt] at (9, 4) {};

		\node (B) at (11, 0.5)  {critical value};	
		\draw [->, thin] (B) -- (A);
	\end{scope}

	\draw (-1, -2) arc (180: 0: 0.5cm);
	\draw (0.5, -1) arc (180: 360: 0.5cm);
	\node (C) at (6, -1.5) {Elementary 1-Dimensional Pieces};
	\draw [->] (C) -- (2, -1.5);
\end{tikzpicture}
\caption{A One-Dimensional Decomposition (of a 1-Manfold) Induced by a Morse Function}
\label{Fig1DDecompFromMorse}
\end{center}
\end{figure}
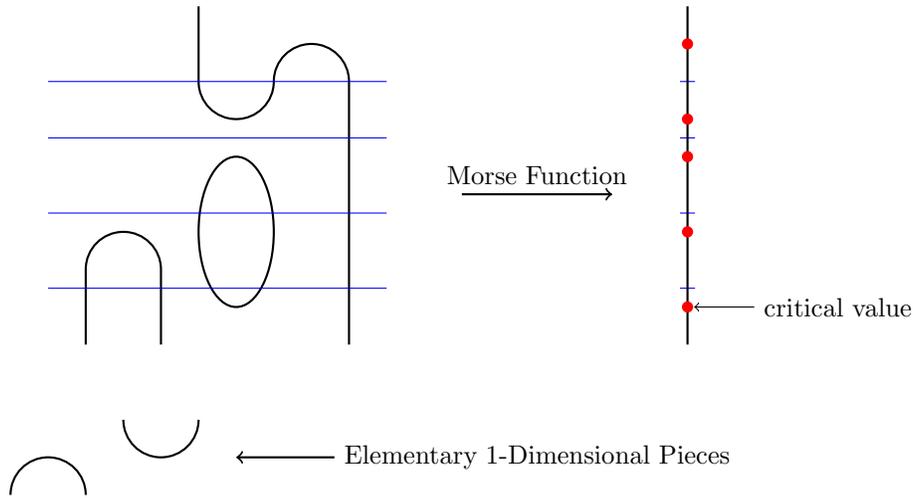

Choosing a generic map to the one-dimensional space $\R$ gives us a template for how to decompose the space $X$ in a one-dimensional fashion. We can decompose $X$ into regions where $f$ contains a single isolated critical point and into regions where $f$ has no critical points at all, see Figure~\ref{Fig1DDecompFromMorse}.
To completely understand this decomposition of $X$ one must further analyze these pieces. In particular, if $\dim X > 1$,  to understand the pieces with isolated critical points, one must choose a generic metric and gradient-like vector field. Only then can one obtain the usual handle structure associated to a Morse function. Moreover, all these pieces are glued together via diffeomorphisms of their respective boundaries, and these diffeomorphisms play an important role in the decomposition. 

Instead, we can specialize to the case $\dim X = m = 1$. This simplifies the analysis in two important ways. First of all, metrics and gradient like-vector fields are unnecessary; they give no new information. Secondly, the relevant gluing morphisms are diffeomorphisms of zero-dimensional manifolds and hence are essentially combinatorial.  Moreover, the function $f$, restricted to the regions which contain no critical points, is a local diffeomorphism. 
\index{Morse theory|)}

\subsection{Rudimentary Cerf Theory} \label{SectRudCerfTheory}

Cerf theory,\index{Cerf theory}\index{parametrized Morse theory}\index{parametrized Morse theory|see{Cerf theory}}\index{Morse theory}\index{Morse theory|see{Cerf theory}}
also called parametrized Morse theory, in its most elementary form is the study of families of Morse functions. There are several ways to motivate the study of such families, but perhaps the simplest is to answer the following question: Given a (compact) manifold $X$, we know that there exist Morse functions on $X$ and that these induce decompositions of $X$, but how are these Morse functions related? Is the space of Morse functions connected? How are the corresponding decompositions related?

The answer to this second question is no, the space of Morse functions is not connected. However, the same techniques which lead to our understanding of Morse functions can be used to answer the remaining questions. We know that the space of all functions $X \to \R$ is connected (in fact, contractible), thus between any two Morse functions we may choose a path of functions. Such a path of functions is the same as a single map $F: X \times I \to \R \times I$, which preserves the projections to $I$, i.e., we have $F(x, t) = ( f_t(x), t)$. By the relative transversality theorem we can deform $F$ to be a generic function keeping it the same on the ends of $I$. 

Moreover, as we will see, this can be done is such a way that the component $\partial_t pF$ never vanishes, (here $p: \R \times I \to I$ is the projection). So while this deformed $F$ is not precisely a path, viz. it does not commute with the projection to $I$, it is quite close and sufficient for our purposes.\footnote{In classical Cerf theory it is proven that we may deform $F$ to be generic, while simultaneously keeping it a path. This is done by stratifying the mapping space $C^\infty(X, Y)$ suitably (although it is not a manifold) and proving that any path between Morse functions may be deformed so as to only pass through the codimension-zero and codimension-one strata.} In fact we can pre-compose $F$ with a diffeomorphism of $X \times I$ to obtain a generic path. We now need only specify  a suitable collection of submanifolds of the jet space, with which to make $F$ generic.

What we will find is that such submanifolds can be chosen so that $f_t$ is a Morse function on the complement of a finite number $t \in I$. At these critical times we can again obtain local coordinates, which allow us to completely understand how the corresponding Morse functions change. When $X$ is one-dimensional this allows us to specify a finite list of ``moves'' or ``relations'' between Morse decompositions. Any two decompositions are related by a finite sequence of these moves.

The stratification of the jet space we will use is essentially the Thom-Boardman stratification,
\index{Thom-Boardman stratification}
 though we use the presentation due to Porteous \cite{Porteous71, Porteous72}, see also \cite{GG73}. Some strata will also coincide with the Morin strata. Recall that for any manifolds $M, N$, we have the identification,
\begin{equation*}
{}_{x}J^1(M, N)_{y} \cong \hom(T_x M, T_yN).
\end{equation*}
Recall also that for vector spaces $V$ and $W$,  $L^r(V, W)$ is a (natural) submanifold of $\hom(V, W)$ of dimension $r^2 + r \cdot | \dim V - \dim W|$. It consists of linear maps of corank $r$. The naturality of this submanifold is such that we can construct canonical submanifolds 
%\index{$S_r$}
$S_r \subseteq J^1(M, N)$ such that the fiber over $(x, y)$ is  
\begin{equation*}
S_r := L^r( T_x M, T_yN).
\end{equation*}

Given a linear map $A \in L^r(V, W)$, we can consider the kernel $K_A$ and cokernel $L_A$ which fit into an exact sequence of vector spaces,
\begin{equation*}
0 \to K_A \to V \stackrel{A}{\to} W \to L_A \to 0.
\end{equation*}
Over $L^r(V, W)$ the $K_A$ and $L_A$ assemble into (canonical) vector bundles and the normal bundle of $L^r(V, W)$ in $\hom(V,W)$ at $A$ is canonically isomorphic to the space $\hom(K_A, L_A)$, as shown in \cite[VI  $\oldS$ 1]{GG73}. This induces corresponding vector bundles $K$ and $L$ over $S_r$. As before, we may obtain a preliminary stratification of $J^2(M,N)$ by pulling back the stratification of $J^1(M, N)$. Thus we obtain submanifolds $S^{(2)}_r$. 
%\index{$S^{(2)}_r$}
As before, we have an induced map of fiber bundles over $S_r$,
\begin{equation*}
	S^{(2)}_r \to \hom( K \odot K, L)
\end{equation*}
which is a submersion. We look at the inclusion
\begin{equation*}
	\hom(K \odot K, L) \to \hom(K, \hom(K,L)).
\end{equation*}
  Now $ \hom(K, \hom(K,L))$ can be stratified into manifolds $L^s(K, \hom(K,L))$ and it can be shown that the intersection with $\hom(K \odot K, L) $ is also a manifold $\hom(K \odot K, L)_s$, see \cite{GG73}. Thus we may take the inverse image of $\hom(K \odot K, L)_s$ in $S^{(2)}_r$ to obtain a submanifold $S_{r,s}$. 
%\index{$S_{r,s}$}
These stratify the 2-jets. 

Let us specialize to the case relevant for Morse theory, $M = X \times I$ and $N = \R \times I$. In this case there are only the three strata $S_0$, $S_1$, and $S_2$ of codimensions $0$, $\dim X$, and $2 \cdot \dim X + 2$, respectively. Transversality to $S_2$, thus means that there are no points where the jet of $F$ is in $S_2$. For points in $S_1$, both $K$ and $L$ are one-dimensional, so that there are only two cases: $S_{1, 1}$ which consists of those elements of $S_1^{(2)}$ lying over the zero section of $\hom ( K \odot K, L)$, and its complement $S_{1,0}$. The normal bundle to $S_{1,1}$  in $S^{(2)}_1$
is thus canonically isomorphic to $\hom(K \odot K, L)$, so that $S_{1,1}$ is codimension one in $S_1^{(2)}$, i.e., codimension $\dim X +1$ in $J^2(X \times I, \R \times I)$. 

We may continue this process further. Let $S^{(3)}_{1,1}$
%\index{$S^{(3)}_{1,1}$}
 denote the inverse image of $S_{1,1}$ under the projection $J^3(X \times I, \R \times I) \to J^2(X \times I, \R \times I)$.  As before, we get a map of fiber bundles over $S_{1,1}$,
\begin{equation*}
	S^{(3)}_{1,1} \to \hom(TX, TJ^2(X \times I, \R \times I))
\end{equation*}
which we can project and restrict to get a map to $\hom( K, \hom(K \odot K, L))$. As we will see shortly, this actually lands in $\hom(K \odot K \odot K, L)$, and the composite,
\begin{equation*}
	S^{(3)}_{1,1} \to \hom(K \odot K \odot K, L),
\end{equation*}
is a surjective map of fiber bundles over $S_{1,1}$. Again we may stratify $S^{(3)}_{1,1} $ by stratifying $\hom(K \odot K \odot K, L)$ and pulling the stratification back to $S^{(3)}_{1,1} $. Since all the bundles involved are 1-dimensional, there is an obvious choice for such a stratification. Let  $\hom(K \odot K \odot K, L)_1$ denote the zero section and let $ \hom(K \odot K \odot K, L)_0$ denote its (open) complement. Thus we obtain two strata $S_{1,1,0}$ and $S_{1,1,1}$. The latter has codimension one in $S^{(3)}_{1,1}$ and hence codimension $\dim X +2$ in $J^3(X \times I, \R \times I)$. 

If a map $f: X \times I \to \R \times I$ is such that its jet at $(x, t)$ lands in one of the manifolds we have just constructed, $S_\alpha$, then we will say that $f$ has a {\em singularity of type $S_\alpha$ at $x$}. 
\index{singularity type}
The Thom transversality theorem now applies and shows us that the set $T_W$ of maps $f: X \times I \to \R \times I$ whose induced sections $j^kf$ are transverse to $S_0, S_1, S_2, S_{1,0}, S_{1,1}, S_{1,1,0}$ and $S_{1,1,1}$ form a residual (hence dense) subset. In particular, because the  codimension of $S_{1,1,1}$ is too high, these kinds of singularities do not occur. Moreover, the points with singularity type $S_1$ form a 1-dimensional submanifold of $X \times I$, of which there is a 0-dimensional submanifold of $S_{1,1}$ singularities. 

We may improve upon this in two ways. First, the submanifolds we have constructed here are compatible with those we used to define Morse functions. Thus if $f_0, f_1: X \to \R$ are two Morse functions and $F: X \times I \to \R \times I$ is any map satisfying condition (1) of the  ~\ref{ThomTransversalityTheoremRelativeCorners}, then $F$ automatically satisfies condition (2), as well. Thus Theorem~\ref{ThomTransversalityTheoremRelativeCorners} applies to this relative case. 

Secondly, an arbitrary (generic) map $X \times I \to  \R \times I$ is very far from being a path of functions. We can make it closer as follows. Consider the projection of jet bundles,
\begin{equation*}
	J^1(X \times I, \R \times I) \to J^1(X \times I, I).
\end{equation*}
Inside  $J^1(X \times I, I)$ there is a distinguished open subset $\cO$, whose inverse image in $J^1(X \times I, \R \times I)$ is also open (we will abuse notation and call this open set $\cO$, as well). On $X \times I$ we have a distinguished vector field, which in local coordinates is given by $\partial_t$, where $\partial_t$ is also the standard vector field on $I$.   The span of this gives us 1-dimensional  subbundle $E \subseteq T(X \times I)$, and thus we may consider the projection,
\begin{equation*}
	J^1(X \times I, I) = T^*(X \times I) \otimes TI \to E^* \otimes TI.
\end{equation*}
The inverse image of the complement of the zero section of $ E^* \otimes TI$ gives the open set $\cO$. Finally, we may consider the subset of $C^\infty(X \times I, \R \times I)$ given by,
\begin{equation*}
	M(\cO) = \{ F \; | \; j^1F(X \times I) \subseteq \cO \}.
\end{equation*}
This is open in the Whitney $C^1$-topology, and hence open in the Whitney $C^\infty$-topology. It consists of those maps $F(x,t) = (f_t(x), g(x,t))$ such that the $t$-derivative of the last coordinate (i.e., $\partial_t g$) is non-zero at all points of $X \times I$. A path of functions can be viewed as a map $F: X \times I \to \R \times I$, and such a  map is in $M(\cO)$, so that $M(\cO)$ is non-empty. Thus we may take our previously constructed residual set $T_W$ and intersect it with $M(\cO)$. A dense subset of maps in $M(\cO)$ will also lie in $T_W$.

To complete our analysis we must derive local coordinates for these singularities and study how these local descriptions combine to tell us about the whole. As before, we will simplify the discussion by considering just the case when $X$ is 1-dimensional. In the process we will also justify our claim that the map
\begin{equation*}
	S^{(3)}_{1,1} \to \hom(K, \hom( K \odot K, L)),
\end{equation*}
lands in $\hom(K \odot K \odot K, L)$, and is surjective onto this subspace. 

Up to this point, we have primarily been using the geometric description of jet bundles. To accomplish our further analysis we will now make contact with the algebraic description of jet space. Recall that the fibers of the jet bundles can be identified as certain spaces of algebra homomorphisms,
\begin{equation*}
	{}_x J^k(X,Y)_y \cong \hom( C^\infty_y / \m^{k+1}_y , C^\infty_x / \m^{k+1}_x ).
\end{equation*}
This algebraic description will allow us to use algebraic techniques to derive normal coordinates around the various singularities we have introduced. 

\begin{definition}
	Let $f: (X, x) \to (Y,y)$ be the germ of a local map taking $x\in X$ to $y\in Y$.  The {\em local ring} of $f$ is the quotient ring: $\cR_f =  C^\infty_x / C^\infty_x \cdot f^* \m_y$. This is a local ring with maximal ideal $\m_f$. If $f: X \to Y$ is an actual function, it induces a map germ at each point $x \in X$, and hence a corresponding local ring, $\cR_f(x)$.  	
\end{definition}

The strata of the jet space that we have introduced can alternatively be defined in terms of the corresponding local rings. 
\index{local ring}
For example, consider those map germs between  $X$ and $Y$ in which the the local ring is $\R$. These are precisely those map germs in which $f^* \m_y = \m_x$, and hence $df: T_xX \cong (\m_x / \m_x^2)^* \to (\m_y / \m_y^2)^* \cong T_yY$ is an isomorphism. These are precisely the map germs whose 1-jet lies in $S_0$, i.e., $S_0$ consists of precisely those jets which have local representatives whose local ring is $\cR_f \cong \R$. There are similar characterizations for the remaining strata, and in particular having a map germ whose jet is in the strata $S_{1,0}$ or $S_{1,1,0}$ in fact determines the local ring.

The $S_r$ strata 
%\index{$S_r$}
have a straightforward description in terms of the local ring. As we have seen, the 1-jets consist of the linear maps:
\begin{equation*}
	(df)^*: \m_y/ \m^2_y \to \m_x/ \m^2_x.
\end{equation*}
The stratification into sets $S_r$ corresponds to  precisely those where the cokernel (or kernel) of this map has dimension $r$. Thus the $S_0$ jets are those where $(df^*)$ is an isomorphism (and hence $\cR_f \cong \R$). The cokernel of this linear map is the vector space,
\begin{equation*}
	\tilde K := \m_x/ (C^\infty \cdot f^* \m_y + \m_x^2) \cong \m_f / \m_f^2,
\end{equation*}
and hence is entirely determined by the local ring $\cR_f$.

Consider the case were we have a local map germ $f$, whose 1-jet lies in $S_1$. We wish to determine whether its 2-jet lies in $S_{1,1}$ or $S_{1,0}$. From our earlier discussion we know that we should be able to extract a linear map,
\begin{equation*}
	K \odot K \to L,
\end{equation*}
and that if $s$ is the dimension of the cokernel (or kernel) of this map, then the 2-jet is in  $S_{1,s}$. 
%\index{$S_{1,r}$}
What we can construct from the local ring of $f$ is a linear map $L^* \to K^* \odot K^*$ which is essentially equivalent to the above map. Since $\tilde K$ is the cokernel of $(df^*): T^*_yY \to T^*_xX$, we have $\tilde K \cong K^*$. Note that the inverse image of $\m_x$ under the map $f^*: C^\infty_y \to C^\infty_x$ is precisely the ideal $\m_y$. Let $\tilde L$ be defined by
\begin{equation*}
	\tilde L := (f^*)^{-1}( \m_x^2) / \m_y^2.
\end{equation*}
Thus $\tilde L$ is precisely the kernel of $(df)^*$, and hence $\tilde L \cong L^*$. 

\begin{proposition}
	The multiplication map induces an isomorphism of vector spaces:
\begin{equation*}
	\tilde K \odot \tilde K  = \frac{\m_x}{(C^\infty \cdot f^* \m_y + \m_x^2)} \odot  \frac{\m_x}{(C^\infty \cdot f^* \m_y + \m_x^2)}  \stackrel{\mu}{\to} \frac{\m_x^2}{(\m_x \cdot f^* \m_y + \m_x^3)}.
\end{equation*}
\end{proposition}

\begin{proof}
Choosing local coordinates, it is easy to check that the following square commutes:
\begin{center}
\begin{tikzpicture}[thick]
	\node (LT) at (0,1.5) 	{$(\m_x/ \m_x^2) \odot (\m_x/ \m_x^2) $ };
	\node (LB) at (0,0) 	{$(C^\infty \cdot f^* \m_y/ \m_x^2) \odot (\m_x/ \m_x^2) $};
	\node (RT) at (6,1.5) 	{$\m^2_x/ \m_x^3$};
	\node (RB) at (6,0)	{${( \m_x \cdot f^* \m_y)}/{(\m_x^2 \cdot f^* \m_y) }$};
	\draw [<-right hook] (LT) --  node [left] {$$} (LB);
	\draw [->] (LT) -- node [above] {$\mu$} (RT);
	\draw [<-right hook] (RT) -- node [right] {$$} (RB);
	\draw [->] (LB) -- node [below] {$\mu$} (RB);
\end{tikzpicture}
\end{center}
Here the symmetrization $ C^\infty \cdot f^*m_y \odot \m_x$ makes sense since $C^\infty \cdot f^* \m_y \subseteq \m_x$.
Moreover, the kernel of the top map is clearly zero. A dimension count then shows that both maps in the square labeled $\mu$ are isomorphisms, and hence the induced map between the cokernels of the vertical maps is also an isomorphism.
\end{proof}

The elements of $\tilde L$ are equivalence classes of elements of $\m_y$ which happen to land in $\m_x^2$. Thus we have a natural map,
\begin{equation}\label{EqnNaturalMapIntrinsicDerivative}
	\tilde L  \to \frac{\m_x^2}{(\m_x \cdot f^* \m_y + \m_x^3)} \cong \tilde K \odot \tilde K.
\end{equation}
In fact this map is the transpose of the map $K \odot K \to L$ considered earlier. We can see this by choosing coordinates $(x_0, x_1)$ around $x$ and $(y_0, y_1)$ around $y$ such that $df$ has the form
\begin{equation*}
df_x  = \left( 
\begin{array}{cc}
	0 & 0 \\ 0 & 1
\end{array} \right).
\end{equation*}
Thus $f_0 = a x_0^2 + b x_0 x_1 + c x_1^2 + O(|x|^3)$. The above map is then given precisely by the coefficient $a$, and the map $K \odot K \to L$ is given by (a scalar multiple of) the same coefficient. 

Thus we may characterize the $S_{1,s}$ jets in terms of the local ring $\cR_f$. They are precisely those in which the map in equation \ref{EqnNaturalMapIntrinsicDerivative} has a cokernel (or kernel) of dimension $s$. The cokernel is the space 
%\index{$S_{1,r}$}
\begin{equation*}
	\frac{\m_x^2}{(\m_x^2 \cap (C^\infty_x \cdot f^* m_y) + \m_x^3)} \cong \m_f^2 / \m_f^3,
\end{equation*}
and hence is completely determined by the local ring. Moreover, we see also that the only local ring which induces an $S_{1,0}$ jet is $\cR_f \cong \R[x] /(x^2)$. 

Now we can analyze those local rings which correspond to $S_{1,1,s}$-jets. To characterize these we will need to extract from the local ring a map:
\begin{equation} \label{EqnSecondIntrinsicDerative}
	\tilde L \to \tilde K \odot \tilde K \odot \tilde K.
\end{equation}
A similar calculation as before shows that
\begin{equation*}
	\tilde K \odot \tilde K \odot \tilde K \cong \frac{\m^3_x}{(\m_x^2 \cdot f^*\m_y + \m_x^4)},
\end{equation*}
 via the multiplication map. Moreover, we are already assuming that our local map germ induces a jet in $S_{1,1}$. Thus the map from $\tilde L$ to $\tilde K \odot \tilde K$ is the zero map, and hence elements of $\tilde L$ are equivalence classes of elements of $\m_y$ whose images lies in $(\m_x \cdot f^* \m_y + \m_x^3)$. This has two consequences. First, since $\tilde L$ is also represented by equivalence classes of elements of $\m_y$ which land in $\m_x^2$, we know that the image of such an element is a linear combination of elements in $\m_x^3$ and $\m_x \cdot f^*\m_y$. In particular those elements of the form $\m_x \cdot f^* \m_y \cap \m_x^3$ actually lie in $(\m_x^2 \cdot f^*\m_y + \m_x^4)$. Secondly we have a naturally induced map:
 \begin{equation*}
	\tilde L \to \frac{\m_x \cdot f^* \m_y + \m_x^3}{(\m_x \cdot f^* \m_y + \m_x^4)} \cong \frac{\m_x^3}{(\m_x \cdot f^* \m_y \cap \m_x^3 + \m_x^4 )} \cong \frac{\m_x^3}{ \m_x^2 \cdot f^* \m_y + \m_x^4}.
\end{equation*}
This is precisely the map we needed in Equation \ref{EqnSecondIntrinsicDerative}. 

Working in local coordinates again, we see that $f_0$ has the form,
\begin{equation*}
	f_0 = a x_0 x_1 + b x_1^2 + c x_0^3 + d x_0^2 x_1 + e x_0 x_1^2 + f x_1^3 + O( |x|^4).
\end{equation*}
The above map (which is a map between 1-dimensional vector spaces) is given by the coefficient $c$ (which must be non-zero for $S_{1,1,0}$). This map does indeed agree (at least up to a non-zero scalar) with the transpose of the map $K \odot K \odot K \to L$ used to define $S_{1,1,s}$. Clearly there exist $f$ obtaining any value of $c$ that one wishes, hence the map 
\begin{equation*}
	S^{(3)}_{1,1} \to \hom( K \odot K \odot K, L)
\end{equation*}
is surjective, as claimed earlier. Thus those map germs with induced jets in $S_{1,1,s}$
%\index{$S_{1,1,s}$}
 are precisely those for which the above map has a cokernel (or kernel) of dimension $s$. The cokernel of this map is the space
\begin{equation*}
	\frac{\m_x^3}{C^\infty_x \cdot f^* \m_y \cap \m_x^3} \cong {\m_f^3}/{\m_f^4},
\end{equation*}
so once again this is completely determined by the local ring $\cR_f$, and moreover the only local ring which induces an $S_{1,1,0}$ singularity is $\R_f \cong \R[x]/(x^3)$.

The reason for passing to the algebraic description of jet space and the singularities $S_\alpha$ is that we now have many powerful algebraic tools at our disposal. This will be extremely helpful in deriving the normal coordinates associated to these singularities, especially in our later applications. The main theorem that we will use is the Generalized Malgrange Preparation Theorem and its corollaries. Proofs of these classical theorems may be found, for example, in \cite[Chapter IV]{GG73}. 
	
\begin{theorem} [Generalized Malgrange Preparation Theorem] \label{GeneralizedMalgrangePreparationTheorem}
\index{generalized Malgrange preparation theorem}
Let $X$ and $Y$ be smooth manifolds and $\phi: X \to Y$ be a smooth mapping with $\phi(p) = q$. Let $A$ be a finitely generated $C_p^\infty(X)$-module. Then $A$ is a finitely generated $C^\infty_q(Y)$-module if and only if $A/ \m_q(Y) \cdot A$ is a finite dimensional vector space over $\R$. 
\end{theorem}	

\begin{corollary} \label{GeneralizedMalgrangePreparationTheoremCorollary}
If the projections of $e_1, \dots, e_k$ form a spanning set of vectors in the vector space $A/ (\m_p^{k+1}(X) \cdot A + \m_q(Y) \cdot A)$ then $e_1, \dots, e_k$ form a set of generators for $A$ as a $C^\infty_q(Y)$-module. 
\end{corollary}

Let  $X$ be a manifold of dimension one and $f: X \times I \to \R \times I$ a generic function (i.e., transverse to $S_{i}$, $S_{i,j}$ and $S_{i,j,k}$ for all admissible $i,j,k$). The singularities of $f$ must be one of the following three types: $S_0$, $S_{1,0}$ and $S_{1,1,0}$. These occur with codimensions 0, 1, and 2 respectively. The singularities $S_0$ are not actually singularities at all; they correspond to points where the differential $df$ is invertible, and so by the inverse function theorem, $f$ is a local diffeomorphism at these points. The derivation of the local coordinates for these singularities is a well known classical result.  We will only review the $S_{1,1,0}$ case as a means to demonstrate our use of these powerful algebraic theorems. The $S_{1,0}$ case can be handled in a  similar fashion. 

The $S_{1,1,0}$ critical points are isolated in $X \times I$. Let $p \in X \times I$, be such a critical point. The differential $df$ has rank one at the point $p$ and so by the implicit function theorem, together with our assumptions on the $t$-derivatives of $\pi \circ f$, imply that there exist local coordinates $(x_0,x_1)$ for $X \times I$ centered at $p$ such that $f$ has the following form
\begin{align*}
	f^* y & = h(x_0,x_1) \\
	f^*t & = x_1
\end{align*}
(here $\pi: \R \times I \to I$ is the projection). Moreover since we are at a $S_{1,1,0}$ singularity, $h$ has an expansion,
\begin{equation*}
	h(x_0, x_1) = \alpha_{01} x_0 x_1 + \alpha_{11} x_1^2 + \beta_{000} x_0^3 + \beta_{001} x_0^2 x_1 + \beta_{011} x_0 x_1^2 + \beta_{111} x_1^3 + O(|x|^4)
\end{equation*}
in which $\beta_{000}$ is non-zero. In fact, by assumption $j^2f$ is transverse to $S_{1,1}$, which is equivalent to the statement that the following map of vector spaces,
\begin{equation*}
	T_pX \to T_{j^2f(p)} J^2(X \times I, \R \times I) / T_{j^2f(p)} S_{1,1} \cong \nu S_{1,1} \cong \hom(K, L) \oplus \hom( K \odot K, L),
\end{equation*}
is an isomorphism. In these coordinates, this map is given by the matrix,
\begin{equation*}
	\left( \begin{array}{c c}
			0 & \alpha_{01} \\
			\beta_{000} & \beta_{001}
		\end{array} \right).
\end{equation*}
In particular the coefficient $\alpha_{01}$ must be non-zero.

The local ring $\cR_f(p)$ is isomorphic to $\R[ x_0] / (x_0^3)$, and so by the Generalized Malgrange Preparation Theorem~\ref{GeneralizedMalgrangePreparationTheorem}, and more specifically Corollary \ref{GeneralizedMalgrangePreparationTheoremCorollary}, we know that there exist smooth functions $a,b,c$ on $\R \times I$ such that,
\begin{equation*} \label{EquationCerfDerivationS110}
	x_0^3 = f^* a + f^*b \cdot x_0 + f^*c \cdot x_0^2.
\end{equation*}
Moreover, by replacing $x_0$ by the coordinate $x_0 + \frac 1 2 f^*c$ (and keeping all other coordinates the same) we may assume that $c=0$. Expanding both sides of Equation \ref{EquationCerfDerivationS110} when $x_1 = t = 0$ and identifying the lowest order terms we find that
\begin{equation*}
	\frac{\partial a}{ \partial y}(p) = \frac{1}{\beta_{000}} \neq 0
\end{equation*}
Thus the following is a valid change of coordinates,
\begin{align*}
	\overline y &= a(y,t) \\
	\overline t &= t
\end{align*}
and in these coordinates Equation \ref{EquationCerfDerivationS110} becomes
\begin{equation*} \label{EquationCerfDerivationS110-Second}
	x_0^3 = f^*y + f^*b \cdot x_0.
\end{equation*}
Now expanding Equation \ref{EquationCerfDerivationS110-Second} in both $x_0$ and $x_1$ and gathering the $x_1 x_0$ terms, we see that
\begin{equation*}
	\frac{\partial b}{\partial t}(q)  = - \frac{\alpha_{01}}{\beta_{000}} \neq 0.
\end{equation*}
Thus the following is a valid pair of coordinate changes,
\begin{align*}
	\overline y_0 &= y \\
	\overline y_1 & = b(y,t) \\
	\overline x_0 &= x_0 \\
	\overline x_1 &= f^*b.
\end{align*}
Note: depending on the signs of $\alpha_{01}$ and $\beta_{000}$ this may be an orientation reversing coordinate change. These results may be summarized by the following proposition.
\begin{proposition}
Let $X$ be a manifold of dimension one, $f: X \times I \to \R \times I$ be a generic function, and $p \in X$ be a point with an $S_{1,1,0}$ singularity. Then there exist local coordinates $(x_0, x_1)$ of $X \times I$ centered at $p$ and coordinates $(y_0, y_1)$ of $Y \times I$ centered at $q = f(p)$ such that $f$ has the following form:
\begin{align*}
	f^* y_0 &= x_0^3 + x_1 x_0 \\
	f^* y_1 &= x_1.
\end{align*}
\end{proposition}

A similar analysis applies and gives local coordinates for the $S_{1,0}$ singularities. One may then use these local coordinates to get a global description of generic functions $F: X \times I \to \R \times I$. Any two Morse functions may be extended to such a generic function $F$ and each of the singularities presented has an interpretation in terms of the critical points\index{critical point}\index{birth death} of the Morse function. The $S_0$ points correspond to the regular values of the Morse function, the $S_{1,0}$ singularities correspond to paths of critical points moving around in $X$, and the $S_{1,1,0}$ singularities correspond to birth-death critical point cancelations. One may use this description to prove that any two Morse functions on a (compact) 1-manifold $X$ are related by a path of functions which consists entirely of Morse functions except at a finite number of times, in which a birth-death occurs. We will not pursue this here. Instead we now turn to higher dimensional decompositions.

\subsection{Strategy of the Planar Decomposition Theorem} \label{SectStratofPlaarDecomp}

We have seen how a generic map from a manifold $X$ to $\R$ induces a 1-dimensional decomposition of the manifold. Similarly a generic map from $\Sigma$ to $\R^2$ will induce a 2-dimensional decomposition of $\Sigma$. Just as studying certain generic maps $X \times I \to \R \times I$ allows us to compare different 1-dimensional decompositions, studying generic maps $\Sigma \times I \to \R^2 \times I$ will allow use to compare different 2-dimensional decompositions. 
	
Later we will want to use these planar decompositions 
\index{planar decomposition theorem}
of $\Sigma$ to extract 2-categorical information. For this purpose it is also useful to control not just the map $\Sigma \to \R^2$, but also the composition with the projection to $\R$, a map $\Sigma \to \R$. This is similar to the requirement that we imposed on our maps $X \times I \to \R \times I$. We required that, after composing with the projection to $I$, the map $X \times I \to I$ satisfies certain properties. When we look at maps $\Sigma \times I \to \R^2 \times I$, we will also want to control the projections to $\R \times I$ and to $I$. Incorporating this control is in principle straightforward. We end up with a stratification of the jet bundles which is slightly more refined than the classical Thom-Boardman stratification 
\index{Thom-Boardman stratification}
considered in the previous sections. 
	
After giving this stratification, we will give an analysis similar to the previous one,  deriving normal coordinates for these singularities. For the case of maps $\Sigma \to \R^2$ with $\Sigma$ a manifold of dimension two, this is nearly identical to the derivation of the last section. In the 3-dimensional case some non-standard decompositions appear, but the method of their derivation is not essentially different than the classical case. We will also deal more carefully with multi-jet issues, which were completely ignored in our previous discussion.

The end result of these considerations will be that a generic map $\Sigma \to \R^2$ from a surface will induce a planar decomposition of $\Sigma$ into elementary pieces. Two such generic maps can then by connected by a generic map $\Sigma \times I \to \R^2 \times I$, hence two such planar decompositions are related if there is a decomposition of $\Sigma \times I$ into a 3-dimesnional arrangement of elementary pieces. These give ``relations'' or ``moves'' connecting any two planar decompositions. Finally we will discuss the minor changes necessary to adapt these results to the case where $\Sigma$ has boundary and corners. For simplicity, in what follows we will restrict our attention to the case where $\Sigma$ is a 2-dimensional manifold. The general situation requires a more nuanced discussion involving the choice of gradient like vector fields and similar accoutrements.
	
\subsection{The Stratification of the Jet Space: 2D case} \label{SectStratofJetSpace2D}

The stratification of the jet space that we will use for maps $\Sigma \to \R^2$ is very similar to the Thom-Boardman stratification we used to study maps $X \times I \to \R \times I$. In fact two of our strata, $S_0$ and $S_2$ of codimensions zero and four, respectively, are defined exactly as they were for the Thom-Boardman stratification. The third Thom-Boardman stratum, $S_1$, will be spilt into two separate strata $S_{[01]}$ and $S_{[11]}$ of codimensions one and two, respectively. This finer stratification results from our consideration of the composition $\Sigma \to \R^2 \to \R^1$. These strata then induce preliminary stratifications of the higher jets, just as before, and, also just as before, we may further refine these strata in a straightforward way. All in all, we have the strata listed in Table~\ref{Table2DStratificationOfJets}, grouped according to which jet space they reside in, together with their respective codimensions.
\begin{table}[ht]
\begin{center}
\begin{tabular}{|c|c| c|} \hline
Strata & Codimension  & Name\\ \hline \hline
$S_0$ & 0 & Local Diffeomorphism\\
$S_{[01]}$ & 1& $\cdot$ \\
$S_{[11]}$ & 2 &  $\cdot$\\
$S_{2}$ & 4 & n/a\\ \hline
$S_{[01], 0}$ & 1 & Fold\\
$S_{[01], 1}$ & 2 & $\cdot$\\
$S_{[11], 0}$ & 2 & 2D Morse \\
$S_{[01],1 }$ & 3 & $\cdot$\\ \hline
$S_{[01], 1, 0}$ & 2 & Cusp \\
$S_{[01],1,1}$ & 3 & $\cdot$ \\ \hline
\end{tabular}
\end{center}
\caption{The Stratification of Jet Space: 2D Case}
\label{Table2DStratificationOfJets}
\end{table}%

Based off these codimensions, a generic map $f: \Sigma \to \R^2$ only admits singularities of types $S_0, S_{[0,1], 0}, S_{[11], 0}$, and  $S_{[0,1], 1, 0}$. We take up the task of deriving normal coordinates for these in the next section. Recall that the Thom-Boardman stratum $S_1$ of $J^1( \Sigma, \R^2)$ was defined as the subset of jets, viewed as elements of $\hom( T_x\Sigma, T_y\R^2)$, having corank one (and hence also rank one). The projection $p:\R^2 \to \R$ induces a projection,
\begin{equation*}
	J^1(\Sigma, \R^2) \to J^1( \Sigma, \R)
\end{equation*}
 which on fibers is the projection $ \hom( T_x\Sigma, T_y\R^2) \to \hom( T_x \Sigma, T_{p(y)}\R)$. This later space can also be stratified by corank and pulling this stratification back induces a refined stratification of $S_1$.\footnote{That $S_{[01]}$ is a manifold is clear. The fact that $S_{[11]}$ is also a manifold can be deduced by looking at local coordinates for a  representative map germ.} There are two such strata corresponding to when $d(pf)$ has corank zero and corank one (rank one and rank zero, respectively). This defines the strata $S_{[s \; t]}$ which consists of jets where $d(pf)_x$ has corank $s$ and $df_x$ has corank $t$. 

We see that $S_{[01]}$ is an open subset of $S_1$. Let $S_{[01]}^{(2)}$ and $S_{[01]}^{(3)}$ denote the inverse image of $S_{[01]}$ in $J^2(\Sigma, \R^2)$ and $J^3(\Sigma, \R^2)$ under the respective projections. We define,
\begin{align*}
		S_{[01], r} &:= S^{(2)}_{[01]} \cap S_{1,r}, \\
		S_{[01], r, s} &:= S^{(3)}_{[01]} \cap S_{1,r, s}.
\end{align*}
These are manifolds since $S_{[01]}$ is open in $S_1$. $S_{[11]}$ has codimension 2 in the jet space $J^1(\Sigma, \R^2)$ and its normal bundle is canonically isomorphic to $\hom( T \Sigma, T\R)$ (where these bundles are interpreted as their pull-back to $S_{[11]}$. Over $S_1$ there are two canonical bundles, $K$ and $L$, which we introduced previously, and these restrict to bundles over $S_{[11]}$. Let $S_{[11]}^{(2)}$ denote the inverse image of $S_{[11]}$ in $J^2(\Sigma, \R^2)$. As before we have a map of fiber bundles over $S_{[11]}$,
\begin{equation*}
	S_{[11]}^{(2)} \to \hom( T\Sigma, TJ^1( \Sigma,  \R^2)) \to \hom(K, \hom(K, L)).
\end{equation*}
This composition lands in $\hom( K \odot K, L)$ and is surjective, as can be checked by choosing local coordinates. Our stratification of $\hom(K \odot K, L)$ into the two strata $\hom(K \odot K, L)_0$ and $\hom (K \odot K, L)_1$ induces a refinement of $S_{[11]}^{(2)}$ into the strata $S_{[11],0}$ and $S_{[11], 1}$ respectively. These have codimensions 2 and 3 in $J^2(\Sigma, \R^2)$, respectively. This completes our stratification of the jet space in the 2-dimensional case.

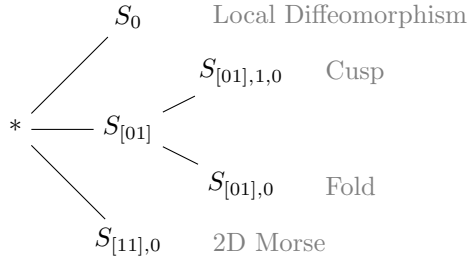
\begin{figure}[ht]
\begin{center}
\begin{tikzpicture}[grow via three points={% 
one child at (0,1) and two children at (-.5,1) and (.5,1)}, grow = right] 
\node at (0,0) { *} 
	child { node (C) {$S_{[11],0}$} node [anchor = west] at +(1,0) { \color{gray} 2D Morse}}	
	child {node {$S_{[01]}$}
		child {node {$S_{[01],0}$} node [anchor = west] at +(1,0) { \color{gray} Fold}}
		child {node {$S_{[01],1,0}$} node [anchor = west] at +(1,0) { \color{gray} Cusp}}
		}
	child {node {$S_0$} node[anchor = west]  at +(1,0) {  \color{gray} Local Diffeomorphism}};
\end{tikzpicture} 
\caption{The 2D Singularity Family Tree }
\index{singularity!3D family tree}
%\label{default}
\end{center}
\end{figure}
	
\subsection{Normal Coordinates: 2D Case} \label{SectNormalCoord2D}

We must derive normal coordinates for the four kinds of singularities of type $S_0$, $S_{[0,1], 0}$, $S_{[11], 0}$, and  $S_{[0,1], 1, 0}$. The $S_0$ points consist of local diffeomorphisms, as already seen. Since we are making use of the projection $\R^2 \to \R$, it is natural for us to find coordinates for $\R^2$ which are compatible with this projection. This is not always possible, however we will always be able to find coordinates $(\overline y_0, \overline y_1)$ in which the kernel of the differential of the projection at $q \in \R^2$ is the span of $\partial_{\overline y_0}$.  
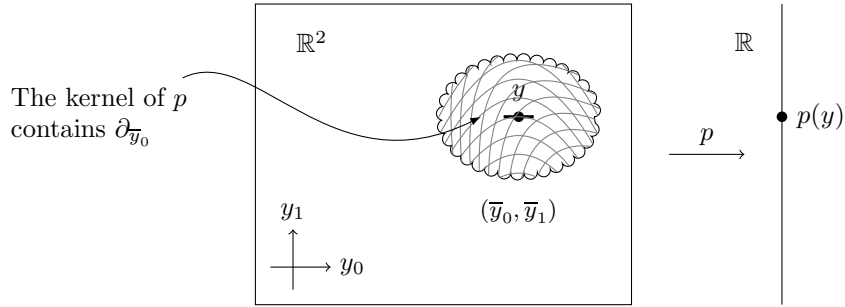
\begin{figure}[ht]
\begin{center}
\begin{tikzpicture}
	\draw (0,0) rectangle (5,4); 
	\draw [->] (.2,.5) -- (1, .5) node [right] {$y_0$};
	\draw [->] (.5, .2) -- (.5, 1) node [above] {$y_1$};
	\node at (.75, 3.5) {$\R^2$};
	\draw [->] (5.5, 2) -- node [above] {$p$} (6.5, 2); 
	\draw (7,0) -- +(0,4);
	\node at (6.5, 3.5) {$\R$};
	\begin{scope}
		\node [circle, fill=black, minimum size=4pt, inner sep=0pt] (q) [label=above:$y$] at (3.5, 2.5) {};
		\draw [clip, decorate,decoration={bumps,mirror}] 
			(q) ellipse (1cm and .75cm); 
		\draw [very thick] (3.3, 2.5) -- (3.7, 2.5);
		\foreach \y in {-.25, 0, ..., 1.75} 
			\draw [gray] (2, \y) parabola[bend pos=0.5] bend +(0,2) +(3,0); 
		
		\begin{scope} [rotate = 30, yshift = -2cm, xshift = 1cm]
			\foreach \y in {-.25, 0, ..., 1.75} 
				\draw [gray] (2, \y) parabola[bend pos=0.5] bend +(0,2) +(3,0); 
		\end{scope}
	\end{scope}	
	\node at (3.5, 1.25) {\small $(\overline  y_0, \overline  y_1)$};
	\node [circle, fill=black, minimum size=4pt, inner sep=0pt]  [label=right:$p(y)$] at (7, 2.5) {};
	
	\node [text width = 2.5cm] (k) at (-2, 2.5) {The kernel of $p$ contains $\partial_{\overline y_0}$};
	\draw   [-latex] (k) .. controls +(2, 1) and +(-2,-1) .. +(5,0);
\end{tikzpicture} 
\caption{The Standard Coordinates for $\R^2$ and its projection to $\R$}
\label{TheStandardCoordR2Fig}
\end{center}
\end{figure}
This is depicted in Figure~\ref{TheStandardCoordR2Fig}. 
We will also make the convention that $(y_0, y_1)$ denote the standard coordinates in $\R^2$ and that the projection is given by $p(y_0, y_1) = y_1$. For later purposes, it will also be convenient to keep track of whether the coordinates $(\overline y_0, \overline y_1)$ preserve the standard orientation on $\R^2$, and whether the orientation on $\R$ induced by the image of $\partial_{\overline y_1}$ agrees with the standard orientation of $\R$. If this is the case, then we say  $(\overline y_0, \overline y_1)$ preserves both the orientations of $\R$ and $\R^2$.

\begin{proposition} [Cusp (see Figure~\ref{CuspGraphicFig} on pg.~\pageref{CuspGraphicFig})] \label{2dCuspProp}
\index{cusp singularity}
\index{singularity!cusp}
Let $\Sigma$ be a surface, $f: \Sigma \to \R^2$ be generic, and $x \in \Sigma$ be a critical point with singularity type $S_{[01], 1, 0}$. Then there exist coordinates $(\overline x_0, \overline x_1)$ of $\Sigma$ centered at $x$ and coordinates $(\overline y_0,\overline  y_1)$ of $\R^2$ centered at $y = f(x)$ such that $\partial_{y_0} \in \ker p$ at $y$, such that in these coordinates $f$ takes the normal form:
\begin{align*}
	f^* \overline y_0 &= \overline x_0^3 \pm \overline x_1 \overline x_0 \\
	f^* \overline y_1 &= \overline x_1.
\end{align*}
Moreover, the coordinates $(\overline y_0,\overline y_1)$ can be chosen so as to agree with the orientations of $\R^2$ and $\R$ at $y$ and $p(y)$. 
\end{proposition}
\begin{proof}
This follows essentially from the calculations of the previous section.  The differential $df$ has rank one at the point $x \in \Sigma$ and the differential $d(pf)$ also has rank one. Thus by the implicit function theorem there exist local coordinates $(x_0,x_1)$ for $\Sigma$ centered at $x$ such that $f$ has the following form:
\begin{align*}
	f^* y_0 & = h(x_0,x_1) \\
	f^*y_1 & = x_1
\end{align*}
Moreover, $h$ has an expansion of the form,
\begin{equation*}
	h(x_0, x_1) =\sum_{0 \leq i \leq j \leq 1} \alpha_{ij} x_i x_j + \sum_{0 \leq i \leq j \leq k \leq 1} \beta_{ijk} x_i x_j x_k +  O(|x|^4).
\end{equation*}
Since $x$ is an $S_{[01], 1,0}$ we know that $\alpha_{00} = 0$ and that $\beta_{000} \neq 0$. Morevoer, since $j^3f \pitchfork S_{[01], 1,0}$, the matrix,
\begin{equation*}
	\left( \begin{array}{c c}
			0 & \alpha_{01} \\
			\beta_{000} & \beta_{001}
		\end{array} \right).
\end{equation*}
is non-degenerate, which is equivalent to $\alpha_{01} \neq 0$ (together with the already satisfied condition,
 $\beta_{000} \neq 0$). As before the local ring is isomorphic to $\cR_f(x) = \R[x_0]/(x_0^3)$, so by the Malgrange Preparation Theorem~\ref{GeneralizedMalgrangePreparationTheorem}, there are functions $a, b,c$ on $\R^2$ such that 
 \begin{equation*}
	x_0^3 = f^*a + f^*b \cdot x_0 + f^*c \cdot x_0^2.
\end{equation*}
The same coordinate change in the domain $\Sigma$ that was used last time allows us to assume that $c \equiv 0$. By expanding these variables and collecting similar terms we find that,
\begin{equation*}
	\frac{\partial a}{\partial y_0}(y) \neq 0, \qquad \frac{\partial a}{\partial y_1}(y) = 0, \qquad \frac{\partial b}{\partial y_1}(y) \neq 0.
\end{equation*}
Thus we find that the following is a valid pair of coordinate changes, and that at $y$ $\partial_{\overline y_0}$ is proportional to $\partial y_0$, 
\begin{align*}
	\overline y_0 &= \pm a(y_0, y_1) \\
	\overline y_1 & = \pm b(y_0,y_1) \\
	\overline x_0 &= \pm x_0 \\
	\overline x_1 &=  \pm f^*b.
\end{align*}
Choosing the first two signs in the above coordinate change correctly we may arrange for the desired property regarding orientations, and adjusting the sign of $\pm x_0$ appropriately ensures that in these coordinates $f$ has the desired normal form. 
\end{proof}

\begin{proposition} [2D Morse (see Figure~\ref{2DMorseGraphicFig} on pg.~\pageref{2DMorseGraphicFig})] \label{2d2dmorseprop}
\index{singularity!2D Morse}
\index{2D Morse singularity}
	Let $\Sigma$ be a surface, $f: \Sigma \to \R^2$ be generic, and $x \in \Sigma$ be a critical point with singularity type $S_{[11], 0}$. Then there exist coordinates $(\overline x_0, \overline x_1)$ of $\Sigma$ centered at $x$, and coordinates $(\overline y_0)$ of $\R$ centered at $p(y) = pf(x)$ (thus $(\overline y_0, y_1)$ form coordinates for $\R^2$, with $y_1$ standard),  such that $f$ takes the following normal form,
 \begin{align*}
	f^* \overline y_0 &= \overline x_0 \\
	f^* y_1 &= \pm \overline x_0^2 \pm \overline x_1^2.
\end{align*}
Moreover, these coordinates preserve the orientations of $\R^2$ and $\R$.   
\end{proposition}
\begin{proof}
The differential $df$ has rank one at the point $x \in \Sigma$ and the differential $d(pf)$ has rank zero. Thus by the implicit function theorem there exist local coordinates $(x_0,x_1)$ for $\Sigma$ centered at $x$ such that $f$ has the following form:
\begin{align*}
	f^* y_0 & =x_0 \\
	f^*y_1 & = h(x_0,x_1).
\end{align*}
Since $j^2f \pitchfork S_{[11], 0}$ we know that the map, 
\begin{equation*}
	T_xX \to \hom(T_xX, T_{p(y)}\R) \cong (\nu S_{[11], 0})_{j^2f(x)} ,
\end{equation*}
is a surjection. This implies that the Hessian of $h$ is non-degenerate. The expansion of $h$ has the form,
\begin{equation*}
	h(x_0, x_1) =\sum_{0 \leq i \leq j \leq 1} \alpha_{ij} x_i x_j +  O(|x|^3),
\end{equation*}
and our conditions that $x$ is an $S_{[11], 0}$ singularity ensure that $\alpha_{11} \neq 0$. 

Let $K={- \alpha_{01}}/{\alpha_{00}}$. The following coordinate change is valid:
\begin{align*}
	\tilde x_0 &= x_0 \\
	\tilde x_1 &= x_1 - K x_0
\end{align*}
and in these coordinates we have (removing the tildes), 
\begin{equation*}
	h = \alpha_{11} \cdot x_1^2 + \beta \cdot x_0^2 + O(|x|^3),
\end{equation*}
where $\beta \neq 0$. Thus for some functions $a(x_0)$, $b(x_0)$, and $c(x_0, x_1)$, such that $a(0) = \beta\neq 0$ and $c(0,0) = \alpha_{11} \neq 0$, we have
\begin{equation*}
	h = a(x_0) \cdot x_0^2 + 2 b(x_0) \cdot x_0^2 x_1 + c(x_0, x_1) \cdot x_1^2.
\end{equation*}
Now we can make the coordinate change,
\begin{align*}
	\tilde x_0 &= x_0 \\
	\tilde x_1 &= x_1 + \frac{b(x_0)}{\alpha_{11}} \cdot x_0^2
\end{align*}
under which we have (again removing tildes),
\begin{equation*}
	h = \tilde a(x_0) \cdot x_0^2 + \tilde c(x_0, x_1) \cdot x_1^2, 
\end{equation*}
for some new functions $\tilde a$ and $\tilde c$ such that $\tilde a(0) = \beta \neq 0$ and $c(0,0) = \alpha_{11} \neq 0$. Finally we form the following valid pair of coordinate changes,
\begin{align*}
	\overline x_0 &= x_0{\sqrt{\tilde a(x_0)} }\\
	\overline x_1 &=  x_1{\sqrt{\tilde c(x_0, x_1)} }\\
	\overline y_0 &= y_0{\sqrt{\tilde a(y_0)} }\\
	\overline y_1 & =y_1 
\end{align*}
at least in a sufficiently small neighborhood of $x$ and $y$. In these coordinates $f$ has the claimed normal form. 
\end{proof}

Let $S_{[01], 0}(f) = Y$ denote the set of elements of $\Sigma$ which have $S_{[01], 0}$ singularities. By transversality, we know that $Y$ is an embedded 1-dimensional submanifold of $\Sigma$. The conditions the we used to define the stratum $S_{[01], 0}$ ensure that both $d(f|_Y)$ and $ d( pf|_Y)$ are non-degenerate. In particular, $pf$ is a local diffeomorphism and  the map $f: Y\to \R^2$ is (locally) an embedding. In fact (locally) its image is the graph $\{ (\gamma(y_1), y_1) \}$ of a function  $\gamma: \R \to \R$. These considerations give the first half of the following proposition.

\begin{proposition}[Folds (see Figure~\ref{FoldGraphicsFig} on pg.~\pageref{FoldGraphicsFig})] \label{2dFoldsProp}
\index{singularity!fold}
\index{fold singularity}
Let $\Sigma$ be a surface, $f: \Sigma \to \R^2$ be generic, and $x \in \Sigma$ be an $S_{[01],0}$ singularity point. Then there exists a submanifold  $Y^1 \subseteq \Sigma$ containing $x$ and consisting entirely of $S_{[01],0}$ singularities,  a neighborhood $U\subseteq \R$ of $p(y) = pf(x)$, and 
 a function $\gamma: U \to \R$ whose graph $\{ (\gamma(y_1), y_1) \; | y_1 \in U \}$ consists precisely of the image of $Y$ under $f$. Moreover there exist coordinates $(x_0, x_1)$ for $\Sigma$ centered at $x$, such that in these coordinates $Y = \{ x_0 = 0 \}$ (the $x_1$-axis) and such that $f$ has the following normal form:
\begin{eqnarray*}
	f^* y_0 &=&  \pm \overline{x_0}^2 + \gamma( x_1)\\
	f^*y_1 &=& x_1.
\end{eqnarray*}
where $(y_0, y_1)$ are the standard coordinated for $\R^2$. 
\end{proposition}

\begin{proof} The existence of $Y$, $U$, and $\gamma$ as in the statement of the proposition have already been established. By the implicit function theorem we may choose coordinates $(x_0, x_1)$ for $\Sigma$ centered at $x$, such that $Y = \{ x_0 = 0 \}$, and such that in these coordinates
\begin{eqnarray*}
	f^* y_0 &=&  h(x_0, x_1)\\
	f^*y_1 &=& x_1.
\end{eqnarray*}
As before, we can expand $h$ as,
\begin{equation*}
	h(x_0, x_1) = \sum_{0 \leq i \leq j \leq 1} \alpha_{ij} \cdot x_i x_j + O(|x|)^3
\end{equation*}
and our assumption that the points of $Y$ are $S_{[01],0}$ singularity points implies that $h(0, x_1) = \gamma(x_1)$ and that
\begin{equation*}
	\frac{\partial h}{\partial x_0} (0, x_1) \equiv 0.
\end{equation*}
We temporarily consider the new function $g(x_0, x_1) = h(x_0, x_1) - \gamma(x_1)$. We have $\frac{\partial g}{\partial x_0} (0, x_1) \equiv 0$ and $g(x_0, x_1) \equiv 0$. Thus we may write
\begin{equation*}
	h(x_0, x_1) = x_0^2 \tilde g(x_0, x_1) + \gamma(x_1)
\end{equation*}
for some new function $\tilde g$. Our transversality assumptions imply that $\tilde g (0, x_1) $ is a non-vanishing function, and so the following is a valid coordinate change:
\begin{eqnarray*}
	\overline x_0 &=& x_0{\sqrt{\tilde g(x_0, x_1)}}\\
	\overline x_1 &=& x_1,
\end{eqnarray*}
at least in a sufficiently small neighborhood of $Y$. In these coordinates $f$ has the desired form. 
\end{proof}

\subsection{The Stratification of the Jet Space: 3D Case} \label{SectStratofJetSpace3D}

We will also want to fully understand generic maps from $\Sigma \times I$ to $\R^2 \times I$ in a similar manner to our understanding of generic maps for $\Sigma$ to $\R^2$. As before, we will accomplish this by judiciously choosing a stratification of the jet spaces $J^k(\Sigma \times I, \R^2 \times I)$, and again we will want to leverage the fact that we have projections $p: \R^2 \times I \to \R \times I$ and $p:\R \times I \to I$. Due to the higher dimensionality, our stratification is necessarily more complicated than the previous stratification, but the method which we use to construct the stratification is, in principle, the same. Ignoring multi-jet issues, we will require that our generic maps induce sections $j^kf$ which are transverse to the 28 strata listed in Table~\ref{Table3DStratificationOfJets}, which lists the strata and their codimensions, grouped according to which jet space they reside in. 
\begin{table}[ht]
\begin{center}
\begin{tabular}{|c|c| c|} \hline
Strata & Codimension & Name \\ \hline \hline
%%%  The 1-Jets
$S_0$ & 0 & Local Diffeomorphism\\
$S_1$ & 1 & $\cdot$ \\
$S_{[001]}$ & 1& $\cdot$ \\
$S_{[011]}$ & 2 &  $\cdot$\\
$S_{[111]}$ & 3 & $\cdot$ \\
$S_{2}$ & 4 & $\cdot$\\ 
$S_3$ & 9 & $\cdot$ \\ \hline
%%% The 2-Jets
$S_{[001], 0}$ & 1& 1D Morse (Fold) \\
$S_{[001], [01]}$ & 2& $\cdot$ \\
$S_{[001], [11]}$ & 3 & $\cdot$ \\
$S_{[011], [000]}$ & 2 &  2D Morse\\
$S_{[011], [001]}$ & 3 &  $\cdot$\\
$S_{[011], [100]}$ & 3 &  $\cdot$\\
$S_{[011], [101]}$ &  4 &  $\cdot$\\
$S_{[011], [111]}$ &  5 &  $\cdot$\\
$S_{[011], [211]}$ &  6 &  $\cdot$\\
$S_{[111], 0}$ & 3 & 3D Morse \\ 
$S_{[111], 1}$ & 4 & n/a \\ \hline
%%% The 3-Jets
$S_{[001], [01], 0}$ & 2& 1D Morse Relation (Cusp) \\
$S_{[001], [01], 1}$ & 3& $\cdot$ \\
$S_{[001], [11], 0}$ & 3 & Cusp Inversion \\
$S_{[001], [11], 1}$ & 4 & $\cdot$ \\
$S_{[011], [001], 0}$ & 3 & Cusp Flip \\
$S_{[011], [001], 1}$ & 4 &  $\cdot$\\
$S_{[011], [100],0}$ & 3 & 2D Morse Relation \\
$S_{[011], [100],1}$ & 4 &  $\cdot$\\  \hline
%%% The 4-Jets
$S_{[001], [01], 1, 0}$ & 3& Swallowtail \\
$S_{[001], [01], 1, 1}$ & 4& $\cdot$ \\
\hline
\end{tabular}
\end{center}
\caption{The Stratification of Jet Space: 3D Case}
\label{Table3DStratificationOfJets}
\end{table}%

Based on this stratification and the respective codimensions, we see that a generic map $f: \Sigma \times I \to \R^2 \times I$ will only have singularities of the nine types labeled with names in Table~\ref{Table3DStratificationOfJets}. Among these, we will be able to eliminate the possibility of ``3D-Morse'' singularities ($S_{[111], 0}$) by restricting ourselves to a suitable open subset of maps $f$. To complete our analysis, we will need to derive local coordinates for these eight remaining singularity types. This is taken up in the following section. 

For these dimensions (a 3-manifold mapping to a 3-manifold) the classical Thom-Boardman stratification divides the first jet space into four submanifolds, $S_r$ with  $r = 0, 1, 2, 3$, based on the corank $r$ of the differential $df$. Using our projections, $p$, we can further refine the strata $S_1$ and $S_2$, however since $S_2$ is already of codimension four, we will only need to refine $S_1$. %\footnote{ Maps $X \times I \to \R^2 \times I$ when $X$ has dimension greater than two can also be studied in this manner. In this case one would need to consider the $S_2$ stratum and also its  stratification [Is this true?]}
 Just as for the 2-dimensional stratification, we will compare the relative coranks of all the differentials $df$, $d(pf)$ and $d(p^2 f)$. In this case we see that there are three possibilities, which are indexed by $[r \; s \; t]$, where $r$ is the corank of $d(p^2 f)$, $s$ is the corank of $d(pf)$, and $t$ is the corank of $df$ (which is 1 in all cases). The proof that these are manifolds is completely analogous to the previous 2-dimensional case. 

Let $S^{(2)}_\alpha$ denote the inverse image of $S_\alpha$ in $J^2( \Sigma \times I , \R^2 \times I)$ under the natural projection to $J^1(\Sigma \times I , \R^2 \times I)$. The classical Thom-Boardman stratification of the second jet space spilts $S^{(2)}_1$ into two strata $S_{1,1}$ and $S_{1,0}$. These are defined by considering the natural bundles $K$ and $L$ over $S_1$ (which are, respectively, the kernel and cokernel of $df$) and looking at the natural map of fiber bundles (over $S_1$):
\begin{equation*}
	S^{(2)}_1 \to \hom( K \odot K, L).
\end{equation*}
As before this is a submersion and the range is stratified by the corank $r$ of the corresponding map in $\hom(K, \hom(K, L)$. Pulling back the stratification to $S^{(2)}_1$ defines the manifolds $S_{1,r}$.

This allows us to define several new strata as follows:
\begin{align*}
S_{[001], 0} &: =  S_{[001]}^{(2)} \cap S_{1,0} \\
S_{[111], 0} &  := S_{[111]}^{(2)} \cap S_{1,0} \\ 
S_{[111], 1} & := S_{[111]}^{(2)} \cap S_{1,1} 
\end{align*}
The first %three
two 
of these are clearly manifolds, since $S_{1,0}$ is open in $S_1^{(2)}$. The last one is also a manifold, which can be seen by choosing local coordinates. This leaves three intersections which give us a preliminary stratification,
\begin{align*}
	S_{[001], 1} & := S_{[011]}^{(2)} \cap S_{1,1} \\
	S_{[011], 0} & := S_{[011]}^{(2)} \cap S_{1,0} \\
	S_{[011], 1} & := S_{[011]}^{(2)} \cap S_{1,1}. 
\end{align*}
The first and second of these will be divided into two strata each, and the last into four. Let us first consider the $S_{[001]} $ stratum. Over this stratum there are three important  vector bundles $K_1$, $K_2$ and $L$, which are the kernel of $df$, the kernel of $d(p^2f)$, and the cokernel of $df$, respectively. Thus $K_2$ is 2-dimensional and contains $K_1$ as a 1-dimensional subspace. The jet space induces a natural bundle map
\begin{equation*}
	S^{(2)}_{[001]} \to \hom (T(\Sigma \times I),  TJ^1(\Sigma \times I, \R^2 \times I))
\end{equation*}
and so we may project to the normal bundle of $S_{[001]}$ (which, recall, is isomorphic to $\hom(K,L)$) and restrict to both of the two kernels. Thus, there exist natural maps of fiber bundles over $S_{[001]}$,
\begin{equation*}
	S^{(2)}_{[001]} \to \hom ( K_2 \odot K_1, L) \to \hom(K_1 \odot K_1, L). 
\end{equation*}
  We may attempt to stratify $S^{(2)}_{[001]}$ by the coranks of the corresponding operators in the spaces 
\begin{align*}
	\hom(K_2, \hom(K_1, L)) &\cong \hom(K_1, \hom(K_2, L)), \\
	\text{and} \quad & \hom(K_1, \hom(K_1,L)).
\end{align*} 
    The possibilities are labeled $[r \; s]$, where $r$ is the corank of the operator in \\  $\hom(K_2, \hom(K_1, L))$ (which agrees with the corank of the operator in \\ $ \hom(K_1, \hom(K_2, L)) $ ) and $s$ is the corank of the operator in $\hom(K_1, \hom(K_1, L))$. Thus there are three possibilities: $[00]$, $[01]$, and $[11]$, corresponding to a stratification of $S^{(2)}_{[001]}$ into three submanifolds $S_{[001], [00]}$, $S_{[001], [01]}$, and $S_{[001], [11]}$. The first of these is exactly the stratum $S_{[001], 0}$, which we have already defined. Proving that these are submanifolds of the jet space of the claimed codimensions is now relatively straightforward.

The $S_{[011]}$ stratum leads to more possibilities, as we shall see. Over the stratum $S_{[011]}$, there are several naturally defined bundles. There are three fundamental bundles $K_1$, $K_2$, and $L$, which correspond to the kernel of $df$, the kernel of $d(pf)$ and the cokernel of $df$, respectively. $S_{[011]}$ sits inside $S_1$, and so we may restrict the normal bundle of $S_1$ to $S_{[011]}$ . This is the bundle $\hom(K_1, L)$. In addition, we may consider the normal bundle of $S_{[011]}$ itself. This is the bundle $\hom(K_2, L)$. This leads to the following sequence of fiber bundle maps over $S_{[011]}$,
\begin{equation*}
S_{[011]}^{(2)} \to \hom( K_2 \odot K_2, L) \to \hom( K_2 \odot K_1, L) \to \hom( K_1 \odot K_1, L).
\end{equation*}

Each of these spaces, $\hom( K_i \odot K_j, L)$ can be stratified according to the corresponding corank in $\hom(K_i, \hom(K_j, L))$. This leads to the following six possibilities: 
\begin{equation*}
	[001], \; [000], \;[100], \;[101], \; [111], \; \text{ and } \; [211]
\end{equation*}
where $[r \; s \; t]$ corresponds to the case where the element in $\hom(K_2, \hom(K_2, L))$ has corank $r$, the element in $\hom(K_2, \hom(K_1, L))$ has corank $s$, and the element in $\hom(K_1, \hom(K_1, L))$ has corank $t$. Again it is straightforward to check that these are submanifolds of the jet space of the claimed codimensions. 

The Thom-Boardman stratification continues and our remaining strata may be defined using the Thom-Boardman stratification together with the strata we have already constructed. The normal bundle to $S_{1,1}$ in $S_{1}^{(2)}$ is isomorphic to the bundle $\hom( K \odot K, L)$ and thus, by the same projection-restirction technique we have been employing, there is a natural map of fiber bundles over $S_{1,1}$,
\begin{equation*}
	S^{(3)}_{1,1} \to \hom( K \odot K \odot K, L) \to \hom(K, \hom( K \odot K, L)).
\end{equation*}
This allows use to stratify $S^{(3)}_{1,1}$ into two strata $S_{1,1,0}$ and $S_{1,1,1}$, based on the corank of the operator in $ \hom(K, \hom( K \odot K, L))$. Similarly, the normal bundle of $S_{1,1,1}$ is isomorphic to,
 \begin{equation*}
 \hom( K \odot K \odot K, L) 
\end{equation*}
and so we get a map of fiber bundles (over $S_{1,1,1}$),
\begin{equation*}
	S^{(4)}_{1,1,1} \to \hom( K \odot K \odot K \odot K, L) \to \hom(K, \hom( K \odot K \odot K, L)).
\end{equation*}
Again we can stratify according to the corank of the final operator, yielding two strata: $S_{1,1,1,0}$ and $S_{1,1,1,1}$. The later has codimension four. We can define the remaining strata as the intersections:
\begin{align*}
S_{[001], [01], 0} & := S_{[001], [01]}^{(3)} \cap S_{1,1,0}  \\
S_{[001], [01], 1} &:= S_{[001], [01]}^{(3)} \cap S_{1,1,1}  \\
S_{[001], [11], 0} &  := S_{[001], [11]}^{(3)} \cap S_{1,1,0} \\
S_{[001], [11], 1} & := S_{[001], [11]}^{(3)} \cap S_{1,1,1} \\
S_{[011], [001], 0} & := S_{[011], [001]}^{(3)} \cap S_{1,1,0} \\
S_{[011], [001], 1} & := S_{[011], [001]}^{(3)} \cap S_{1,1,1} \\
S_{[011], [100],0} & := S_{[011], [100]}^{(3)} \cap S_{1,1,0} \\
S_{[011], [100],1} & := S_{[011], [100]}^{(3)} \cap S_{1,1,1} \\
%%% The 4-Jets
S_{[001], [01], 1, 0}& := S_{[001], [01]}^{(4)} \cap S_{1,1,1,0} \\
S_{[001], [01], 1, 1} & := S_{[001], [01]}^{(4)} \cap S_{1,1,1,1} 
\end{align*}

Just as the Cerf theory strata were compatible with the Morse theory strata, the 3D strata listed in Table~\ref{Table3DStratificationOfJets} are compatible with 2D strata considered in the previous two sections. In particular if
$f: \Sigma \to \R^2$ is generic with respect to the 2D stratification, then the product map $F= (f, id): \Sigma \times I \to \R^2 \times I$ is generic with respect to the 3D stratification. Moreover, given a pair of generic maps, $f_0, f_1: \Sigma \to \R^2$, we may choose a path $F: \Sigma \times I \to \R^2 \times I$ from $f_0$ to $f_1$ which is a product map near the two ends of the interval. This map $F$ satisfies the conditions of the Relative Transversality Theorem~\ref{ThomTransversalityTheoremRelativeCorners}. Let $p:\R^2 \to I$ be the projection onto the time variable. Recall that $F$ is a path precisely if $pF(x,t) = t$. 

The Relative Transversality Theorem~\ref{ThomTransversalityTheoremRelativeCorners} ensures there is residual subset of maps whose jet sections are transverse to the above strata, and which agree with $F$ on a neighborhood of the ends $\Sigma \times \partial I$. Hence any non-empty subset of maps open in the Whitney $C^\infty$-topology contains such maps. The set of paths, however, does not form an open set in the Whitney $C^\infty$-topology. Instead we will be interested in studying those $F$ for which $\partial_t(pF)$ is positive. These form an open neighborhood of the paths, which can be seen by following an analogous argument to one given in Section~\ref{SectRudCerfTheory}.
First we consider the projection of jet bundles:
\begin{equation*}
	J^1(\Sigma \times I, \R^2 \times I) \to J^1(\Sigma \times I, I).
\end{equation*}
On  $\Sigma \times I$ there is a distinguished vector field which in local coordinates is given by $\partial_t$. This induces a function on the 1-jets $J^1(\Sigma \times I, I) \cong T^*(\Sigma \times I) \otimes TI$ as follows. Given a jet we may pair it with the vector field $\partial_t$ to obtain an element of $TI \cong I \times \R$. Projecting to $\R$ gives the desired function. The inverse image of the positive real numbers, an open set, consists of all those 1-jets in which the $t$-partial derivative is positive. The inverse image in $J^1(X \times I, \R \times I)$ is also open, and we denote this open set $\cO$.
% $\cO$ is defined on considering, on  $\Sigma \times I$, the distinguished vector field, which in local coordinates is given by $\partial_t$, where $\partial_t$ is the standard vector field on $I$.   The span of this gives us 1-dimensional  subbundle $E \subseteq T(\Sigma \times I)$, and thus we may consider the projection,
%\begin{equation*}
%	J^1(X\Sigma \times I, I) = T^*(\Sigma \times I) \otimes TI \to E^* \otimes TI.
%\end{equation*}
%The inverse image of the complement of the zero section of $ E^* \otimes TI$ gives the open set $\cO$. 
From this we may consider the subset of $C^\infty(X \times I, \R \times I)$ given by
\begin{equation*}
	M(\cO) = \{ F \; | \; j^1F(\Sigma \times I) \subseteq \cO \}.
\end{equation*}
This is open in the Whitney $C^1$-topology, and hence open in the Whitney $C^\infty$-topology. It consists of those maps $F$, as above, for which $\partial_t(pF)$ never vanishes.
% $F(x_0, x_1,t) = (f_t(x_0, x_1), g(x_0, x_1,t))$ such that the $t$-derivative of the last coordinate (i.e., $\partial_t g$) is non-zero at all points of $\Sigma \times I$. A path of functions, $F: \Sigma \times I \to \R \times I$, is a map in $M(\cO)$, so that we see that $M(\cO)$ is non-empty. 
Fix two maps $f_0, f_1: \Sigma \to \R^2$ which are generic with respect to the 2D stratification. 
We will restrict our attention to only those maps $F \in M(\cO)$ which agree with $(f_0, id)$ and $(f_1, id)$ on a neighborhood of $\Sigma \times \partial I$. 

\begin{definition}
Let $\Sigma$ be a surface and let $f_0, f_1: \Sigma \to \R^2$ be two maps which are generic with respect to the 2D stratification. 
A map $F: \Sigma \times I \to \R \times I$ such that $F \in M(\cO)$ and which agrees with $(f_0, id)$ and $(f_1, id)$ on a neighborhood of $\Sigma \times \partial I$  
 will be called {\em generic} if it lies in the residual (hence dense) subset of those $F$ whose induced jet sections $j^k F$ are transverse to the strata listed in Table~\ref{Table3DStratificationOfJets}. We will say that the generic map $F$ {\em connects} the generic maps $f_0$ and $f_1$.
\end{definition}

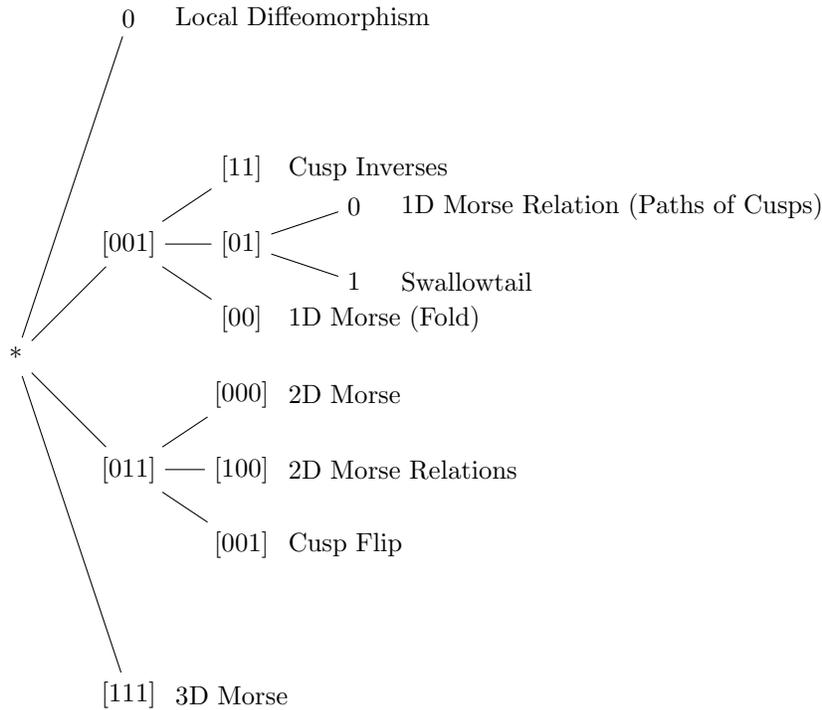
\begin{figure}[htb]
\begin{center}
\begin{tikzpicture}[
%grow via three points={% 
%one child at (0,1) and two children at (-1,1) and (1,1)}, grow = right] 
level distance={1.5cm}, 
level 1/.style={sibling distance=3cm}, 
level 2/.style={sibling distance=1cm}, 
level 3/.style={sibling distance=1cm}, grow = right] 

\node at (0,0) { *} 
	child { node {[111]} node  [anchor = west] at +(0.5,0) { \color{black} 3D Morse}}	
	child {node {[011]}
		child {node {[001]} node  [anchor = west] at +(0.5,0) {  \color{black} Cusp Flip}}
		child {node {[100]} node [anchor = west] at +(0.5,0) {  \color{black} 2D Morse Relations }}
		child {node {[000]} node  [anchor = west] at +(0.5,0) {  \color{black}  2D Morse  }}
		}
	child {node {[001]} 
		child {node {[00]} node  [anchor = west] at +(0.5,0) {  \color{black}  1D Morse (Fold) }}
		child {node {[01]}
			child {node {1} node  [anchor = west] at +(0.5,0) {  \color{black} Swallowtail}}
			child {node {0} node  [anchor = west] at  +(0.5,0)
				 {  \color{black} 1D Morse Relation (Paths of Cusps)}}
			}
		child {node {[11]} node  [anchor = west] at +(0.5,0) {  \color{black} Cusp Inverses}}
		}
	child {node {0} node [anchor = west] at +(0.5,0) {  \color{black} Local Diffeomorphism}};
	
\end{tikzpicture} 
\caption{The 3D Singularity Family Tree }
\index{singularity!3D family tree}
\label{fig:3DSingFamilyTree}
\end{center}
\end{figure}

\subsection{Derivation of Local Coordinates: 3D Case} \label{SectNormalCoord3D}

As before, let $\Sigma$ be a surface.
In the last section we have seen that a generic map  $F: \Sigma \times I \to \R \times I$ (connecting generic maps  $f_0, f_1: \Sigma \to \R^2$) has singularities of only eight types. Using the Thom-Boardman stratification of jet space, we may determine the local ring of each of these singularities, which are listed in Table~\ref{TableLocalRingof3DSingularities}, along with their codimensions. In this section, under the hypothesis that $F$ is generic, we will derive local coordinates around each of these singularities. 

\begin{table}[ht]
\begin{center}
\begin{tabular}{|c|c| c| c|} \hline
Singularity Stratum & Name & Codim. & Local Ring \\ \hline \hline
$S_0$ & Local Diffeomorphism & 0 & $\R$ \\ \hline 
$S_{[001], 0}$  & 1D Morse (Fold) & 1 & $\R[ x] / (x^2)$ \\
$S_{[011], 0}$  & 2D Morse & 2 & $\R[ x] / (x^2)$\\
$S_{[011], [100],0}$  & 2D Morse Relation  & 3 &$\R[ x] / (x^2)$ \\ \hline
$S_{[001], [11], 0}$  & Cusp Inversion & 3 &$\R[ x] / (x^3)$ \\
$S_{[001], [01], 0}$  & 1D Morse Rel. (Paths of Cusps) & 2 & $\R[ x] / (x^3)$ \\
$S_{[011], [001], 0}$  & Cusp Flip & 3 & $\R[ x] / (x^3)$ \\ \hline
$S_{[001], [01], 1, 0}$  & Swallowtail & 3 & $\R[ x] / (x^4)$ \\ \hline
\end{tabular}
\end{center}
\caption{The Local Rings of the 3D Singularities}
\label{TableLocalRingof3DSingularities}
\end{table}%

The stratification of jet space that we use in the 3-dimensional case considered here was derived by incorporating information about the composition of $F$ with the projections $\R^2 \times I \to \R \times I $ and $\R \times I \to I$. It is natural for us to try to consider those coordinates changes which preserve some of this relevant structure. For all but the $S_0$ singularity in Table~\ref{TableLocalRingof3DSingularities} (pg.~\pageref{TableLocalRingof3DSingularities}), there are two associated subspace $K_1 \subseteq T_x( \Sigma \times I)$ and $K_2 \subseteq T_x( \Sigma \times I) $, where $x \in \Sigma \times I$ is a critical point of the given singularity type. The coordinate changes that we will consider will preserve these spaces in a sense that we will make precise later. 

The singularities of type $S_0$ are not actually singularities at all. In a neighborhood of such a point $F$ is a diffeomorphism, and so we call these points {\em local diffeomorphism points}. In this case $F$ itself provides standard coordinates around such a point. The remaining singularities fall into two kinds, $S_{[001]}$-singularities and $S_{[011]}$-singularities, and our initial coordinates will only depend on this division.

For the remaining cases our derivation of local coordinates will generally consist of four stages. First we use the implicit function theorem to obtain initial coordinates and we analyze which conditions $F$ must satisfy in order that the sections $j^kF$ be transverse to the strata of Table~\ref{Table3DStratificationOfJets}. Second, using Lemma \ref{LmaCoordinateChange} below, we will obtain a new, simplified coordinate system in which the higher partial derivatives of $F$ take convenient values. Next, using the local ring structure induced by $F$, we will employ the Generalized Malgrange Preparation Theorem~\ref{GeneralizedMalgrangePreparationTheorem} and Corollary \ref{GeneralizedMalgrangePreparationTheoremCorollary} to obtain a functional equation, just as we have done previously in our discussion of Cerf theory. Finally, using this functional equation and the conditions on the partial derivatives of $F$ that we obtained in step two, we will be able to derive a final coordinate system in which $F$ has a simple form.  

The results are summarized in the seven Propositions below: Prop.~\ref{PropCuspInversionCoord}, Prop.~\ref{PropPathofCuspCoordinates}, Prop.~\ref{PropSwallowtailCoord},  Prop.~\ref{FoldCoordProp}, Prop.~\ref{Prop2DMorseCoord}, Prop.~\ref{Prop2DMorseRelationCoord}, and Prop.~\ref{PropCuspFlipCoord}. The proofs of these propositions share many commonalities. The extent to which the proofs are similar is precisely reflected by the proximity of the relevant strata in the stratification of the previous section. This can be visualized by the `singularity family tree' of Figure~\ref{fig:3DSingFamilyTree} (pg.~\pageref{fig:3DSingFamilyTree}). For example the four proofs of Prop.~\ref{PropCuspInversionCoord}, Prop.~\ref{PropSwallowtailCoord},  Prop.~\ref{PropPathofCuspCoordinates}, and Prop.~\ref{FoldCoordProp} (i.e. those corresponding to $S_{[001]}$-singularities) all begin in the same way. Rather than repeat the beginnings of these proofs four times over, we will first provide preliminaries for the $S_{[001]}$-singularities. Similarly we provide a preliminary discussion of the $S_{[011]}$-singularities. 

The singularities of codimension strictly less than three (cf. Table~\ref{TableLocalRingof3DSingularities} on pg.~\pageref{TableLocalRingof3DSingularities}) also deserve a few remarks before we state the propositions.  
%We now turn to the lower codimension singularities. 
Let $S_{[001], [01], 0}(F)$, $S_{[011], 0}(F)$, and $S_{[001], 0}(F)$ denote those points in $\Sigma \times I$ in which $F$ obtains an $S_{[001], [01], 0}$-, $S_{[011], 0}$-, or $S_{[001], 0}$-singularity, respectively. By transversality, each of these is a smooth submanifold of $\Sigma \times I$ of dimension one, one, or two, respectively. Moreover, in the $S_{[001], [01], 0}$ and $S_{[011], 0}$ cases, the projections
\begin{align*}
	F: &S_{[001], [01], 0}(F) \to \R^2 \times I \stackrel{p^2}{\to} I \\
	F: & S_{[011], 0}(F)  \to \R^2 \times I \stackrel{p^2}{\to} I, 
\end{align*} 
are local diffeomorphisms. Thus in both of these cases the image of this submanifold under $F$ is (locally) the image of the graph of some function $\gamma= (\gamma_0, \gamma_1) : I \to \R^2$. Similarly, in the  $S_{[001], 0}$ case the projection,
\begin{equation*}
		F: S_{[001], 0}(F) \to \R^2 \times I \stackrel{p}{\to} \R\times  I 
\end{equation*}
is a local diffeomorphism and hence the image under $F$ of this submanifold is (locally) the graph of a function $\gamma: \R \times I \to \R$. 
%The $S_{[001], 0}$ case is the easiest of the three and we begin there.

We now state the results of this section.
\begin{proposition}[Cusp Inversion (see Figure~\ref{CuspInversionGraphicFig} on pg.~\pageref{CuspInversionGraphicFig})] \label{PropCuspInversionCoord}
\index{cusp inversion singularity}
\index{singularity!cusp inversion}
Let $x \in \Sigma \times I$ be a point where $F$ obtains an $S_{[001], [11], 0}$ (i.e., Cusp Inversion) singularity. Then there exist local coordinates $(\overline x_0, \overline x_1, \overline x_2)$ for $\Sigma \times I$ centered at $x$, and local coordinates $(\overline y_0, \overline y_1, y_2)$ for $\R^2 \times I$ centered at $y= F(x)$ (here $y_2$ is the standard coordinate of $I$), such that in these coordinates $F$ takes the form:
\begin{align*}
	F^* \overline y_0 & = \overline x_0^3 \pm \overline x_2 \overline x_0 \pm \overline x_1^2  \overline x_0 \\
	F^* \overline y_1 & = \overline x_1 \\
	F^* y_2 & = \overline x_2.
\end{align*}
Moreover these coordinate changes preserve the orientations of $\R^2 \times I$, $\R \times I$ and $I$.
\end{proposition}

\begin{proposition}[1D Morse Relation (Paths of Cusps) (see Figure~\ref{CuspAnd2DMorseIn3DFig} on pg.~\pageref{CuspAnd2DMorseIn3DFig})] \label{PropPathofCuspCoordinates}
\index{cusp singularity}
\index{singularity!cusp}
\index{1D Morse relation singularity}
\index{1D Morse relation singularity|see{cusp singularity}}
	Let $\Sigma$ be a surface, $F: \Sigma \times I \to \R^2 \times I$ be generic, and $x \in \Sigma \times I$ be an  $S_{[001],[01], 0}$ singularity point. Then there exists a 1-dimensional submanifold  $Y \subseteq \Sigma \times I $ containing $x$ and consisting entirely of $S_{[001],[01], 0}$ singularities,  a neighborhood $U\subseteq I$ of $p^2(y) = p^2F(x)$, and 
 a function $\gamma = (\gamma_0, \gamma_1)  : U \to \R^2$ whose graph $\{ (\gamma_0(y_2), \gamma_1(y_2), y_2) \; | \; y_2 \in U \}$ consists precisely of the image of $Y$ under $F$. Moreover there exist coordinates $(\overline x_0, \overline x_1, x_2)$ for $\Sigma \times I$ centered at $x$, (here $x_2$ is the standard coordinate for $I$) and coordinates $(\overline y_0, \overline y_1, y_2)$ for $\R^2 \times I$, centered at $y$ (here $y_2$ is the standard coordinate for $I$) and  such that in these coordinates $Y = \{ x_0 = 0, x_1 = 0 \}$ (the $x_2$-axis) and  $F$ has the following normal form:
\begin{align*}
	F^* y_0 &=  \overline{x_0}^3  \pm \overline x_0 \overline x_1 +  \gamma_0( x_2)\\
	F^*y_1 &= x_1 + \gamma_1(x_2)\\
	F^*y_2 &= x_2 \\
\end{align*}
These coordinate changes can be taken to be compatible with the orientations of $\R^2 \times I$, $\R \times I$ and $I$. 
\end{proposition}

\begin{proposition}[Swallowtail (see Figure~\ref{SwallowtailGraphicFig} on pg.~\pageref{SwallowtailGraphicFig})] \label{PropSwallowtailCoord}
\index{Swallowtail singularity}
\index{singularity!swallowtail}
Let $x \in \Sigma \times I$ be a point where $F$ obtains an $S_{[001], [01], 1, 0}$ (i.e., Swallowtail) singularity. Then there exist local coordinates $(\overline x_0, \overline x_1, \overline x_2)$ for $\Sigma \times I$ centered at $x$, and local coordinates $(\overline y_0, \overline y_1, \overline y_2)$ for $\R^2 \times I$ centered at $y= F(x)$, such that in these coordinates $F$ takes the form:
\begin{align*}
	F^* \overline y_0 & = \pm \overline x_0^4 +  \overline x_1 \overline x_0 \pm \overline x_2  \overline x_0^2 \\
	F^* \overline y_1 & = \overline x_1 \\
	F^* \overline y_2 & = \overline x_2.
\end{align*}
Moreover these coordinate changes preserve the orientations of $\R^2 \times I$, $\R \times I$ and $I$.
\end{proposition}

\begin{proposition}[1D Morse (Fold) (see Figures~\ref{FoldSingularitiesIn3DFig} and \ref{FoldGraphicsFig} on pg.~\pageref{FoldSingularitiesIn3DFig} and pg.~\pageref{FoldGraphicsFig}, respectively)] \label{FoldCoordProp}
\index{fold singularity}
\index{singularity!fold}
\index{1D Morse singularity}
\index{1D Morse singularity|see{fold singualrity}}
	Let $\Sigma$ be a surface, $F: \Sigma \times I \to \R^2 \times I$ be generic, and $x \in \Sigma \times I$ be an  $S_{[001], 0}$ singularity point. Then there exists a 2-dimensional submanifold  $Y \subseteq \Sigma \times I $ containing $x$ and consisting entirely of $S_{[001], 0}$ singularities,  a neighborhood $U\subseteq \R \times I$ of $p(y) = pF(x)$, and 
 a function $\gamma  : U \to \R$ whose graph $\{ (\gamma(y_1,y_2), y_1, y_2) \; | \; (y_1,y_2) \in U \}$ consists precisely of the image of $Y$ under $F$. Moreover there exist coordinates $(\overline x_0, \overline x_1, x_2)$ for $\Sigma \times I$ centered at $x$, (here $x_2$ is the standard coordinate for $I$)  such that in these coordinates $Y = \{ \overline x_0 = 0 \}$ (the $\overline x_1$-$x_2$-plane) and  $F$ has the following normal form:
\begin{align*}
	F^* y_0 &=  \pm \overline{x_0}^2  +  \gamma(\overline x_1, \overline x_2)\\
	F^*y_1 &= \overline x_1\\
	F^*y_2 &= x_2 \\
\end{align*}
Here $(y_0, y_1, y_2)$ are the standard coordinates for $\R^2 \times I$ centered at $y = F(x)$. 
\end{proposition}

\begin{proposition}[2D Morse (see Figure~\ref{CuspAnd2DMorseIn3DFig} on pg.~\pageref{CuspAnd2DMorseIn3DFig})] \label{Prop2DMorseCoord}
\index{2D Morse singularity}
\index{singularity!2D Morse}
	Let $\Sigma$ be a surface, $F: \Sigma \times I \to \R^2 \times I$ be generic, and  $x \in \Sigma \times I$ be an  $S_{[011], [000]}$ singularity point. Then there exists a one-dimensional submanifold  $Y \subseteq \Sigma \times I $ containing $x$ and consisting entirely of $S_{[011], [000]}$ singularities,  a neighborhood $U\subseteq I$ of $p^2(y) = p^2F(x)$, and 
 a function $\gamma = (\gamma_0, \gamma_1)  : U \to \R^2$ whose graph $\{ (\gamma_0(y_2), \gamma_1(y_2), y_2) \; | \; y_2 \in U \}$ consists precisely of the image of $Y$ under $F$. Moreover there exist coordinates $(\overline x_0, \overline x_1, x_2)$ for $\Sigma \times I$ centered at $x$, (here $x_2$ is the standard coordinate for $I$) and coordinates $(\overline y_0, \overline y_1, y_2)$ for $\R^2 \times I$, centered at $y$ (here $y_2$ is the standard coordinate for $I$) and  such that in these coordinates $Y = \{ x_0 = 0, x_1 = 0 \}$ (the $x_2$-axis) and  $F$ has the following normal form:
\begin{align*}
	F^* y_0 &=  x_0 +  \gamma_0( x_2) \\
	F^*y_1 &= \pm x_0^2 \pm x_1^2 + \gamma_1(x_2)\\
	F^*y_2 &= x_2 \\
\end{align*}
These coordinate changes can be taken to be compatible with the orientations of $\R^2 \times I$, $\R \times I$ and $I$. 
\end{proposition}

\begin{proposition}[2D Morse Relation (see Figure~\ref{2DMorseRelationFig} on pg.~\pageref{2DMorseRelationFig})] \label{Prop2DMorseRelationCoord}
\index{2D Morse relation singularity}
\index{singularity!2D Morse relation}
Let $x \in \Sigma \times I$ be a point where $F$ obtains an $S_{[011], [100], 0}$ (i.e., 2D Morse Relation) singularity. Then there exist local coordinates $(\overline x_0, \overline x_1, \overline x_2)$ for $\Sigma \times I$ centered at $x$, and local coordinates $(\overline y_0, \overline y_1, \overline y_2)$ for $\R^2 \times I$ centered at $y= F(x)$, such that in these coordinates $F$ takes the form:
\begin{align*}
	F^* \overline y_0 & =  \overline x_0 \\
	F^* \overline y_1 & = \pm \overline x_1^2 \pm \overline x_0 \overline x_2 \pm \overline x_0^3  \\
	F^* \overline y_2 & = \overline x_2.
\end{align*}
Moreover these coordinate changes preserve the orientations of $\R^2 \times I$, $\R \times I$ and $I$.
\end{proposition}

\begin{proposition}[Cusp Flip (see Figure~\ref{CuspFlipsGraphicFig} on pg.~\pageref{CuspFlipsGraphicFig})] \label{PropCuspFlipCoord}
\index{cusp flip singularity}
\index{singularity!cusp flip}
Let $x \in \Sigma \times I$ be a point where $F$ obtains an $S_{[011], [001], 0}$ (i.e., Cusp Flip) singularity. Then there exist local coordinates $(\overline x_0, \overline x_1, \overline x_2)$ for $\Sigma \times I$ centered at $x$, and local coordinates $(\overline y_0, \overline y_1, \overline y_2)$ for $\R^2 \times I$ centered at $y= F(x)$, such that in these coordinates $F$ takes the form:
\begin{align*}
	F^* \overline y_0 & = \overline x_0 \\
	F^* \overline y_1 & = \overline x_1^3 \pm \overline x_1 \overline x_0 \pm \overline x_1^2  \overline x_2 \\
	F^* \overline y_2 & = \overline x_2.
\end{align*}
Moreover these coordinate changes preserve the orientations of $\R^2 \times I$, $\R \times I$ and $I$.
\end{proposition}

Let $(y_0, y_1, y_2)$ denote the standard coordinates of $\R^2 \times I$ (here $y_2$ is the standard coordinate of $I$). Let $x \in \Sigma \times I$ be a point where $F$ obtains one of the seven remaining singularities. The implicit function theorem then ensures that there exist coordinates $(x_0, x_1, x_2)$ of $\Sigma \times I$, centered at $x$, in which $F$ has the form listed in Table~\ref{InitialChoiceofCoordTable}. In both cases $x_2$ will coincide with the standard coordinate on $I$. Moreover in each of these coordinates the function $h$, defined in Table~\ref{InitialChoiceofCoordTable}, satisfies $dh = 0$.

\begin{table}[ht]
\begin{center}
\begin{tabular}{|c|c|c|c| c|} \hline
Type & $F$ & $dF$ & $K_1$ & $K_2$ \\ \hline \hline
$S_{[001]}$ & 
$\begin{array}{r c l}
	& & \\
	F^*y_0  & = & h(x_0, x_1, x_2) \\
	F^*y_1 & = & x_1 \\
	F^* y_2 & = & x_2 \\
	&& 
\end{array}$
& $\left( \begin{array}{c c c}
	0 & 0 & 0 \\ 0 & 1 & 0 \\ 0 & 0 & 1
\end{array} \right)$ & $\text{Span} \{ \partial_{x_0} \}$ & $\text{Span} \{ \partial_{x_0}, \partial_{x_1} \}$ \\ \hline
$S_{[011]}$ & 
$\begin{array}{r c l}
	& & \\
	F^*y_0  & = & x_0 \\
	F^*y_1 & = & h(x_0, x_1, x_2) \\
	F^* y_2 & = & x_2 \\
	&&
\end{array}$
& $\left( \begin{array}{c c c}
	1 & 0 & 0 \\ 0 & 0 & 0 \\ 0 & 0 & 1
\end{array} \right)$ & $\text{Span} \{ \partial_{x_1} \}$ & $\text{Span} \{ \partial_{x_0}, \partial_{x_1} \}$ \\ \hline
\end{tabular}
\end{center}
\caption{Initial Choice of Coordinates}
\label{InitialChoiceofCoordTable}
\end{table}%

Recall that over the strata $S_{[001]}$ and $S_{[011]}$ there are vector bundles $K_1 \subseteq K_2 \subseteq T(\Sigma \times I)$, defined as the kernel of $dF$ and $d(pF)$, respectively. At the point $x$ these give a pair of subspaces $K_1$ and $K_2$ of $T_x(\Sigma \times I)$, which are also listed in Table~\ref{InitialChoiceofCoordTable}. As we saw in the previous section, these subspaces play a key role in the further stratification of the jet space. All the coordinate changes that we will consider will preserve these spaces in the sense that, in the $S_{[001]}$-case, if $(\overline x_1, \overline x_2, \overline x_3)$ are our final coordinates, then we have $K_1 = \text{Span} \{ \partial_{\overline x_0} \}$ and  $K_2 = \text{Span} \{ \partial_{\overline x_0}, \partial_{\overline x_1} \}$ at $x$. Similarly, in the $S_{[011]}$-case, if $(\overline x_1, \overline x_2, \overline x_3)$ are our final coordinates, then we have $K_1 = \text{Span} \{ \partial_{\overline x_1} \}$ and  $K_2 = \text{Span} \{ \partial_{\overline x_0}, \partial_{\overline x_1} \}$ at $x$. %We will precede with the codimension 3 singularities first. 

We omit the proofs of the following two trivial lemmas. 

\begin{lemma} \label{LmaCoordinateChange}
Let $h: \R^3 \to \R$ be a function such that $dh =0$ at the origin. Let $(x,y,z)$ denote the standard coordinates on $\R^3$, and consider the following coordinate change:
\begin{align*}
	u & = x + A y + Bz\\
	v & = y + Cz \\
	w & =  z
\end{align*}
Then the following partial derivatives of $h$ at the origin with respect to $u, v, w$, are given by:
\begin{align*}
	h_{uu} &= h_{xx} \\
	h_{uv} &= h_{xy} - A h_{xx} \\
	h_{uw} &= h_{xz} + (- B + AC)h_{xx} - Ch_{xy} \\
	h_{uuv} &= h_{xxy} - A h_{xxx} \\
 	h_{uuw} &= h_{xxz} + (- B + AC)h_{xxx} - Ch_{xxy} \\
	h_{uuuv} &= h_{xxxy} - A h_{xxxx} \\
	h_{uuuw} &= h_{xxxz}  + (- B + AC)h_{xxxx} - Ch_{xxxy}
\end{align*}
where $h_{\alpha \beta}$,  $h_{\alpha \beta \gamma}$, and $h_{\alpha \beta \gamma \delta}$ with $\alpha, \beta, \gamma, \delta \in \{ x, y, z\}$ denote the respective partial derivatives of $h$ at the origin. Moreover, we have equality of the following tangent vectors at the origin:
\begin{align*}
	\partial_u & =  \partial_x \\
	\partial_v & =  \partial_y -A \partial_x \\
	\partial_w & = \partial_z + (- B + AC) \partial_x -C \partial_y  
\end{align*}
\end{lemma}

\begin{lemma} \label{LmaCoordinateChange2}
Let $h: \R^3 \to \R$ be a function such that $dh =0$ at the origin. Let $(x,y,z)$ denote the standard coordinates on $\R^3$, and consider the following coordinate change:
\begin{align*}
	u & = x + A x^2\\
	v & = y\\
	w & =  z
\end{align*}
Then the following partial derivatives of $h$ at the origin with respect to $u, v, w$, are given by:
\begin{align*}
	h_{uu} &= h_{xx} -2 A h_x &	h_{uuw} &= h_{xxz} - 2A h_{xz}\\
	h_{uv} &= h_{xy}  &  	h_{uuuv} &= h_{xxxy} -4A h_{xxy} + 12 h_{xy} \\
	h_{uw} &= h_{xz}  & 	h_{uuuw} &= h_{xxxz} -4A h_{xxz} + 12 h_{xz}  \\
	h_{uuv} &= h_{xxy} - 2A h_{xy} 
\end{align*}
where $h_{\alpha \beta}$,  $h_{\alpha \beta \gamma}$, and $h_{\alpha \beta \gamma \delta}$ with $\alpha, \beta, \gamma, \delta \in \{ x, y, z\}$ denote the respective partial derivatives of $h$ at the origin. Moreover, we have equality of the following tangent vectors:
\begin{align*}
	\partial_u & =  \frac{1}{1 + 2A x}\partial_x \\
	\partial_v & =  \partial_y  \\
	\partial_w & = \partial_z   
\end{align*}
\end{lemma}

\begin{proof}[Preliminaries for Propositions \ref{PropCuspInversionCoord}, \ref{PropPathofCuspCoordinates}, \ref{PropSwallowtailCoord},  and \ref{FoldCoordProp}] The implicit function theorem guarantees that there exist coordinates with the properties listed in Table~\ref{InitialChoiceofCoordTable}, i.e., 
	\begin{align*}
		F^* y_0 & = h(x_0, x_1, x_2) \\
		F^* y_1 & = x_1 \\
		F^* y_2 & = x_2.
	\end{align*}
	In the cases corresponding to singularities of codimension three (Prop.~\ref{PropCuspInversionCoord} and Prop.~\ref{PropSwallowtailCoord}), the equation $dh=0$ holds at the origin and we may expand $h$ in the variables $x_0, x_1, x_2$. We have,
	\begin{align*}
		h(x_1, x_2, x_3) = &\sum_{0 \leq i \leq j \leq 2} \alpha_{ij} x_i x_j + \sum_{0 \leq i \leq j \leq k \leq 2 } \beta_{ijk} x_i x_j x_k \\
		& + \sum_{0 \leq i \leq j \leq k \leq \ell 2} \gamma_{ijk\ell} x_ i x_j x_k x_\ell + O(|x|^5).
	\end{align*}
	For the cases corresponding to singularities of codimension less than three (Prop.~\ref{PropPathofCuspCoordinates} and Prop.~\ref{FoldCoordProp}) then the equation $dh=0$ holds along the entire manifold $Y$. In the case of Prop.~\ref{PropPathofCuspCoordinates}, we may arrange so that in these coordinates $Y = \{ (0, \gamma_1(x_2), x_2) \}$.  In the case of Prop.~\ref{FoldCoordProp}) we may arrange for $Y = \{ x_0 = 0 \}$. 
	
	In all cases $j^1F$ is transverse to the stratum $S_{[001]}$, which means that the following map is surjective:
		\begin{align*}
			d(j^1F): T_x(\Sigma \times I) \to &\;  T_{j^1F(x)}J^1(\Sigma \times I, \R^2 \times I) / T_{j^1F(x)}S_{[001]}  \\
			& \quad \cong \nu_{j^1F(x)} S_{[001]} \cong \hom(K_1, L) 
		\end{align*}
	For Prop.~\ref{PropPathofCuspCoordinates} and Prop.~\ref{FoldCoordProp} this surjectivity holds along the entire manifold $Y$, not just at the origin. At this stage the proofs of Prop.~\ref{PropCuspInversionCoord}, Prop.~\ref{PropPathofCuspCoordinates}, Prop.~\ref{PropSwallowtailCoord},  and Prop~.\ref{FoldCoordProp} begin to diverge. 
\renewcommand{\qedsymbol}{}  % remove the QED box symbol for this proof.	
\end{proof}

\begin{proof}[Proof of Prop.~\ref{PropCuspInversionCoord}] 
%Again, the implicit function theorem guarantees that there exist coordinates with the properties listed in Table~\ref{InitialChoiceofCoordTable}, i.e., 
%\begin{align*}
%	F^* y_0 & = h(x_0, x_1, x_2) \\
%	F^* y_1 & = x_1 \\
%	F^* y_2 & = x_2.
%\end{align*}
%where $dh = 0$. Again we can expand $h$ in the variables $x_0, x_1, x_2$ to get,
%\begin{align*}
%	h(x_1, x_2, x_3) = &\sum_{0 \leq i \leq j \leq 2} \alpha_{ij} x_i x_j + \sum_{0 \leq i \leq j \leq k \leq 2 } \beta_{ijk} x_i x_j x_k  + O(|x|^4).
%\end{align*}
Since $x$ is an $S_{[001], [11], 0}$-singularity, we know that $\alpha_{00} = 0$, $\alpha_{01}= 0$, and $\beta_{000} \neq 0$.
In addition to $j^1F$ being transverse to $S_{[001]}$, we have $j^2F$ transverse to $S_{[001], [11]}$. This later means the following map is surjective at $x$:
%By assumption, $j^kF$ is transverse to the strata $S_{[001]}$ and $S_{[001], [11]}$ for $k = 1$ and 2, respectively. Thus the following maps are surjective at $x$:
\begin{align*}
%	d(j^1F): T_x(\Sigma \times I) \to &\;  T_{j^1F(x)}J^1(\Sigma \times I, \R^2 \times I) / T_{j^1F(x)}S_{[001]}  \\
%	& \quad \cong \nu_{j^1F(x)} S_{[001]} \cong \hom(K_1, L) \\
	d(j^2F): T_x(\Sigma \times I) \to & \; T_{j^2F(x)}J^2(\Sigma \times I, \R^2 \times I) / T_{j^2F(x)}S_{[001], [11]}  \\
	&\quad \cong \nu_{j^2F(x)} S_{[001], [11]} \cong \hom(K_1, L) \oplus \hom(K_1 \cdot K_2, L)
\end{align*}
Given that $\alpha_{00} = \alpha_{01} = 0$ and $\beta_{000} = 0$, the transversality assumptions are equivalent to the condition that $\alpha_{02} \neq 0$ and that the following matrix is non-degenerate:
\begin{equation*}
\left( \begin{array}{ccc}
0& 0& \alpha_{02} \\
\beta_{000} & \beta_{001} & \beta_{002} \\
\beta_{001} & \beta_{011} & \beta_{012}
\end{array} \right)
\end{equation*}
 Using Lemma \ref{LmaCoordinateChange} 
 %%%% (1) eliminate \beta_001 and \beta_002
 %%%% (2) eliminate \beta_012
 we may find acceptable coordinates $(x_0', x_1', x_2')$ and $(y_0, y_1', y_2)$ in which the above matrix becomes:
 \begin{equation*}
\left( \begin{array}{ccc}
0& 0& \pm r \\
\pm s & 0 & 0 \\
 0 & \pm t & 0
\end{array} \right)
\end{equation*}
for some positive numbers $r,s,t$. By scaling $x_0$ appropriately, we may assume that $s= 1$.

The local ring of this singularity is $\cR_F \cong \R[x_0']/( (x_0')^3)$ and so by the generalized Malgrange Preparation Theorem there exist functions smooth $a,b,c$ on $\R^2 \times I$ such that 
\begin{equation*}
	(x_0')^3 = F^*a + F^*b \cdot x'_0 + F^*c \cdot (x'_0)^2.
\end{equation*}
By suitably changing the coordinate $x_0'$, we may assume that $c \equiv 0$. Moreover, we have $\frac{\partial a}{\partial y_0}(0) \neq 0$, $\frac{\partial a}{\partial y_1}(0) = 0$, and 
\begin{equation*}
	b(y_0, y_1', y_2) = \pm y_2 + g(y_0, y_1, y_2), 
\end{equation*}
where the function $g$ has the form,
\begin{equation*}
	g(y_0, y_1, y_2) = \pm t \cdot y_1^2 + O(|y|^3).
\end{equation*}
Thus the following is a valid pair of coordinate changes in  sufficiently small neighborhoods of $x$ and $y$: 
\begin{align*}
	\overline x_0 &= \pm x_0' & 	\overline y_0 &= a(y_0, y_1, y_2) \\
	\overline x_1 &= \sqrt{F^*g} & 	\overline y_1 &=  \sqrt{g(y_0, y_1, y_2) } \\
	\overline x_2 &= x_2' & \overline y_2 &= y_2
\end{align*}
In these coordinates $F$ has the desired form.
\end{proof}

\begin{proof}[Proof of Prop.~\ref{PropPathofCuspCoordinates}]
%Again, the implicit function theorem guarantees that there exist coordinates around $x$ with the properties listed in Table~\ref{InitialChoiceofCoordTable}, i.e., 
%\begin{align*}
%	F^* y_0 & = h(x_0, x_1, x_2) \\
%	F^* y_1 & = x_1 \\
%	F^* y_2 & = x_2.
%\end{align*}
%where $dh = 0$ along $Y$. Moreover, we may arrange so that in these coordinates $Y = \{ (0, \gamma_1(x_2), x_2) \}$. 
 Since this is an  $S_{[001], [01], 0}$ singularity, we have $h_{x_0 x_0} = 0$, $h_{x_0 x_1} \neq 0$, and $h_{x_0 x_0 x_0} \neq 0$ along $Y$. 
Next, we preform the following coordinate change:
\begin{align*}
	x_0' &= x_0 & 					y_0' &= y_0 - \gamma_0(y_2)  \\
	x_1' &=x_1 - \gamma_1(x_2) & 	y_1' &=  y_1 - \gamma_1(y_2)  \\
	x_2' &= x_2 &					y_2' &= y_2
\end{align*}
In these coordinates, $Y$ consists of the $x_2$-axis and its image under $F$ consists of the $y_2$-axis. We introduce the following new function, 
\begin{equation*}
	g(x_0', x_1', x_2') := h(x_0', x_1' + \gamma_1(x_2'), x_2') - \gamma_0( x_2').
\end{equation*}
Notice that $g = F^* y_0'$. Our assumptions on $h$ guarantee that
\begin{align*}
g(0, 0, x_2') &\equiv 0 & \frac{\partial g}{\partial x_0'}(0, 0, x_2') &\equiv 0 &  \frac{\partial^2 g}{\partial (x_0')^2}(0, 0, x_2') & \equiv 0 \\
\frac{\partial^2 g}{\partial x_0' \partial x_1'}(0, 0, x_2') & \neq 0 & & & \frac{\partial^3 g}{\partial (x_0')^3}(0, 0, x_2') & \neq 0
\end{align*}
The local ring of this singularity is $\R[x_0'] / ((x_0')^3)$ and so by the Malgrange Preparation Theorem, there exist smooth functions $a, b, c$ on $\R^2 \times I$ such that
\begin{equation*}
	(x_0')^3 = F^*a + F^*b \cdot x_0' + F^*c \cdot (x_0')^2.
\end{equation*}
We may rewrite this as
\begin{equation*}
	(x_0 - \frac 1 3 F^* c)^3 + F^* \tilde b \cdot ( x_0 - \frac 1 3 F^* c) = F^*\tilde a
\end{equation*}
for some new functions $\tilde a, \tilde b$ on $\R^2 \times I$. Collecting terms we see that,
\begin{align*}
	\frac{\partial \tilde a}{\partial y_0 }(0,0, y_2) & \neq 0 & \frac{\partial \tilde b}{\partial y_0 }(0,0, y_2) &= 0 & \frac{\partial \tilde b}{\partial y_1 }(0,0, y_2) & \neq 0. 
\end{align*}
Thus the following is a valid pair of coordinate changes:
\begin{align*}
	\overline x_0 &= \pm(x_0' - \frac 1 3 F^*c)  & \overline y_0 
		&= \tilde a( y_0' + \gamma(y_2'), y_1', y_2') + \gamma(y_2') \\
	\overline x_1 &= F^* \tilde b& 	\overline y_1 &= \tilde b(y_0', y_1', y_2') \\
	\overline x_2 &=x_2'        & \overline y_2 &= y_2'
\end{align*}
Notice that $\overline y_2 = y_2' = y_2$. In these coordinates $F$ has the desired form. 
\end{proof}

\begin{proof}[Proof of Prop.~\ref{PropSwallowtailCoord}]
%As mentioned earlier, the implicit function theorem guarantees that there exist coordinates with the properties listed in Table~\ref{InitialChoiceofCoordTable}, i.e., 
%\begin{align*}
%	F^* y_0 & = h(x_0, x_1, x_2) \\
%	F^* y_1 & = x_1 \\
%	F^* y_2 & = x_2.
%\end{align*}
%where $dh = 0$ at the origin. Let us expand $h$ in the variables $x_0, x_1, x_2$. We have,
%\begin{align*}
%	h(x_1, x_2, x_3) = &\sum_{0 \leq i \leq j \leq 2} \alpha_{ij} x_i x_j + \sum_{0 \leq i \leq j \leq k \leq 2 } \beta_{ijk} x_i x_j x_k \\
%	& + \sum_{0 \leq i \leq j \leq k \leq \ell 2} \gamma_{ijk\ell} x_ i x_j x_k x_\ell + O(|x|^5).
%\end{align*}
Since $x$ is an $S_{[001], [01], 1,0}$-singularity, we know that $\alpha_{00} = 0$, $\alpha_{01} \neq 0$, $\beta_{000} = 0$, and $\gamma_{0000} \neq 0$. 
By assumption, we additionally have $j^2F$ and $j^3F$ transverse to the strata  $S_{[001], [01]}$ and $S_{[001], [01], 1}$, respectively. Thus we also have that the following maps are surjective at $x$:
%By assumption, $j^kF$ is transverse to the strata $S_{[001]}$, $S_{[001], [01]}$, and $S_{[001], [01], 1}$ for $k = 1$, 2, and 3, respectively. Thus the following maps are surjective at $x$:
\begin{align*}
%	d(j^1F): T_x(\Sigma \times I) \to &\;  T_{j^1F(x)}J^1(\Sigma \times I, \R^2 \times I) / T_{j^1F(x)}S_{[001]}  \\
%	& \quad \cong \nu_{j^1F(x)} S_{[001]} \cong \hom(K_1, L) \\
	d(j^2F): T_x(\Sigma \times I) \to & \; T_{j^2F(x)}J^2(\Sigma \times I, \R^2 \times I) / T_{j^2F(x)}S_{[001], [01]}  \\
	&\quad \cong \nu_{j^2F(x)} S_{[001], [01]} \cong \hom(K_1, L) \oplus \hom(K_1 \cdot K_1, L)\\
	d(j^3F): T_x(\Sigma \times I) \to & \; T_{j^3F(x)}J^3(\Sigma \times I, \R^2 \times I) / T_{j^3F(x)}S_{[001], [01],1}  \\
	&\quad \cong \nu_{j^3F(x)} S_{[001], [01],1} \\
	&\quad \cong \hom(K_1, L) \oplus \hom(K_1 \cdot K_1, L) \oplus \hom(K_1 \cdot K_1 \cdot K_1, L)
\end{align*}
Given that $\alpha_{00} = 0$, $\alpha_{01} \neq 0$, $\beta_{000} = 0$, and $\gamma_{0000} \neq 0$, our transversality assumptions are equivalent to the condition that the $2 \times 2$ matrix,
\begin{equation*}
	\left( \begin{array}{cc}
		\alpha_{01} & \alpha_{02} \\
		\beta_{001} & \beta_{002} \\
\end{array} \right)
\end{equation*}
is non-degenerate. 

Using a combination of Lemmas~\ref{LmaCoordinateChange} and  \ref{LmaCoordinateChange2} 
%%%% (1) eliminate alpha_{02} using the first lemma (requires y_1 to change, too).
%%%% (2) eliminate beta_{001} using the second lemma.
%%%% (3) eliminate gamma_{0001} and gamma_{0002} using the 1st lemma.
we may find new coordinates $(x'_0, x'_1, x'_2)$ centered at $x$ and $(y_0, y'_1, y_2)$ centered at $y =F(x)$, (here $y_0$ and $y_2$ are standard coordinates) such that in these coordinates we have,
\begin{equation} \label{EqnSwallowtailAssumption}
	\left( \begin{array}{ccc}
		\alpha_{00} & \alpha_{01} & \alpha_{02} \\
		\beta_{000} & \beta_{001} & \beta_{002} \\
		\gamma_{0000} & \gamma_{0001} & \gamma_{002}
\end{array} \right) = \left( \begin{array}{ccc}
	0 & \pm r & 0 \\
	0 & 0 & \pm s \\
	\pm t & 0 & 0
\end{array} \right).
\end{equation}
for some strictly positive constants $r,s,t$. 

The local ring is given by $\cR_F = \R[x'_0] /((x'_0)^4)$ and hence, by the Generalized Malgrange Preparation Theorem~\ref{GeneralizedMalgrangePreparationTheorem}, there exist functions $a,b,c,d$ on $\R^2 \times I$ such that,
\begin{equation*}
(x'_0)^4 = F^*a + F^*b \cdot x'_0 + F^*c \cdot (x'_0)^2 + F^*d \cdot (x'_0)^3. 
\end{equation*}
By our conditions \ref{EqnSwallowtailAssumption} and by collecting terms we see that $\frac{\partial d}{\partial y_1}(0) = 0$ and $\frac{\partial d}{\partial y_2}(0) = 0$. Thus we may replace $x_0'$ by $x_0' - \frac 1 3 F^*d$, without effecting our conditions \ref{EqnSwallowtailAssumption}. In otherwords, we may assume without loss of generality that $d \equiv 0$. 

By using our conditions (\ref{EqnSwallowtailAssumption}) and again collecting terms, we see that
\begin{align*}
	\frac{\partial a}{\partial y_0} (0) &= \pm t & \frac{\partial a}{\partial y_1} (0) &= 0 & \frac{\partial a}{\partial y_2} (0) &= 0 \\
	 \frac{\partial b}{\partial y_0} (0) &= 0 & \frac{\partial b}{\partial y_1} (0) &= \pm r & \frac{\partial b}{\partial y_2} (0) &= 0 \\
	 \frac{\partial b}{\partial y_0} (0) &= 0 & \frac{\partial c}{\partial y_1} (0) &= 0 & \frac{\partial c}{\partial y_2} (0) &= \pm s 
\end{align*}
Thus the following pair of coordinate changes are valid in neighborhoods of $x$ and $y$:
\begin{align*}
	\overline x_0 &= \pm x_0 & 	\overline y_0 &= a(y_0, y_1, y_2) \\
	\overline x_1 &= F^*b & 	\overline y_1 &= b(y_0, y_1, y_2) \\
	\overline x_2 &= F^*c & 	\overline y_2 &= c(y_0, y_1, y_2)
\end{align*}
In these coordinates $F$ has the desired form. 
\end{proof}

\begin{proof}[Proof of Prop.~\ref{FoldCoordProp}]
%Again, the implicit function theorem guarantees that there exist coordinates around $x$ with the properties listed in Table~\ref{InitialChoiceofCoordTable}, i.e., 
%\begin{align*}
%	F^* y_0 & = h(x_0, x_1, x_2) \\
%	F^* y_1 & = x_1 \\
%	F^* y_2 & = x_2.
%\end{align*}
%where $dh = 0$ along $Y$. We may arrange for $Y = \{ x_0 = 0 \}$. 
Since this is an  $S_{[001], 0}$ singularity, we know that $h_{x_0x_0} \neq 0$ along $Y$. %, and moreover $F$ satisfies the transversality assumption that $j^1F$ is transverse to the strata $S_{[001]}$, i.e., the following map is surjective along $Y$:
%\begin{align*}
%	d(j^1F): T_x(\Sigma \times I) \to &\;  T_{j^1F(x)}J^1(\Sigma \times I, \R^2 \times I) / T_{j^1F(x)}S_{[001]}  \\
%	& \quad \cong \nu_{j^1F(x)} S_{[001]} \cong \hom(K_1, L). 
%\end{align*}
We introduce a new function,
\begin{equation*}
	g(x_0, x_1, x_2) := h(x_0, x_1, x_2) - \gamma( x_1, x_2).
\end{equation*}
Our assumptions on $h$ guarantee that
\begin{align*}
g(0, x_1, x_2) &\equiv 0 & \frac{\partial g}{\partial x_0}(0, x_1, x_2) &\equiv 0 &  \frac{\partial^2 g}{\partial x_0^2}(0, x_1, x_2) & \neq 0 
\end{align*}
and thus we have:
\begin{equation*}
	h(x_0, x_1, x_2) = x_0^2 \tilde g( x_0, x_1, x_2) + \gamma( x_1, x_2).
\end{equation*}
for some function $\tilde g$, such that $ \tilde g( 0, x_1, x_2) \neq 0$. Thus in a possibly smaller neighborhood of $Y$, the following is a valid coordinate change:
\begin{align*}
	\overline x_0 & = x \sqrt{|\tilde g (x_0, x_1, x_2)|} \\
	\overline x_1 &= x_1 \\
	\overline x_2 &= x_2. 
\end{align*}
In these coordinates $F$ has the desired normal form. 
\end{proof}

\begin{proof}[Preliminaries for Propositions \ref{Prop2DMorseRelationCoord}, \ref{Prop2DMorseCoord}, and \ref{PropCuspFlipCoord}] In these cases the implicit function theorem implies that there are coordinates $(x_0, x_1, x_2)$ centered at $x$ with the properties listed in Table~\ref{InitialChoiceofCoordTable}, i.e., 
\begin{align*}
	F^* y_0 & = x_0 \\
	F^* y_1 & = h(x_0, x_1, x_2)  \\
	F^* y_2 & = x_2.
\end{align*}
When the codimension is exactly three (Prop.~\ref{Prop2DMorseRelationCoord} and Prop.~\ref{PropCuspFlipCoord}), the equation $dh=0$ holds at the origin. We may expand $h$ in the variables $x_0, x_1, x_2$ to obtain
\begin{align*}
	h(x_1, x_2, x_3) = &\sum_{0 \leq i \leq j \leq 2} \alpha_{ij} x_i x_j + \sum_{0 \leq i \leq j \leq k \leq 2 } \beta_{ijk} x_i x_j x_k  + O(|x|^4).
\end{align*}

In all case we have the assumption that $j^1F$ is transverse to the stratum $S_{[011]}$, and hence the following map is surjective at $x$:
\begin{align*}
	d(j^1F): T_x(\Sigma \times I) \to &\;  T_{j^1F(x)}J^1(\Sigma \times I, \R^2 \times I) / T_{j^1F(x)}S_{[011]}  \\
	& \quad \cong \nu_{j^1F(x)} S_{[011]} \cong \hom(K_2, L) 
%	d(j^2F): T_x(\Sigma \times I) \to & \; T_{j^2F(x)}J^2(\Sigma \times I, \R^2 \times I) / T_{j^2F(x)}S_{[011], [001]}  \\
%	&\quad \cong \nu_{j^2F(x)} S_{[011], [001]} \cong \hom(K_2, L) \oplus \hom(K_1 \cdot K_1, L)
\end{align*}
In the codimension-three cases there will be additional transversality conditions. 
	\renewcommand{\qedsymbol}{}
\end{proof}

\begin{proof}[Proof of Prop.~\ref{Prop2DMorseCoord}]
%The implicit function theorem guarantees that there exist coordinates around $x$ with the properties listed in Table~\ref{InitialChoiceofCoordTable}, i.e., 
%\begin{align*}
%	F^* y_0 & = x_0 \\
%	F^* y_1 & = h(x_0, x_1, x_2) \\
%	F^* y_2 & = x_2.
%\end{align*}
In this case we further have that $Y = \{ (\gamma_0(x_2), 0 , x_2) \}$ in our chosen coordinates.  Since this is an  $S_{[011], [000]}$ singularity, along $Y$ we have $dh \equiv 0$, $\frac{\partial^2 h}{\partial x_1^2} \neq 0$, and 
\begin{equation*}
\left( \begin{array}{cc}
	\frac{\partial^2 h}{\partial x_0 \partial x_1} &  \frac{\partial^2 h}{\partial x_1^2} \\
	\frac{\partial^2 h}{\partial x_0^2} &  \frac{\partial^2 h}{\partial x_0 \partial x_1}
\end{array} \right)
\end{equation*}
is non-singular. 
Next we preform the following coordinate change, which is analogous to the one preformed in the proof or Prop.~\ref{PropPathofCuspCoordinates}:
\begin{align*}
	x_0' &= x_0  - \gamma_0(x_2)  & 					y_0' &= y_0 - \gamma_0(y_2)  \\
	x_1' &=x_1& 	y_1' &=  y_1 - \gamma_1(y_2)  \\
	x_2' &= x_2 &					y_2' &= y_2
\end{align*}
In these coordinates, $Y$ consists of the $x_2$-axis and its image consists of the $y_2$-axis. 
Let $g = F^* y_1'$. We have,
\begin{equation*}
	F^* y_1' = g(x_0', x_1', x_2') = h( x_0' + \gamma_0(x_2'), x_1', x_2') - \gamma_1(x_2').
\end{equation*}
Our conditions on $F$ assure
\begin{align*}
	g(0,0, x_2') &\equiv 0		& dg (0,0, x_2) & \equiv 0 \\
	\left( \begin{array}{cc}
		\frac{\partial^2 g}{\partial x_0' \partial x_1'} &  \frac{\partial^2 g}{\partial (x_1')^2} \\
		\frac{\partial^2 g}{\partial (x_0')^2} &  \frac{\partial^2 g}{\partial x_0' \partial x_1'}
	\end{array} \right)(0,0,x_2') & \text{ is non-singular} &   \frac{\partial^2 g}{\partial (x_1')^2} (0,0, x_2) & \neq 0.
\end{align*}
%%%%
%%%%  If we just use the Malgrange Prep. Theorem then we are tempted to get coordinates in which 
%%%%	F^* y_0 = x_0
%%%%	F^*y_1 = x_1^2	
%%%%	F^* y_2 = x_2
%%%% This is fine by the criteria we set out at the begining: it preserves the kernels of df and dpf along Y.
%%%% However it doesn't preserve these kernels along the singularities in this neighborhood 
%%%% (which are fold singularities in this case).
The local ring of this singularity is given by $\cR_F = \R[x_1'] / ((x_1')^2)$ and so by the Malgrange Preparation Theorem there exist functions $a, b$ on $\R^2 \times I$ such that 
\begin{equation*}
	(x_1')^2 = F^*a + F^*b \cdot x_1'.
\end{equation*}
By setting $x_1'' = x_1' - \frac 1 2 F^*b$, we may assume that $b \equiv 0$. Our conditions on $F$ (and hence $g$) imply that,
\begin{equation*}
a(y'_0, y'_1, y'_2) = A \cdot (y_0')^2 + \tilde a( y'_0, y'_1, y'_2)
\end{equation*}
for some non-zero constant $A$ and where,
\begin{align*}
	 \frac{\partial a}{\partial y'_0}(0,0,y'_2) &= 0 &  \frac{\partial a}{\partial y_1'}(0,0,y'_2)  & \neq 0 \\
	 \frac{\partial a}{\partial y'_2}(0,0,y'_2)  & = 0 & \frac{\partial^2 a}{\partial (y_0')^2}(0,0,y'_2)  & \neq 0 
\end{align*}
Thus the following is a valid pair of coordinate changes:
\begin{align*}
	\overline x_0 &= \sqrt{|A|} x_0' & 			\overline	y_0 &=  \sqrt{|A|}  y_0 +  \gamma_0(y_2')  \\
	\overline x_1 &=x_1'' & 					\overline y_1 &=  \tilde a(y_0', y_1', y_2') + \gamma_1(y_2)  \\
	\overline x_2 &= x_2' &					\overline y_2 &= y_2'
\end{align*}
In these coordinates, $F$ has the desired form. 
\end{proof}

\begin{proof}[Proof of Prop.~\ref{PropCuspFlipCoord}] 
%The implicit function theorem implies that there are coordinates $(x_0, x_1, x_2)$ centered at $x$ with the properties listed in Table~\ref{InitialChoiceofCoordTable}, i.e., 
%\begin{align*}
%	F^* y_0 & = x_0 \\
%	F^* y_1 & = h(x_0, x_1, x_2)  \\
%	F^* y_2 & = x_2.
%\end{align*}
%where $dh = 0$ at the origin. Expanding $h$ in the variables $x_0, x_1, x_2$ we get
%\begin{align*}
%	h(x_1, x_2, x_3) = &\sum_{0 \leq i \leq j \leq 2} \alpha_{ij} x_i x_j + \sum_{0 \leq i \leq j \leq k \leq 2 } \beta_{ijk} x_i x_j x_k  + O(|x|^4).
%\end{align*}

In this case $x$ is an $S_{[011], [001], 0}$-singularity and we have that $\alpha_{01} \neq 0$, $\alpha_{11} = 0$, and $\beta_{111} \neq 0$. 
By assumption, we have the further condition that $j^2F$ is transverse to the stratum $S_{[011], [001]}$, which means that the following map is surjective at $x$:
\begin{align*}
%	d(j^1F): T_x(\Sigma \times I) \to &\;  T_{j^1F(x)}J^1(\Sigma \times I, \R^2 \times I) / T_{j^1F(x)}S_{[011]}  \\
%	& \quad \cong \nu_{j^1F(x)} S_{[011]} \cong \hom(K_2, L) \\
	d(j^2F): T_x(\Sigma \times I) \to & \; T_{j^2F(x)}J^2(\Sigma \times I, \R^2 \times I) / T_{j^2F(x)}S_{[011], [001]}  \\
	&\quad \cong \nu_{j^2F(x)} S_{[011], [001]} \cong \hom(K_2, L) \oplus \hom(K_1 \cdot K_1, L)
\end{align*}
Given that $\alpha_{11} = 0$, $\alpha_{01} \neq 0$, and $\beta_{111} \neq 0$, our transversality assumptions are equivalent to the condition that the $3 \times 3$ matrix,
\begin{equation*}
	\left( \begin{array}{ccc}
		\alpha_{01} & 0 & \alpha_{12} \\
		\alpha_{00} & \alpha_{01} & \alpha_{02} \\
		\beta_{011} & \beta_{111} & \beta_{112}
\end{array} \right)
\end{equation*}
is non-degenerate. Using Lemma \ref{LmaCoordinateChange} and Lemma \ref{LmaCoordinateChange2} 
 %%%% (1) eliminate \alpha_{12}
 %%%% (2) eliminate \alpha_{02}
 %%%% (3) eliminate \alpha_{00}
 %%%% (4) using 2nd lma, eliminate \beta_{011}
 we may find coordinates $(x_0', x_1', x_2')$ and $(y'_0, y_1', y'_2)$ in which the above matrix becomes:
 \begin{equation*}
\left( \begin{array}{ccc}
\pm r& 0& 0 \\
0& \pm r  & 0 \\
 0 & \pm s & \pm t
\end{array} \right)
\end{equation*}
for some positive numbers $r,s,t$. %By scaling $x_0$ appropriately, we may assume that $s= 1$.

The local ring of this singularity is $\cR_F \cong \R[x_1']( (x_1')^3)$, and so by the Malgrange Preparation Theorem there exist functions $a,b,c$ on $\R^2 \times I$ such that,
\begin{equation*}
	(x_1')^3 = F^* a +  F^*b \cdot x'_1 + F^*c \cdot (x'_1)^2.
\end{equation*}
Collecting terms we see that,
\begin{align*}
	\frac{\partial a}{\partial y_0} (0) &= 0 & \frac{\partial a}{\partial y_1} (0) & \neq 0 & \frac{\partial a}{\partial y_2} (0) &= 0 \\
	 \frac{\partial b}{\partial y_0} (0) & \neq 0 & \frac{\partial b}{\partial y_1} (0) &= 0 & \frac{\partial b}{\partial y_2} (0) &= 0 \\
	 \frac{\partial b}{\partial y_0} (0) &= 0 & \frac{\partial c}{\partial y_1} (0) &= 0 & \frac{\partial c}{\partial y_2} (0) & \neq 0 
\end{align*}
Thus the following pair of coordinate changes are valid in neighborhoods of $x$ and $y$:
\begin{align*}
	\overline x_0 &= F^* b & 	\overline y_0 &= b(y_0', y_1', y_2') \\
	\overline x_1 &= \pm x_1'  &	\overline y_1 &= a(y_0', y_1', y_2') \\
	\overline x_2 &= F^*c & \overline y_2 &= c(y_0', y_1', y_2')
\end{align*}
In these coordinates $F$ has the desired form. 
\end{proof}

\begin{proof}[Proof of Prop.~\ref{Prop2DMorseRelationCoord}] 
%The implicit function theorem implies that there are coordinates $(x_0, x_1, x_2)$ centered at $x$ with the properties listed in Table~\ref{InitialChoiceofCoordTable}, i.e., 
%\begin{align*}
%	F^* y_0 & = x_0 \\
%	F^* y_1 & = h(x_0, x_1, x_2)  \\
%	F^* y_2 & = x_2.
%\end{align*}
%where $dh = 0$ at the origin. Expanding $h$ in the variables $x_0, x_1, x_2$ we get
%\begin{align*}
%	h(x_1, x_2, x_3) = &\sum_{0 \leq i \leq j \leq 2} \alpha_{ij} x_i x_j + \sum_{0 \leq i \leq j \leq k \leq 2 } \beta_{ijk} x_i x_j x_k  + O(|x|^4).
%\end{align*}
Since $x$ is an $S_{[011], [100], 0}$-singularity we know that $\alpha_{11} \neq 0$, and that 
\begin{equation*}
		\alpha_{01}^2 - \alpha_{11} \alpha_{00} = 0.
\end{equation*}

By assumption, we further know that $j^2F$ is transverse to the stratum  $S_{[011], [100]}$, which means that the following map are surjective at $x$:
\begin{align*}
%	d(j^1F): T_x(\Sigma \times I) \to &\;  T_{j^1F(x)}J^1(\Sigma \times I, \R^2 \times I) / T_{j^1F(x)}S_{[011]}  \\
%	& \quad \cong \nu_{j^1F(x)} S_{[011]} \cong \hom(K_2, L) \\
	d(j^2F): T_x(\Sigma \times I) \to & \; T_{j^2F(x)}J^2(\Sigma \times I, \R^2 \times I) / T_{j^2F(x)}S_{[011], [100]}   \\
	&\cong \nu_{j^2F(x)} S_{[011], [001]}
\end{align*}
Our transversality assumptions are now equivalent to the conditions that the $2 \times 3$ matrix,
\begin{equation*}
\left( \begin{array}{ccc}
	\alpha_{01} & \alpha_{11} & \alpha_{12} \\
	\alpha_{00} &\alpha_{01} & \alpha_{02}
\end{array} \right)
\end{equation*}
has rank two and that the $3 \times 3$-matrix,
\begin{equation*}
\left( \begin{array}{ccc}
	\alpha_{01} & \alpha_{11} & \alpha_{12} \\
	\alpha_{00} &\alpha_{01} & \alpha_{02} \\
	& \nabla (\alpha_{01}^2 - \alpha_{11} \alpha_{00}) & 
\end{array} \right)
\end{equation*}
is non-degenerate. Using Lemma \ref{LmaCoordinateChange} and Lemma \ref{LmaCoordinateChange2} 
 %%%% (1) eliminate \alpha_{01} (hence \alpha_{00} = 0 and \alpha_{02} \neq 0)
 %%%% (2) eliminate \alpha_{12} (hence \beta_{000} \neq 0)
 %%%% (3) using 2nd lma, eliminate \beta_{001}
 %%%% (4) eliminate \beta_{002}
 we may find acceptable coordinates $(x_0', x_1', x_2')$ and $(y'_0, y_1', y'_2)$ in which, for some strictly positive constants $r,s,t$, the following equation holds:
\begin{equation*} %\label{EqnSwallowtailAssumption}
	\left( \begin{array}{ccc}
		\alpha_{01} & \alpha_{11} & \alpha_{12} \\
		\alpha_{00} & \alpha_{01} & \alpha_{02} \\
		\beta_{000} & \beta_{001} & \beta_{002}
\end{array} \right) = \left( \begin{array}{ccc}
	0 & \pm r & 0 \\
	0 & 0 & \pm s \\
	\pm t & 0 & 0
\end{array} \right).
\end{equation*}

The local ring of this singularity is $\cR_F \cong \R[x_1']( (x_1')^2)$, and so by the Malgrange Preparation Theorem there exist functions $a,b,c$ on $\R^2 \times I$ such that,
\begin{equation*}
	(x_1')^2 = F^* a +  F^*b \cdot x'_1, 
\end{equation*}
and making a suitable change to the $x_1'$ coordinate we may assume that $b \equiv 0$. Expanding both sides we see that,
\begin{equation*}
	a(y'_0, y'_1, y'_2) = \tilde a ( y'_0, y'_1, y'_2) \pm s \cdot y_0' y_2'  \pm t \cdot y_0^3,
\end{equation*}
and moreover,
\begin{align*}
	\frac{\partial^3 \tilde a}{\partial y_0^3} (0) = 0, \qquad
	\frac{\partial^2 \tilde a}{\partial y_0 \partial y_2} (0) = 0, \textrm{ and } \qquad
	\frac{\partial^3 \tilde a}{\partial y_0^3} (0) = \pm \frac{1}{t} y_1 + O(|y|^2).
\end{align*}
Thus the following pair of coordinate changes are valid in neighborhoods of $x$ and $y$:
\begin{align*}
	\overline x_0 &= x_0' \sqrt[3]{t} 	& \overline y_0 &= y_0' \sqrt[3]{t}\\
	\overline x_1 &= x_1' & 	\overline y_1 &= \tilde a(y_0', y_1', y_2') \\
	\overline x_2 &=x_2' \cdot s \frac{1}{\sqrt[3]{t}} & \overline y_2 &= y_2' \cdot s \frac{1}{\sqrt[3]{t}}
\end{align*}
In these coordinates $F$ has the desired form. 
\end{proof}

\begin{remark}
The Coordinates for singularities that we have derived often  occur in multiple forms. For example there are two forms of the fold singularities of Proposition~\ref{FoldCoordProp}, according to the sign of $x_0$-coordinate. Viewing a fold singularity as a path of 1-manifold Morse singularities (which we elaborate on in the following section), we see that these two forms  correspond to the two possible indices of such Morse functions. In general we will call these various forms of singularities the {\em indices} of the singularity. Table~\ref{TableNumIndicesof3DSingularities} gives the number of indices in the 3-dimensional case. 
\end{remark}

\begin{table}[ht]
\begin{center}
\begin{tabular}{|c|c| c| } \hline
Singularity Stratum & Name & Number of Indices \\ \hline \hline
$S_0$ & Local Diffeomorphism & 0 \\ \hline 
$S_{[001], 0}$  & 1D Morse (Fold) & 2   \\
$S_{[011], 0}$  & 2D Morse & 4 \\
$S_{[011], [100],0}$  & 2D Morse Relation  & 8  \\ \hline
$S_{[001], [11], 0}$  & Cusp Inversion & 4 \\
$S_{[001], [01], 0}$  & 1D Morse Relation (Paths of Cusps) & 2 \\
$S_{[011], [001], 0}$  & Cusp Flip & 4  \\ \hline
$S_{[001], [01], 1, 0}$  & Swallowtail &  4 \\ \hline
\end{tabular}
\end{center}
\caption{The Number of Indices of the 3D Singularities}
\label{TableNumIndicesof3DSingularities}
\end{table}%

\section{The Geometry of the Singularities} \label{SectGeomofSing}

In the previous sections we derived normal coordinates for each of the types of singularities we encounter in generic maps $\Sigma \to \R^2$ and $\Sigma \times I \to \R^2 \times I$, when $\Sigma$ is a surface. These normal coordinates allow us to understand the local geometric nature of such generic maps. In this section we explain qualitatively what each of these singularities looks like and how this determines the local structure of $\Sigma$. 

In Morse theory, 
\index{Morse theory}
one of the pieces of data that one can recover from a Morse function is a picture in the target $\R$ of the singularities. In this 1-dimensional setting, this is not much data, just a finite sequence of isolated critical values. In the higher dimensional setting we are concerned with here, the image in the target of the singular locus carries far more information. Most of the singularities are non-isolated, for example a cusp singularity $S_{[01],1,0}$ occurs at a single point, but it is part of a whole arc of $S_{[01]}$ singularities. The configurations of these singularities amongst each other provides a great deal of local information about the source manifold. Moreover, their image in the target provides a partial picture of the global structure of this manifold, and in subsequent sections this will form the basis of a more elaborate diagram, from which we can recover the source manifold completely. 

\begin{figure}[htb]
\begin{center}
\begin{tikzpicture}[thick, decoration={snake, segment length = 6mm, amplitude = 0.5mm, pre length = 2mm, post length = 1mm}]
\draw (0,0) rectangle (6, 3)  rectangle (12, 0);
%%% The negative fold
\begin{scope}[xshift = 1cm, yshift = 0.5 cm]
\node (A) at (.1,.2) {};
\draw (A.center) -- +(-0.1, 0) -- +(-0.1,1.6) -- +(.2,1.6) arc (90: -90:  0.3cm and 0.1cm) -- +(-0.2, 0)
	-- +(-0.2,-1.6) -- +(0, -1.6) arc (-90: 0:  0.3cm and 0.1cm) decorate { -- +(0, 1.6)};
\draw [densely dashed] (A.center) -- +(.2, 0) arc (90: 0: 0.3cm and 0.1cm);
\end{scope}
%%% Neg. Fold Graphic
\draw (3,0.5) rectangle (5,2.5); \draw[red, decorate] (4.5, 2.5) -- +(0, -2);
%%%% The positive fold
\begin{scope}[xshift = 8cm, yshift = 0.5 cm, xscale = -1]
\node (A) at (.1,.2) {};
\draw (A.center) -- +(-0.1, 0) -- +(-0.1,1.6) -- +(.2,1.6) arc (90: -90:  0.3cm and 0.1cm) -- +(-0.2, 0)
	-- +(-0.2,-1.6) -- +(0, -1.6) arc (-90: 0:  0.3cm and 0.1cm) decorate {-- +(0, 1.6)};
\draw [densely dashed] (A.center) -- +(.2, 0) arc (90: 0: 0.3cm and 0.1cm);
\end{scope}
%%% Pos. Fold Graphic
\draw (9,0.5) rectangle (11,2.5); \draw[red, decorate] (9.5, 2.5) -- +(0, -2);
\end{tikzpicture}
\caption{Fold Singularities and Their Graphics}
\label{FoldGraphicsFig}
\index{fold singularity}
\index{singularity!fold}
\index{singularity!graphic}
\index{graphic}
\end{center}
\end{figure}
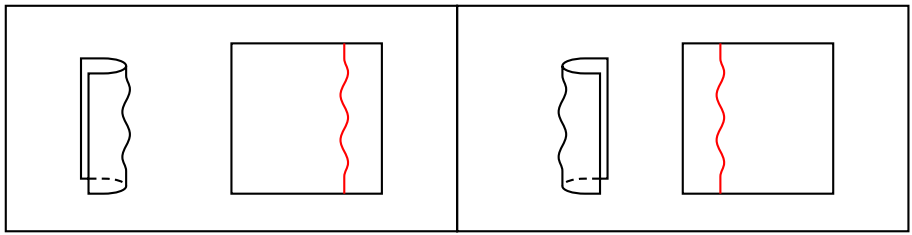

Given the importance of the image of the singular locus, it is not surprising that it has a name. We call this the {\em graphic} 
\index{graphic}
\index{graphic|see{singularity}}
\index{singularity!graphic}
of the generic map. We will see many examples in what follows. Let us look at our first examples. The fold singularities of Proposition~\ref{2dFoldsProp} and their corresponding graphics are depicted in Figure~\ref{FoldGraphicsFig}.
To the left of each diagram  in the figure we have depicted a 3-dimensional embedded surface. The projection of this surface into the plane of the page is a map to $\R^2$ which is generic and exhibits $S_{[01],0}$-singularities. 

Notice that these singularities are not isolated and that in the coordinates of Proposition~\ref{2dFoldsProp} the singular locus is parametrized by a curve, $\gamma$. The graphic of this singularity consists of a smoothly embedded arc in $\R^2$, whose projection to $\R$ is a diffeomorphism. In Figure~\ref{FoldGraphicsFig} we have depicted this as the wiggly red line. This is how we will typically draw our graphics.  The fold singularity can be considered as a vertical path of Morse functions for a 1-manifold. 
\begin{figure}[htb]
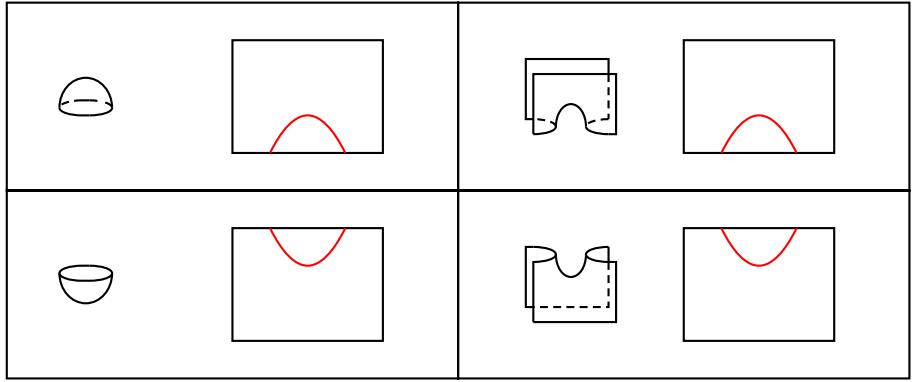

\begin{center}

\caption{2D Morse Singularities and Their Graphics}
\label{2DMorseGraphicFig}
\index{2D Morse singularity}
\index{singularity!2D Morse}
\index{singularity!graphic}
\index{graphic}
\end{center}
\end{figure}

\begin{figure}[htb]
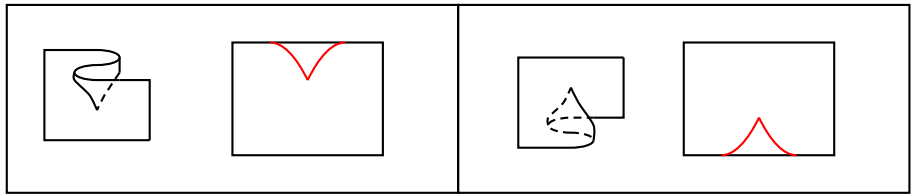

\begin{center}
%
\caption{Cusp Singularities and Their Graphics}
\label{CuspGraphicFig}
\index{cusp singularity}
\index{singularity!cusp}
\index{singularity!graphic}
\index{graphic}
\end{center}
\end{figure}

The four 2D Morse ($S_{[11]}, 0$) Singularities of Proposition~\ref{2d2dmorseprop} are depicted in Figure~\ref{2DMorseGraphicFig} while the two Cusp ($S_{[01], 1, 0}$) Singularities of Proposition~\ref{2dCuspProp} are depicted in Figure~\ref{CuspGraphicFig}. Our naming convention should now be clear. The fold singularities are so named because the correspond to singularities which occur when one sheet of the surface $\Sigma$ is ``folded'' over another.  

The 2D Morse singularities are similar to the fold singularities, but are distinguished by their projection to $\R$ (which is  projection onto the $y$-axis in the depictions in Figures~\ref{2DMorseGraphicFig} and \ref{FoldGraphicsFig}). For the folds, the projection to $\R$ has no critical points, while for the 2D Morse, there is a single isolated Morse singularity. If we were just concerned with the Morse function given by  projecting to $\R$, then we would not be able to distinguish between the two ``saddles''. These would both correspond to Morse singularities of index one. However we have more than a map to $\R$, we have a map to $\R^2$, and this additional coordinate allows us to distinguish these singularities. 

The Cusp singularities are equally well known. It is a singularity that arises in the study of Morse functions via Cerf theory. Our definition of the $S_{[01], 1, 0}$ stratum ensures that around these singularities the projection to $\R$ can be viewed as a time parameter and that the map to $\R^2$ can be viewed as a path of Morse functions on a 1-dimensional manifold. This path consists entirely of Morse functions except at the Cusp singularity itself, where a birth/death is occurring. 

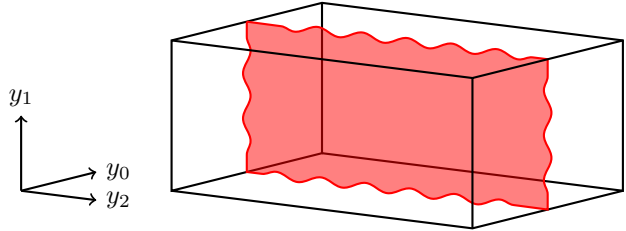
\begin{figure}[htb]
\begin{center}
\begin{tikzpicture}[thick, decoration={snake, segment length = 6mm, amplitude = 0.5mm, pre length = 2mm, post length = 2mm}]
\draw [->] (-2,0) -- (-1, 0.25) node [right] {$y_0$};
\draw [->] (-2,0) -- (-2, 1) node [above] {$y_1$};
\draw [->] (-2,0) -- (-1, -0.125) node [right] {$y_2$};
\draw (2, 2.5) -- (2, 0.5) -- (0,0) -- (0, 2) -- (2, 2.5) -- (6, 2) -- (6, 0) -- (2, 0.5) (6,0) -- (4, -0.5) -- (0,0); 
\fill[red, semitransparent] (1, 2.25) decorate {-- (1, 0.25)} decorate { -- (5, -0.25)}  decorate {-- (5, 1.75)} decorate {-- (1, 2.25)};
\draw [red] (1, 2.25) decorate {-- (1, 0.25)} decorate {-- (5, -0.25)} decorate {-- (5, 1.75)} decorate {-- (1, 2.25)};
\draw (6,2) -- (4, 1.5) -- (4, -0.5) (4, 1.5) -- (0, 2);
\end{tikzpicture}
\caption{Fold Singularities in Three Dimensions}
\label{FoldSingularitiesIn3DFig}
\index{fold singularity}
\index{singularity!fold}
\index{singularity!graphic}
\index{graphic}
\end{center}
\end{figure}

%\draw [red] (1, 2.25) decorate {-- (1, 0.25)} -- (5, -0.25) decorate {-- (5, 1.75)} -- (1, 2.25);

This brings us to the 3-dimensional singularities. These singularities can be divided into those of codimension three and those of lower codimension. The lower codimensional singularities correspond precisely to paths of the previous 2-dimensional singularities. Thus the $S_{[001], 0}$-singularities correspond to paths of the fold singularities. The image in $\R^2 \times I$ consists of a 2-dimensional sheet which is swept out, and whose projection to $\R \times I$ is  a diffeomorphism, see Figure~\ref{FoldSingularitiesIn3DFig}.

\begin{figure}[htb]
\begin{center}
\begin{tikzpicture}[thick, decoration={snake, segment length = 6mm, amplitude = 0.5mm, pre length = 2mm, post length = 1mm}, scale=0.9]
%\draw [->] (-2,0) -- (-1, 0.25) node [right] {$y_0$};
%\draw [->] (-2,0) -- (-2, 1) node [above] {$y_1$};
%\draw [->] (-2,0) -- (-1, -0.125) node [right] {$y_2$};
\draw (2, 2.5) -- (2, 0.5) -- (0,0) -- (0, 2) -- (2, 2.5) -- (6, 2) -- (6, 0) -- (2, 0.5) (6,0) -- (4, -0.5) -- (0,0); 
\draw [red] (0.5, 2.125) parabola (1, 1)  decorate {-- (5, 0.5)} parabola [bend at end] (4.5, 1.625) -- (0.5, 2.125);
\draw [red] (1.5, 2.375) parabola (1, 1)  decorate {-- (5, 0.5)} parabola [bend at end] (5.5, 1.875) -- (1.5, 2.375);
\fill[red, nearly transparent]  (0.5, 2.125) parabola (1, 1)  decorate {-- (5, 0.5)} parabola [bend at end] (4.5, 1.625) -- (0.5, 2.125);
\fill[red, nearly transparent]  (1.5, 2.375) parabola (1, 1)  decorate {-- (5, 0.5)} parabola [bend at end] (5.5, 1.875) -- (1.5, 2.375);
\draw (6,2) -- (4, 1.5) -- (4, -0.5) (4, 1.5) -- (0, 2);
\begin{scope}[xshift = 7cm]
\draw (2, 2.5) -- (2, 0.5) -- (0,0) -- (0, 2) -- (2, 2.5) -- (6, 2) -- (6, 0) -- (2, 0.5) (6,0) -- (4, -0.5) -- (0,0); 
\draw [red] (0.5, 2.125) parabola [bend at end] (1, 1)  decorate {-- (5, 0.5)} parabola (4.5, 1.625) -- (0.5, 2.125);
\draw [red] (1.5, 2.375) parabola [bend at end]  (1, 1)  decorate {-- (5, 0.5)} parabola  (5.5, 1.875) -- (1.5, 2.375);
\fill[red, nearly transparent]  (0.5, 2.125) parabola [bend at end]  (1, 1)  decorate {-- (5, 0.5)} parabola (4.5, 1.625) -- (0.5, 2.125);
\fill[red, nearly transparent]  (1.5, 2.375) parabola [bend at end]  (1, 1)  decorate {-- (5, 0.5)} parabola(5.5, 1.875) -- (1.5, 2.375);
\draw (6,2) -- (4, 1.5) -- (4, -0.5) (4, 1.5) -- (0, 2);
\end{scope}
\end{tikzpicture}
\caption{Cusp and 2D Morse Singularities in Three Dimensions}
\label{CuspAnd2DMorseIn3DFig}
\index{cusp singularity}
\index{singularity!cusp}
\index{2D Morse singularity}
\index{singularity!2D Morse}
\index{singularity!graphic}
\index{graphic}
\end{center}
\end{figure}
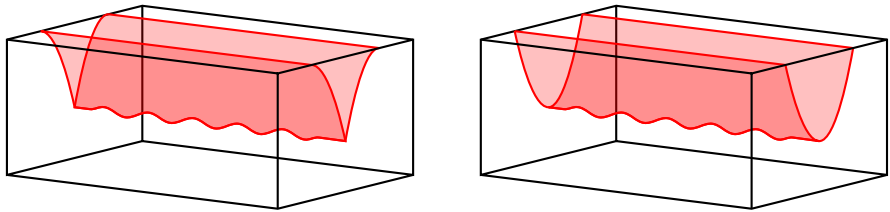

Similarly the Cusp singularities $S_{[001], [01], 0}$ and the 2D Morse singularities $S_{[011], 0}$ have graphics which consist of paths of the previous 2-dimensional graphics. Some of these are depicted in Figure~\ref{CuspAnd2DMorseIn3DFig}, with the others obtained by the obvious permutations of these graphics.

\begin{figure}[htb]
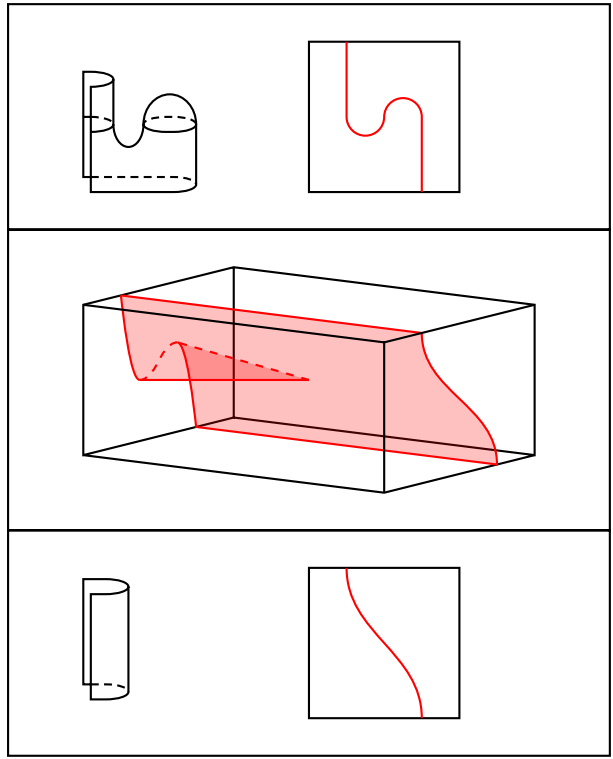

\begin{center}

\caption{2D Morse Relations and Their Graphics}
\label{2DMorseRelationFig}
\index{2D Morse relation singularity}
\index{singularity!2D Morse relation}
\index{singularity!graphic}
\index{graphic}
\end{center}
\end{figure}

The codimension three singularities are more interesting. For example the 2D Morse Relation singularity, $S_{[011], [100], 0}$, occurs,  as the name suggests, around a path of maps to $\R^2$ in which two 2D Morse singularities form a birth/death. There are a total of eight possibilities for this singularity, and Figure~\ref{2DMorseRelationFig} shows one of these. This figure shows the initial and final surfaces, together with their 2-dimensional graphics. In the middle pane it shows the full 3-dimensional graphic. The remaining seven versions of this singularity are obtained as the seven obvious reflections of this diagram.

\begin{figure}[htb]
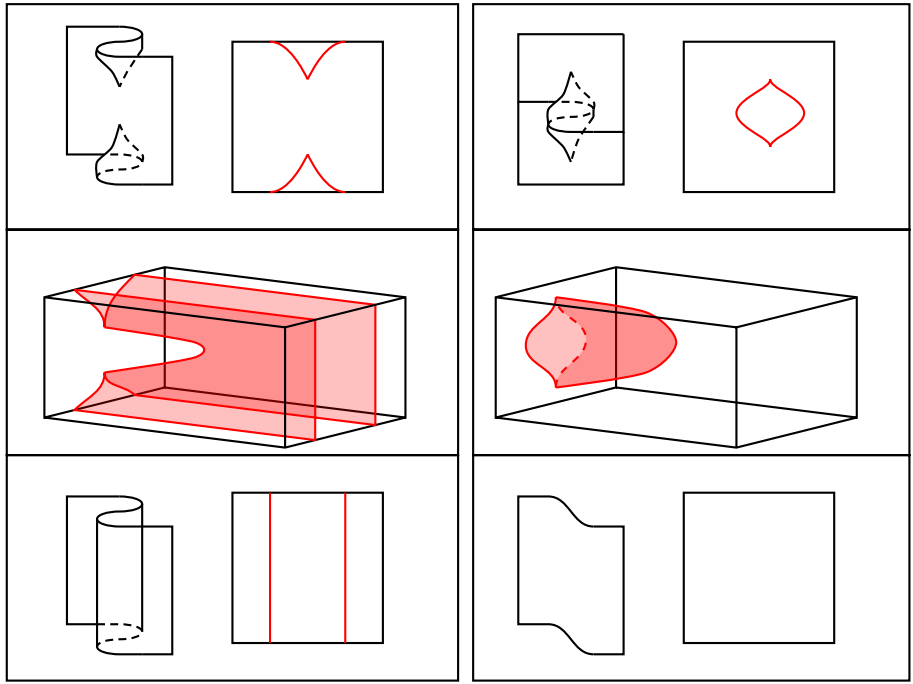

\begin{center}

\caption{Cusp Inversions and Their Graphics}
\label{CuspInversionGraphicFig}
\index{cusp inversion singularity}
\index{singularity!cusp inversion}
\index{singularity!graphic}
\index{graphic}
\end{center}
\end{figure}

The Cusp Inversion singularities, $S_{[001], [11], 0}$, occur in four types, of which two are depicted in Figure~\ref{CuspInversionGraphicFig}. Again we have show the initial and final surfaces, their 2-dimensional graphics, and the 3-dimensional graphic of the singularity. The remaining two forms of this singularity are given by reflecting these. These singularities are important for the bordism bicategory because, as we shall see, they witness the fact that the cusp bordisms are inverse to one another.

\begin{figure}[htb]
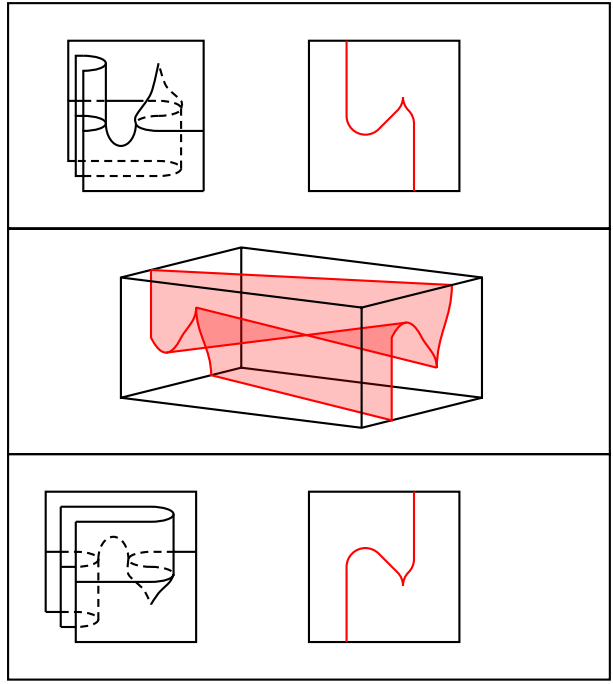

\begin{center}

\caption{Cusp Flips and Their Graphics}
\label{CuspFlipsGraphicFig}
\index{cusp flip singularity}
\index{singularity!cusp flip}
\index{singularity!graphic}
\index{graphic}
\end{center}
\end{figure}

Perhaps the most interesting of the singularities  are the Cusp Flip singularities, $S_{[011], [001], 0}$. They involve a non-trivial interaction between Cusp and 2D Morse singularities. They appear in four types, one of which is depicted in Figure~\ref{CuspFlipsGraphicFig}. The remaining types are again given by permutation. In terms of the graphic, they are easy enough to visualize; the point of the cusp simply ``flips over,'' producing a new graphic. We will have more to say about this singularity and its role in the bordism bicategory in Chapter \ref{ClassificationChapter}. 

\begin{figure}[htb]
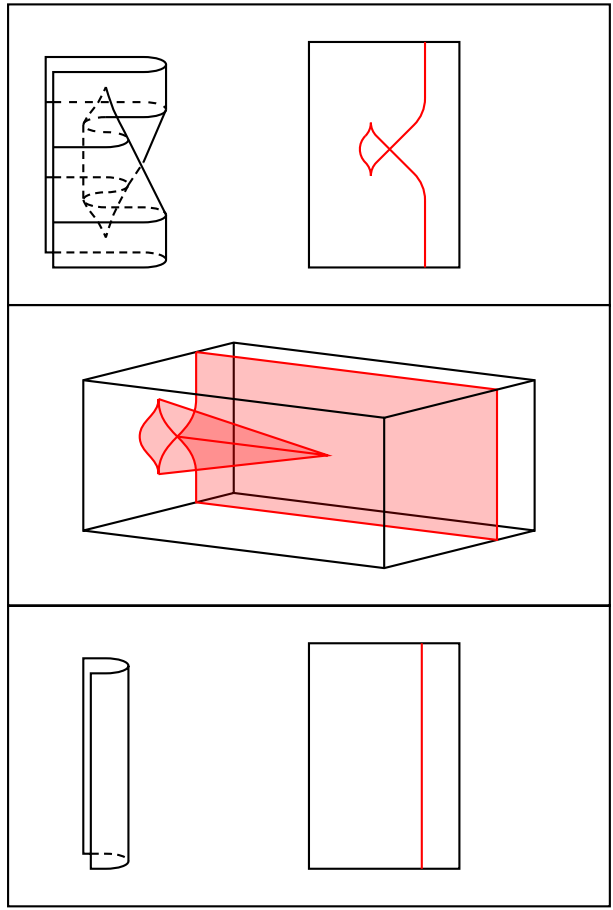

\begin{center}

\caption{Swallowtails and Their Graphics}
\label{SwallowtailGraphicFig}
\index{swallowtail singularity}
\index{singularity!swallowtail}
\index{singularity!graphic}
\index{graphic}
\end{center}
\end{figure}

The final singularity, the Swallowtail Singularity, $S_{[001], [01], 1,0}$, is a standard singularity arising from the two parameter version of Cerf theory. One of the four types is depicted in Figure~\ref{SwallowtailGraphicFig}, and the remaining are given by permuting this diagram. We invite the reader to spend some time visualizing this singularity. 

\section{Multijet Considerations} \label{SectMultijetConsiderations}
	
In the previous sections we have analyzed the local behavior of generic maps $\Sigma \to \R^2$ and $\Sigma \times I \to \R^2 \times I$ and we found that a convenient tool for understanding the structure of $\Sigma$ in terms of the singularities of these maps was the {\em graphic} of the map, i.e., the image in the target of the singular locus. So far,  we have only paid attention to the local behavior of this graphic, but the graphic itself makes sense globally. In fact,  the global graphic will play key role in the next section. In this section we make the first steps toward understanding its structure.

By jet transversality, we know that each of the singular loci occurs on a submanifold of $\Sigma$ or $\Sigma \times I$ and that these submanifolds map into the target via an immersion.\footnote{This is true, at least if we restrict to the Local Diffeomorphism, Fold, Cusp, and 2D Morse singularities, in the 2-dimensional case,  and restrict to the eight singularities listed in Table~\ref{TableLocalRingof3DSingularities} in the 3-dimesnional case.} 
 They are not typically embedded. In fact, we saw that in any neighborhood of the image of a Swallowtail singularity, there must be an arc where the images of two Fold loci intersect. In this section we want to better understand how these immersed manifolds intersect. For example, consider the following three scenarios:
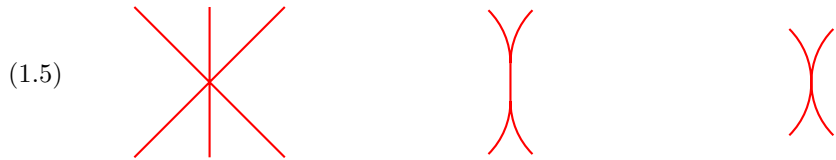
\begin{equation} \label{ThreeSingularitiesEqn}
\begin{tikzpicture}[thick]
\pgfsetbaseline{1cm}
\draw [red] (0,2) -- (2,0) (0, 0) -- (2,2) (1,2) -- (1,0);
\node (B) at (5,1.25) {}; 
\node (C) at (5, 0.75) {};
\draw [red] (B.center) arc (0: 45 : 1cm) (C.center) arc (0: -45 : 1cm) (B.center) arc (180: 135 : 1cm) (C.center) arc (180: 225 : 1cm);
\draw [red] (B.center) -- (C.center);
\node (A) at (9,1) {}; 
\draw [red] (A.center) arc (0: 45 : 1cm) (A.center) arc (0: -45 : 1cm) (A.center) arc (180: 135 : 1cm) (A.center) arc (180: 225 : 1cm);
\end{tikzpicture}
\end{equation}
These curves are not in general position, but from the analysis we have done so far, such configurations cannot be ruled out. 

However, we can force our singular loci to be in general position by enlarging our collection of singularity strata and requiring transversality with respect to these strata. The new strata will live in the multi-jet spaces, and transversality to these strata will imply, among other things, that the images of the singular loci meet in general position. After introducing these strata, the Multi-jet Transversality Theorem~\ref{MultiJetTransversalityThm} ensures that a generic set of maps satisfy both these new transversality conditions, as well as the previous ones. 

Recall that the multi-jet space 
\index{multi-jet bundle}
is defined as the pull-back of $J^{\bf k}(X, Y) = J^{k_1}(X, Y) \times \cdots \times J^{k_s}(X, Y)$ to  $X^{(s)} = \{ (x_1, \dots, x_s) \in X^s \; | \; x_i \neq x_j \; \textrm{ for } 1 \leq i < j \leq s\}$, where  ${\bf k} = (k_1, \dots, k_s)$ is a multi-index of natural numbers. The multi-jet of a map $f: X \to Y$ is then given by,
\begin{equation*}
	j^{\bf k} f (x_1, \dots, x_s) = ( j^{k_1} f( x_1), \dots, j^{k_s}f( x_s) ). 
\end{equation*}
Let $\Delta Y \subseteq Y^s$ denote the diagonal, and let $\beta:  J^{\bf k}(X, Y) \to Y^s$ be the canonical projection. Notice that $\beta^{-1}(\Delta Y)$ is a sub-manifold of $J^{\bf k}(X, Y)$.

Specializing to the case $X = \Sigma$ is a surface and $Y = \R^2$, 
consider the following submanifolds of the multi-jet spaces:
\begin{align*}
	(S_{[11]} \times S_{[11]}) \cap \beta^{-1}(\Delta Y) & \subseteq J^{(1,1)}(\Sigma, \R^2) \\
	(S_{[11]} \times S_{[01]}) \cap \beta^{-1}(\Delta Y) & \subseteq J^{(1,1)}(\Sigma, \R^2) \\
	(S_{[01]} \times S_{[01], 1}) \cap \beta^{-1}(\Delta Y) & \subseteq J^{(1,2)}(\Sigma, \R^2) \\
	(S_{[01],0} \times S_{[01], 0}) \cap \beta^{-1}(\Delta Y)& \subseteq J^{(2,2)}(\Sigma, \R^2)  
\end{align*}
These have codimensions 6, 5, 5, and 4, respectively and since the dimension of $\Sigma^{(2)}$ is $2 + 2 = 4$, only the last singularity occurs. Moreover we can consider the stratum,
\begin{equation*}
	(S_{[01]} \times S_{[01]} \times  S_{[01]}) \cap \beta^{-1}(\Delta Y) \subseteq J^{(1,1,1)}(\Sigma, \R^2). 
\end{equation*}
The codimension of this submanifold is 7, while the dimension of $\Sigma^{(3)}$ is 6, so this singularity does not occur. 

What does this mean? It means that for a generic map $f: \Sigma \to \R^2$, whose multi-jets are transverse to all these strata, the only singularities that can have the same critical value (i.e., image in $\R^2$) are fold singularities, and that these must occur in pairs and not in triples. Moreover (if $\Sigma$ is compact, which we will tacitly assume from now on) the critical values corresponding to multiple fold singularities must be isolated, and  only a finite number occur. This rules out the first two scenarios from \ref{ThreeSingularitiesEqn}, but the third still remains.

The characteristic that distinguishes this final  scenario from general position is that the tangent vectors to both curves coincide at the intersection point. This is in contrast to the intersection point which occurs, for example, in a neighborhood of a swallowtail critical value, in which the tangent vectors are transverse, see Figure~\ref{SwallowtailGraphicFig}. In order to rule out this scenario, we will have to understand how to read off these tangent vectors from the singularity strata. 

Recall that the singularity stratum $S_{[01]}$ was defined as an open subset of the Thom-Boardmann stratum $S_1 \subseteq J^1(\Sigma, \R^2)$. It will be enough to work with the simpler Thom-Boardman stratum $S_1$.  Recall that over $S_1$, there existed two canonical line bundles, $K = \ker df$ and $L = \coker df$ (where $f$ is any local map representing the jet in $S_1$). Recall also, that $S_1^{(2)}$ denotes the inverse image of $S_1$ in the second jet space $ J^2(\Sigma, \R^2)$, and that there is a map of bundles over $S_1$,
\begin{equation*}
	S_1^{(2)} \to \hom(T\Sigma, TJ^1(\Sigma, \R^2)) \to \hom(T\Sigma, TJ^1(\Sigma, \R^2) / TS_1)
\end{equation*}
where these vector bundles should be interpreted as their pullbacks to $S_1$. 

Let $S_1(f)$ denote those points in $\Sigma$ in which $f$ obtains an $S_1$-singularity. The tangent space to $S_1(f) = (j^1f)^{-1}( S_1)$ in $\Sigma$ will be the kernel of the induced map,
\begin{equation} \label{EqnTSigToNuS1}
	T\Sigma \to TJ^1(\Sigma, \R^2) / TS_1 = \nu S_1
\end{equation}
which can be read off from the second jet, $j^2f$. We would then like to use these vector spaces to obtain a stratification of $(S_1 \times S_1) \cap \beta^{-1}(\Delta Y)$, in particular by stratifying it according to whether the images of these vector spaces under $df$ coincide or not. The problem is that the kernel to (\ref{EqnTSigToNuS1}) does not form a vector bundle over $S^{(2)}_1$.  The dimensions of these vector spaces jump. 

Fortunately there is an open subset,   $S^{(2), 0}_1 \subseteq S^{(2)}_1$, which has the minimal corank of zero (i.e., where the map (\ref{EqnTSigToNuS1}) is surjective). Over this subspace, the kernel, $C$, forms a line bundle and we may try to proceed with our construction as planned. Note that the condition that $j^1f \pitchfork S_1$ at $x \in \Sigma$ is precisely that the above map is surjective, i.e., that $j^2f(x) \in S^{(2), 0}_1$. This is why it is  sufficient just to work with $S^{(2), 0}_1$ and to ignore the rest of the stratum $S_1^{(2)}$. 

Thus the space $(S^{(2), 0}_1 \times S^{(2), 0}_1) \cap \beta^{-1}(\Delta Y)$ has a canonical pair of line bundles over it. We would like to map these line bundles to $TY$ and form a further stratification according whether their images coincide.  We run into a similar problem in that the images in $TY$ do not form a vector bundle over $(S^{(2), 0}_1 \times S^{(2), 0}_1) \cap \beta^{-1}(\Delta Y)$, again because their dimensions may jump. 

Determining the dimension of the image of $C$ under $df$ is equivalent to determining the dimension of the intersection of $C$ with $K = \ker df$. Fortunately we have already stratified $S_1^{(2)}$ according to this dimension. This is precisely the corank of the induced map,
\begin{equation*}
	K \to TJ^1(\Sigma, \R^2) / TS_1 = \nu S_1 \cong \hom(K, L).
\end{equation*}
The only relevant multi-jet stratum will be,
\begin{equation*}
	S_{1,0 \times 1,0} := [(S^{(2), 0}_1 \cap S_{1,0}) \times (S^{(2), 0}_1 \cap S_{1,0}) ) ] \cap \beta^{-1}(\Delta Y)
\end{equation*}
over which there exists a pair of lines $C_1$ and $C_2$, whose images under $df$ form another pair of lines $df(C_1), df(C_2) \subseteq TY$. The map from $S_{1,0 \times 1,0}$ to the Grassmannian of pairs of lines is a submersion, which can be seen after choosing local coordinates. This Grassmannian is stratified into a codimension zero stratum and a codimension one stratum according to whether these lines coincide or not. 
Hence we get an induced stratification of  $S_{1,0 \times 1,0}$ by pulling-back the stratification of this Grassmannian, and hence a pair of submanifolds of $J^{(2,2)}( \Sigma, \R^2)$ of total codimensions 4 and 5 respectively. The later type of singularity does not occur as $\Sigma^{(2)}$ only has dimension 4. This rules out the last scenario of \ref{ThreeSingularitiesEqn} for generic maps $f: \Sigma \to \R^2$. We summarize these results in the following theorem.

\begin{theorem} \label{ThmGraphicOfMaptoR2}
The graphic 
\index{singularity!graphic}
\index{graphic}
of a generic map $f: \Sigma \to \R^2$, with $\Sigma$ a compact surface, consists of a finite number of embedded curves (the $S_{[01], 0}$-singularities) whose projections to $\R$ are diffeomorphisms, and two finite collections of isolated points (corresponding to $S_{[11],0}$- and $S_{[01],1,0}$-singularities), such that in a neighborhood of each of these points there exist coordinates in which the graphic has the form of a 2D Morse graphic (Figure~\ref{2DMorseGraphicFig}) or the form of a Cusp graphic (Figure~\ref{CuspGraphicFig}), respectively. Moreover, the embedded curves may intersect, but are in general position. 
\end{theorem}

It will be useful if we codify these sorts of diagrams in an abstract definition. We will choose to call this new sort of diagram a {\em graphic}, which generalizes our previous use of terminology.  Then the above theorem may be expressed by saying that the images of the critical points of a generic function form an abstract graphic. %With an eye towards our future uses, we will now also introduce the version of graphic based on the square $I^2$ rather than the plane $\R^2$. This is suitable for the kinds of manifolds with corners that will be introduced in Chapter~\ref{ClassificationChapter}.

\begin{definition} \label{2dGraphicDefn}
A {\em 2-dimensional graphic} 
\index{singularity!graphic}
\index{graphic}
is a diagram in $\R^2$ consisting of a finite number of embedded  closed curves whose projections to $\R$ ($y$-axis) are local diffeomorphisms, and two finite collections of isolated points. These curves and points are labeled. The curves are labeled by $S_{[01], 0}$ together with one of the two possible indices of the fold singularities. The points are  labeled by either $S_{[11],0}$ or $S_{[01],1,0}$, together with the possible indices of these singularities. We require that in a neighborhood of each of the points there exist coordinates in which the graphic has the form of a 2D Morse graphic (Figure~\ref{2DMorseGraphicFig}) or the form of a Cusp graphic (Figure~\ref{CuspGraphicFig}), respectively. Furthermore, the embedded curves are required to be in general position. %Moreover the points are disjoint from $\partial I^2$, the curves are disjoint from the left and right sides of the square ($\{0\} \times I$ and $\{1\} \times I$), and the curves intersect the rest of the boundary transversally. Finally the boundary of each curve is consists precisely of its intersection with the points and the boundary of $I^2$.  

For later convenience we also define a lower dimensional case. A {\em 1-dimensional graphic} consists of a finite collection of isolated points in $\R$. 
\end{definition}

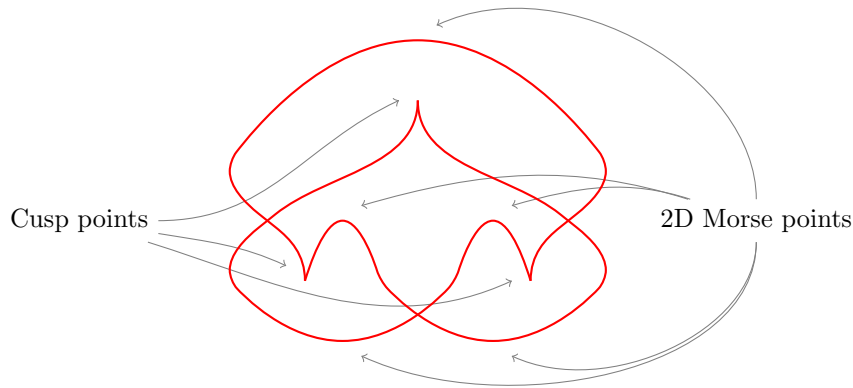
\begin{figure}[htbp]
	\begin{center}
	\begin{tikzpicture}[thick, yscale=0.8]
		\node (C) at (-2, 2) {Cusp points};
		\node (M) at (7,2) {2D Morse points};
		\draw [->, gray, thin] (C) to [out=0, in=210] (2.25, 4);
		\draw [->, gray, thin] (C) to [out=350, in=150] (0.75, 1.25);
		\draw [->, gray, thin] (C) to [out=340, in=210] (3.75, 1);	
		\draw [->, gray, thin] (M) to [out=270, in=330] (1.75, -0.25);
		\draw [->, gray, thin] (M) to [out=270, in=330] (3.75, -0.25);
		\draw [->, gray, thin] (M) to [out=90, in=30] (2.75, 5.25);
		\draw [->, gray, thin] (M) to [out=160, in=20] (1.75, 2.25);
		\draw [->, gray, thin] (M) to [out=160, in=20] (3.75, 2.25);
		\draw [red, rounded corners] (2.5, 4) to [out = -90, in = 45] (0.5, 2) to [out = 225, in = 90]  (0,1) parabola bend (1.5, 0) (3,1) parabola bend (3.5, 2) (4,1);
		\draw [red, rounded corners] (2.5, 4) to [out = -90, in = 135] (4.5, 2) to [out = -45, in = 90] (5,1) parabola bend (3.5, 0) (2, 1) parabola bend (1.5, 2) (1,1);
		\draw [red, rounded corners] (1,1) to [out = 90, in = -45] (0.5, 2) to [out = 135, in = -90] (0, 3) parabola bend (2.5, 5) (5, 3) to [out = -90, in = 45] (4.5, 2) to [out = 225, in = 90] (4, 1);
	\end{tikzpicture}
	\end{center}
	\caption{The following is one possible graphic for a generic map from $\Sigma = \R \P^2$ to $\R^2$. }
	\label{fig:RP2graphic}
\end{figure}

The intersections of the 3-dimensional singular loci can similarly be controlled by providing an appropriate multi-jet stratification. In this case the stratification and resulting phenomena are more interesting. We will again begin with a coarse analysis that will rule out many possible situations. Then we will provide a finer stratification of the remaining cases, taking into account tangential information, as we did in the 2-dimensional case. Recall the seven non-tirival singularities listed in Table~\ref{TableLocalRingof3DSingularities}, with their codimensions. 

The codimension of $S_a \times S_b \cap \beta^{-1}( \R^2 \times I)$ is given by,
\begin{equation*}
	\codim S_a + \codim S_b + 3
\end{equation*}
while the dimension of $(\Sigma \times I)^{(2)}$ is 6. By comparing the codimensions of the relevant strata we see that there are only the following three possibilities with codimension $\leq 6$,
\begin{align*}
	&S_{[001], 0} \times S_{[001], [01], 0} \cap \beta^{-1}( \R^2 \times I)&& \codim = 6 \\
	&S_{[001], 0} \times S_{[011], 0} \cap \beta^{-1}( \R^2 \times I)&& \codim = 6 \\
	&S_{[001], 0} \times S_{[001], 0} \cap \beta^{-1}( \R^2 \times I)&& \codim = 5. 
\end{align*}
Moreover, the dimension of $(\Sigma \times I)^{(3)}$ is 9, and hence $S_{[001], 0} \times S_{[001], 0} \times S_{[001], 0} \cap \beta^{-1}( \R^2 \times I)$, of codimension 9 in $J^{(2,2,2)}( \Sigma \times I, \R^2 \times I)$, is the only triple intersection which can occur. 

We must now analyze the tangential structure of all these singularity loci. As before, for each of the singularities $S_\alpha$, where $\alpha = [001]$ or $[011]$, we may look at it's pre-image $S_\alpha^{(2)}$ in $J^2(\Sigma \times I , \R^2 \times I)$, which forms a bundle over $S_\alpha$. Again there are maps of fiber bundles over $S_\alpha$,
\begin{equation*}
S_\alpha^{(2)} \to \hom (T(\Sigma \times I), TJ^1(\Sigma \times I, \R^2 \times I)) \to \hom (T(\Sigma \times I, \nu S_\alpha),
\end{equation*}
where each vector bundle is considered a bundle over $S_\alpha$ via pulling it back. 

We can consider the open subset of those maps $TX \to \nu S_\alpha$ consisting of minimal corank (i.e., maximal rank). The inverse image defines an open subset $S_\alpha^{(2), 0} \subseteq S_\alpha^{(2)}$. Over $ S_\alpha^{(2), 0}$ there is a vector bundle $C_\alpha$ consisting of the kernel of the map, $TX \to \nu S_\alpha$. If $f: \Sigma \times I \to \R^2 \times I$ is generic and takes an $S_\alpha$-singularity at $x \in X$, then $j^1f \pitchfork S_\alpha$ at $x$, which is equivalent to $j^1f(x) \in  S_\alpha^{(2), 0} \subseteq S_\alpha^{(2)}$. If $S_\alpha(f)$ denotes the submanifold of $\Sigma  \times I$ consisting of $S_\alpha$-singularities, then the tangent space to $S_\alpha(f)$ at $x$ is precisely the vector space $C_{\alpha, j^2f(x)} \subseteq T_x(\Sigma \times I)$. 

We will want to consider the image of $C_\alpha$ under the map $df: TX \to TY$. Again this is not a bundle on all of $S_\alpha^{(2), 0} $, but we have already stratified $S_\alpha^{(2)}$ into substrata on which $df$ has a fixed dimensional intersection with $C_\alpha$, and on the intersections of this stratification with the open set $  S_\alpha^{(2), 0} $, we get vector bundles $df(C_\alpha)$. The possible cases are listed in Table~\ref{TangetSpaceTable}.

\begin{table}[ht]
\begin{center}
\begin{tabular}{|c|c| c | c|} \hline
Singularity ($\alpha$) & Dim. of $C_\alpha$ & Dim. of $df(C_\alpha)$  & Properties \\ \hline
$S_{[001], 0}$ & 2& 2 & Proj. to $T(\R \times I)$ \\
&&& is surjective. \\
$S_{[001], [01], 0}$ & 1 & 1 &  \\
$S_{[011], 0}$ & 1 & 1 & Does not contain $\partial_{y_0}$ \\ \hline
\end{tabular}
\end{center}
\caption{Tangent Spaces to the Images of 3D Singularities}
\label{TangetSpaceTable}
\end{table}%

Now we may use these bundles on the multi-jet strata, together with the subbundles $\text{Span}\{ \partial_{y_0}\}$ and $\text{Span}\{ \partial_{y_0}, \partial_{y_1}\}$, to construct a further refinement of this stratification. Transversality with respect to this stratification will then imply that the images of the singular loci are in general position with respect to one another and with respect to the subbundles 
$\text{Span}\{ \partial_{y_0}\}$ and $\text{Span}\{ \partial_{y_0}, \partial_{y_1}\}$.

On the multi-jet stratum $S_{[001], 0} \times S_{[001],0}$ (intersected with $S_{[001]}^{(2),0} \times S_{[001]}^{(2),0}$), we have the two vector bundles denoted by $df(C_{[001]})$. We will further stratify $S_{[001], 0} \times S_{[001],0}$ by the dimension of the intersection of these 2-dimensional bundles, as well as their intersection with the plane $\ker\left(dp^2: T(\R^2 \times I) \to T( I)\right)$. There is a codimension 2 substratum, where these vector bundles coincide in $T(\R^2 \times I)$. The codimension zero substratum consists of where they only have a 1-dimensional intersection. Inside this stratum, there is an additional sub-stratification according to whether this intersection is contained in $\text{Span}\{ \partial_{y_0}, \partial_{y_1}\}$, or not. See Table~\ref{CodimOfMultijetStrata}. Note that $\dim df(C_{\alpha_1}) \cap df(C_{\alpha_2}) \cap \text{Span}\{ \partial_{y_0}\} = (0)$ automatically. 

Next consider the multi-jet stratum $S_{[001], 0} \times S_{[011],0}$. Over this stratum,   (intersected with $S_{[001]}^{(2),0} \times S_{[011]}^{(2),0}$), we have a line bundle and a rank two vector bundle $C_{[011]}$ and $C_{[001]}$, respectively. We will stratify this according to the whether their images under $df$ intersect each other. % and whether they intersect  $\text{Span}\{ \partial_{y_0}, \partial_{y_1}\}$, non-trivially. 
There is a an open stratum in which they do not intersect and a codimension one stratum in which they do. See Table~\ref{CodimOfMultijetStrata}.

\begin{table}[ht]
\begin{center}
\begin{tabular}{|c|c| c| c|c|} \hline
Singularity  & $A$ & $B$ & $C$ & Total Codimension  \\ \hline
$S_{[001], 0} \times S_{[001],0} \cap \beta^{-1}(\R^2 \times I)$ & 2 &0&  & 5+2 = 7 \\
$S_{[001], 0} \times S_{[001],0} \cap \beta^{-1}(\R^2 \times I)$& 1 & 0&1 & 5 + 1 = 6 \\
$S_{[001], 0} \times S_{[001],0} \cap \beta^{-1}(\R^2 \times I)$ & 1 &0& 0 & 5  \\ \hline
%%%%
$S_{[001], 0} \times S_{[011],0} \cap \beta^{-1}(\R^2 \times I)$ & 0 & 0 & 0  & 6  \\
$S_{[001], 0} \times S_{[011],0} \cap \beta^{-1}(\R^2 \times I)$ & 1 & &   & 6 + 1 = 7  \\ \hline
%%%%
$S_{[001], 0} \times S_{[001],[01],0} \cap \beta^{-1}(\R^2 \times I)$ & 0 & 0 & 0  & 6  \\
$S_{[001], 0} \times S_{[001],[01],0} \cap \beta^{-1}(\R^2 \times I)$ & 1 & &   & 6 + 1 = 7  \\ \hline
\end{tabular}
\vspace{0.25cm}

$A = \dim df(C_{\alpha_1}) \cap df(C_{\alpha_2})$, \\
$B = \dim df(C_{\alpha_1}) \cap df(C_{\alpha_2}) \cap \text{Span}\{ \partial_{y_0} \}$, and \\
$C = \dim df(C_{\alpha_1}) \cap df(C_{\alpha_2}) \cap \text{Span}\{ \partial_{y_0}, \partial_{y_1}\}$.

\end{center}
\caption{Codimensions of 3D Multi-jet Strata}
\label{CodimOfMultijetStrata}
\end{table}%

Now consider the multi-jet stratum $S_{[001], 0} \times S_{[001],[01],0}$ intersected with the  stratum $S_{[001]}^{(2),0} \times S_{[011]}^{(2),0}$. Over this intersection we have  two rank two vector bundles  $C_1$ and $C_2$, whose images $df(C_1)$ and $df(C_2)$ form a rank two vector bundle and a line bundle, respectively. We will further stratify according to the dimension of the intersection of $df(C_1)$ and $df(C_2)$. There is an open substratum where they intersect trivially, and a codimension one stratum where they intersect non-trivially. See Table~\ref{CodimOfMultijetStrata}.

Finally over the intersection of $S_{[001]}^{(2),0} \times S_{[001]}^{(2),0} \times S_{[001]}^{(2),0}$ with $S_{[001], 0} \times S_{[001], 0}\times S_{[001], 0} $ there are three rank two vector bundles $df(C_1)$, $df(C_2)$, and $df(C_3)$. We will further stratify according to the dimension of their intersection in $TY$. There is a open substratum in which none of these 2-planes intersect. There is a codimension two substratum where two of the 2-planes coincide and there is a codimension 4 substratum in which all the 2-planes coincide. The total codimension of these strata are 9, 11, and 13, respectively, and since the dimension of $(\Sigma \times I)^{(3)}$ is 9, only the first occurs. 

Motivated by the above analysis of 3-dimensional multi-jet singularities, we introduce the following abstract definition of a 3-dimensional graphic. A generic map $f: \Sigma \times I \to \R^2 \times I$ then induces such a graphic.

\begin{definition} \label{3dGraphicDefn}
A {\em 3-dimensional graphic} 
\index{singularity!graphic}
\index{graphic}
is a diagram in $\R^2 \times I$ (with standard coordinates $(y_0, y_1, y_2)$, the coordinate $y_2$ being thought of as the `time' coordinate) consisting of a finite collection of closed labeled embedded points, arcs, and surfaces, such that the following conditions are met:
\begin{itemize}
\item The points are labeled by one of the following four labels: $S_{[011],[100],0}$,  $S_{[001], [11], 0}$, $S_{[011], [001], 0}$, and $S_{[001], [01], 1, 0}$, together with an appropriate index. 
\item The  arcs are labeled by either $S_{[011], 0}$ or $S_{[001], [01], 0}$, an appropriate index,  and for  each arc the projection to $I$ is a local diffeomorphism.
\item The surfaces are labeled by $S_{[001], 0}$, an appropriate index,  and for each surface the projection to $I \times I$ (the $(y_1, y_2)$-plane) is a local diffeomorphism. 
\item In a neighborhood of each point in each labeled arc, there exist coordinates in which the arrangement of arcs and surfaces matches the arrangement of arcs and surfaces induced by the images of singular loci associated to $S_{[011], 0}$ and $S_{[001], [01], 0}$ singularities. See Propositions 
\ref{Prop2DMorseCoord}
and \ref{PropPathofCuspCoordinates} and Figure~\ref{CuspAnd2DMorseIn3DFig}.
\item In a neighborhood of each labeled point there exist coordinated in which the arrangement of arcs and surfaces matches the images of the singular loci associated to  $S_{[011],[100],0}$,  $S_{[001], [11], 0}$, $S_{[011], [001], 0}$, and $S_{[001], [01], 1, 0}$ singularities. See Propositions \ref{Prop2DMorseRelationCoord}
\ref{PropCuspInversionCoord}, \ref{PropCuspFlipCoord}, and \ref{PropSwallowtailCoord} and Figures \ref{2DMorseRelationFig}, \ref{CuspInversionGraphicFig}, \ref{CuspFlipsGraphicFig}, and \ref{SwallowtailGraphicFig}.
\item The restrictions of the 3D graphic to the $y_2=0$ and $y_2=1$ planes consist of 2D graphics in the sense of Def.~\ref{2dGraphicDefn}.		
\item The arcs, points, and surfaces are in general position with respect to eachother and $\R^2 \times \partial I$. 
\item In particular when two surfaces intersect, it is along a 1-dimensional submanifold. We further require the the tangent space to this 1-dimensional intersection locus only intersect $\text{Span}\{ \partial_{y_0}, \partial_{y_1}\}$ non-trivially at a finite number of isolated points. We call these points {\em inversion points}. %The rest of this 1-manifold will also be considered an arc.  
%\item The points are disjoint from $\partial I^3$, the arcs are disjoint all of the boundary faces of the cube except the  $y_2=0$ and $y_2=1$ planes, the surfaces are disjoint from the $y_0=0$ and $y_0=1$ sides. 
%\item The boundary of the arcs is their intersection with the points and $\partial I^3$, and the boundary of the surfaces is their intersection with $\partial I^3$ and the arcs. 
\item The restriction of the graphic to $\R^2 \times \{0\}$ and $\R^2 \times \{1\}$ is a 2-dimensional graphic in the sense of Def.~\ref{2dGraphicDefn}
\end{itemize}
\end{definition}

\begin{theorem}
Let $\Sigma$ be a compact surface, and $f: \Sigma \times I \to \R^2 \times I$  a generic map. Then the images of the singular loci of $f$ form a {\em graphic}, as in the above definition.
\end{theorem}

\begin{remark}
	Several variants of the above notions are possible. First our analysis of the singularities of maps into the plane are equally valid with $\R^2$ replaced by any surface with a regular foliation. The leaves of the foliation play the role of the parallel translates of the $x$-axis. Similarly our 3-dimensional analysis is valid for targets consisting of 3-manifolds equipped with a complete flag of regular foliations. 
	
	Another important variation, which we will make use of in later chapters, is that we may allow the target surface to have boundary and to specify the restriction of the graphic to the boundary. This yields a relative notion of graphic. This allows us to replace $\R^2$ with $I^2$ and $\R^2 \times I$ with $I^3$. See the discussion leading up to Theorem~\ref{ThomTransversalityTheoremRelativeCorners}.
\end{remark}

A 3-dimensional graphic can be viewed as a sort of generalized path or `movie' of 2-dimensional graphics. For almost all times, the slice $\R^2 \times \{t\}$ of a 3-dimensional graphic will be a 2-dimensional graphic. As we progress through time this graphic will change by isotopy, except at a few critical time values. At these values the 2-dimensional graphic undergoes a topological change. These `movie moves' come from codimension 3 singularities, depicted in Table~\ref{LocalMovesForPlanarDiagramsTable1}, or from codimension 3 multijet singularities, depicted in Table~\ref{LocalMovesForPlanarDiagramsTable2}.  These table also list the number of indices for these singularities which dictates the number of normal forms for each of these kinds of operations.

\begin{table}[h]
\begin{center}
\begin{tabular}{|m{1.75in} | c |} \hline
Image of the Graphic & Number of Indices \\ \hline
\vspace{5pt}
\begin{tikzpicture} [scale = 0.75]
%% Cusp Flip A
\begin{scope}[xshift = 3cm, yshift = 6cm]
\draw (4,-5.5) rectangle (6,-3.5); 
\draw[red] (4.5, -5.5) -- +(0, 1) arc (180: 45: 0.25cm) -- +(0.25,-0.25) arc(45: 0: 0.25cm) arc (180: 135: 0.25cm) arc (-45: 0: 0.25cm) node (B) {} (B.center) |- (B.center |- 5.5, -3.5);
\end{scope}
\begin{scope}
\draw (4,0.5) rectangle (6,2.5); \draw[red] (4.5, 2.5) -- +(0, -1) arc (180: 315: 0.25cm) -- +(0.25,0.25) arc(-45: 0: 0.25cm) arc (180: 225: 0.25cm) arc (45: 0: 0.25cm) node (A) {} (A.center) |- (A.center |- 5.5, 0.5);
\end{scope}
\node at (6.5, 1.5) {$\leftrightarrow$};
\end{tikzpicture}
& 1 Index \\
 
 %Cusp flip B
\begin{tikzpicture}[xscale = -1,scale = 0.75]
\begin{scope}[xshift = 3cm, yshift = 6cm]
\draw (4,-5.5) rectangle (6,-3.5); 
\draw[red] (4.5, -5.5) -- +(0, 1) arc (180: 45: 0.25cm) -- +(0.25,-0.25) arc(45: 0: 0.25cm) arc (180: 135: 0.25cm) arc (-45: 0: 0.25cm) node (B) {} (B.center) |- (B.center |- 5.5, -3.5);
\end{scope}
\begin{scope}
\draw (4,0.5) rectangle (6,2.5); \draw[red] (4.5, 2.5) -- +(0, -1) arc (180: 315: 0.25cm) -- +(0.25,0.25) arc(-45: 0: 0.25cm) arc (180: 225: 0.25cm) arc (45: 0: 0.25cm) node (A) {} (A.center) |- (A.center |- 5.5, 0.5);
\end{scope}
\node at (6.5, 1.5) {$\leftrightarrow$};
\end{tikzpicture}
& 1 Index \\
\hline
\end{tabular}
\end{center}
\caption{Movie moves for 2D graphics coming from codim 3 singularities of 3D graphics.}
\label{LocalMovesForPlanarDiagramsTable1}
\end{table}%

\begin{table}[h]
\begin{center}
\begin{tabular}{|m{1.75in} | c |} \hline
Image of the Graphic & Number of Indices \\ \hline
\vspace{5pt}

%% Cusp intersect Fold A
\begin{tikzpicture} [scale = 0.75]
\draw (0,0) rectangle +(2,2); 
\draw [red] (0.5, 2) parabola (1, 1) parabola [bend at end] (1.5, 2);
\draw [red] (0.25, 2) -- (0.25, 1) to [out = -90, in = 90] (1.75, 0); 
\draw (3,0) rectangle +(2,2);
\draw [red] (3.5, 2) parabola (4, 0.5) parabola [bend at end] (4.5, 2);
\draw [red] (3.25, 2) to [out = -90, in = 90] (4.75, 1) -- (4.75, 0); 
\node at (2.5, 1) {$\leftrightarrow$};
\end{tikzpicture}
& 2 Indices \\

%% Cusp intersect Fold B
\begin{tikzpicture} [xscale = -1, scale = 0.75]
\draw (0,0) rectangle +(2,2); 
\draw [red] (0.5, 2) parabola (1, 1) parabola [bend at end] (1.5, 2);
\draw [red] (0.25, 2) -- (0.25, 1) to [out = -90, in = 90] (1.75, 0); 
\draw (3,0) rectangle +(2,2);
\draw [red] (3.5, 2) parabola (4, 0.5) parabola [bend at end] (4.5, 2);
\draw [red] (3.25, 2) to [out = -90, in = 90] (4.75, 1) -- (4.75, 0); 
\node at (2.5, 1) {$\leftrightarrow$};
\end{tikzpicture}
& 2 Indices \\

%% Cusp intersect Fold C
\begin{tikzpicture} [yscale = -1, scale = 0.75]
\draw (0,0) rectangle +(2,2); 
\draw [red] (0.5, 2) parabola (1, 1) parabola [bend at end] (1.5, 2);
\draw [red] (0.25, 2) -- (0.25, 1) to [out = -90, in = 90] (1.75, 0); 
\draw (3,0) rectangle +(2,2);
\draw [red] (3.5, 2) parabola (4, 0.5) parabola [bend at end] (4.5, 2);
\draw [red] (3.25, 2) to [out = -90, in = 90] (4.75, 1) -- (4.75, 0); 
\node at (2.5, 1) {$\leftrightarrow$};
\end{tikzpicture}
& 2 Indices \\
\hline
\end{tabular}
\end{center}
\caption{More movie moves for 2D graphics coming from codim 3 multijet singularities of 3D graphics.}
\label{LocalMovesForPlanarDiagramsTable2}
\end{table}%

\section{The Planar Decomposition Theorem} \label{SectPlanarDecompThm}
\index{planar decomposition theorem|(}

In this section we prove the main theorems of this chapter. Building on the results already obtained, we introduce the notion of a planar diagram. A planar diagram is essentially an arrangement of arcs and points in $\R^2$, together with a certain amount of combinatorial data. A generic map from a surface $\Sigma$ to $\R^2$ leads to a planar diagram, and moreover the surface can be completely recovered from the diagram, up to diffeomorphism. Similarly, we introduce the notion of a spatial diagram, which plays the same role for maps from $\Sigma \times I$ to $\R^2 \times I$. Two planar diagrams are defined to be equivalent if there is a spatial diagram extending the planar diagrams. We show, among other results,  that equivalence classes of planar diagrams are  in natural bijection with diffeomorphism classes of surfaces.

By the Thom transversality theorem, we know that there exist generic maps from any surface, $\Sigma$, into $\R^2$. What we have been studying so far has essentially been the local structure of these maps, and this gives us a complete, but local, understanding of our surface $\Sigma$. In the last section we introduced graphics, and in this section we will see how we can use the graphic to assemble these local descriptions into a global description of $\Sigma$. 

\subsection{Warm-up: linear decompositions of 1-manifolds}

\begin{figure}[htb]
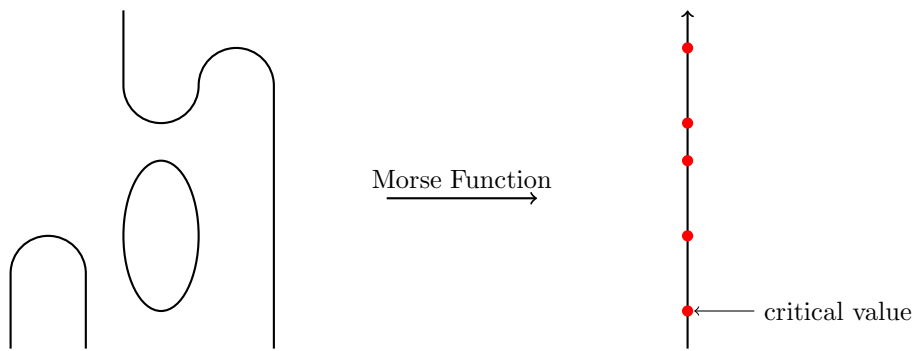

\begin{center}

\caption{A Morse Function for a 1-Manfold}
\label{FigMorseFunctionFor1Manifold}
\end{center}
\end{figure}

As we have done before, we will begin by describing the analogous situation one dimension lower. Suppose that $Y$ is a 1-dimensional manifold, equipped with a Morse function, $f$, as in Figure~\ref{FigMorseFunctionFor1Manifold}. For every point $r \in \R$ and every point $y \in Y$, such that $f(y) =r$, we know that there exist neighborhoods around $y$ and $r$ in which $f$ takes a particularly nice form. In particular $f$ will be a local diffeomorphism if $r$ is not a critical value for $f$, and if $r$ is a critical value for $f$, then there is precisely one critical point $y \in Y$ whose critical value is $r$. Around this point there are standard Morse coordinates in which $f$ takes the form of Equation \ref{MorseLemmaCoordEqn}. Around any other point in $f^{-1}(r)$, $f$ is a local diffeomorphism. 

\begin{figure}[htb]
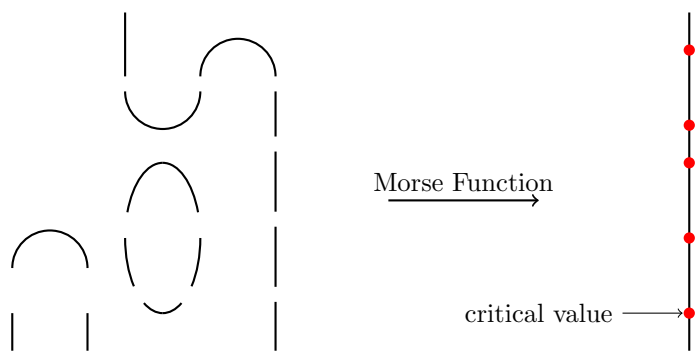

\begin{center}
%

\caption{A Linear Decomposition of a 1-Manifold}
\label{FigMorseFunctionDecomposition}
\end{center}
\end{figure}

Thus we may choose open neighborhoods covering $\R$ in such a way that the inverse images of these open sets decompose $Y$ into pieces, $U_i$, which consist of either several sheets mapping by a local diffeomorphism onto $Y$, or the disjoint union of such sheets with a single {\em elementary bordism}.
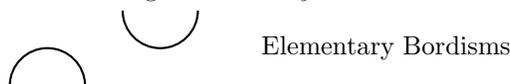
\begin{center}
%
\end{center}
In this manner we may write $Y$ as the union $ \cup U_i$, compare Figures \ref{FigMorseFunctionFor1Manifold} and \ref{FigMorseFunctionDecomposition}.

We may choose these open sets in $\R$ to be connected and to be as small and as well-shaped as we like. Those pieces which contain no critical points, and hence are simply a collection of sheets, may be (canonically) identified with $S_i \times f(U_i)$, where $S_i$ is the set of sheets and $f(U_i)$ is the image of $U_i$ in $\R$. The pieces which contain an isolated critical point have a similar description, with one factor being a copy of standard Morse coordinates. 

Given the decomposition of $\R$, together with sheets $S_i$ for each region, we would like to be able to reconstruct $Y$. For this we need to know how to glue the resulting $U_i$ together. The overlaps of any two $U_i$ are regions which contain no critical values, and hence give rise to a transition function $S_i \cong S_j$, which is an isomorphism of sets. Since $\R$ has covering dimension one, we may assume that all triple intersections are empty. 

For convenience, we can trivialize these sets $S_i$ as follows. 
Let ${\bf N}$ denote the ordered set $(1 \; 2 \; \dots \; N)$. We trivialize the sets $S_i$ by choosing an isomorphism $S_i \cong {\bf N}$, for some $N$. The elementary bordisms corresponding to critical points involve precisely two sheets, and so are determined by a set $S_i$ with a distinguished pair of elements and a labeling of  ``\tikz{\draw (0,0) -- +(-3pt, 0) arc (270:90: 3pt) -- +(3pt, 0); }'' or ``\tikz{\draw (0,0) -- +(3pt, 0) arc (-90:90: 3pt) -- +(-3pt, 0); }'' depending on the two possible indices of the critical point. %We can then choose an identification  $S_i \cong {\bf N}$. 
 %in which the distinguished pair of elements is $\{ N, N-1\}$. Thus 
 The two boundaries of such a region have sheets which are identified with either ${\bf N}$ or ${\bf N-2}$.  

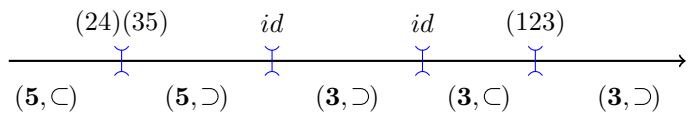
\begin{figure}[htb]
\begin{center}
\begin{tikzpicture}[thick]
	\draw [->]  (0, 1) -- +(9,0);
	\draw  [thin, blue, )-(] (1.5, 1.2) -- +(0, -.4);
	\draw  [thin, blue, )-(] (3.5, 1.2) -- +(0, -.4);
	\draw  [thin, blue, )-(] (5.5, 1.2) -- +(0, -.4);
	\draw [thin, blue, )-(] (7, 1.2) -- +(0, -.4);

	\node at (0.5,0.5) {$({\bf 5}, \tikz{\draw (0,0) -- +(-3pt, 0) arc (270:90: 3pt) -- +(3pt, 0); })$};
	\node at (2.5,0.5) {$({\bf 5}, \tikz{\draw (0,0) -- +(3pt, 0) arc (-90:90: 3pt) -- +(-3pt, 0); })$};
	\node at (4.5,0.5) {$({\bf 3}, \tikz{\draw (0,0) -- +(3pt, 0) arc (-90:90: 3pt) -- +(-3pt, 0); })$};
	\node at (6.25,0.5) {$({\bf 3}, \tikz{\draw (0,0) -- +(-3pt, 0) arc (270:90: 3pt) -- +(3pt, 0); })$};
	\node at (8.25,0.5) {$({\bf 3}, \tikz{\draw (0,0) -- +(3pt, 0) arc (-90:90: 3pt) -- +(-3pt, 0); })$};

	\node (A) at (1.5, 1.5) {$(24)(35)$};
	\node (B) at (3.5, 1.5) {$id$};
	\node (C) at (5.5, 1.5) {$id$};
	\node (D) at (7, 1.5) {$(123)$};
\end{tikzpicture}
\caption{A Diagram giving a Combinatorial Description of the 1-manifold from Figure~\ref{FigMorseFunctionFor1Manifold}.
}
\label{LinearDiagramFig}
\end{center}
\end{figure}

Together with these choices, a Morse function $f: Y \to \R$ gives us a combinatorial description of $Y$, which consists of the 1-dimensional graphic, the cover $\{U_i\}$, and the sheet data. Figure~\ref{LinearDiagramFig} gives an example of such a  diagram based off of the Morse function depicted in Figure~\ref{FigMorseFunctionFor1Manifold}. For simplicity in the regions labeled  {$({\bf N}, \tikz{\draw (0,0) -- +(-3pt, 0) arc (270:90: 3pt) -- +(3pt, 0); })$} and  {$({\bf N}, \tikz{\draw (0,0) -- +(3pt, 0) arc (-90:90: 3pt) -- +(-3pt, 0); })$} in Figure~\ref{LinearDiagramFig} it is assumed that the `distinguished pair of elements' corresponds to $\{ N, N-1\}$.
%We get a diagram on $\R$ by dividing it into regions labeled by either $({\bf N})$, $({\bf N}, \tikz{\draw (0,0) -- +(-3pt, 0) arc (270:90: 3pt) -- +(3pt, 0); })$, or  $({\bf N}, \tikz{\draw (0,0) -- +(3pt, 0) arc (-90:90: 3pt) -- +(-3pt, 0); })$. The boundaries of these regions are labeled by permutations of the appropriate sets. 

Conversely, given a graphic, a cover $\{U_i\}$, and compatible sheet data one can reconstruct a 1-manifold together with a Morse function inducing the original diagram. This 1-manifold can be constructed by gluing the basic pieces together exactly as the sheet data dictates. We invite the reader to check that Figure~\ref{LinearDiagramFig} gives a prescription which exactly reproduces the 1-manifold of Figure~\ref{FigMorseFunctionFor1Manifold}, up to diffeomorphism. In general we  have:

\begin{lemma}%\label{lem:}
	Let $(M,f)$ be a 1-dimensional compact manifold equipped with a Morse function. Let $\Psi$ be the 1-dimensional graphic of $f$ (see Def.~\ref{2dGraphicDefn}). Let $\{U_i\}$ be a cover of $\R$ (with no triple intersections) and let $\cS$ be a choice of sheet data with respect to $\Psi$ and $\{U_i\}$, as above. Let $M'$ be the manifold obtained from $(\Psi, \{U_i\}, \cS)$ by the above gluing prescription. Then there exists a diffeomorphism $M \cong M'$ over $\R$. \qed 
\end{lemma}

Thus we can recover $M$, up to diffeomorphism, from the combinatorial data $(\Psi, \{U_i\}, \cS)$. To complete our description of 1-manifolds in terms of this data, we must consider when two such collections correspond to the same 1-manifold. This is where Cerf theory provides an essential tool. Any two Morse functions on the same manifold $Y$ may be connected by a generic path of functions, which we consider as a map $Y \times I \to \R \times I$, commuting with the projection to $I$. We have seen how to derive normal coordinates for these sorts of functions, as in Propositions \ref{2dCuspProp}, \ref{2d2dmorseprop}, and \ref{2dFoldsProp}, and in this way we obtain a completely local description of such functions. 

Just as a Morse function provides a 1-dimensional graphic in $\R$, so too a generic map $Y \times I \to \R \times I$ provides a 2-dimensional graphic in $\R \times I$ (this kind of graphic will differ slightly from the kind considered in Def.~\ref{2dGraphicDefn} in that only cusp and crossing singularities are permitted, not 2D Morse singularities.)
%together with a suitable choice of cover $\{U_i\}$ and sheet data $\cS$, provided a diagram in the target $\R$, so too a generic map $Y \times I \to \R \times I$ provides a diagram in the target $\R \times I$, which we termed the {graphic} of the map. We build a diagram based on the graphic. 
 Our local description ensures that we may choose a covering of $\R \times I$ by open sets so that in each open set $f$ has a particular normal form, and again the inverse images over each such region are determined up to canonical isomorphism by a certain amount of combinatorial sheet data. 

Thus as a first approximation we could try to define an equivalence relation on triples $(\Psi, \{U_i\}, \cS)$ by declaring two triples to be equivalent if they extend to a common 2-dimensional graphic on $\R \times I$ with a cover and compatible sheet data. By construction such equivalence classes will be in natural bijection with diffeomorphism classes of 1-manifolds. However we would like to also understand this equivalence relation combinatorially. Unfortunately arbitrary open sets in $\R \times I$ can intersect in fairly complicated ways, and unless we restrict the allowed open sets of the cover of $\R \times I$, we will have difficulty making the equivalence relation a combinatorial one.   

One approach is to note that the covering dimension of $\R \times I$ is two, and so we may choose the covering so there are at most triple intersections. Moreover we may restrict ourselves to only considering convex open sets in our covering. One must also take some care to ensure that these intersections are arranged well with respect to the natural foliation of $\R \times I$ coming from the projection to $I$. In the end, with a little work, one obtains a simple finite collection of rules which generate the above equivalence relation on triples $(\Psi, \{U_i\}, \cS)$.

\subsection{Reconstruction via sheet data}\label{sec:planar:reconstruction}

A similar strategy can be employed in higher dimensions. 
In dimension two we are concerned not with Morse functions, but with maps to the plane, $\R^2$. We have already analyzed the singularities of these maps. Fix a surface $\Sigma$ and a generic map $f:\Sigma \to \R^2$. If $p \in \R^2$ is a point which is not a critical value of $f$, then $f$ is a local diffeomorphism around $p$. Propositions \ref{2dCuspProp}, \ref{2d2dmorseprop}, and \ref{2dFoldsProp} allow us to choose normal coordinates around any point which is a critical value. Theorem~\ref{ThmGraphicOfMaptoR2} ensures that the images of the singular loci will form a 2-dimensional graphic in $\R^2$ in the sense of Definition~\ref{2dGraphicDefn}. Thus we can proceed analogously to above, by choosing an open covering of $\R^2$ over which we have local models, and compatible sheet data explaining how to glue these local models together. 

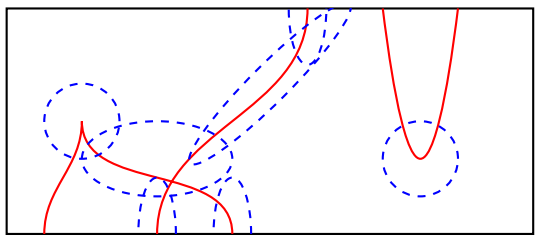
\begin{figure}[htb]
\begin{center}
\begin{tikzpicture}[thick]
\draw[clip] (0,0) rectangle (7,3);
\draw [red] (0.5, 0) to [out = 90, in = -90] (1,1.5) to [out = -90, in = 90] (3, 0)
	(2,0) to  [out = 90, in = -90]  (4,3)
	(5,3) parabola bend (5.5, 1) (6,3) ;
	 \draw[dashed, blue] (1, 1.5) circle (0.5 cm);
	\draw[dashed, blue] (5.5, 1) circle (0.5 cm);
	\draw[dashed, blue] (2, 1) ellipse (1cm and 0.5cm);
	\begin{scope}[xshift = 3.5cm, yshift = 2cm, rotate = 45] 
		\draw[dashed, blue] (0,0) ellipse (1.5cm and 0.25cm); 
	\end{scope}
		\draw[dashed, blue] (2, 0) ellipse (0.25cm and 0.75cm);
		\draw[dashed, blue] (4, 3) ellipse (0.25cm and 0.75cm);
		\draw[dashed, blue] (3, 0) ellipse (0.25cm and 0.75cm);
\end{tikzpicture}
\caption{Choice of Local Neighborhoods Decomposes $\R^2$}
\label{ChoiceofLocalNbhdsFig}
\end{center}
\end{figure}

Furthermore we can conduct a similar analysis one dimension higher. Two maps to $\R^2$ can be connected by a path of maps, which we consider as a single map $\Sigma \times I \to \R^2 \times I$, and which we can replace by a generic map. Again we have analyzed  the singularities of such maps. From a triple consisting of a 2-dimensional graphic, an open cover, and a compatible collection of sheet data we can construct a surface (with a generic map to the plane inducing the graphic). Two such triples produce diffeomorphic manifolds when they extend to a 3-dimensional graphic with a cover and compatible sheet data. In this section we will describe more precisely the sheet data and make precise these reconstruction results. 

Consider first the 2-dimensional case of a closed surface $\Sigma$ with a generic map to the plane $f: \Sigma \to \R^2$. Suppose that $U \subseteq \R^2$ is an open set and that the intersection of $U$ with the singular locus of $f$ is empty. Then $f$ restricts to a local diffeomorphism from $f^{-1}(U)$ to $U$. Suppose that $U$ is contractible, then we further have
\begin{equation*}
	f^{-1}(U) \cong S \times U
\end{equation*}
as smooth manifolds over $U$, where $S$ is the discrete set of fibers over any point in $U$, the `sheets' over $U$. Thus the manifold $f^{-1}(U)$ is completely determined up to diffeomorphism over $U$ by the combinatorial data of the set $S$. 

Now suppose that $U$ intersects the singular locus of $f$ non-trivially. For each point in the singular locus of the manifold we may choose normal coordinates as in Propositions \ref{2dCuspProp}, \ref{2d2dmorseprop}, and \ref{2dFoldsProp}. These normal coordinates identify a neighborhood of the given singular point in $M$ with a standard manifold (in fact in each case it is a disk with a specific map to $\R^2$) and this identification is over $\R^2$. By shrinking $U$ if necessary, we have that the inverse image $f^{-1}(U)$ consists of exactly this normal neighborhood and some disjoint `sheets' where $f$ restricts to a local diffeomorphism. We will now describe the combinatorial data that we will need to specify in order to make this identification a canonical diffeomorphism over $\R^2$. We assume that $U$ is contractible to avoid any `monodromy' issues. 

\begin{figure}[htb]
	\begin{center}
	\begin{tikzpicture}[thick]
	\node (A) [draw, ellipse, rotate = -30, dashed, blue, minimum height=3cm, minimum width=2cm]  at (0,0) {};
	\draw [red] (A.260) to [out = 30, in = -90] (A.90);
	\node at (1,-1) {$S$};
	\node at (-1.5,1) {$S \sqcup T $};

	\node (B) [draw, ellipse, rotate = -30, dashed, blue, minimum height=3cm, minimum width=2cm]  at (5,0) {};
	\draw [red] (B.300) to [out = 120, in = -90] (B.90);
	\draw [red] (B.180) to (B.0);
	\node at (5,2) {$S \sqcup T \sqcup T'$};
	\node at (3,-1) {$S \sqcup T$};
	\node at (7,0) {$S  \sqcup T'$};
	\node at (6,-1) {$S $};
	\end{tikzpicture}
	\end{center}
	\caption{Fold and Fold Crossing Sheet Data}
	\label{FoldSheetDataFig}
\end{figure}
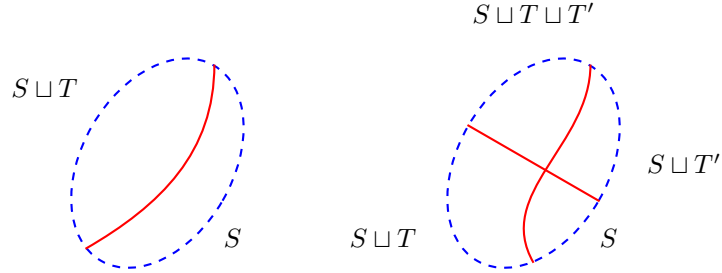

If $U$ is a region containing a portion of a single arc of the graphic, then that arc divides $U$ into two components, $U_0$ and $U_1$.  This in turn divides $f^{-1}(U)$ into two components, and away from the singular fold arc, $f$ restricts on $f^{-1}(U)$ to a local diffeomorphism. Thus, on each component we have  a canonical identification via $f$ with $S_i \times U_i$, where $S_i$ is the set of sheets over $U_i$. $S_0$ and $S_1$ differ by exactly two elements, in the sense that there is a two element set $T$ and a canonical identification $S_0 \cong  S_1 \sqcup T$ or a canonical identification $S_0 \sqcup T \cong S_1$. See Figure~\ref{FoldSheetDataFig}. 

If $U$ is a region which contains a point where two fold arcs cross, then the fold arcs divide $U$ into four regions. A similar analysis of sheet data applies to this situation. In this case there exist two 2-element sets, $T$ and $T'$, and four regions which have sheet data that can be identified with the sets $S$, $S \sqcup T$,  $S \sqcup T'$, and $S \sqcup T \sqcup T'$, respectively. The sheet data associated to a 2D Morse singularity is essentially the same as the data associated to the fold singularity. See Figure~\ref{FoldSheetDataFig}.

\begin{figure}[htb]
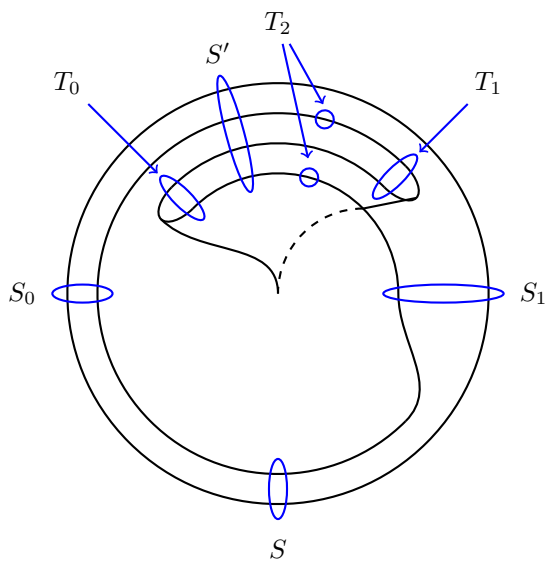

\begin{center}
%
\caption{Cusp Sheet Data}
\label{CuspSheetDataFig}
\end{center}
\end{figure}

If $U$ is a region containing a cusp singularity, then there exist two arcs of fold singularities within $U$. These divide $U$ into two regions, which have their associated sheet data $S$ and $S'$. Each of the fold singularities also has sheet data $(S_0, S_0 \sqcup T_0)$ and $(S_1, S_1 \sqcup T_1)$, where $T_0$ and $T_1$ are two element sets. Moreover we have isomorphisms $S \cong S_0$, $S \cong S_1$, and $S_0 \sqcup T_0 \cong S' \cong S_1 \sqcup T_1$, and under these identifications, $T_0$ and $T_1$ share one element in common, $\{ p\}$. Let $T_2 = T_0 \cup T_1 \setminus \{ p\}$. Then the following diagram commutes,
\begin{center}

\end{center}
where $\sigma: T_2 \to T_2$ is the non-trivial automorphism of this 2-element set. 
We call such data {\em cusp sheet data}, 
\index{cusp sheet data}
see Figure~\ref{CuspSheetDataFig}. It is determined by the sets $S, S'$ and the two maps $S \rightrightarrows S'$.

If the intersection of the graphic in $\R^2$ with $U$ is labeled by the above describe combinatorial data (the {\em sheet data}), then we can canonically recover $f^{-1}(U)$ up to diffeomorphism over $U$ from this data. To get a global identification of the surface $\Sigma$ over $\R^2$ we need to specify how to glue these local models together. This requires additional data ({\em gluing data}) on the intersections of the open sets.

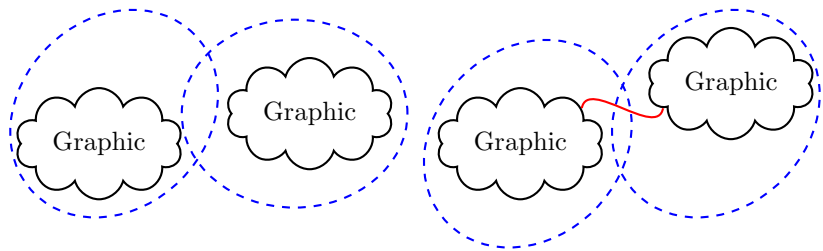
\begin{figure}[htb]
\begin{center}

\begin{center}
\begin{tikzpicture}[thick, scale=0.8]
	\node [rotate = 45, draw, ellipse, blue, dashed, minimum width = 3cm, minimum height = 2.5cm] at (0,0) {};
	\node [ draw, ellipse, blue, dashed, minimum width = 3cm, minimum height = 2.5cm] at (3,0) {};

	\node  [cloud, draw, aspect =2] at (-0.25, -0.5) {Graphic};
	\node  [cloud, draw, aspect =2] at (3.25, 0) {Graphic};	

	\node [xshift = 5.5cm, rotate = 45, draw, ellipse, blue, dashed, minimum width = 3cm, minimum height = 2.5cm] at (0,-0.5) {};
	\node [xshift = 8cm, rotate = 45, draw, ellipse, blue, dashed, minimum width = 3cm, minimum height = 2.5cm] at (0,0) {};
	\node (A) [cloud, draw, aspect =2] at (6.75, -0.5) {Graphic};
	\node (B) [cloud, draw, aspect =2] at (10.25, 0.5) {Graphic};

	\draw [red] (A.30) to [out = 80, in = -100] (B.200);
	
\end{tikzpicture}
\end{center}

\caption{Possible Normal Coordinate Intersections}
\label{PossibleNormalCoordinateIntersectionsFig}
\end{center}
\end{figure}

The covering dimension of $\R^2$ is two, and hence we may assume that these normal coordinate neighborhoods have at most triple intersections. Moreover, we may choose the open sets so that in these triple intersections there are are no critical values. Similarly, we may assume that the double intersections either contain no other critical points, or contain a portion of a single arc in the graphic, see Figure~\ref{PossibleNormalCoordinateIntersectionsFig}.

In the case that there are no critical points in the intersection, the gluing data amounts to an identification of the sheet data, and is thus given by an isomorphism of sets $S_0 \cong S_1$. In the second case we need to understand how to glue together elementary pieces which correspond to single fold singularities. Let the normal coordinate neighborhoods be $U_0$ and $U_1$. Recall that the sheet data for such singularities consists of a set of sheets $S_i$, and an additional pair of sheets $T_i$. Thus the fold locus splits, say, $U_0$ into two halves and the sheets over each half are identified with $S_1$ and $S_1 \sqcup T_1$, respectively. A straightforward calculation shows that a diffeomorphism
$f^{-1}(U_0) \cong f^{-1}(U_1)$, respecting the map to $\R^2$ is equivalent to an isomorphism of pairs of sets $(S_0, T_0) \cong (S_1, T_1)$. We call such data \emph{fold gluing data}. 

\begin{table}[ht]
\begin{center}
\begin{tabular}{|c|c|} \hline
Normal Coordinate Neighborhood & Sheet Data \\ \hline
 No Crit. Points & a set $S$ \\ \hline
 Fold & $(S, T)$,  $|T| =2$ \\ \hline
  2D Morse & $(S, T)$,  $|T| =2$ \\ \hline
 Fold Crossing & $(S, T_0, T_1)$,  $|T_i| =2$ \\ \hline
Cusp & $(S \rightrightarrows S')$, cusp sheet data \\ \hline \hline
Gluing Region & Gluing Data \\ \hline
 No Crit. Points & $S_0 \cong S_1$ \\ \hline
 Fold & $(S_0, T_0) \cong (S_1, T_1)$ \\ \hline
\end{tabular}
\end{center}
\caption{2D covering Sheet Data and Gluing Data}
\label{PlanarSheetDataAndGluingDataTable}
\end{table}%

Triple intersections of the normal coordinate neighborhoods only occur away from any critical points. On such triple intersections, the gluing data corresponding to the three adjacent double intersections must satisfy the obvious cocycle identity.

\begin{definition}\label{def:2Dsheetdata}
	Let $\Psi$ be a 2-dimensional graphic in $\R^2$ and let $\{U_i\}$ be a covering of $\R^2$ such that there are at most triple intersections, that the triple intersections contain no isolated points of $\Psi$ and such that the double intersections contain either no critical values of the graphic or contain a single arc.  A choice of {\em sheet data} $\cS$ with respect to $(\Psi, \{U_i\})$ is an assignment of sets and maps of sets to each open set and double intersection as in Table~\ref{PlanarSheetDataAndGluingDataTable}, and satifying the cocycle condition on triple intersections. 
\end{definition}

Given a triple $(\Psi, \{U_i\}, \cS)$ consisting of a 2-dimensional graphic $\Psi$, a cover $\{U_i\}$ as above, and a choice of sheet data $\cS$, then there exists a canonical surface $\Sigma$ with a generic map $f: \Sigma \to \R^2$ inducing the graphic $\Psi$. This surface is obtained by gluing together the local models determined by the sheet data $\cS$ over each $U_i$. The gluing is preformed exactly as dictated by the gluing data contained in $\cS$. We will say that $\Sigma$ is {\em constructed} from the triple $(\Psi, \{U_i\}, \cS)$. 

Our analysis up to this point yields:

\begin{lemma}\label{lem:Reconstruction2D}
	Let $\Sigma$ be a compact surface equipped with a generic map $f: \Sigma \to \R^2$, and let $\Psi$ be the 2-dimensional graphic corresponding to $f$. Then there exists a choice of cover $\{U_i\}$ and a choice of sheet data $\cS$ with respect to $(\Psi, \{U_i\})$, such that $\Sigma$ is diffeomorphic over $\R^2$ to the surface constructed from $(\Psi, \{U_i\}, \cS)$. \qed
\end{lemma} 

A similar analysis applies in the 3-dimensional setting. 
Any two generic  maps from a surface $\Sigma$ to $\R^2$ may be connected by a generic map $\Sigma \times I \to \R^2 \times I$. The results of the previous sections imply that we can choose coordinate neighborhoods in a manner which is compatible with the graphic of this map. We must determine the additional sheet data needed to identify the inverse images of these coordinate neighborhoods with standard elementary pieces. By the results of the last section, we may assume that each normal coordinate neighborhood is one of thirteen possible types, as listed in Table~\ref{3DSheetDataTable}. The analysis of the sheet data then occurs just as it did in the 2-dimensional setting.  

The no critical point, fold, 2D Morse, fold crossing, fold crossing inversion, cusp, and cusp inversion normal coordinate neighborhoods have sheet data exactly identical to their 2-dimensional analogs. The fold triple crossing sheet data is analogous to the fold crossing sheet data. The three sheets of folds divide the normal coordinate neighborhood into eight regions, and these regions have sheet assignments $S$, $ S \sqcup T_i$, $S \sqcup T_i \sqcup T_j$, or $S \sqcup T_0 \sqcup T_1 \sqcup T_2$.  Similarly the sheet data for normal coordinate neighborhoods containing an intersection point  of a 2D Morse path and fold surface or an intersection point  of a Cusp path and Fold surface are equally easy to derive. We leave the details to the reader. The sheet data for 2D Morse relation singularities is precisely the same as for a fold singularity, and likewise the sheet data for a cusp flip singularity is precisely the same as for a Cusp singularity.%\footnote{Note that the sheet data only depends on the class of the singularity in the multi-jet Thom-Boardman classification \cite{GG73}. } 

\begin{figure}[htb]
\begin{center}
\begin{tikzpicture}[thick, scale = 0.8]
	\node (A) at (0,0) {}; 
	\node (B) at (2,3) {}; 
	\node (C) at (2,-3) {}; 
	\node (D) at (4,0) {}; 
	
	\draw [red] (6, 4) to [out = -90, in = 45] (D.center) to [out = 225, in = 90] (C.center) to [out = 90, in = -90] (A.center) to [out = 90, in = -90] (B.center) to [out = -90, in = 135] (D.center) to [out = -45, in = 90] (6, -4); 
	\node [blue] (E) at (2, 0) {$S'$};
	\node [blue] (F) at (4, 2) {$S$};
	\node [blue] (G) at (4, -2) {$S$};
	\node [blue] (H) at (-1, 0) {$S$};
	\draw [blue, right hook->] (H) -- node [above] {$i_1$} (E);
	\draw [blue, left hook->] (G) -- node [below] {$i_2$} (E);
	\draw [blue, left hook->] (F) -- node [above] {$i_0$} (E);

	\node [blue] at (6, 0) {$i_0(S) \cap i_2(S) = S''$};
\end{tikzpicture}
\caption{Swallowtail Sheet Data}
\label{SwallowtailSheetDataFig}
\end{center}
\end{figure}
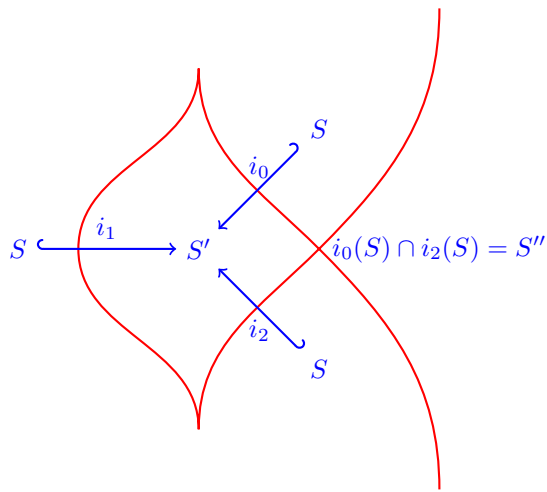

The remaining singularity, the swallowtail singularity, deserves more attention. The sheet data is determined by restricting attention to the graphic of Figure~\ref{SwallowtailSheetDataFig}, which exists on a hyperplane slice in a small neighborhood of the swallowtail singularity. There are three distinct regions, and so a priori there exist three sheet data $S, S', S''$. There are also five fold arcs which relate these regions via inclusions. Three of these are drawn in Figure~\ref{SwallowtailSheetDataFig}, and are labeled $i_o, i_1$, and $i_2$. Be virtue of the fold crossing singularity, we see that $i_0(S) \cap i_2(S) = S''$, so that this set and the remaining inclusions are already determined by the three maps $i_0, i_1, i_2: S \to S'$. By virtue of the cusp singularities, we see that the pairs $(i_0, i_1)$ and $(i_1, i_2)$ form cusp sheet data. Finally, by examining the geometry of the swallowtail singularity we see that $i_0(S) \cap i_2(S) = i_0(S) \cap i_1(S) \cap i_2(S)$. These conditions are equivalent to the following three conditions on the triple of maps $(i_0, i_1, i_2): S \to S'$:
\begin{enumerate}
\item $i_0(S) \cap i_2(S) = i_0(S) \cap i_1(S) \cap i_2(S)$,
\item $(i_0, i_1)$ is valid cusp sheet data,
\item $(i_1, i_2)$ is valid cusp sheet data.
\end{enumerate}
We call such a triple {\em swallowtail sheet data},
\index{swallowtail sheet data}
 $(S \triplearrow S')$.
 
\begin{table}[ht]
\begin{center}
\begin{tabular}{|c|c|} \hline
Normal Coordinate Neighborhood & Sheet Data \\ \hline \hline
No Crit. Points & $S$, set \\ \hline
Fold & $(S, T)$, sets, $|T|=2$ \\ \hline
Fold Crossing & $(S, T_0, T_1)$, sets, $|T_i|=2$ \\ \hline
Fold Crossing Inversion & $(S, T_0, T_1)$, sets, $|T_i|=2$ \\ \hline
Fold Triple Crossing & $(S, T_0, T_1, T_2)$, sets, $|T_i|=2$ \\ \hline
2D Morse & $(S, T)$, sets, $|T|=2$ \\ \hline
2D Morse $\cap$ Fold & $(S, T_0, T_1)$, sets, $|T_i|=2$ \\ \hline
Cusp & $(S \rightrightarrows S')$, Cusp Data \\ \hline
Cusp $\cap$ Fold &  $(S \rightrightarrows S', T)$, Cusp Data and set, $|T| =2$\\ \hline
2D Morse Relation & $(S, T)$, sets, $|T|=2$ \\ \hline
Cusp Inversion &  $(S \rightrightarrows S')$, Cusp Data \\ \hline
Cusp Flip &  $(S \rightrightarrows S')$, Cusp Data \\ \hline
Swallowtail & $(S \triplearrow S')$, Swallowtail Data \\ \hline
\end{tabular}
\end{center}
\caption{3D covering Sheet Data}
\label{3DSheetDataTable}
\end{table}%

Given the sheet data listed in Table~\ref{3DSheetDataTable}, we can construct over each open set of our cover a standard elementary piece. We need to understand how to glue these together, and this requires additional data associated to the intersections of the covering sets. The covering dimension of $\R^2 \times I$ is three, so we may assume our cover has at most 4-fold intersections. Moreover the triple intersections come in two types: those with no singularities and those with a single fold surface. The double intersections come in five kinds. There may be no critical points, a single fold surface, a single fold crossing curve with its accompanying fold surfaces, a single 2D Morse curve and its accompanying fold surfaces, or a single cusp curve and its accompanying surfaces. Determining what this gluing data amounts to in terms of the sheet data is again straightforward. The only non-trivial case are the intersections of the cover which contain curves of cusp singularities. Suppose that we are considering such an open set of our cover, and suppose that we are given two sets of sheet data $(i_0, i_1): S \rightrightarrows S'$ and $(\overline i_0, \overline i_1): \overline S \rightrightarrows \overline S'$. Thus there exist two elementary pieces over this open set. The isomorphisms of these elementary pieces, relative to the map to $\R^2 \times I$, are given precisely by the pairs of maps $(f, f'): (S, S') \to \overline S, \overline S')$ such that the following pair of diagrams commutes,
\begin{center}
\begin{tikzpicture}[thick]
	\node (LT) at (0,1.5) 	{$S$ };
	\node (LB) at (0,0) 	{$\overline S$};
	\node (RT) at (2,1.5) 	{$S'$};
	\node (RB) at (2,0)	{$\overline S'$};
	\draw [->] (LT) --  node [left] {$f$} (LB);
	\draw [->] (LT) -- node [above] {$i_0$} (RT);
	\draw [->] (RT) -- node [right] {$f'$} (RB);
	\draw [->] (LB) -- node [below] {$\overline i_0$} (RB);
	
	\node (LTA) at (4,1.5) 	{$S$ };
	\node (LBA) at (4,0) 	{$\overline S$};
	\node (RTA) at (6,1.5) 	{$S'$};
	\node (RBA) at (6,0)	{$\overline S'$};
	\draw [->] (LTA) --  node [left] {$f$} (LBA);
	\draw [->] (LTA) -- node [above] {$i_1$} (RTA);
	\draw [->] (RTA) -- node [right] {$f'$} (RBA);
	\draw [->] (LBA) -- node [below] {$\overline i_1$} (RBA);
\end{tikzpicture}
\end{center}
We call  such a pair $(f, f')$ {\em cusp gluing data}. The gluing data for all the double intersections is listed in Table~\ref{Table3DGluingData}

\begin{table}[ht]
\renewcommand{\arraystretch}{1.5}
\begin{center}
\begin{tabular}{|c|c|} \hline
Normal Coordinate Intersection & Gluing Data \\ \hline \hline
No Crit. Point & $S \cong \overline S$, set isom. \\ \hline
Fold & $(S, T) \cong (\overline S, \overline T)$, pair of set isom. \\ \hline
2D Morse & $(S, T) \cong (\overline S, \overline T)$, pair of set isom. \\ \hline
Fold Crossing & $(S, T_0, T_1) \cong (\overline S, \overline T_0, \overline T_1)$, triple of set isom. \\ \hline
Cusp & $(f, f'):  (S, S') \to \overline S, \overline S')$, cusp gluing data \\ \hline
\end{tabular}
\end{center}
\caption{3D Gluing Data}
\label{Table3DGluingData}
\end{table}%

In order for the gluing to be well defined, on triple intersection the gluing data must satisfy the obvious cocycle conditions. On 4-fold intersections there are no additional conditions.

\begin{definition}
	Let $\Psi$ be a 3-dimensional graphic in $\R^2 \times I$ and let $\{U_i\}$ be an open cover of $\R^2 \times I$ such that there are at most 4-fold intersections. Moreover we suppose that
	\begin{itemize}
		\item the 4-fold intersections are disjoint from the singular locus of $\Psi$;
		\item the triple intersections are disjoint from the isolated points and arcs of the singular locus of $\Psi$ and contain at most a single component of the of the 2-dimensional singular locus of $\Psi$; and
		\item the double intersections are disjoint from the isolated points of $\Psi$.
	\end{itemize}
A choice of \emph{sheet data} $\cS$ with respect to $(\Psi, \{U_i\})$ is an assignment of sets and maps of sets to each open set as per Table~\ref{3DSheetDataTable}, and an assignment of isomorphisms of these collections of sets for each double intersection are per Table~\ref{Table3DGluingData}, subject to the cocycle condition on triple intersections. 	 
\end{definition}

Just as in the 2-dimensional case the above date is enough to reconstruct $\Sigma \times I$. Given a triple $(\Psi, \{U_i\}, \cS)$ consisting of a 3-dimensional graphic in $\R^2 \times I$, a cover $\{U_i\}$ as above, and a choice of sheet data $\cS$ with respect to $(\Psi, \{U_i\})$, these exists a canonical 3-manifold of the form $\Sigma \times I$ together with a generic map $\Sigma \times I \to \R^2 \times I$ inducing the given graphic $\Psi$. As before $\Sigma \times I$ is obtained by gluing the local models over each $U_i$ together according to the gluing data of $\cS$. We will say that $\Sigma \times I$ is \emph{constructed} from  $(\Psi, \{U_i\}, \cS)$. We have:

\begin{lemma}\label{lem:Reconstruction3D}
	Let $\Sigma$ be a compact surface with a generic map $f: \Sigma \times I \to \R^2 \times I$ inducing the 3-dimensional graphic $\Psi$. Then there exists a choice of cover $\{U_i\}$ and a choice of sheet data $\cS$ with respect $(\Psi, \{U_i\})$, such that $\Sigma \times I$ is diffeomorphic over $\R^2 \times I$ to the 3-manifold constructed from $(\Psi, \{U_i\}, \cS)$. Moreover these choices of open cover and sheet data $(\Psi, \{U_i\}, \cS)$ may be made to extend any such choices already provided on the boundary $\R^2 \times I$, provide these boundary conditions satisfy the conclusion of Lemma~\ref{lem:Reconstruction2D}. \qed 
\end{lemma}

\begin{corollary}\label{cor:Reconstruction3Drel}
	Suppose we are given two triples $(\Psi, \{U_i\}, \cS)$ and $(\Psi', \{U'_i\}, \cS')$ each consisting of a 2-dimensional graphic, an open cover as in Definition~\ref{def:2Dsheetdata}, and sheet data with respect to these. Let $\Sigma$ and $\Sigma'$, respectively, be the surfaces constructed from these triples. Then $\Sigma$ is diffeomorphic to $\Sigma'$ if and only if there exists a triple  $(\Psi'', \{U''_i\}, \cS'')$ consisting of a 3-dimensional graphic in $\R^2 \times I$, a cover as above, and a choice of sheet data with respect to these which restricts to $(\Psi, \{U_i\}, \cS)$ and $(\Psi', \{U'_i\}, \cS')$ on the boundary. \qed
\end{corollary}

Thus we have partially succeeded in providing a 2-dimensional combinatorial model of surfaces. Any surface can be encoded in the largely combinatorial data of the triple $(\Psi, \{U_i\}, \cS)$ consisting of the graphic, the cover,  and sheet data. Moreover the previous corollary tells us when two such triples give equivalent surfaces. We can think of this a providing an equivalence relation on 2-dimensional triples $(\Psi, \{U_i\}, \cS)$, two such triples are equivalent if they extend to a common 3-dimensional triple.  

The difficulty with this approach is in expressing this equivalence relation combinatorially. Potentially, this could be accomplished by requiring the open sets to be convex and also requiring certain regularity conditions for the intersections of the open sets with respect to the foliation of $\R^2 \times I$. However in what follows we will employ a slightly different tactic. We do this for two reasons (1) making precise the nature of the regularity conditions is somewhat cumbersome, and (2) we would like to establish a simple description that meshes well with the higher categorical framework developed in the next two chapters. 

\subsection{Chambering graphs and foams}

Here we will introduce the notion of a {\em chambering graph} for a 2-dimensional graphic. Briefly, this is special sort of graph embedded into $\R^2$ which is transverse to the given graphic $\Psi$.  A chambering graph will cut the plane into regions, the {\em chambers}, and hence gives rise to a cover (though no longer by open sets). Rather than considering general open covers, as we did in the last section, we will effectively be limiting ourselves to just covers coming from a chambering graph. We will call the combined package consisting of a 2-dimensional graphic, a chambering graph, and compatible sheet data a {\em planar diagram}. These will be the topic of the next section. As expected, from a planar diagram we will be able to construct a surface and a map to $\R^2$ reproducing the graphic. 

Furthermore, we will have an analogous construction in the 3-dimensional case. In this section we will introduce the notion of a {\em chambering foam} which is a sort of higher dimensional version of a graph. It cuts $\R^2 \times I$ into regions (again called {\em chambers}). A 3-dimensional graphic, plus a chambering foam, plus compatible sheet data will constitute a {\em spatial diagram}. These will be defined in the next section as well. There are also relative versions of these notions. Two planar diagrams will give rise to diffeomorphic surfaces if and only if they extend to a common spatial diagram. Moreover, in contrast to the previous section, a spatial diagram always decomposes as a finite sequence of simple elementary operations. This gives us the desired 2-dimensional combinatorial model for surfaces which we have been after. 

As mentioned, a chambering foam is a higher dimensional version of a graph. It will consist of vertices, edges, and surfaces between the edges. In short it is a stratified space.  There are several competing notions of stratified space available in the literature, which have varying degrees of generality. Here, however, we will only need the most simplistic sorts of stratified spaces, those which are \emph{locally conical}. A locally conical stratified space is built from a collection of local models, much like a manifold is built from Euclidean space. The local models consist of products of cones on the lower dimensional locally conical stratified spaces. Thus the notion is inductive. Moreover it is possible and quite useful to restrict to only specific specified local models. 

\begin{definition}
	A \emph{stratified space of dimension $n$} is a topological space $X$ equipped with a filtration by closed subsets
	\begin{equation*}
		\emptyset = X_{-1} \subseteq X_0 \subseteq X_1 \subseteq X_2 \subseteq \cdots \subseteq X_n = X 
	\end{equation*}
	such that for each $0 \leq k \leq n$ the space $X_k \setminus X_{k-1}$, the $k^\text{th}$ \emph{stratum}, is given the structure of a smooth manifold (without boundary) of dimension $k$. An \emph{isomorphism} of stratified spaces is a homeomorphism of topological spaces preserving the filtration which restricts to a diffeomorphism on each stratum. We will say $X$ is of \emph{compact type} if each filtration $X_i$ is compact. 
\end{definition}

When we say that a map from a stratified space to a manifold has a certain property (smooth, transverse to a given submanifold, etc) we will mean that its restriction to each stratum has that property.
 
A smooth manifold is a stratified space in which all the lower dimensional strata are empty. A manifold with boundary or more generally a manifold with corners is a stratified space in an obvious way. Any open subset of a stratified space is a stratified space, and the product of two stratified spaces is a stratified space.  
 
\begin{example}
	 Suppose that $X$ is a non-empty stratified space of dimension $n$. Then the {\em open cone} on $X$ is the quotient space $CX = X \times [0,1) / X \times \{0\}$. This is a stratified space of dimension $n+1$, where $(CX)_0$ is the cone point and  for $1 \leq k \leq n+1$ $(CX)_k = X_{k-1} \times (0,1)$.  By convention we declare $C(\emptyset) = pt$. 
\end{example}

\begin{definition}
	A stratified space $X$ of dimension $n$ is \emph{locally conical} if for every point $x \in X$
there exists an integer $m \geq 0$ and a set $\{P_i\}$ of locally conical stratified spaces for $0\leq i \leq m$, where $P_i$ has dimension $n_i < n$, and there exists a neighborhood of $x$ which is isomorphic as a stratified space to: 
	\begin{equation*}
		\R^k \times CP_0 \times \cdots CP_m.
	\end{equation*}
	Here $k$ is determined by the equality 
	\begin{equation*}
		n = k  + \sum_{i=0}^m (n_i+1).
	\end{equation*}
(this puts restrictions on the possible $n_i$). We may also restrict the allowed local models. 	

A \emph{system of local models} for a locally conical stratified space of dimension $n$ is a set $\sU = \{ \cU_i\}$ where each element is itself a set $\cU_i = \{P_{i,j}\}$ of locally conical stratified spaces as above. We may equivalent specify $\sU$ by the local models corresponding the each $\cU_i$. 
We say is $X$ is \emph{locally conical with respect to $\sU$} if the local neighborhoods around each point of $x$ are chosen from $\sU$. By convention we will always assume that $\sU$ contains the set $\{ \emptyset\}$, and hence that $\R^n$ is always an allowed local model. 
\end{definition}

Part of the structure of a chambering graph or a chambering foam will consist of a locally conical stratified space which is embedded into $\R^2$ or $\R^2 \times I$, respectively. Moreover we will require that these stratified spaces are locally conical with respect to a very particular system of local models. This system is carefully chosen so that the combinatorial description we obtain is both simple and closely aligned with the higher categorical structures we will encounter in the subsequent chapters. 

\begin{figure}[htbp]
	\begin{center}
		\begin{tikzpicture}
			\node at (1,0) {$\R^1$};
			\draw [chambering] (1,0.5) -- (0.5,2.5);
			
			\node at (3,0) {$C_3$};
			\node at (3,-0.5) {cone on 3 points};
			\draw [chambering, fill] (3,1.5) circle (2pt);
			\draw [chambering] (3, 1.5) -- (2.5,2.5) (3, 1.5) -- (3.5,2.5) (3, 1.5) -- (3,0.5);
			\node at (6,0) {$C_1$};
			\node at (6, -0.5) {cone on 1 point};
			\draw [chambering, fill] (6,1.5) circle (2pt);
			\draw [chambering] (6,1.5) -- (6, 0.5);
		\end{tikzpicture}		
		\caption{The local models for chambering graphs.}
		\label{fig:LocalmodelsChamberingGraph}
	\end{center}
\end{figure}
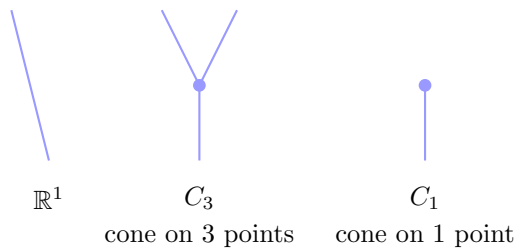

We will now describe the local models which will be used to construct chambering graphs and foams, but first we start with a few observations. 
A stratified space of dimension zero is just a 0-dimensional manifold. A locally conical stratified space of dimension 1 is a topological graph (a.k.a. a 1-complex). It is locally modeled on $\R$ and on the cones over 0-dimensional manifolds. For chambering graphs we will require that they only have trivalent and univalent vertices. Thus the local models consist of $\R$ and the cones on one point  and three point sets ($C_1$ and $C_3$ respectively), as depicted in Figure~\ref{fig:LocalmodelsChamberingGraph}.

\begin{figure}[htbp]
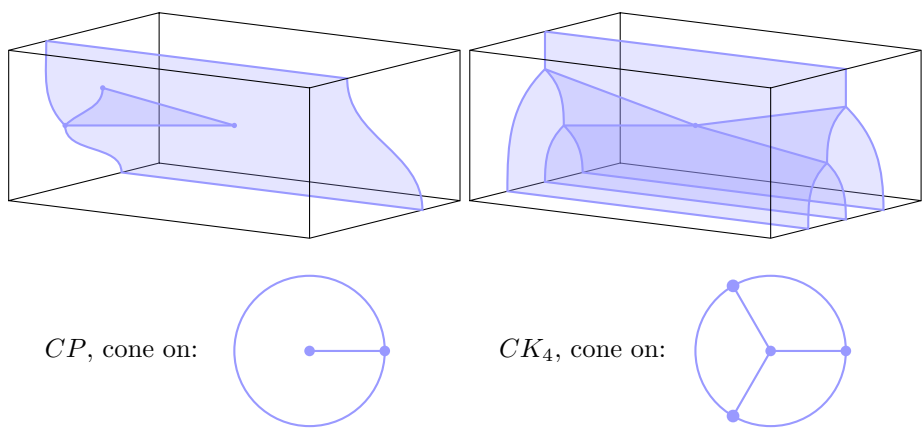

	\begin{center}

	\end{center}
	\caption{The local models for chambering foams.}
	\label{fig:LocalModelChamberingFoam}
\end{figure}

The local models for chambering foams are given by products of the previous models with $\R$, i.e., $\R^2$, $\R^1 \times C_1$, and $\R^1 \times C_3$, together with two additional models, depicted in Figure~\ref{fig:LocalModelChamberingFoam}. These are cones on two specific locally conical stratified spaces of dimension 1: the \tikz[baseline=0.075cm]{\fill (0,0) circle (1pt) (0.1,0.2) circle (1pt); \draw (0,0) -- (0.1,0.2) (0.2,0.2) circle (0.1cm);} graph (`$P$') and the complete graph on four vertices $K_4$. In the categorical language that we will develop in later chapters these local models will be responsible for unitality and associativity relations, respectively.   

We are now ready to give the definition of chambering graphs and foams. 
 
\begin{definition} \label{def:chamberinggraph}
	Let $\Psi$ be a 2-dimensional graphic in $\R^2$. A {\em chambering graph} $\Gamma$ for $\Psi$ consists of a 1-dimensional stratified space $\Gamma$ of compact type which is locally conical with respect to the system of local models comprised of $\R$, $C_1$, and $C_3$ (i.e. $\Gamma$ is a graph\footnote{Strictly speaking, $\Gamma$ could fail to be a graph in the usual sense in that it could have loops as well as components which are circles (with no vertices). However these phenomena are ruled out by the additional conditions we impose on $\Gamma$, and so there is no harm in thinking of $\Gamma$ as a graph in the traditional sense.} with only univalent and trivalent vertices), together with a closed embedding of $\Gamma$ as a subspace of $\R^2$ such that:	
	\begin{enumerate}
		\item The embedding of $\Gamma$ is smooth and in general position with respect to the arcs and vertices of the graphic $\Psi$. Recall that this means this condition is satisfied on each of the strata of $\Gamma$. In particular the vertices of $\Gamma$ are disjoint from the arcs and vertices of the graphic and the edges of $\Gamma$ intersect only the arcs of $\Psi$ transversely.		
		\item The following \emph{arrangement conditions} are satisfied:
		\begin{enumerate}
			\item For each edge in $\Gamma$ (i.e. connected component of the 1-dimensional stratum), the projection from that edge to the $y$-axis of $\R^2$ is a local diffeomorphism. In otherwords the projection from the 1-dimensional stratum of $\Gamma$ to the $y$-axis of $\R^2$ is a map without singularities. 
			\item Let $p \in \Gamma$ be a $C_3$ point (i.e. trivalent vertex). By the previous condition each of the three edges leaving $p$ leaves either upward or downward with respect to the $y$-coordinate of $\R^2$. We require that for each $C_3$ point (trivalent vertex) there is at least one of each, an upward and a downward arc.  
		\end{enumerate}
	\end{enumerate}
\end{definition}

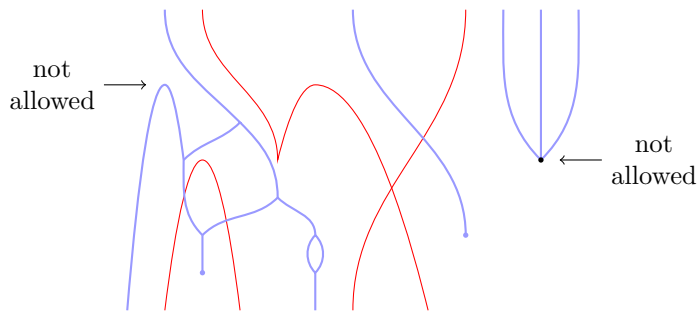
\begin{figure}[htbp]
	\begin{center}
	\begin{tikzpicture}[align=center]
		% immagine a box from (0,0) to (8,4)
		\node (A) at (-1, 3) {not \\ allowed};
		\node (B) at (7, 2) {not \\ allowed};
		\draw [->] (A) -- ++(1.25,0);
		\draw [->] (B) -- ++(-1.25,0);

		\draw [red] (0.5, 0) parabola bend (1,2) (1.5,0);
		\draw [red] (1,4) to [out = -90, in = 90] (2,2) parabola bend (2.5, 3) (4,0);
		\draw [red] (3,0) to [out = 90, in = -90] (4.5,4);
		
		\draw [chambering] (0,0) parabola bend (0.5,3) (0.75,2) to [out = -90, in = 135] (1,1); \draw [chambering] (1,1) -- (1, 0.5); \fill [chambering] (1, 0.5) circle (1pt);
		
		\draw [chambering] (5,4) to [out = -90, in = 135] (5.5,2) (5.5,4) -- (5.5,2) (6,4) to [out = -90, in = 45] (5.5,2); \fill (5.5,2) circle (1pt);
		\draw [chambering] (0.5, 4) to [out = -90, in =135] (1.5, 2.5) to [out = -135, in = 45] (0.75,2)
			(1.5, 2.5) to [out = -45, in = 90] (2, 1.5) to [out = -135, in = 45] (1,1)
			(2, 1.5) to [out = -45, in = 90] (2.5, 1) to [out = -135, in = 135] (2.5, 0.5)
			(2.5, 1) to [out = -45, in = 45] (2.5, 0.5) -- (2.5, 0);
		\draw [chambering] (3,4) to [out = -90, in = 90] (4.5, 1); \fill [chambering] (4.5,1) circle (1pt);

	\end{tikzpicture}
	\end{center}
	\caption{A chambering graph with two illegal points.}
	\label{fig:chamberinggraphdemo}
\end{figure}

Figure~\ref{fig:chamberinggraphdemo} shows a typical chambering graph with two illegal points. One of the illegal points violates condition (2)(a), while the other violates condition (2)(b). 

\begin{definition}\label{def:chamberingfoam}
	Let $\Psi$ be a 3-dimensional graphic in $\R^2 \times I$. A {\em chambering foam} $\Gamma$ for $\Psi$ consists of a 2-dimensional stratified space $\Gamma$ of compact type which is locally conical with respect to the system of local models comprised of $\R^2$, $\R \times C_1$, $\R \times C_3$, $CP$, and $CK_4$ (see Figure~\ref{fig:LocalModelChamberingFoam} for these last two), together with an embedding as a closed subspace of $\R^2 \times I$ such that:
	\begin{enumerate}
		\item The embedding of $\Gamma$ is smooth and in general position with respect to the vertices, arcs, and surfaces of the graphic $\Psi$, and well as the foliation of $\R^2 \times I$ induced by its projections to $\R \times I$ (the $yt$-coodinate plane) and to $I$. In particular this implies that:
		\begin{itemize}
			\item The 0-dimensional stratum of $\Gamma$ is disjoint from the vertices, arcs, and surfaces of the graphic $\Psi$.
			\item The 1-dimensional stratum of $\Gamma$ is disjoint from the vertices and arcs of $\Psi$, and intersects the surfaces of $\Psi$ transversely.
			\item the 2-diemsnional stratum of $\Gamma$ is dijoint from the vertices of $\Psi$ and intersects the arcs and surfaces of $\Psi$ transversely.  
		\end{itemize}
		\item The following \emph{arrangement conditions} are satisfied:
		\begin{enumerate}
			\item The projection from $\Gamma$ to $\R \times I$ (the $yt$-coodinate plane) has no singularities (which means that restricted to each strata of $\Gamma$ there are no singularities); The projection from the 2-dimensional stratum of $\Gamma$ to $I$ has no singularities. 
			\item  For each time $t \in I$ such that the slice $\R^2 \times \{t\}$ is disjoint from the vertices of the graphic $\Psi$, disjoint from the 0-dimensional stratum of $\Gamma$, and such that $t$ is not a critical value of the projection (from the strata) of $\Gamma$ to $I$, 		
			 we have that $\Gamma \cap \R^2 \times \{t\}$ is a chambering graph in $\R^2$ (with respect to the 2-dimensional graphic $\Psi \cap \R^2 \times \{t\}$). 
			\item Let $p\in \Gamma$ be a $CK_4$-point of the 0-dimensional stratum. In a neighborhood of $p$ there exist exactly four components, arcs, of the 1-dimensional stratum which converge to $p$ (these are the four arcs of $C_3$-points in the stratified space $CK_4$).  Under the above conditions, each of these arcs approaches $p$ from either the positive or negative time direction (i.e. the $I$-coordinate direction). We require that there is at least one arc for each direction, positive or negative. 
			\item Let $p \in \Gamma$ be a $CP$-point of the 0-dimensional stratum. In a neighborhood of $p$ there exist exactly two components, arcs, of the 1-dimensional stratum which converge to $p$ (these are the $C_3$ and $C_1$ ares in the stratified space $CP$). Under the above conditions, each of these arcs approaches $p$ from either the positive or negative time direction, as above. In this case we require that both arcs approach from the \emph{same} direction, either both positive or both negative. 
		\end{enumerate}
	\end{enumerate}	
\end{definition}

Figure~\ref{fig:LocalModelChamberingFoam} shows the typical allowed behavior near the two types of points of the 0-dimensional stratum of a chambering graph. For a $CK_4$-point there is another allowed configuration in which three of the $C_3$-arcs  point in the same direction. 

A 3-dimensional graphic $\Psi$ in $\R^2 \times I$ together with a chambering foam $\Gamma$ can be read along the time direction. For generic values of $t$ the slice $\R^2 \times \{t\}$ will contain a 2-dimensional graphic and a chambering graph. Thus the 3-dimensional graphic and foam can be viewed as a kind of generalized path (or `movie') between  2-dimensional graphics with chambering graphs. 

\begin{figure}[htbp]
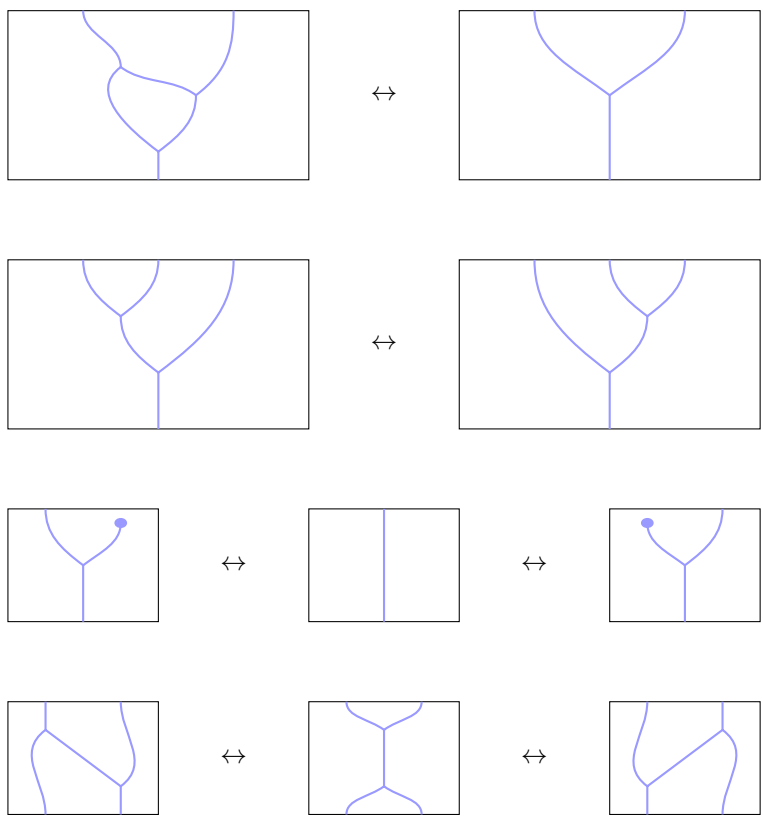

	 \begin{center}

	 \end{center}
	\caption{Movie moves for chambering graphs coming from vertices of a chambering foam.}
	\label{fig:0stratumChamberMoves}
\end{figure}

For most of the `movie' the 2-dimensional graphic and chambering graph only change by an isotopy. However at a certain times the graphic or the chambering graph will undergo a topological change. These `movie moves' happen locally and they occur in several varieties. First there are the changes which happen purely in the graphic. These have been analyzed carefully in previous sections and the results are tabulated in Table~\ref{LocalMovesForPlanarDiagramsTable1} and Table~\ref{LocalMovesForPlanarDiagramsTable2}. 

\begin{figure}[htbp]
	\begin{center}
		\begin{tikzpicture}[yscale = 0.75]
			\draw (0,8) rectangle (1.5,10) (3,8) rectangle (4.5,10) (5.5,8) rectangle (7,10) (8.5,8) rectangle (10,10);
			\draw [chambering] (0.75,10) -- (0.75,9.5) (0.75,8.5) -- (0.75,8);
			\draw [chambering, fill] (0.75,9.5) circle (2pt) (0.75,8.5) circle (2pt);
			
			\draw [chambering] (3.75,10) -- (3.75,8);
			
			\draw [chambering] (6.25,9.5) -- (6.25,8.5);
			\draw [chambering, fill] (6.25,9.5) circle (2pt) (6.25,8.5) circle (2pt);
			
			\draw (0,4) rectangle (1.5,6) (3,4) rectangle (4.5,6) (5.5,4) rectangle (7,6) (8.5,4) rectangle (10,6);
			\draw [chambering] (0.5,6) to [out = -90, in = 135] (0.75,5.5)
				(1,6) to [out = -90, in = 45] (0.75,5.5)
				(0.75, 5.5) -- (0.75, 4.5)
				(0.75,4.5) to [out = -135, in = 90] (0.5,4)
				(0.75,4.5) to [out = -45, in = 90] (1,4);
			\draw [chambering] (3.5,6) -- (3.5,4) (4,6) -- (4,4);
			\draw [chambering] (6.25,6) -- (6.25, 5.5) to [out = -135, in = 135] (6.25,4.5) -- (6.25,4)
				(6.25, 5.5) to [out = -45, in = 45] (6.25, 4.5);
			\draw [chambering] (9.25, 6) -- (9.25, 4);

			\node at (2.25,9) {$\leftrightarrow$};
			\node at (7.75,9) {$\leftrightarrow$};
			\node at (2.25,5) {$\leftrightarrow$};
			\node at (7.75,5) {$\leftrightarrow$};
			
		\end{tikzpicture}
	\end{center}
	\caption{Movie moves for chambering graphs coming from time-like critical points of a chambering foam.}
	\label{fig:localmoveschamberinggraph}	
\end{figure}
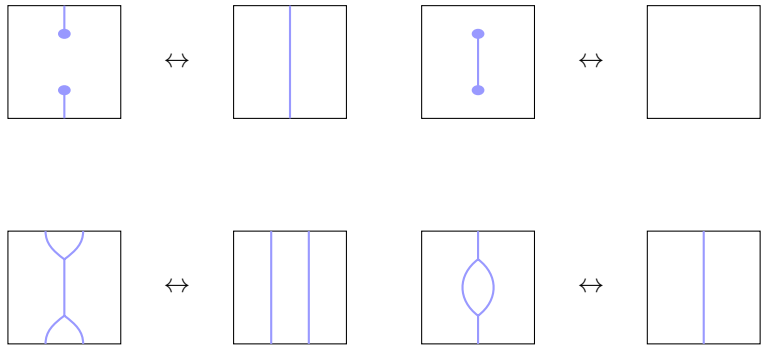

Then there are the changes which happen purely in the chambering graph. These arise from two sources, passing a vertex of the chambering foam (depicted up to reflection in Figure~\ref{fig:0stratumChamberMoves}) and passing a time-like critical point of the foam $\Gamma$, i.e. a critical point of the projection from $\Gamma$ to $I$. Such a critical point can only occur on the 1-dimensional stratum of $\Gamma$ and the possible movie moves are depicted, up to reflection, in Fugure~\ref{fig:localmoveschamberinggraph}. 

\begin{figure}[htbp]
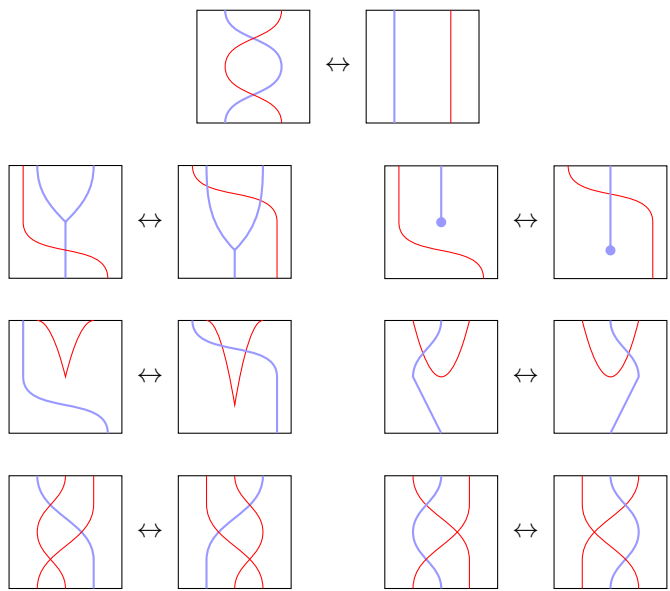

	\begin{center}
		%% Fold Crossing Inversion

	\end{center}
	\caption{Movie moves for intersection points between chambering graphs and 2D graphics.}
	\label{fig:chambergraphicmoviemoves}
\end{figure}

Finally there are movie moves which involve both the chambering graph and the graphic. These moves occur in three situation: where the interstection of the 2-dimensional strata of the 3D graphic and the chambering foam (a 1-dimensional intersection) obtains a time-like critical point, where the 1-dimensional stratum of the foam intersects the 2-dimensional surfaces of the 3D graphic, and where the 2-dimensional stratum of the foam intersects the either the 1-dimensional arcs of the 3D graphic or the 1-dimension intersection fold-crossing locus of the 3D graphic. These moves are depicted in Figure~\ref{fig:chambergraphicmoviemoves}. 

\begin{lemma}%\label{lem:}
	Let $\Psi$ be a 3-dimensional graphic in $\R^2 \times I$ which restricts to 2-dimensional graphics $\Psi_0$ and $\Psi_1$ on $\R^2 \times \partial I$. Let $\Gamma_0$ and $\Gamma_1$ be chambering graphs for $\Psi_0$ and $\Psi_1$ respectively. Then there exists a chambering foam $\Gamma$ for $\Psi$ whose restriction to $\R^2 \times I$ yields $\Psi_0$ and $\Psi_1$ respectively. 
\end{lemma}

\begin{proof} In a neighborhood of $\R^2 \times \{0\}$ (respectively $\R^2 \times \{0\}$) $\Psi$ is isotopic to a trivial product graphic `$\Psi_0 \times I$'. This it will suffice to assume that $\Psi$ is such a product graphic and to show that $\Gamma_0$ extends to a chambering foam whose restriction to $\R^2 \times \{1\}$ is empty. This can easily be achieved by applying the above movie moves. For example the following algorithm produces the desired chambering foam:
	\begin{enumerate}
		\item The first movie move in Figure~\ref{fig:chambergraphicmoviemoves} allows us to `cut' any edge of the chambering graph by introducing a pair of $C_1$-points. In the neighborhood of each intersection of an edge of the chambering graph with the graphic, we may introduce such a cut. Then we may slide on of the $C_1$-points past the arc of the graphic (this is the second movie move on the second row of Figure~\ref{fig:chambergraphicmoviemoves}). The result is a chambering graph which is disjoint from the 2-dimensional graphic. 
		\item Next we take each edge of the chambering graphic and we again cut each one by introducing a pair of $C_1$-points. What remains is a chambering graph consisting of a disjoint union of graphs of the following three types:
		
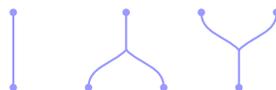
\begin{center}
		\begin{tikzpicture}
				\draw [chambering,fill] (0,1) circle (1pt) -- (0,0) circle (1pt);
				\draw [chambering] (1.5,1) -- (1.5,0.5) to [out = -135, in = 90] (1,0) (1.5,0.5) to [out = -45, in = 90] (2,0) ;
				\draw [chambering,fill] (1.5,1) circle (1pt) (1,0) circle (1pt)(2,0) circle (1pt);
				\draw [chambering] (3,0) -- (3,0.5) to [out = 135, in = -90] (2.5,1)  (3,0.5) to [out = 45, in = -90] (3.5,1);	
				\draw [chambering,fill] (3,0) circle (1pt) (2.5,1) circle (1pt) (3.5,1) circle (1pt);
		\end{tikzpicture}
		\end{center}  
		\item Each of these component graphs may be eliminated by applying the second movie move in Figure~\ref{fig:chambergraphicmoviemoves} and the third movie move in Figure~\ref{fig:0stratumChamberMoves}. \qedhere
	\end{enumerate}
\end{proof}

\begin{definition} \label{def:chambers}
	Let $\Psi$ be a 2-dimensional graphic in $\R^2$ (respectively a 3-dimensional graphic is $\R^2 \times I$) and let $\Gamma$ be a chambering graph (resp. chambering foam) with respect to $\Psi$. We define the \emph{chambers} of $(\Psi,\Gamma)$ to be the connected components of the compliment of $\Gamma$ and the points, arcs, (and surfaces) of $\Psi$ in $\R^2$ (resp. $\R^2 \times I$). Let $\cU = \{U_i\}$ be an open covering of $\R^2$ (resp. $\R^2 \times I$). We will say that $\Gamma$ is \emph{subordinate} to $\cU$ if each chamber is a subset of at least one $U_i$.  
\end{definition}

\begin{proposition}\label{prop:existencesubgraphicfoam}
	Let $\Psi$ be a 2-dimensional graphic in $\R^2$ (respectively a 3-dimensional graphic is $\R^2 \times I$) and let $\cU = \{U_i\}$ be an open covering of $\R^2$ (resp. $\R^2 \times I$). Then there exists a chambering graph (resp. chambering foam) $\Gamma$ with respect to $\Psi$ which is subordinate to $\cU$.  
\end{proposition}

\begin{proof}
	The proofs in the 2-dimensional and 3-dimensional cases are very similar.  We will explain the slightly more complicated 3-dimensional case. Moreover, without loss of generality we may assume that the open cover is locally finite. However we will only prove the result in the case that the open cover is actually finite. This is no real loss, first because for our applications in later chapters we will only need the finite cover case, and secondly because the proof below does also go through in the locally finite case, but with more careful bookkeeping.  We will try to sketch the necessary changes along the way.

	Let $H$ be a hyperplane in $\R^3 \supseteq \R^2 \times I$ and fix $\epsilon >0$ a small positive number. The hyperplane divides $\R^3$ into two half spaces and by taking $\epsilon$-neighborhoods we obtain a two set open cover $\cV^{H, \epsilon} = \{ V_+, V_-\}$. The intersection $V_+ \cap V_-$ is an $\epsilon$-neighborhood of the hyperplane $H$.  Given an open cover $\cW = \{W_i\}$, we may obtain a new open cover $\cW^H = \{ W_i \cap V_+, W_i \cap V_-\}$, with twice as many open sets, obtained by intersecting the open sets of $\cU$ with $V_+$ and $V_-$. This new open cover is subordinate to $\cV^{H, \epsilon}$ and $\cW$.

What we will show is that if $\Gamma$ is a chambering foam for $\Psi$ which is subordinate to $\cW$, and $H$ is suitably generic hyperplane, then we may alter $\Gamma$ to obtain $\Gamma^H$ such that $\Gamma^H$ is subordinate to $\cW^H$. Thus by repeating this operation finitely many times (for varying $H$ and $\epsilon$) we may obtain a chambering foam which subordinate to any finite cover. In the locally finite case we must preform this operation countably many times, and so an argument must be made that shows that, locally, the foam is still only altered finitely many times. In this way the limit will still be a chambering foam with the desired properties.  
	
Thus assume that $\Gamma$ is a chambering foam for $\Psi$ which is subordinate to a fixed initial cover $\cW$. Let $H$ be a hyperplane which is  (1) in general position with respect to the projection $\R^2 \times I \to \R \times I$  (by which we mean that it does not contain a line parallel to the $x$-axis), (2) in general position with respect to the points, arcs, and surfaces of the graphic $\Psi$, and (3) in general position with respect to the strata of the chambering foam $\Gamma$. The idea is that $\Gamma^H$ should be obtained by taking the union of $\Gamma$ and $H$. The components of the compliment of $\Gamma$, the points, arcs, and surfaces of $\Psi$, and $H$ (i.e. the chambers of this hypothetical $\Gamma^H$) are certainly subordinate to $\cW^H$ (for any choice of $\epsilon$). 

Unfortunately the union of $\Gamma$ and $H$ is \emph{not} a chambering foam. Specifically it fails to be a chambering foam precisely at the points of intersection between $H$ and $\Gamma$. To obtain the true $\Gamma^H$ we will have to modify/alter these intersection points so as to replace the union of $\Gamma$ and $H$ with an honest chambering foam. These modifications are local and can be preformed in an arbitrarily small neighborhood of of the intersection locus. By making these alterations in a suitably small neighborhood (depending on $\epsilon$) we will obtain a chambering graph $\Gamma^H$ which is subordinate to $\cW^H$ (for that fixed parameter $\epsilon$). 

There are three kinds of intersection points that we must discuss, corresponding to the intersection of $H$ with the 2-dimensional stratum of $\Gamma$, and the two kinds of 1-dimensional strata (types $C_1$ and $C_3$). These intersect $H$ along a 1-dimensional submanifold of $H$ and at isolated points, respectively. (In fact $\Gamma \cap H$ will be a 1-dimensional stratified space which is locally conical with respect to $C_1$ and $C_3$, just like a chambering graph). 

We will first discuss how to modify points of intersection away from the $C_1$ and $C_3$ intersection points. A time-like cross sectionof such a point  (i.e. intersection with $\R^2 \times \{t\}$ for some fixed $t\in I$) looks like the intersection of two transverse lines. It is a 4-valent vertex with two upward point arcs (positive $y$-direction) and two downward arcs (negative $y$-direction). We must resolve it into a trivalent graph, and in fact there is a canonical minimal way to do so. We resolve it so that the two upward arcs meet at a trivalent vertex and the two lower arcs meet at their own trivalent vertex, with a new edge passing between them. We can make this modification in a sufficiently small neighborhood so that the modified chambering foam is with $\epsilon$ distance of the hyperplane and the original foam. A cross section of this modification is depicted in the first line of Figure~\ref{fig:modifyingchamber+hyperplane}.

\begin{figure}[htbp]
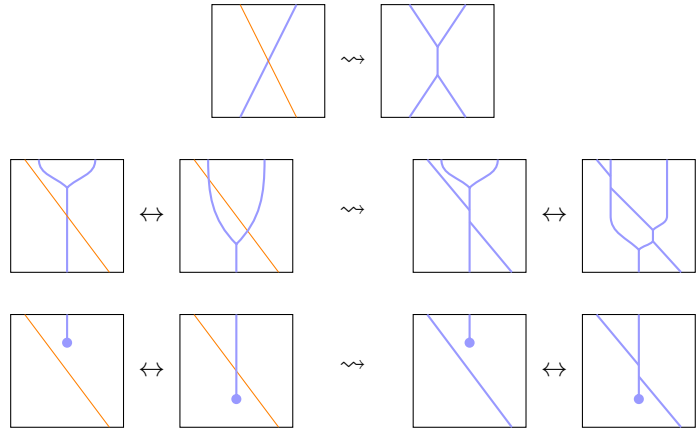

	\begin{center}

	\end{center}

	\caption{Modifying the intersection of a hyperplane and a chambering foam.}
	\label{fig:modifyingchamber+hyperplane}
\end{figure}
 
For the $C_3$ and $C_1$ intersection points we need a more complicated modification. A neighborhood of these points can be viewed as a movie which modifies the configuration of the hyperplane and the $C_3$- or $C_1$-points of a chambering graph. This movie depicts such a point sliding through the hyperplane, and these moves are depicted one the left-hand-side of the second and third row of Figure~\ref{fig:modifyingchamber+hyperplane}. After we resolve the non-$C_3$ and $C_1$ intersection points we obtain the chambering graphs depicted on the right-hand-side of Figure~\ref{fig:modifyingchamber+hyperplane}. Thus we must choose chambering foams which implement these `movie moves'. 

For the $C_3$ case this is easily done using the movie moves depicted in the last rows of 
Figure~\ref{fig:0stratumChamberMoves} and Figure~\ref{fig:localmoveschamberinggraph}. The $C_1$ case is equally easy, using the movie moves depicted in the last two rows of those figures. The precise sequences of moves is not unique and we leave it to the reader to decide the precise sequence, depending on taste. The sequence itself is not important, provided that it can be chosen so these operations stay with $\epsilon$-distance of the original chambering foam and hyperplane. This is clearly the case. 
\end{proof}

\begin{corollary}[of the proof]
		Let $\Psi$ be a 3-dimensional graphic in $\R^2 \times I$ with chambering foam $\Gamma$ and let $\cU = \{U_i\}$ be an open covering of $\R^2 \times I$. Let $(\Psi_0, \Gamma_0)$ and $(\Psi_1, \Gamma_1)$ be the 2-dimensional graphics and chambering graphs obtained by restricting $(\Psi, \Gamma)$ to $\R^2 \times \{0\}$ and $\R^2 \times \{1\}$ respectively. Similarly, let $\cU_0$ and $\cU_1$ be the restrictions of $\cU$ to $\R^2 \times \{0\}$ and $\R^2 \times \{1\}$ respectively. Suppose that $\Gamma_0$ and $\Gamma_1$ are subordinate to $\cU_0$ and $\cU_1$, respectively. Then there exists new  chambering foam $\Gamma'$ with respect to $\Psi$ which is subordinate to $\cU$, and whose restrictions to   $\R^2 \times \{0\}$ and $\R^2 \times \{1\}$ still yield $\Gamma_0$ and $\Gamma_1$, respectively. \qed
\end{corollary}

\subsection{Planar and Spatial Diagrams} \label{sec:planarandspacialdiagrams}
 With these considerations we introduce the following definition.

\begin{definition} \label{def:planardiagram}
	A {\em planar diagram} 
	\index{planar diagram} is a triple $(\Psi, \Gamma, \cS)$ consisting of a graphic $\Psi$ in $\R^2$ (as in Definition~\ref{2dGraphicDefn}), a chambering graph $\Gamma$ with respect to $\Psi$ (as in Definition~\ref{def:chamberinggraph}), and a collection of {\em sheet data} $\cS$. The sheet data consists of: 
	\begin{itemize}
		\item An assignment of a set to each \emph{chamber} (see Definition~\ref{def:chambers});
		 \item An assignment of sheet data to the arcs and isolated points of the graphic as in Table~\ref{PlanarSheetDataAndGluingDataTable};  
		\item An assignment of, for each component of the compliment of $\Psi$ in the 1-dimensional stratum of $\Gamma$, an isomorphism of sets between the two adjacent regions; and 
		\item An assignment of, for each intersection point between an edge of $\Gamma$ and $\Psi$, an isomorphism of pairs of sets $(S_0, T_0) \cong (S_1, T_1)$, as listed in the fold-gluing data entry of Table~\ref{PlanarSheetDataAndGluingDataTable}.
	\end{itemize}
For each of the points of the zero dimensional stratum of $\Gamma$ we further require a cocycle (or monodromy) condition: For a $C_1$-point we require the edge leaving the point is labeled with the identity isomorphism of sets, for a $C_3$-point we require that the circular composite of the three isomorphisms labeling the three edges leaving the point yields the identity isomorphism of sets.  	
\end{definition}

We also have a 3-dimensional analog of this notion. 

\begin{definition}
	A {\em spatial diagram} 
	\index{spatial diagram} is a triple $(\Psi, \Gamma, \cS)$
	consisting of a 3-dimensional graphic $\Psi$ (in the sense of Definition \ref{3dGraphicDefn}), a  chambering foam $\Gamma$ (in the sense of Definition~\ref{def:chamberingfoam}), and a collection of {\em sheet data} $\cS$. The sheet data consists of: 
	\begin{itemize}
		\item An assignment of a set to each \emph{chamber} (see Definition~\ref{def:chambers});
		 \item An assignment of sheet data to the surfaces, arcs, and isolated points of the graphic as in Table~\ref{3DSheetDataTable};  
		\item An assignment of, for each component of the compliment of $\Psi$ in the 2-dimensional stratum of $\Gamma$, an isomorphism of sets between the two adjacent regions; and 
		\item An assignment of sheet data to the intersection points of $\Gamma$ with $\Psi$ as in Table~\ref{Table3DGluingData}. 
	\end{itemize}
Along the 1-dimensional strata of $\Gamma$ we require the same cocycle conditions as in the definition of planar diagram. 	 
\end{definition}

A spatial diagram can be restricted to either of its two boundaries, and the last condition ensures that this yields a planar diagram.  This leads to a natural notion of equivalence between planar diagrams. 

\begin{definition}
Two planar diagrams are {\em equivalent} if there exists a spatial diagram whose restriction to the boundary agrees with the pair of planar diagrams. 
\end{definition} 
 
The notion of equivalence of planar diagrams is clearly reflexive and is symmetric by definition. It is transitive since we may concatenate spatial diagrams\footnote{To make the glued spatial diagram smooth, we may need to adjust the original spatial diagrams slightly. This is not a serious issue since we only need to show the relation is transitive.}, and thus forms an equivalence relation. The detailed analysis that we have preformed now pays off. We can understand very precisely when two planar diagrams are equivalent.

\begin{theorem}	\label{PlanarDecompMovesThm}
	Two planar diagrams are equivalent if and only if they can be related by a finite sequence of the following local moves: 
	\begin{enumerate}
		\item Isotopy;
		\item The changes in graphic shown in Table~\ref{LocalMovesForPlanarDiagramsTable1} and
		Table~\ref{LocalMovesForPlanarDiagramsTable2}, provided there exist compatible choices of sheet data as in Table~\ref{3DSheetDataTable};
		\item The changes in chambering graph depicted in Figure~\ref{fig:0stratumChamberMoves}, Figure~\ref{fig:localmoveschamberinggraph}, and Figure~\ref{fig:chambergraphicmoviemoves}, provided both sides of the move may be equipped with compatible choices of sheet data. \qed
	\end{enumerate}
\end{theorem}

In complete parallel to our discussion in Section~\ref{sec:planar:reconstruction}, from a planar diagram $(\Psi, \Gamma, \cS)$ we may obtain a canonical surface $\Sigma$ together with a generic map $\Sigma \to \R^2$ whose associated graphic is $\Psi$. The surface $\Sigma$ is obtained by gluing the local models, dictated by the graphic, together as dictated by the chambering graph and the sheet data. We will say that $\Sigma$ is {\em constructed} from the planar diagram $(\Psi, \Gamma, \cS)$. Similarly from a spatial diagram there is a canonical 3-manifold of the form $\Sigma \times I$ with a generic map to $\R^2 \times I$ (commuting with the projection to $I$) which reproduces the underlying 3-dimensional graphic. We will also say $\Sigma \times I$ is \emph{constructed} from the spatial diagram. 

\begin{theorem}[Planar Decomposition Theorem] \label{PlanarDecompositionTheorem} 
	\index{planar decomposition theorem}
	Let $\Sigma$ be a compact surface.
	\begin{enumerate}
		\item Let $f: \Sigma \to \R^2$ be a generic map to the plane, with 2-dimensional graphic $\Psi$. Then there exists a  chambering graph $\Gamma$ with respect to $\Psi$ and a choice of sheet data $\cS$, making the triple $(\Psi, \Gamma, \cS)$ a planar diagram. Moreover the surface constructed from this planar diagram is diffeomorphic to $\Sigma$ over $\R^2$. 
		\item Let $f: \Sigma \times I \to \R^2 \times I$ be a generic map (commuting with projection to $I$), with $3$-dimensional graphic $\Psi$. Then there exists a chambering foam $\Gamma$ and a choice of sheet data $\cS$, making the triple $(\Psi, \Gamma, \cS)$ a spatial diagram. Moreover the 3-manifold $\Sigma' \times I$ constructed from this spatial diagram is diffeomorphic to $\Sigma \times I$ over $\R^2 \times I$. 
		\item In the previous situation, if we are given planar diagrams on the boundary $\R^2 \times \partial I$ and isomorphisms of $\Sigma$ over $\R^2$ with the surfaces constructed from these planar diagrams, then the above spatial diagram may be chosen so as to extend these two planar diagrams. 
	\end{enumerate}
Consequently, diffeomorphism classes of compact surfaces are in natural bijection with 	equivalence classes of planar diagrams. 
\end{theorem}

\begin{proof}
	
	Suppose that $\Psi$ is a 2-dimensional graphic in $\R^2$, that $\cU = \{U_\alpha\}$ is an open cover of $\R^2$, and that $\overline{\cS}$ is a compatible collection of sheet data, and we considered in Section~\ref{sec:planar:reconstruction}. Suppose that $\Gamma$ is a chambering graph with respect to $\Psi$. The key observation is that if $\Gamma$ is subordinate to $\cU$, then we may obtain an induced planar diagram $(\Psi, \Gamma, \cS)$. This is done by choosing for each chamber a specific member of the cover $\cU$ containing that chamber. Moreover the surfaces constructed from these two sets of data, the triple $(\Psi, \cU, \overline{\cS})$ and the planar diagram, are naturally diffeomorphic over $\R^2$. A similar observation holds in the 3-dimensional case. 

Now the three enumerated claims follow from our previous considerations, specifically from Lemma~\ref{lem:Reconstruction2D}, Lemma~\ref{lem:Reconstruction3D}, Corollary~\ref{cor:Reconstruction3Drel}, and Proposition~\ref{prop:existencesubgraphicfoam}.
The last statement follows immediately from the three enumerated claims. 
\end{proof}

\index{planar decomposition theorem|)}

 %Planar Decomposition Theorem
\chapter{Symmetric Monoidal Bicategories} \label{SymMonBicatChapt}

The goal of this chapter is to introduce symmetric monoidal bicategories and to establish a collection of fundamental theorems about them. These results provide a package which makes it possible to manipulate symmetric monoidal bicategories effectively and to preform a variety of calculations. Throughout we will endeavor to emphasize those aspects of the theory which are analogous to results in topology, as the theory of symmetric monoidal bicategories shares many common elements with the theory of $E_\infty$-spaces (or more precisely $E_\infty$-$2$-types). In Section~\ref{sec:historysymmonbicat} we review the history of symmetric monoidal bicategories, explaining some of what has and what has not been done previously in the literature. In particular the 
definition of symmetric monoidal bicategory has historically been scattered in several pieces throughout the literature, with a few minor gaps.\footnote{After the original version of this work appeared in 2009, a complete account of the definition of symmetric monoidal bicategory did surface in the literature \cite{Stay:2013aa}.} 

In Section \ref{SectSmyMonBicat} we collect together these result and present the fully-weak definition of symmetric monoidal bicategory, as well as the accompanying definitions of symmetric monoidal homomorphism, symmetric monoidal transformation, symmetric monoidal modification, and their various compositions.  In Section \ref{SectSymMonWhisker} we introduce the operation of symmetric monoidal whiskering. To our knowledge this has never appeared in the literature, but an analogous operation for (non-symmetric) monoidal bicategories can be extracted from \cite{GPS95}. This is the key operation which allows one to construct a tricategory of symmetric monoidal bicategories (in the sense of   \cite{GPS95}).  See also Definition \ref{DefnWhiskering} in the appendix for a lower categorical version of this operation. 

In Section \ref{SectWhiteheadsTheoremforSymmetricMonoidalBicategories} we prove the first of three key theorems about symmetric monoidal bicategories, which we dub ``Whitehead's theorem for symmetric monoidal bicategories''. Whitehead's theorem in topology provides a recognition principle for when a map of (nice) topological spaces is a homotopy equivalence. Theorem~\ref{WhiteheadforSymMonBicats} provides a similar recognition principle for determining when a symmetric monoidal homomorphism is an equivalence of symmetric monoidal bicategories. This theorem is surely ``known to experts'', 
but as we were unable to find it in the literature 
%but as neither its proof nor its statement exists in the literature %in the form of Theorem~\ref{WhiteheadforSymMonBicats} 
we provide a complete proof here. To reduce the amount of coherence data one needs to check, we break the proof of Theorem~\ref{WhiteheadforSymMonBicats} into a sequence of six lemmas, and introduce the notion of skeletal symmetric monoidal bicategories, mimicking the well known definition of skeletal categories. 

Then we move on to the important topic of presentations for symmetric monoidal bicategories. These enjoy a duplicitous existence being, on the one hand, similar to the theory of presentations of groups (or any other algebraic structure) via generators and relations, and being, on the other hand, similar to CW-structures for topological spaces. The precise definition of presentation for a symmetric monoidal bicategory that we introduce is reminiscent of the definition of $n$-computad, first introduced in dimension two by R. Street in \cite{MR0401868} to formalize the notion of presentation for strict $2$-categories. The notion of $n$-computad was rediscovered by A. Burroni \cite{MR1224519} some 20 years later under the moniker `$n$-polygraph'. Batanin  generalized the notion of $n$-computad to any categorical structure controlled by a globular operad \cite{MR1664991}, which, although quite general, does not include the theory of symmetric monoidal bicategories. More recently Batanin has generalized the notion of computad to arbitrary finitary monads on globular sets \cite{MR1664991}, which includes both the symmetric monoidal bicategory monad and the ordinary bicategory monad.
% This allows us to develop the basics of the theory, while avoiding a general treatment of the auxiliary theory of globular operads. 
We will describe a variation on Batanin's construction in the next section, Section~\ref{sec:computads}.

In this generality a {\em symmetric monoidal 3-computad} corresponds to the data of a presentation for a symmetric monoidal bicategory. There is an adjunction between our symmetric monoidal computads and symmetric monoidal bicategories, and thus every presentation gives rise to a symmetric monoidal bicategory, the free symmetric monoidal bicategory on that presentation. The symmetric monoidal bicategories which arise this way enjoy many special properties. We call them {\em computadic}.
The counit of this adjunction is an equivalence of symmetric monoidal bicategories, and thus every symmetric monoidal bicategory has a functorial {\em computadic replacement}. In fact there is a Quillen model structure on symmetric monoidal bicategories (and strict functors) and this also serves as a functorial cofibrant replacement in terms of that structure. 

In Sections~\ref{SectFreelyGenSymMonBicats} and \ref{sec:cofibrancyThm} we establish our second key result, the Cofibrancy Theorem (Thm.~\ref{thm:cofibrancythm}). We show that computadic symmetric monoidal bicategories satisfy an important universal property: the bicategory of symmetric monoidal homomorphisms out of a cofibrant symmetric monoidal bicategory is equivalent to the bicategory given by specifying the value of the homomorphism (or transformation of modification) on each generator, subject to the relations. In particular every weak symmetric monoidal functor out of a cofibrant symmetric monoidal bicategory is equivalent to a strict one, and every symmetric monoidal transformation between strict symmetric monoidal functors out of a cofibrant symmetric monoidal bicategory is equivalent to a strict transformation. The precise details involve a lengthly inductive analysis.

Our third and final theorem for symmetric monoidal bicategories is a very general coherence theorem. Coherence theorems, which often take the equivalent form of strictification theorems, have been well studied in the literature. They usually come in two varieties. The first, most common, variety asserts that every weak categorical object is weakly equivalent to a stricter version of that categorical object. A typical example of this kind is the coherence theorem of Gordon-Powers-Street \cite{GPS95} which asserts that every tricategory is equivalent to a Gray-category, a much stricter version of tricategory. In these cases one also has a comparison functor which witnesses the equivalence. 

The second kind of coherence theorem compares the free weak categorical object generated by a particular kind of datum to the free stricter categorical object generated by that same datum. This kind of coherence theorem asserts that the natural comparison map is an equivalence. Some examples of this sort are \cite[Thm.~1.2]{JS93} for monoidal categories freely generated by a given category and \cite[Cor.~1.9]{GO13} for symmetric monoidal bicategories freely generated by a set. 

In higher category theory one often needs to compare diagrams involving complicated composites of higher morphisms, and in practice one would like to have an efficient notation for these composites. For example the composites of 2-morphisms in bicategories are often described using pasting diagrams (see Appendix~\ref{PastingsStringsAdjointsAndMatesSection}). At first inspection, pasting diagrams appear ill-founded and without a well-defined value. To give meaning to a pasting diagram one must choose ways of parenthesizing the morphisms involved in the pasting diagram and must also insert various associators from the bicategory to make the 2-morphisms composable. It is a consequence of MacLane's coherence theorem (Thm.~\ref{MacLanesCoherenceTheorem}) that these choices don't effect the value of the composite 2-morphism, and hence pasting diagrams have a well-defined meaning. 

Coherence theorems of the second sort have immediate corollaries of the above kind. Any diagram constructed using the stricter categorical calculus can be lifted in an essentially unique way to a diagram in the calculus of the weaker categorical structure. However this second kind of coherence theorem is also limited in so far as it only applies to diagrams coming from the kinds of free objects under consideration. This is a mild annoyance for the theory of bicategories, but is a serious problem for the theory of symmetric monoidal bicategories. 

The current state-of-the-art in coherence theorems for symmetric monoidal bicategories are either of the first kind above, or, in \cite{GO13}, of this second kind for symmetric monoidal bicategories freely generated by a bicategory. Unfortunately most symmetric monoidal bicategories are very far from symmetric monoidal bicategories of this kind, and most diagrams one would like to consider don't arise from symmetric monoidal bicategories freely generated by a bicategory. For example such a free symmentric monoidal bicategory has no non-trivial morphisms to or from the unit object. In the case of the bordism category, an example relevant to topological field theories, this rules out, for example, diagrams involving closed manifolds. 

We fill this gap. In Section~\ref{sec:strictsymbicats} we introduce the notion of {\em quasistrict} symmetric monoidal bicategory, in which most, but not all, of the coherence data is trivial. This is expected as it is not possible to fully strictify symmetric monoidal bicategories. This notion is algebraic and there is a corresponding monad on 2-globular sets. The material of Section~\ref{sec:computads} applies to this notion as well to give a corresponding notion of quasistrict symmetric monoidal computad. Every 
%semi-strict symmetric monoidal bicategory is, in particular, a symmetric monoidal bicategory, which may equivalently be expressed as a map of monads. Hence, every 
symmetric monoidal computad induces a quasistrict symmetric monoidal computad, and there is a comparison homomorphism between the corresponding symmetric monoidal bicategories: 
\begin{equation*}
	\sF(P) \to \sF_\text{q.s.}(P).
\end{equation*}
In Section~\ref{sec:coherence} we show that for every symmetric monoidal 3-computad (i.e. symmetric monoidal presentation) this comparison map is an equivalence of symmetric monoidal bicategories. 

This has several important consequences. For example, as every symmetric monoidal bicategory is equivalent to a computadic symmetric monoidal bicategory, we learn that every symmetric monoidal bicategory is equivalent to a quasistrict symmetric monoidal bicategory. This is a coherence theorem of this first kind. The functorial computadic replacement of a arbitrary symmetric monoidal bicategory is quite large and corresponds to a presentation in which every object, 1-morphism, and 2-morphism of the original bicategory is used as a generator. The functorial quasistrictification is correspondingly large. However our Coherence Theorem provides a more nuanced result. If one happens to have a smaller cofibrant replacement (i.e. symmetric monoidal equivalence $\sF(P) \simeq \cC$ for some smaller presentation $P$), then one may obtain a correspondingly smaller strictification. We also give an example showing that these results are optimal. It is not possible to find a stronger general strictification theorem.

Let us remark that this more nuanced coherence result is absolutely essential for our applications to the classification of topological field theories in the next chapter. Any coherence theorem which asserts that symmetric monoidal bicategories are equivalent to a functorially constructed strictification will not be adequate as any {\em functorial} strictification is necessarily large and useless for the purpose of classifying topological field theories. The result we establish here permits us to find a handicrafted, but not functorial, strictification of the bordism bicategory which leads to a classification of topological field theories. The situation is in precise analogy with topological spaces. Any functorial CW-approximation of topological spaces (such as the geometric realization of the singular simplicial set) is necessarily extremely large and mostly of theoretical importance. For doing computations one usually selects an explicit small CW-structure, which typically cannot be constructed functorially.

These results, the theory of symmetric monodidal computads, Whitehead's theorem, the Cofibrancy theorem, and the Coherence theorem, provide a useful toolkit for manipulating symmetric monoidal bicategories and proving results about them. Analogous results for the theory of bicategories have long been available and are used implicitly throughout the literature. For example, MacLane's coherence theorem implies that a well-typed pasting diagram has a unique value once the boundaries of the pasting diagram have been bracketed. However this alone is not enough to interpret the pasting diagrams used by Gordon-Powers-Street \cite{GPS95} in the definition of tricategory. This requires a stronger coherence result, such as an analog of the cofibrancy theorem. The same issues arise in the definition of symmetric monoidal bicategory and so we discuss this is some detail in Section~\ref{sec:coherenceinaction}.

In the next chapter we will use these results in a different way to classify extended 2D topological field theories. We will describe the symmetric monoidal bicategory $\bord_2$ of 2-dimension bordisms, and topological field theories are given by symmetric monoidal functors out of this. A presentation of the bordism bicategory, that is a symmetric monoidal 3-computad $P$ and an equivalence $\sF(P) \to \bord$ (i.e., a cofibrant replacement), gives a way to classify functors out of the bordism bicategory. From the Cofibrancy theorem of Section~\ref{sec:cofibrancyThm} such functors are given by specifying the value of each generator, subject to the relations. This is analogous to topological fact that maps out of CW-complexes can be easily describe cell by cell. However if we try to use the canonical functorial cofibrant replacement of $\bord$, then we gain nothing as this presentation would use every manifold and bordism as generators. Similarly the functorial CW-approximation of a space is much too large for computations. Its importance is of a more theoretical nature. 

Instead we will provide a handcrafted and much smaller presentation, $P$. There will be an obvious symmetric monoidal functor $\sF(P) \to \bord$, and we will be tasked with showing this is an equivalence of symmetric monoidal bicategories. A direct approach, using say Whitehead's theorem for symmetric monoidal bicategories, is not practical, since the symmetric monoidal bicategory $\sF(P)$ is quite complicated. However we will be able to show this indirectly by comparing two equivalent symmetric monoidal bicategories. The Planar Decomposition Theorem~\ref{PlanarDecompositionTheorem} of the last chapter will enable us to replace the bordism bicategory with an equivalent bicategory whose 2-morphisms are given by certain planar string diagrams. This replacement will then be directly comparable to a strictified version of $\sF(P)$.

\section{The Tortuous History of Symmetric Monoidal Bicategories} \label{sec:historysymmonbicat}

\index{symmetric monoidal bicategory!history|(}

The mathematical axiomatization of symmetric monoidal bicategories goes back to the work of Kapranov and Voevodsky \cite{KV94, KV94-2} on braided monoidal 2-categories. They proposed a definition of braiding for strictly associative monoidal 2-categories. There was a minor omission in this original definition of braided (strict) monoidal 2-category, which was repaired in the work of Baez and Neuchl \cite{BN96}, where the definition is also simplified and put into a more conceptual context. 

It was further clarified in the work of Day and Street \cite{DS97}. They explain how the categorification of monoidal, braided monoidal, and symmetric monoidal categories gains an additional layer. There are monoidal 2-categories, braided monoidal 2-categories, {\em sylleptic} monoidal 2-categories and, finally, symmetric monoidal 2-categories. The definition was further modified by Crans \cite{MR1626844}, who added a normalization condition which corrected a small error in one of the proofs of \cite{BN96} (see \cite{BL03}).

Just as symmetric monoidal categories are braided monoidal categories which satisfy additional axioms, symmetric monoidal bicategories are sylleptic monoidal bicategories satisfying additional axioms. This is in contrast to the relationship between braided and sylleptic monoidal bicategories: a sylleptic monoidal bicategory not only satisfies additional axioms compared to a braided monoidal bicategory, but it is equipped with additional structure, a {\em syllepsis}. 
Again, all this was carried out using a partially strict notion of monoidal 2-category called a {\em Gray monoid}. The authors justify this by invoking the coherence theorem of Gordon, Powers and Street. % \cite{GPS95}. 

Gordon, Powers, and Street, in \cite{GPS95}, introduced what is essentially a fully weak notion of tricategory and accompanying notions of trihomomorphism, tritransformation, and higher morphisms. Mirroring the well known definition of a (weak) monoidal category as a bicategory with one object, they define a monoidal bicategory to be a tricategory with one object. This is the most common form encountered in examples, and it is essentially the definition we use below. These authors also prove a coherence theorem, which ensures that any such monoidal bicategory is equivalent to a Gray monoid. 
The notion of tricategory introduced by Gordon, Powers, and Street is {\em non-algebraic} in the sense that the forgetful functor to 3-globular sets does not admit a left adjoint. An essentially equivalent algebraic version of tricategory was given in N. Gurski's Ph.D. dissertation \cite{MR2717302}. 

In light of the Gordon-Powers-Street coherence theorem, Day and Street defined braided, sylleptic and symmetric monoidal bicategories only in the context of Gray monoids. To quote them,
\begin{quote}
Examples naturally occur as monoidal bicategories rather than Gray monoids. However, the coherence theorem of \cite{GPS95} allows us to transfer our definitions and results. ...
\end{quote}
This transference was carried out explicitly for braided and sylleptic monoidal bicategories in part of the thesis work of P. McCrudden and occurs in the appendices of \cite{McCrudden00}. The symmetric monoidal bicategories of primary interest to the current work are, naturally, of this fully weak kind. In \cite{Shu1004} Shulman gives a method which easily produces many examples of symmetric monoidal bicategories. Further examples appear in Stay \cite{Stay:2013aa}, where the definition of symmetric monoidal bicategory is also spelled out in full detail and in one place. 

The Gordon-Powers-Street coherence theorem for tricategories also gives a coherence theorem for monoidal bicategories \cite{GPS95}. It concludes that every monoidal bicategory is equivalent to a Gray monoid, and is a coherence theorem of the `first kind' in the terminology of the introduction to this chapter. More refined coherence theorems for tricategories appear in Gurski's work \cite{MR2717302, MR3076451}. Gurski has also proven coherence theorems for braided monoidal bicategories \cite{MR2770448}, and, with Osorno, for symmetric monoidal bicategories \cite{GO13}. These include coherence results of the `second kind' for free symmetric (resp. braided) monoidal bicategories generated by ordinary bicategories. 

When the original version of the current work appeared in 2009, there still remained two gaps, which were relatively straightforward to fill, but which were missing from the literature. First, the axiom that a sylleptic monoidal bicategory must satisfy in order to be symmetric must be transfered from the partially strict version presented in \cite{DS97} to the fully weak context. The simplicity of this axiom makes this completely straightforward. Second, while monoidal bicategories can be defined as tricategories with a single object this does not give the correct notion of monoidal transformation and monoidal modification. Indeed, if one were to adopt these as the correct notion of higher morphisms, then there would be an additional categorical layer: permutations between the modifications, as given in \cite{GPS95}. This situation mirrors that for monoidal categories.

Monoidal categories can be defined as single object bicategories and homomorphisms between these result in monoidal functors. However transformations between these homomorphism are not the same as monoidal natural transformations between monoidal functors, and the modifications provide an additional categorical layer not usually discussed in the context of monoidal categories (cf. \cite{MR2342826, MR2839900}).  

	These two points of view can be reconciled by observing that the single object bicategory corresponding to a monoidal category is canonically {\em pointed}, and this should be considered as part of the structure. Homomorphisms between single object bicategories are automatically pointed homomorphisms, as expected. However, pointed transformations precisely reproduce monoidal natural transformations, and moreover all pointed modifications between single object bicategories are necessarily trivial. In this way monoidal categories form a precise subset of the theory of bicategories, see Appendix  \ref{MonoidalCatsAsBicats}. 
	
A similar discussion applies to monoidal bicategories, which should be considered as canonically pointed single object tricategories. Again all homomorphisms between these are automatically pointed and there are no non-identity pointed permutations. We take as a definition that monoidal transformations and monoidal modifications should not be given as those in \cite{GPS95}, but rather by their pointed analogs.  We will only be considering pointed homomorphism, transformations and modifications which occur between single object tricategories and this considerably simplifies the relevant diagrams presented in \cite{GPS95}. We provide these simplified diagrams below.  Additionally, the following tables should help in the translation between the various structures defined in the aforementioned works. 

\begin{table}[ht]
\begin{center}
\begin{tabular}{|c|c| c | c|} \hline
 Author & Monoidal $( \clubsuit)$ & Braided $( \diamondsuit)$ & Sylleptic $( \heartsuit)$  \\ \hline
 %KV  &  &  & n/a  \\
 BN  &  $(\cC, \otimes, I, id, id,id, id, id, id, id )$  & $(R, \tilde R^{-1}_{(-|-,-)},  \tilde R^{-1}_{(-, -|-)}     )$  & n/a  \\
 DS  & $(\cM, \otimes, I, id, id,id, id, id, id, id )$ & $(\rho, \bar \omega^{-1}_{- - I -}, \bar \omega^{-1}_{- I - -} )$  &   $v$ \\
 M & $(\cK, \otimes, I, {\sf a}, {\sf l}, {\sf r}, \pi, \nu, \lambda, \rho)$ & $(\rho, R, S)$ & $v$ \\
 GPS  & $(T, I, \otimes, \mathfrak{a}, \mathfrak{l}, \mathfrak{r}, \pi, \mu, \lambda, \rho )$ & n/a  & n/a  \\ \hline
 Here &$(\sM, \otimes, 1, \alpha, \ell, r, \pi, \mu, \lambda, \rho)$ & $(\beta, R, S)$ & $\sigma$  \\ \hline
\end{tabular}
\end{center}
\caption{Translating Definitions of Monoidal Bicategories}
\label{transdefofmonoidalbicats}
\end{table}%

\begin{table}[ht]
\begin{center}
\begin{tabular}{|c|c| c | } \hline
 Author & Monoidal $(\clubsuit)$ & Braided $(\diamondsuit)$   \\ \hline
 %KV  &  &    \\
 BN  &   -- & --    \\
 DS  & $(T,\chi, \iota, \omega, \kappa, \xi)$  &  $u$    \\
 M & $(T,\chi, \iota, \omega, \gamma, \delta)$ & $u$   \\
 GPS  & $(H, \chi, \iota, \omega, \gamma, \delta)$ & n/a   \\ \hline
 Here & $(H, \chi, \iota, \omega, \gamma, \delta)$ & $u$ \\ \hline
\end{tabular}
\end{center}
\caption{Translating Definitions of Homomorphisms of Monoidal Bicategories}
\label{transdefofhomsofmonoidalbicats}
\end{table}%

\begin{table}[ht]
\begin{center}
\begin{tabular}{|c| c  c c c | c|} \hline
 Authors  & BN & DS & M & GPS & Here   \\ \hline
 Monoidal Trans. $(\clubsuit)$    & -- & $(\theta, \theta_2, \theta_0)$ & $(\theta, \Pi, M)$ &  $(\theta, \Pi, M)$ &  $(\theta, \Pi, M)$  \\ \hline
Monoidal Mod.  $(\clubsuit)$ &   -- & $s$ &$s$ &    $m$  &    $m$ \\ \hline
\end{tabular}
\end{center}
\caption{Translating Definitions of Transformations and Modifications of Monoidal Bicategories}
\label{transdefofmodsofmonoidalbicats}
\end{table}%

In what follows, we will only need to consider symmetric monoidal bicategories, and not the related notions of monoidal, braided monoidal, and sylleptic monoidal. However these related notions might be of interest to the reader, and so we have grouped the various data and axioms accordingly. We have also given each a designated symbol. Those axioms and data which are designated with the symbol $\clubsuit$ correspond to the definition of monoidal bicategories. Those designated with the symbol $\diamondsuit$ correspond to braided monoidal bicategories, and similarly $\heartsuit$ corresponds to sylleptic and $\spadesuit$ to symmetric monoidal bicategories.

\index{symmetric monoidal bicategory!history|)}

\section{Coherence for Bicategories, in action} \label{sec:coherenceinaction}

The purpose of this section is two-fold. On the one hand, the definition of symmetric monoidal bicategory makes use of certain pasting diagrams which have an apparent ambiguity beyond which can be easily resolved by MacLance's coherence theorem (Theorem~\ref{MacLanesCoherenceTheorem}). In the literature which acknowledges this fact (some appear to overlook it), the issue is quickly brushed aside by appealing to a coherence theorem for functors (see for example \cite[Remark~4.5]{MR3076451}). We would first like to spell this out in more detail as a service to the reader. 

Furthermore, in this chapter we establish a collection of theorems about symmetric monoidal bicategories which together form a toolkit of sorts for working with them. There is an analogous and well established toolkit available for ordinary (non-symmetric monoidal) bicategories, and this gives us an opportunity to demonstrate the potential utility of such theorems, but in the more familiar context of bicategories.   

The definition of symmetric monoidal bicategory, as with many other kinds of higher categorical structure, will consist of a certain collection of data consisting of homomorphisms, natural transformations between various composites of these functors, and modifications between various composites of these transformations. These are then required to satisfy certain relations, which take the form of equations between composite of the modifications. Pasting diagrams are supposed to encode these composites. Consider the following situation.  

\begin{figure}[ht]
\begin{center}
\begin{tikzpicture}[yscale=0.8]
		\node (LT) at (0, 4) {$F_1 G_1 H_1$};
		\node (LB) at (0, 0) {$F_1 G_1 H_2$};
		\node (RT) at (5, 4) {$F_1 G_2 H_2$};
		\node (RB) at (5, 0) {$F_2 G_1 H_2$};
		\draw [->] (LT) -- node [left] {$F_1 \alpha_1$} (LB);
		\draw [->] (LT) -- node [above] {$F_1 G_1 \theta_1$} (RT);
		\draw [<-] (RT) -- node [right] {$\theta_3 G_1 H_2$} (RB);
		\draw [->] (LB) -- node [below] {$\alpha_2 H_2$} (RB);
		\draw [->] (LB) -- node [above left] {$F_1 \theta_2 H_2$} (RT);
		\node at (1.5, 3) {$\Downarrow F_1 \lambda_1$};
		\node at (3.5, 1) {$\Downarrow \lambda_2 H_2$};
		%\node at (0.5, 1) {$\ulcorner$};
		%\node at (1.5, 0.5) {$\lrcorner$};
\end{tikzpicture}
\end{center}
	\caption{A pasting diagram. Does it make sense?}
	\label{fig:pastingdiagramambiguity}
\end{figure}
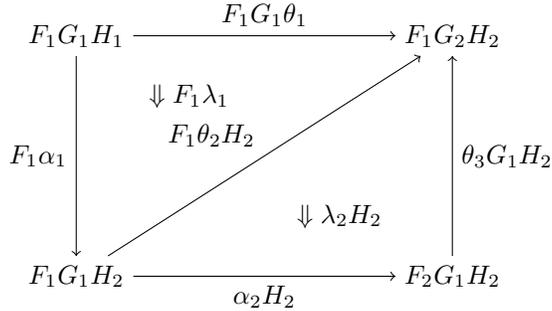

Suppose that we have four bicategories $\sB_1, \sB_2, \sB_3$, and $\sB_4$ (these could, for example, be various products of a single bicategory) and suppose that for each $i=1,2$ we have homomorphisms $H_i: \sB_1 \to \sB_2$, $G_i: \sB_2 \to \sB_3$, and $F_i: \sB_3 \to \sB_4$. Further suppose that we have transformations:
\begin{align*}
	\theta_1:  H_1 &\to H_2 \\
	\theta_2: G_2 &\to G_1 \\
	\theta_3: F_2 &\to F_1 \\
	\alpha_1: G_1 &\circ H_1 \to G_2 \circ H_2 \\
	\alpha_2: F_1 &\circ G_2 \to F_1 \circ G_1
\end{align*}
and that we have modifications:
\begin{align*}
	\lambda_1: G_1 \theta_1 &\Rightarrow (\theta_2 H_2) \circ \alpha_1 \\
	\lambda_2: F_1 \theta_2 &\Rightarrow (\theta_3 G_1) \circ \alpha_2.
\end{align*}
We would like to be able to make sense of the pasting diagram in Figure~\ref{fig:pastingdiagramambiguity}.

This is a pasting diagram in the functor bicategory $\bicat(\sB_1, \sB_4)$.  The familiar calculus of pasting diagrams, which relies on MacLance's coherence theorem, permits us to add parentheses and to insert associativity and unit constraints into the pasting diagram in order to interpret it. Assuming that the pasting diagram is well-typed, which means that by adding these parentheses and extra morphisms we can make the sources and targets match up, we are assured that there is a unique value for the pasting diagram, that it does not in fact depend on the choices of bracketing or which associativity constraints we use. 

However the pasting diagram in Figure~\ref{fig:pastingdiagramambiguity} is not well-typed in this way. The pasting diagram is supposed to correspond to a composite between $(F_1 \lambda_1)$ and $(\lambda_2 H_2) * Id_{(F_1 \alpha_1)}$. The target of the first morphism is $F_1 ((\theta_2 H_2) \circ \alpha_1)$, while the source of the second is $(F_1 \theta_2 H_2) \circ (F_1 \alpha_1)$. These are generally not equal unless the functor $F_1$ is strict (see Remark~\ref{rmk:strictnessofwhisker} in the appendix), and they cannot be made equal by simply using the associators and unitors of the ambient bicategory. The identification requires using the coherence constraints of the functor $F_1$. 

In this example with such a very small pasting diagram there is essentially one clear way to use the coherence constraints of the functors to make the diagram well-typed. However for larger and more complex diagrams this is less clear. Usually we will know that there is at least one way to use the functor coherence constraints to make the pasting diagram well-typed, but how can we be assured that the result doesn't somehow depend on our choices? Before addressing that we will make brief digression. 

Let $\bicat_s$ denote the category of bicategories and strict homomorphisms. Lack has shown that this category supports the structure of a {\em Quillen model category} \cite{MR2138540}, a powerful structure for doing abstract homotopy theory, in this case a homotopy theory of bicategories. Part of this structure defines the notion of {\em cofibrant} bicategories, and produces a comonad $(-)^c: \bicat_s \to \bicat_s$, a functorial {\em cofibrant replacement}. For each bicategory $\sB$, $\sB^c$ is a cofibrant bicategory and the counit $\varepsilon: \sB^c \to \sB$ is an equivalence of bicategories. In fact this comonad is induced from an adjunction 
\begin{equation*}
	\sF_3:\comp_3^B \rightleftarrows \bicat_s: V_3
\end{equation*}
between bicategories and the category of `3-computads' with respect to the bicategory monad $B$ on 2-globular sets; this is a special case of the theory described in Section~\ref{sec:computads}. There is a canonical inverse equivalence $\sB \to (\sB)^c$, but we won't describe that here. We will need the following analog of our cofibrancy theorem:

\begin{lemma}[Cofibrancy Theorem for Bicategories {\cite[Lemma 9]{MR2138540}}]%\label{lem:}
	If $\sA$ is a cofibrant bicategory, and $f: \sA \to \sB$ is a homomorphism of bicategories, then there exists a strict homomorphism $f': \sA \to \sB$ which is equivalent to $f$. \qed
\end{lemma}

Moreover a careful reading of the proof of \cite[Lemma 9]{MR2138540} shows that the construction of $f'$ and the equivalence $f \simeq f'$ is natural in $f$. Transformations between strict functors can also be replaced by equivalent strict transformations, but we will not need that fact either. 

The next property that we will use concerns the adjunction:
\begin{equation*}
	L:\bicat_s \leftrightarrows 2\cat: i.
\end{equation*}
Here $L$ is the left adjoint to the inclusion of the category of strict 2-categories and strict homomorphisms into all bicategories and strict homomorphisms. The functor $L$ turns a bicategory into a strict 2-category by brutally forcing the coherence cells to become identities. In general we would expect $iL(\sB)$ to be a very different bicategory from $\sB$. However for cofibrant bicategories this is not the case.

\begin{lemma}[{\cite[Lemma 10]{MR2138540}}] %\label{lem:}
	Let $\sB$ be a cofibrant bicategory. Then the unit $\sB \to iL(\sB)$ is an equivalence of bicategories. \qed
\end{lemma}

In other words the brutal act of declaring all coherence morphisms to be identities is perfectly cromulent, provided the bicategory is cofibrant. Moreover for each bicategory $\sB$ we have a zig-zag of equivalences (each of which is a strict homomorphism):
\begin{equation*}
	\sB \leftarrow (\sB)^c \to iL(\sB^c) = \st(\sB).
\end{equation*}
Here $\st(\sB)$ is the classical strictification of $\sB$ (which is described in explicit concrete terms in Appendix~\ref{sec:app2-coherencestrict}). Moreover if $f: \sA \to \sB$ is a homomorphism we get a diagram:
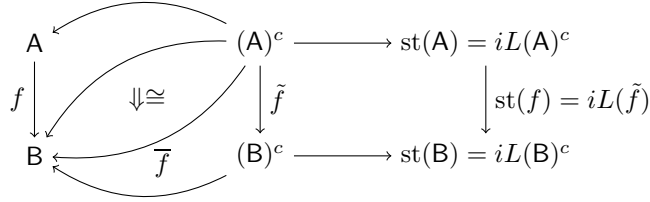
\begin{center}
\begin{tikzpicture}
		\node (LT) at (0, 1.5) {$\sA$};
		\node (LB) at (0, 0) {$\sB$};
		\node (MT) at (3, 1.5) {$(\sA)^c$};
		\node (MB) at (3, 0) {$(\sB)^c$};
		\node (RT) at (6, 1.5) {$\st(\sA) = iL(\sA)^c$};
		\node (RB) at (6, 0) {$\st(\sB) = iL(\sB)^c$};
		\draw [->] (LT) -- node [left] {$f$} (LB);
		\draw [<-] (LT) to [bend left] node [above] {$$} (MT);
		\draw [->] (MT) -- node [right] {$\tilde{f}$} (MB);
		\draw [<-] (LB) to [bend right] node [below] {$$} (MB);
		\draw [->] (MT) to [bend left] node [below] {$\overline{f}$} (LB);
		\draw [->] (MT) to [bend right] node [above] {$$} (LB);
		\draw [->] (MT) -- (RT);
		\draw [->] (MB) -- (RB);
		\draw [->] (RT) -- node [right] {$\st(f) = iL(\tilde{f})$} (RB);
		\node at (1.5, 0.75) {$\Downarrow \cong$};
		%\node at (0.5, 1) {$\ulcorner$};
		%\node at (1.5, 0.5) {$\lrcorner$};
\end{tikzpicture}
\end{center}
Each of the triangles and the square commute strictly, the indicated 2-morphism is an equivalence, and all the homomorphisms are strict except $f$ and the composite $(A)^c \to \sB$. The strict homomorpism $\tilde{f}: \sF_3 V_3(\sA) = (\sA)^c \to (\sB)^c = \sF_3 V_3(\sB)$ is obtained by taking the mate of $\overline{f}:\sF_3 V_3(\sA) \to \sB$, which is a map $V_3(\sA) \to V_3(\sB)$, and then applying $\sF_3$. 
 
In short, not only can we strictify bicategories, but, up to equivalence, we can functorially strictify functors as well. This is (a version of) the coherence theorem for functors. A more detailed analysis shows that this construction is strictly compatible with composition, $\st(f \circ g) = \st(f) \circ \st(g)$ (see \cite{JS93} \cite{MR3076451}). 

Now let us return to the problem of interpreting the pasting diagram of Figure~\ref{fig:pastingdiagramambiguity}.  
After applying the strictification functor, each component of the pasting diagram corresponds to a transformation between strict homomorphisms. Because the bicategories and homomorphisms are now strict, all the relevant associativity constraints are identities and the process of whiskering strictly commutes with vertical composition. The resulting diagram is now well-typed and so gives a well-defined 2-cell of the functor category $\bicat(\st(\sB_1), \st(\sB_2))$. A choice of bracketing of the boundary of the original pasting diagram corresponds to a commutative diagram:
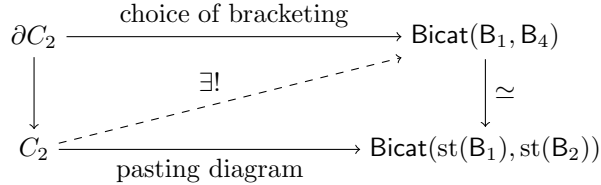
\begin{center}
\begin{tikzpicture}
	\node (LT) at (0, 1.5) {$\partial C_2$};
	\node (LB) at (0, 0) {$C_2$};
	\node (RT) at (6, 1.5) {$\bicat(\sB_1, \sB_4)$};
	\node (RB) at (6, 0) {$\bicat(\st(\sB_1), \st(\sB_2))$};
	\draw [->] (LT) -- node [left] {$$} (LB);
	\draw [->] (LT) -- node [above] {choice of bracketing} (RT);
	\draw [->] (RT) -- node [right] {$\simeq$} (RB);
	\draw [->] (LB) -- node [below] {pasting diagram} (RB);
	\draw [dashed, ->] (LB) -- node [above left] {$\exists !$} (RT);
	%\node at (0.5, 1) {$\ulcorner$};
	%\node at (1.5, 0.5) {$\lrcorner$};
\end{tikzpicture}
\end{center}
Since the right-hand homomorphism is an equivalence, there exists a unique lift as indicated by the dashed arrow. Thus we do in fact have a canonical and unambiguous way to interpret each such pasting diagram, once parentheses have been added to the boundary.

\section{Symmetric Monoidal Bicategories}	\label{SectSmyMonBicat}

\begin{definition}\label{DefnSymMonBicat}
A {\em symmetric monoidal bicategory} 
\index{bicategory!symmetric monoidal}
\index{bicategory!symmetric monoidal|see{symmetric monoidal bicategory}}
\index{symmetric monoidal bicategory}
consists of a bicategory $\sM$ together with the following data: 
\begin{enumerate}
\item [$\clubsuit$] a distinguished object $1 \in \sM$,
\item [$\clubsuit$] a homomorphism
	\begin{equation*}
		\otimes = (\otimes, \phi_{(f,f'), (g, g')}^\otimes, \phi^\otimes_{(a,a')} ): \sM \times \sM \to \sM
\end{equation*}
\item [] transformations:
\begin{align*} 
&\clubsuit
\begin{cases}
\alpha = ( \alpha_{abc}, \alpha_{fgh})  : (a \otimes b) \otimes c \to a \otimes (b \otimes c) \\
	\ell = (\ell_a, \ell_f)  : 1 \otimes a \to a \\
	r = (r_a, r_f) : a \to a \otimes 1 
\end{cases} \\
& \diamondsuit \begin{cases}
\beta = (\beta_{ab}, \beta_{fg}) : a \otimes b \to b \otimes a
\end{cases}
\end{align*}
which are equivalence transformations. We also choose adjoint inverses $\alpha^*$, $\ell^*$, and  $r^*$ and their associated adjunction data, which we will not name. 
%\item [$\diamondsuit$] an equivalence transformation:
%\begin{equation*}
%	\beta = (\beta_{ab}, \beta_{fg}) : a \otimes b \to b \otimes a
%\end{equation*}
\item [$\clubsuit$] invertible modifications $\pi$, $\mu$, $\lambda$, and $\rho$ as in Figure~\ref{fig:DefnSymMonBicat1}.

\begin{figure}[ht]
	  \begin{center}

	\end{center}
	
	\caption{Modifications for Monoidal Bicategories ($\clubsuit$)}
	\label{fig:DefnSymMonBicat1}
\end{figure}
\item [$\diamondsuit$] invertible modifications $R$ and $S$, as in Figure~\ref{fig:DefnSymMonBicat2}.
\begin{figure}[ht]
	
	\begin{center}
	%
	\end{center}
	\caption{Modifications for Braided Monoidal Bicategories ($\diamondsuit$)}
	\label{fig:DefnSymMonBicat2}
\end{figure}

\item [$\heartsuit$] an invertible modification $\sigma$, as in Figure~\ref{fig:DefnSymMonBicat3}.
\begin{figure}[ht]
	\begin{center}
	%
	\end{center}
	\caption{Modifications for Sylleptic Monoidal Bicategories ($\heartsuit$)}
	\label{fig:DefnSymMonBicat3}
\end{figure}
\end{enumerate}
such that the following axioms are satisfied:
\begin{enumerate}
\item [$\clubsuit$] The equations (SM1), (SM2.i), and (SM2.ii) from Appendix \ref{app:defnsymbicat} (depicted in Figures \ref{fig:SM1a}) through \ref{fig:SM2ii}) are satisfied, using the data
$(\sM, 1, \otimes, \alpha, \ell, r, \pi, \mu, \lambda, \rho)$.
\item [$\diamondsuit$] The equations (SM3.i), (SM3.ii), (SM4), and (SM5) from Appendix \ref{app:defnsymbicat} (depicted in Figures \ref{fig:SM3ia}) through \ref{fig:SM5b}) are satisfied, using the above data and also $(\beta, R, S)$. 
\item [$\heartsuit$] The equations (SM6.i) and (SM6.ii) from Appendix \ref{app:defnsymbicat} (depicted in Figures \ref{fig:SM6i}) and \ref{fig:SM6ii}) are satisfied, using the above data and $\sigma$.
\item [$\spadesuit$] Furthermore, equation (SM7) from Appendix \ref{app:defnsymbicat} (depicted in Figure \ref{fig:SM7}) is satisfied.
\end{enumerate}
A bicategory equipped with the data $\clubsuit, \diamondsuit, \heartsuit$, but only satisfying axioms $\clubsuit, \diamondsuit$ and $\heartsuit$ is a {\em sylleptic monoidal bicategory}. A bicategory equipped with the data $\clubsuit$ and $\diamondsuit$ and satisfying the axioms $\clubsuit$ and $\diamondsuit$ is a {\em braided monoidal bicategory}, and lastly a bicategory equipped only with data $\clubsuit$ satisfying the axioms $\clubsuit$ is a {\em monoidal bicategory}.
\index{bicategory!monoidal}
\index{bicategory!braided monoidal}
\index{bicategory!sylleptic monoidal}
\index{sylleptic monoidal bicategory}
\index{braided monoidal bicategory}
\index{monoidal bicategory}
\end{definition}

\begin{remark}
	The purpose of specifying the inverses $\alpha^*$, $\ell^*$, and  $r^*$ and their associated adjunction data, is to make the notion of symmetric monoidal bicategory {\em algebraic}, meaning that the forgetful functor from symmetric monoidal bicategories to 2-dimensional globular sets is the right adjoint of a monadic adjunction.  We will refer to the corresponding monad on 2-dimensional globular sets as the {\em symmetric monoidal bicategory monad}. 
	
	The original definition of tricategory (hence also monoidal bicategory) from Gordon-Powers-Street \cite{GPS95} merely required the morphisms $\alpha, \ell$ and $r$ to be equivalences, without specifying there inverses. This yields a non-algebraic notion, and \cite{MR2717302} modified this notion to obtain an equivalent algebraic definition of tricategory. By making the inverse equivalences into adjoint equivalences the notion becomes algebraic and remains essentially the same. 

Note however that in a {\em symmetric} monoidal bicategory we do not need to specify an additional inverse adjoint equivalence for the transformation $\beta$. The morphism $\beta$ is its own inverse adjoint equivalence with the syllepsis modification $\sigma$ (and its inverse) providing the unit and counit of the adjunction equivalence. The zig-zag identities required of an adjunction are equivalent to axiom (SM7) (See Figure~\ref{fig:SM7}).
\end{remark}

\begin{definition} \label{def:symmonhom}
Let $\sM$ and $\sM'$ be two symmetric monoidal bicategories. 
 A {\em symmetric monoidal homomorphism} 
\index{homomorphism}
\index{bicategory!symmetric monoidal homomorphism}
\index{symmetric monoidal bicateogry!homomorphism}
$H: \sM \to \sM'$ consists of the following data:
 \begin{enumerate}
\item [$\clubsuit$]  a homomorphism $H: \sM \to \sM'$
\item  [$\clubsuit$] transformations:
\begin{align*}
&\chi = (\chi_{ab}, \chi_{fg}): H(a) \otimes H(b) \to H(a \otimes b) \\
&\iota: 1' \to H(1)
\end{align*}
\item [] adjoint equivalence transformations $\chi^*$ and $\iota^*$, together with their associated adjunction data, which we leave unnamed. 
\item  [$\clubsuit$] invertible modifications $\omega$, $\gamma$, and $\delta$, as in Figure~\ref{fig:MonoidalHom1}. 
\begin{figure}[ht]
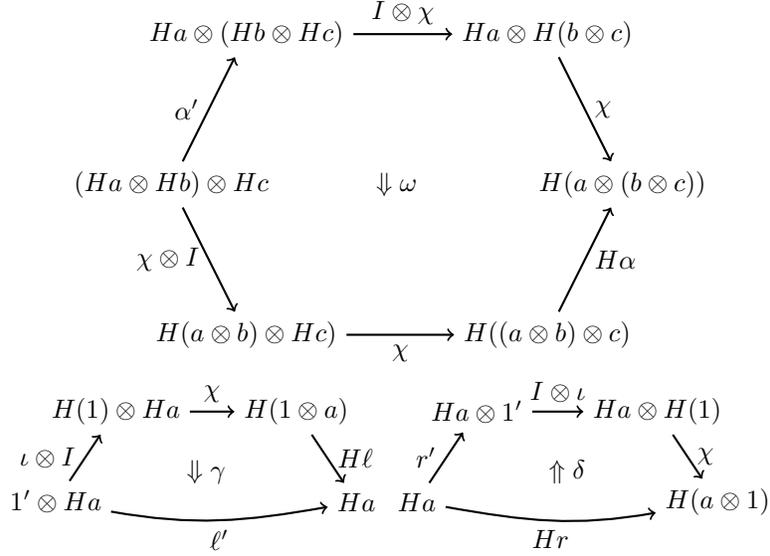

	\begin{center}

	\end{center}
	\caption{Modifications for Monoidal Homomorphism ($\clubsuit$)}
	\label{fig:MonoidalHom1}
\end{figure}

\item [$\diamondsuit$] an invertible modification $u$, as in Figure~\ref{fig:MonoidalHom2}.
\begin{figure}[ht]
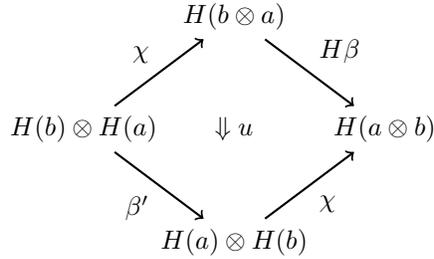

	\begin{center}
	%
	\end{center}
	\caption{Modifications for Braided Monoidal Homomorphism ($\diamondsuit$)}
	\label{fig:MonoidalHom2}
\end{figure}
\end{enumerate}
such that the following axioms hold:
\begin{enumerate}
\item [$\clubsuit$] Equations (HTA1) and (HTA2) of \cite{GPS95} hold using the data $(H, \chi, \iota, \omega, \gamma, \delta)$. 
\item [$\diamondsuit$] Equations (BHA1) and (BHA2) of \cite{McCrudden00} hold using the previous data and $u$. 
\item [$\heartsuit, \spadesuit$] Equation (SHA1) of \cite{McCrudden00} holds.  
\end{enumerate}
A symmetric monoidal homomorphism will be called {\em strict} if $H$ is a strict homomorphism of bicategories, the coherence transformations $\chi$ and $\iota$ are identities, the modifications $\omega$, $\gamma$, $\delta$ are the canonical modifications given by the coherence isomorphisms of the functor $\otimes$ or the hom bicategory (see \cite[Def.3.3.3]{MR2717302}), and the modification $u$ is an identity.  

If $\sM$ and $\sM$ are sylleptic monoidal bicategories, then the identical data and axioms define a {\em sylleptic monoidal homomorphism}. If $\sM$ and $\sM$ are braided monoidal bicategories, then the data $\clubsuit$ and $\diamondsuit$ subject to axioms $\clubsuit$ and $\diamondsuit$ define a {\em braided monoidal homomorphism}. If $\sM$ and $\sM'$ are merely monoidal bicategories, then the data and axioms $\clubsuit$ define a {\em monoidal homomorphism}. 
\end{definition}

\begin{remark}
	Strict homomorphisms between symmetric monoidal bicategories may be identified with the morphisms of $T$-algebras where $T$ is the symmetric monoidal bicategory monad on 2-dimensional globular sets. In this way symmetric monoidal bicategories and strict homomorphism form a 1-category. The reader should be warned however that defined this way the composition of strict homomorphisms does not precisely agree with the composition of general symmetric monoidal homomorphisms as defined below in Section~\ref{SectSymMonWhisker}. With the composition defined there, the composite of strict homomorphisms is not necessarily strict. A completely analogous problem exists for strict functors of tricategories, see \cite[Chap.~4.1]{MR2717302}.
\end{remark}

\begin{definition}
Let $\sM$ and $\sM'$ be symmetric monoidal bicategories and let $H, \bar H: \sM \to \sM'$ be two symmetric monoidal homomorphisms. Then a {\em symmetric monoidal transformation} 
\index{transformation}
\index{bicategory!symmetric monoidal transformation}
\index{symmetric monoidal bicateogry!transformation}
$\theta: H \to \bar H$ consists of:
\begin{enumerate}
\item [$\clubsuit$] A transformation: $\theta = (\theta_a, \theta_f): H \to \bar H$
\item [$\clubsuit$] Modifications $\Pi$ and $M$, as in Figure~\ref{fig:MonoidalTrans}. 

\begin{figure}[ht]
	 \begin{center}
	\begin{tikzpicture}[thick]
		%\node at (0, 4.5) {\underline{The Pentagonator}};
		\node (LT) at (0,2) 	{$Ha \otimes Hb$ };
		\node (LB) at (1,0) 	{$\bar Ha \otimes Hb$};
		\node (RT) at (8,2) 	{$\bar H ( a \otimes b)$};
		\node (MT) at (4,3.5) {$H(a \otimes b)$};
		\node (RB) at (7,0)	{$ \bar Ha \otimes \bar H b$};
		\draw [->] (LT) --  node [left] {$ \theta \otimes I$} (LB);
		\draw [->] (LT) -- node [above left] {$\chi $} (MT);
		\draw [->] (MT) -- node [above right] {$\theta $} (RT);
		\draw [->] (RB) -- node [right] {$\bar \chi $} (RT);
		\draw [->, in = 200, out = 340] (LB) to node [below] {$I \otimes \theta $} (RB);
		\node at (4,1 ) {$\Uparrow \Pi$};
	\end{tikzpicture}
	\end{center}
	\begin{center}
	\begin{tikzpicture}[thick]
		\node (A) at (0, 1.5) {$1'$};
		\node (B) at (3, 1.5) {$\bar H (1)$};
		\node (C) at (1.5, 0) {$H(1)$};
		\node at (1.5, .75) {$\Uparrow M$};

		\draw [->] (A) -- node [above] {$ \bar \iota$} (B);
		\draw [->, out = -90, in = 180] (A) to node [below left] {$\iota$} (C);
		\draw [->, out = 0, in = -90] (C) to node [below right] {$\theta$} (B);
	\end{tikzpicture}
	\end{center}
	\caption{Modifications for Monoidal Natural Transformation}
	\label{fig:MonoidalTrans}
\end{figure}
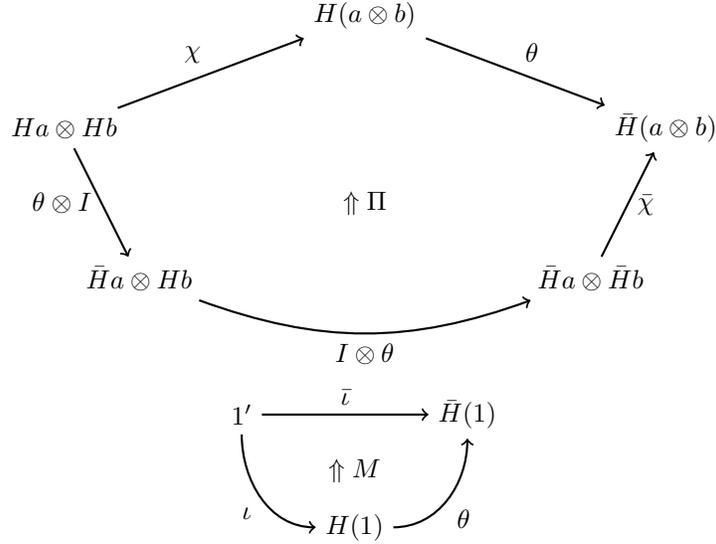
\end{enumerate}
such that the following axioms are satisfied:
\begin{enumerate}
\item [$\clubsuit$] The following equations hold, where these designations refer to the 2-morphisms in Figures \ref{fig:MBTA1.a}, \ref{fig:MBTA1.b}, \ref{fig:MBTA2.a}, \ref{fig:MBTA2.b}, \ref{fig:MBTA3.a}, and \ref{fig:MBTA3.b}.
\begin{align*}
(MBTA1.a)& = (MBTA1.b), \\
(MBTA2.a) &= (MBTA2.b), \\
(MBTA3.a)& = (MBTA3.b).
\end{align*}

\item [$\diamondsuit, \heartsuit, \spadesuit$] Equation (BTA1) of \cite{McCrudden00} holds. 
\end{enumerate}
%A symmetric monoidal transformation is {\em strict} if $\theta$ is a strict transformation and the coherence modifications $\Pi$ and $M$ are identities. 

{\em Braided monoidal transformations} and {\em sylleptic monoidal transformations} are defined by the same data and axioms. If $\sM$ and $\sM'$ are merely monoidal bicategories, and $H$ and $\bar H$ are monoidal homomorphisms, then the data and axioms $\clubsuit$ define a {\em monoidal transformation}.
\end{definition}

\begin{figure}[ht]
%	\begin{center}

	\caption{Axiom MBTA1.a}
	\label{fig:MBTA1.a}
\end{figure}

\begin{figure}[ht]
%
	\caption{Axiom MBTA1.b}
	\label{fig:MBTA1.b}
\end{figure}

\begin{figure}[ht]
%
	\caption{Axiom MBTA2.a}
	\label{fig:MBTA2.a}
\end{figure}

\begin{figure}[ht]
%
	\caption{Axiom MBTA2.b}
	\label{fig:MBTA2.b}
\end{figure}

\begin{figure}[ht]
%
	\caption{Axiom MBTA3.a}
	\label{fig:MBTA3.a}
\end{figure}

\begin{figure}[ht]
%
	\caption{Axiom MBTA3.b}
	\label{fig:MBTA3.b}
\end{figure}
%\end{center}

\begin{definition} \label{DefnSymMonoidalModifiaction}
 ($\clubsuit, \diamondsuit, \heartsuit$ and $\spadesuit$) Let $\sM$ and $\sM'$ be symmetric monoidal bicategories, $H, \bar H: \sM \to \sM'$ be symmetric monoidal homomorphisms and $\theta, \tilde \theta: H \to \bar H$ be symmetric monoidal transformations. Then a {\em symmetric monoidal modification} 
\index{modification}
\index{bicategory!symmetric monoidal modification}
\index{symmetric monoidal bicategory!modification}
$m: \theta \to \tilde \theta$ consists of a modification $m: \theta \to \tilde \theta$ such that the  two equations in Figures \ref{fig:BMBM1} and \ref{fig:BMBM2} hold:
\begin{figure}[ht]

	\caption{Axiom BMBM1}
	\label{fig:BMBM1}	
\end{figure}

\begin{figure}[ht]
%
	\caption{Axiom BMBM2}
	\label{fig:BMBM2}	
\end{figure}
\end{definition}

\begin{definition}
Let $\sM$ and $\sM'$ be two symmetric monoidal bicategories, let $H, \overline{H}: \sM \to \sM'$ be symmetric monoidal homomorphisms, let $\theta, \tilde \theta, \tilde{\tilde{\theta}}: H \to \overline{H}$ be symmetric monoidal transformations, and let $m: \theta \to \tilde \theta$ and $\tilde m: \tilde \theta \to \tilde{\tilde{ \theta}}$ be symmetric monoidal modifications. Then the modification $\tilde m \circ m: \theta \to \tilde{\tilde{\theta}}$ is a symmetric monoidal modification called the {\em (vertical) composition} of the symmetric monoidal modifications $m, \tilde m$. 
\end{definition}

\begin{definition}
Let $\sM$ and $\sM'$ be two symmetric monoidal bicategories, let $H, \overline{H}, \overline{\overline H}: \sM \to \sM'$ be symmetric monoidal homomorphisms, and let $\theta: H \to \overline H$, $\tilde \theta: \overline H \to \overline{\overline H}$ be symmetric monoidal transformations. Then the {\em composition} of $\theta$ and $\tilde \theta$ is the symmetric monoidal transformation given by the transformation $\tilde \theta \circ \theta$ and the modifications whose components are given by the pasting diagrams in Figures \ref{fig:PastingDiagComp1} and \ref{fig:PastingDiagComp2}.

\begin{figure}[ht]
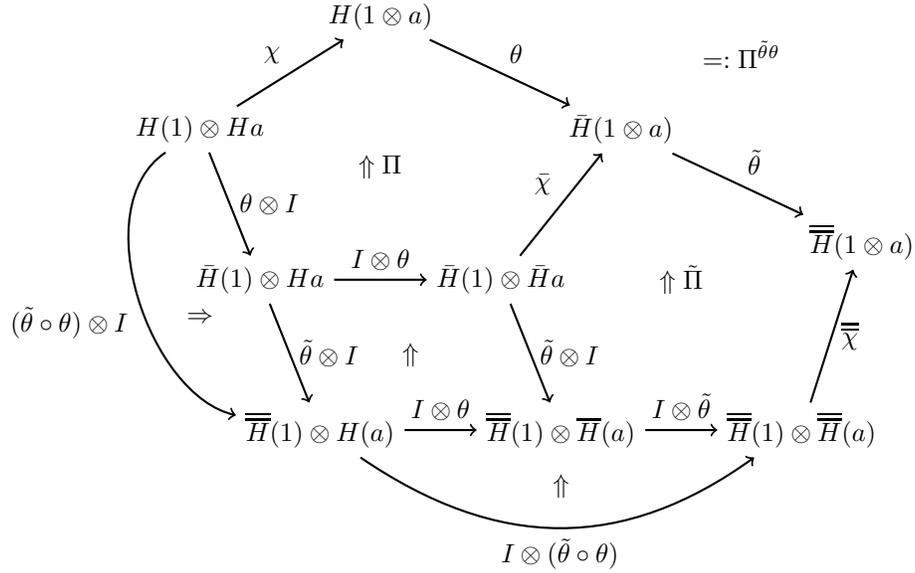


	\caption{Pasting Diagram for $\Pi^{\tilde \theta \theta}$ }
	\label{fig:PastingDiagComp1}
\end{figure}

\begin{figure}[ht]
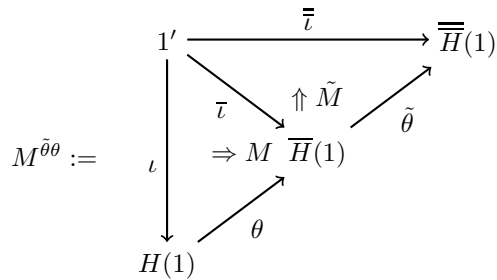

%
	\caption{Pasting Diagram for $M^{\tilde \theta \theta}$ }
	\label{fig:PastingDiagComp2}
\end{figure}
\end{definition}

\begin{definition}
Let $\sM$ and $\sM'$ be two symmetric monoidal bicategories, let $H, \overline{H}, \overline{\overline H}: \sM \to \sM'$ be symmetric monoidal homomorphisms, let $\theta, \xi: H \to \overline H$, $\tilde \theta, \tilde \xi: \overline H \to \overline{\overline H}$ be symmetric monoidal transformations and let $m: \theta \to \xi$ and $\tilde m: \tilde \theta \to \tilde \xi$ be symmetric monoidal modifications. Then the {\em horizontal composition} of $m$ and $\tilde m$ is the symmetric monoidal modification $\tilde m * m$. 
\end{definition}

Together with these compositions and the obvious associators and unitors, we have a bicategory $\symbicat( \sM, \sM')$ whose objects are the symmetric monoidal homomorphisms, the 1-morphisms are the symmetric monoidal transformations, and whose 2-morphisms are the symmetric monoidal modifications.

\section{Symmetric Monoidal Whiskering} \label{SectSymMonWhisker}

Now we will begin to develop the theory of symmetric monoidal bicategories. 

\begin{definition}
Let $\sM, \sM'$ and $\sM''$ be symmetric monoidal bicategories. Let $H: \sM \to \sM'$ and $\overline H : \sM' \to \sM''$ be symmetric monoidal homomorphisms. Then the {\em composition} of $H$ and $\overline H$ is the symmetric monoidal homomorphism given by the homomorphism $\overline H H : \sM \to \sM''$, the following transformations:
\begin{align*}
\chi^{\overline H H} &: \overline HH (a) \otimes'' \overline HH(b) \stackrel{\overline \chi}{\to} \overline H (H a \otimes' Hb) \stackrel{\overline H \chi}{\to} \overline HH( a \otimes b) \\
\iota^{\overline HH} & : 1'' \stackrel{\overline \iota}{\to} \overline H ( 1') \stackrel{\overline H \iota}{\to} \overline H H (1)
\end{align*}
and the modifications whose components are given by the pasting diagrams in Figures \ref{fig:PastingDiagHomComp1}, \ref{fig:PastingDiagHomComp3}, \ref{fig:PastingDiagHomComp2}, and \ref{fig:PastingDiagHomComp4}. Here the unlabeled 2-morphisms are the components of the transformations $\chi, \iota, \overline \chi$ and $\overline \iota$ and the canonical coherence 2-morphisms from $\sM''$. 
%\begin{landscape}
%	\begin{figure}[p]
\begin{sidewaysfigure}
		\vspace{5in}
		{%\scriptsize
 }
	\caption{Pasting Diagram for $\omega^{\overline HH}$ }
	\label{fig:PastingDiagHomComp1}	
\end{sidewaysfigure}
%	\end{figure}
%\end{landscape}

\begin{figure}[ht]
%
	
\caption{Pasting Diagram for $\gamma^{\overline H H}$ }
\label{fig:PastingDiagHomComp2}
\end{figure}

\begin{figure}[ht]
%
\caption{Pasting Diagram for $\delta^{\overline H H}$ }
\label{fig:PastingDiagHomComp3}
\end{figure}

\begin{figure}[ht]
%
\caption{Pasting Diagram for $u^{\overline H H}$ }
\label{fig:PastingDiagHomComp4}
\end{figure}
\end{definition}

\begin{definition}
Let $\sA, \sB, \sC, \sD$ be symmetric monoidal bicategories, $F: \sA \to \sB$, $G, \overline G: \sB \to \sC$, and $H: \sC \to \sD$ be symmetric monoidal homomorphisms, and $\theta: G \to \overline G$ be a symmetric monoidal transformation. Then the {\em whiskering} 
\index{whiskering}
of $F$ and $\theta$ is the symmetric monoidal homomorphism from $GF$ to $\overline G F$ defined by the transformation $\theta F: GF \to \overline GF$ and the modification whose components are given by the pasting diagrams in Figures \ref{fig:PastingDiagWisk1} and \ref{fig:PastingDiagWisk2}.

\begin{figure}[ht]

	\caption{Pasting Diagram for $\Pi^{\theta F}$}
	\label{fig:PastingDiagWisk1}
\end{figure}

\begin{figure}[ht]
	%
	\caption{Pasting Diagram for $M^{\theta F}$}
	\label{fig:PastingDiagWisk2}
\end{figure}

The {\em whiskering} of $\theta$ and $H$ is the symmetric monoidal transformation from $HG$ to $H \overline G$ given by the transformation $H\theta: HG \to H \overline G$ and the modification whose components are given by the pasting diagrams in Figures \ref{fig:PastingDiagWisk3} and \ref{fig:PastingDiagWisk4}.

\begin{figure}[ht]
	%
	\caption{Pasting Diagram for $M^{H \theta }$}
	\label{fig:PastingDiagWisk3}
\end{figure}

\begin{figure}[ht]
	%
	\caption{Pasting Diagram for $\Pi^{H \theta} $}
	\label{fig:PastingDiagWisk4}
\end{figure}

Here the unlabeled arrows are the structure 2-morphisms from the transformation $\chi^H$ and the homomorphism $H$. 

If $\tilde \theta: G \to \overline G$ is another symmetric monoidal homomorphism, and $m: \theta \to \tilde \theta$ is a symmetric monoidal modification, then the {\em whiskering} of $m$ by $H$ and $F$ is the symmetric monoidal modification consisting of the modification $HmF: H \theta F \to H \tilde \theta F$.  
\end{definition}

\begin{definition}
A symmetric monoidal homomorphism $H: \sM \to \sM'$ is an {\em equivalence}
\index{symmetric monoidal bicategory!equivalence}
\index{bicategory!symmetric monoidal equivalence}
 of symmetric monoidal bicategories if there exists a symmetric monoidal homomorphism $F: \sM' \to \sM$ such that $FH \simeq id_\sM$ in the bicategory $\symbicat(\sM, \sM)$ and $HF \simeq id_{\sM'}$ in the bicategory $\symbicat(\sM', \sM')$. 
\end{definition}

Using whiskering, it is straightforward to construct the necessary compositions which realize symmetric monoidal bicategories as a tricategory. This is done exactly as it has been done for ordinary bicategories, \cite{GPS95}, and for monoidal bicategories \cite{MR2839900}. 
The hom bicategories are precisely the bicategories $\symbicat(\sM, \sM')$ considered previously. Equivalence in any tricategory is defined exactly as above, and we moreover have the following result, an analog of which is true in any tricategory. 

\begin{lemma}
Let $\sM$ and $\sM'$ be symmetric monoidal bicategories, and let $H: \sM \to \sM'$ be a symmetric monoidal homomorphism, which is an equivalence. Then for all symmetric monoidal bicategories $\sB$, the canonical homomorphisms (induced by whiskering) 
\begin{align*}
	H^*: & \symbicat( \sM', \sB) \to  \symbicat( \sM, \sB) \\
	H_*: & \symbicat( \sB, \sM) \to \symbicat( \sB, \sM')
\end{align*}
are equivalences of bicategories. 
\end{lemma}

\section[Whitehead's Thm for Symmetric Monoidal Bicategories]{Whitehead's Theorem for Symmetric Monoidal Bicategories} \label{SectWhiteheadsTheoremforSymmetricMonoidalBicategories}
\index{Whitehead's theorem!for symmetric monoidal bicategories|(}
\index{symmetric monoidal bicategory!Whitehead's theorem|(}
\index{bicategory!Whitehead's theorem for symmetric monoidal|(}

Whitehead's theorem for topological spaces states that a map between sufficiently nice\footnote{Here ``sufficiently nice'' can be taken to mean spaces with the homotopy type of a CW-complex.} topological spaces is a homotopy equivalence if and only if it induces an isomorphism of all homotopy groups (at all base points). A direct analog of this is the statement that a functor between groupoids $F:G \to G'$ is an equivalence of categories if and only if it induces a bijection on the set of isomorphisms classes of objects (i.e., is an isomorphism on $\pi_0$) and for each object $x \in G$, $F$ induces an isomorphism $F: G(x,x) \to G'(Fx,Fx)$, i.e., is an isomorphism on $\pi_1(G, x)$ for each object $x \in G$. 

A similar well known statement holds for functors between categories: a functor $F: C \to C'$ is an equivalence of categories if and only if it induces a surjection of the set of isomorphism classes of objects (i.e., $F$ is essentially surjective) and induces a bijection $C(x,y) \to C'(Fx, Fy)$ for all pairs of objects $x,y \in C$ (i.e., $F$ is fully-faithful). 
In this section we prove an analogous statement for symmetric monoidal bicategories, which allows us to easily verify when a homomorphism between symmetric monoidal bicategories is an equivalence.% This will later be used to prove that the small model of the bordism category that we will construct is equivalent to the geometric bordism bicategory. 

%We begin with a small combinatorial digression on binary trees. Although, as stated, Lemma \ref{AnyPathsofTreesareEquivalent} is evidently true and it's proof relatively straightforward, we will see that it implies a sort of  ``coherence theorem'' which will simplify later proofs and constructions. This should be compared with the construction given in \cite[VII.2]{MacLane71}, which is in fact more complicated then the one presented here. 

\begin{definition}
A {\em binary tree} 
\index{binary tree}
is defined recursively as follows. The symbol $(-)$ is a binary tree. If $u$ and $v$ are binary trees, then $u \boxdot v = (u) \boxdot (v)$ is a binary tree. Similarly the {\em edges} of a binary tree are defined recursively as a set of subtrees. The edges of the binary tree $(-)$ consist of the set $\{ (-)\}$. The edges of a binary tree $u \boxdot v$ consist of the disjoint union of the edges of $u$, the edges of $v$, and the set $\{ u \boxdot v\}$. A {\em marked} binary tree $(t, S)$ consists of a pair where $t$ is a binary tree and $S$ is a subset of the edges of $t$. 
\end{definition}

\begin{figure}[ht]
\begin{center}

\begin{tikzpicture} [thick, scale=0.8]
	\begin{scope}[scale = .75]
		\draw (3.5, 4.5) -- (3.5, 3.5) -- (0,0) (.5, .5) -- (1,0) (1.5, 1.5) -- (3,0) ;
		\draw (3.5, 3.5) -- node [above right] {$e$} (5.5, 1.5) -- (7,0) (5.5, 1.5) -- (4, 0) (4.5, .5) -- (5,0) (6.5, .5) -- (6, 0);
	\end{scope}		
	
	\begin{scope}[xshift = 8cm, scale = .75]
		\draw (3.5, 4.5) -- (3.5, 3.5) -- (0,0) (.5, .5) -- (1,0) (1.5, 1.5) -- (3,0) ;
		\draw [dashed] (3.5, 3.5) --  (7,0) (5.5, 1.5) -- (4, 0) (4.5, .5) -- (5,0) (6.5, .5) -- (6, 0);
	\end{scope}
	
	\node at (7, -2) {$t = (((-) \boxdot (-)) \boxdot  (-) ) \boxdot \underbrace{(( (-) \boxdot (-) ) \boxdot ( (-) \boxdot (-)  ))}_e$};
%	\node at (7, -3) {$e = (( (-) \boxdot (-) ) \boxdot ( (-) \boxdot (-)  ))$};

	\draw [<->] (6, 2) -- (7, 2);
	\draw [<->] (4, -1.2) -- (3.5, -.8);
	\draw [<->] (10, -1.2) -- (10.5, -.8);
	
\end{tikzpicture}

\caption{Equivalent Descriptions of Binary Trees and Edges}
\label{BinTreesandEdges}
\end{center}
\end{figure}
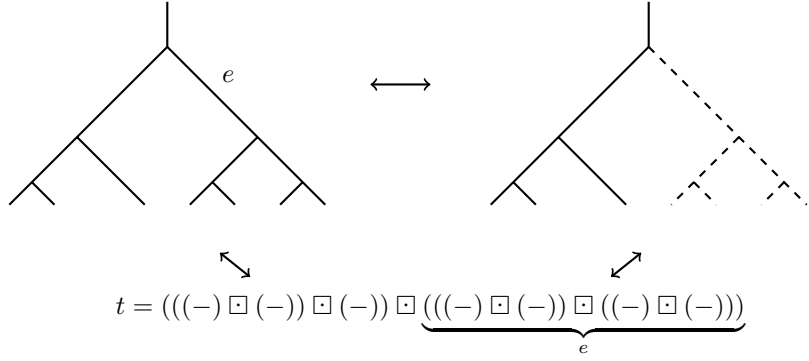

The definition given above is logically equivalent to the usual definition of a proper binary planar rooted tree. Every such tree has a top-most edge (the root) and so every edge $e$ can be equivalently described by specifying the maximal subtree with $e$ as its root, see Figure~\ref{BinTreesandEdges}. The edges of a tree form a partially ordered set with partial order given by inclusion of subtrees. We refer to this partial order as the {\em height} of the edge and thus may speak of edges which have comparable height and edges which have incomparable height.  

\begin{definition}
	Let $(t, S)$ and $(t, S')$ be marked binary trees with the same underlying tree $t$, and such that $S = S' \sqcup S''$ for some collection of edges $S''$. A {\em path} from $(t, S)$ to $(t, S')$ is a word whose letters consist precisely of the edges $S''$ with no omissions or repetitions.
\end{definition}

Paths have an evident composition given by concatenation, and we also consider the trivial/empty word to be a path from the marked tree $(t,S)$ to itself. Any two paths can be related by permuting the letters of the given word. For later reference it is useful to split this process into the following elementary operations:
 
\begin{lemma} \label{AnyPathsofTreesareEquivalent}
	Let $\sim$ be an equivalence relation on paths between marked binary trees such that:
	\begin{enumerate}
	\item [(i)] $\sim$ is closed under composition of paths,
	%For any paths $u: (t, S) \to (t, S')$, $w, \hat w: (t, S') \to (t, S'')$ and $v: (t, S'') \to (t, S''')$, if $w \sim \hat w$, then $vw u \sim v \hat w u$. (Compatible with composition)
	\item [(ii)] Edges of incomparable heights commute: If $e_0, e_1$ are two edges of $t$ which have incomparable heights, then 
	\begin{equation*}
		(e_0 e_1) \sim (e_1 e_0): (t, S \cup \{ e_0, e_1 \}) \to (t,S),
	\end{equation*}
	\item [(iii)]  Edges of comparable heights commute: If $e_0, e_1$ are two edges of $t$ such that $e_0 \leq e_1$, then 
	\begin{equation*}
		(e_0 e_1) \sim (e_1 e_0): (t, S \cup \{ e_0, e_1 \}) \to (t,S).
	\end{equation*}
	\end{enumerate}
	Then any two paths $P$ and $Q$ from $(t, S)$ to $(t, S')$ are equivalent, $P \sim Q$.   \qed
\end{lemma}

We now turn our attention to symmetric monoidal bicategories. The significance of the previous lemma will become clear momentarily. The goal of this section is to prove Theorem~\ref{WhiteheadforSymMonBicats}, which states that a symmetric monoidal homomorphisms between symmetric monoidal bicategories is a symmetric monoidal equivalence if and only if it is an equivalence of the underlying bicategories, which in turn is true if an only if the homomorphism is essentially surjective on objects and induces equivalence of hom categories. While it is possible to prove this theorem in a single direct calculation (by building an inverse symmetric monoidal homomorphism explicitly using the axiom of choice) the magnitude of such a calculation would soon become overwhelming. 

Instead we prefer to break the problem into more calculationally manageable pieces. We first introduce notions of 0- and 1-skeletal symmetric monoidal bicategories and prove that any symmetric monoidal bicategory is equivalent to one of these. This permits us to reduce the general Whitehead's theorem to the case when both bicategories are skeletal, a case which is more tractable. Recall that a category $C$ is {\em skeletal} if $x \cong y$ implies $x = y$ for any objects $x,y \in C$.
\begin{definition}
A bicategory $\sB$ is {\em 1-skeletal}
%\index{bicategory!1-skeletal}
%\index{bicategory!0-skeletal}
\index{bicategory!skeletal}
\index{skeletal bicategory}
 if the hom bicategories $\sB(x, y)$ are skeletal for each pair of objects $x,y \in \sB$. $\sB$ is {\em 0-skeletal} when it satisfies the condition: $x$ is equivalent to $y$ if and only if $x = y$ in $\sB$. Finally, a bicategory $\sB$ is {\em skeletal} if it is both 0-skeletal and 1-skeletal. A symmetric monoidal bicategory is $k$-skeletal if it underlying bicategory is $k$-skeletal.%\footnote{In any good inductive definition of $n$-category one should be able to similarly define the notion of $n$-skeletal. In this case a bicategory is always $k$-skeletal for $k \geq 2$, and a skeletal bicategory is $k$-skeletal for all $k$.} 
\end{definition}

The following lemma is well known, but we reproduce it here since we will need to use the details of its proof in what follows. % It also serves as a means to illustrate how we will employ Lemma \ref{AnyPathsofTreesareEquivalent} in later calculations.

\begin{lemma}
Every bicategory is equivalent to a 1-skeletal bicategory.
\end{lemma}

\begin{proof} Let $\sB$ be a given bicategory.
We will build a 1-skeletal bicategory $\sB^s$ and an equivalence $\sB^s \simeq \sB$. The objects of $\sB^s$ will coincide with those of $\sB$.  Choose for each hom bicategory $\sB(x, y)$  representatives for each isomorphism class of 1-morphism. Let $\sB^s(x,y)$ be the full subcategory with these 1-morphisms as objects. $\sB^s(x,y)$ is skeletal, and the  fully-faithful inclusion functor $i: \sB^s(x,y) \hookrightarrow \sB(x,y)$ is an equivalence of categories. We construct an explicit inverse functor by choosing, for each $f \in \sB(x,y)$ an isomorphism $\eta_f: \bar f \to f$, where $\bar f$ is the unique representative in the isomorphism class of $f$. We impose the convention that if $f = \bar f$ happens to be our chosen representative, then $\eta_{\bar f} = id_{\bar f}$. 
 
One choice of inverse functor $P:  \sB(x,y) \to  \sB^s(x,y)$ is defined on objects as the identity, on 1-morphism by $P(f) =: \overline f$ and on 2-morphisms as the unique 2-morphism that makes the following diagram commute:
\begin{center}
\begin{tikzpicture}[thick]
	\node (LT) at (0,1.5) 	{$\bar f$ };
	\node (LB) at (0,0) 	{$f$};
	\node (RT) at (3,1.5) 	{$\bar g$};
	\node (RB) at (3,0)	{$g$};
	\draw [->] (LT) --  node [left] {$\eta_f$} (LB);
	\draw [->, dashed] (LT) -- node [above] {$P(\alpha) = : \overline{\alpha}$} (RT);
	\draw [->] (RT) -- node [right] {$\eta_g$} (RB);
	\draw [->] (LB) -- node [below] {$\alpha$} (RB);
\end{tikzpicture}
\end{center}
In particular, the composition  $P \circ i = id_{\sB^s(x,y)}$ (strict equality) and the $\eta_f$ provide the components of a natural isomorphism $\eta: i \circ P \to id_{\sB(x,y)}$. Strictly speaking, both $P$ and $\eta$ depend on $x$ and $y$, but we will omit this dependence from our notation.

We equip $\sB^s$ with identity 1-morphism $P(I_x) \in \sB^s(x,x)$ and with compositions:
\begin{equation*}
\sB^s(y,z) \times \sB^s( x,y) \stackrel{i \times i}{\hookrightarrow} \sB(y,z) \times \sB( x,y) \stackrel{c}{\to} \sB(x,z) \stackrel{P}{\to} \sB^s(x,z)
\end{equation*}
i.e., the composition of 1-morphisms $f$ and $g$ is given by $P(i(f) \circ i(g))$. The fully-faithful functor $i$ is automatically injective, and so we may identify $f \in \sB^s(x,y)$ with its image $i(f) \in \sB(x,y)$, and thereby simply write $P( f \circ g)$ for the composition of $f$ and $g$ in $\sB^s$. Finally, we must specify the associators and unitors. These are given on components by:
\begin{equation*}
a^{\sB^s}: P(P(f \circ g) \circ h) \stackrel{P( \eta * 1) }{\to} P( (f \circ g) \circ h) \stackrel{P(a^\sB)}{\to} P(f \circ (g \circ h)) \stackrel{P(1 * \eta^{-1})}{\to} P( f \circ P(g \circ h)) 
\end{equation*}
\begin{equation*}
r^{\sB^s}:	P( f \circ P(I)) \stackrel{P(1 * \eta)}{\to} P( f \circ I) \stackrel{P(r^\sB)}{\to} P(f) = f
\end{equation*}
\begin{equation*}
\ell^{\sB^s}:	P( P(I) \circ f) \stackrel{P( \eta * 1)}{\to} P( I \circ f) \stackrel{P(\ell^\sB)}{\to} P(f) = f
\end{equation*}
These are easily verified to satisfy the pentagon and triangle identities, and hence give $\sB^s$ the structure of a bicategory. We will return to these identities after finishing the rest of this proof. 

The functors $i:  \sB^s(x,y) \to \sB(x,y)$ and $P: \sB(x,y) \to \sB^s(x,y)$ assemble into a pair of homomorphisms: $i: \sB^s \to \sB$ and $P: \sB \to \sB^s$ which are given as follows:
The functor $i = (id, i, \phi^i_{g,f}, \phi^i_x)$ is the identity on the objects and on hom categories is given by the component functors $i: \sB^s(x,y) \to \sB(x,y)$. The natural transformations $\phi^i$ are given by $ \phi^i_{g,f} = \eta^{-1}  : g \circ f \to P(g \circ f)$ and $\phi^i_x = \eta^{-1}: I_x \to P(I_x)$.  

Similarly the homomorphism $P = (id, P, \phi^P_{g,f}, \phi^P_x)$ is the identity on objects and on hom categories is given by $P: \sB(x,y) \to \sB^s(x,y)$. The natural transformations $\phi^P$ are given by $ \phi^P_{g,f} = P(\eta * \eta)  : P( P(g) \circ P(f)) \to P(g \circ f)$ and $\phi^i_x = id: P(I_x) \to P(I_x)$.  With these choices the composite homomorphism $P \circ i: \sB^s \to \sB^s$ is precisely the identity homomorphism, as the reader can readily check. In particular the morphisms,
\begin{equation*}
P(I_x) \stackrel{id}{\to} P(I_x) \stackrel{P( \eta^{-1}) }{\to} P(I_x)
\end{equation*}
\begin{equation*}
P( g \circ f) = P(P(g) \circ P(f))  \stackrel{P( \eta * \eta) }{\to} P( g \circ f) \stackrel{P( \eta^{-1})}{\to} P( g \circ f))
\end{equation*}
are identities. 

The reverse composition is the homomorphism $i \circ P = (id, i \circ P, \phi^{i \circ P}_{f,g}, \phi^{i \circ P}_x)$ where these natural transformations are given on components by:
\begin{equation*}
\phi^{i \circ P}_{f,g}: \bar f \circ \bar g \stackrel{\eta^{-1}}{\to} \overline{\bar f \circ \bar g} \stackrel{\overline{\eta * \eta} }{\to } \overline{f \circ g},
\end{equation*}
\begin{equation*}
\phi^{i \circ P}_x: I_x \stackrel{\eta^{-1}}{\to} \overline{I}_x.
\end{equation*}
Recall that $\overline{ f\circ g}$ means $P(f \circ g)$. This is not the identity homomorphism, but it is equivalent to the identity by a transformation whose components are $\sigma = (I_x, \eta_f)$. Thus the bicategories $\sB$ and $\sB^s$ are equivalent. 

Now let us return to the verification of the pentagon and triangle identities and the choice of associator $a^{\sB^s}$ and unitors, above. The associator for $\sB^s$ is a morphism from $\overline{(\overline{f \circ g}) \circ h}$ to $\overline{ f \circ ( \overline{g \circ h})}$. The definition of $a^{\sB^s}$ we gave above used a combination of the $\eta$ 2-morphisms to construct isomorphisms:
\begin{align*}
	\overline{(\overline{f \circ g}) \circ h} & \cong (f \circ g) \circ h \\
	\overline{ f \circ ( \overline{g \circ h})} & \cong f \circ (g \circ h)
\end{align*}
and then used these to pullback the associator of $\sB$,  $a: (f \circ g) \circ h \to f \circ (g \circ h)$. One might ask the following question: why did we choose these particular isomorphisms? For example the first isomorphism, as specified in the above proof, is the composite: 
\begin{equation*}
	\overline{(\overline{f \circ g}) \circ h} \to \overline{(f \circ g) \circ h} \to  (f \circ g) \circ h,
\end{equation*}
but we could have just as easily insisted on using the composite:
\begin{equation*}
	\overline{(\overline{f \circ g}) \circ h} \to (\overline{f \circ g}) \circ h \to  (f \circ g) \circ h.
\end{equation*}
In fact, these isomorphisms are the same, as can be seen from the naturality of $\eta$. However more is true. A marked bracketed expression, such as  $\overline{(\overline{f \circ g}) \circ h} $ determines a marked binary tree. In this case it is the tree $((-) \boxdot (-)) \boxdot (-)$ and the marked edges are precisely those corresponding to $f \circ g$ and $(f \circ g) \circ h$. 

Each marked edge corresponds to a single application of the functor $P$, which may be removed by composing with a single instance of the $\eta$ 2-morphisms. Thus a path between marked binary trees corresponds to a particular choice of composite, such as the two considered above. These various compositions are often equal, and we may impose an equivalence relation on the paths between marked binary trees such the two paths are equivalent precisely when the corresponding composites of 2-morphisms are equal. The naturality of the natural transformation $\eta$ implies that property (iii) of Lemma \ref{AnyPathsofTreesareEquivalent} is satisfied by this equivalence relation. Similarly, the functoriality of $\circ$ implies both properties (i) and (ii). Thus, by Lemma \ref{AnyPathsofTreesareEquivalent} we see that any two paths are equivalent, and hence that any two corresponding composites of the $\eta$'s coincide.

This implies that given any marked bracketed expression, such as $\overline{(\overline{(\overline{f \circ g} ) \circ h} ) \circ j}$ used in the verification of the pentagon identity, there is a {\em canonical} isomorphism with the same bracketed expression, equipped with fewer markings, and moreover these canonical isomorphisms are coherent. With this observation at hand the verification of the pentagon axiom is straightforward. 
\end{proof}

\begin{lemma}
Every symmetric monoidal bicategory is equivalent to a 1-skeletal symmetric monoidal bicategory.
\end{lemma}

\begin{proof}
Let $\sB$ be a symmetric monoidal bicategory, and construct the 1-skeletal bicategory $\sB^s$ and equivalence of bicategories $\sB^s \simeq \sB$ as in the previous lemma. We must transfer the symmetric monoidal structure from $\sB$ to $\sB^s$ and promote the equivalence $\sB^s \simeq \sB$ to a symmetric monoidal equivalence. 

$\sB^s$ becomes a symmetric monoidal bicategory when equipped with the following additional structure. The unit $1 \in \sB$ is the identical object of $\sB^s$ (recall $\sB$ and $\sB^s$ have the same objects). 
The tensor product $\otimes^{\sB^s}$ on $\sB^s$ is defined to be the composite homomorphism:
\begin{equation*}
	\sB^s \times \sB^s \stackrel{i \times i}{\hookrightarrow} \sB \times \sB \stackrel{\otimes}{\to} \sB \stackrel{P}{\to} \sB^s.
\end{equation*}
Thus on objects, $\otimes^{\sB^s}$ agrees with $\otimes^\sB$, and on 1-morphisms and 2-morphisms it is given by $f \otimes^{\sB^s} g := \overline{f \otimes g}$. Lemma \ref{AnyPathsofTreesareEquivalent} again implies the existence of canonical coherent isomorphisms between marked bracketed expressions involving any mixture of $\circ$ and $\otimes$ (property (iii) follows from the naturality of $\eta$ and properties (i) and (ii) follow from the functoriality of $\circ$ and $\otimes$. Thus we may define $\alpha^{\sB^s}$ (and all the remaining structure) by conjugating the structure on $\sB$ by these canonical isomorphisms. The relevant coherence diagrams are then automatically satisfied. 
Explicitly, we have the following structure transformations:
\begin{align*}
& \alpha^{\sB^s} = ( \bar \alpha_{abc}, \eta^{-1} \circ (1 \otimes \eta^{-1}) \circ \alpha_{fgh} \circ (\eta \otimes 1) \circ \eta ) = (  \bar \alpha_{abc}, \textrm{can}^{-1} \circ \alpha_{fgh} \circ \textrm{can})  \\
& \ell^{\sB^s} = ( \bar \ell_a, \ell_f \circ \eta) \\
& r^{\sB^s} = ( \bar r, r_f \circ \eta) \\
& \beta^{\sB^s} = (\bar \beta, \eta^{-1} \circ \beta_{f,g} \circ \eta)
\end{align*}
together with structure modifications defined by the pasting diagrams in Figures \ref{fig:PastingDiagSkelitalization1}, \ref{fig:PastingDiagSkelitalization2}, \ref{fig:PastingDiagSkelitalization3}, \ref{fig:PastingDiagSkelitalization4}, \ref{fig:PastingDiagSkelitalization5}, and \ref{fig:PastingDiagSkelitalization6}.
Here the unlabeled 2-morphisms are the natural 2-morphisms $\eta$. With these structures $\sB^s$ becomes a symmetric monoidal bicategory. 

\begin{figure}[ht]
	 \begin{center}

	\end{center}
	\caption{Pasting Diagram for $\pi^{\sB^s}$}
	\label{fig:PastingDiagSkelitalization1}
\end{figure}

\begin{figure}[ht]
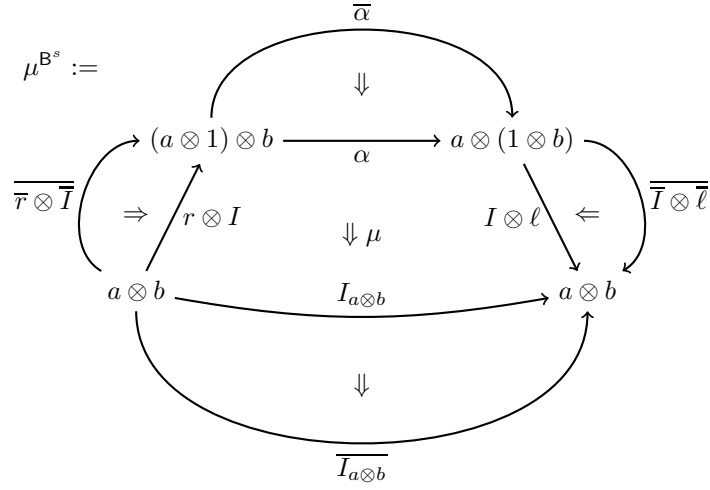

	\begin{center}
	%
	\end{center}
	\caption{Pasting Diagram for $\mu^{\sB^s}$}
	\label{fig:PastingDiagSkelitalization2}
\end{figure}

\begin{figure}[ht]
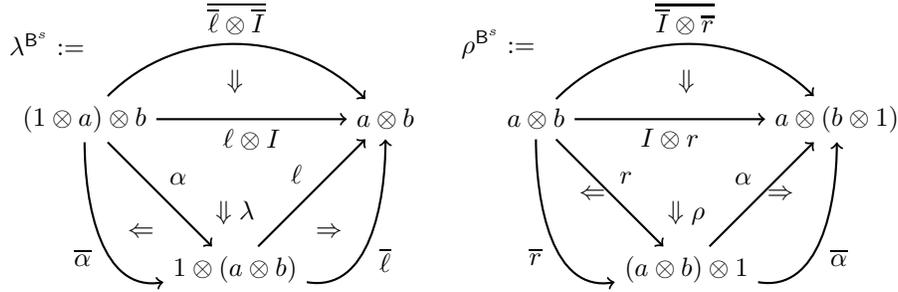

	
	\begin{center}
	%
	\end{center}
	\caption{Pasting Diagram for $\lambda^{\sB^s}$ and $\rho^{\sB^s}$}
	\label{fig:PastingDiagSkelitalization3}
\end{figure}

\begin{figure}[ht]
	\begin{center}
	%
	\end{center}	
	\caption{Pasting Diagram for $R^{\sB^s}$}
	\label{fig:PastingDiagSkelitalization4}
\end{figure}

\begin{figure}[ht]
	\begin{center}
	%
	\end{center}
	\caption{Pasting Diagram for $S^{\sB^s} $}
	\label{fig:PastingDiagSkelitalization5}
\end{figure}

\begin{figure}[ht]
	\begin{center}
	%
	\end{center}
	\caption{Pasting Diagram for $\sigma^{\sB^s}$ }
	\label{fig:PastingDiagSkelitalization6}
\end{figure}

The inclusion homomorphism $i: \sB^s \to \sB$ becomes a symmetric monoidal homomorphism when equipped with the following structure. The structure transformations are given by:
\begin{align*}
\chi^i & = (I, f \otimes g \stackrel{\textrm{can}^{-1}}{\to} \overline{f \otimes g}) \\
\iota^i &= (I, id)
\end{align*}
and the structure modifications $\omega^i, \gamma^i, \delta^i, u^i$ are given by the canonical 2-morphism such as $\omega^i$ depicted in Figure~\ref{fig:PastingDiagInclusHom}.
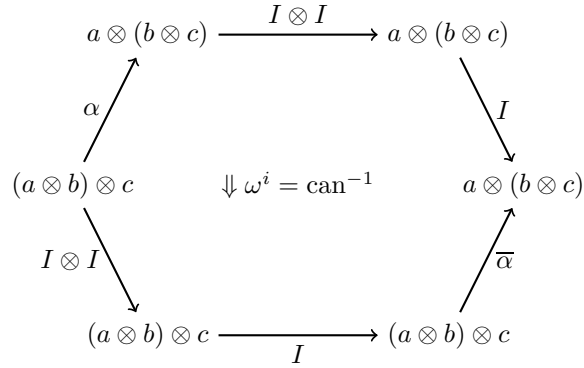
\begin{figure}[ht]
	\begin{center}
	\begin{tikzpicture}[thick]
		\node (LT) at (0,4) 	{$a \otimes (b \otimes c)$ };
		\node (LM) at (-1,2) 	{$(a \otimes b) \otimes c$ };
		\node (LB) at (0,0) 	{$( a\otimes b) \otimes c$};

		\node (RT) at (4,4) 	{$a \otimes (b \otimes c)$};
		\node (RM) at (5,2) 	{$ a \otimes (b \otimes c)$ };
		\node (RB) at (4,0)	{$(a \otimes b) \otimes c$};
		\draw [->] (LM) --  node [left] {$I \otimes I$} (LB);
		\draw [->] (LM) --  node [left] {$\alpha$} (LT);	
		\draw [->] (LT) -- node [above] {$I \otimes I$} (RT);
		\draw [->] (LB) -- node [below] {$I$} (RB);
		\draw [->] (RT) -- node [right] {$I$} (RM);
		\draw [->] (RB) -- node [right] {$\overline{\alpha}$} (RM);

		\node at (2,2) {$ \Downarrow \omega^{i} = \text{can}^{-1} $};	
	\end{tikzpicture}
	\end{center}
	\caption{Pasting Diagram for $\omega^i$}
	\label{fig:PastingDiagInclusHom}
\end{figure}

Similarly, $P: \sB \to \sB^s$ becomes a symmetric monoidal homomorphism when equipped with the structure transformations:
\begin{align*}
\chi^i & = (I, \overline{\overline f \otimes \overline g} \stackrel{\textrm{can}}{\to} \overline{f \otimes g}) \\
\iota^i &= (I, id)
\end{align*}
together with the canonical structure 2-morphisms $\omega^P, \gamma^P, \delta^P, u^P$ such as $\omega^P$ depicted in Figure~\ref{fig:PastingDiagPHom}.
\begin{figure}[ht]
\begin{center}
\begin{tikzpicture}[thick]
	\node (LT) at (0,4) 	{$a \otimes (b \otimes c)$ };
	\node (LM) at (-1,2) 	{$(a \otimes b) \otimes c$ };
	\node (LB) at (0,0) 	{$( a\otimes b) \otimes c$};
	
	\node (RT) at (4,4) 	{$a \otimes (b \otimes c)$};
	\node (RM) at (5,2) 	{$ a \otimes (b \otimes c)$ };
	\node (RB) at (4,0)	{$(a \otimes b) \otimes c$};
	\draw [->] (LM) --  node [left] {$I \otimes I$} (LB);
	\draw [->] (LM) --  node [left] {$\overline{\alpha}$} (LT);	
	\draw [->] (LT) -- node [above] {$I \otimes I$} (RT);
	\draw [->] (LB) -- node [below] {$I$} (RB);
	\draw [->] (RT) -- node [right] {$I$} (RM);
	\draw [->] (RB) -- node [right] {$\overline{\alpha}$} (RM);
	
	\node at (2,2) {$ \Downarrow \omega^{P} = \text{can} $};	
\end{tikzpicture}
\end{center}
	\caption{Pasting Diagram for $\omega^P$}
	\label{fig:PastingDiagPHom}
\end{figure}

The composition of $P$ and $i$, equipped with these symmetric monoidal structures, can now be readily computed. We have, in both cases, that the resulting compositions are equivalent to the corresponding identity homomorphism, and thus $i$ is an equivalence of symmetric monoidal bicategories. 
\end{proof}

\begin{lemma} \label{EverySymBicatIsEquivToSkeletal}
Every symmetric monoidal bicategory is equivalent to a skeletal symmetric monoidal bicategory. 
\end{lemma}

\begin{proof}
Let $\sB$ be a symmetric monoidal bicategory. By the previous lemma, we may assume, without loss of generality, that $\sB$ is 1-skeletal. Choose for each equivalence class of objects in $\sB$ a representative object and let $\sB^{sk}$ denote the full sub-bicategory of $\sB$ spanned by these objects. $\sB^{sk}$ is a skeletal bicategory. We will equip $\sB^{sk}$ with the structure of a symmetric monoidal bicategory and promote the inclusion functor $j: \sB^{sk} \to \sB$ to a symmetric monoidal equivalence. Note that a similar transference-of-structure has been shown to hold for monoidal bicategories along general equivalences of bicategories \cite[Thm.~5.1]{MR2972968}.

Choose for each object $x \in \sB$ an equivalence $\xi_x: \overline{x} \to x$ in $\sB$, such that $\xi_{\overline x} = I_{\overline x}$. Since $\sB$ is 1-skeletal, there exists a unique 1-morphism $\xi^{-1}: x \to \overline x$ such that $\xi \xi^{-1}$ and $\xi^{-1} \xi$ are identity 1-morphisms (strict equality). $j: \sB^{sk} \to \sB$ is an equivalence of bicategories and an inverse equivalence is given by $Q: \sB \to \sB^{sk}$ which is defined as follows: 
\begin{itemize}
\item For objects $x \in \sB$, $Q(x) = \overline x$.
\item For 1-morphisms $f: x \to y$, we have $Q(f) =  \overline f := \overline x \stackrel{\xi}{\to} x \stackrel{f}{\to} y \stackrel{\xi^{-1}}{\to} \overline y$. Note that because $\sB$ is 1-skeletal the order of this composition is irrelevant. 
\item For 2-morphisms $\alpha: f \to g$, we define $Q(\alpha) = \overline \alpha$ to be the horizontal composition $[ id_\xi * \alpha] * id_{\xi^{-1}}$. 
\end{itemize}
Just as before, we have that $Q \circ j$ is the identity homomorphism on $\sB^{sk}$ and that $j \circ Q$ is equivalent to the identity homomorphism on $\sB$, with $\xi$ providing the structure of that equivalence.

$\sB^{sk}$ becomes a symmetric monoidal bicategory when equipped with the following additional structure. The unit is given by $\overline 1 \in \sB^{sk}$. The tensor product $\otimes^{\sB^{sk}}$ on $\sB^{sk}$ is defined to be the composite homomorphism:
\begin{equation*}
	\sB^{sk} \times \sB^{sk} \stackrel{j \times j}{\hookrightarrow} \sB \times \sB \stackrel{\otimes}{\to} \sB \stackrel{Q}{\to} \sB^{sk}.
\end{equation*}
Thus on objects, $\otimes^{\sB^{sk}}$ is given by $x \otimes^{\sB^{sk}} y = \overline { x \otimes y}$, and similarly on 1-morphisms it is given by $\overline {f \otimes g}$ and on 2-morphisms by $\overline{ \zeta \otimes \kappa}$. 
 Lemma \ref{AnyPathsofTreesareEquivalent} implies the existence of canonical coherent 1-morphisms
   between marked bracketed expressions involving the objects of $\sB^{sk}$ and $\otimes$. 
   Thus we may define the symmetric monoidal structure on $\sB^{sk}$ by conjugating the corresponding structures on $\sB$ by these canonical 1-morphisms. For example the pentagonator $\pi^{\sB^{sk}}$ is given by the following pasting diagram in Figure~\ref{fig:PastingDiagSkelPentagonator}.
Here the symbol ``$\#$'' means the identity 2-morphism, i.e., that the diagram commutes one the nose. It is straightforward to check that the relevant coherence diagrams are automatically satisfied. 

\begin{figure}[ht]
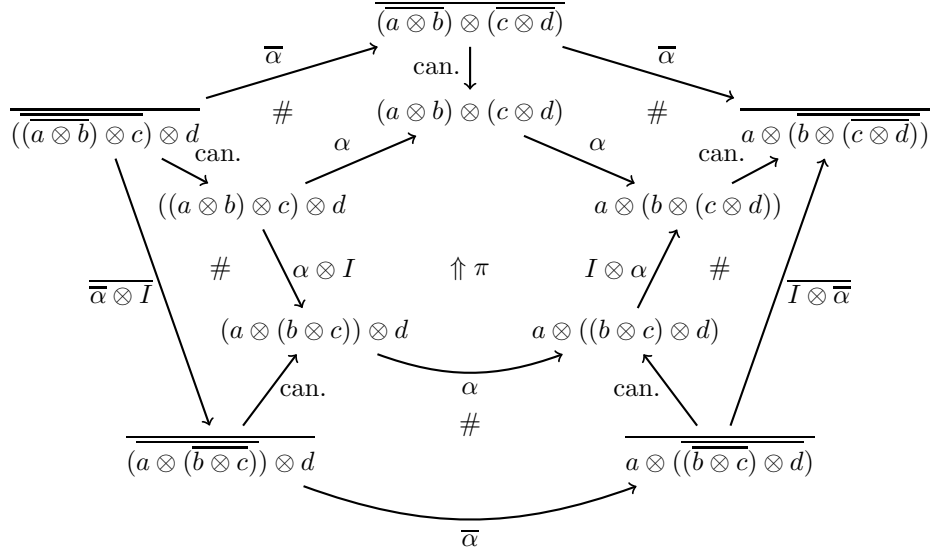

	\begin{center}

	\end{center}
	\caption{Pasting Diagram for the Pentagonator $\pi^{\sB^{sk}}$}
	\label{fig:PastingDiagSkelPentagonator}
\end{figure}

We can promote $j$ to a symmetric monoidal homomorphism as well. We equip it with the following structure transformations and modifications: $\chi^j$ is the transformation which in components is given by the canonical 1-morphism:
\begin{equation*}
	\chi^j_{a,b}: a \otimes b \to \overline{ a \otimes b}.
\end{equation*}
There is also a component 2-morphism for $\chi$. This 2-morphism can be taken to be the identity, which fills the pasting diagram in Figure~\ref{fig:PastingDiagComponentForChi} (when appropriately bracketed)\footnote{Actually the bracketing which allows this to be filled by the identity 2-morphism is not compatible with the bracketing needed to define $\chi^j_{f,g}$, however there is a canonical $\chi^j_{f,g}$ given by rebracketing the identity 2-morphism appropriately.}. Similarly $\iota^j: 1 \to \overline 1$ is the canonical 1-morphism and we can take $\omega^j, \gamma^j, \delta^j$ and $u^j$ to be identity 2-morphisms.  
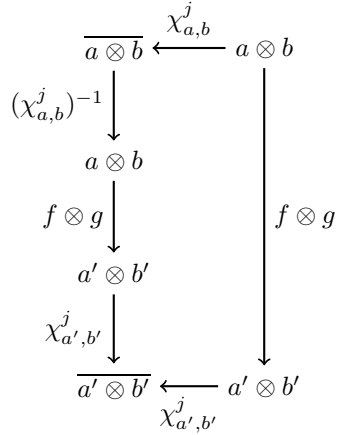
\begin{figure}[ht]
	\begin{center}
	\begin{tikzpicture}[thick]
		\node (LT) at (0,4.5) 	{$\overline{ a \otimes b}$ };
		\node (LMT) at (0,3) 	{${ a \otimes b}$ };
		\node (LMB) at (0,1.5) 	{${ a' \otimes b'}$ };
		\node (LB) at (0,0) 	{$\overline{ a' \otimes b'}$};
		\node (RT) at (2,4.5) 	{${ a \otimes b}$};
		\node (RB) at (2,0)	{${ a' \otimes b'}$};
		\draw [->] (LT) --  node [left] {$(\chi^j_{a,b})^{-1}$} (LMT);
		\draw [->] (LMT) --  node [left] {$f \otimes g$} (LMB);
		\draw [->] (LMB) --  node [left] {$\chi^j_{a',b'}$} (LB);
		\draw [<-] (LT) -- node [above] {$\chi^j_{a,b}$} (RT);
		\draw [->] (RT) -- node [right] {$f \otimes g$} (RB);
		\draw [<-] (LB) -- node [below] {$\chi^j_{a',b'}$} (RB);
	\end{tikzpicture}
	\end{center}
	\caption{Pasting Diagram for component 2-morphism for $\chi$}
	\label{fig:PastingDiagComponentForChi}
\end{figure}

 $Q$ can be made into a symmetric monoidal homomorphism in the same manner by equipping it with the transformations $\chi^Q_{a,b} = \textrm{can.} : \overline{\overline a \otimes \overline b} \to \overline{ a \otimes b}$ and $\iota^Q = I_{\overline 1}$. All the additional structure 2-morphisms can be taken to be identities. $Q$ and $j$ are inverse equivalences of symmetric monoidal bicategories. 
\end{proof}

\begin{lemma} \label{caninverseequivforskeletalbicats}
Let $F : \sB \to \sC$ be an equivalence between skeletal bicategories. Then there exists a canonical inverse homomorphism $F^{-1}: \sC \to \sB$ such that  $F^{-1} \circ F = id_\sB$ and $F \circ F^{-1} = id_\sC$ {\em on the nose}.
\end{lemma}

\begin{proof}
Since $\sB$ and $\sC$ are skeletal and since $F$ is an equivalence, the components $F_0, F_1, F_2$ are all bijections with inverses $F^{-1}_i: \sC_i \to \sB_i$. Let the components of $\eta$ be given by
\begin{align*}
	 \eta_{f,g} = F^{-1}(({\phi^F}_{F^{-1}(f), F^{-1}(g)})^{-1})  \\
	\eta_a = F^{-1}((\phi^F_{F^{-1}(a)})^{-1}):  I_{F^{-1}(a)} \to F^{-1}( I_a)
\end{align*}
 for all 1-morphisms $f,g \in \sC$ and all objects $a \in \sC$.  
Applying the maps $F^{-1}_i$ to the coherence diagrams of $(F, \phi)$ yield precisely the necessary coherence diagrams for $(F^{-1}, \eta)$ to be a homomorphism. The compositions of $F$ and $F^{-1}$ can easily be checked to be the identity homomorphisms.
\end{proof}

Up to this point we have mainly been describing symmetric monoidal bicategories with the language of pasting diagrams. The proof of the next lemma, however, will need to make use of {\em mates}, as described in Appendix \ref{PastingsStringsAdjointsAndMatesSection}. The language of string diagrams is far better suited to this task, and so we will use it below. The unfamiliar reader should consult  Appendix \ref{Chap:Bicats}. 

\begin{lemma} \label{SymMonHomBetweenSkeletalBicats}
A symmetric monoidal homomorphism $F: \sM \to \sM'$ between skeletal symmetric monoidal bicategories is an equivalence if and only if it is an equivalence of underlying bicategories, i.e., the following properties are satisfied:
\begin{enumerate}
\item (essentially surjective on objects) It is surjective on $\pi_0$. 
\item (essentially full on 1-morphism) $F: \pi_0\sM(x, y) \to \pi_0\sM'(Fx, Fy)$ is surjective.
\item (fully-faithful on 2-morphisms) $F: \sM(x, y) \to \sM'(Fx, Fy)$ is fully-faithful.   
\end{enumerate}
\end{lemma}

\begin{proof}
The necessity of $F$ to be an equivalence of underlying bicategories is clear. We construct the inverse symmetric monoidal equivalence as follows: First construct the canonical inverse equivalence $(F^{-1}, \eta)$ of bicategories from Lemma \ref{caninverseequivforskeletalbicats}. Then we construct the remaining coherence data for $F^{-1}$ by applying $F^{-1}$ to certain mates of the coherence data of $F$. Specifically, $\chi^{F^{-1}} $ is given by whiskering $(\chi^F)^{-1}$ by $F^{-1}$ and $F^{-1} \times F^{-1}$. Note that ``$(\chi^F)^{-1}$'' makes sense since we are dealing with {\em skeletal} bicategories. 
 On objects we have $\chi^{F^{-1}}_{x,y} := F^{-1}( (\chi_{F^{-1}x, F^{-1}y}^F)^{-1} )$. We can similarly, define $\iota^{F^{-1}}$ as $F^{-1}((\iota^F)^{-1})$. 

The higher coherence data is constructed as follows. Applying $u^F$ to the objects $F^{-1}(x)$ and $F^{-1}(y)$ yields the 2-morphism in Figure~\ref{fig:SymMonHomBetweenSkeletalBicats1}. In string diagrams we can represent this 2-morphism as in Figure~\ref{fig:SymMonHomBetweenSkeletalBicats2}.
We can then form the mate which is depicted in Figure~\ref{fig:SymMonHomBetweenSkeletalBicats3}. Applying $F^{-1}$ to the mate $u^*$ yields $u_{x,y}^{F^{-1}}$. 

\begin{figure}[ht]
	\begin{center}

	\end{center}
	\caption{Constructing Coherence Data}
	\label{fig:SymMonHomBetweenSkeletalBicats1}
\end{figure}

\begin{figure}[ht]
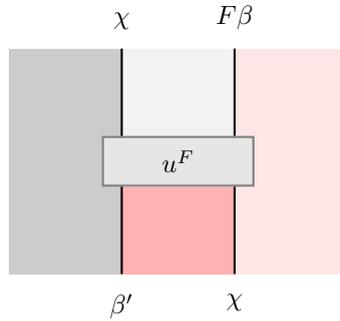

	\begin{center}
	%
	\end{center}
	\caption{A String Diagram Version of Figure~\ref{fig:SymMonHomBetweenSkeletalBicats1}}
	\label{fig:SymMonHomBetweenSkeletalBicats2}
\end{figure}

\begin{figure}[ht]
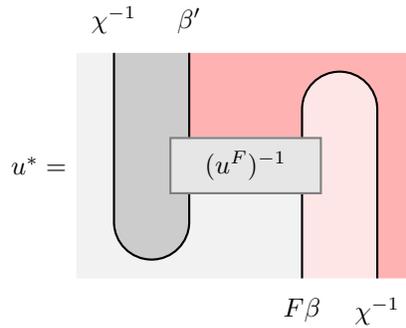

	\begin{center}
	%
	\end{center}
	\caption{The Mate of the 2-Morphism of Figure~\ref{fig:SymMonHomBetweenSkeletalBicats2}}
	\label{fig:SymMonHomBetweenSkeletalBicats3}
\end{figure}

%%%%%%%%%

Similarly we may construct $\gamma^{F^{-1}}, \delta^{F^{-1}},$ and $\omega^{F^{-1}}$. Applying $\gamma^F$ to the object $F^{-1}x$ yields the 2-morphism of Figure~\ref{fig:SymMonHomBetweenSkeletalBicats4}, which we can represent by the string diagram in Figure~\ref{fig:SymMonHomBetweenSkeletalBicats5}. 
\begin{figure}[ht]
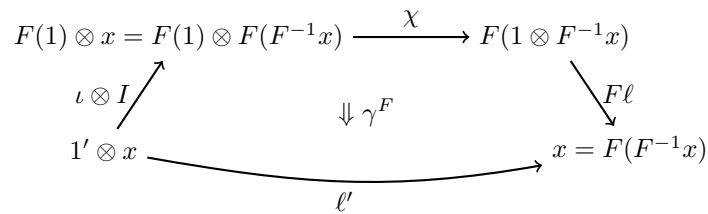

	\begin{center}
	%
	\end{center}
	\caption{More Coherence Data}
	\label{fig:SymMonHomBetweenSkeletalBicats4}
\end{figure}
\begin{figure}[ht]
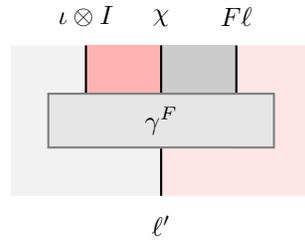

	\begin{center}
	%
	\end{center}
	\caption{The String Diagram of the 2-Morphism of Figure~\ref{fig:SymMonHomBetweenSkeletalBicats4}}
	\label{fig:SymMonHomBetweenSkeletalBicats5}
\end{figure}
Now we can form the 2-morphism depicted in Figure~\ref{fig:SymMonHomBetweenSkeletalBicats6},
\begin{figure}[ht]
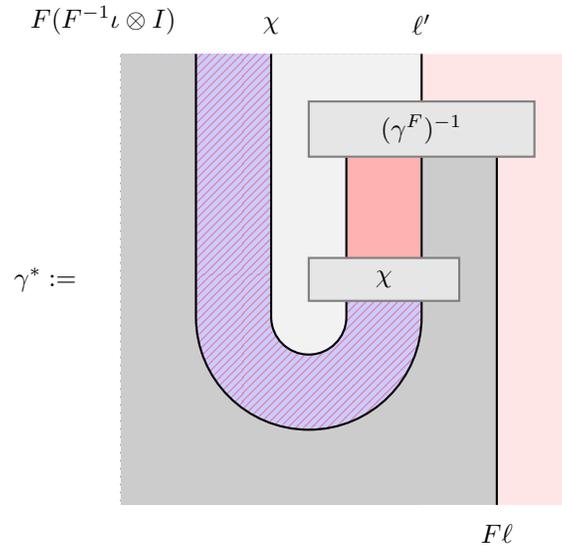

	\begin{center}
	%
	\end{center}
	\caption{A 2-Morphisms built from the 2-Morphism of Figure~\ref{fig:SymMonHomBetweenSkeletalBicats5}}
	\label{fig:SymMonHomBetweenSkeletalBicats6}
\end{figure}
where the 2-morphism labeled ``$\chi$'' is part of the structure making $\chi$ a transformation. It fills the following diagram:
\begin{center}
%
\end{center}
where we have identified $\iota \otimes I = F(F^{-1}(\iota)) \otimes F(F^{-1}(I_x))$. The 2-morphism $\gamma^{F^{-1}}$ is obtained by applying $F^{-1}$ to $\gamma^*$. The construction of $\delta^{F^{-1}}$ and $\omega^{F^{-1}}$ are completely analogous. 

Using string diagrams it is relatively easy, but quite tedious, to verify that these do indeed satisfy the necessary axioms to make $F^{-1}$ into a symmetric monoidal homomorphism. Moreover a direct calculation shows that the compositions $F \circ F^{-1}$ and $F^{-1} \circ F$ are exactly (with equality on the nose) the respective identity symmetric monoidal homomorphisms. 
\end{proof}

\begin{theorem}[Whitehead's Theorem for Symmetric Monoidal Bicategories] \label{WhiteheadforSymMonBicats}
Let $\sB$ and $\sC$ be symmetric monoidal bicategories. A symmetric monoidal homomorphism $F:\sB \to \sC$ is a symmetric monoidal equivalence if and only if it is an equivalence of underlying bicategories, that is $F$ is essentially surjective on objects, essentially full on 1-morphisms, and fully-faithful on 2-morphisms. 
\end{theorem}

\begin{proof}
The necessity of $F$ to be an equivalence of underlying bicategories is clear. Conversely, let $F$ be a symmetric monoidal homomorphism, which is an equivalence of underlying bicategories. By Lemma \ref{EverySymBicatIsEquivToSkeletal}, $\sB$ is equivalent as a symmetric monoidal bicategory to a skeletal symmetric monoidal bicategory $\sB^{sk}$, and similarly $\sC$ is equivalent to a skeletal symmetric monoidal bicategory $\sC^{sk}$. The symmetric monoidal homomorphism $F$ is a symmetric monoidal equivalence if and only if the induced homomorphism $F^{sk}: \sB^{sk} \to \sC^{sk}$ is a symmetric monoidal equivalence. Our assumption that $F$ is an equivalence of underlying bicategories, together with Lemma \ref{SymMonHomBetweenSkeletalBicats}, ensures this is the case. 
\end{proof}

\index{Whitehead's theorem!for symmetric monoidal bicategories|)}
\index{symmetric monoidal bicategory!Whitehead's theorem|)}
\index{bicategory!Whitehead's theorem for symmetric monoidal|)}

\section{Stricter versions of Symmetric Monoidal Bicategories} \label{sec:strictsymbicats}

As mentioned previously in Section~\ref{sec:historysymmonbicat}, there have been several versions of {\em semistrict symmetric monoidal 2-categories} in the literature. These are variants of symmetric monoidal bicategories which are more rigid than the fully weak notion. In these theories roughly half the coherence isomorphisms are identities. More precisely the Gordon-Powers-Street coherence theorem \cite{GPS95} dictates which coherence cells of the monoidal bicategory structure may be taken to be trivial; the remaining coherence constraints are either left unconstrained or required to satisfy a mild normalization condition. 

In this section we introduce three variants on the theme of symmetric monoidal bicategories, two of which are, to the best of our knowledge, new. The most rigid of these is our notion of {\em quasistrict symmetric monoidal 2-category}. This version of symmetric monoidal bicategory is almost, but not quite, entirely strict. All of the coherence cells are identities except the braiding and interchanger, and moreover there are certain additional normalizations on these remaining coherence cells, which we explain below. The significance of this notion will become apparent in Section~\ref{sec:coherence} in which we establish our coherence theorem for symmetric monoidal bicategories. We will show, among other things, that every symmetric monoidal bicategory is equivalent to a quasistrict symmetric monoidal 2-category. 

In a {\em semistrict} symmetric monoidal 2-category we will allow some of the coherence cells to be non-trivial, roughly half of them. The literature has many variants semistrict symmetric monoidal 2-category, which differ in small ways in precisely which coherence cells are required to be identities. The notion we review is essentially the symmetric monoidal variant of Crans' semistrict braided monoidal 2-categories \cite{MR1626844}.

The axioms of a symmetric monoidal bicategory are, naturally, built on the axioms of braided monoidal bicategories. These axioms favor the elementary exchanges: $\beta_{a,b}: a \otimes b \to b \otimes a$. In a {\em symmetric} monoidal structure one might equally well expect operations for every permutation in the symmetric group. It is possible to construct such operations out of the elementary exchanges, but not uniquely so due to the fact that elements of the symmetric group do not have unique representations in terms of elementary (braided) exchanges. 
%\CSP{Are there canonical representatives? Maybe using positive braids? cf GO.}

The third notion that we introduce is {\em unbiased semistrict symmetric monoidal 2-categories} which, like other semistrict notions, uses the Gordon-Powers-Street coherence theorem as a starting place. However rather than favoring the elementary exchanges, operations are defined for every permutation in the symmetric group. In terms of the amount of coherence data, this notion is on par with previous versions of semistrict symmetric monoidal 2-categories.
%, however the coherence axioms, the equations that the coherence data must satisfy, are more transparent, especially when expressed in terms of string diagrams. 
This notion will figure prominently in the next chapter when we analyze the bordism bicategory. Using the planar decomposition theorem from the Chapter~\ref{ChapPlanarDecomp} we will be able to replace the bordism bicategory with an equivalent 2-category which is naturally an unbiased semistrict symmetric monoidal 2-category.  
%\CSP{IS it really more transparent?}

Finally, we will also give an example that demonstrates that the notion of quasistrict symmetric monoidal 2-category is minimal in the sense that it is not possible to model all symmetric monoidal bicategories if any of the last remaining coherence cells are required to be identities.  

\begin{definition}\label{def:graymonoid}
	A {\em Gray monoid} is a monoidal bicategory $\sM$ such that: 
	\begin{enumerate}
		\item $\sM$ is a strict 2-category;
		\item The transformations $\alpha$, $\ell$, $r$, $\pi$, $\mu$, $\lambda$, and $\rho$ are identities. Moreover the inverse adjoint equivalences $\alpha^*$, $\ell^*$, and $r^*$ are also identities with trivial adjunction data. 
		\item The functor $\otimes = (\otimes, \phi_{(f,f'), (g, g')}^\otimes, \phi^\otimes_{(a,a')} )$ is {\em cubical}. This means that the 2-morphism witnessing naturality:
		\begin{equation*}
			\phi_{(f,f'), (g, g')}^\otimes: (f \otimes f') \circ (g \otimes g') \to (f \circ g) \otimes (f' \circ g')
		\end{equation*}
		 is an identity if either $f$ or $g'$ is an identity morphism. 
	\end{enumerate}
\end{definition}

\noindent In particular in a Gray monoid the operation $\otimes$ is strictly associative and unital, though the {\em interchanger $\phi^\otimes$} may still be non-trivial. Moreover the normalization axioms 
imply that $\phi^\otimes_{(a,a')}$ is also an identity. 

The coherence theorem of Gordon-Powers-Street \cite{GPS95} implies that any monoidal bicategory is equivalent to a Gray monoid. One may then transfer the symmetric monodial structure along such an equivalence, thereby obtaining a partial strictification theorem. In the braided monoidal case, this is essentially the notion considered in \cite{BN96}, which was modified by Crans \cite{MR1626844} with additional normalization conditions for the braiding. The following is the symmetric monoidal version of Crans' notion (c.f. \cite{GO13}).

\begin{definition}\label{def:CransSemiStrict}
	A {\em Crans semistrict symmetric monoidal 2-category} is a symmetric monoidal bicategory $\sM$ such that:
	\begin{enumerate}
		\item[(CSS.1)] The underlying monoidal bicategory is a Gray monoid;
		\item[(CSS.2)] The following additional normalization conditions apply: 
		\begin{enumerate}
			\item The 1-morphisms $\beta_{1,x}$ and $\beta_{x,1}$ are identity morphisms on $x$. 
			\item The isomorphisms $R_{1,x|y}$, $R_{x,1|y}$, $S_{x| 1,y}$, and $S_{x| y,1}$ are the identity 2-isomorphism of $\beta_{x,y}$.
			\item The isomorphisms $R_{x,y|1}$ and $S_{1|x,y}$ are the identity 2-isomorphisms of $I_{x\otimes y}$. 
		\end{enumerate}
	\end{enumerate}
\end{definition}
Note that in a Crans semistrict symmetric monoidal 2-category the seven Axioms (SM3.i), (SM3.ii), (SM4), (SM5), (SM6.i), (SM6.ii), and (SM7) of the definition of symmetric monoidal bicategory are not automatic, and must still be satisfied. Their form, though, is simplified due to the fact that $\sM$ is a Gray monoid.

In \cite{MR2770448} it was shown that every braided monoidal bicategory is equivalent to a semistrict braided monoidal bicategory satisfying the normalization conditions of Crans and in \cite{GO13} this was extended to the symmetric monoidal case as well. 

We will now introduce an even more strict notion of symmetric monoidal 2-category. 

\begin{definition} \label{def:quasistrict}
	A  {\em quasistrict symmetric monoidal 2-category} is a Crans semistrict symmetric monoidal 2-category $\sM$ such that:
	\begin{enumerate}
		\item[(QS.1)] The modifications $R$, $S$, and $\sigma$ are identities. (Which implies the normalizations (2)(b) and (2)(c) of Def.~\ref{def:CransSemiStrict}).
%		$\sM$ is a strict 2-category;
%		\item The following transformations and modifiacations are identities: $\alpha$, $\ell$, $r$, $\pi$, $\mu$, $\lambda$, $\rho$, $R$, $S$, and $\sigma$. Moreover the inverse adjoint equivalences $\alpha^*$, $\ell^*$, and $r^*$ are also identities with trivial adjunction data.  
%		\item The functor $\otimes = (\otimes, \phi_{(f,f'), (g, g')}^\otimes, \phi^\otimes_{(a,a')} )$ is {\em cubical}. This means that the 2-morphism witnessing naturality:
%		\begin{equation*}
%			\phi_{(f,f'), (g, g')}^\otimes: (f \otimes f') \circ (g \otimes g') \to (f \circ g) \otimes (f' \circ g')
%		\end{equation*}
%		 is an identity if either $f$ or $g'$ is an identity morphism. Note: this implies that $\phi^\otimes_{(a,a')}$ is also an identity. 
		\item[(QS.2)] For the transformation $\beta = (\beta_{(a,b)}, \beta_{(f,g)})$, the component $\beta_{(f,g)}$ is the identity if either $f$ or $g$ is an identity morphism. 
		%Furthermore, $\beta_{1b} = id_b$ when $a =1$ is the unit object of $\sC$; similarly $\beta_{a1} = id_a$. 
		\item[(QS.3)]  The 2-morphism witnessing naturality:
		\begin{equation*}
			\phi_{(f,f'), (g, g')}^\otimes
		\end{equation*}
		 is an identity if either $f$ or $g'$ is a component of $\beta$, i.e. if $f= \beta_{a,b}$ or $g' = \beta_{a,b}$.
	\end{enumerate}
\end{definition}

The notion of quasistrict symmetric monoidal 2-category has very little coherence data, and in fact much of that is redundant, as we will see. They also satisfy an nice `wire diagram' calculus \cite{B14}. 

\begin{lemma} \label{lma:QSbraidinggivenbyinterchange}
	In a quasistrict symmetric monoidal 2-category the coherence cells $\phi^\otimes_{(f, id), (id, g' )}$ determine the coherence cells $\phi^\otimes_{(f, f'), (g, g')}$, and together with the components $\beta_{a,b}$ determine the coherence cells $\beta_{f,f'}$. In particular $\beta_{f,f'}$ is an identity if either $f$ or $f'$ is a component of $\beta$. 
\end{lemma}

\begin{proof}
	Let $g: a \to b$, $g': a' \to b'$, $f:b \to c$, and $f': b' \to c'$.
	Consider the following square of 2-morphisms, which is commutative since $\otimes$ is a homomorphism:
	\begin{center}
	\begin{tikzpicture}
		\node (LT) at (0, 1.5) {$( id_{c} \otimes f') \circ (f \otimes id_{b'}) \circ (g \otimes g')$};
		\node (LB) at (0, 0) {$( id_{c} \otimes f') \circ [ (f \circ g) \otimes (id_{b'} \circ g')]$};
		\node (RT) at (8, 1.5) {$[(id_c \circ f) \otimes ( f' \circ id_{b'})] \circ (g \otimes g')$};
		\node (RB) at (8, 0) {$(f \circ g) \otimes (f' \circ g')$};
		\draw [->] (LT) -- node [left] {$1 * \phi_{(f, id_{b'}), (g,g')}$} (LB);
		\draw [->] (LT) -- node [above] {$\phi_{(id_{c}, f'), ( f, id_{b'})} * 1$} node [below] {$=$} (RT);
		\draw [->] (RT) -- node [right] {$\phi_{(f,f'),(g,g')}$} (RB);
		\draw [->] (LB) -- node [below] {$\phi_{(id_{c}, f'), (fg, g')}$} node [above] {$=$} (RB);
		%\node at (0.5, 1) {$\ulcorner$};
		%\node at (1.5, 0.5) {$\lrcorner$};
	\end{tikzpicture}
	\end{center}
This establishes that $\phi_{(f,id_{b'}), (g,g')}$ determines $\phi_{(f,f'),(g,g')}$. In a similar manner $\phi_{(f, id_{b'}), (g,g')}$ is determined by $\phi_{(f, id_{b'}), (id_b,g')}$.

Similarly, since $\beta$ is a natural transformation, it follows from our assumption that the symmetric monoidal bicategory is quasistrict that the composite,
\begin{center}
\begin{tikzpicture}
	\node (L) at (0, 0) {$(f' \otimes id_{c} ) \circ (id_{b'} \otimes f) \circ \beta_{(b,b')}$};
	\node (M) at (6, 0) {$(f' \otimes f) \circ \beta_{(b,b')}$};
	\node (R) at (10, 0) {$\beta_{c,c'} \circ (f \otimes f')$};
	\draw [->] (L) -- node [above] {$\phi_{(f',id_{c}),(id_{b'}, f)} * 1$} (M);
	\draw [->] (M) -- node [above] {$\beta_{(f,f')}$} (R);
\end{tikzpicture}
\end{center}
is the identity morphism. Hence $\beta_{(f,f')}$ is indeed determined by $\phi_{(f',id_{c}),(id_{b'}, f)}$. The finial statement of the lemma is also clear from this formula. 
\end{proof}

\begin{example}
	The following is a very important example of a quasistrict symmetric monoidal 2-category $\sS$, which cannot be strictified any further. It is skeletal. This 2-category has only one object (the unit object 1) and only two 1-morphisms $id_1$ and $f$. These satisfy the equation $f \circ f = id_1$, which implies that $f \otimes f = id_1$. Note that the 1-morphism are invertible. There are exactly two endomorphisms of $id_1$, which are both invertible: $id_{id_1}$ and $\beta_{f,f}$. This, plus the structural equations of a quasistrict symmetric monoidal 2-category fully determine $\sS$.
	The only non-trivial coherence cell is $\beta_{f,f}$ or equivalently $\phi_{(f,id), (id, f)}$.
	
	 This quasistrict symmetric monoidal 2-category corresponds to a stable 1-type, specifically the one for which both homotopy groups are $\Z/2\Z$ and the $k$-invariant is non-trivial, see \cite[App.B2]{MR2192936} and \cite{MR2981952}. The underlying monoidal 2-groupoid models a space which cannot be modeled by a strict monoid in 2-groupoids, i.e. by a one object strict 3-groupoid (see the argument from \cite[Sect.~2.7]{MR2883823}).
\end{example}

Recall that the symmetric group on $n$-letters has the following presentation: It is generated by $\sigma_1, \dots, \sigma_{n-1}$ subject to the relations
\begin{align*}
	\sigma_i^2 &= 1  \\
	\sigma_i \sigma_j &= \sigma_j \sigma_i & j \neq i, i \pm 1 \\
	\sigma_i \sigma_{i+1} \sigma_i &= \sigma_{i+1} \sigma_i \sigma_{i+1}.
\end{align*}
Here $\sigma_i$ corresponds to the permutation of the $i^\textrm{th}$ and $(i+1)^\textrm{st}$ elements. We will call this the standard presentation of $\Sigma_n$.  
Recall also the operadic composition of symmetric groups:
\begin{align*}
	\Sigma_n \times \Sigma_{m_1} \times \cdots \times \Sigma_{m_n} & \to \Sigma_{m_1 + \cdots + m_n} \\
	(\sigma, (\tau_1, \dots, \tau_n)) & \mapsto \sigma \circ (\tau_1, \dots, \tau_n). %= \tau_{\sigma(1)} \sqcup \cdots \sqcup \tau_{\sigma(n)}. 
\end{align*}
Finally, we have the $i^\textrm{th}$ deletion homomorphism $d_i: \Sigma_n \to \Sigma_{n-1}$, which erases the $i^\textrm{th}$ letter. 

\begin{lemma} \label{lem:qsisuss}
	Let $\sC$ be a quasistrict symmetric monoidal 2-category. 
	Then for each $n$ and each element $\sigma \in \Sigma_n$ of the symmetric group on $n$-letters there exists a natural transformation $\beta^\sigma: x_1 \otimes \cdots \otimes x_n \to x_{\sigma(1)} \otimes \cdots \otimes x_{\sigma(n)}$. These transformations are characterized by the following properties:
	\begin{enumerate}
		\item If $\sigma \in \Sigma_2$ is the non-trivial element, then 
		\begin{equation*}
			\beta^\sigma = \beta_{x_1,x_2}: x_1 \otimes x_2 \to x_2 \otimes x_1
		\end{equation*}
		is given by the braiding transformation of $\sC$.
		\item This association is compatible with composition and units: if $\sigma, \sigma' \in \Sigma_n$, then
		\begin{equation*}
			\beta^\sigma \circ \beta^{\sigma'} = \beta^{\sigma \sigma'}.
		\end{equation*}
		and $\beta^e = id$. 
		\item This association is compatible with the operadic composition of symmetric groups. For each fixed $\sigma \in \Sigma_n$ and $\tau_i \in \Sigma_{m_i}$, $1 \leq i \leq n$, we set 
		 $\tilde{\sigma} = \sigma \circ (\tau_1, \dots, \tau_n)$. Then we have:
		\begin{equation*}
			\beta^{\tilde{\sigma}} = \beta^\sigma \circ (\beta^{\tau_1} \otimes \cdots \otimes \beta^{\tau_n}). 
		\end{equation*} 
	\end{enumerate}
\end{lemma}

\begin{proof}
	It is clear that if such transformations exist, then they are uniquely characterized by the above properties. Moreover the above properties give a prescription for defining these transformations. Fix $n$ and let $\sigma_i \in \Sigma_n$ be the permutation of the $i^\textrm{th}$ and $(i+1)^\textrm{st}$ elements. The above properties require that $\beta^{\sigma_i}$ is given by the formula
	\begin{equation*}
		\beta^{\sigma_i} \cong \underbrace{id \otimes \cdots \otimes id}_{\textrm{first } (i-1) \textrm{ elements}} \otimes \beta_{x_i, x_{i+1}} \otimes \underbrace{id \otimes \cdots \otimes id}_{\textrm{last } (n-i-1) \textrm{ elements}}.
	\end{equation*}
If $\sigma = \sigma_{i_1} \cdots \sigma_{i_k}$ is a word in the standard presentation for $\Sigma_n$, then  $\beta^\sigma$ must be given as the composite $\beta^{\sigma_1} \circ \cdots \circ \beta^{\sigma_{i_k}}$. For this assignment to be well-defined we must verify that the relations of the symmetric group hold. 
	
The first relation holds because in a quasistrict symmetric monoidal 2-category the syllepsis transformation $u$ is the identity transformation. The second relation holds because in a 	quasistrict symmetric monoidal 2-category the interchanger coherence cell $\phi^\otimes$ is the identity when applied to identities and the components of the braiding transformations $\beta_{x_i, x_{i+1}}$. Finally the last relation follows from the fact that the transformations $R$ and $S$ are identities together with the fact that the `naturator' cell $\beta_{f,g}$ is an identity if either $f$ or $g$ is an identity.  

Thus the above prescription for the transformation $\beta_\sigma$ is well-defined, and satisfies the first two properties required by the lemma. The final property, compatibility with the operadic composition, also follows from the fact that the interchanger coherence cell $\phi^\otimes$ is the identity when applied to identities and the components of the braiding transformations $\beta_{x_i, x_{i+1}}$.
\end{proof}

Quasistrict symmetric monoidal 2-categories are extremely rigid, and as we have just seen, you get well defined operations for each element of the symmetric group.
We will now turn to our unbiased semistrict symmetric monoidal 2-categories, which weakens the axioms of quasistrict symmetric monoidal 2-categories, while keeping operations for each element of the symmetric group. In a Gray monoid $\sM$ the tensor product operation is strictly associative and gives rise to a well-defined operation 
	\begin{equation*}
		\otimes^n: \underbrace{\sM \times \cdots \times \sM}_{n \textrm{ times}} \to \sM.
	\end{equation*}
which is a cubical homomorphism in the sense of \cite{GPS95}. Note that by convention $\otimes^1: M \to M$ is the identity homomorphism and $\otimes^0: pt \to \sM$ is the unit element. 

\begin{definition}
	An {\em unbiased semistrict symmetric monoidal 2-category} consists of:
	\begin{itemize}
		\item a Gray monoid $\sM$;
		\item transformations $\beta_{\sigma}: \otimes^n \to \otimes^n \circ \sigma$ for each $\sigma \in \Sigma_n$, $n \geq 0$; 
		\item invertible modifications $X^{\sigma, \sigma'}: (\beta^\sigma  * 1) \circ \beta^{\sigma'} \to \beta^{\sigma \sigma'}$ for each pair of elements $\sigma, \sigma' \in \Sigma_n$; and for the identity element $e \in \Sigma_n$, we have an invertible modification $X^e: id \to \beta^e$;
	\end{itemize}
	These are required to satisfy the following three coherence conditions:
	\begin{enumerate}
		\item [(USS.1.a)] (Associativity)
		\begin{equation*}
			X^{\sigma \sigma', \sigma''} \circ (X^{\sigma, \sigma'} * 1)  = X^{\sigma, \sigma' \sigma''} \circ (1 * X^{\sigma', \sigma''})
		\end{equation*}
		\item [(USS.1.b)]  (Right Unitality)
		\begin{equation*}
			X^{\sigma, e} \circ (1 * X^e)  = id_{\beta^\sigma}
		\end{equation*}
		\item [(USS.1.c)] (Left Unitality)
		\begin{equation*}
			X^{e, \sigma} \circ (X^e * 1)  = id_{\beta^\sigma}
		\end{equation*}
	\end{enumerate}
	where $\sigma, \sigma', \sigma'' \in \Sigma_n$. Here the first equation is an identity between 2-morphisms from $\beta^{\sigma} \circ \beta^{\sigma'} \circ \beta^{\sigma''}$ to $\beta^{\sigma \sigma' \sigma''}$.
	
In addition these are required to satisfy the following four compatibilities with the operadic product: 	
		\begin{enumerate}

			\item [(USS.2)] We require that $\beta^{id \sqcup \sigma} = id \otimes \beta^\sigma$ and $\beta^{\sigma \sqcup id} = \beta^\sigma \otimes id$. 
			Moreover these identifications are to carry $X^{(id \sqcup \sigma),(id \sqcup \sigma') }$ to $id * X^{\sigma, \sigma'}$, and  $X^{(\sigma \sqcup id), (\sigma' \sqcup id)}$ to $X^{\sigma, \sigma'} * id$. 
			This implies that $\beta^e$ is tirvial in the case $n=0$  (i.e. the identity morphism of the unit in $\sM$). 
			
			\item [(USS.3)] For each pair $\sigma$ and $\sigma'$ we obtain two 2-morphisms from $\otimes^{n + n'}$ to $\otimes^{n + n'} \circ (\sigma \sqcup \sigma')$:
			\begin{align*}
				X^{id \sqcup \sigma', \sigma \sqcup id} :  \beta^\sigma \otimes \beta^{\sigma'} &= (id \otimes \beta^{\sigma'}) \circ (\beta^\sigma \otimes id) \\ & = \beta^{id \sqcup \sigma'} \circ \beta^{\sigma \sqcup id} \to \beta^{\sigma \sqcup \sigma'} \\
				X^{\sigma \sqcup id, id \sqcup \sigma'} \circ \phi^{\otimes}_{(\beta^{\sigma}, id), (id, \beta^{\sigma'})} : \beta^\sigma \otimes \beta^{\sigma'} &\to (\beta^\sigma \otimes id) \circ (id \otimes \beta^{\sigma'}) \\
				& = \beta^{\sigma \sqcup id} \circ  \beta^{id \sqcup \sigma'} \to \beta^{\sigma \sqcup \sigma'} 
			\end{align*}
			We require that they coincide. This implies that there is a canonical isomorphism $\otimes^n \beta^{\sigma_i} \cong \beta^{\sqcup_n \sigma_i}$ for each collection $\sigma_i \in \Sigma_{n_i}$, with $1 \leq i \leq n$. 
			
			\item [(USS.4)]
			 Fix $n$ and a collection of natural numbers $n_i \geq 0$ for $1 \leq i \leq n$. Set $N = \sum_i n_i$. Let $\sigma \in \Sigma_n$ and let $\tilde{\sigma} \in \Sigma_{N}$ be the element obtained from the operadic product $\sigma \circ (\tau_i)$ with each $\tau_i = e \in \Sigma_{n_i}$ trivial. Then we have: $\beta^{\tilde\sigma} = \beta^\sigma$, and moreover these identifications are to carry $X^e$ to $X^{\tilde e}$ and $X^{\sigma, \sigma'}$ to $X^{\tilde \sigma, \tilde \sigma'}$. 
			
			This means that for $(x_i)$ with $1 \leq i \leq N$ a collection of objects of $\sM$, we may set
				\begin{equation*}
					y_k = x_{\sum_{i=1}^{k-1} n_i + 1} \otimes x_{\sum_{i=1}^{k-1} n_i + 2} \otimes \cdots \otimes x_{\sum_{i=1}^{k} n_i}
				\end{equation*}	
				and then we have $\beta^{\tilde \sigma}_{(x_i)} = \beta^\sigma_{(y_k)}$. 
				
			This and the previous compatibility imply that if the $k^\textrm{th}$ element of $(x_i)$ is the unit object in $\sM$, then $\beta^\sigma = \beta^{\overline{\sigma}}$ where $\overline{\sigma} = d_k(\sigma) \in \Sigma_{n-1}$ is $\sigma$ with the $k^\textrm{th}$ letter deleted. 
			\item [(USS.5)] Fix $n$ and a collection of natural numbers $n_i \geq 0$ for $1 \leq i \leq n$. Set $N = \sum_i n_i$. Let $\sigma \in \Sigma_n$ and let $\tilde{\sigma} \in \Sigma_{N}$ be the element obtained from the operadic product $\sigma \circ (e)$, as in the previous axiom. In this case the $\tau_i$ may be non-trivial. We require that the coherence cell 
			\begin{equation*}
				\beta^\sigma_{(\beta^{\tau_i})}: (\otimes^n \beta^{\tau_{\sigma(i)}}) \circ \beta^\sigma \to \beta^\sigma \circ (\otimes^n \beta^{\tau_i})
			\end{equation*}
is given as a composite: 
			\begin{align*}
				(\otimes^n \beta^{\tau_{\sigma(i)}}) \circ \beta^\sigma & \cong \beta^{\sqcup_n \tau_{\sigma(i)}} \circ \beta^{\tilde{\sigma}} \\
				 &\cong \beta^{(\sqcup_n \tau_{\sigma(i)}) \circ \tilde{\sigma}} = \beta^{\tilde{\sigma} \circ \sqcup_n \tau_{i}} \\
				 & \cong \beta^{\tilde{\sigma}} \circ \beta^{\sqcup_n \tau_i} \\
				 & \cong \beta^\sigma \circ (\otimes^n \beta^{\tau_i})
			\end{align*}
			Where the first an last isomorphisms are given by our previous requirements and the middle two are induced from the modifications $X$. Moreover we further require that $\beta^\sigma_{(f_i)} = id$ if each $f_i = id_{x_i}$ is an identity. 
			\end{enumerate}
\end{definition}

\begin{remark}
	The notion of unbiased semistrict symmetric monoidal 2-category can be re-expressed in a more compact and natural form using the sophisticated language of Quillen model categories and enriched symmetric operads. There is a symmetric monoidal Quillen model category structure on the category of 2-categories, which was introduced by Lack \cite{MR1931220, MR2138540}. The monoidal structure, however, is not the Cartesian product (which is not compatible with the Quillen model structure) but the {\em Gray tensor product} \cite{MR0371990, GPS95}. This is an alternative symmetric monoidal closed structure which has many nice properties. It is intimately connected to cubical homomorphisms: for 2-categories $A$ and $B$ there exists a universal cubical homomorphism $A \times B \to A \otimes_G B$, where $A \otimes_G B$ is the Gray tensor product. Thus strict homomorphisms out of $A \otimes_G B$ are in natural bijection with cubical homomorphisms out of $A \times B$. A Gray monoid (Def~\ref{def:graymonoid}) is equivalent to a (strict) monoid with respect to this monoidal structure. 
	
	We may now consider the theory of symmetric operads (and their algebras) with respect to the symmetric monoidal category $\Gray= (2\cat, \otimes_G)$ of strict 2-categories with the Gray tensor product (this is made slightly easier by the fact that the terminal category is also the unit for the gray tensor product). For example the terminal symmetric operad is the operad $Comm$ whose algebras are strictly commutative monoids. Unbiased semistrict symmetric monoidal 2-categories are the algebras, not for the operad $Comm$, but for the operad $E_\infty$ which is a particular cofibrant replacement of $Comm$. Specifically the $n^\textrm{th}$ 2-categories $E_\infty(n)$ are cofibrant 2-categories (in the Lack model structure) which are equivalent to the terminal 2-category $pt$ but admit a free $\Sigma_n$-action. This is a 2-categorical analog of the topologists $E_\infty$-operad.
	
This can be seen directly from the definition of unbiased symmetric monoidal 2-category. The compatibility with the operadic product and the normalization conditions arise because these are {\em strict} algebras for the operad $E_\infty$. The weakness of the notion comes only from the cofibrant replacement in passing from $Comm$ to $E_\infty$. For example comparing the axioms that the modifications $X$ are required to satisfy to the notion of functor of bicategories, we see that these $X$ essentially amount to the data of certain weak functors from the pair groupoid on the set $\Sigma_n$ to certain functor categories, which are strictly $\Sigma_n$-equivariant.    

Needless to say, giving a complete and precise account along these lines would be quite involved and is somewhat orthogonal to our main task. Thus we invite the interested reader to pursue this point of view, but we will not pursue it any further here. This notion has also been considered by Gurski, Johnson, and Osorno \cite{GJO14}. 
\end{remark}

\subsection{A graphical calculus for unbiased semistrict symmetric monoidal 2-categories} \label{sec:graphcalculusUSS}
Now we will address how unbiased semistrict symmetric monoidal 2-categories relate to the other notions of symmetric monoidal bicategory we have seen.

\begin{lemma}
	Every quasistrict symmetric monoidal 2-categories gives rise to an unbiased semistrict symmetric monoidal 2-categories in which the modifications $X$ are the identity.
\end{lemma}

\begin{proof}
	 This is largely consequence of Lemma~\ref{lem:qsisuss}. In particular that lemma estabilishes that there exist transformations $\beta^\sigma$, as required in the definition of unbiased semistrict symmetric monoidal 2-category, and moreover that they are strictly compatible with the composition in the symmetric group and the operadic composition. Thus we may take the modifications $X^{\sigma, \sigma'}$ and $X^e$ to be identities. We must show that these prescriptions satisfy the axioms of a unbiased semistrict symmetric monoidal 2-category. Because these transformations are identities, axioms (USS.1) and (USS.2) are automatic. Moreover axiom (USS.3) also holds automaticallly because these modifcations are identities and because of axiom (QS.3) in the definition of quasistrict symmetric monoidal bicategory. Axiom (USS.4) holds because of Lemma~\ref{lem:qsisuss}, and the final axiom (USS.5) holds because of axiom (QS.2) in the definition of quasistrict symmetric monoidal 2-category, in particular because of Lemma~\ref{lma:QSbraidinggivenbyinterchange}.  
\end{proof}

\begin{figure}[htbp]
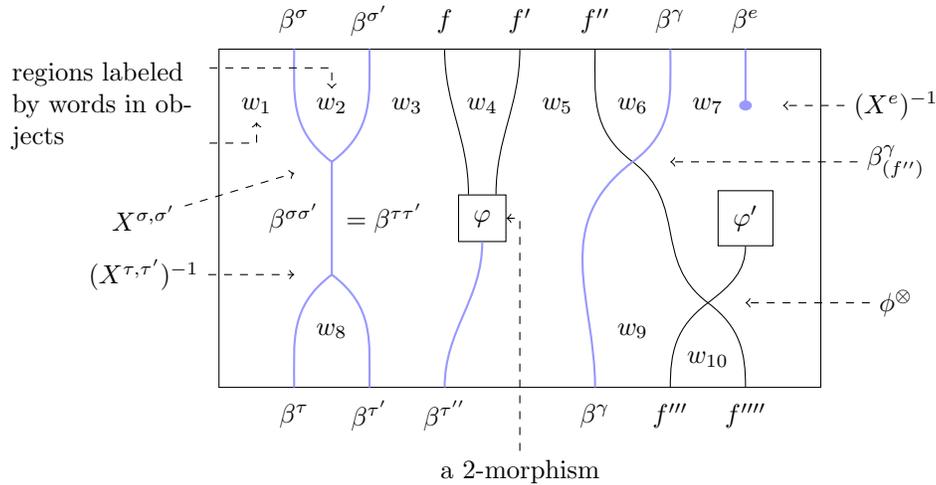

	\begin{center}

	\end{center}
	\caption{A Typical String Diagram in an Unbiased Semistrict Symmetric Monoidal 2-Category}
	\label{fig:USS-stringdiagram}
\end{figure}

Next we will show that every unbiased symmetric monoidal 2-category naturally inherits the structure of a Crans semistrict symmetric monoidal 2-category, however before we do it will be helpful to introduce a diagrammatic calculus for unbiased semistrict symmetric monoidal bicategories. There are two well-established diagrammatic calculi for 2-categories, the theory of pasting diagrams and the equivalent theory of string diagrams, both of which are reviewed in Appendix~\ref{PastingsStringsAdjointsAndMatesSection}. Unbiased semistrict symmetric monoidal 2-categories enjoy an enhanced version of the string diagram calculus. The typical sort of diagram on which we can evaluate is shown in Figure~\ref{fig:USS-stringdiagram}. \emph{We will always read these diagrams from right to left}. 

In a standard string diagram for a 2-category the regions are labeled with objects of the category, the arcs (or `strings') are labeled by the 1-morphisms and there are {\em coupons} (small `boxes') which are labeled by the 2-morphisms. The string diagrams for unbiased semistrict symmetric monoidal 2-categories are similar, except now all the regions are labeled by formal words in the objects. We will distinguish some of the arcs (colored chambering in Figure~\ref{fig:USS-stringdiagram}) which are labeled by the $\beta$ morphisms. We also distinguish several coupons: those coming from the modifications $X^{\sigma, \sigma'}$, $X^e$, their inverses, the coherence morphisms of the $\beta$ and $\phi^\otimes$ transformations. Some examples are displayed in the figure. 

As with ordinary string diagrams, the labels are subject to the obvious source-target matching conditions. Moreover since an unbiased symmetric monoidal string diagram gives a string diagram in the usual sense for the 2-category $\sM$, it has a well defined evaluation in $\sM$. 

However we have more operations with these string diagrams than with ordinary string diagrams. Given such a string diagram and an object $x$ of $\sM$ we may `tensor' these together to obtain a new string diagram. For each of the words we simply concatenate on the new object. We can do this on either side of the word and, indeed, there are left and right versions of this tensor operation. We will describe the operation of tensor by an object on the left. For the ordinary string-arcs and coupons we replace each label $f$ with $id_x \otimes f$ or $\varphi$ with $id_{id_x} * \varphi$. For the $\beta$-arcs we do the same, but then we may reinterpret this as $\beta^{id \sqcup \sigma}$, coming from a permutation from $\Sigma_{n+1}$. A similar reinterpretation happens for each of the distinguished kinds of coupons. The reason this is a sensible operation is due to axiom (USS.2) in the definition of unbiased semistrict symmetric monoidal 2-category. This axiom ensures that these identifications/reinterpretations are valid.

\begin{figure}[ht]
\begin{center}
% associative-like relation	
\begin{tikzpicture}[thick,scale = 0.5]
	\draw (0,0) rectangle (3,2);
	\draw (1,2) -- (1,1) 
		(0.5,0.75) rectangle (1.5,1) 
		(1,0.75) -- (1,0)
		(2,2) to [out=270, in=90] (2.5,0);
		
	\node at (4,1) {$\otimes$};

	\draw (5,0) rectangle (8,2);
	\draw (6,2) to [out=270, in=90] (5.5,0) 
		(6.75,2) to [out=270, in=90] (7,1.25)
		(6.5,1) rectangle (7.5,1.25)
		(6.75,1) to [out=270, in=90] (6.5,0)
		(7.25,1) to [out=270, in=90] (7.75,0);

	\draw (10,1.5) rectangle (16,3.5);
		\draw (11,3.5) to [out=270, in=90] (10.5,1.5) 
			(11.75,3.5) to [out=270, in=90] (12,2.75)
			(11.5,2.5) rectangle (12.5,2.75)
			(11.75,2.5) to [out=270, in=90] (11.5,1.5)
			(12.25,2.5) to [out=270, in=90] (12.75,1.5);
	\draw (10,-1.5) rectangle (16,0.5);
		\draw (14,0.5) -- (14,-0.5) 
			(13.5,-0.75) rectangle (14.5,-0.5) 
			(14,-0.75) -- (14,-1.5)
			(15,0.5) to [out=270, in=90] (15.5,-1.5);
	
	\draw (18,0) rectangle (24,2);
		\draw (19,2) to [out=270, in=90] (18.5,0) 
			(19.75,2) to [out=270, in=90] (20,1.25)
			(19.5,1) rectangle (20.5,1.25)
			(19.75,1) to [out=270, in=90] (19.5,0)
			(20.25,1) to [out=270, in=90] (20.75,0);
		\draw (22,2) -- (22,1) 
				(21.5,0.75) rectangle (22.5,1) 
				(22,0.75) -- (22,0)
				(23,2) to [out=270, in=90] (23.5,0);
	\draw [thin, ->] (1,-0.5) to [out=-90, in = 180] (9.5,-0.5);
	\draw [thin, ->] (6,2.5) to [out=90, in = 180] (9.5,2.5);
	\draw [thin] (16.5,2.5) to [out=0, in = 180] (17.5,1);
	\draw [thin] (16.5,-0.5) to [out=0, in = 180] (17.5,1);
	\draw [thin, ->] (17.5,1) -- (17.75,1);
	\node at (3,3) {`stretch'};
	\node [text width=3cm] at (20,3) {`merge' via tensor with objects};
\end{tikzpicture}
\caption{The monoidal operation for unbiased semistrict symmetric monoidal 2-categories in terms of string diagrams.}
\label{Fig:MonoidalStructureUSS-StringDiagram}
\end{center}
\end{figure}
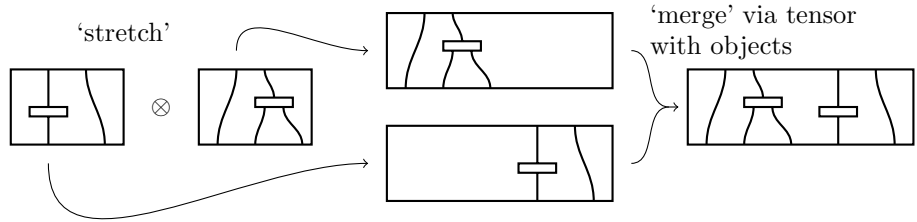

We can tensor 2-morphisms and this operation can be meaningfully expressed in terms of these enhanced string diagrams. Because the monoidal structure of an unbiased symmetric monoidal 2-category is cubical we have that
\begin{equation*}
	 f_1 \otimes f_2 = (id_{x} \circ f_1) \otimes (f_2 \circ id_{y}) = (id_{x} \otimes f_2) \circ (f_1 \otimes id_{y})
\end{equation*}
and similarly for 2-morphisms (the second equality holds only because the functor $\otimes$ is cubical). Thus we may rewrite the tensor operations in terms of the horizontal composition and the operation of tensoring with objects. In terms of string diagrams this process is depicted in Figure~\ref{Fig:MonoidalStructureUSS-StringDiagram}.

\begin{figure}[htbp]
	\begin{center}
		\begin{tikzpicture}[yscale = 0.75]
			\draw (0,0) rectangle (2,2) (4,0) rectangle (6,2) (8,0) rectangle (10,2);
			\draw (0,4) rectangle (4,7) (6,4) rectangle (10,7);
			
			\node at (3,1) {$=$}; \node at (7,1) {$=$}; \node at (5,5.5) {$=$};
			
			\draw [chambering] (1,7) to [out = -90, in = 135] (1.5, 6) 
				(2,7) to [out = -90, in = 45] (1.5,6)
				(1.5,6) to [out = -90, in = 135] (2,5)
				(3,7) to [out = -90, in = 45] (2,5)
				(2,5) -- (2,4);
				
			\node at (1,7.5) {$\beta^\sigma$};
			\node at (2,7.5) {$\beta^{\sigma'}$};
			\node at (3,7.5) {$\beta^{\sigma''}$};
			\node at (1,5.5) {$\beta^{\sigma \sigma'}$};
			\node at (2,3.5) {$\beta^{\sigma \sigma' \sigma''}$};

			\draw [chambering] (8,7) to [out = -90, in = 135] (8.5, 6) 
				(9,7) to [out = -90, in = 45] (8.5,6)
				(8.5,6) to [out = -90, in = 45] (8,5)
				(7,7) to [out = -90, in = 135] (8,5)
				(8,5) -- (8,4);
				
			\node at (7,7.5) {$\beta^\sigma$};
			\node at (8,7.5) {$\beta^{\sigma'}$};
			\node at (9,7.5) {$\beta^{\sigma''}$};
			\node at (9,5.5) {$\beta^{\sigma' \sigma''}$};
			\node at (8,3.5) {$\beta^{\sigma \sigma' \sigma''}$};	
				
			\draw [chambering] (0.5, 2) to [out = -90, in = 135] (1,1)
				 (1.5,1.75) to [out = -90, in = 45] (1,1)
				(1,1) -- (1,0);
			\draw [chambering,fill] (1.5, 1.75) circle (2pt);
			
			\draw [chambering] (5,2) -- (5,0);
			
			\draw [chambering, fill] (8.5,1.75) circle (2pt);
			\draw [chambering]  (8.5,1.75) to [out = -90, in = 135] (9,1)
			(9.5, 2) to [out = -90, in = 45] (9,1)
			(9,1) -- (9,0);
			
			\node at (0.5,2.5) {$\beta^\sigma$};
			\node at (1,-0.5) {$\beta^\sigma$};
			
			\node at (5,2.5) {$\beta^\sigma$};
			
			\node at (9.5,2.5) {$\beta^\sigma$};
			\node at (9,-0.5) {$\beta^\sigma$};
		\end{tikzpicture}
	\end{center}
	\caption{The axioms (USS.1) of unbiased semistrict symmetric monoidal 2-categories.}
	\label{fig:graphicaxiomUSS1}		
\end{figure}
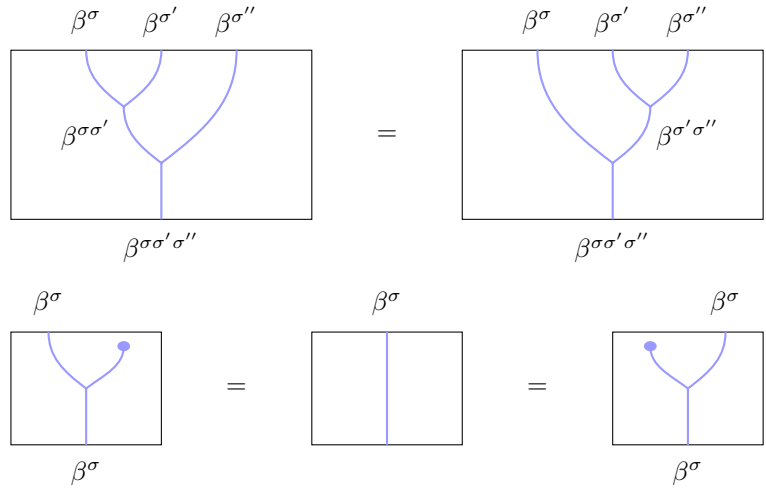

\begin{figure}[htpb]
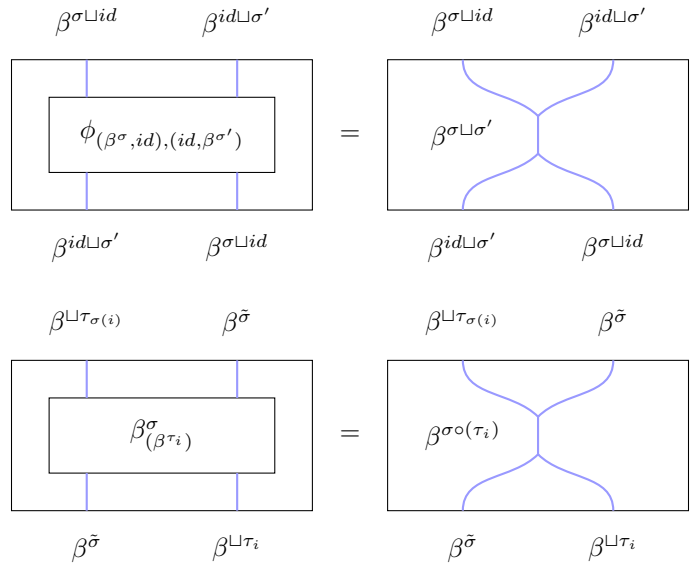


	\caption{The axioms (USS.3) and (USS.5) of unbiased semistrict symmetric monoidal 2-categories.}
	\label{fig:graphicaxiomUSS3+5}		
\end{figure}

The axioms of an unbiased semistrict symmetric monoidal 2-category have a clean interpretation in terms of these string diagrams. The axioms (USS.1) are depicted in Figure~\ref{fig:graphicaxiomUSS1}. Of course the inverse relations, which correspond to flipping the string diagrams along the $x$-axis, also hold. Axioms (USS.2) and (USS.4) are used in reinterpreting the $\beta$ and $X$ morphisms as coming from different symmetric groups, such as in the operation of tensoring with an object, describe above. Axioms (USS.3) and (USS.5) are shown in Figure~\ref{fig:graphicaxiomUSS3+5}.

\begin{figure}[htbp]
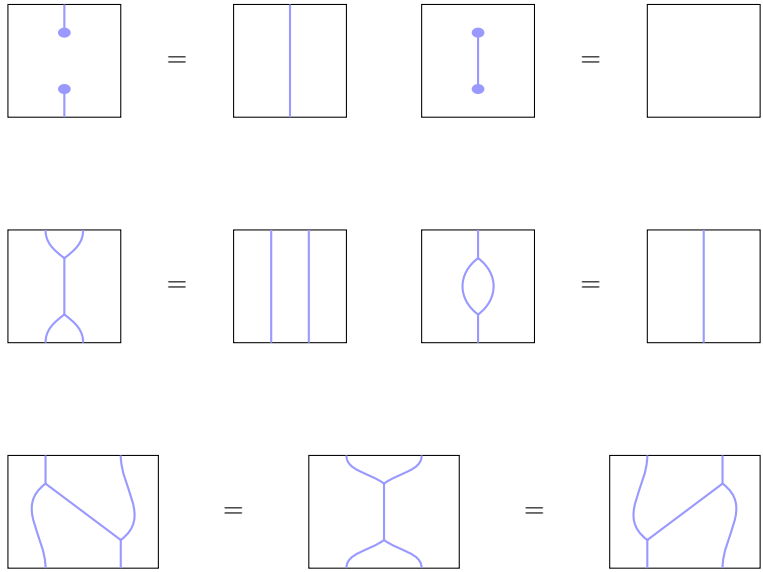

	\begin{center}
		%
	\end{center}
	\caption{Relations for string diagrams for unbiased semistrict symmetric monoidal 2-categories.}
	\label{fig:graphicUSSrelations}	
\end{figure}

\begin{lemma}\label{lem:basicRelationsforUSSstringDiagrams}
	Each of the relations depicted in Figure~\ref{fig:graphicUSSrelations} holds for string diagrams for unbiased semistrict symmetric monoidal bicategories, whenever consistent labelings of each side of the equations is possible. 
\end{lemma}

\begin{figure}[htpb]
	\begin{tikzpicture}
		\draw (0,0) rectangle (2,4) (3,0) rectangle (5,4) (6,0) rectangle (8,4) (9,0) rectangle (11,4);
		\draw [chambering] (0.5,4) -- (0.5,3.5) to [out = -135, in = 90] (0.5, 0)
			(0.5,3.5) to [out = -45, in = 135] (1.5, 2.5) -- (1.5, 0)
			(1.5,2.5) to [out = 45, in = -90] (1.5,4);
		
		\draw [chambering] (3.5,4) -- (3.5,3.5) to [out = -135, in = 135] 
			(4, 1.5) -- (4,1) to [out = -135, in = 90] (3.5,0)
			(3.5,3.5) to [out = -45, in = 135] (4.5, 2.5) to [out = -90, in = 45] (4,1.5)
			(4.5,2.5) to [out = 45, in = -90] (4.5,4)
			(4,1) to [out = -45, in = 90] (4.5,0);
		
		\draw [chambering] (6.5,4) -- (6.5,3.5) to [out = -135, in = 135] (6.5,2.5)
			(6.5,3.5) to [out = -45, in = 45] (6.5, 2.5) to [out = -90, in = 135] (7,1.5)
			(7, 1.5) -- (7,1) to [out = -135, in = 90] (6.5,0) 
			(7,1.5) to [out = 45, in = -90] (7.5,4)
			(7,1) to [out = -45, in = 90] (7.5,0);	
		
		\draw [chambering] (9.5,4) to [out = -90, in = 135] (10,1.5)
				(10, 1.5) -- (10,1) to [out = -135, in = 90] (9.5,0) 
				(10,1.5) to [out = 45, in = -90] (10.5,4)
				(10,1) to [out = -45, in = 90] (10.5,0);
		
		\node at (2.5,2) {$=$};
		\node at (5.5,2) {$=$};
		\node at (8.5,2) {$=$};
			
	\end{tikzpicture}
	\caption{Proof of `Frobenius' relation.}
	\label{fig:USSFrobRelnproof}	
\end{figure}
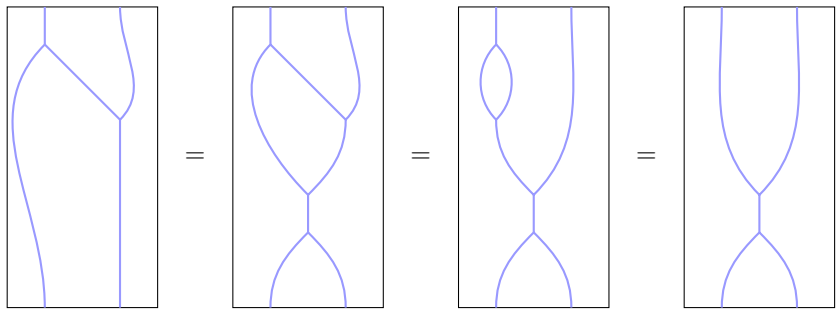

\begin{proof} The first four equations from Figure~\ref{fig:graphicUSSrelations} simply express the invertibility of the $X^{\sigma, \sigma'}$ and $X^e$ 2-morphisms. The last two equalities are consequences of these and the relations in Figure~\ref{fig:graphicaxiomUSS1}. For example the first of these two equalities can be seen by the sequence of moves depicted in Figure~\ref{fig:USSFrobRelnproof}.  In the first and last equality of this figure we are using the invertibility of the $X^{\sigma, \sigma'}$s, and in the middle equality we are using (the inverse of) the `associativity' relation from Figure~\ref{fig:graphicaxiomUSS1}. Figure~\ref{fig:USSFrobRelnproof} is unlabeled, but there is a unique way to label it consistent with the starting configuration and the operations. The final equality of the lemma is similar. 	
\end{proof}

\begin{lemma}
	An unbiased symmetric monoidal 2-category naturally inherits the structure of a Crans semistrict symmetric monoidal 2-category.  
\end{lemma}

\begin{proof}
	An unbiased semistrict symmetric monoidal 2-category has an underlying Gray monoid, just as a Crans semistrict symmetric monoidal 2-category. We need to provide the braiding transofrmation $\beta$, the modifications $R$, $S$, and $\sigma$ and demonstrate that they satisfy the axioms of a symmetric monoidal bicategory and the normalization axioms of a Crans semistrict symmetric monoidal 2-category. Each of these morphisms has a clear candidate coming from the structure of a unbiased semistrict symmetric monoidal 2-category:
	\begin{itemize}
		\item The braiding is given by the transformation $\beta^\sigma$, where $\sigma \in \Sigma_2$ is the non-trivial element;
		\item The modifications $R$ and $S$ are given by transformations $(X^{\sigma, \sigma'})^{-1}$, with $\sigma$ and $\sigma'$ once stabilized braidings;
		\item The syllepsis modification $\sigma$ is given by a composite $(X^{\tau, \tau})^{-1} \circ X^e$, where $\tau \in \Sigma_2$ is the non-trivial element.
	\end{itemize}
With these choices the Crans normalization (CSS.2a) follows directly from (USS.2) and the Crans normalizations (CSS.2b) and (CSS.2c) follow from (USS.4). 
We must now verify the seven remaining nontrivial axioms from the definition of symmetric monoidal bicategory: Axioms (SM3.i), (SM3.ii), (SM4), (SM5), (SM6.i), (SM6.ii), and (SM7). 

The axioms (SM3.i) and (SM3.ii) are both special cases of (USS.1a), the `associativity' relation. 
Figure~\ref{fig:USS-AxiomSM4} shows, in terms of the string diagrams valid for unbiased semistrict symmetric monoidal 2-categories,  both sides of the equation stating axiom (SM4) (we leave it to the reader to add the correct labels to this diagram, which is easy to read off from the pasting diagram (SM4)). 
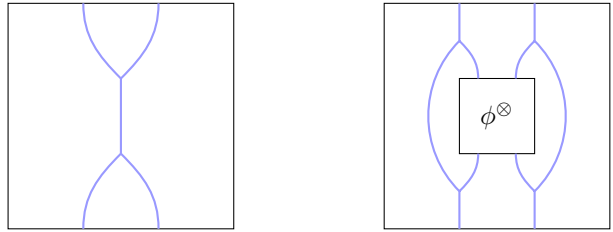
\begin{figure}[htpb]
	\begin{tikzpicture}
		\draw (0,0) rectangle (3,3) (5,0) rectangle (8,3);
		\draw (6,1) rectangle (7,2);
		\node at (6.5,1.5) {$\phi^\otimes$};
		%\node at (4,1.5) {$=$};
		\draw [chambering] (1,3) to [out = -90, in = 135] (1.5, 2) -- (1.5, 1) to [out = -135, in = 90] (1,0)
			(2,3) to [out = -90, in = 45] (1.5, 2) (1.5,1) to [out = -45, in = 90] (2,0);
		\draw [chambering] (6,3) -- (6,2.5) to [out = -135, in = 135] (6, 0.5) -- (6, 0)
			(6,2.5) to [out = -45, in = 90] (6.25, 2) (6.25, 1) to [out = -90, in = 45] (6, 0.5)
			(7,3) -- (7,2.5) to [out = -45, in = 45] (7, 0.5) -- (7,0)
			(7, 2.5) to [out = -135, in = 90] (6.75, 2) (6.75, 1) to [out = -90, in = 135] (7,0.5);
	\end{tikzpicture}
	\caption{Unbiased semistrict string diagrams for axiom (SM4).}
	\label{fig:USS-AxiomSM4}
\end{figure}
Notice that the right-hand diagram involves a non-trivial interchanger $\phi^\otimes$ coming from the cubical functor $\otimes$. Axiom (USS.3) allows us to replace this interchanger with a composite of $X$ 2-morphisms (see also Figure~\ref{fig:graphicaxiomUSS3+5}). After we do this replacement, it is straightforward to change one side of the (SM4) equation into the other using the relations from (USS.1) and Lemma~\ref{lem:basicRelationsforUSSstringDiagrams} (see also Figures~\ref{fig:graphicaxiomUSS1} and \ref{fig:graphicUSSrelations}). Axiom (SM5) (the `Breen polytope') is proven in exactly the same manor, except now each side of the equation involves a non-trivial instance of the braiding coherence cell $\beta^\sigma_{(\beta^{\tau_i})}$. These must first be replaced with $X$ 2-morphisms, which permitted due to Axiom (USS.5), see Figure~\ref{fig:graphicaxiomUSS3+5}. 

Now we consider the sylleptic monoidal axioms, which connect the $R$ and $S$ modifications. Axioms (SM6.i) and (SM6.ii) each take the form of the first equality depicted in Figure~\ref{fig:USS-SyllepticAxioms} in terms of unbiased semistrict string diagrams. They differ slightly in the precise labelings of the chambering string arcs (which we have suppressed). Axiom (SM7) is given in the same figure as the second equality. All are easily verified using the relations we have established for such diagrams, see Figures~\ref{fig:graphicaxiomUSS1} and \ref{fig:graphicUSSrelations}. This completes our verifications of the axioms of a Crans semistrict symmetric monoidal 2-category. 
\begin{figure}[htbp]
	\begin{tikzpicture}
		\draw (0,0) rectangle (2,5) (3,0) rectangle (7,5);
		
		\draw [chambering] (1,5) -- (1,2.5) to [out = -135, in = 90] (0.5,0) (1,2.5) to [out = -45, in = 90] (1.5, 0);
		
		\draw [chambering] (3.5,5) to [out = -90, in = 135] (4,1) -- (4,0.5) 
		 	(4,1) to [out = 45, in = -90] (5,2) to [out = 45, in = -135] (5.5,3) -- (5.5,3.5)
		 	(5.5,3) to [out = -45, in = 90] (6,0)
		 	(5,2) to [out = 135, in = -135] (5.5, 4) -- (5.5,4.5)
			(5.5,4) to [out = -45, in= 90] (6.5, 0);
		\draw [chambering,fill] (4,0.5) circle (2pt) (5.5,3.5) circle (2pt) (5.5, 4.5) circle (2pt);
		\node at (2.5, 2.5) {$=$};
		
		\draw (0,-4) rectangle (2,-1) (3,-4) rectangle (7,-1);
		\draw [chambering] (1,-1) -- (1,-4);
		\draw [chambering] (4, -1) to [out = -90, in = 135] (4.5, -3) -- (4.5, -3.5)
			(4.5, -3) to [out = 45, in = -135] (5.5, -2) -- (5.5, -1.5)
			(5.5, -2) to [out = -45, in = 90] (6,-4);
		\draw [chambering, fill] (4.5, -3.5) circle (2pt) (5.5, -1.5) circle (2pt);
		\node at (2.5, -2.5) {$=$};
		
	\end{tikzpicture}
	\caption{Axioms (SM6.i) (or (SM6.ii)) and (SM7) in string diagrams.}
	\label{fig:USS-SyllepticAxioms}
\end{figure}
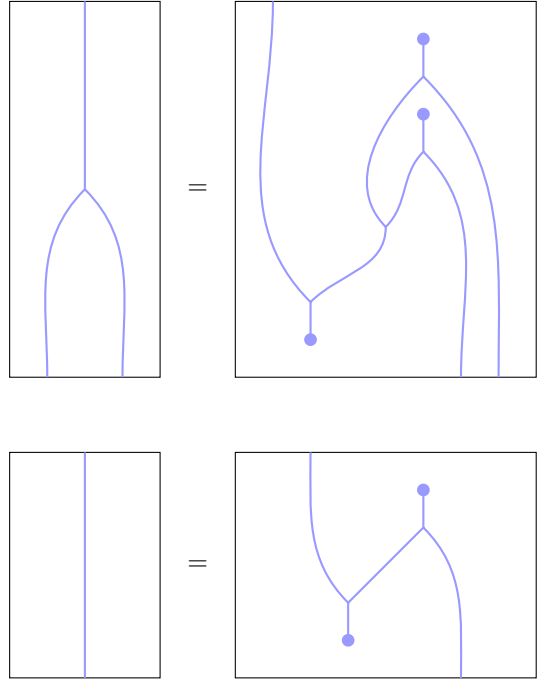
\end{proof}

\section{Generalized Computads} \label{sec:computads}

Here we will describe a very general method of defining and constructing presentations of higher categorical objects. These presentations will also be called {\em computads}, following the literature. Our treatment shares many similarities with Batanin's construction \cite{MR1664991}, but differs in a few small details as we will explain below. The theory is very general and best expressed with the language of monads.

	A {\em monad} $(T, \mu, \eta)$ on a category $\sC$ is a monoid in the monoidal category of endofunctors $\End(\sC)$. A morphism of monads is a monoid homomorphism (though many authors use a more general notion of morphism). A monad is {\em finitary} if the underlying endofunctor $T$ commutes with filtered colimits. 
	
	An {\em algebra} for the monad $T$, also called a {\em $T$-algebra}, consists of an object $X \in \sC$ and an {\em action map} $\alpha: TX \to X$ such that $X \stackrel{\eta_X}{\to} TX \stackrel{\alpha}{\to} X$ is the identity and the square:
	\begin{center}
	\begin{tikzpicture}
			\node (LT) at (0, 1.5) {$T^2X$};
			\node (LB) at (0, 0) {$TX$};
			\node (RT) at (2, 1.5) {$TX$};
			\node (RB) at (2, 0) {$X$};
			\draw [->] (LT) -- node [left] {$T\alpha$} (LB);
			\draw [->] (LT) -- node [above] {$\mu_X$} (RT);
			\draw [->] (RT) -- node [right] {$\alpha$} (RB);
			\draw [->] (LB) -- node [below] {$\alpha$} (RB);
			%\node at (0.5, 1) {$\ulcorner$};
			%\node at (1.5, 0.5) {$\lrcorner$};
	\end{tikzpicture}
	\end{center}
	commutes. 
	
	A {\em morphism} of $T$-algebras $f: (X, \alpha^X) \to (Y, \alpha^Y)$ is a morphism $f: X \to Y$ in $\sC$ such that the following square commutes:
	\begin{center}
	\begin{tikzpicture}
			\node (LT) at (0, 1.5) {$TX$};
			\node (LB) at (0, 0) {$TY$};
			\node (RT) at (2, 1.5) {$X$};
			\node (RB) at (2, 0) {$Y$};
			\draw [->] (LT) -- node [left] {$Tf$} (LB);
			\draw [->] (LT) -- node [above] {$\alpha^X$} (RT);
			\draw [->] (RT) -- node [right] {$f$} (RB);
			\draw [->] (LB) -- node [below] {$\alpha^Y$} (RB);
			%\node at (0.5, 1) {$\ulcorner$};
			%\node at (1.5, 0.5) {$\lrcorner$};
	\end{tikzpicture}
	\end{center}

	 The $T$-algebras and their morphisms form a category $\alg_T$. The forgetful functor from $T$-algebras to $\sC$ is part of an adjunction
	 \begin{equation*}
	 	F^T: \sC \leftrightarrows \alg_T: U.
	 \end{equation*}
	 There is an isomorphism of monads $UF^T\cong T$ and so the monad $T$ may be recovered from the adjunction. We will say that an adjunction equivalent to one of this form is a {\em monadic adjunction}.

In many approaches to higher categories the basic underlying structure consists of a collection of {\em cells}. These are the objects, 1-morphisms, 2-morphisms, etc of the higher categorical object. A given higher categorical structure will also include various maps (such as compositions) between these cells, but in good circumstances these can be regarded as placing additional algebraic structure on the underlying collection of cells. Such collections are called {\em globular sets}, which we define precisely below. For these `algebraic' notions of higher categorical object, they will be the algebras for a monad on globular sets.

Let $\G_n$ denote the category of {\em globes}, which is the category generated by the graph:
\begin{center}
	\begin{tikzpicture}
		\node (A) at (0,0) {$0$};
		\node (B) at (1,0) {$1$};
		\node (C) at (2,0) {$2$};
		\node (D) at (3,0) {$\cdots$};
		\node (E) at (4,0) {$n$};
			\draw [transform canvas={yshift=0.5ex},->]  (A) -- node [above] {s} (B);
			\draw [transform canvas={yshift=-0.5ex},->] (A) -- node [below] {t} (B);
			\draw [transform canvas={yshift=0.5ex},->]  (B) -- node [above] {s} (C);
			\draw [transform canvas={yshift=-0.5ex},->] (B) -- node [below] {t} (C);
			\draw [transform canvas={yshift=0.5ex},->]  (C) -- node [above] {s} (D);
			\draw [transform canvas={yshift=-0.5ex},->] (C) -- node [below] {t} (D);
			\draw [transform canvas={yshift=0.5ex},->]  (D) -- node [above] {s} (E);
			\draw [transform canvas={yshift=-0.5ex},->] (D) -- node [below] {t} (E);
	\end{tikzpicture}
\end{center}
and subject to the relations $ss=st$ and $tt = ts$. The category of presheaves on $\G_n$, denoted $\gset_n := \fun(\G_n^\textrm{op}, \set)$, is the category of {\em (non-reflexive) globular sets} (of dimension $n$). Thus a globular set consists exactly of the data of a collection of {\em cells}: 0-cells, 1-cells, 2-cells, up through $n$-cells, and each cell of positive dimension has a source and a target, which are cells of one lower dimension. 

\begin{example}
	The {\em $k$-cell} $C_k$ is the representable gobular set represented by the object $k$. This is not to be confused with the cells of a general globular set $X$; these are the maps $C_k \to X$. The {\em boundary of the $k$-cell} $\partial C_k$ is defined inductively: $\partial C_0 = \emptyset$, and $\partial C_k = C_{k-1} \cup_{\partial C_{k-1}} C_{k-1}$. There is a canonical inclusion $\partial C_k \to C_k$ given by the union $s \cup t$. This map is the one used to define the pushout, which completes the induction.  Note that $\partial C_{n+1}$ is a well-defined $n$-dimensional globular set. The map $id \cup id: \partial C_{k+1} = C_k \cup_{\partial C_{k-1}} C_k \to C_k$ will be called the {\em collapse map}.  
\end{example}

Many types of categorical objects can be described as algebras for a monad on globular set. These include the notions of strict $n$-categories, bicategories, tricategories\footnote{This is not true for the original definition of tricategory \cite{GPS95}, but is true for the essentially equivalent modification provided in \cite{MR2717302}.}, symmetric monoidal bicategories, etc. In these cases the $T$-algebra morphisms give a notion of strict homomorphism making this a category. 

For example there is the strict $n$-category monad, whose algebras are globular sets equipped with the structure of a strict $n$-category, $\alg_T = n\cat_s$; there is the bicategory monad on $\gset_2$, whose algebras are bicategories, $\alg_T = \bicat_s$; and, as we will see, there are monads on $\gset_2$ whose algebras are symmetric monoidal bicategories, and monads whose algebras are the more strict versions of symmetric monoidal bicategories we encountered in Section~\ref{sec:strictsymbicats}.

If the monad is {\em finitary} (meaning the endofunctor on globular sets commutes with filtered colimits) then the category of algebras is cocomplete, see \cite[Thm.~9.3.9]{MR2178101}. % This thm states that Alg_T has coequalizers, but since it also has filteredt colimits (and U creates them) it is cocomplete.
In that case we can define the notion of {\em computad} which gives a way of presenting the given categorical objects. Let $T$ be a finitary monad on the category $\gset_n$ of $n$-dimensional globular sets. We will inductively define the category $\comp_k^T$ of {\em $k$-computads} with respect to the monad $T$. This category will come equipped with an adjunction
\begin{equation*}
	\sF_k:\comp_k^T \leftrightarrows \alg_T: V_k
\end{equation*} 
which takes a $k$-computad to the $T$-algebra it presents. For $k=-1$ we set $\comp^T_k = pt$ and $\sF_{-1}(pt)$ is the initial $T$-algebra. The objects of the category $\comp_k^T$ consist of triples $(C, X, x)$ where $C \in \comp_{k-1}^T$ is an $(k-1)$-computad, $X \in \set$ is a set, and 
\begin{equation*}
	x: X \times \partial C_{k} \to U(\sF_{k-1}C)
\end{equation*}
a map of globular sets. Here $U: \alg_T \to \gset_n$ is the forgetful functor from $T$-algebras to globular sets. The map $x$ is equivalent to a map of $T$-algebras
\begin{equation*}
	x: T(X \times \partial C_{k}) \to \sF_{k-1}(C)
\end{equation*}
out of the free $T$-algebra generated by the globular set $X \times \partial C_{k}$. A morphism of $k$-computads is a pair morphisms $C \to C'$, $X \to X'$, such that the following square commutes:
\begin{center}
\begin{tikzpicture}
	\node (LT) at (0, 1.5) {$X\times \partial C_{k}$};
	\node (LB) at (0, 0) {$U(\sF_{k-1}C)$};
	\node (RT) at (3, 1.5) {$X'\times \partial C_{k}$};
	\node (RB) at (3, 0) {$U(\sF_{k-1}C')$};
	\draw [->] (LT) -- node [left] {$$} (LB);
	\draw [->] (LT) -- node [above] {$$} (RT);
	\draw [->] (RT) -- node [right] {$$} (RB);
	\draw [->] (LB) -- node [below] {$$} (RB);
	%\node at (0.5, 1) {$\ulcorner$};
	%\node at (1.5, 0.5) {$\lrcorner$};
\end{tikzpicture}.
\end{center}

To complete the induction we must also define the adjunction $\sF_k \dashv V_k$. The functor $V_k$ is given inductively, for a $T$-algebra $\sC \in \alg_T$, by 
\begin{equation*}
	V_k(\sC) = (V_{k-1}\sC, X_{(\sC,k)}, x_{(\sC,k)})
\end{equation*}
where $X_{(\sC,k)}$ and $x_{(\sC,k)}$ are defined via the pullback square:
\begin{center}
\begin{tikzpicture}
	\node (LT) at (0, 1.5) {$X_{(\sC,k)}$};
	\node (LB) at (0, 0) {$\gset(\partial C_k, U\sF_{k-1}V_{k-1}\sC)$};
	\node (RT) at (4, 1.5) {$\gset(C_k, U\sC)$};
	\node (RB) at (4, 0) {$\gset(\partial C_k, U\sC)$};
	\draw [->] (LT) -- node [left] {$x_{(\cC,k)}$} (LB);
	\draw [->] (LT) -- node [above] {$$} (RT);
	\draw [->] (RT) -- node [right] {$$} (RB);
	\draw [->] (LB) -- node [below] {$$} (RB);
	\node at (0.5, 1) {$\ulcorner$};
	%\node at (1.5, 0.5) {$\lrcorner$};
\end{tikzpicture}.
\end{center}
The functor $\sF_k$ is defined on the $k$-computad $D=(C,X,x)$ as the pushout of $T$-algebras:
\begin{center}
\begin{tikzpicture}
	\node (LT) at (0, 1.5) {$T(X \times \partial C_k)$};
	\node (LB) at (0, 0) {$\sF_{k-1}(C)$};
	\node (RT) at (3, 1.5) {$T(X \times C_k)$};
	\node (RB) at (3, 0) {$\sF_k(D)$};
	\draw [->] (LT) -- node [left] {$$} (LB);
	\draw [->] (LT) -- node [above] {$$} (RT);
	\draw [->] (RT) -- node [right] {$$} (RB);
	\draw [->] (LB) -- node [below] {$$} (RB);
	%\node at (0.5, 1) {$\ulcorner$};
	\node at (2.5, 0.5) {$\lrcorner$};
\end{tikzpicture}.
\end{center}
Here we use the fact that the monad $T$ is finitary to ensure that the above pushout exists. By convention, we will extend the above construction to $k=n+1$ by replacing the inclusion $\partial C_{n+1} \to C_{n+1}$ (which is not a map of $n$-dimensional globular sets) with the collapse map $\partial C_{n+1} \to C_n$. We will call the $T$-algebras in the image of $\sF_{k+1}$ as {\em computadic $T$-algebras}.

\begin{proposition}\label{pro:computadsexist}
	Let $T$ be a finitary monad on the category of $n$-dimensional globular sets. Then the above inductive construction is well-defined and gives rise to an adjunction 
	\begin{equation*}
		\sF_k:\comp_k^T \leftrightarrows \alg_T: V_k
	\end{equation*}
	for all $0 \leq k \leq n+1$. 
\end{proposition}

\begin{proof} All that remains to complete the induction is to show that $\sF_k$ and $V_k$ are adjoint. By construction this holds for $k=-1$. We will focus on the case $k\leq n$, the case $k = n+1$ is similar. 
The counit map $\varepsilon_k:\sF_k V_k (\sC) \to \sC$ may be written down directly. It is induced by the natural commutative diagram of $T$-algebras:
	\begin{center}
		\begin{tikzpicture}
			\node (LT) at (0, 1.5) {$T(X_{(\sC,n)} \times \partial C_k)$};
			\node (LB) at (0, 0) {$\sF_{k-1}V_{k-1} \sC$};
			\node (RT) at (4, 1.5) {$T(X_{(\sC,n)} \times C_k)$};
			\node (RB) at (4, 0) {$\sF_k V_k \sC$};
			\node (X) at (8, -1.5) {$\sC$};
			\draw [->] (LT) -- node [left] {$$} (LB);
			\draw [->] (LT) -- node [above] {$$} (RT);
			\draw [->] (RT) -- node [right] {$$} (RB);
			\draw [->] (LB) -- node [below] {$$} (RB);
			\draw [->] (RT) to [bend left] node [above right] {$$} (X);
			\draw [->] (LB) to [bend right] node [below left] {$\varepsilon_{k-1}$} (X);
			\draw [->, dashed] (RB) -- node [above right] {$\varepsilon_k$} (X);
			%\node at (0.5, 1) {$\ulcorner$};
			\node at (3.5, 0.5) {$\lrcorner$};
		\end{tikzpicture}
	\end{center}
where the morphism $T(X_{(\sC,n)} \times C_k) \to \sC$ is the mate of natural map $X_{(\sC,n)} \to \gset_n(C_k, U \sC)$.

The map $\varepsilon_k$ is the counit of an adjunction between $\sF_k$ and $V_k$ precisely if the composite
\begin{equation*}
	\comp_k^T( D, V_k(\sC)) \stackrel{\sF_k}{\to} \alg_T( \sF_k(D), \sF_k V_k(\sC)) \stackrel{(\varepsilon_k)_*}{\to} \alg_T( \sF_k(D), \sC)
\end{equation*}
 is a bijection. A direct calculation shows that this is indeed the case. 
\end{proof}

A morphism of monads $f: T \to T'$ induces a functor $f^*: \alg_{T'} \to \alg_T$ such that $U^{T} f^* = U^{T'}$. In the context above we have:

\begin{lemma}\label{lem:morphismofCompatads}
	Let $f: T \to T'$ be a morphism between finitary monads on the category of $n$-dimensional globular sets. Then for each $0 \leq k \leq n+1$ there is an induced functor on categories of computads $f_*: \comp^T_k \to \comp^{T'}_k$. Moreover we have an adjunction between the composites: $\sF^{T'}_k \circ f_* \dashv V^T_k \circ f^*$, and this induces a canonical homomorphism of $T$-algebras 
	\begin{equation*}
		\phi: \sF^T_k(D) \to f^* \sF^{T'}_k(f_* D) 
	\end{equation*}
	between the $T$-algebra generated by the $k$-computad $D$ and the $T$-algebra underlying the $T'$-algebra generated by $f_*(D)$.
\end{lemma}

\begin{proof}
	This is again constructed inductively. The base case $k=-1$ has $f_*$ the identity functor. To prove the lemma we need to construct the functor $f_*$ and demonstrate the adjunction. This later is determined by the unit of the adjunction
	\begin{equation*}
		id \to V^T_k \circ f^* \circ \sF^{T'}_k \circ f_*,
	\end{equation*}
which itself is determined by its mate $\phi: \sF_k^T \to f^* \circ \sF^{T'}_k \circ f_*$ (see Appendix~\ref{PastingsStringsAdjointsAndMatesSection} for a review of the theory of mates). Thus we first aim to construct $f_* = f_*^{(k)}$ and $\phi = \phi^{(k)}$, assuming we have done so for lower $k$. 
	
Let $D = (C, X,x) \in \comp^T_k$ be a $k$-computad with respect to the monad $T$. The natural transformation $\phi^{(k-1)}: \sF_{k-1}^T \to f^* \circ \sF^{T'}_{k-1} \circ f^{(k-1)}_*$ induces
\begin{align*}
	X \stackrel{x}{\to} \gset_n(\partial C_k, U \sF_{k-1}^T C) \stackrel{\phi}{\to} & \gset_n(\partial C_k, U f^*\sF^{T'}_{k-1} f_*^{(k-1)} C ) \\ & = \gset_n(\partial C_k, U \sF^{T'}_{k-1} f_*^{(k-1)} C ). 
\end{align*}
Thus we may define $f_*^{(k)} D := ( f_*^{(k-1)}C, X, \phi \circ x)$, where `$\phi \circ x$' denotes the above composite. 	

The natural transformation $\phi^{(k)}$ is induced from the following natural commutative diagram of $T$-algebras:	
\begin{center}
	\begin{tikzpicture}
		\node (LT) at (0, 1.5) {$T(X \times \partial C_k)$};
		\node (LB) at (0, 0) {$\sF_{k-1}^T C$};
		\node (RT) at (4, 1.5) {$T(X \times C_k)$};
		\node (RB) at (4,0) {$\sF^T D$};
		
		\node (XLT) at (2.5, -1) {$f^*T'(X \times \partial C_k)$};
		\node (XLB) at (2.5, -2.5) {$f^*\sF_{k-1}^{T'} f_*^{(k-1)} C$};
		\node (XRT) at (6.5, -1) {$f^*T'(X \times C_k)$};
		\node (XRB) at (6.5, -2.5) {$f^* \sF_k^{T'} f_* D$};
	
		\draw [->] (LT) -- node [left] {$$} (LB);
		\draw [->] (LT) -- node [above] {$$} (RT);
		\draw [->] (RT) -- node [right] {$$} (RB);
		\draw [->] (LB) -- node [below] {$$} (RB);
		
		\draw [->] (XLT) -- node [left] {$$} (XLB);
		\draw [->] (XLT) -- node [above] {$$} (XRT);
		\draw [->] (XRT) -- node [right] {$$} (XRB);
		\draw [->] (XLB) -- node [below] {$$} (XRB);
		
		\draw [->] (RT) to node [above right] {$f$} (XRT);
		\draw [->] (LT) to  node [above right] {$f$} (XLT);
		\draw [->] (LB) to  node [below left] {$\phi^{(k-1)}$} (XLB);
		\draw [->, dashed] (RB) to  node [below left] {$\phi^{(k)}$} (XRB);	
		%\node at (0.5, 1) {$\ulcorner$};
		\node at (3.5, 0.5) {$\lrcorner$};
	\end{tikzpicture}
\end{center}
To complete the induction and verify that $\phi$ induces an adjunction $\sF^{T'}_k \circ f_* \dashv V^T_k \circ f^*$, it is enough to check that for any $T'$-algebra $\sA$, the induced map
\begin{equation*}
	\alg_{T'}(\sF^{T'} f_*D, \sA) \stackrel{f^*}{\to} \alg_T(f^* \sF^{T'} f_* D, f^* \sA) \stackrel{\phi^*}{\to} \alg_T(\sF^T D, f^* \sA)
\end{equation*}
is an isomorphism. This follows from a direct calculation, similar to the one from the previous proposition. 	
\end{proof}

\begin{corollary}\label{cor:leftadjointmorphismoftheories}
	If the functor $f^*: \alg_{T'} \to \alg_T$ admits a left adjoint $L$, then there is a canonical natural isomorphism $L \circ \sF^T \cong \sF^{T'} \circ f_*$, and the homomorphism $\phi$ is induced by the unit map
	\begin{equation*}
		\sF^T(D) \to f^* L \sF^T(D) \cong f^* \sF^{T'} (f_* D). \qed
	\end{equation*}
\end{corollary}

\begin{remark}
	Batanin's approach to generalized computads in \cite{MR1664991} has some minor differences with ours. Rather than construct an adjunction from $k$-computads to $T$-algebras, Batanin constructs a family of monads $T_k$ on $k$-dimensional globular sets from the single monad $T$. This is accomplished by conjugating the monad $T$ by the inclusion and truncation functors between $k$-dimensional globular sets and $n$-dimensional globular sets. Then Batanin constructs an adjunction between $T_k$-computads and $T_k$-algebras, rather than with $T$-algebras. Nevertheless our category $\comp_k^T$ and Batanin's category of $T_k$-computads are equivalent. Batanin doesn't consider morphisms of monads in \cite{MR1664991}. 
\end{remark}

\section[Algebraic properties of symmetric monoidal bicats]{Algebraic Properties of the theories of Symmetric Monoidal Bicategories}\label{sec:algpropsSymBicats}

The results of Section~\ref{sec:strictsymbicats} give us a very pleasing situation with regards to the four theories of symmetric monoidal bicategory we have encountered. Every quasistrict symmetric monoidal 2-category gives rise to an unbiased semistrict symmetric monoidal 2-category; Every unbiased semistrict symmetric monoidal 2-category gives rise to a Crans semistrict symmetric monoidal 2-category; and every Crans semistrict symmetric monoidal 2-category gives rise to a symmetric monoidal bicategory. Moreover each of these four notions has an underlying bicategory and hence an underlying globular set. In each case the forgetful functor is the right adjoint of a monadic adjunction, thus we have a chain of monadic transformations:
\begin{equation*}
	T \to T_{css} \to T_{uss} \to T_{qs}
\end{equation*}
where $T$, $T_{css}$, $T_{uss}$, and $T_{qs}$ correspond to the monads on 2-dimensional globular sets for the theories of symmetric monoidal bicategories, Crans semistrict symmetric monoidal 2-categories, unbiased semistrict symmetric monoidal 2-categories, and quasistrict symmetric monoidal 2-categories, respectively. We will now show that each of these monads is finitary.

The category of 2-dimensional globular sets is a presheaf category on a finite category, hence locally finitely presentable. Thus we may apply the general machinery of finite equational theories for finitary endofunctors on locally finitely presentable categories. This is a very general context which applies to algebraic theories of many kinds. Roughly the idea is to describe a given endofunctor in terms a presentation. This presentation will consist of a {\em signature} and a set of {\em flat equations}. The signature, defined formally below, is a set parametrizing the possible finitary operations of the theory and the flat equations are the equations that these operations must satisfy. The main result in this general setting is:

\begin{proposition}[{\cite[Sect.~5]{MR1239558} \cite[Pr.~4.22]{MR3003200}}]  \label{pro:}
	An endofunctor on a locally finitely presentable category is finitary if and only if it has a presentation by a {\em signature} and a set of {\em flat equations}. \qed
\end{proposition}

We will specialize to the case of endofunctors on 2-dimensional globular sets, and our treatment closely follows \cite{MR3003200}. Let $\gset_2$ denote the category of 2-dimensional globular sets, which is the presheaf category on $\G_2$, introduced in  Section~\ref{sec:computads}. Let $\cF \subseteq \gset_2$ denote the full-subcategory consisting of the {\em compact} presheaves, those presheaves which are finite colimits of representable presheaves.  These are exactly the presheaves $f$ such that the functor $\gset_2(f, -)$ commutes with filtered colimits. 

\begin{definition}\label{def:signature}
		A {\em signature} is a collection $\Sigma = (\Sigma_f)_{f \in \cF}$ of 2-dimensional globular sets indexed by (representatives of) the compact globular sets $\cF$.  
		
		A {\em $\Sigma$-algebra} is a globular set $A \in \gset_2$, together with a rule which assigns to every morphism $\phi:f \to A$, a morphism $\hat{\phi}:\Sigma_f \to A$. Given another $\Sigma$-algebra $B$, a {\em homomorphism} $h: A \to B$ is a map of globular sets such that $h \hat{\phi} = \widehat{h \phi}$.  
\end{definition}

Let $C_0$, $C_1$, and $C_2$ denote the representable globular sets. We will write $\sigma: f \leadsto C_k$ to mean an operation of arity-$f$ which produces a $k$-cell. Taking all of these together gives the global operation $f \leadsto \Sigma_f$. We now give some examples of the kinds of operations that arise in the above theories of symmetric monoidal bicategory. 

\begin{example}
	$\emptyset \leadsto C_0$. A nullary operation of this kind declares the existence of a distinguished 0-cell, for example the unit object $1 \in \sM$ in a symmetric monoidal bicategory. There are also unitary operations $C_0 \leadsto C_1$ and $C_1 \leadsto C_2$ which, for example, can produce the identity morphisms.  
\end{example}

\begin{example}
	$C_0 \sqcup C_0 \leadsto C_0$. An operation of this kind is a binary operation on the 0-cells, producing a new 0-cell. An example of this kind is the zeroth component of the functor $\otimes$ in a symmetric monoidal bicategory. 
\end{example}

\begin{example}
	$C_1 \cup^{C_0} C_1 \leadsto C_1$. Here the pushout is formed using the two inclusions of the 0-cell. An operation of this kind takes a pair of `composable' 1-cells and produces a new 1-cell. An example of this type is the composition operation for 1-morphisms. 
\end{example}

\begin{example}
	$C_0 \cup C_0 \leadsto C_1$. We can also have operations which take, say, a pair of objects and produce a 1-cell. For example the component of the braiding transformation $\beta_{x,y}$ is an operation of this kind. Each of the structural natural transformations and modifications are variations on this sort of operation. 
\end{example}

\begin{example}
	$(C_1 \cup^{C_0} C_1) \sqcup (C_1 \cup^{C_0} C_1) \leadsto C_2$. The `interchanger' $\phi^\otimes_{(f,f'),(g,g')}$ is an operation which takes two pairs of composable 1-morphisms and produces a 2-morphism. Similarly all the other coherence cells can be interpreted as coming from finitary operations. 
\end{example}

Thus we see that every bit of data that makes up the definition of symmetric monoidal bicategory (or Crans semistrict symmetric monoidal 2-category, unbiased semistrict symmetric monoidal 2-category, or quasistrict symmetric monoidal 2-category) can be viewed as coming from a finitary operation belonging to a signature. 

Each signature gives rise to a {\em polynomial endofunctor} $H_\Sigma$ on globular sets given by:
\begin{equation*}
	H_\Sigma(X) := \coprod_{f \in \cF} \Sigma_f \times \gset_2(f, X).
\end{equation*}
This endofunctor is finitary (it commutes with filtered colimits). The category of $\Sigma$-algebras and homomorphisms is equivalent to the category of algebras for the endofunctor $H_\Sigma$. The four theories of symmetric monoidal bicategory we are concerned with are not theories of this kind; the operations are subject to a variety of equations. 

\begin{definition}\label{def:flateqn}
	A {\em flat equation} for a signature $\Sigma$ means a pair of parallel arrows
	\begin{equation*}
		u,u': f \to H_\Sigma(g)
	\end{equation*}
	for some $f,g \in \cF$. 
\end{definition}

Given a set $I$ of flat equations $u_i, u'_i: f_i \to H_\Sigma(g_i)$, $i \in I$, we may define a new signature $\overline{\Sigma}$ as follows: 
\begin{equation*}
	\overline{\Sigma}_g : = \coprod_{i \in I, g = g_i} f_i.
\end{equation*}
The $u_i$ and $u'_i$ assemble into natural transformations:
\begin{equation*}
	u, u': H_{\overline{\Sigma}} \to H_\Sigma.
\end{equation*}

\begin{definition}\label{def:presentationofendofunctor}
	An endofunctor $T$ is presented by a signature $\Sigma$ and a set $I$ of flat equations for $\Sigma$ if it is the coequalizer
	\begin{equation*}
		H_{\overline{\Sigma}} \rightrightarrows H_{\Sigma} \to T
	\end{equation*}
	of the natural transformations $u$ and $u'$. 
\end{definition}

The $T$-algebras are exactly the $\Sigma$-algebras which {\em satisfy} the flat equations in the set $I$, meaning that the structure map $a: H_\Sigma(A) \to A$ coequalizes $H_\Sigma(phi) \circ u_i$ and $H_\Sigma(\phi) \circ u'_i$ for all $\phi: g_i \to A$.

\begin{example}
	Consider a signature which contains the operations $(-)\otimes (-)$ and $\beta_{-,-}$ on 0-cells. We may define the flat equation of the from $u, u': C_0 \to H_{\Sigma}(C_0 \sqcup C_0)$ by sending the unique  0-cell $z$ of $C_0$ to:
\begin{align*}
	u(z) & =  x_0 \otimes x_1 \\
	u'(z) &= s(\beta_{x_0,x_1})
\end{align*}	
where $x_0$ and $x_1$ denote the two 0-cells of $C_0 \sqcup C_0$, and $s$ denoted the source map. A $\Sigma$-algebra satisfying this equation satisfies the condition that $\beta_{x,y}$ is a 1-cell whose source is $x \otimes y$, for all 0-cells $x,y \in A$. By adding additional flat equations one may establish all the source/target matching of the coherence cells. 
\end{example}

\begin{example}
	Other, more pedestrian examples include: 
	\begin{itemize}
		\item an equation of 1-cells establishing the equality of $(x \otimes y) \otimes z$ and $x \otimes (y \otimes z)$ (which holds in the semistrict monoidal theories);
		\item the naturality equations which express that $(\otimes, \phi^\otimes)$ is a functor of bicategories; similar equations for the other structural transformations and modifications. 
		\item the pentagon equation for the associator transformations, and all manner of other axioms in the theory of symmetric monoidal bicategories.
	\end{itemize}
\end{example}

In short all of the axioms which the four theories of symmetric monoidal bicategory must satisfy arise as flat equations for a signature. We have established that the four monads $T$, $T_{css}$, $T_{uss}$, and $T_{qs}$ are each presented by a signature and a set of flat equations for it, though at the end of the day this signature and the number of flat equations will be quite large (but finite). From this analysis we have the following corollary.

\begin{corollary}\label{cor:finitary-monads}
	Each of the four monads $T$, $T_{css}$, $T_{uss}$, and $T_{qs}$ is finitary.\qed
\end{corollary}

Thus theory of $T$-computads developed in Section~\ref{sec:computads} applies to each of these four theories. 

There are further consequences of this analysis. Each of these four theories is built upon the theory of bicategories, meaning that the monads for these theories can be obtained by taking the signature and flat equations which define the bicategory monad and adding additional operations and equations. (Note that strict 2-categories, which underly all the semistrict versions of symmetric monoidal 2-category, are also obtained from the theory of bicategories by imposing additional flat equations). 

\begin{corollary}
	Let $S= T$, $T_{css}$, $T_{uss}$, or $T_{qs}$ be one of the monads corresponding to the four theories of symmetric monoidal bicategories considered previously. Then the frgetful functor $\alg_S \to \bicat_s$ to bicategories
	 is part of a monadic adjunction. \qed
\end{corollary}

Similarly, the semistrict theories can be obtained from the theory of symmetric monoidal bicategories by adding additional operations and equations (in the Crans semistrict and the quasistrict cases we merely need to add extra flat equations). Thus we have:

\begin{corollary} \label{cor:fstarleftadjoint}
	The forgetful functors fit into a chain of adjunctions:
	\begin{equation*}
		\symbicat_{qs} \leftrightarrows \symbicat_{uss} \leftrightarrows \symbicat_{css} \leftrightarrows \symbicat
	\end{equation*}
	where in each case the morphisms of the relevant categories are the strict homomorphisms. \qed
\end{corollary}

Thus the conditions of Cor.~\ref{cor:leftadjointmorphismoftheories} are met.

\section{Further properties of Computads}\label{sec:furthercomputads}

Let us return to our discussion of general computads for monads on globular sets. The cells which are used to construct any globular set have a notion of dimension. In many cases of interest the given monad on globular sets respect this notion to some degree. We will say that the monad is {\em monotone}. Monotone monads and the theory of computads over them enjoy additional properties which we will describe below. 

To make this precise it will be useful to compare monads on globular sets of differing dimensions. Let $k < n$, and consider the canonical inclusion $i_k:\G_k \to \G_n$. This induces the fully-faithful {\em restriction functor} $i^*_k: \gset_n \to \gset_k$. This functor admits both a left adjoint $i_{k,!}: \gset_k \to \gset_n$ and a right adjoint $i_{k, *}: \gset_k \to \gset_n$. Thus given a monad $T$ on $n$-dimensional globular sets, we get an induced monad $T_k: i_k^* T i_{k,!}$ on $k$-dimensional globular sets. 

\begin{definition}\label{def:monotonemonad}
	A monad $T$ on $n$-dimensional globular sets is {\em monotone} if for each $k < n$, the canonical natural transformation $T_k i_k^* \to i_k^* T$ is an isomorphism. 
\end{definition}

In words this says that the $k$-dimensional part of a free $T$-algebra generated by a globular set $X$ only depends on the $k$-dimensional part of $X$. The examples we are interested in satisfy this property.

\begin{lemma}\label{lem:symmontheoriesaremonotone}
	Each of the monads corresponding to the theories of symmetric monoidal bicategories, Crans semistrict symmetric monoidal 2-categories, unbiased semistrict symmetric monoidal 2-categories, and quasistrict symmetric monoidal 2-categories are monotone. 
\end{lemma}

\begin{proof}
	As we saw in Section~\ref{sec:algpropsSymBicats}, each of these theories is presented by a signature and a set of flat equations. The signature gives us a collection of operations that are used to define the relevant theory. In each case these operations satisfy equations which make them dimension increasing, which means that an operation which takes as input cells of dimension $\geq k$ will only produce new cells of dimension $\geq k$. For example, the interchanger $\phi^\otimes$ is an operation which takes as input two sets of composable cells and produces a 2-cell. A 2-cell, of course, has source and target 1-cells and 0-cells.  As a bare operation from just a signature, this would produce new 0-cells. However, in the four theories we are considering, the interchanger is required to satisfy equations which identify its source and target with previously constructed 0-cells. The source and target 0-cells of the interchanger morphism $\phi^\otimes_{(f,f'), (g,g')}$ are given by the tensor product of the sources and targets of $f$ and $f'$, and $g$ and $g'$, respectively. In a similar manner all of the operations which make up the theories symmetric monoidal bicategories, Crans semistrict symmetric monoidal 2-categories, unbiased semistrict symmetric monoidal 2-categories, and quasistrict symmetric monoidal 2-categories satisfy this dimension increasing property, which directly implies that the corresponding monads are monotone. 
\end{proof}

\begin{lemma}\label{lem:monotonecolimit}
	Let $T$ be a monotone monad on $n$-dimensional globular sets. The the restriction functor induces an adjunction $(i_k^*, R)$
	\begin{equation*}
		i_k^*: \alg_T \leftrightarrows \alg_{T_k}: R
	\end{equation*}
	between the theories of $T$-algebras on $n$-dimensional globular sets and $T_k$-algebras on $k$-dimensional globular sets. In particular $i_k^*$ admits a right adjoint. This right adjoint is fully-faithful. Moverover if $T$ is finitary, then $i^*_k: \alg_T \to \alg_{T_k}$ admits a left-adjoint. 
\end{lemma}

\begin{proof}
	First, because $T$ is monotone restriction induces a functor $\alg_T \to \alg_{T_k}$, which we also denote by $i^*_k$, which is given by sending a $T$-algebra $(X,\alpha)$ to the $T_k$-algebra $(i_k^* X, i_k^* \alpha)$. Thus we have a diagram as follows:
\begin{center}
\begin{tikzpicture}
		\node (LT) at (0, 2) {$\alg_T$};
		\node (LB) at (0, 0) {$\gset_n$};
		\node (RT) at (3, 2) {$\gset_k$};
		\node (RB) at (3, 0) {$\alg_{T_k}$};
		\draw [->] ([xshift=1ex]LT.south) -- node [right] {$U$} ([xshift=1ex]LB.north);
		\draw [<-] ([xshift=-1ex]LT.south) -- node [left] {$F$} ([xshift=-1ex]LB.north);
		
		\draw [->] (LT) -- node [above] {$i_k^*$} (RT);
		\draw [->] ([xshift=1ex]RT.south) -- node [right] {$U_k$} ([xshift=1ex]RB.north);
		\draw [<-] ([xshift=-1ex]RT.south) -- node [left] {$F_k$} ([xshift=-1ex]RB.north);
		\draw [->] ([yshift=1ex]LB.east) -- node [above] {$i^*_k$} ([yshift=1ex]RB.west);
		\draw [<-] ([yshift=-1ex]LB.east) -- node [below] {$i_{k,*}$} ([yshift=-1ex]RB.west);
		%\node at (0.5, 1) {$\ulcorner$};
		%\node at (1.5, 0.5) {$\lrcorner$};
\end{tikzpicture}
\end{center}
Each of the three sides of double arrows represents an adjunction, and we wish to know that the top arrow admits a right adjoint. As it happens this is precisely the context of \cite[Thm.~3.7.4(a)]{MR2178101} which gives a precise set of six criteria, which if satisfied, ensure that the top arrow, indeed, admits a right adjoint. This right adjoint will satisfy $U \circ R = i_{k,*} \circ U_k$, but this implies that $i^*_k R \cong id$, and so $R$ is fully-faithful. 

Of these criteria, (i), (ii), (iv), and (v) are obvious. Criterion (iii) is also easily established using the fact that $T$ is monotone. This leaves the last criterion which asks that $\alg_T$ has coequalizers of $U$-contractible coequalizer pairs, and that $i^*_k: \alg_T \to \alg_{T_k}$ preserves them. 

In this case, because the adjunction $(F,U)$ is monadic, $\alg_T$ does admit coequalizers of $U$-contractible pairs and in fact $U$ creates them \cite[Prop.~3.3.5]{MR2178101}. Let $A \rightrightarrows B$ be such a coequalizer pair, and let $C$ be the coequalizer. Since $i_k^*: \gset_n \to \gset_k$ preserves colimits, $i^*_k UA \rightrightarrows i^*_k UB$ is also admits a coequalizer in $\gset_k$ (and is also a contractible coequalizer pair). Specifically the coequalizer is $i_k U C$. However, since $i^*_k U = U_k i_k^*$ (since $T$ is monotone), this implies that $i^*_k A \rightrightarrows i^*_k B$ is a $U_k$-contractible coequalizer pair in $\alg_{T_k}$. Furthermore since $(F_k, U_k)$ is a monadic adjunction, the functor $U_k$ creates coequalizers of $U_k$-contractible coequalizer pairs. It follows that $i_k^* C$ is the coequalizer of $i^*_k A \rightrightarrows i^*_k B$, which is exactly the criterion we wished to show. Thus $i^*_k: \alg_T \to \alg_{T_k}$ admits a right adjoint.

The last statement follows from a similar application of \cite[Thm.~3.7.4(b)]{MR2178101}, where in this case the only non-trivial criterion is that $\alg_T$ has coequalizers. This follows from the fact that $T$ is assumed to be finitary, and hence $\alg_T$ is cocomplete.  
	% G is fully-faithful.
\end{proof}

One upshot of the previous lemma is that for finitary monotone monads $T$, the functor $i^*_k: \alg_T \to \alg_{T_k}$ commutes with all limits and colimits of $T$-algebras that exist. This will be useful in what follows. 

If $\cP$ is a class of maps between globular sets, then we will say that a monad $T$ preserves $\cP$ if $Tf$ is in $\cP$ whenever $f$ is in $\cP$. Another consequence which we will explore is that monotone monads will automatically preserve certain kinds of special epimorphisms between globular sets. This will eventually allow us to show, among other things, that every symmetric monoidal bicategory is equivalent to a computadic symmetric monoidal bicategory, and similar results for the related semistrict theories. 

Before describing what holds true for monotone monads, it is instructive to understand what can go wrong. An epimorphism of globular sets is a map which is surjective on each sets of cells (surjective on 0-cells, surjective on 1-cells, etc.). In general monads on globular sets, even monotone monads, will fail to preserve epimorphisms. 

\begin{example}
	Consider the free category monad $T$ on 1-dimensional globular sets (the category of 1-dimensional globular sets is the same as the category of graphs). Let $X = C_1 \sqcup C_1$ be the globular set consisting of two disjoint edges (it has four 0-cells). Let $Y$ be the globular set $C_1 \cup^{C_0} C_1$ obtained by gluing the edges together, one tail to one tip. $Y$ has three $0$-cells, and may be obtained as a quotient of $X$. The quotient map $q:X \to Y$ is an epimorphism of globular sets. 
	
	However the map $Tq: TX \to TY$ is not an epimorphism of gobular sets. $TX$ is the globular set underlying the free category generated by $X$, while $TY$ is gloular set underlying the free category generated by $Y$. The map $Tq$ is not surjective on 1-cells. In particular $TY$ contains a {\em composite} 1-cell, which fails to be in the image of $TX$ (the relevant arrows are not composable in $TX$). 
\end{example}

I am confident that the reader can imagine how this problem becomes even more serious for the theory of symmetric monoidal bicategories and its cousins. Nonetheless, there are interesting classes of epimorphisms which are preserved by all monotone monads. Let us introduce some terminology.
\begin{definition}\label{def:acyclicfib}
	A map of $n$-globular sets $f:X \to Y$ is {\em $k$-bijective} if it induces an isomorphism on sets of cells $X_m \to Y_m$ for $0 \leq m \leq k$.  We will say that it is {\em $k$-surjective} if it is $(k-1)$-bijective and the induced map on $k$-cells $X_k \to Y_k$ is surjective. We will say that it is {\em $k$-injective} if it is $(k-1)$-bijective and the induced map $X_k \to Y_k$ is injective. 
	
		It is {\em $k$-acyclic} if it has the right-lifting property with respect to the inclusions $\partial C_m \to C_m$ for each $0 \leq m \leq k$.	
%	the induced map of sets $X_m \to Y_m$ is surjective for each $0 \leq m \leq k$.  
And finally, a map will be called an {\em acyclic fibration} if it is $n$-acyclic and also has the right-lifting property with respect to the collapse map $\partial C_{n+1} \to C_n$.
\end{definition}

Thus a $k$-bijection is precisely a morphism $f: X \to Y$ such that $i_k^*(f)$ is an isomorphism. The $n$-surjective morphisms between $n$-dimensional globular sets are characterized as exactly those morphisms which are $(n-1)$-bijective and admit a section. Furthermore, a map is $k$-surjective if and only if $i^*_k(f)$ is $k$-surjective. Every $k$-surjection (hence 
also every $k$-bijection) is in particular a $k$-acyclic morphism. 

Recall the following standard notion:
\begin{definition}
	Let $\cC$ be a category. An {\em orthogonal factorization system} consists of two classes of arrows of $\cC$ $(E,M)$ such that
	\begin{enumerate}
		\item every morphism in $\cC$ factors as a composite of a morphism in $E$ followed by a morphism in $M$;
		\item $E$ is the class of morphisms which have the left-lifting property with respect to $M$;
		\item $M$ is the class of morphisms which have the right-lifting property with respect to $E$; and
		\item For each lifting problem
		
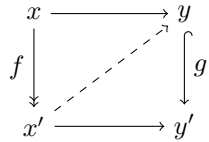
\begin{center}
		\begin{tikzpicture}
			\node (LT) at (0, 1.5) {$x$};
			\node (LB) at (0, 0) {$x'$};
			\node (RT) at (2, 1.5) {$y$};
			\node (RB) at (2, 0) {$y'$};
			\draw [->>] (LT) -- node [left] {$f$} (LB);
			\draw [->] (LT) -- node [above] {$$} (RT);
			\draw [right hook->] (RT) -- node [right] {$g$} (RB);
			\draw [->] (LB) -- node [below] {$$} (RB);
			\draw [->, dashed] (LB) -- (RT);
			%\node at (0.5, 1) {$\ulcorner$};
			%\node at (1.5, 0.5) {$\lrcorner$};
		\end{tikzpicture}
		\end{center}
with $f \in E$ and $g \in M$, the solution to the lifting problem (the dashed arrow) is unique. 
	\end{enumerate}
\end{definition}

We have the following basic fact concerning $n$-surjections and $n$-injections:

\begin{lemma}
	The $n$-surjections and the $n$-injections form an orthogonal factorization system for the category of $n$-globular sets and $(n-1)$-bijections. \qed
\end{lemma}

Let $\cP$ be a class of maps between $n$-dimensional globular sets, such as the classes of $k$-bijective maps, $k$-acyclic maps, $k$-surjective maps, or acyclic fibrations. If $T$ is a monad on $n$-dimensional globular sets then a morphism of $T$-algebras will be said to be of class $\cP$ if its underlying map of globular sets is so.  

Let $\sC$ be a category and $\cP$ a class of maps. We will say that $\cP$ is {\em closed under pullbacks} if for each pullback square (depicted on the left-hand side)
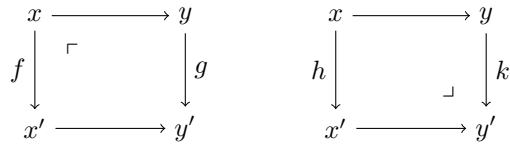
\begin{center}
\begin{tikzpicture}
	\node (LT) at (0, 1.5) {$x$};
	\node (LB) at (0, 0) {$x'$};
	\node (RT) at (2, 1.5) {$y$};
	\node (RB) at (2, 0) {$y'$};
	\draw [->] (LT) -- node [left] {$f$} (LB);
	\draw [->] (LT) -- node [above] {$$} (RT);
	\draw [->] (RT) -- node [right] {$g$} (RB);
	\draw [->] (LB) -- node [below] {$$} (RB);
	\node at (0.5, 1) {$\ulcorner$};
	%\node at (1.5, 0.5) {$\lrcorner$};
	\node (LTA) at (4, 1.5) {$x$};
	\node (LBA) at (4, 0) {$x'$};
	\node (RTA) at (6, 1.5) {$y$};
	\node (RBA) at (6, 0) {$y'$};
	\draw [->] (LTA) -- node [left] {$h$} (LBA);
	\draw [->] (LTA) -- node [above] {$$} (RTA);
	\draw [->] (RTA) -- node [right] {$k$} (RBA);
	\draw [->] (LBA) -- node [below] {$$} (RBA);
	\node at (5.5, 0.5) {$\lrcorner$};	
\end{tikzpicture}
\end{center}
such that $g$ is in $\cP$, we have that $f$ is in $\cP$. Dually we will say $\cP$ is {\em closed under pushouts} if for each pushout square (depicted on the right-hand side) such that $k$ is in $\cP$, we have that $h$ is in $\cP$. 

\begin{lemma} \label{lem:pullbackofclassesofmaps}
	In the category $\gset_n$ of $n$-dimensional globular sets, the classes of $k$-bijective maps, $k$-surjective maps, and $k$-injective maps, are preserved under both pullbacks and pushouts. The classes of $k$-acyclic maps and acyclic fibrations are preserved under pullbacks. 
	
	For every monad $T$ on $n$-dimensional globular sets, the category $\alg_T$ is complete and each of the classes of $k$-bijective maps, $k$-surjective maps, $k$-injective maps, $k$-acyclic maps, and acyclic fibrations are preserved under pullbacks. 
\end{lemma}

\begin{proof}
	In the category of sets the isomorphisms, surjections, and injections are preserved under both pushouts and pullbacks. The category of globular sets is a presheaf category and so limits and colimits are computed `objectwise'. From these observations it immediately follows that the classes of $k$-bijective maps, $k$-surjective maps, and $k$-injective maps between globular sets are also preserved under both pullbacks and pushouts. 
	
	The classes of  $k$-acyclic maps and acyclic fibrations are defined as those maps with the right-lifting property with respect to a fixed set of maps. Such classes are automatically preserved under pullbacks. 
	
	For any monad, the forgetful functor $U: \alg_T \to \gset_n$ creates limits \cite[Thm.~3.4.2]{MR2178101}, thus since $\gset_n$ is complete it follows that $\alg_T$ is complete. Moreover it follows that any classes of maps between globular sets which are preserved under pullbacks, are also preserved under pullbacks in the category of $T$-algebras. 
\end{proof}

\begin{lemma} \label{lem:imageofbijisTalg}
	Let $T$ be a monad on $n$-dimensional globular sets which preserves $n$-surjections (we will see this is automatically the case for monotone monads). Then the image of an $(n-1)$-bijective $T$-algebra map is a subalgebra of its codomain. 
\end{lemma}

\begin{proof}
	What this means is that if $f: A \to C$ is an $(n-1)$-bijective $T$-algebra map, and 
	\begin{equation*}
		A \twoheadrightarrow C \hookrightarrow B
	\end{equation*}
	is the unique factorization of $f$ as a map of globular sets into an $n$-surjection followed by an $n$-injection, then $C$ admits a (unique)  $T$-algebra structure and moreover each of these maps is in fact a $T$-algebra homomorphism. The proof of this is exactly as the proof of \cite[Prop.~9.3.7]{MR2178101}. 
\end{proof}

\begin{proposition}
	Let $T$ be a monotone monad on $n$-dimensional globular sets. Then $T$ preserves $k$-bijections and $k$-surjections. Moreover if $T$ is finitary (so that the theory of $T$-algebras is cocomplete) then in the category of $T$-algebras $k$-bijective and $k$-surjective maps are preserved along pushouts.
\end{proposition}

\begin{proof}
	Let $T$ be a monotone monad on $n$-dimensional globular sets. Let $f: X \to Y$ be a $k$-bijective map of globular sets. Hence $i^*k(f)$ is an isomorphism. Thus $T_k i^*_k(f) = i^*_k T(f)$ is an isomorphism, and hence $T(f)$ is $k$-bijective. Thus $T$ preserves $k$-bijective maps. 
	
	Now suppose that $f: X \to Y$ is $k$-surjective. Since $T$ is monotone, $T(f)$ is $k$-surjective if and only if $T_k(i^*_k(f))$ is $k$-surjective. Thus it suffices to consider the case $k=n$. The $n$-surjective maps between $n$-dimensional globular sets are characterized by the property that they are $(n-1)$-bijective, and that they admit a section. Since $T$ is a functor, $T(f)$ will also admit a section. Since $T$ is monotone, $T(f)$ is $(n-1)$-bijective, hence $T(f)$ is $n$-surjective. Thus $T$ preserves $k$-surjections. 
	
Now suppose that $T$ is finitary and consider a pushout square of $T$-algebras. 
\begin{center}
\begin{tikzpicture}
	\node (LT) at (0, 1.5) {$Y$};
	\node (LB) at (0, 0) {$X$};
	\node (RT) at (2, 1.5) {$Z$};
	\node (RB) at (2, 0) {$P$};
	\draw [->] (LT) -- node [left] {$f$} (LB);
	\draw [->] (LT) -- node [above] {$$} (RT);
	\draw [->] (RT) -- node [right] {$g$} (RB);
	\draw [->] (LB) -- node [below] {$$} (RB);
	%\node at (0.5, 1) {$\ulcorner$};
	\node at (1.5, 0.5) {$\lrcorner$};
\end{tikzpicture}
\end{center}
Suppose that $f$ is $k$-bijective. By Lemma~\ref{lem:monotonecolimit}, $i_k^*$ preserves all colimits of $T$-algebras, in particular pushouts. So this square remains a pushout square after applying $i^*_k$. However $i^*_k(f)$ is an isomorphism by assumption, and hence $i^*_k(g)$ is also an isomorphism, and so $g$ is $k$-bijective. Thus in the category of $T$-algebras, $k$-bijective maps are preserved along pushouts. 

Now consider the same pushout square, but assume only that $f$ is $k$-surjective. As before (and due to Lemma~\ref{lem:monotonecolimit}) it suffices to consider the case $k=n$. Since $f$ is $(n-1)$-bijective, it follows that $g$ is $(n-1)$-bijective. 

To show that $g$ is actually $n$-surjective, we will adapt the argument of \cite[Prop.~9.3.8]{MR2178101} to produce an explicit construction of this pushout. From the explicit construction it will follow that $g$ is $n$-surjective. To this end consider the set of all commutative squares in which $Z \to P_i$ is an $n$-surjective map.
\begin{center}
\begin{tikzpicture}
	\node (LT) at (0, 1.5) {$Y$};
	\node (LB) at (0, 0) {$X$};
	\node (RT) at (2, 1.5) {$Z$};
	\node (RB) at (2, 0) {$P_i$};
	\draw [->>] (LT) -- node [left] {$f$} (LB);
	\draw [->] (LT) -- node [above] {$$} (RT);
	\draw [->>] (RT) -- node [right] {$g_i$} (RB);
	\draw [->] (LB) -- node [below] {$$} (RB);
	%\node at (0.5, 1) {$\ulcorner$};
	%\node at (1.5, 0.5) {$\lrcorner$};
\end{tikzpicture}
\end{center}
We have an induced map $Z \to \prod P_i$, however it will no longer be $(n-1)$-bijective. To remedy this note that we have canonical isomorphisms of $T$-algebras 
\begin{equation*}
	R i^*_{n-1}(\prod P_i) \cong \prod R i^*_{n-1}(\prod P_i) \cong \prod_i R i^*_{n-1} Z
\end{equation*}
where $R$ denotes the right adjoint to $i^*_{n-1}: \alg_T \to \alg_{T_{n-1}}$. 
This is due to the fact that both $R$ and $i^*_{n-1}$ commute with products of $T$-algebras. Moreover the unit map $\prod P_i \to R i^*_{n-1}(\prod P_i)$ is, by construction, $(n-1)$-bijective. Thus we may form the pullback of $T$-algebras:
\begin{center}
\begin{tikzpicture}
	\node (LT) at (0, 1.5) {$\tilde{P}$};
	\node (LB) at (0, 0) {$R i^*_{n-1} X$};
	\node (RT) at (4, 1.5) {$\prod P_i$};
	\node (RB) at (4, 0) {$R i^*_{n-1}(\prod P_i)$};
	\draw [->] (LT) -- node [left] {$$} (LB);
	\draw [->] (LT) -- node [above] {$$} (RT);
	\draw [->] (RT) -- node [right] {$$} (RB);
	\draw [->] (LB) -- node [above] {$\Delta$} (RB);
	\node at (0.5, 1) {$\ulcorner$};
	%\node at (1.5, 0.5) {$\lrcorner$};
\end{tikzpicture}
\end{center}
where the bottom map is the diagonal map. By Lemma~\ref{lem:pullbackofclassesofmaps}, $Z \to \tilde{P}$ is $(n-1)$-bijective. Moreover since $Z \to P$ is $(n-1)$-bijective, the map $P \to \prod P_i$ automatically factors through $\tilde{P}$.

Next we may apply Lemma~\ref{lem:imageofbijisTalg}. We may factor the previous map into an $n$-surjective $T$-algebra homomorphism followed by an $n$-injective one:
\begin{equation*}
	Z \twoheadrightarrow P_0 \hookrightarrow \tilde{P}.
\end{equation*}
Since $Y \twoheadrightarrow X$ is assumed to be $n$-surjective, there exists a unique lift $X \to P_0$ making the following square commute:
\begin{center}
\begin{tikzpicture}
	\node (LT) at (0, 1.5) {$Y$};
	\node (LB) at (0, 0) {$X$};
	\node (RT) at (2, 1.5) {$Z$};
	\node (RB) at (2, 0) {$P_0$};
	\draw [->] (LT) -- node [left] {$$} (LB);
	\draw [->] (LT) -- node [above] {$$} (RT);
	\draw [->] (RT) -- node [right] {$$} (RB);
	\draw [->] (LB) -- node [below] {$$} (RB);
	%\node at (0.5, 1) {$\ulcorner$};
	%\node at (1.5, 0.5) {$\lrcorner$};
\end{tikzpicture}
\end{center}
We claim that this square realizes $P_0$ as the pushout $P_0 \cong P = X \cup^Y Z$. By construction $g:Z \to P_0$ is an $n$-surjective map, and so if the claim holds we are done. 

We must show that $P_0$ satisfies the universal property of the pushout. Consider a commutative square of $T$-algebras
\begin{center}
\begin{tikzpicture}
	\node (LT) at (0, 1.5) {$Y$};
	\node (LB) at (0, 0) {$X$};
	\node (RT) at (2, 1.5) {$Z$};
	\node (RB) at (2, 0) {$W$};
	\draw [->] (LT) -- node [left] {$$} (LB);
	\draw [->] (LT) -- node [above] {$$} (RT);
	\draw [->] (RT) -- node [right] {$$} (RB);
	\draw [->] (LB) -- node [below] {$$} (RB);
	%\node at (0.5, 1) {$\ulcorner$};
	%\node at (1.5, 0.5) {$\lrcorner$};
\end{tikzpicture}
\end{center}
The map $P \to W$ factors through $W'$, where $W'$ is defined as the pullback
\begin{center}
\begin{tikzpicture}
	\node (LT) at (0, 1.5) {$W'$};
	\node (LB) at (0, 0) {$Ri^*_{n-1} P = R i^*_{n-1} Z$};
	\node (RT) at (4, 1.5) {$W$};
	\node (RB) at (4, 0) {$Ri^*_{n-1} W$};
	\draw [->] (LT) -- node [left] {$$} (LB);
	\draw [->] (LT) -- node [above] {$$} (RT);
	\draw [->] (RT) -- node [right] {$$} (RB);
	\draw [->] (LB) -- node [below] {$$} (RB);
	\node at (0.5, 1) {$\ulcorner$};
	%\node at (1.5, 0.5) {$\lrcorner$};
\end{tikzpicture}
\end{center}
Since the maps $P \to W'$ and $Z \to P$ are $(n-1)$-bijective, to demonstrate the universal property of $P_0$ it suffices to consider only those commutative squares in which the map $Z \to W$ is $(n-1)$-bijective. So without loss of generality suppose this is the case. 

Applying Lemma~\ref{lem:imageofbijisTalg} again, the map $Z \to W$ factors
\begin{equation*}
	Z \twoheadrightarrow P_i \hookrightarrow W
\end{equation*}
as an $n$-surjective $T$-algebra map (which is necessarily isomorphic to one of our original $g_i: Z \to P_i$), followed by an $n$-injection. Because $Y \to X$ is an $n$-surjection, there exists a unique lift $X \to P_i$, and thus the map $P \to W$ factors uniquely through the square defining $P_i$. By construction this factors uniquely through $P_0$, and thus we have shown that $P_0$ satisfies the universal property of the pushout. 
\end{proof}

\begin{lemma} \label{lem:acyclicfibrationlemma}
	Let $T$ be a finitary monotone monad on $n$-globular sets.
	% Suppose that the underling endo functor of $T$ preserves $k$-bijections and moreover in the category of $T$-algebras $k$-bijections are preserved along pushouts. 
	Then for $0\leq k \leq n$ we have:
\begin{enumerate}
	\item for any $k$-computad $D= (C,X,x) \in \comp^T_k$  the canonical map $\sF_{k-1}C \to \sF_k D$ is $(k-1)$-bijective; 
	\item for any $n$-computad $D=(C,X,x)\in \comp^T_n$,  the canonical map  $\sF_n(C) \to \sF_{n+1} D$ is $n$-surjective; and
	\item the counit $\sF_{k} V_{k} \sC \to \sC$ is $k$-acyclic. 
\end{enumerate}	
Moreover, 
%if $T$ preserves $n$-surjective maps and in the category of $T$-algebras $n$-surjective maps are preserved along pushouts, then 
the counit $\sF_{n+1} V_{n+1} \sC \to \sC$ is an acyclic fibration.
\end{lemma}

\begin{proof}
For the first statement we have that $\partial C_{k} \to C_k$ is $(k-1)$-bijective, therefore by our assumptions on $T$, so is the map of free $T$-algebras $T(\partial C_k) \to T(C_k)$, and therefore also the pushout $\sF_{k-1}C \to \sF_k D$.	

The second statement can be seen as follows. First note that by the first statement the canonical map  $\sF_{k-1} V_{k-1} \sC \to \sF_k V_k \sC$ is $(k-1)$-bijective. Now suppose that we are given a lifting problem:
\begin{center}

\end{center}
Since $\sF_{k-1} V_{k-1} \sC \to \sF_k V_k \sC$ is $(k-1)$-bijective there exists a unique lift making the following diagram commute:
\begin{center}
%
\end{center}
Now the commutative square
\begin{center}
%
\end{center}
defines an element $x_0 \in X_{\sC, k}$. This element gives rise to a commutative diagram:
\begin{center}
%
\end{center}
and by construction the composite $C_k \to U\sF_k V_k \sC \to U\sC$ is the original map. Hence this defines the desired lift. 

For the final statement we note that an identical argument shows that $\sF_n V_n \sC \to \sF_{n+1} V_{n+1} \sC$ is $n$-surjective. Hence, given a lifting square as below, we can find the dashed arrow making the diagram commute.  
\begin{center}
%
\end{center}
This lift defines an element of $x_0 \in X_{(\sC,n+1)}$, and the remainder of the argument proceeds exactly as in the previous case. 
\end{proof}	

\begin{corollary} \label{cor:everysymbicatequivtocomputadic}
	Every symmetric monoidal bicategory is equivalent, as a symmetric monoidal bicategory, to a computadic one via the strict symmetric monoidal homomorphism $\sF_3 V_3 \sM \to \sM$. 
\end{corollary}

\begin{proof}
	By the previous lemma, for every symmetric monoidal bicategory $\sM$ the strict homomorphism of symmetric monoidal bicategories $\sF_3 V_3 \sM \to \sM$ is an acyclic fibration with source a computadic symmetric monoidal bicategory. Since an acyclic fibration is surjective on onjects and 1-morphisms, and fully-faitful on 2-morphisms it follows from 	Theorem~\ref{WhiteheadforSymMonBicats} (Whiteheads theorem for symmetric monoidal bicategories) that this homomorphism is an equivalence of symmetric monoidal bicategories. 
\end{proof}

\subsection{A brief digression on Quillen model categories}

 We will now focus on the case where the given higher categorical structure is built upon the theory of bicategories. More precisely we will consider a finitary monad on $2$-dimensional globular sets $T$, together with a map of monads $f: B \to T$ from the bicategory monad $B$. As we saw previously this gives rise to a functor
\begin{equation*}
	f^*:\alg_T \to \bicat_s = \alg_B
\end{equation*}
Important examples are when $T$ is the strict $2$-category monad (in which case $f^*$ is the inclusion of strict 2-categories into bicategories), when $T$ is the symmetric monoidal bicategory monad (in which case $f^*$ is a forgetful functor), and also the three semistrict theories: Crans semistrict symmetric monoidal 2-categories, unbiased semistrict symmetric monoidal 2-categories, and quasistrict symmetric monoidal 2-categories. In each of these cases, and in many others of interest, the functor $f^*$ will admit a left adjoint and the resulting adjunction $(L, f^*)$ will be monadic. 

The category $\bicat_s$ has an important Quillen model structure on it constructed by Lack \cite{MR2138540}. It is interesting to know when this model structure can be transferred to a model structure on $\alg_T$. 

\begin{proposition}\label{pro:TalgModelStr}
	Let $T$ be a finitary monad on $2$-globular sets with a monad morphism $f: B \to T$ from the bicategory monad. Suppose that
	\begin{enumerate}
		\item The functor $f^*: \alg_T \to \bicat_s$ is the right-adjoint of a monadic adjunction; and
		\item There exists path objects in $\alg_T$ in the sense of \cite{MR1734325}.
	\end{enumerate}
	Then Lack's model structure on $\bicat_s$ \cite{MR2138540} induces a {\em transferred} model structure on $\alg_T$; this is a model structure in which a morphism of $T$-algebras $X \to Y$ is a fibration, respectively weak equivalence, if and only if the induced morphism of bicategories $f^*X \to f^* Y$ is so (hence the acyclic fibrations of this model structure agree with the acyclic fibrations from Definition~\ref{def:acyclicfib}).  It is combinatorial and every object is fibrant (hence it is a right proper model structure). 
		
	Moreover if $T$ is monotone, then the comonad $\sF_{3}V_{3}$ induced by the adjunction with $T$-computads is a functorial cofibrant replacement for this model structure. In this case the cofibrant objects are characterized as the retracts of computadic $T$-algebras.
	\qed
\end{proposition}

We won't need the above proposition in what follows and so we will only indicate its proof briefly. All but the last statements are a direct application of \cite[Lma~2.3(2)]{MR1734325}, which gives criteria for the existence of a transferred model structure. There are four generating cofibrations of the resulting transferred model structure: $T(\emptyset) \to T(C_0)$, $T(\partial C_1) \to T(C_1)$, $T(\partial C_2) \to T(C_2)$, and $T(\partial C_3) \to T(C_2)$.
By construction $\sF_{3}V_{3} (\sC)$ is an iterated pushout of these basic maps, and hence cofibrant. In fact these pushouts are preformed in order of dimension; the computadic $T$-algebras are direct analogs of CW-complexes, where the cells are attached in the order of dimension. The counit map $\sF_{3}V_{3} \sC \to \sC$ is an acyclic fibration by Lemma~\ref{lem:acyclicfibrationlemma}. Thus if $\sC$ is cofibrant in the transferred model structure, we get a diagram
\begin{center}
\begin{tikzpicture}
	\node (LT) at (0, 1.5) {$T(\emptyset)$};
	\node (LB) at (0, 0) {$\sC$};
	\node (RT) at (2, 1.5) {$\sF_3 V_3 \sC$};
	\node (RB) at (2, 0) {$\sC$};
	\draw [right hook->] (LT) -- node [left] {$$} (LB);
	\draw [->] (LT) -- node [above] {$$} (RT);
	\draw [->>] (RT) -- node [right] {$$} (RB);
	\draw [->] (LB) -- node [below] {$=$} (RB);
	\draw [->, dashed] (LB) -- (RT);
	%\node at (0.5, 1) {$\ulcorner$};
	%\node at (1.5, 0.5) {$\lrcorner$};
\end{tikzpicture}
\end{center}
where the left-most vertical map is a cofibration and the right-most vertical map is an acyclic fibration. Hence the dashed lift exists making the diagram commute. In particular $\sC$ is a retract of the computadic $T$-algebra $\sF_3 V_3(\sC)$.

\begin{remark} \label{rmk:symbicatmodelstr}
	Although we won't pursue it here, let us remark that the above proposition applies to the case where $T$ is the symmetric monoidal bicategory monad, defined in later sections. The least straightforward property is the existence of a path object. However this can be constructed as the functor category $P \sC = \bicat(E, f^*\sC)$ of weak homomorphisms, transformations, and modifications into the underlying bicategory of $\sC$, where $E$ is the {\em free-walking isomorphism}, the contractible category on two objects. The bicategory $P \sC$ inherits a `pointwise' symmetric monoidal structure from $\sC$.  
\end{remark}

\section{Presentations of Symmetric Monoidal Bicategories} \label{SectFreelyGenSymMonBicats}

In this section we give a detailed account of {\em symmetric monoidal computads}, which are computads as defined in Section~\ref{sec:computads} for the symmetric monoidal bicategory monad on 2-dimensional globular sets. These computads define the data of a presentation for a symmetric monoidal bicategory. They consist of a collection of generating objects, 1-morphisms, 2-morphisms, and relations together with appropriate maps. Given such data $P$, we will be able to produce a symmetric monoidal bicategory, $\sF(P)$. The details considered here will then be used in the next section to prove the cofibrancy theorem which states that $\sF(P)$ has the property that the bicategory of symmetric monoidal homomorphism out of $\sF(P)$ into $\sM$ is equivalent to a bicategory of ``$P$-shaped data'' in the target symmetric monoidal bicategory $\sM$.

 %introduce freely generated symmetric monoidal bicategories. Roughly, the idea is that given certain {\em generating data}, $G$, we will be able to produce a symmetric monoidal bicategory, $\sF(G)$, which has the property that the bicategory of symmetric monoidal homomorphism out of $\sF(G)$ into $\sM$ is equivalent to a bicategory of ``$G$-shaped'' data in the target symmetric monoidal bicategory $\sM$. The generating data will consist of, among other things, three sets $G_0$, $G_1$ and $G_2$ called {\em generating objects}, {\em generating 1-morphism}, and {\em generating 2-morphisms}, respectively, and $\sF(G)$ should be regarded as the free symmetric monoidal bicategory whose objects are generated by $G_0$, whose 1-morphisms are generated by $G_1$ and whose 2-morphisms are generated by $G_2$. 

If we only cared about symmetric monoidal bicategories generated by a bicategory $G$, 
 then the symmetric monoidal bicategory $\sF(G)$ would be relatively straightforward to construct.\footnote{For example we could consider $\sF$ as the (weak) left adjoint to the forgetful (tri-)homomorphism $\symbicat \to \bicat$.} However, we wish to allow more general generating data which may not be easily organized into a bicategory. For example, we will wish to allow our generating 1-morphisms to have sources and targets which are {\em not} part of the generating objects, but merely consequences, (bracketed) words in the generating objects like ``$(w \otimes x) \otimes (y \otimes z)$''. Thus to give the precise definition of generating morphisms, we must already have some rudimentary knowledge of what the objects of $\sF(G)$ will be. A similar problem persists for generating 2-morphisms. Consequently, we are lead to the following sequence of recursive definitions, which simply spell out the details of  Section~\ref{sec:computads} in the case at hand.

\begin{definition}
A {\em symmetric monoidal $0$-computad} consists of a set $G_0$. Given such a set, a {\em binary word} in $G_0$ is a binary tree with leaves labeled by elements of $G_0$, i.e.,  we have the following recursive definition: the word $(1)$ and the words $(a)$ for $a \in G_0$ are binary words in $G_0$. If $u$ and $v$ are binary words in $G_0$ then so is $u \otimes v = (u) \otimes (v)$.
\end{definition}

\begin{definition}
	A {\em symmetric monoidal $1$-computad} consists of a set $G_0$, a set $G_1$, and a pair of maps $s,t: G_1 \to BW(G_0)$, where $BW(G_0)$ denotes the set of binary words in $G_0$. Given a 1-truncated  datum $G_1$,  then  a {\em binary word} in $G_1$, which we will define momentarily, has a {\em source} and a {\em target}, which are both binary words in $G_0$.  {\em Binary words} in $G_1$ are defined recursively as follows: \begin{itemize}

\item If $f \in G_1$, then $(f)$ is a binary word in $G_1$, with source $s(f)$ and target $t(f)$. 
\item If $u, v, w$ are binary words in $G_0$ and $x$ is a symbol from Table~\ref{BinWordsin1TruncatedGenDataTable}, then $(x)$ is a binary word in $G_1$ with source and target as listed in Table~\ref{BinWordsin1TruncatedGenDataTable}. 
\item If $u, v, w$ are binary words in $G_0$ and $x$ is a symbol from Table~\ref{BinWordsin1TruncatedGenDataTable}, then $(x^*)$ is a binary word in $G_1$ with source and target the opposite of those listed in Table~\ref{BinWordsin1TruncatedGenDataTable}. 
%\item If $f$ and $g$ are two binary words in $G_1$, then $f \otimes g$ is a binary word in $G_1$ with source $s(f) \otimes s(g)$ and target $t(f) \otimes t(g)$.  
\end{itemize}
\begin{table}[h]
\begin{center}
\begin{tabular}{|c|c| c|} \hline
Symbol & Source & Target \\  \hline
$I_u$ & $u$ &$u$ \\
  $\alpha_{u,v,w}$ & $(u \otimes v) \otimes w$ & $u \otimes ( v \otimes w)$ \\
  $\ell_u$ & $1 \otimes u$ & $u$ \\
  $r_u$ &$u$& $u \otimes 1$ \\
  $\beta_{u,v}$ & $u \otimes v$& $v \otimes u$ \\ \hline
\end{tabular}
\end{center}
\caption{Binary Words in Symmetric Monoidal 1-Computad.}
\label{BinWordsin1TruncatedGenDataTable}
\end{table}%
We also define {\em binary sentences} in $G_1$, which again have sources and targets which are binary words in $G_0$. These are defined recursively as follows:
\begin{itemize}
\item If $f$ is a binary word in $G_1$, then $(f)$ is a binary sentence in $G_1$ with source $s(f)$ and target $t(f)$. 
\item If $f$ and $g$ are binary sentences in $G_1$ such that $s(f) = t(g)$, then $(f) \circ (g)$ is a binary sentence in $G_1$ with source $s(g)$ and target $t(f)$. 
\item If $f$ and $g$ are binary sentences in $G_1$, then $f \otimes g$ is a binary sentence in $G_1$ with source $s(f) \otimes s(g)$ and target $t(f) \otimes t(g)$.  
\end{itemize}

\end{definition}

\begin{definition}
A {\em symmetric monoidal 2-computad} $G$, consists of a symmetric monoidal 1-computad $(G_0, G_1, s,t)$ and a set $G_2$ together with a pair of maps $s,t: G_2 \to BS(G_0, G_1, s,t)$ where $BS(G_0, G_1, s,t)$ denotes the set of binary sentences in $G_1$. These maps are required to satisfy:
\begin{align*}
& s(s(\epsilon)) = s(t(\epsilon)), \\
& t(s(\epsilon)) = t(t(\epsilon)),
\end{align*}
for all $\epsilon \in G_2$.  Given a symmetric monoidal 2-computad $(G_2,G_1, s,t)$, then a {\em binary word} in $G_2$ is defined recursively as follows. A binary word has a source and a target which are binary sentences in $G_1$.
\begin{itemize}
\item If $\epsilon \in G_2$ then $(\epsilon)$ is a binary word in $G_2$ with source $s(\epsilon)$ and target $t(\epsilon)$.
\item If $a, a', b, c,d$ are binary words in $G_0$, $f:a \to b, g, h$, $f'$ and $g'$ are binary sentences in $G_1$ such that $s(f) = t(f')$, $s(f') = t(f'')$ and $s(g) = t(g')$ and
$x$ is a symbol from Table~\ref{BinWordsin2TruncatedGenDataTable}, then $(x)$ is a binary word in $G_2$ with source and target as listed in Table~\ref{BinWordsin2TruncatedGenDataTable}. 
\item If $a, a', b, c,d$ are binary words in $G_0$, $f:a \to b, g, h$, $f'$ and $g'$ are binary sentences in $G_1$ such that $s(f) = t(f')$, $s(f') = t(f'')$ and $s(g) = t(g')$ and $x$ is a symbol from Table~\ref{BinWordsin2TruncatedGenDataTable}, then $(x^{-1})$ is a binary word in $G_2$ with source and target the opposite of those listed in Table~\ref{BinWordsin2TruncatedGenDataTable}. 
%\item If $u$ and $v$ are binary words in $G_2$, then $u \otimes v$ is a binary word in $G_2$ with source $s(u) \otimes s(v)$ and target $t(u) \otimes t(v)$. 
\end{itemize}
\begin{table}[ht]
\begin{center}
\begin{tabular}{|c|c| c|} \hline
Symbol & Source & Target \\  \hline
$id_f$ & $f$ &$f$ \\
$a^c_{f,f',f''}$ & $(f \circ f') \circ f''$ & $f \circ (f' \circ f'')$ \\
$r^c_f$ & $f \circ I_a$ & $f$ \\
$\ell_f^c$& $I_b \circ f$ & $f$ \\ \hline
$\eta_f$ & $I_a$ & $f^* \circ f$ \\
$\varepsilon_f$ & $f \circ f^*$ & $I_b$ \\ \hline
$\phi^\otimes_{(f,g), (f', g')}$ &$(f \otimes g) \circ (f' \otimes g')$& $(f \circ f') \otimes (g \circ g')$ \\
 $\phi^\otimes_{a, a'}$ &$ I_{a \otimes a'}$ & $I_a \otimes I_{a'}$ \\
  $\alpha_{f,g,h}$ & $[f \otimes ( g \otimes h)]\circ \alpha_{s(f), s(g), s(h)}$   & $\alpha_{t(f), t(g), t(h)} \circ[(f \otimes g) \otimes h]$ \\
  $\ell_f$ & $f\circ \ell_{s(f)}$ & $\ell_{t(f)} \circ (I_1 \otimes f)$ \\
  $r_f$ &  $(f \otimes I_1) \circ r_{s(f)}$ &  $ r_{t(f)} \circ f$   \\
  $\beta_{f,g}$ & $(g \otimes f) \circ \beta_{s(f), s(g)} $ & $\beta_{t(f), t(g)} \circ (f \otimes g)$  \\ \hline
$\pi_{a,b,c,d}$& $[( I \otimes \alpha) \circ \alpha ] \circ ( \alpha \otimes I)$ &$\alpha \circ \alpha$ \\
$\mu_{a,b}$ &$ [ I_a \otimes \ell_b) \circ \alpha_{a, 1, b}] \circ (r_a \otimes I_b)$ & $I_{a \otimes b}$ \\
$\lambda_{a,b}$ & $\ell_a \otimes I_b$ & $\ell_{a \otimes b} \circ \alpha_{1, a,b}$ \\
$\rho_{a,b}$ & $I_a \otimes r_b$ & $\alpha_{a,b,1} \circ r_{a \otimes b}$ \\
$R_{a,b,c}$ &$[\alpha_{b,c,a} \circ \beta_{a, b \otimes c}] \circ \alpha_{a,b,c}$ & $[(I_b \otimes \beta_{a, c}) \circ \alpha_{b,a,c} ] \circ (\beta_{a, b} \otimes I_c)$ \\
$S_{a,b,c}$ & $[\alpha_{c, a,b}^* \circ \beta_{a \otimes b, c}] \circ \alpha_{a,b,c}^*$ & $[( \beta_{a, c} \otimes I_b) \circ \alpha^*_{a, c,b}] \circ (I_a \otimes \beta_{b,c})$ \\
$\sigma_{a,b}$ & $I_{a \otimes b}$ &$ \beta_{b,a} \circ \beta_{a,b}$ \\
 \hline
\end{tabular}
\end{center}
\caption{Binary Words in Symmetric Monoidal 2-Computad.}
\label{BinWordsin2TruncatedGenDataTable}
\end{table}%
%Binary Sentences
%{\em Binary sentences} in $G_2$ are defined recursively as follows:
%\begin{itemize}
%\item If $\epsilon$ is a binary word in $G_2$ then $(\epsilon)$ is a binary sentence in $G_2$ with source $s(\epsilon)$ and target $t(\epsilon)$. 
%\item If $u$ and $v$ are binary sentences in $G_2$ such that $t(t(v)) = s(s(u))$, then $u * v$ is a binary word in $G_2$ with source $s(u) \circ s(v)$ and target $t(u) \circ t(v)$.  
%\item If $u$ and $v$ are binary sentences in $G_2$, then $u \otimes v$ is a binary sentence in $G_2$ with source $s(u) \otimes s(v)$ and target $t(u) \otimes t(v)$.
%\end{itemize}
%Paragraphs
{\em Paragraphs} in $G_2$ are defined recursively as follows:
\begin{itemize}
\item If $\epsilon$ is a binary word in $G_2$ then $(\epsilon)$ is a paragraph in $G_2$ with source $s(\epsilon)$ and target $t(\epsilon)$. 
\item  If $p$ and $p'$ are paragraphs in $G_2$, such that $t(t(p')) = s(s(p))$, then $(p) * ( p')$ is a paragraph in $G_2$ with source $s(p) \circ s(p')$ and target $t(p) \circ t(p')$. 
\item If $p$ and $p'$ are paragraphs in $G_2$,  then $(p) \otimes ( p')$ is a paragraph in $G_2$ with source $s(p) \otimes s(p')$ and target $t(p) \otimes t(p')$. 

\item If $p_0, p_1, \cdots p_k$ is a composable sequence of paragraphs in $G_2$, i.e., $s(p_{i-1}) = t(p_{i})$ for $1 \leq i \leq k$, then the (non-binary) word $p_0 p_1 \cdots p_k$ is a paragraph.  We also write this paragraph as $p_0 \circ p_1 \circ \cdots \circ p_k$.
\end{itemize}
\end{definition}

\begin{definition} \label{DefnFreeSymmetricMonBicat}
 Let $G = ( G_2, G_1, G_0)$ be a symmetric monoidal 2-computad. Then define $\sF(G) = \sF_2(G)$ as the following symmetric monoidal bicategory. The objects of $\sF(G)$ consist of the binary words in $G_0$. The 1-morphisms consist of the binary sentences in $G_1$ and the 2-morphisms consist of equivalence classes of paragraphs in $G_2$. The equivalence relation on paragraphs is the finest such that:
 \begin{itemize}
\item If $p$ is a paragraph with source $f$ and target $g$, and $p^{-1}$ is also a paragraph, then $p p^{-1} \sim id_g$ and $p^{-1} p \sim id_f$. 
\item $a^c, r^c, \ell^c, \phi^\otimes_{(f,g), (f', g')},  \phi^\otimes_{a, a'}, \alpha_{f,g,h}, \ell_f, r_f,$ and $\beta_{f,g}$ are natural. This means that for each of these morphisms a certain square commutes. For example, consider binary sentences $f, \tilde f: a \to b$ and any paragraph $\xi: f \to \tilde f$. If $r^c$ is natural then the following square commutes:
\begin{center}
\begin{tikzpicture}[thick]
	\node (LT) at (0,1.5) 	{$f \circ I_a$ };
	\node (LB) at (0,0) 	{$\tilde f \circ I_a$};
	\node (RT) at (2,1.5) 	{$f$};
	\node (RB) at (2,0)	{$\tilde f$};
	\draw [->] (LT) --  node [left] {$\xi * id_{I_a}$} (LB);
	\draw [->] (LT) -- node [above] {$r^c_f$} (RT);
	\draw [->] (RT) -- node [right] {$\xi$} (RB);
	\draw [->] (LB) -- node [below] {$r^c_{\tilde f}$} (RB);
\end{tikzpicture}
\end{center}
which is to say that the paragraphs $\xi \circ r^c_f$ and $r^c_{\tilde f} \circ (\xi * id_{I_a})$ are equivalent. Naturality for $a^c, \ell^c, \phi^\otimes_{(f,g), (f', g')},  \phi^\otimes_{a, a'}, \alpha_{f,g,h}, \ell_f, r_f,$ and $\beta_{f,g}$ is similar. 
\item $a^c, r^c, \ell^c$ satisfy the pentagon and triangle axioms. 
\item $\otimes, \alpha, \ell, r, \beta, \pi, \mu, \lambda, \rho, R, S, \sigma$ satisfy the axioms for symmetric monoidal bicategories.

That is the equations (TA1), (TA2) and (TA3) of \cite{GPS95} are satisfied, the equations (BA1), (BA2), (BA3), (BA4), (SA1) and (SA2)  of  \cite{McCrudden00} are satisfied and equation (SMA) is satisfied, see Definition \ref{DefnSymMonBicat}.
\item The equivalence relation is closed under $\otimes, *$ and concatenation $\circ$ (i.e., composition).
\end{itemize}
 Vertical composition of 2-morphisms is given by concatenation. The composition of 1-morphisms is given by $\circ$ and the horizontal composition of 2-morphisms is given by $*$.  
\end{definition}

In this work we only consider relations which occur at the level of 2-morphisms, although more general constructions do exist. Let $\sM$ be a symmetric monoidal bicategory equipped with a set  $\cR$ of pairs of 2-morphism $(\alpha,\beta) \in \cR$ such that the sources and targets of $\alpha$ and $\beta$ agree. $\cR$ is called the set of {\em generating relations}.
Let $\sim$ denote the finest equivalence relation on the 2-morphisms of $\sM$ such that:
\begin{enumerate}
\item $\alpha \sim \beta$ if $(\alpha,\beta) \in \cR$,
\item $f \circ \alpha \circ g \sim f \circ \beta \circ g$ whenever these compositions exist and $\alpha \sim \beta$,
\item $\alpha * g \sim \beta *g$  and $f * \alpha \sim f* \beta$ whenever these compositions exist and $\alpha \sim \beta$,
\item $\alpha \otimes g \sim \beta  \otimes g$  and $f \otimes \alpha \sim f \otimes \beta$ whenever $\alpha \sim \beta$.
\end{enumerate}
Let $x, y \in \sM$ be objects and $f,g: x \to y$ be 1-morphisms in $\sM$. Denote by
$(\sM/ \cR)(x,y)(f,g)$ the set $\sM(x,y)(f,g) / \sim$, and by $(\sM/ \cR)_2$ the totality $\sM_2 / \sim$. We define the following partial operations on the sets $(\sM/ \cR)(x,y)(f,g)$:
\begin{enumerate}
\item vertical composition: $[\alpha] \circ [ \beta] := [ \alpha \circ \beta]$,
\item horizontal composition:  $[\alpha] * [ \beta] := [ \alpha * \beta]$,
\item tensor product:  $[\alpha] \otimes [ \beta] := [ \alpha \otimes \beta]$.
\end{enumerate}
These operations are only given when the corresponding operation is defined on representatives. The above required properties  of $\sim$ ensure that these operations are well defined. The vertical composition makes $(\sM/ \cR)( x,y)$ into a category with the same objects as $\sM(x,y)$. Let $\sM / \cR$ denote the collection of the objects of $\sM$ together with these categories $(\sM/ \cR)( x,y)$. We equip $\sM/\cR$ with  associator and unitor natural isomorphisms induced by their images under the evident quotient map $\sM_2 \to (\sM/\cR)_2$. In this way $\sM/ \cR$ becomes a bicategory. Similarly $\otimes$, together with the images of the unit and coherence data of $\sM$, equip $\sM/ \cR$ with the structure of a symmetric monoidal bicategory. The quotient map $\sM \to \sM/ \cR$ is a strict symmetric monoidal (strict) homomorphism. 

\begin{proposition} \label{RelationsForSMBicatsProp}
The quotient homomorphism $\sM \to \sM/ \cR$ satisfies the following universal property:
\begin{equation*}
\symbicat(\sM/ \cR, \sS) \cong \symbicat_\cR(\sM, \sS) 
\end{equation*}
for all symmetric monoidal bicategories $\sS$. Here $\symbicat_\cR$ denotes the full sub-bicategory consisting of those symmetric monoidal homomorphisms $H$, such that $H(\alpha) = H(\beta)$ whenever $(\alpha, \beta) \in \cR$.  
\end{proposition}

\begin{proof}
Pre-composition with $\sM \to \sM/\cR$ clearly lands in $\symbicat_\cR(\sM,\sS)$. Given a symmetric monoidal homomorphism $H: \sM \to \sS$ in  $\symbicat_\cR(\sM,\sS)$, the properties of $\sim$ and $H$ ensure that $H(\alpha) = H(\beta)$ whenever $\alpha \sim \beta$. Hence there is, in fact, a unique extension of $H$ to $\sM/ \cR$. It is also follows, since $\sM/ \cR$ has the same objects and 1-morphisms as $\sM$, that symmetric monoidal transformations and modifications extend uniquely as well, so that we indeed have an isomorphism of bicategories. 
\end{proof}

\begin{definition}
	A {\em symmetric monoidal 3-computad} $P$ (a.k.a. a {\em symmetric monoidal presentation}) consists of a symmetric monoidal 2-computad $G$ together with a set of generating relations $\cR$ for the symetric monoidal bicategory $\sF(G)$. The symmetric monoiadal bicategory $\sF(P)$ is given by $\sF(G)/\cR$. A symmetric monoidal bicategory strictly isomorphic to $\sF(P)$ for some symmetric monoidal 3-computad will be called {\em computadic}.
\end{definition}

%
%\begin{remark} \label{rem:filtration}
%	Let $G = ( G_2 \rightrightarrows G_1 \rightrightarrows G_0)$ be a generating datum. By discarding the generating 2-morphisms or generating 2- and 1-morphisms we obtain two new generating data
%	\begin{align*}
%		G_{\leq 1} &= ( \emptyset \rightrightarrows G_1 \rightrightarrows G_0) \\
%		G_{\leq 0} &= ( \emptyset \rightrightarrows \emptyset \rightrightarrows G_0).
%	\end{align*} 
%	We have a corresponding sequence of strict functors $\sF(G_{\leq0}) \to \sF(G_{\leq 1}) \to \sF(G)$, which we can regard as a filtration. This will play an important role in section~\ref{sec:coherence}.
%\end{remark}

\section{The Cofibrancy Theorem} \label{sec:cofibrancyThm}

In this section we prove that computadic symmetric monoidal bicategories enjoy a cofibrancy property which amounts to a sort of coherence theorem for homomorphisms. In particular if $P = (G_0, G_1, G_2, \cR)$ is a symmetric monoidal 3-computad (a.k.a. presentation) and $\sM$ is an arbitrary symmetric monoidal bicategory, then every (weak) symmetric monoidal homomorphism from $\sF(P)$ to $\sM$ is naturally equivalent in $\symbicat(\sF(P), \sM)$ to a strict homomorphism. In fact we show that $\symbicat(\sF(P), \sM)$ is equivalent to an auxiliary bicategory $P(\sM)$ of ``$P$-data'' in $\sM$. 

 Roughly, the objects of $P(\sM)$ consist of an assignment of an object of $\sM$ for each element of $G_0$, a 1-morphisms in $\sM$ for each element of $G_1$, and a 2-morphism in $\sM$ for each element of $G_2$. These are required to be compatible with the source and target maps. For example suppose that $G_0 = \{ x, y, z\}$, $G_1 = \{f \}$, and that the source of $f$ is $x \otimes (y \otimes z)$. If we are given an assignment of objects of $\sM$ for each element $G_0$ (i.e., a map $G_0 \to \sM_0$), then we get a canonical extension to all binary words in $G_0$ (i.e., an extension to $BW(G_0) \to \sM_0$). This extension is given by taking the binary word in $G_0$ and evaluating it in the symmetric monoidal bicategory $\sM$. The 1-morphism that we assign to $f$ must now have a source which is exactly the value of $x \otimes (y \otimes z)$, evaluated in $\sM$. A similar statement holds for the targets of generating 1-morphisms and for 2-morphisms as well. Of course the relations $\cR$ must be satisfied. 
The morphisms and 2-morphisms of $P(\sM)$ can similarly be defined by simply assigning values to the generators of $P$.

Let $G$ be the underlying symmetric monoidal 2-computad of $P$ and $\cR$ the set of generating relations. As we will see, the bicategory $P(\sM)$ will be identified with a sub-bicategory of $\symbicat(\sF(P), \sM)$. In light of Proposition~\ref{RelationsForSMBicatsProp}, $\symbicat(\sF(P), \sM)$ is a full sub-bicategory of $\symbicat(\sF(G), \sM)$.  Namely it consists of those homomorphisms which satisfy the relations. The bicategory $P(\sM)$ will similarly be defined as the full-sub-bicategory of $G(\sM)$ which such that the relations $\cR$ are satisfied. Thus we will focus on the case that $P = G$ is a 2-computad without relations.

\begin{definition} \label{DefnBicatOfGShapedData}
Let $G = ( G_2, G_1, G_0)$ be a symmetric monoidal 2-computad and let $\sM$ be a symmetric monoidal bicategory. Define the bicategory $G(\sM)$ as follows:
\begin{itemize}
\item The objects of $G(\sM)$ consist of a triple of maps: 
\begin{equation*}
\phi: (G_0, G_1, G_2)  \to (\sM_0, \sM_1, \sM_2).
\end{equation*}
 Given such a triple, we have a canonical extension to a triple of maps: 
 \begin{equation*}
\tilde \phi :(BW(G_0), BS(G_1), P(G_2)) \to  (\sM_0, \sM_1, \sM_2),
\end{equation*}
where $BW(G_0)$ denotes the set of binary words in $G_0$, $BS(G_1)$ denotes the set of binary sentences in $G_1$, and $P(G_2)$ denotes the paragraphs in $G_2$. This extension map is given by using $\phi$ to obtain a corresponding expression in $\sM$ and evaluating it in $\sM$ using the symmetric monoidal structure on $\sM$. The triple $\phi$ is required to satisfy the conditions:
\begin{enumerate}
\item $\tilde \phi (s(f)) = s( \phi(f))$ and $\tilde \phi (t(f)) = t( \phi(f))$  for all $f \in G_1$,
\item $\tilde \phi (s(\epsilon)) = s( \phi(\epsilon))$ and $\tilde \phi (t(\epsilon)) = t( \phi(\epsilon))$  for all $\epsilon \in G_2$.
\end{enumerate}
\item The 1-morphisms of $G(\sM)$ (from $\phi_0$ to $\phi_1$, say) consist of pairs of maps, 
\begin{equation*}
\psi: (G_0, G_1) \to (\sM_1, \sM_2),
\end{equation*}
These maps are required to satisfy a number of conditions. Given such a pair, we have a canonical extension $\tilde \psi:BW(G_0) \to \sM_1$ given by evaluation in $\sM$. We require for all $a \in G_0$ and for all $f\in G_1$ that, 
\begin{align*}
s( \psi(a)) &= \phi_0(a)\\  
t(\psi(a)) &= \phi_1(a) \\
s(\psi(f)) &= \phi_1(f) \circ \tilde \psi (s(f)) \\
  t(\psi(f)) &= \tilde \psi (t(f)) \circ \phi_0(f). 
\end{align*}
 We can also extend $\psi$ to a map $\tilde \psi: BS(G_1) \to \sM_2$, also satisfying this property. This extension is defined inductively as follows:
\begin{itemize}
\item $\tilde \psi = \psi$ on the elements of  $G_1 \subset BS(G_1)$, 
\item $\tilde \psi(I_u)$ is the canonical 1-morphism: \begin{equation*}
(r^c)^{-1} \circ \ell: I_{\tilde \phi_1(u)} \circ \psi(u) \to \psi(u) \to \psi(u) \circ  I_{\tilde \phi_0(u)}. 
\end{equation*}
\item For $\ell, r, \alpha$ and $\beta$, $\tilde \psi$ is defined to be the structure 2-morphisms which realize $\ell, r, \alpha$ and $\beta$ as transformations. For example, $\tilde \psi(\alpha)$ is defined to be $\alpha_{\tilde \psi(u), \tilde \psi(v), \tilde \psi(w)}^{-1}$, which fills the following diagram:
\begin{center}
\begin{tikzpicture}[thick]
	\node (LT) at (0,1.5) 	{$(\phi_0u \otimes \phi_0v) \otimes \phi_0 w$ };
	\node (LB) at (0,0) 	{$(\phi_1u \otimes \phi_1v) \otimes \phi_1 w$};
	\node (RT) at (3.5,1.5) 	{$\phi_0u \otimes (\phi_0v \otimes \phi_0 w)$};
	\node (RB) at (3.5,0)	{$\phi_0u \otimes (\phi_0v \otimes \phi_0 w)$};
	\draw [->] (LT) --  node [left] {$\stackrel{\psi_{(u \otimes v) \otimes w} =}{ (\psi(u) \otimes \psi(v)) \otimes \psi(w)}$} (LB);
	\draw [->] (LT) -- node [above] {$\alpha$} (RT);
	\draw [->] (RT) -- node [right] {$ \psi(u) \otimes (\psi(v) \otimes \psi(w))$} (RB);
	\draw [->] (LB) -- node [below] {$\alpha$} (RB);
	\node at (2, .75) {$\Rightarrow$};
\end{tikzpicture}
\end{center}
The values of $\tilde \psi (\ell), \tilde \psi (r)$, and $\tilde \psi (\beta)$ are defined similarly, and this determines $\tilde \psi$ on binary words in $G_1$ of length one.  
\item We inductively define $\tilde \psi$ on binary words of longer length, by defining $\tilde \psi (f \otimes f')$ as the composition:
\begin{align*}
 [\phi_1f \otimes \phi_1f'] \circ [ \psi a \otimes \psi a'] &\to [\phi_1f \circ \psi a] \otimes [ \phi_1f' \circ \psi a'] \\
 & \stackrel{\tilde \psi f \otimes \tilde \psi f'}{\to}  [\psi b \circ \phi_0f] \otimes [  \psi b' \circ \phi_1f']  \\
 & \to  [\psi b \otimes  \psi b'] \circ  [   \phi_0 f \otimes \phi_1f']
\end{align*}
\item  Finally, to define $\tilde \psi$ on all binary sentences in $G_1$, we must define it on compositions.   If $f$ and $g$ are binary sentences in $G_1$ such that $s(f) = t(g)$, then $\tilde \psi( f \circ g)$ is defined by the diagram
	\begin{center}

		\end{center}
\end{itemize}
Thus we have formed an extension $\tilde \psi: BS(G_2) \to \sM_2$. For $\psi$ to be a 1-morphism of $G(\sM)$ we require a further axiom. For all $\epsilon \in G_2$ (with source binary sentence $f$ and target $g$) we require that the naturality square commutes:
\begin{center}
%
\end{center}

\item The 2-morphism of $G(\sM)$ (from $\psi: \phi_0 \to \phi_1$ to $\theta: \phi_0 \to \phi_1$) are given by certain maps $m:G_0 \to \sM_2$. These maps are required to satisfy several properties. First, we require that $s(m(a)) = \psi(a)$ and that $t(m(a)) = \theta(a)$. Such a map has a canonical extension $m: BW(G_0) \to \sM_2$ (again given by evaluation in $\sM$), and we further require that for all $f \in G_1$, the following diagram commutes: 
\begin{center}
%
\end{center}

\end{itemize}
The vertical and horizontal compositions of 2-morphisms in $G(\sM)$ are given point wise in $\sM$, and the 
composition of 2-morphisms in $G(\sM)$ is given by the composition of transformations as in Definition \ref{defncompositonoftransformations}. The associators and unitors of $\sM$ induce associators and unitors for $G(\sM)$.  
\end{definition}

The canonical extensions considered in the above definitions give the data of the homomorphisms, transformations, and modifications of $\symbicat(\sF(G), \sM)$. Indeed $\tilde \phi$, together with the trivial coherence structure $(\chi, \iota, \omega, \gamma, \delta, u)$ gives a strict symmetric monoidal homomorphism, $\tilde \psi$ (with trivial $\Pi$ and $M$) gives a strict symmetric monoidal transformation between the $\tilde \phi_i$, and $m$ gives a similar symmetric monoidal modification. In constructing $G(\sM)$ we have simultaneously built a canonical homomorphism: 
\begin{equation*}
j: G(\sM) \to \symbicat(\sF(G), \sM).
\end{equation*}
In fact, the composition in $G(\sM)$ is given by precisely the composition in $\symbicat(\sF(G), \sM)$, so that this homomorphism is strict.

\begin{proposition} \label{PropJEssentiallySurjOnObjects}
Every symmetric monoidal homomorphism $H \in  \symbicat(\sF(G), \sM)$ is equivalent to one in the image of the canonical homomorphism:  \begin{equation*}
j: G(\sM) \to \symbicat(\sF(G), \sM).
\end{equation*} 
\end{proposition}
\begin{proof} Given $H = (H, \chi, \iota, \omega, \gamma, \delta u)$ we construct an object $\res(H)$ of $G(\sM)$ together with  a canonical equivalence transformation $\theta$ between $H$ and $j \circ \res (H)$. The object $\res H$ is essentially the {\em restriction} of $H$ to the generating sets $G_0, G_1$, and $G_2$.  
 An object of $G(\sM)$ consists of three maps, and the first of these is precisely given by restricting $H$ to $G_0 \subseteq \text{ob } \sF(G)$. This gives a map: 
\begin{equation*}
	\res H : G_0 \to \sM_0.
\end{equation*}

By our previous considerations, there is a canonical extension of $\res H$ to $BW(G_0)$, which usually will not agree with $H$. Thus we cannot simply take $\res H$  to be the actual restriction of $H$ to $G_1$ and $G_2$; the sources and targets won't match up correctly, as required in Definition \ref{DefnBicatOfGShapedData}. 
However there is a canonical 1-morphism in $\sM$, $\theta_w: H(w) \to \res H(w)$, which is defined inductively as follows. On the elements $a \in G_0$, $\theta_a = I_a$ the identity 1-morphism, and $\theta_1 = \iota^*: H(1) \to 1$. This defines $\theta$ for words of length less of length one. For longer binary words, $\theta$ is defined inductively as the composition:
\begin{equation*}
\theta_{u \otimes v}: H( u \otimes v) \stackrel{\chi^*}{\to} H(u) \otimes H(v) \stackrel{\theta_u \otimes \theta_v}{\to} \res H(u) \otimes \res H(v) = \res H(u \otimes v).
\end{equation*}
Similarly $\theta_u^*: \res H (u) \to H(u)$ is defined using $\iota$ and $\chi$ in place of $\iota^*$ and $\chi^*$. 
We define $\res H: G_1 \to \sM_1$ by setting $\res H (f)$ to be the composition $\res H (f) = [\theta_{t(f)} \circ H(f)] \circ \theta_{s(f)}^*$. We also define $\res H: G_2 \to \sM_2$ by setting $\res H( \alpha) = [ id_{\theta_{t(f)}} * H(\alpha) ] * id_{\theta_{s(f)}^*}$, i.e., by whiskering. The data $\res H: (G_0, G_1, G_1) \to (\sM_0, \sM_1, \sM_2)$ now defines an object of $G(\sM)$. 

The components, $\theta_u$, are part of what will be a symmetric monoidal transformation $\theta: H \to j ( \res H)$. We must define $\theta$ on the 1-morphisms of $\sF(G)$. This is again done inductively. For any binary sentence $s: u \to v$, we must have a 2-morphism $\theta_s: \res H (s) \circ \theta_u \to \theta_v \circ H(s)$. For $f \in G_1$, we have $\res H(s) = [ \theta_v \circ H(f) ] \circ \theta^*_u$, thus we may take $\theta_s$ to be the obvious composition of associators, unitors, and the adjunction 2-isomorphism $\theta^*_u \circ \theta_u \to  I_u$. Similarly, we define $\theta$ on $I_u$ by the pasting diagram in Figure~\ref{fig:PropJEssentiallySurjOnObjects1}.
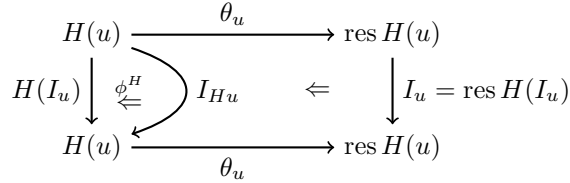
\begin{figure}[ht]
	\begin{center}
	\begin{tikzpicture}[thick]
		\node (LT) at (0,1.5) 	{$H(u)$ };
		\node (LB) at (0,0) 	{$H(u)$};
		\node (RT) at (4,1.5) 	{$\res H (u)$};
		\node (RB) at (4,0)	{$\res H (u)$};
		\draw [->] (LT) --  node [left] {$H(I_u)$} (LB);
		\draw [->] (LT) -- node [above] {$\theta_u$} (RT);
		\draw [->] (RT) -- node [right] {$I_u = \res H(I_u)$} (RB);
		\draw [->] (LB) -- node [below] {$\theta_u$} (RB);
		\draw [->] (LT).. controls (1.5,1) and (1.5, 0.5) .. node [right] {$I_{Hu}$} (LB);
		\node at (.5, .75) {$\stackrel{\phi^H}{\Leftarrow}$};
		\node at (3, .75) {$\Leftarrow$};
	\end{tikzpicture}
	\end{center}
	\caption{A Pasting Diagram}
	\label{fig:PropJEssentiallySurjOnObjects1}
\end{figure}
Here the unlabeled 2-morphism is the obvious composition of unitors. With this definition $\theta$ automatically satisfies the second commutative diagram in Figure~\ref{FigBicatTransformationAxioms}.
 
For $\alpha$, $\ell$, $r$, and $\beta$ we can use similar pasting diagrams, involving the coherence morphisms $\omega, \gamma, \delta, u$ and their mates. For example $\theta_\ell$ is defined by the pasting diagram in Figure~\ref{fig:PropJEssentiallySurjOnObjects2}.
\begin{figure}[ht]
	\begin{center}

	\end{center}
	\caption{A Pasting Diagram for $\theta_\ell$}
	\label{fig:PropJEssentiallySurjOnObjects2}
\end{figure}
The pasting diagrams for $\alpha, r$, and $\beta$ are similar and use mates of $\omega, \delta$ and $u$, respectively. 

This defines $\theta$ on 1-morphisms of length and width one. We extend it to all binary sentence by induction. If $f \circ g$ is a binary sentence in $G_1$, then we define $\theta_{g \circ f}$ to be the unique 2-morphism filling the diagram of Figure~\ref{fig:PropJEssentiallySurjOnObjects3}.
\begin{figure}[ht]
	\begin{center}
		%
		\end{center}
	\caption{A diagram defining $\theta_{g \circ f}$}
	\label{fig:PropJEssentiallySurjOnObjects3}
\end{figure}

It is instructive to compare this with Definition \ref{DefnBicatTransformation} of transformations of bicategories. Defining $\theta$ on compositions as above ensures that $\theta$ satisfies the first diagram of Figure~\ref{FigBicatTransformationAxioms}.
Finally, we must inductively define $\theta$ on tensor products $ f \otimes g$. This is given by the  pasting diagram in Figure~\ref{fig:PropJEssentiallySurjOnObjects4}.
\begin{figure}[ht]
	\begin{center}

	\end{center}
	\caption{A Pasting Diagram defining $\theta$ on $ f \otimes g$}
	\label{fig:PropJEssentiallySurjOnObjects4}
\end{figure}

With these definitions $\theta$ becomes a transformation $\theta: H \to \res H$, and similarly $\theta^*: \res H \to H$ becomes an inverse transformation. We will now promote $\theta$ to be a symmetric monoidal transformation by equipping it with modifications $\Pi$ and $M$. 
For $M$, we must have a 2-morphism which fills the diagram in Figure~\ref{fig:PropJEssentiallySurjOnObjects5}.
\begin{figure}[ht]
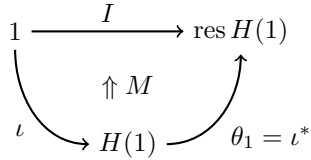

	\begin{center}
	%
	\end{center}
	\caption{A Pasting Diagram for $M$ }
	\label{fig:PropJEssentiallySurjOnObjects5}
\end{figure}
The obvious choice is the adjunction 2-morphism for $\iota$, and this is what we will take. Similarly we use the adjunction 2-morphisms for $\chi$ to define $\Pi$, which is given precisely by the  pasting diagram in Figure~\ref{fig:PropJEssentiallySurjOnObjects6}.
\begin{figure}[ht]
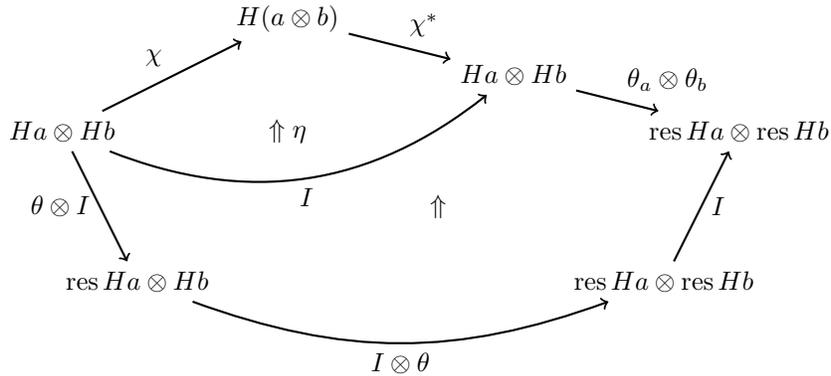

	 \begin{center}
	%
	\end{center}
	\caption{A Pasting Diagram}
	\label{fig:PropJEssentiallySurjOnObjects6}
\end{figure}
Here $\eta$ is an adjunction 2-morphism for $\chi$ and the unlabeled 2-morphism is the canonical one given by the homomorphism structures of $\otimes$. 

With these choices of $\Pi$ and $M$, $\theta$ becomes a symmetric monoidal transformation. The necessary coherence diagrams commute precisely because of our choices for $\theta$ and the fact that $\omega, \gamma, \delta$, and $u$ come from a symmetric monoidal homomorphism. A similar statement holds for $\theta^*$, which is easily seen to be an inverse equivalence transformation.
\end{proof}

In the proof of the last proposition, we constructed an object $\res H \in G(\sM)$ for every homomorphism $H \in \symbicat(\sF(G), \sM)$.  In fact this fits into a  homomorphism:
\begin{equation*}
\res : \symbicat(\sF(G), \sM) \to G(\sM),
\end{equation*}
 which satisfies $\res  \circ j = id_{G_{\sM}}$, so that we may regard $j$ as an inclusion. This follows from the following proposition.

\begin{theorem}[The Cofibrancy Theorem] \label{thm:cofibrancythm}
The homomorphisms 
\begin{equation*}
	j:G(\sM) \leftrightarrows \symbicat(\sF(G), \sM): \res
\end{equation*}
 are inverse equivalences of bicategories; The composite $\res \circ j = id$ is the identity and we have a natural equivalence $\theta: id \to j \circ \res$ as defined above. 
\end{theorem}

\begin{proof} Proposition~\ref{PropJEssentiallySurjOnObjects} asserts that the homomorphism $j$ is essentially surjective on objects. By Theorem~\ref{Whiteheadforbicats} (see also Theorem~\ref{WhiteheadforSymMonBicats}), we must show that it is also fully-faithful on 2-morphisms and essentially full on 1-morphisms. 

	Recall that the objects and 1-morphisms in the image of $j$ are strict,  i.e., the higher coherence data $(\chi, \iota, \omega, \gamma, \delta, u)$ and $(M, \Pi)$ are trivial. Let  $H$ and $\overline H$ be two homomorphisms in the image of $j$ and let $\theta$ and $\tilde \theta$ be two transformations between these, also in the image of $j$. Then axioms (BMBM1) and (BMBM2) of Definition \ref{DefnSymMonoidalModifiaction} ensure that any modification $m: \theta \to \tilde \theta$ is in fact determined by its restriction to $G_0 \subseteq \sF(G)$, and agrees with $j$ applied to this restriction. Thus $j$ is fully-faithful on 2-morphisms. 
	
	It remains to show, given $H$ and $\overline H$ in the image of $j$, that any symmetric monoidal transformation $(\theta, \Pi, M)$ between them is equivalent to one in the image of $j$. Since $H$ and $\overline H$ are in the image of $j$, the restriction of $\theta$ to $G_0$ and $G_1$ is an actual 1-morphism of $G(\sM)$ (between the corresponding objects of $G(\sM)$). Call this $\res \theta$. 
We will show $j \res \theta$ is equivalent to $\theta$. 

We must construct an invertible symmetric monoidal modification $m:j \res \theta \to \theta$. This consists of 2-morphisms, 
\begin{equation*}
	m_u:  j\res \theta (u) \to \theta(u)
\end{equation*}
for each object $u \in BW(G_0)$, which must satisfy certain coherence conditions. In fact $\Pi$ and $M$ provide a canonical choice for $m_u$. 

Recall that $j \res \theta(u)$ is an iterated binary tensor product of $\theta$ restricted to $G_0$. Since $H$ and $\overline H$ are strict, $\Pi$ provides 2-morphisms $ \theta_a \otimes \theta_b \to \theta_{a \otimes b}$ for any pair of objects $a, b \in \F(G)$, and $M$ provides a 2-morphism $\theta_1 \to I_1$. Thus given any binary word $u$ in $G_0$ we may apply a sequence of these 2-morphisms $\Pi$ and $M$ to $j \res \theta(u)$ to obtain $ \theta (u)$. Any two choices of  sequence result in the same 2-morphism, as can be proven by appealing to  Lemma \ref{AnyPathsofTreesareEquivalent}, (an argument  we have demonstrated several times before). The fact that this is an (invertible) symmetric monoidal modification then follows readily from the coherence axioms of $\Pi$ and $M$. 
\end{proof}

The interpretation of the above proposition is that the symmetric monoidal bicategory $\sF(G)$ is the {\em free symmetric monoidal bicategory generated by $G$}. The bicategory of symmetric monoidal homomorphisms  from $\sF(G)$ into any other symmetric monoidal bicategory, $\sM$, is equivalent to the bicategory of $G$-data in $\sM$. 

%\section{Relations for Symmetric Monoidal Bicategories} \label{SectRelnForSymMonBicat}

\begin{corollary}
	If $\sC$ is a symmetric monoidal bicategory which is cofibrant in the model structure on symmetric monoidal bicategories (See Rmk.~\ref{rmk:symbicatmodelstr} and Prop.~\ref{pro:TalgModelStr}), and $\sM$ is an arbitrary symmetric monoidal bicategory, then every symmetric monoidal homomorphism from $\sC$ to $\sM$ is equivalent to a strict symmetric monoidal homomorphism. 
\end{corollary}  

\begin{proof}
	The proof sketch following Prop.~\ref{pro:TalgModelStr} shows that every cofibrant symmetric monoidal bicategory is a strict retract of a computadic symmetric monoidal bicategory, and so the corollary follows directly from the above cofibrancy theorem. 
\end{proof}

\begin{remark}
	In light of Corollary~\ref{cor:leftadjointmorphismoftheories}, the above theorem holds also when $\sC$ and $\sM$ are replaced by unbiased semistrict symmetric monoidal 2-categories or by quasistrict symmetric monoidal 2-categories. 
\end{remark}

\section{Presentations of Stricter Symmetric Monoidal 2-Categories} \label{sec:PresSemistrict}

In Sections~\ref{SectSmyMonBicat} and \ref{sec:strictsymbicats} we introduced four different notions of symmetric monoidal bicategory of varying degrees of strictness. The most weak theory introduced was that of symmetric monoidal bicategories and the most strict that of quasistrict symmetric monoidal 2-categories. Moreover each of these theories was shown to be algebraic and to correspond to a finitary monotone monad on 2-dimensional globular sets (see Cor.~\ref{cor:finitary-monads} and Lem.~\ref{lem:symmontheoriesaremonotone}). The general theory of {\em computads} explained in Section~\ref{sec:computads} applies in each case, and gives a theory of presentations for each notion of symmetric monoidal structure. 

In Section~\ref{SectFreelyGenSymMonBicats} we gave a completely explicit and detailed account of the theory of computads/presentations in the fully-weak case of symmetric monoidal bicategories. A completely analogous treatment applies in each of the remaining cases, except that now the operations and equivalence relations are dictated by the corresponding theory of Crans semistrict symmetric monoidal 2-categories, unbiased semistrict symmetric monoidal 2-categories, or quasistrict symmetric monoidal 2-categories. 
 In this section we will describe explicitly the resulting computadic versions of these 2-categories, paying special attention to the quasistrict and unbiased semistrict cases.

\subsection{Crans Semistrict Symmetric Monoidal 2-Categories} \label{sec:CSScomputad}

Computadic Crans semistrict symmetric monoidal 2-categories are constructed analogously to computadic symmetric monoidal bicategories, which were described in Section~\ref{SectFreelyGenSymMonBicats}. Since that section went into complete (and perhaps even painful) detail in that case, and since this case bares much similarity to that case, we will allow ourselves to be more brief, and trust the reader to fill in the necessary details. Just as before the objects of a computadic Crans semistrict symmetric monoidal 2-category will be words in the generating objects, the 1-morphisms will be sentences formed using the generating 1-morphisms, and the 2-morphisms will be paragraphs built using the generating 2-morphisms.

However there are a few very important differences. The first notable difference occurs because the tensor operation $\otimes$ and the horizontal composition operations are both strictly associative. Thus as before a 0-computad is a set $G_0$ of generating objects, but the objects of the corresponding computadic Crans semistrict symmetric monoidal 2-category are not binary words in $G_0$, as they were for computadic symmetric monoidal bicategories, but just ordinary words in $G_0$ with no choice of parentheses. 

The sentences which make up the 1-morphisms of a computadic Crans semistrict symmetric monoidal 2-category are defined just as they are for computadic symmetric monoidal bicategories, except now when we form words and sentences we do not keep track of the parentheses. Moreover we omit the generating morphisms from Table~\ref{BinWordsin1TruncatedGenDataTable} which correspond to the associator and unitors of the $\otimes$ operation. 

Another important difference occurs because the operation $\otimes$ arrises from a functor which is {\em cubical}. This can be interpreted in two equivalent ways. On the one hand, this can be interpreted as identifying certain sentences, and so, in contrast to the symmetric monoidal bicategorical case, the 1-morphisms of a computadic Crans semistrict symmetric monoidal 2-category are not sentences, but rather {\em equivalence classes} of sentences. Part of the equivalence relation allows us to insert identity 1-morphisms as we wish, or to tensor by the identity object as we wish. An equivalent interpretation is that we would only permit sentences of a particular form, which we explain momentarily. 

The key fact is that the interchanger coherence cell $\phi^\otimes$ is sometimes an identity 2-morphism, which implies the following lemma, valid in any Gray monoid:
\begin{lemma} \label{lem:coherencegraymonoid}
	In a Gray-monoid $(\sM, \otimes, 1)$ the tensor product functor
	\begin{equation*}
		\otimes: \sM(a,b) \times \sM(a', b') \to \sM(a \otimes a', b \otimes b')
	\end{equation*}
	is determined by the horizontal composition of $\sM$ together with the operation of tensoring by identity morphisms, according to the formula
	\begin{equation*}
		f \otimes f'  =  (id_b \otimes f') \circ (f \otimes id_{a'}).
	\end{equation*}
	for 1-morphisms $f: a \to b$ and $f': a' \to b'$.
\end{lemma}
\begin{proof}
	For 1-morphisms $f: a \to b$ and $f': a' \to b'$ we have 
	\begin{equation*}
		f \otimes f'  = (id_b \circ f) \otimes (f' \circ id_{a'}) = (id_b \otimes f') \circ (f \otimes id_{a'}).
	\end{equation*}
	where the second equality follows from the fact that $\otimes$ is cubical (and hence $\phi^\otimes_{id, f, f', id}$ is an identity morphism).  
	There is a similar formula for 2-morphisms. 
\end{proof}

\begin{corollary}
	In a computadic Crans semistrict symmetric monoidal 2-category, the objects are given by words in the generating set $G_0$ and  every 1-morphism is either an identity or is given as finite (horizontal) composite of elementary 1-morphisms of the form $(id_a \otimes f \otimes id_b)$ for words $a$ and $b$. Here $f$ is either a generating 1-moprhism from the set $G_1$ or $f$ is a 1-morphism of the form $\beta_{x,y}$ for words $x$ and $y$. \qed
\end{corollary}

The 2-morphisms of a computadic Crans semistrict symmetric monoidal 2-category are defined similarly using paragraphs without parentheses and using only the structural morphisms of Crans semistrict symmetric monoidal 2-categories, and subject to the relations thereof and any additional relations $\cR$ specified by the computad. 

\subsection{Quasistrict Symmetric Monoidal 2-Categories} \label{sec:qsymcomputad}

We will now give a more detailed analysis of what happens in the case of quasistrict symmetric monoidal 2-categories, which are the most strict version we have considered. 
Let $G_0$ be a set. We define a (permutative) symmetric monoidal category $\Sigma_{G_0}$ as follows: The objects of $\Sigma_{G_0}$ consist of finite words (possibly empty) with letters being the elements of $G_0$. The morphisms consist of permutations taking one word to another. More precisely the morphisms are pairs $(w_1 w_2 \dots w_n, \sigma)$ of a word in $G_0$ of length $n$ and an element $\sigma \in \Sigma_n$ of the symmetric group. The target of such a morphism is the word $w_{\sigma(1)} w_{\sigma(2)} \dots w_{\sigma(n)}$, and composition is the obvious operation. Concatenation makes $\Sigma_{G_0}$ into a strict monoidal category, which is symmetric with the obvious braiding. 

\begin{lemma}
	Let $G_0$ be a set, then $\sF_{qs}(G_0)$ the free quasistrict symmetric monoidal 2-category generated by $G_0$ is isomorphic to the symmetric monoidal category $\Sigma_{G_0}$, viewed as a quasistrict symmetric monoidal 2-category with only identity 2-morphisms.
\end{lemma}

\begin{proof}
	In the algebraic theory of quasistrict symmetric monoidal 2-categories the only non-identity 1-morphisms that are postulated are tensor products and composites of existing 1-cells and the components of $\beta$ on objects. Similarly the only non-identity 2-cells that are postulated are tensor products and composites of already exiting 2-cells and the coherence cells of the functor $\otimes$ and the transformation $\beta$.  Thus the lemma is established if we can show that for the free quasistrict symmetric monoidal 2-category generated by a set, the only 2-morphisms are identity 2-cells. Moreover it suffices to consider just the coherence cells of the functor $\otimes$ and the transformation $\beta$.
	
	The axioms of quasistrict symmetric monoidal 2-categories ensure that in this case these coherence morphisms are also identities. For the coherences cells of $\otimes$ this follows from the fact that $\otimes$ satisfies properties (3) and (5) of Definition~\ref{def:quasistrict}. Moreover Lemma~\ref{lma:QSbraidinggivenbyinterchange} states that the braiding coherence cell is determined by the interchanger. Since the later is an identity morphism it follows that the naturator of the braiding must also be an identity. 
\end{proof}

Let $(G_1, G_0)$ be a quasistrict symmetric monoidal 1-computad. This consists of a set $G_0$, a set $G_1$ and a function assigning to each element of $G_1$ a source and target object of $\sF_{qs}(G_0) = \Sigma_{G_0}$. These objects are words in $G_0$. From this data we obtain a quasistrict symmetric monoidal 2-category $\sF_{qs}(G_1,G_0)$ following the prescription of Section~\ref{sec:computads}. The objects are the same as those of $\Sigma_{G_0}$: words in the elements of $G_0$. The 1-morphisms are generated by the 1-morphisms of $\Sigma_{G_0}$ and the elements of $G_1$ under the operations of horizontal composition $\circ$ and tensor product $\otimes$. There are also 2-morphisms which are generated under the operations of tensor product, vertical, and horizontal composition  by identities on the 1-morphisms together with 2-morphisms $\phi^\otimes_{(f,f'), (g,g')}$ and $\beta_{f,g}$ coming from the coherence morphisms of the functor $\otimes$ and transformation $\beta$. 

The situation is completely analogous to the construction of the free computadic symmetric monoidal bicategory constructed previously in Section~\ref{SectFreelyGenSymMonBicats}, except that now these structures are required to satisfy the axioms of a quasistrict symmetric monoidal 2-category. 
This simplifies the situation significantly. 

Our next goal is to give a concise description of the quasistrict symmetric monoidal 2-category $\sF_{qs}(G_1,G_0)$. By Lemma~\ref{lem:coherencegraymonoid} we can always write tensor products of morphisms in terms of simpler horizontal composites of the individual factors tensored with identity morphisms on the left and right. Furthermore we have:

\begin{lemma}\label{lem:braidelementgen}
	In a quasistrict symmetric monoidal 2-category the following identity holds for any object $x$ and 1-morphism $f$: 
	\begin{equation*}
		(f \otimes id_x) = \sigma_1 \circ (id_x \otimes f) \circ \sigma_0
	\end{equation*}
	where $\sigma_0: s(f) \otimes x \to x \otimes s(f)$ and $\sigma_1: x \otimes t(f) \to t(f) \otimes x$ are the braiding isomorphisms. 
\end{lemma}
\begin{proof}
	This is a direct consequence of the fact that $\beta_{id_x,f}$ is the identity morphism. 
\end{proof}

\begin{corollary} \label{cor:QSequivreln1-mor}
	Let $(G_1, G_0)$ be a quasistrict symmetric monoidal 1-computad. Then every 1-morphism in $\sF(G_1,G_0)$ may be written as a composite
	\begin{equation*}
		\alpha_1 \circ (id_{u_n} \otimes f_n) \circ \sigma_{n-1} \circ \cdots \circ \sigma_1 \circ (id_{u_1} \otimes f_1) \circ \alpha_0
	\end{equation*}
	where the $u_i$ are objects of $\Sigma_{G_0}$ (words in $G_0$), the $\sigma_i$ and $\alpha_i$ are morphisms of $\Sigma_{G_0}$, and the $f_i$ are elements of $G_1$. We will refer to the number of $f_i$ factors as the {\em length} of the morphism. 
	
Two such decompositions 
	\begin{align*}
		\alpha_1 \circ (id_{u_n} \otimes f_n) \circ \sigma_{n-1} \circ \cdots \circ \sigma_1 \circ (id_{u_1} \otimes f_1) \circ \alpha_0 \\
		\beta_1 \circ (id_{v_n} \otimes g_n) \circ \tau_{n-1} \circ \cdots \circ \tau_1 \circ (id_{v_1} \otimes g_1) \circ \beta_0
	\end{align*}
	give rise to the same 1-morphisms if and only if they have the same length, if $f_i = g_i$, and if there exist $\phi_i: u_i \to v_i$ morphisms in $\Sigma_{G_0}$ such that
	\begin{align*}
		(id_{s(f_1)} \cdot \phi_1 ) \circ \alpha_0 &= \beta_0 \\
		\beta_1 \circ (id_{s(f_1)} \cdot \phi_n) & = \alpha_1 \\
		(id_{s(f_{i+1})} \cdot \phi_{i+1}) \circ \sigma_i &= \tau_i \circ (id_{t(f_i)} \cdot \phi_i).
	\end{align*}
\end{corollary}

\begin{proof}
	The existence of such a decomposition follows immediately by Lemma~\ref{lem:coherencegraymonoid}, Lemma~\ref{lem:braidelementgen}, and Lemma~\ref{lem:qsisuss}. Moreover the `if' part of the second claim is also clear, by applying Lemma~\ref{lem:qsisuss}. To see the `only if' part, we simply have to observe that the only relations between 1-morphism come from identifying certain coherence cells in the notion of Crans symmetric monoidal bicategory with identity cells. We can check that in each case the above relations are sufficient to relate each side of the coherence cell (in particular the order of the elementary morphisms $f_i$ is never altered).  
\end{proof}

The above corollary gives a precise normal form for all the 1-morphisms of $\sF_{qs}(G_1,G_0)$. Moreover it is straightforward to write the effect of both horizontal composition and tensor product in terms of this normal form. The former is given by concatenation and simplification, the later is given by the former plus the operation of tensoring by objects. 
 We will now move on to the 2-morphisms. 

The 2-morphisms of $\sF(G_1,G_0)$ are generated by the coherence cells  $\phi^\otimes_{(f,f'), (g,g')}$ and $\beta_{f,g}$ and subject to the relations of a quasistrict symmetric monoidal 2-category. In light of Lemma~\ref{lma:QSbraidinggivenbyinterchange}, it is sufficient to provide just the coherence cells $\phi^\otimes_{(f, id), (id, g)}$, and by the naturality of these cells and the normalization axioms of quasistrict symmetric monoidal 2-categories, it is sufficient to provide these cells just in case $f$ and $g$ have length one. All other 2-morphisms are obtained as horizontal and vertical composites of these cells and identity 2-cells. Moreover the only defining relation for symmetric monoidal bicategories which involves these coherence cells are the naturality hexagons implicit in the assertion that $\otimes$ is a homomorphism. Thus we have established:

\begin{lemma}
	Let $(G_1,G_0)$ be a quasistrict symmetric monoidal 1-computad, then the 2-morphisms of $\sF_{qs}(G_1,G_0)$ are generated under the operations of horizontal and vertical composition in a 2-category, and the operation of tensoring with identity 2-morphisms, by coeherence cells
	\begin{equation*}
		\phi^\otimes_{(f, id), (id,g)}: (id_{b'} \otimes f) \circ \beta_{a, b'} \circ (id_a \otimes g) \circ \beta_{b,a} \to \beta_{a',b'} \circ (id_{a'} \otimes g) \circ \beta_{b,a'} \circ (id_b \otimes f)
	\end{equation*}
	where $f: a \to a'$ and $g: b \to b'$ are in $G_1$, and the $\beta$ are the braiding 1-morphisms. Moreover these are only subject to the relations implicit in strict 2-categories (interchange, associativity, and units) together with the hexagon equations asserting that $\phi^\otimes$ forms the coherence cell of a homomorphism $\otimes$. \qed 
\end{lemma}

We may now describe the general case. Let $P=(\cR, G_2, G_1, G_0)$ be a quasistrict symmetric monoidal 3-computad. Then the quasistrict symmetric monoidal 2-category $\sF_{qs}(P)$ presented by $P$ has the same objects and 1-morphisms as $\sF_{qs}(G_1,G_0)$. The 2-morphisms are generated under the horizontal and vertical composition by the 2-morphisms  $\phi^\otimes_{(f, id), (id,g)}$ from above together with 2-morphisms $(\id_u \otimes \alpha)$ for $u$ an object and $\alpha \in G_2$. These are subject to the relations defining quasistrict symmetric monoidal 2-categories and the additional relations $\cR$.

\subsection{Unbiased Semistrict Symmetric Monoidal 2-Categories} \label{sec:presentationUSS}

In this case we could proceed as in the previous two cases, mimicking the constructions from Section~\ref{SectFreelyGenSymMonBicats}, however the theory of computads based on unbiased semistrict symmetric monoidal 2-categories has pleasant geometrical interpretation. In fact the string diagram calculus for unbiased semistrict symmetric monoidal 2-categories, which was developed in Section~\ref{sec:strictsymbicats}, is directly related to the structure of the computads for this theory. 

Let $P= (G_0, G_1, G_2, \cR)$ be an unbiased semistrict symmetric monoidal 3-computad. The objects and 1-morphism of $\sF(P)$ are similar to those for computadic Crans semistrict symmetric monoidal 2-categories. Specifically the objects are words in the set of generating objects $G_0$. Thus each object has a length. The 1-morphisms are either identities on the objects or are given by finite composable sequences of elementary 1-morphisms. The elementary 1-morphisms are of two varieties:
\begin{itemize}
	\item 1-morphisms of the kind $\beta^\sigma: w \to w'$. Hence $w = (w_i)$ is a word in $G_0$ of length $n$, $\sigma \in \Sigma_n$ is a permutation, and $w'= (w_{\sigma(i)})$ is the result of applying the permutation $\sigma$ to $w$. 
	\item 1-morphisms of the kind $(id_w \otimes f \otimes id_{w'})$ where $f \in G_1$ and $w,w'$ are objects. The source of this 1-morphism is $w s(f) w'$ and the target is $w t(f) w'$. 
\end{itemize}

These 1-morphisms have an obvious strictly associative and unital (horizontal) composition, given by concatenation. We can tensor these 1-morphisms by identity 1-morphisms on either the left or right sides. For the first kind of elementary 1-morphism this occurs by stabilizing the permutation, while this operation is obvious for the second kind of elementary 1-morphism. We tensor a composable chain by an identity morphism by tensoring each factor separately. Finally the tensor product of arbitrary 1-morphisms is determined by these operations, tensoring with identities and horizontal composition, in light of Lemma~\ref{lem:coherencegraymonoid}. This gives a complete description of the 1-morphisms. 

\begin{figure}[htbp]
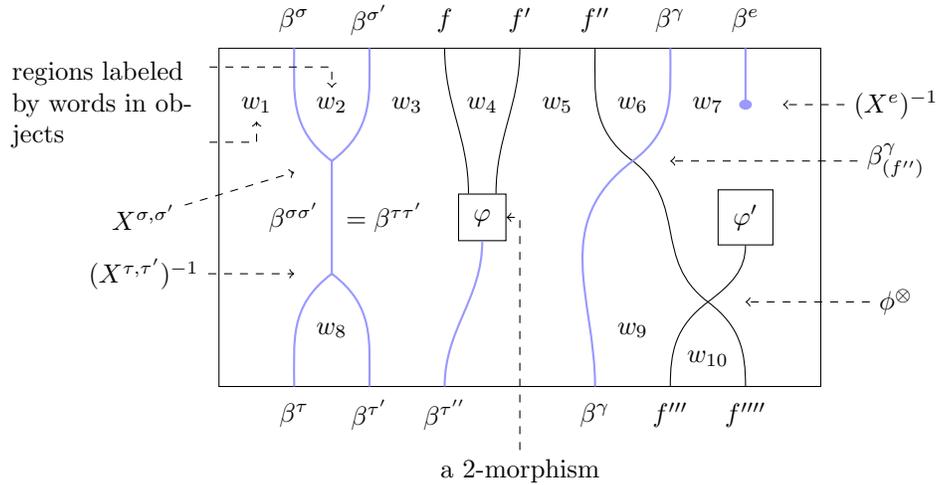

	\begin{center}

	\end{center}
	\caption{The 2-morphisms of computadic unbiased semistrict symmetric monoidal 2-categories are given by string diagrams.}
	\label{fig:USScomputad2Mor}
\end{figure}

The 2-morphisms are given by equivalence classes of string diagrams. These string diagrams are precisely the string diagrams for unbiased semistrict symmetric monoidal 2-categories that was introduced in Section~\ref{sec:strictsymbicats}. For convenience we reproduce Figure~\ref{fig:USS-stringdiagram} in Figure~\ref{fig:USScomputad2Mor}, which shows a typical string diagram of this sort. 

These string diagrams consist of a number of labeled  coupons, strings, and regions. The regions are labeled with the objects of the computadic unbiased semistrict symmetric monoidal 2-category, that is words in the elements $G_0$. There are two kinds of strings corresponding to the two kinds of elementary 1-morphisms. In Figure~\ref{fig:USScomputad2Mor} these are distinguished as the chambering and black strings. The chambering strings are labeled by the 1-morphisms of the form $\beta^\sigma$, and the black strings are labeled by the 1-morphisms of the form $(id_w \otimes f \otimes id_{w'})$. 

The coupons that are allowed in the string diagram are of three kinds:
\begin{itemize}
	\item coupons coming from the structure 2-morphisms $X^{\sigma, \sigma'}$, $X^e$, or there inverses (these are the trivalant and univalent vertices among the chambering strings); 
	\item coupons coming from the coherence cells of the $\beta^\sigma$ morphisms with respect to the $(id_w \otimes f \otimes id_{w'})$ morphisms, and their inverses (these are crossing points between chambering and black strings);
	\item coupons labeled by $id_w \otimes u \otimes id_{w'}$ where $u \in G_2$ is a generating 2-morphism. 
\end{itemize}
The vertical composition of such diagrams is given by vertical concatenation, and the horizontal composition by horizontal concatenation. The operation of tensoring with the identity extends the operation defined on 1-morphisms in the obvious way, and the tensor product of diagrams is determined by these other operations, see Figure~\ref{Fig:MonoidalStructureUSS-StringDiagram} and the discussion in Section~\ref{sec:strictsymbicats}.

These string diagrams are taken up to equivalence which is generated by local operations. These operations are precisely the axioms governing the theory of unbiased semistrict symmetric monoidal 2-categories, which express several properties, such as the naturalily of the $\beta^\sigma$ morphisms and the axioms (USS.1) for the $X$ 2-morphisms. Many of these local relations were discussed previously in Section~\ref{sec:strictsymbicats} while describing operations on general string diagrams for unbiased semistrict symmetric monoidal 2-categories. In addition there may be relations $\cR$ which are part of the 3-computad.

\section{Coherence Theorems} \label{sec:coherence}

Coherence theorems come in many varieties. MacLane's original coherence theorem states that every diagram in a bicategory of a certain form commutes. Other kinds of coherence theorems relate instance of a weak categorical structure to instances of a stricter version of that categorical structure. We will primarily be interested in two kinds coherence theorems. A coherence theorem of the first kind states for every instance of a weak categorical structure there exists an equivalent stricter categorical structure. Many theorems of this sort also specify a functorial construction of the stricter categorical object. A coherence theorem of the second kind considers categorical objects which are constructed from some sort of presentation data. These coherence theorems state that the comparison functor from the weak categorical object generated by a given datum to the stricter categorical object generated by the same datum is an equivalence. 

In this section we will discuss various coherence results of both of these kinds. The first result, a sort of extended warm-up, will be a coherence theorem of the second kind relating unbiased semistrict symmetric monoidal bicategories and quasistrict symmetric monoidal bicategories. Later, we will prove a coherence theorem of the second kind for symmetric monoidal bicategories, where the presentation data is precisely the data of a symmetric monoidal 3-computad $P = (G_0, G_1, G_2, \cR)$ as in Section~\ref{SectFreelyGenSymMonBicats}, and the stricter notion is that of quasistrict symmetric monoidal 2-category (this was the strictest notion introduced in Section~\ref{sec:strictsymbicats}). This later result will build upon the prior coherence results of Gurski \cite{MR2770448} and Gurski-Osorno \cite{GO13} and implies several other coherence theorems.

Our first theorem is the following:

\begin{theorem}\label{thm:USSCoherence}
	Let $P$ be an unbiased semistrict symmetric monoidal 3-computad. Then the canonical functor (see Section~\ref{sec:computads})
		\begin{equation*}
			\phi: \sF_{uss}(P) \to f^* \sF_{qs}(f_* P),
		\end{equation*}
		which is a strict homomorphism of unbiased semistrict symmetric monoidal 2-categories, induces an equivalence of underlying bicategories. Here  $\sF_{uss}(P)$ is the computadic semistrict symmetric monoidal 2-category generated by $P$, and $\sF_{qs}(f_* P)$ is the computadic quasistrict symmetric monoidal 2-category generated by $P$.  
\end{theorem}

\begin{corollary}%\label{cor:}
	The homomorphism $\phi: \sF_{uss}(P) \to f^* \sF_{qs}(f_* P)$ induces an equivalence of underlying symmetric monoidal bicategories. 
\end{corollary}

\begin{proof}
	The homomorphism $\phi$ induces a (strict) symmetric monoidal homomorphism between underlying symmetric monoidal bicategories. By Theorem~\ref{thm:USSCoherence} the underlying homomorphism of bicategories is an equivalence, and hence it is a symmetric monoidal equivalence by Theorem~\ref{WhiteheadforSymMonBicats} (Whitehead's Theorem for Symmetric Monoidal Bicategories).
\end{proof}

\begin{proof}[Proof of Theorem~\ref{thm:USSCoherence}]
	Let $P = (G_0, G_1, G_2, \cR)$ be an unbiased symmetric monoidal 3-computad.
	We will construct a homomorphism $H: \sF_{qs}(f_* P) \to \sF_{uss}(P)$ and a transformation $\eta: id \to H \phi$, such that the composite $\phi H = id$ is the identity (on the nose) and $\eta$ is a natural equivalence. We can describe this functor explicitly using the descriptions of $\sF_{uss}(P)$ and $\sF_{qs}(f_*P)$ given in Section~\ref{sec:PresSemistrict}. 
To construct $H$ and the transformation $\eta$, we must assign certain data in $\sF_{uss}(P)$ to the objects, 1-morphisms, and 2-morphisms of $\sF_{qs}(f_*P)$, which we will also describe explicitly.

First let us briefly recall, from Section~\ref{sec:PresSemistrict}, the description of the computadic 2-categories $\sF_{uss}(P)$ and $\sF_{qs}(f_*P)$. Both of these 2-categories have the same objects, which are the words in the set of generating objects $G_0$. The morphisms of $\sF_{uss}(P)$ consists of composable sequences of elementary 1-morphisms. These elementary 1-morphisms are either of the kind $\beta^\sigma$, corresponding to a permutation $\sigma \in \Sigma_n$, or they are of the form $(id_x \otimes f \otimes id_y)$ for some objects $x$ and $y$ and a single generating 1-morphism $f \in G_1$. 

The 1-morphisms of $\sF_{qs}(f_* P)$ are similar, except that for the second kind of elementary morphism only those of the form $(id_x \otimes f)$ are used. Moreover, the only sequences which occur are those without repetition of the permutation morphisms. This is because in $\sF_{qs}(f_*P)$ we have strict equality $\beta^\sigma \circ \beta^{\sigma'} = \beta^{\sigma \sigma'}$. Finally, these sequences are subject to an equivalence relation described in Cor.~\ref{cor:QSequivreln1-mor}, and which we will review below.

 The 2-morphisms of $\sF_{uss}(P)$ are given precisely by the equivalence classes of string diagrams, which we described at the end of Section~\ref{sec:PresSemistrict}. The elementary components of these string diagrams consist of the structure morphisms $X^{\sigma, \sigma'}$, $X^e$, and their inverses, the coherence morphisms $\beta^\sigma_{(f)}$,  $\phi^{\otimes}_{(f,id)(id, g)}$, and their inverses, and the generating 2-morphisms $\varphi \in G_2$.  The 2-morphisms of $\sF_{qs}(f_* P)$ are generated from the coherence cells $\phi^{\otimes}_{(f,id)(id, g)}$ and the generating 2-morphisms $\varphi \in G_2$.

Now we may describe the functor $\phi$ explicitly. On objects, it is the identity map. On 1-morphisms it preforms a kind of `evaluation'. Specifically, since $\phi$ is a strict functor it is sufficient to specify the image of each of the generating morphisms. Morphisms corresponding to permutations are sent to the corresponding permuation; morphisms of the kind $(id_x \otimes f \otimes id_y)$ are sent to $(id_x \sqcup \sigma_1) \circ (id_{xy} \otimes f) \circ (id_x \sqcup \sigma_0)$, as in Lemma~\ref{lem:braidelementgen}. Globally this has the effect of preforming this second replacement on each generating morphism
of this sort and then collecting together and composing the terms corresponding to permutations. On 2-morphisms there is a similar evaluation. The 2-morphisms $X^{\sigma, \sigma'}$, $X^e$, and their inverses are sent to identities, as are the coherence cells for the $\beta$ morphisms which only involve a single generator $f \in G_1$. The $\phi^\otimes$ morphisms are mapped to each other in the obvious way. 

Now we turn to the task of defining the homomorphism $H$. On objects the homomorphism $H$ is defined to also be the identity. On 1-morphisms the homomorphism $H$ is slightly more difficult, and we must make some non-canonical choices. The 1-morphisms of $\sF_{qs}(f_* P)$ are given by equivalence classes of certain composable sequences of elementary 1-morphisms. If we were working directly with these sequences, and not the equivalence classes, then it would be clear what to do. These sequences are particular instances of the ones which define 1-morphisms in  $\sF_{uss}(P)$, and so we would simply declare the value of $H$ on such a sequence to be the corresponding sequence in $\sF_{uss}(P)$. 
However the morphisms in $\sF_{qs}(f_* P)$ are not the sequences themselves, but only equivalence classes thereof. To define the homomorphism $H$ we must choose, once and for all, a representative of each equivalence class defining each 1-morphism in the 2-category $\sF_{qs}(f_* P)$. Having made this choice, we simply declare the value of $H$ on a given 1-morphism to be precisely this sequence, viewed as a 1-morphism in $\sF_{uss}(P)$. It is clear that defined this way $\phi \circ H$ is the identity on objects and 1-morphisms. 

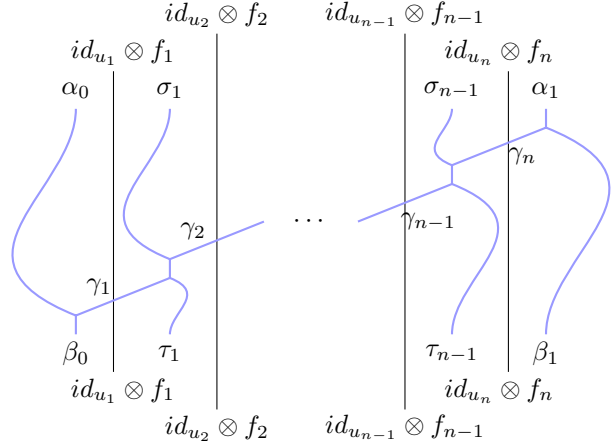
\begin{figure}[htbp]
	\begin{center}
	\begin{tikzpicture}[yscale=0.5, xscale = 1.25]
			%\draw (0,0) rectangle (6,4);
			\node at (3,3) {$\cdots$};
			
			\draw (0.9,7) -- (0.9,-1);
			\draw (2,8) -- (2,-2);
			\draw (4,8) -- (4,-2);
			\draw (5.1,7) -- (5.1,-1);
			
			\draw [chambering] (0.5,6) to [out = -90, in = 135] (0.5, 0.5) -- (0.5,0)
				(0.5,0.5) to (1.5,1.5) -- (1.5, 2) to [out = 135, in = -90] (1.5,6)
				(1.5,1.5) to [out = -45, in = 90] (1.5,0)
				(1.5,2) -- (2.5, 3);
				
			\draw [chambering] (5.5,6) -- (5.5,5.5) to [out=-45, in = 90] (5.5,0)
				(5.5, 5.5) -- (4.5, 4.5) to [out = 135, in = -90] (4.5, 6)
				(4.5, 4.5) -- (4.5, 4) to [out = -45, in = 90] (4.5,0)
				(4.5, 4) -- (3.5,3);
			
			\node at (0.5, 6.5) {$\alpha_0$};
			\node at (1, 7.5) {$id_{u_1} \otimes f_1$};
			\node at (1.5, 6.5) {$\sigma_1$};
			\node at (2, 8.5) {$id_{u_2} \otimes f_2$};
			\node at (4, 8.5) {$id_{u_{n-1}} \otimes f_{n-1}$};
			\node at (4.5, 6.5) {$\sigma_{n-1}$};
			\node at (5, 7.5) {$id_{u_n} \otimes f_n$};
			\node at (5.5, 6.5) {$\alpha_1$};
			
			\node at (0.5, -0.5) {$\beta_0$};
			\node at (1, -1.5) {$id_{u_1} \otimes f_1$};
			\node at (1.5, -0.5) {$\tau_1$};
			\node at (2, -2.5) {$id_{u_2} \otimes f_2$};
			\node at (4, -2.5) {$id_{u_{n-1}} \otimes f_{n-1}$};
			\node at (4.5, -0.5) {$\tau_{n-1}$};
			\node at (5, -1.5) {$id_{u_n} \otimes f_n$};
			\node at (5.5, -0.5) {$\beta_1$};
			
			\node at (0.75, 1.25) {$\gamma_1$};
			\node at (1.75, 2.75) {$\gamma_2$};
			\node at (4.25, 3) {$\gamma_{n-1}$};
			\node at (5.25, 4.75) {$\gamma_n$};

	\end{tikzpicture}
	\end{center}
	\caption{A canonical isomorphism between potential values of $H$.}
	\label{fig:Hnaturality}
\end{figure}  

As it happens, we in fact have slightly more structure than the above description lets on. The equivalence relation on sequences, described in Cor.~\ref{cor:QSequivreln1-mor}, is as follows. Two sequences 
\begin{align*}
	\alpha_0 \circ (id_{u_1} \otimes f_1) \circ \sigma_{1} \circ \cdots \circ \sigma_{n-1} \circ (id_{u_n} \otimes f_n) \circ \alpha_1 \\
	\beta_0 \circ (id_{v_1} \otimes g_1) \circ \tau_{1} \circ \cdots \circ \tau_{n-1} \circ (id_{v_n} \otimes g_n) \circ \beta_1
\end{align*}
represent the same 1-morphisms if and only if they have the same {\em length} (i.e. number of $f_i$ factors, see Cor.~\ref{cor:QSequivreln1-mor}), if $f_i = g_i$, and if there exist $\gamma_i: u_i \to v_i$ morphisms in $\Sigma_{G_0}$ (i.e. permutations) such that
\begin{align*}
	(id_{s(f_n)} \cdot \gamma_n ) \circ \beta_1 &= \alpha_1 \\
	\alpha_0 \circ (id_{t(f_1)} \cdot \gamma_1) & = \beta_0 \\
	 \sigma_i  \circ (id_{s(f_{i+1})} \cdot \gamma_{i+1})&=  (id_{t(f_i)} \cdot \gamma_i) \circ \tau_i .
\end{align*}
Suppose that we have two different choices of representatives for the same 1-morphism of  $\sF_{qs}(f_* P)$. Then for each choice of sequence of permutations $\gamma_i: u_i \to v_i$ which are used to witness the equivalence between these two representatives, there exists a canonical 2-isomorphisms in $\sF_{uss}(P)$ between the corresponding 1-morphisms of $\sF_{uss}(P)$. This 2-morphism is given by the string diagram depicted in Figure~\ref{fig:Hnaturality}. A simple exercise in the string diagram calculus of Section~\ref{sec:strictsymbicats} show that this morphism only depends on the source and target morphisms and does not depend on the choice of permutations $\gamma_i: u_i \to v_i$ used to witness the equivalence relation. 

Before defining $H$ on 2-morphisms, we will first describe the remaining coherence structure and the components of the natural transformation $\eta$. The homomorphism $H = (H, \phi^H_{gf}, \phi^H_a)$ is  weak, and so we must provide the coherence cells $\phi^H_{fg}$ and $\phi_a$. Since $H$ preserves identity 1-morphisms on the nose, we will choose to take $\phi_a$ to be trivial. However the coherence cells $\phi^H_{fg}$ must be non-trivial in general. 

There are two reasons that $H$ fails to commute strictly with composition. The first is due to the fact that we must make arbitrary choices of representatives of the sequences which define the 1-morphism of $\sF_{qs}(f_* P)$. The second is that even after making these choices, the homomorphism still fails to commute with the composition of sequences due to the fact that the permutation morphisms compose strictly. The first issue is easily solved by conjugating with the canonical isomorphisms shown in Figure~\ref{fig:Hnaturality}. Thus suppose that we are given representative sequences of composable 1-morphisms. 
\begin{align*}
	f = \alpha_1 \circ (id_{u_n} \otimes f_n) \circ \sigma_{n-1} \circ \cdots \circ \sigma_1 \circ (id_{u_1} \otimes f_1) \circ \alpha_0 \\
	g = \alpha_3 \circ (id_{v_m} \otimes g_m) \circ \tau_{m-1} \circ \cdots \circ \tau_1 \circ (id_{v_1} \otimes g_1) \circ \alpha_2
\end{align*}
The sequences $H(g) \circ H(f)$ and $H(g \circ f)$ in $\sF_{uss}(P)$ are nearly identical, and differ only at the location of composition:
\begin{align*}
	H(g) \circ H(f) & = \cdots \circ (id_{v_1} \otimes g_1) \circ \beta^{\alpha_2} \circ \beta^{\alpha_1} \circ (id_{u_n} \otimes f_n) \circ \cdots \\
	H(g \circ f) & = \cdots \circ (id_{v_1} \otimes g_1) \circ \beta^\omega \circ (id_{u_n} \otimes f_n) \circ \cdots 
\end{align*}
where $\omega$ is the permutation corresponding to the composition of $\alpha_1$ and $\alpha_2$. The coherence cell $\phi^H$ is given by applying the 2-morphism $X^{\alpha_2, \alpha_1}$. With these choices it is a simple matter to verify that $(H, \phi^H_{fg}, \phi^H_a)$ satisfies the hexagon and unitor squares necessary to define a homomorphism of bicategories. Note also that the homomorphism $\phi$ sends these coherence cells to identities. 

Next we construct the components of the natural equivalence $\eta: id \to H \phi$. Since the identity functor and $H \phi$ agree on objects, we may take the object-wise components of $\eta$ to be identities. In other words $\eta$ is an invertible {\em icon} in the terminology of \cite{Lack10}. It is given by an assignment of an isomorphism $\eta_f: f \to H\phi(f)$ for each 1-morphism $f$. In words, this morphism gives a string diagram which transforms an arbitrary 1-morphism in $\sF_{uss}(P)$, that is an arbitrary sequence of the elementary 1-morphism $(id_x \otimes f \otimes id_y)$ and $\beta^\sigma$ into one of the special form that arises in the image of $H$. This 2-morphism is built by preforming three operations. First we transform each elementary 1-morphism $(id_x \otimes f \otimes id_y)$ into a 1-morphism of the form $(id_z \otimes f)$.  This introduces additional permutation terms. Next all the adjacent permutation terms are `merged' together. Finally, if necessary, we compose with the canonical transformation depicted in Figure~\ref{fig:Hnaturality} to change to the correct reduced sequence corresponding to our arbitrary choice of representative. A typical example of the effect of these first two operations is given by the string diagram depicted schematically in Figure~\ref{fig:etastringdiagram}.

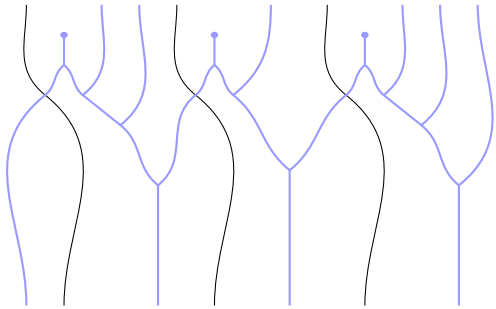
\begin{figure}[htpb]
	\begin{center}
		\begin{tikzpicture}[yscale=0.8]
			%\draw (0,0) rectangle (8,5);
			
			\draw (1,5) to [out = -90, in = 135] (1.25, 3.5) to [out = -45, in = 90] (1.5, 0);
			\draw (3,5) to [out = -90, in = 135] (3.25, 3.5) to [out = -45, in = 90] (3.5, 0);
			\draw (5,5) to [out = -90, in = 135] (5.25, 3.5) to [out = -45, in = 90] (5.5, 0);
			
			\draw [chambering, fill] (1.5, 4.5) circle (1pt);
			\draw [chambering, fill] (3.5, 4.5) circle (1pt);
			\draw [chambering, fill] (5.5, 4.5) circle (1pt);
			
			\draw [chambering] (1.5, 4.5) -- (1.5, 4) to [out = -135, in =45] 
				(1.25, 3.5) to [out = -135, in = 90] (1,0);
			\draw [chambering] (1.5, 4) to [out = -45, in = 135] (1.75, 3.5)
				to [out = -45, in = 135] (2.25, 3)
				to [out = -45, in = 135] (2.75, 2) -- (2.75, 0);
			\draw [chambering] (2,5) to [out = -90, in = 45] (1.75, 3.5);	
			\draw [chambering] (2.5,5) to [out = -90, in = 45] (2.25, 3);

			\draw [chambering] (3.5, 4.5) -- (3.5, 4) to [out = -135, in =45] 
				(3.25, 3.5) to [out = -135, in = 45] (2.75, 2);
			\draw [chambering] (3.5, 4) to [out = -45, in = 135] (3.75, 3.5) to 
				[out = -45, in = 135] (4.5, 2.25) -- (4.5, 0);
			\draw [chambering] (4.25, 5) to [out = -90, in = 45] (3.75, 3.5);
				
			\draw [chambering] (5.5, 4.5) -- (5.5, 4) to [out = -135, in =45] 
				(5.25, 3.5) to [out = -135, in = 45] (4.5, 2.25);
			\draw [chambering] (5.5, 4) to [out = -45, in = 135] (5.75, 3.5)
				to [out = -45, in = 135] (6.25, 3)
				to [out = -45, in = 135] (6.75, 2) -- (6.75,0);
				
			\draw [chambering] (6,5) to [out = -90, in = 45] (5.75, 3.5);
			\draw [chambering] (6.5,5) to [out = -90, in = 45] (6.25, 3);
			\draw [chambering] (7,5) to [out = -90, in = 45] (6.75, 2);
			
		\end{tikzpicture}
	\end{center}
	
	\caption{The components of $\eta: id \to H \phi$.}
	\label{fig:etastringdiagram}
\end{figure}

We may now, finally, define $H$ on 2-morphisms. Recall that the morphism 2-morphism of $\sF_{qs}(f_*P)$ are generated under the operations of horizontal and vertical composition, and by the operation of tensoring with identity 2-morphisms, by the elementary 2-morphisms $\phi^\otimes_{(f,id), (id, g)}$ with $f, g \in G_1$, and the generating 2-morphisms $\varphi \in G_2$. These are subject to the following relations: the relations from definition of 2-category (associativity, interchange, strict units), relations expressing the naturality of $\phi^\otimes$, the hexagon relation expressing that $\phi^\otimes$ is the coherence cell making $\otimes$ a homomorphism, and finally the relations $\cR$ from the 3-computad $P$. 

Each of these generating 2-morphisms naturally corresponds to a 2-morphism in $\sF_{uss}(P)$, however they will generally have the wrong source an target 1-morphisms to be the 2-cell part of a homomorphisms. Instead the homomorphism $H$ is defined on each of these generators to be the corresponding 2-morphism in $\sF_{uss}(P)$ conjugated by the components of the  2-isomorphisms $\eta$. Thus, for example, the image of $\phi^\otimes_{(f,id), (id, g)}$ under $H$ is given by the string diagram shown schematically in Figure~\ref{fig:H2cellimage}.

\begin{figure}[htbp]
	\begin{center}
	\begin{tikzpicture}
		%\draw (0,0) rectangle (3 ,3);
		\draw (1,3) to [out = -90, in = 135] (1.5, 2) to [out = -45, in = 90] (2,0);
		\draw (2,3) to [out = -90, in = 45] (1.75, 2.5) to [out = -135, in = 45] (1.5, 2) to [out = -135, in = 135] (1, 0.5) to [out = -45, in = 90] (1.25, 0);
		\draw [chambering] (1.5, 3) to [out = -90, in = 135] (1.75, 2.5) to [out = -45, in = 135] (2.25, 2.25) -- (2.25, 2)
			(2.25, 2.25) to [out = 45, in = -90] (2.5, 3);
		\draw [chambering, fill] (2.25, 2) circle (2pt);
		\draw [chambering] (0.75, 0) to [out = 90, in = -135] (1, 0.5) to [out = 45, in = -135] (1.25, 0.75) -- (1.25, 1)
			(1.25, 0.75) to [out = -45, in = 90] (1.5, 0);
		\draw [chambering, fill] (1.25, 1) circle (2pt);
		
	\end{tikzpicture}
	\end{center}
	\caption{$H$ applied to $\phi^\otimes_{(f,id), (id, g)}$. }
	\label{fig:H2cellimage}
\end{figure}

To see that this is well defined we must check that $H$ is compatible with the relations on the 2-morphisms of $\sF_{qs}(f_*P)$. The relations $\cR$ are imposed on both sides and thus hold automatically. The remaining relations are, in each case, easy to verify directly, though the lead to somewhat tedious string diagram calculations. We will omit the details. Note that with this choice the 2-morphism $\eta$ becomes a natural transformation. Specifically for each 2-morphism $\varphi$ of $\sF_{uss}(P)$ we have the naturality equation: $H \phi( \varphi) = \eta_{t(\varphi)} \circ \varphi \circ \eta_{s(\varphi)}$.
It is sufficient to check this on the generating 2-morphisms, for which again follows from a direct string diagram calculation. 
\end{proof}

Now we will turn to coherence results for general symmetric monoidal bicategories. Gurski-Osorno have several results in this direction. 

\begin{theorem}[{\cite[Thm~27]{MR2770448} and \cite[Thm~1.13]{GO13}}] \label{thm:GOCransCoherencetype1}
	Every symmetric monoidal bicategory is equivalent to a Crans semistrict symmetric monoidal 2-category. \qed.
\end{theorem}

The second coherence theorem of Gurski and Osorno is a coherence theorem of the second kind for symmetric monoidal bicategories freely generated by a set of objects. Let $G_0$ be a set and let $\Sigma_{G_0}$ be the symmetric monoidal category generated by that set, introduced in Section~\ref{sec:qsymcomputad}. This is a symmetric monoidal 1-category whose objects are words in the letters $G_0$. Viewed as a 2-category with only identities, this is also the free quasistrict symmetric monoidal 2-category $\sF_{qs}(G_0)$ generated by $G_0$. Gurski and Osorno's result is as follows:

\begin{proposition}[{\cite{GO13}}] \label{pro:GO-coherence} % I don't include a theorem number because the precise statement below is only stated for $G_0 = {pt}$, but it follows directly from their results. 
	The canonical strict homomorphism $\phi:\sF(G_{0}) \to \sF_{qs} = \Sigma_{/G_0}$ is an equivalence of symmetric monoidal bicategories, where $\sF(G_0)$ and $\sF_{qs}(G_0)$ denote the free symmetric monoidal bicategory, respectively the free quasistrict symmetric monoidal 2-category, generated by the 0-computad $G_0$. \qed
\end{proposition}

Gurski's and Osorno's coherence theorem is a bit more general than the one stated above, applying to symmetric monoidal bicategories freely generated by bicategories, however it does not immediately cover general symmetric monoidal computads. 

Our starting point will be to extend an argument of Lack \cite{MR1931220} which shows that in the presence of a cofibrancy theorem a coherence theorem of the first kind can sometimes be promoted to one of the second kind.

\begin{proposition}\label{pro:CrassCoherence}
	For every symmetric monoidal 3-computad $P$, the canonical (strict) homomorphism $\phi:\sF(P) \to f^*\sF_{css}f_*(P)$ is an equivalence of symmetric monoidal bicategories, where $\sF(P)$ and $\sF_{css}(f_*P)$ denote the free symmetric monoidal bicategory, respectively the free Crans semistrict symmetric monoidal 2-category, generated by the computad $P$.
\end{proposition}

\begin{proof}
	The first observation is that this homomorphism is essentially surjective on objects, essentially full on 1-morphisms, and full on 2-morphisms. This will be proven in Cor.~\ref{cor:phiis2-acyclic}, after a few technical lemmas on taking quotients in bicategories. Assuming this result, we will continue with the proof of Prop.~\ref{pro:CrassCoherence}. By Theorem~\ref{WhiteheadforSymMonBicats} (Whitehead's Theorem for Symmetric Monoidal Bicategories) all that remains is to verify that $\phi$ is faithful on 2-morphisms. 
	
By the Gurski-Osorno coherence result, Theorem~\ref{thm:GOCransCoherencetype1}, there is some symmetric monoidal equivalence $\sF(P) \to \sM$, with $\sM$ a Crans semistrict symmetric monoidal 2-category. By Theorem~\ref{thm:cofibrancythm} (the Cofibrancy Theorem) we can replace this with an equivalent strict homomorphism $f: \sF(P) \to \sM$, which of course is still an equivalence. In particular $f$ is faithful on 2-morphisms. 

By Cor~\ref{cor:leftadjointmorphismoftheories} and Cor~\ref{cor:fstarleftadjoint}, the canonical strict homomorphism $\phi$ can be identified with the unit  of the adjunction $L: \symbicat \leftrightarrows \symbicat_{css}: f^*$ between the categories of symmetric monoidal bicategories and Crans semistrict symmetric monoidal 2-categories. 
\begin{equation*}
	\phi = \eta: \sF(P) \to f^* L (\sF(P)) \cong f^* \sF_{css}f_*(P)
\end{equation*}
Here each is a category where the morphisms are the strict homomorphisms. By the universal property of this adjunction, there exists a strict symmetric monoidal homomorphism
\begin{equation*}
	h: f^* L (\sF(P)) \cong f^* \sF_{css}f_*(P) \to \sM
\end{equation*}
which factors the strict symmetric monoidal equivalence $f: \sF(P) \to \sM$, i.e., $f = h \phi$. However since $f$ is faithful on 2-morphisms, it follows that $\phi$ must also be faithful on 2-morphisms, and hence is an equivalence. 
\end{proof}

To prove Cor.~\ref{cor:phiis2-acyclic} below, we will need two lemmas.
Now we will return the the theories of symmetric monoidal bicategory that we have been studying and we will look at the comparison homomorphism between the theories of computadic $T$-algebras. Our goal is to show that for each of these theories the comparison homomorphism $\phi: \sF(P) \to f^* \sF' (f_* P)$ is $2$-acyclic, but first we will need some general but technical lemmas about forming quotient bicategories.

\begin{lemma}\label{lem:bicatquot2mor}
	Let $\sB$ be a bicategory and let $W$ be a set of 2-morphisms of $\sB$. Let $q: \sB \to \sB'$ denote the quotient map in the category $\bicat_s$ of bicategories and strict functors obtained by declaring the 2-morphisms in $W$ to be identity cells. Suppose that $W$ arrises as the collection of 2-morphisms in a sub-bicategory $\sW \subseteq \sB$ which contains all objects and 1-morphism, and such that each 2-morphism of $\sW$ is invertible in $\sW$. Then the quotient map $q: \sB \to \sB'$ is a bijection on objects, and surjective on both 1-morphisms and 2-morphisms. Moreover $\sB'$ is a strict 2-category. 
\end{lemma}

\begin{proof}
	Declaring certain 2-morphisms to be identities has the effect of identifying certain 1-morphisms and 2-morphisms. In general this could result in the formation of new 2-morphisms as composites of these newly formed quotients. However this has no effect on the objects, and so the quotient map $q: \sB \to \sB'$ will always be a bijection on objects. This observation doesn't require any conditions on the class $W$. 
	 
	The fact that $\sB'$ is a strict 2-category follows from the fact that $\sW$, as a subbicategory, contains all the associators and unitors. 
	
	Consider the following equivalence relations on 1-morphism and 2-morphisms of $\sB$:
	\begin{itemize}
		\item Two 1-morphisms $f,f' \in \sB$ are equivalent $f \sim f'$ if and only if there exists a 2-morphism $w \in \sW$ such that $w: f \to f'$. (note: this is symmetric because $w$ is invertible).
		\item Two 2-morphisms $\alpha, \alpha'$ in $\sB$ are equivalent if and only if there exists a sequence of 2-morphisms $\alpha_1, \alpha_2, \dots, \alpha_n$, and two sequences of 2-morphisms in $\sW$,  $w_0, w_1, \dots, w_n$ and $w_0', w_1', \dots, w_n'$ such that
		\begin{align*}
			\alpha & = w_0 \circ \alpha_1 \circ w_1 \circ \alpha_2 \circ \dots \circ \alpha_n \circ w_n \\
			\alpha' & = w_0' \circ \alpha_1 \circ w_1' \circ \alpha_2 \circ \dots \circ \alpha_n \circ w_n'. 
		\end{align*}
	\end{itemize}
Declaring the 2-morphisms of $\sW$ to be identities certainly enforces these equivalence relations to hold, and hence the functor $q: \sB \to \sB'$ must pass through the equivalence classes. However these relations are easily seen to be compatible with both the horizontal and vertical composition and hence we get a bicategory by passing to equivalence classes of both 1-morphisms and 2-morphisms. This bicategory is necessarily the quotient $\sB'$. By construction the map to equivalence classes is surjective on both 1-morphisms and 2-morphisms. 
\end{proof}

\begin{remark}
	Let $W$ be any class of invertible 2-morphisms. Then we can form a class which meets the conditions of the lemma by taking $\sW$ to be the smallest sub-bicategory which contains all the objects and 1-morphisms, as well as the morphisms of $W$. $\sW$ also necessarily contains the associator and unitor morphisms of $\sB$, but this is essentially the only restriction. If $W$ contains all the associator and unitor 2-morphisms, then identifying the morphisms of $\sW$ with identities makes no additional identifications beyond what is implied by identifying the 2-morphisms of $W$ with identities. To see this we let $q: \sB \to \sB'$ be the quotient map. We let $\sB'_0$ conisist of all objects and 1-morphisms of $\sB'$ and only the identity 2-morphisms. Since $W$ contains all associators and unitors, $\sB'$ is a strict 2-category and hence $\sB'_0$ is a sub-2-category. Then $\sW' = q^{-1}(\sB'_0)$ is a sub-bicategory of $\sB$ consisting of all objects and all 1-morphisms, together with all 2-morphisms of $\sB$ which become identities after applying $q$. It follows that $\sW \subseteq \sW'$.
\end{remark}

\begin{lemma}\label{lem:bicatquot1mor}
	Let $\sB$ be a strict 2-category and let $W$ be a set of 1-morphism of $\sB$. Let $q: \sB \to \sB'$ denote the quotient map in the category $\bicat_s$ of bicategories and strict functors obtained by declaring the 1-morphisms in $W$ to be identity cells. Suppose that $W$ satisfies the properties that it is obtained as the 1-morphisms of a sub-2-category which contains all the objects of $\sB$ and in which every 1-morphism in $W$ is invertible in this subcategory, and  
 that for any two morphisms $f, f' \in W$ which are parallel, i.e. the source and target of $f$ and $f'$ agree, then we have that $f \cong f'$ in $\sB$. Under these assumptions the quotient map $q: \sB \to \sB'$ is surjective on objects, 1-morphisms, and 2-morphisms. 
\end{lemma}

\begin{proof}
	 The assumption that $\sB$ be a strict 2-category is probably not necessary, but is convenient and  sufficient for our purposes. The proof is this is similar to the last lemma. Since the 1-morphisms of $W$ are invertible and obtained from a subcategory we can describe the 1-morphisms and objects of the quotient $\sB'$ explicitly as the quotient by an equivalence relation. Namely two objects are declared equivalent if there exists an morphism of $W$ connecting them, and two 1-morphisms $f, f'$ are declared equivalent if there exists a sequence of 1-morphisms $f_1, f_2, \dots, f_n$, and two sequences of 1-morphisms in $W$,  $w_0, w_1, \dots, w_n$ and $w_0', w_1', \dots, w_n'$ such that
		\begin{align*}
			f & = w_0 \circ f_1 \circ w_1 \circ f_2 \circ \dots \circ f_n \circ w_n \\
			f' & = w_0' \circ f_1 \circ w_1' \circ f_2 \circ \dots \circ f_n \circ w_n'. 
		\end{align*}
As before this equivalence relation is compatible with composition. It follows that the quotient $q: \sB \to \sB'$ is surjective on objects and 1-morphisms. However we must also see that it is surjective on 2-morphisms. This follows from the fact that parallel 1-morphisms which are identified where already isomorphic, and so no new 2-cells are created by this quotient process. 
\end{proof}

\begin{corollary}\label{cor:phiis2-acyclic}
	Let $P$ be a symmetric monoidal 3-computad. Then the symmetric monoidal homomorphisms $\phi:\sF(P) \to f^*\sF_{css}(f_*P)$ and $\phi:\sF(P) \to f^*\sF_{qs}(f_*P)$ are surjective on objects, 1-morphisms, and 2-morphisms, that is they are 2-acyclic in the sense of Def.~\ref{def:acyclicfib}.  
\end{corollary}

\begin{proof}
	In each case the homomorphism $\phi$ is obtained exactly by declaring certain invertible 1-morphisms and 1-morphisms to be identities. 
We can form this quotient in two stages, first forming the quotient by declaring the relevant 2-cells to be identities and then the relevant 1-cells. For the first stage of this process, Lemma~\ref{lem:bicatquot2mor} applies directly because the class of 1-morphism includes the associator and unitors of the underlying bicategorical structure. Because of this the class of morphisms that we intend to declare to be identities comes from such a sub-bicategory $\sW$.  
	
To apply Lemma~\ref{lem:bicatquot1mor} for the second stage of the quotient procedure, we need to know that any two parallel 2-morphisms which we are declaring to be identities are already equivalent. However these 1-morphisms are in the image of the inclusion functor $\sF(G_0) \to \sF(P)$ obtained by including just the underlying $0$-computad (i.e. just the generating objects). The necessary closure property then follows, for example, from the Gurski-Osorno coherence result Prop.~\ref{pro:GO-coherence}.
\end{proof}

We are now in a position to prove the main theorem of this section:

\begin{theorem}[The Coherence Theorem for Symmetric Monoidal Bicategories]\label{thm:CoherenceThm}
	Let $P = (G_0, G_1, G_2, \cR)$ be any symmetric monoidal 3-computad. Then the canonical strict homomorphism
	\begin{equation*}
		\phi: \sF(P) \to f^* \sF_{qs}(f_*P)
	\end{equation*}
	is an equivalence of symmetric monoidal bicategories. Here $\sF(P)$ and $\sF_{qs}(f_*P)$ denote the computadic symmetric monoidal bicategory, respectively the computadic quasistrict symmetric monoidal 2-category, generated by $P$. In particular, every symmetric monoidal bicategory is equivalent to a quasistrict symmetric monoidal 2-category.
\end{theorem}

\begin{proof}
	The final statement follows because of Cor.~\ref{cor:everysymbicatequivtocomputadic}. Specifically the unit map (a strict symmetric monoidal homomorphism) $\sF V_3 \sM \to \sM$ is a symmetric monoidal equivalence between an arbitrary symmetric monodial bicategory $\sM$ to a computadic one (presented by the 3-computad $V_3(\sM)$). 

In light of Prop.~\ref{pro:CrassCoherence}, it is sufficient to prove that the canonical map
\begin{equation*}
	\phi: f_{css}^* \sF_{css} (f_{css,*} P) \to f_{qs}^* \sF_{qs} (f_{qs,*} P)
\end{equation*} 
is an equivalence, where $\sF_{css}( f_{css,*} P)$, respectively $\sF_{qs}( f_{qs,*} P)$, is the computadic Crans semistrict symmetric monoidal 2-category, respectively the computadic quasistrict symmetric monoidal 2-category, generated by $P$.

Our strategy for proving this theorem will be to mimic our previous proof of Theorem~\ref{thm:USSCoherence}. As in that case, in light of Theorem~\ref{WhiteheadforSymMonBicats} (Whitehead's Theorem for Symmetric Monoidal Bicategories) it is sufficient to show that $\phi$ is an equivalence of underlying bicategories. As before, we will demonstrate this by explicitly constructing an inverse equivalence of bicategories. Specifically we will construct a homomorphism of bicategories:
\begin{equation*}
	H: \sF_{qs} (f_{qs,*} P) \to \sF_{css} (f_{css,*} P)
\end{equation*}
which, like $\phi$, is the identity on objects and which is a section of $\phi$ in the sense that $\phi H = id$. We will simultaneously construct a natural equivalence (in fact an invertible {\em icon} \cite{Lack10}) $\eta: id \to H \phi$. To alleviate our cluttered notation, we will stop writing the push forward functor $f_*$ and leave it implicit. 

The 3-computad $P$ is built inductively starting with a 0-computad $G_0$, and progressing through a 1-computad, 2-computad, and finally the 3-computad $P$. This naturally filters the problem of constructing the inverse equivalence $H$. 
\begin{center}
\begin{tikzpicture}
		\node (LT) at (-1, 3) {$\sF(G_0)$};
		\node (LM) at (-1, 1.5) {$\sF_{css}(G_0)$};
		\node (LB) at (-1, 0) {$\sF_{qs}(G_0)$};
		
		\node (MT) at (2, 3) {$\sF(G_0, G_1)$};
		\node (MM) at (2, 1.5) {$\sF_{css}(G_0, G_1)$};
		\node (MB) at (2, 0) {$\sF_{qs}(G_0, G_1)$};

		\node (RT) at (5, 3) {$\sF(G_0, G_1, G_2)$};
		\node (RM) at (5, 1.5) {$\sF_{css}(G_0, G_1, G_2)$};
		\node (RB) at (5, 0) {$\sF_{qs}(G_0, G_1, G_2)$};
		
		\node (RRT) at (8, 3) {$\sF(P)$};
		\node (RRM) at (8, 1.5) {$\sF_{css}(P)$};
		\node (RRB) at (8, 0) {$\sF_{qs}(P)$};		
		
		\draw [->] (LT) -- node [left] {$\simeq$} (LM);
		\draw [->] (LM) -- node [left] {$\simeq$} (LB);
		
		\draw [->] (MT) -- node [left] {$\simeq$} (MM);
		\draw [->] (MM) -- node [left] {$$} (MB);
		
		\draw [->] (RT) -- node [left] {$\simeq$} (RM);
		\draw [->] (RM) -- node [left] {$$} (RB);
		
		\draw [->] (RRT) -- node [left] {$\simeq$} (RRM);
		\draw [->] (RRM) -- node [left] {$$} (RRB);
		
		\draw [->] (LT) -- node [above] {$$} (MT);
		\draw [->] (MT) -- node [above] {$$} (RT);
		\draw [->] (RT) -- node [above] {$$} (RRT);
		
		\draw [->] (LM) -- node [above] {$$} (MM);
		\draw [->] (MM) -- node [above] {$$} (RM);
		\draw [->] (RM) -- node [above] {$$} (RRM);
		
		\draw [->] (LB) -- node [below=2ex] {$0$-bijective} (MB);
		\draw [->] (MB) -- node [below=2ex] {$1$-bijective} (RB);
		\draw [->] (RB) -- node [below=2ex] {$2$-surjective} (RRB);
		
		\draw [->] (LB) to [bend right] node [right] {$H_0$} (LM);
\end{tikzpicture}
\end{center}
The first, and probably most difficult, stage of this filtration is already dealt with by Gurski-Osorno's second coherence result, Prop.~\ref{thm:GOCransCoherencetype1}, which implies that 
\begin{equation*}
	\phi: \sF_{css}(G_0) \to \sF_{qs}(G_0)
\end{equation*}
is an equivalence of symmetric monoidal bicategories. Thus we may choose an inverse equivalence $H_0$, and we could make this an inverse symmetric monoidal equivalence, if we desired to do so. 

We will arrange so that $H_0$ is an equivalence we the desired properties. This is eased by two facts, first that $\sF_{qs}(G_0)$ is an ordinary symmetric monoidal 1-category, and second that, via Prop.~\ref{thm:GOCransCoherencetype1}, we already know that $\phi$ is an equivalence. From this later observation it suffices to construct any section $H_0$ of $\phi$; it will automatically be an inverse equivalence. 

To construct this section, first notice that since $\phi$ is the identity on objects between $\sF_{css}(G_0)$ and $\sF_{qs}(G_0)$, we can choose $H_0$ to also be the identity on objects. Next, recall from Cor.~\ref{cor:phiis2-acyclic} that $\phi$ is surjective on 1-morphisms. Thus for every 1-morphism in $\sF_{qs}(G_0)$ we choose, arbitrarily, a 1-morphism of $\sF_{css}(G_0)$ lying above it via $\phi$. An important observation is that an any two such lifts are isomorphic via a unique 2-isomorphism. These 2-isomorphisms provide the necessary coherence cells. Since  $\sF_{qs}(G_0)$ has only identity 2-cells, there is nothing to do in 2-morphisms, and so we have constructed the desired section $H$. The same observations construct the invertible transformation (invertible icon) $\eta: id \to H \phi$, whose object-wise components are identities. 

Now we will extend the section $H_0$ to $H$ on all of $\sF_{qs}(P)$. At this point the argument becomes very similar to the proof of Theorem~\ref{thm:USSCoherence}. The 1-morphisms of $\sF_{qs}(P)$ are equivalence classes of sequences of the form
\begin{equation*}
	\alpha_1 \circ (id_{u_n} \otimes f_n) \circ \sigma_{n-1} \circ \cdots \circ \sigma_1 \circ (id_{u_1} \otimes f_1) \circ \alpha_0
\end{equation*}
where each $f_i \in G_1$ is a generating 1-morphism and each $\alpha_i$, $\sigma_i$ is a morphism from $\sF_{qs}(G_0)$. Given such a sequence, we may `evaluate' it in $\sF_{css}(G_0)$ as the sequence:
\begin{equation*}
	H_0(\alpha_1) \circ (id_{u_n} \otimes f_n) \circ H_0(\sigma_{n-1}) \circ \cdots \circ H_0(\sigma_1) \circ (id_{u_1} \otimes f_1) \circ H_0(\alpha_0)
\end{equation*}
In short we simply apply $H_0$ to all the morphisms coming from $\sF_{qs}(G_0)$ and just lift the elementary generating 1-morphisms from $G_1$ in the obvious way. As before this is almost strictly compatible with composition. If we have two composable sequences,
\begin{align*}
	f = \alpha_1 \circ (id_{u_n} \otimes f_n) \circ \sigma_{n-1} \circ \cdots \circ \sigma_1 \circ (id_{u_1} \otimes f_1) \circ \alpha_0 \\
	g = \alpha_3 \circ (id_{v_m} \otimes g_m) \circ \tau_{m-1} \circ \cdots \circ \tau_1 \circ (id_{v_1} \otimes g_1) \circ \alpha_2
\end{align*}
then the sequences $H(g) \circ H(f)$ and $H(g \circ f)$ in $\sF_{css}(P)$ differ only at the location of composition:
\begin{align*}
	H(g) \circ H(f) & = \cdots \circ (id_{v_1} \otimes g_1) \circ H_0({\alpha_2}) \circ H_0({\alpha_1}) \circ (id_{u_n} \otimes f_n) \circ \cdots \\
	H(g \circ f) & = \cdots \circ (id_{v_1} \otimes g_1) \circ H_0(\alpha_2 \alpha_1) \circ (id_{u_n} \otimes f_n) \circ \cdots 
\end{align*}
and thus the coherence morphism $\phi^{H_0}$ for $H_0$ also yields a coherence morphism for $H$. 

Of course the 1-morphisms of $\sF_{qs}(P)$ are not actually such sequences, but instead are equivalence classes of such sequences. Thus to define $H$ we must again make a choice of representative of each equivalence class for each 1-morphism in $\sF_{qs}(P)$. However, just as in the proof of Theorem~\ref{thm:USSCoherence}, there is a canonical isomorphisms between any two such choices. In fact given a collection $\gamma_i$ of morphisms in $\sF_{qs}(G_0)$, which are used to witness the equivalence between 2-representatives, then the string diagram depicted in Figure~\ref{fig:Hnaturality} actually gives a prescription for the induced 2-isomorphism. The string diagram is built from 2-morphisms which merge permutations (these come from coherence cells $\phi^{H_0}$) and `crossings' between permuations and the generators $f_i$. This later case corresponds to using an application of an interchanger $\phi^\otimes$ applied to a permuation and the generating 1-morphism. In short the string diagram depicted in Figure~\ref{fig:Hnaturality} gives a prescription for the a 2-isomorphism between different potential values of $H$ for a given 1-morphism of $\sF_{qs}(P)$. As tedious exercise in applying the axioms of Crans semistrict symmetric monoidal 2-categories shows that this 2-morphism is independent of the choice of morphisms $\{\gamma_i\}$ witnessing the equivalence. Thus in this case, as in the proof of Theorem~\ref{thm:USSCoherence}, there is a canonical 2-isomorphism between the values of different sequences representing the same 1-morphism in $\sF_{qs}(P)$. 

Next we turn to the value of the transformation $\eta: id \to H\phi$. The morphisms of $\sF_{css}(P)$ are give by compositions of elementary 1-morphisms of the from $id_w \otimes f \otimes id_{w'}$ where $w, w'$ are arbitrary objects and either $f \in G_1$ is a generating 1-morphism or $f$ is the component of a braiding transformation $f = \beta_{x,y}$. The transformation $\eta$ is a 2-morphisms which transforms this into an equivalent 1-morphism lying in the image of $H$. The reader should compare with the construction of $\eta$ in the proof of Theorem~\ref{thm:USSCoherence}, and in particular with Figure~\ref{fig:etastringdiagram}, which depicts that $\eta$ in terms of string diagrams. The current $\eta$ is built analogously. We first transform all the elementary 1-morphisms coming from $f \in G_1$ into ones of the form $(id_z \otimes f)$. Instead of an application of the transformations $X^e$ and $(X^{\sigma, \sigma'})^{-1}$, we use the syllepsis, otherwise the process is virtually identical. Next all the permutation terms are composed together. In current situation we may use the coherence morphisms of $H_0$ to accomplish this, rather than the $X$ transformations. 

Thus once again, the string diagram we used previously gives us a prescription for constructing the desired morphisms. In fact this is not at all surprising. In light of Theorem~\ref{thm:USSCoherence} and Prop.~\ref{pro:GO-coherence} we have a chain of equivalences
\begin{equation*}
	\sF_{css}(G_0) \stackrel{\phi \simeq}{\longrightarrow} \sF_{uss}(G_0) \stackrel{\phi \simeq}{\longrightarrow}\sF_{qs}(G_0).
\end{equation*}
Thus any string diagram calculation which occurs away from the dark lines (labeled by generating 1-morphisms) is valid in the 2-category $\sF_{css}(G_0)$. Thus our string diagram calculations from before are much more than just a heuristic for constructing analogous 2-morphisms in this new context. 
 
The final task is to define $H$ on 2-morphisms, and this, too, is accomplished just as for Theorem~\ref{thm:USSCoherence}, by conjugating the naive choices of lifts by the components of $\eta$. To verify this is a well-defined choice we must verify the same relations as we did previously for Theorem~\ref{thm:USSCoherence}. The proofs that these relations hold are similar. 	
\end{proof}

In the next chapter we will see how the coherence results of this section can be applied and used to help classify extended topological field theories. Before ending, though, we would like to record one last coherence theorem, a coherence theorem for functors, which is a consequence of Theorem~\ref{thm:CoherenceThm} and our previous considerations. 

\begin{theorem}[Cohernece for Symmetric Monoidal Homomorphisms]\label{thm:CoherenceFunctors}
	Let $f: \sA \to \sB$ be a symmetric monoidal homomorphism between symmetric monoidal bicategories. Then $f$ fits into a functorial diagram of the form:	
	
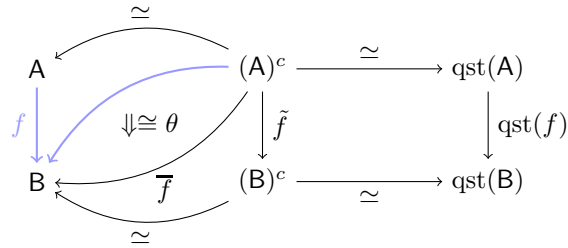
\begin{center}
	\begin{tikzpicture}
			\node (LT) at (0, 1.5) {$\sA$};
			\node (LB) at (0, 0) {$\sB$};
			\node (MT) at (3, 1.5) {$(\sA)^c$};
			\node (MB) at (3, 0) {$(\sB)^c$};
			\node (RT) at (6, 1.5) {$\qst(\sA) $};
			\node (RB) at (6, 0) {$\qst(\sB)$};
			\draw [->, thick, chambering] (LT) -- node [left] {$f$} (LB);
			\draw [<-] (LT) to [bend left] node [above] {$\simeq$} (MT);
			\draw [->] (MT) -- node [right] {$\tilde{f}$} (MB);
			\draw [<-] (LB) to [bend right] node [below] {$\simeq$} (MB);
			\draw [->] (MT) to [bend left] node [below] {$\overline{f}$} (LB);
			\draw [->, thick, chambering] (MT) to [bend right] node [above] {$$} (LB);
			\draw [->] (MT) -- node [above] {$\simeq$} (RT);
			\draw [->] (MB) -- node [below] {$\simeq$}(RB);
			\draw [->] (RT) -- node [right] {$\qst(f)$} (RB);
			\node at (1.5, 0.75) {$\Downarrow \cong \theta$};
			%\node at (0.5, 1) {$\ulcorner$};
			%\node at (1.5, 0.5) {$\lrcorner$};
	\end{tikzpicture}
	\end{center}
	In which every morphism is a strict homomorphism of symmetric monoidal bicategories, except the two indicated with thicker chambering arrows; in which every unlabeled region of the diagram commutes strictly; in which the symmetric monoidal bicategories $(\sA)^c$ and $(\sB)^c$ are computadic; in which the $\qst(\sA)$ and $\qst(\sB)$ are quasistrict symmetric monoidal bicategories (depending only on $\sA$ and $\sB$, respectively); and  in which the indicated arrows are symmetric monoidal equivalences. 
\end{theorem}

\begin{proof}
	For a symmetric monoidal bicatgory $\sM$, we take $(\sM)^c = \sF V_3(\sM)$ to be its canonical computadic replacement. Theorem~\ref{thm:cofibrancythm} (The Cofibrancy Theorem) then constructs a strict replacement for the composite
	\begin{equation*}
		f_1: (\sA)^c \to \sA \to \sB,
	\end{equation*}
	specifically the replacement is given by $j \circ \res(f_1) = \overline{f}$, and we have a natural symmetric monoidal 2-isomorphism $\theta: f_1 \to j \circ \res(f_1)$. The morphism $\tilde{f}$ arises as the mate of $\overline{f}$, using the universal property of the adjunction $(\sF, V_3)$. The final square then arises from the natural transformation $\phi$. The indicated morphisms are equivalences because of Corollary~\ref{cor:everysymbicatequivtocomputadic} and Theorem~\ref{thm:CoherenceThm}.
\end{proof}

 % Symmetric Monoidal Bicategories

\chapter[The Classification of 2D Extended TFTs ]{The Classification of 2-Dimensional Extended Topological Field Theories} \label{ClassificationChapter}

In this chapter we assemble the results of the previous two chapters to prove our main classification theorems. The majority of the results of this chapter are entirely new. There are four primary goals for this chapter. First we must introduce the symmetric monoidal bicategory of bordisms with structure. 
%Previous attempts at defining the bordism bicategory have been made in \cite{KL01} and \cite{Morton07}, but neither of these is adequate. First, they essentially ignore the {\em symmetric monoidal} structure of the bordism bicategory, which is essential to define extended topological field theories. This is perhaps understandable, since they did not have available the results of Chapter \ref{SymMonBicatChapt}. 
%
%Additionally, neither  \cite{KL01} nor \cite{Morton07}, discuss the issue of gluing bordisms in sufficient detail to adapt their definition to a bicategory of bordisms with structure. The main issue is that 
When one glues two smooth manifolds along a common boundary, the resulting space, $W \cup_Y W'$, which is a topological manifold, does not have a unique smooth structure extending those on $W$ and $W'$. While any two choices of such smooth structure result in diffeomorphic manifolds, this diffeomorphism is not {\em canonical}. In the 1-categorical bordism category, this issue is easily swept under the rug, but in the higher categorical setting it must be dealt with.   %Kerler-Lyubashenko \cite{KL01} and Morton  \cite{Morton07}, narrowly avoid this issue. 

In Section \ref{SectionBordBicat1} we discuss these gluing issues at length and give a preliminary, naive  version of the $d$-dimensional bordism bicategory as a symmetric monoidal bicategory.  In Section \ref{SectBordBicat2} we improve upon this definition and give the definition of the symmetric monoidal bicategory of bordisms with $\cF$-structures. The allowed structures are very general and include structures which can be described by sheaves, such as orientations, or even stacks, such as super manifold structures, spin structures, or principal bundles. 

In Section \ref{SectUnorientedClassification} we address our second primary goal, which is to prove our main classification theorem (Theorem~\ref{TheoremUnorientedClassification}) for 2-dimensional extended unoriented topological field theories. This theorem gives a generators and relations presentation for the extended 2-dimensional unoriented bordism bicategory as a symmetric monoidal bicategory. Using the results of Chapeters \ref{ChapPlanarDecomp} and \ref{SymMonBicatChapt}, the proof is relatively straightforward. This effectively classifies extended 2-dimensional unoriented topological field theories with {\em any} target symmetric monoidal bicategory, in the sense that such topological field theories are completely determined by a explicit list of data, the data coming from the presentation.  %This is the first time such a result has been rigorously proven in the literature. 

In Section \ref{SectOrientedClassification} we prove our third primary goal, and show how the techniques developed here can be generalized and adapted to classify bordism bicategories with structure. We demonstrate this by proving the analogous classification theorem (Theorem~\ref{OrientedBordismClassification}) for the 2-dimensional {\em oriented} extended bordism bicategory. 

Finally, in Section \ref{SectExamplesandApps}, we achieve our fourth primary goal. We use the classification results of the previous section to classify extended topological field theories with values in the symmetric monoidal bicategory of algebras, bimodules, and intertwiners. In particular, over a perfect field we prove that extended oriented topological field theories are the same as semi-simple (non-commutative) Frobenius algebras and that in the unoriented case they are essentially the same as semi-simple Frobenius $*$-algebras.
		
\section{The Bordism Bicategory I: Gluing} \label{SectionBordBicat1}

The modern interpretation of Atiyah's axiomatic definition of topological quantum field theories hinges upon the existence of the bordism category. The objects and morphism of this category are clear. The objects (for the $d$-dimensional unoriented category) are closed $(d-1)$-dimensional manifolds and the morphisms from $W_0$ to $W_1$ are $d$-dimensional (1-)bordisms from $W_0$ to $W_1$, i.e., compact $d$-manifolds $S$ equipped with a decomposition and isomorphism of their boundaries $\partial S = \partial_\text{in} S \sqcup \partial_\text{out}S \cong W_0 \sqcup W_1$. These bordisms are taken up to isomorphism of bordism, i.e., isomorphisms $S \cong S'$, such that the induced map,
\begin{equation*}
 W_0 \sqcup W_1 \cong \partial_\text{in} S \sqcup \partial_\text{out}S  = \partial S \cong \partial S' = \partial_\text{in} S' \sqcup \partial_\text{out}S' \cong W_0 \sqcup W_1,
\end{equation*}
is the identity.

In the purely topological setting (i.e., where these manifolds are {\em topological} manifolds) composition is clear and is given simply by gluing bordisms in the obvious way. In the smooth setting, gluing is mildly more complicated. A smooth structure on the topological manifold $S \cup_W S'$ is not uniquely determined by smooth structures on $S$ and $S'$. In addition one needs to specify a tubular neighborhood $W \times \R \hookrightarrow S \cup_W S'$ which is smooth when restricted to the collars $W_+ = W \times [0, \infty) \hookrightarrow S $ and $W_- = W \times (-\infty, 0] \hookrightarrow S' $. Different choices of collars generally induce different smooth structures on $S \cup_W S'$. Nevertheless, the resulting smooth manifolds are diffeomorphic, but by a non-canonical diffeomorphism (as we shall see shortly). This is sufficient to provide the bordism category with a well-defined composition and one easily verifies that it is, indeed, a category. 

In the higher categorical setting these mild difficulties are only amplified. In the bicategorical case what one is after is a bicategory whose objects are closed $(d-2)$-manifolds, $Y$, whose 1-morphism are compact $(d-1)$-dimensional 1-bordisms, $W$, between these and whose 2-morphisms are compact bordisms, $S$, between the 1-bordisms. Such a setup almost invariably requires the use of manifolds with corners. More seriously, in such a bicategory one must be able to glue the 1-morphisms ( 1-bordisms) horizontally. This is problematic because, as we have just seen, this gluing is not well defined on the nose, but only up to {\em non-canonical} diffeomorphism. One may be tempted to take diffeomorphism classes 1-bordisms as the 1-morphisms, but such an approach is doomed to failure since one will be unable to construct an acceptable notion of 2-bordism (there will be no well-defined composition of 2-bordisms).

Alternatively one may try to skirt this issue by using the axiom of choice to choose for each pair of composable 1-bordisms, $W$ and $W'$, a smooth composite $W \circ W'$ (as suggested in \cite[\oldS 5.1 ]{Morton07}) or equivalently by choosing a skeleton (or finite set of representatives) for each diffeomorphism class of 1-bordism (as done in \cite{KL01}). Both of these approaches result in well defined composition functors for pairs of composable 1-bordisms, and moreover both approaches construct a natural diffeomorphism:
\begin{equation*}
	a: (W \circ W') \circ W'' \to W \circ (W' \circ W'').
\end{equation*}
However this natural diffeomorphism will {\em fail} to satisfy the pentagon axiom! This problem is narrowly avoided in \cite{KL01} and \cite{Morton07}, but for subtle reasons. 

%There doesn't appear to be an adequate discussion elucidating these issues in the literature and so, i
In this section,  we begin by discussing the gluing of bordisms and provide a solution to the above problem. This also provides us an opportunity to freshen the mind of the reader and to fix terminology.
In the next section we will give an alternative solution to this problem which is more compatible with equipping our bordisms with additional structures. 

\subsection{Manifolds with Corners}

	Let $\R_+ = \{ y \in \R \;  | \; y \geq 0 \} = [0, \infty)$ and let $X$ be a topological $m$-manifold with boundary. Let $x \in X$. A {\em chart at $x$} (of {\em index} $p$) is a continuous map $\varphi: U \to \R^m_+$, where $U$ is an open subset of $X$ containing $x$, $\varphi$ is a homeomorphism onto its image, and in which $\varphi(x)$ has $p$-many coordinates which take the value zero. 
	Two charts $\varphi: U \to \R^m_+$ and $\varphi': U' \to \R^m_+$ are {\em compatible} if the composite,
	\begin{equation*}
		 \varphi' \circ \varphi^{-1}: V \to \R^m_+
		 \end{equation*}
is a diffeomorphism onto its image, where $V =  \varphi(U \cap U')$. Here if $S \subset  \R^m_+$ is a subset, then  $f: S \to   \R^m_+$ is {\em smooth} if there exists an open neighborhood $W \subset \R^m$ of $S$ and a smooth extension of $f$ to $W$. Note that if $\varphi$ and $\varphi'$ are two charts at $x$ which are compatible, then the index of $\varphi$ coincides with the index of $\varphi'$.

If $X$ is a topological manifold with boundary then an {\em atlas} for $X$ is a compatible family of charts, containing a chart at $x$ every point $x \in X$. Given an atlas for $X$, the {\em index} of a point $x \in X$ is the index of any chart at $x$. Atlases are ordered by inclusion. A {\em manifold with corners} \index{manifold with corners} is a topological manifold with boundary equipped with a maximal atlas. Given an atlas $\cA$ for $X$, there exists a unique maximal atlas containing $\cA$.

\begin{figure}[ht]
\begin{center}
\begin{tikzpicture} [thick]
\filldraw [fill=green!20!white, draw=green!50!black] (0,0) -- (0,1) arc (180: -90: 1 cm) -- (0,0);
\end{tikzpicture}
\caption{A manifold with corners without a collar.}
\label{CirclewithCornerFig}
\end{center}
\end{figure}
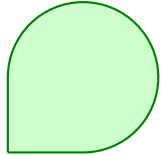

A submanifold of a manifold with corners $X$ is a manifold with corners $Y$ with an embedding $i: Y \hookrightarrow X$, such that in any charts at $y$ for $Y$  and at $i(y)$ for $X$, $i$ is smooth. Let $Y$ be a submanifold of $X$ of codimension one such that $Y \subset \partial X$. Then a {\em collar} of $Y$ is a map $f_Y: U_Y \to Y \times \R_+$ such that $U_Y$ is an open neighborhood of $Y$ in $X$, $f_Y(y) = (y, 0)$ for all $y \in Y$, and $f_Y$ is a diffeomorphism of manifolds with corners. Not every submanifold $Y \subset \partial X$ admits a collar, as the example in Figure~\ref{CirclewithCornerFig} demonstrates.

A {\em connected face} of a manifold with corners $X$ is the closure of a component of $\{ x \in X \; | \; \text{ index}(x) = 1\}$. A manifold with corners is a {\em manifold with faces} if each $x \in X$ belongs to $\text{index}(x)$ different connected faces. For such a manifold a {\em face} is a disjoint union of connected faces, and is itself a manifold with faces (and hence a manifold with corners). If $X$ is a manifold with faces, then so is $\R_+ \times X$.

\begin{proposition} \label{GluingProp}
Let $X$ and $X'$ be manifolds with faces and $Y$ a manifold with corners, together with inclusions $f:Y \to X$ and $f': Y \to X'$ realizing $Y$ as a face of each. Then $X \cup_Y X'$ is a topological manifold with boundary. If in addition we are given collars $f_+: Y \times \R_+ \to X$ and $f_-: Y \times \R_+ \to X'$, there exists a canonical smooth structure on $X \cup_Y X'$ which is compatible with the smooth structures on $X$ and $X'$. 
\end{proposition}

\begin{proof}
	The  property that $X \cup_Y X'$ is a topological manifold with boundary is a completely local question and may be readily checked. To provide $X \cup_Y X'$ with the structure of a manifold with corners we must provide a maximal smooth atlas. It is sufficient to provide a covering by open smooth manifolds with corners, such that the transition functions are smooth on overlaps. Let $\breve{X}$ denote $X \setminus Y$ and similarly let $\breve{X}'$ denote $X' \setminus Y$. Then  $\breve{X}$, $\breve{X}'$, and $\R \times Y$ provide such a cover, where the first two maps to $X \cup_Y X'$ are the obvious inclusions and the third, from $\R \times Y$, is given by the union of the collars $f_+$ and $f_-$.
\end{proof}

The next lemma and theorem are taken from \cite[\oldS 6]{Munkres61}, the proofs of which are readily adapted to cover the case of manifolds with faces. 
We will delay the proof of the following lemma until after that of the following theorem.

\begin{lemma} \label{UniquenessOfGluingHelperLma}
	Let $M$ be a compact manifold with faces and $N$ an (open) neighborhood of $M \times \{0\}$ in $M \times \R_+$. Let $f: N \hookrightarrow M \times \R_+$ be smooth and equal to the identity on $M\times \{0 \}$. Then there exists a diffeomorphism $\tilde f: N \to f(N)$ such that
\begin{enumerate}
\item $\tilde f = f$ in a neighborhood of the complement of $N$,
\item $\tilde f = id$ in a neighborhood of $M\times \{0\}$. 
\end{enumerate}
\end{lemma}

\begin{theorem} \label{UniquenessOfGluingThm}
	Let $X$ and $X'$ be manifolds with faces and $Y$ a manifold with corners, together with inclusions $f:Y \to X$ and $f': Y \to X'$ realizing $Y$ as a face of each. Suppose that $Y$ admits collars  $f_+: Y \times \R_+ \to X$ and $f_- : Y \times \R_+ \to X'$, but that these are not specified. \index{bordism!gluing without a choice of collars}
Then the glued manifold with faces, $X \cup_Y X'$, is unique up to diffeomorphism fixing $Y$, and which coincides with the identity map outside a neighborhood of $Y$.   
\end{theorem}

\begin{proof}
Let $P =(f_+, f_-)$ and $P' = (f_+', f_-')$ be two choices of pairs of collars for $Y$. Then $P$ is equivalent to a map
%\begin{equation*}
$	P: Y \times \R \to X \cup_Y X'
$%\end{equation*}
which is a homeomorphism onto its image, $V$, and smooth on the complement of $Y \times \{0\}$. Similarly, $P': Y \times \R \to X \cup_Y X'$ is a map which a is a homeomorphism onto its image, $V'$, and smooth on the complement of $Y \times \{0\}$.

Let $U = P^{-1}( V \cap V') \subset Y \times \R$. Then the composite $(P')^{-1} \circ P$ is defined on $U$ and is a homeomorhism onto $(P')^{-1} \circ P(U)$. Moreover it is the identity on $Y \times \{0\}$ and is smooth when restricted to $U_+ = U \cap Y \times [0, \infty)$ and $U_- = U \cap Y \times (-\infty, 0]$, separately. 

By Lemma \ref{UniquenessOfGluingHelperLma} there exists a homeomorphism $g: U \to (P')^{-1} P(U)$ such that
\begin{enumerate}
\item $g$ equals $ (P')^{-1} P$ near the boundary of $U$,
\item $g$ is a diffeomorphism when restricted to $U_+$ and $U_-$,
\item $g$ coincides with the identity map in a neighborhood of $Y \times \{0\}$. 
\end{enumerate}
Hence $P' g P^{-1}$ is defined on the neighborhood $V \cap V'$ of $Y$ in $X \cup_Y X'$ and coincides with the identity map near the boundary of this neighborhood. Extending by the identity map yields a homoeomorphism $f$ of $X \cup_Y X'$. One may readily check that $f$ is a diffeomorphism.
\end{proof}

\begin{proof}[Proof of Lemma \ref{UniquenessOfGluingHelperLma}]
There exists a smooth function $\beta: M \to [0, 1)$ such that $N' = \{ (x, t) \in M \times \R_+ \; | \; t \leq \beta(x) \} \subseteq N$. Proving the lemma for $N'$ proves it for $N$, hence without loss of generality we may assume that $N = N'$. 

Let $f(x, t) = ( X(x,t), T(x,t))$. Then $\frac{\partial T}{\partial t} ( x, 0) > 0$. Hence we may choose $\varepsilon: M \to (0,1)$, such that 
\begin{enumerate}
\item $\frac{\partial T}{\partial t} ( x, t) > 0$ for all $0 \leq t \leq \varepsilon(x)$, 
\item $\varepsilon(x) |\frac{\partial T}{\partial t} ( x, t) | < 1$ for all $0 \leq t \leq \beta(x)$.
\end{enumerate}
 Choose a smooth monotone function $\alpha: \R \to \R$, such that $\alpha(t) = 0$ for all $t \leq \frac 1 3$ and $\alpha(t) = 1$ for all $t \geq \frac 2 3$. Define the function on $N$,
 \begin{equation*}
	\psi(x,t) = ( 1- \alpha(\frac{t}{\beta(x)} )) \varepsilon(x) t + \alpha(\frac{t}{\beta(x)} ) t.
\end{equation*}
Then $g(x,t) = (x, \psi(x,t))$ defines a diffeomorphism $g: N \to N$. 

Setting $f_1 = f \circ g: N \to f(N)$, we see that $f_1$ is also a diffeomorphism. Let $f_1(x,t) = (X_1(x,t), T_1(x,t))$.  Define the function on $N$,
\begin{equation*}
	\phi(x,t) = \alpha( \frac{2t}{\beta(x)}) t + [1 - \alpha( \frac{2t}{\beta(x)})] T_1(x,t).
\end{equation*}
Then $h(x,t) = (x, \phi(x,t))$ defines a diffeomorphism $h: N \to N$. Set $f_2 = f_1 \circ h^{-1}$, and let $f_2(x,t) = (X_2(x,t), T_2(x,t))$.

We have that $T_2 \equiv t$ on a neighborhood $L$ of $M \times \{0\}$, and since $M$ is compact, there exists a $\delta >0$ such that $M \times [0, \delta] \subseteq L$. Finally we let
\begin{equation*}
	X_3(x,t) = X_2( x, \alpha( \frac t \delta) t),
\end{equation*}
and define $\tilde f (x, t) = (X_3(x,t), T_2(x,t))$. One may check that $\tilde f$ has the desired properties. 
\end{proof}

\begin{remark} \label{UniqueIsotopyClassRmk}
	The diffeomorphism constructed in Lemma \ref{UniquenessOfGluingHelperLma} depends on the choice of functions $\varepsilon$, $\beta$ and $\alpha$. Any other choice of such data $(\tilde \varepsilon, \tilde \beta, \tilde \alpha)$ is isotopic to the original choice of data through valid choices of this data. In fact the space of such data is a contractible space. \index{bordism!contractible space of gluings} Tracing through the proof of Lemma \ref{UniquenessOfGluingHelperLma} and Theorem~\ref{UniquenessOfGluingThm} we see that this implies that, while the diffeomorphism constructed in the proof of Theorem~\ref{UniquenessOfGluingThm} is not canonical, it does have a {\em canonical isotopy class}.
\end{remark}

Following \cite{Janich68, Laures00, LP05, Morton07} we will define our higher bordisms as manifolds with faces equipped with additional structure.

\begin{definition}
	An {\em $\langle n \rangle$-manifold} \index{$\langle n \rangle$-manifold} is a manfold with faces $X$ together with an ordered $n$-tuple $(\partial_0 X, \dots, \partial_{n-1}X)$ of faces of $X$ which satisfy the following conditions:
	\begin{enumerate}
\item $\partial_0 X \cup \dots \cup \partial_{n-1}X = \partial X$, and
\item $\partial_i X \cap \partial_j X$ is a (possibly empty) face of $\partial_i X$ and $\partial_j X$ for all $i \neq j$. 
\end{enumerate}
\end{definition}

Let $[1]$ denote the ordered set $\{0, 1\}$, viewed as a category. Let $\langle n \rangle$ denote the category $[1]^{\times n}$. An $\langle n \rangle$-manifold induces an $\langle n \rangle$-diagram of spaces, i.e., a functor $X: \langle n \rangle \to \Top$. If $X$ is an $\langle n \rangle$-manfold and $a= (a_0, \dots, a_{n-1})$ is an object of $\langle n \rangle$, then
\begin{equation*}
	X(a) = \cap_{a_i \neq 0} \partial_i X
\end{equation*}
 with morphisms being sent to the obvious inclusions. By convention $X(0,0, \dots, 0) = X$. 
 
 \begin{example}
 $X= \R^m_+$ is a manifold with faces in the obvious manner. The faces consist of the $m$-coordinate hyperplanes. Order these hyperplanes $H_0, H_1, \dots, H_m$. Then $X$ becomes an $\langle m \rangle$-manifold where $\partial_i X = H_i$. Moreover, if $e_i$, $0 \leq i \leq m-1$ denote the standard basis vectors for $\R^m$, then
 \begin{equation*}
	X(a_0, \dots, a_{m-1}) = \text{Span} \{ e_i \; | \; a_i = 0 \}.
\end{equation*}
 \end{example}

 As first recognized in \cite{LP05}, the combinatorics of $\langle n \rangle$-manifolds is particularly well suited to extended bordism categories. In addition, the faces of $\langle n \rangle$-manifolds admit collars in the sense of the following lemma.
 
\begin{lemma}[ {\cite[Lemma 2.1.6]{Laures00}} ] \label{CollaringManifoldsWithCornersLemma}
	Any $\langle n \rangle$-manifold $X$ admits and $\langle n \rangle$-diagram $C$ of embeddings
\begin{equation*}
	C( a < b): \R^n_+(a') \times X(a) \hookrightarrow \R^n_+(b') \times X(b)
\end{equation*}
with the property that $C(a<b)$ restricted to $\R^n_+(b') \times X(a)$ is the inclusion map $id \times X(a < b)$. Here $a' = (1, \dots, 1) - a$ and $b' = (1, \dots, 1) - b$. \qed
\end{lemma}

\subsection{The cobordism bicategory} \label{sec:cobbicat}

With these definitions at hand we are ready to give the first definition of the $d$-dimensional bordism bicategory $\cob_d$. The objects are smooth closed $(d-2)$-dimensional manifolds $Y$. The 1-morphisms and 2-morphisms are 1-bordisms and isomorphism classes of 2-bordisms, which we define shortly. 

\begin{definition} \label{def:1bordism}
	Let $Y_0$ and $Y_1$ be smooth closed $(d-2)$-manifolds. A {\em 1-bordism} \index{bordism} from $Y_0$ to $Y_1$ is a smooth compact $(d-1)$-manifold with boundary, $W$, equipped with a decomposition and isomorphism of its boundary $\partial W = \partial_\text{in}W \sqcup \partial_\text{out} W \cong Y_0 \sqcup Y_1$. 
\end{definition}

\begin{definition} \label{Def2BordVer1}
	Let $Y_0$ and $Y_1$ be smooth closed $(d-2)$-manifolds and let $W_0$ and $W_1$ be two 1-bordisms from $Y_0$ to $Y_1$. A {\em 2-bordism} consists of a compact $d$-dimensional $\langle 2 \rangle$-manifold, $S$, equipped with the following additional structures:
	\begin{itemize}
\item A decomposition and isomorphism:
\begin{equation*}
	\partial_0 S = \partial_{0, \text{in}} S \sqcup \partial_{0, \text{out}}S \stackrel{g}{\to} W_0 \sqcup W_1.
\end{equation*}
\item A decomposition and isomorphism:
\begin{equation*}
	\partial_1 S = \partial_{1, \text{in}} S \sqcup \partial_{1, \text{out}}S \stackrel{f}{\to} Y_0 \times I \sqcup Y_1 \times I.
\end{equation*}
\end{itemize}
These are required to induce isomorphisms 
\begin{align*}
	f^{-1}g: & \partial_\text{in}W_0 \sqcup \partial_\text{out} W_0 
		\to Y_0 \times \{0\} \sqcup Y_1 \times \{0\}\\ 
	f^{-1}g: & \partial_\text{in}W_1 \sqcup \partial_\text{out} W_1 \to Y_0 \times \{1\} \sqcup Y_1 \times \{1\},
\end{align*}
which coincide with the structure isomorphisms of $W_0$ and $W_1$. 

Two 2-bordisms,  $S$ and $S'$, are isomorphic if there is a diffeomorphism $h: S \to S'$, which restricts to diffeomorphisms,
\begin{align*}
	\partial_{0, \text{in}} S & \to \partial_{0, \text{in}} S'  \\
	\partial_{0, \text{out}} S & \to \partial_{0, \text{out}} S' \\
	\partial_{1, \text{in}} S & \to \partial_{1, \text{in}} S' \\
	\partial_{1, \text{out}} S & \to \partial_{1, \text{out}} S' 
\end{align*}
and such that $f' \circ h = f$ and $g' \circ h = g$.
\end{definition}

Thus a 2-bordism consists of a diagram of manifolds as in Figure~\ref{2BordDiagramFig}. There are several equivalence relations on diffeomorphism that will be useful for us. 
\begin{figure}[ht]
\begin{center}
	\begin{tikzpicture}[thick, scale = 2]
		% top
		\draw (1, 3.5) arc (90: 270: 0.3cm and 0.1cm) -- (1.1, 3.3) ;
		\draw (-.1, 3.5) -- (0,3.5) arc (90: -90: 0.3cm and 0.1cm);
		\node at (0.5, 3.75) {$W_0$};
		\node at (0.5, 3.25) {$\downarrow$};

		% main
		\draw (0, 3) arc (90: -90: 0.3cm and 0.1cm) -- (0, 2);
		\draw (1, 3) arc (90: 270: 0.3cm and 0.1cm) -- (1.1, 2.8) -- (1.1, 2) -- (0, 2);
		\draw (1,3) -- (1, 2.8); \draw [densely dashed] (1, 2.8) -- (1, 2.2) -- (0,2.2);
		\draw (.3, 2.9) arc (180: 360: 0.2cm and 0.3cm);
		\draw (0, 3) -- (-.1, 3) -- (-.1, 2.2) -- (0, 2.2);
		\node at (0.5, 2.35) {$S$};
		
		% bottom
		\draw (-0.1, 1.5) -- (1,1.5) (0,1.2) -- (1.1, 1.2);
		\node at (0.5, 0.75) {$W_1$};
		\node at (0.5, 1.75) {$\uparrow$};
		
		% left
		\draw (-0.5, 2.8) -- (-0.5, 2) (-0.6, 3) -- (-0.6, 2.2);
		\filldraw (-1, 2.4) circle (1pt) (-1.1, 2.6) circle (1pt); 
		\node at (-1.5, 2.5) {$Y_1$};
		\node at (-0.5,1.75) {$Y_1 \times I$};
		\node at (-0.25, 2.5) {$\rightarrow$};
		\draw [->] (-1.05, 2.75) to [out=90, in = 180] (-0.5, 3.4);
		\draw [->] (-1.05, 2.25) to [out=270, in = 180] (-0.5, 1.35);
		
		% right
		\draw (1.6, 2.8) -- (1.6, 2) (1.5, 3) -- (1.5, 2.2);
		\filldraw (2.1, 2.4) circle (1pt) (2, 2.6) circle (1pt); 
		\node at (2.5, 2.5) {$Y_0$};
		\node at (1.5,1.75) {$Y_0 \times I$};
		\node at (1.25, 2.5) {$\leftarrow$};
		\draw [->] (2.05, 2.75) to [out=90, in=0] (1.5,3.4);
		\draw [->] (2.05, 2.25) to [out=270, in=0] (1.5,1.35);
	\end{tikzpicture}

\caption{The structure of a 2-bordism}
\label{2BordDiagramFig}
\end{center}
\end{figure}
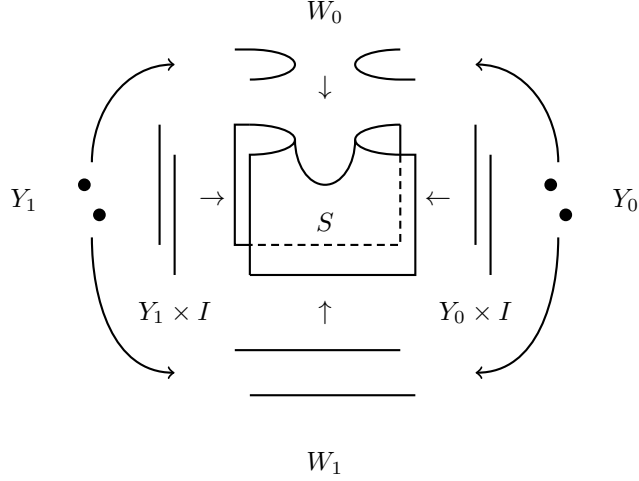
\begin{definition}
	Let $Y$ and $Z$ be smooth manifolds (possibly with corners). Let $f_0, f_1: Y \to Z$ be diffeomorphisms such that $f_0 |_{\partial Y } = f_1|_{\partial Y}$. Then $f_0$ is {\em pseudo-isotopic} \index{pseudo-isotopy} to $f_1$ if there exists a diffeomorphism $F: Y \times I  \to Z \times I$ such that,
	\begin{equation*}
F|_{\partial Y \times I} = f_0|_{\partial Y} \times id: \partial Y \times I \to \partial Z \times I,
\end{equation*}
$F(x,t) = (f_0(x), t)$ in a neighborhood of $Y \times \{0\}$, and $F(x,t) = (f_1(x), t)$ in a neighborhood of $Y \times \{1\}$.
	\end{definition}
Note that an isotopy is a particular kind of pseudo-isotopy, in which $F(y,t) = (f_t(y), t)$ is level preserving. The importance of  pseudo-isotopy for the bordism bicategory arises from the following lemma.
\begin{lemma} \label{PI2-morlemma}
	Let $(S, f,g)$ be a 2-bordism. If we replace $f$ by a pseudo-isotopic map, then the resulting 2-bordism is in the same isomorphism class. \index{bordism!pseuo-isotopy and,}
\end{lemma}
  
\begin{figure}[h]
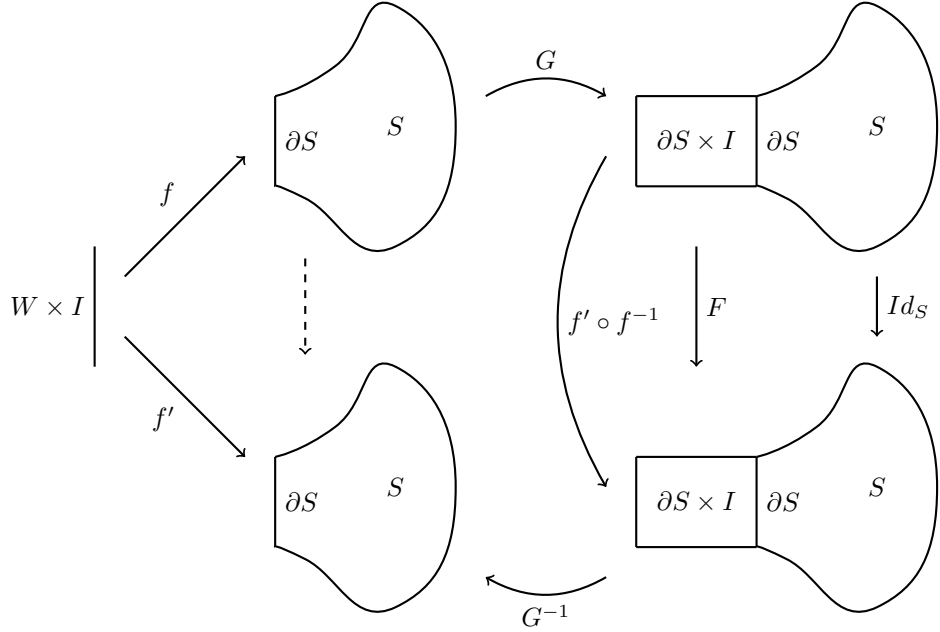

\begin{center}

\caption{Pseudo-isotopy doesn't affect 2-morphisms.}
\label{PI2-morphism}
\end{center}
\end{figure}

\begin{proof}
Suppose that $f$ and $f'$ are pseudo-isotopic. This is equivalent to $f' \circ f^{-1}$ being pseudo-isotopic to the identity map on $\partial S$. Suppose that $F$ is a diffeomorphism of $\partial S \times I$ realizing this pseudo-isotopy (which we may assume is the constant isotopy on $\partial_0S$).  We have that $F \cup id_S$ is a self-diffeomorphism of $(\partial S \times I ) \cup S$ which is $f' \circ f^{-1}$ on the ($\partial_1$) boundary. Choose a diffeomorphism $G: S \to (\partial S \times I ) \cup S$ which is the identity on the boundary.  
Then $h = G^{-1} \circ (F \cup Id_S) \circ G$ is an isomorphism between the 2-bordisms $(S,f,g)$ and $(S,f', g)$. See  Figure~\ref{PI2-morphism}.
% We will construct a morphism of $W$ which is the identity on the boundary not equal to $\partial_iW$, and which restricts to $f_i' \circ (f_i)^{-1}$ on $\partial_iW$. Thus there is a diffeomorphism of $W$, which makes everything commute and the two different classes are equal $[W, f_i, \dots] = [W, f_i', \dots]$. To get the diffeomorphism $W \to W$, we first apply $G$, then we apply $F \cup Id_W$ to $(\partial_i W \times I ) \cup W$. Finally we apply $G^{-1}$ to return to $W$. This is a diffeomorphism of $W$ which is the identity on all the boundaries of $W$, except $\partial_iW$, where it is $f_i' \circ (f_i)^{-1}$, exactly as claimed.
\end{proof}

We can now define vertical composition of 2-bordisms. Given 1-bordisms $W_0$, $W_1$, $W_2$ from $Y_0$ to $Y_1$, and 2-bordisms $S$ from $W_0$ to $W_1$ and $S'$ from $W_1$ to $W_2$, we would like to define their vertical composition by choosing collars for $W_1$ in $S$ and $S'$, as per  Lemma \ref{CollaringManifoldsWithCornersLemma}, and use these to glue $S$ and $S'$ together as in Proposition~\ref{GluingProp}. The result is again a $\langle 2 \rangle$-manifold $S' \circ S$, which has a decomposition and isomorphism,
\begin{equation*}
	\partial_0 (S' \circ S) = \partial_{0, \text{in}}S \sqcup \partial_{0, \text{out}}S' \stackrel{g_\text{in} \sqcup g'_\text{out}}{\longrightarrow} W_0 \sqcup W_2.
\end{equation*}
However the second decomposition and isomorphism is problematic. It is as follows,
\begin{equation*}
	\partial_0 (S' \circ S) = \partial_{1, \text{in}}S \cup_{Y_0}  \partial_{1, \text{in}}S' 
	 \sqcup \partial_{1, \text{out}}S \cup_{Y_1}  \partial_{1, \text{out}}S'  \stackrel{f \cup f'}{\longrightarrow}
	 Y_0 \times (I \cup_{pt} I) \sqcup Y_1 \times (I \cup_{pt} I). 
\end{equation*}
In order to make $S' \circ S$ into a 2-bordism, we must choose isomorphisms,
\begin{equation*}
	Y_i \times (I \cup_{pt} I) \cong Y_i \times I. 
\end{equation*}
This can be done uniformly by choosing once and for all a diffeomorphism $I \cup_{pt} I \cong I$. Any two such choices are pseudo-isotopic, in fact isotopic.  In this way, $S' \circ S$ becomes a 2-bordism from $W_0$ to $W_2$. There were several choices involved, but by Theorem~\ref{UniquenessOfGluingThm}, Remark \ref{UniqueIsotopyClassRmk}, and Lemma \ref{PI2-morlemma} the isomorphism class of this 2-bordism is well defined and doesn't depend on the choices we've made. Moreover this vertical composition defines an operation which is associative on isomorphism classes of 2-bordisms, which again follows from Lemma  \ref{PI2-morlemma}.

\begin{proposition}
	Fix closed $(d-2)$-manifolds $Y_0, Y_1$. Then $\cob_d(Y_0, Y_1)$ is a category, where the objects of $\cob_d(Y_0, Y_1)$ consist of the 1-bordisms from from $Y_0$ to $Y_1$, the morphisms consist of isomorphism classes of 2-bordisms between these and composition is given by vertical composition of 2-bordism.
\end{proposition}

\begin{proof}
The only detail left to be checked is that there exist identity 2-bordisms. We leave it to the reader to verify that $W \times I$ gives such an identity.
\end{proof}

Horizontal composition of 2-bordisms and 1-bordisms is defined similarly, but we will need to use the axiom of choice. Given a 1-bordism $W$ from $Y_0$ to $Y_1$ and a 1-bordism $W'$ from $Y_1$ to $Y_2$, we would like to obtain a new 1-bordism, $W' \circ W$, from $Y_0$ to $Y_1$ by gluing together $W$ and $W'$. Since our 1-bordisms are not taken up to isomorphism, the choices we make in defining this gluing are more important. Propositon \ref{GluingProp} ensures that once we have chosen collars for $Y_1$, that we can form a canonical gluing. So in order to define the horizontal composition of 1-morphisms,  we must use the axiom of choice to choose these  collars for each 1-bordism. 

There is no way to ensure that these choices are compatible, and so this gluing is associative only up to non-canonical isomorphism of 1-bordism. However, following Remark \ref{UniqueIsotopyClassRmk}, there is a {\em canonical isotopy class} of diffeomorphisms and hence a canonical isomorphism class of 2-bordisms realizing the associativity of horizontal composition. Thus even though these diffeomorphisms fail to satisfy the pentagon axioms, the induced bordisms do satisfy the pentagon axioms. This allows one to construct an actually bicategory of bordisms. 

Horizontal composition of 2-morphisms can now be defined just as vertical composition. Given composable 2-bordisms 
\begin{equation*}
	S \in \cob_d(Y_0, Y_1)(W, W') \quad \text{ and } \quad S' \in \cob_d(Y_1, Y_2)(W'', W''')
\end{equation*}
 we define $S'*S \in \cob_d(Y_0, Y_2)(W'' \circ W, W''' \circ W')$ by choosing collars of $Y_1 \times I$ and gluing. These collars can be chosen to be compatible with our choices for $W$, $W'$, $W''$ and $W'''$ and the resulting isomorphism class of 2-bordism is well-defined. Similar considerations involving identities apply (identities are given by the 1-bordisms $Y_0 \times I$) and yield the following proposition.

\begin{proposition}
	 $\cob_d$ is a bicategory in which the objects are closed $(d-2)$-manifolds, the morphism categories are the categories  $\cob_d(Y_0, Y_1)$, and where composition is given by horizontal composition together with its canonical associators and unitors. 
\end{proposition}

Given a diffeomorphism $f: Y_0 \to Y_1$ of closed $(d-2)$-manifolds, we may promote it to a 1-bordism as follows. The $(d-1)$ manifold is $W = Y_0 \times I$, and the decomposition/isomorphism is given by,
\begin{equation*}
	\partial W = Y_0 \times \{0\}  \sqcup Y_0 \times \{1\} \stackrel{id \sqcup f}{\longrightarrow} Y_0 \sqcup Y_1. 
\end{equation*}
This assignment respects composition, at least up to canonical (isomorphism class of) 2-bordism, and similarly isomorphisms of 1-bordisms may be turned into 2-bordisms. In particular we may use this to promote the symmetric monoidal structure on the categories of $(d-2)$-manifolds, 1-bordisms, and 2-bordisms given by disjoint union into a fully fledged symmetric monoidal structure on $\cob_d$. 

\subsection{Symmetric monoidal pseudo-double categories} \label{sec:symmonpseudo2cat}

In this section we review a method for constructing symmetric monoidal bicategories developed by Shulman \cite{Shu1004}.
Roughly speaking, the idea is to first construct the bordism bicategory not as a bicategory, but as a symmetric monoidal pseudo-double category. Such a gadget consists of a pair of (ordinary) symmetric monoidal categories, an `object' symmetric monoidal category, and a `morphism' symmetric monoidal category together with structure maps making it a weak category object in symmetric monoidal categories. 

Bordisms very naturally assemble into such a structure. The category of objects will be the usual cobordism category of $(d-2)$-manifolds and diffeomorphism classes of $(d-1)$-cobordisms. The category of morphisms will have objects $(d-1)$-cobordisms (not taken up to diffeomorphism), and the morphisms will be 2-dimensional cobordisms between cobordisms. As before these are certain $\langle 2\rangle$-manifolds modeled on the square. Only now both the vertical and horizontal sides are allowed to be non-trivial. This is similar to the structure considered by  \cite{Morton07} and \cite{Feshbach:2011aa}. 
From this one may extract a bicategory by discarding the vertical 1-morphisms and restricting to those 2-morphisms whose vertical sides are trivial.   

Recall the following definition from \cite{Shu1004}:
\begin{definition}%\label{def:}
	A \emph{symmetric monoidal pseudo-double category} $\D$ consists of a `symmetric monoidal category of objects' $\D_0$ and a `symmetric monoidal category of arrows' $\D_1$, with structure symmetric monoidal functors
	\begin{align*}
		U: &\D_0 \to \D_1 \\
		S,T: &\D_1 \to \D_0 \\
		\odot: &\D_1 \times_{\D_0} \D_1 \to \D_1
	\end{align*}
	(where the pullback is over $\D_1 \stackrel{T}{\to} \D_0 \stackrel{S}{\leftarrow} \D_1$) such that $S$ and $T$ are strict and
	\begin{align*}
		S(U_A) &= A \\
		T(U_A) &= A \\
		S(M \odot N) &= SN \\
		T(M \odot N) &= TM 
	\end{align*}	
	equipped with natural isomorphisms 
	\begin{align*}
		\alpha: & (M \odot N) \odot P \cong M \odot (N \odot P) \\
		\lambda: & U_B \odot M  \cong M \\
		\rho: & M \odot U_A \cong M
 	\end{align*}
	such that $S(\alpha)$, $T(\alpha)$, $S(\lambda)$, $T(\lambda)$, $S(\rho)$, and $T(\rho)$ are identities, and such that the standard axioms for a monoidal category/bicategory are satisfied (specifically the pentagon and triangle identities). We will further assume that $U$ preserves unit objects strictly. 
\end{definition}

A pseudo double category has \emph{objects} (the objects of $\D_0$), \emph{vertical 1-morphisms} (the morphisms of $\D_0$), \emph{horizontal 1-morphisms} (the objects of $\D_1$), and \emph{2-morphisms} (the morphisms of $\D_1$) which can be thought of as squares. If $g$ is a horizontal 1-morphism, we will sometimes write $g: S(g) \nrightarrow T(g)$ to denote its source and target in $\D$. If $(f: a \to b) \in \D_1$, we picture it as a square as follows:
\begin{center}
\begin{tikzpicture}
		\node (LT) at (0, 2) {$S(a)$};
		\node (LB) at (0, 0) {$S(b)$};
		\node (RT) at (2, 2) {$T(a)$};
		\node (RB) at (2, 0) {$T(b)$};
		\draw [->] (LT) -- node [left] {$S(f)$} (LB);
		\draw [->] (LT) -- node [above] {$a$} (RT);
		\draw [->] (RT) -- node [right] {$T(f)$} (RB);
		\draw [->] (LB) -- node [below] {$b$} (RB);
		\node at (1,1) {$f$};
		%\node at (0.5, 1) {$\ulcorner$};
		%\node at (1.5, 0.5) {$\lrcorner$};
\end{tikzpicture}
\end{center} 
Such squares have a strictly associative vertical composition (coming from the composition in $\D_1$) and a coherently associative horizontal composition induced by the functor $\odot$. 

A 2-morphism $f: a \to b$ whose sides $S(f)$ and $T(f)$ are identities is called \emph{globular}, and the collection consisting of the objects, horizontal 1-morphisms, and globular 2-morphisms naturally inherits the structure of a bicategory, the \emph{horizontal bicategory} $H(\D)$, see \cite{Shu1004}. 

\begin{definition} \label{def:companion}
	Let $\D$ be a symmetric monoidal pseudo double category and $f: a \to b$ a vertical 1-morphisms. A \emph{companion} of $f$ is a horizontal 1-morphism $\hat{f}: a \nrightarrow b$ together with 2-morphisms
	\begin{center}

	\end{center}
	such that the following two equations hold:
	\begin{center}
	%
	\end{center}
	(note that in the second equation we are implicitly using the unit constraint cells to interpret this equation). 
	
	A \emph{conjoint} of $f$, denoted $\check{f}: a \nrightarrow b$, is a companion of $f$ in the double category $\D^\textrm{h-op}$ obtained by reversing the horizontal 1-morphisms, but not the vertical 1-morphisms of $\D$. 
\end{definition}

\begin{definition}
	A symmetric monoidal pseudo double category $\D$ is \emph{fibrant} if every vertical 1-morphism of $\D$ has both a companion and a conjoint. 
\end{definition}

The main result from \cite{Shu1004} which we will utilize is the following:

\begin{theorem}[{\cite[Thm~5.1]{Shu1004} }]\label{thm:ShulmanThm}
	If $\D$ is a fibrant symmetric monoidal pseudo double category, then $H(\D)$ is a symmetric monoidal bicategory. \qed
\end{theorem}

\subsection{The symmetric monoidal cobordism bicategory}\label{sec:symmonbicatcobord}

Now we will implement that above construction in the special case of the cobordism bicategory. We will introduce a symmetric monoidal pseudo double category of cobordisms $\COB_d = ((\COB_d)_0, (\COB_d)_1)$. The symmetric monoidal category of objects, $(\COB_d)_0$, is the symmetric monoidal category of closed compact $(d-2)$-manifolds and their diffeomorphisms. 

% Specifically, the objects of $(\COB_d)_0$ are smooth closed $(d-2)$-manifolds, and the morphisms (the vertical morphisms of $\COB_d$) are equivalence classes of 1-bordisms, as in Definition~\ref{def:1bordism}. Here two such 1-bordisms are equivalent if they are diffeomorphic relative to the boundary. From what we have seen already it is clear that composition makes this into a category and it is easily seen that disjoint union makes it symmetric monoidal. A 1-bordism $W$ from $Y_0$ to $Y_1$ will be denoted $(W; Y_0, Y_1)$.

The symmetric monoidal category of morphisms is more interesting. The objects of $(\COB_d)_1$ (the horizontal morphisms) are 1-bordisms, as in Definition \ref{def:1bordism}. The morphisms of $(\COB_d)_1$ are very similar to the 2-bordisms of Definition~\ref{Def2BordVer1}:

\begin{definition}
	A morphism of $(\COB_d)_1$ from $(W_\textrm{in}; Y_{\textrm{in},0}, Y_{\textrm{in},1})$ to $(W_\textrm{out}; Y_{\textrm{out},0}, Y_{\textrm{out},1})$ consists of a pair of diffeomorphisms 
	\begin{equation*}
		h_0: Y_{\textrm{in},0} \to Y_{\textrm{out},0} \quad \textrm{ and } \quad h_1 :  Y_{\textrm{in},1} \to Y_{\textrm{out},1}
	\end{equation*}
	%(which induce 1-bordisms $(Y_{\textrm{in},0} \times I; Y_{\textrm{in},0} , Y_{\textrm{out},0} )$ and $(Y_{\textrm{in},1} \times I; Y_{\textrm{in},1} , Y_{\textrm{out},1} )$), 
	together with an equivalence class of compact $d$-dimensional $\langle 2 \rangle$-manifold $S$, equipped with a decomposition and isomorphism
	\begin{equation*}
		\partial_0 S = \partial_{0, \text{in}} S \sqcup \partial_{0, \text{out}}S \stackrel{g}{\to} W_\textrm{in} \sqcup W_\textrm{out},
	\end{equation*}
of the `horizontal boundary' $\partial_0S$, and together with a decomposition and isomorphism of the `vertical boundary'
\begin{equation*}
	\partial_1 S = \partial_{1, \text{s}} S \sqcup \partial_{1, \text{t}}S \stackrel{f}{\to} Y_{\textrm{in},0} \times I \sqcup Y_{\textrm{in},1} \times I
\end{equation*}
such that $f^{-1} \circ g$ restricts to the diffeomorphisms $h_0$ and $h_1$ on $Y_{\textrm{in},0} \times \{1\}  \subseteq Y_{\textrm{in},0} \times I$ and $Y_{\textrm{in},1} \times \{1\} \subseteq Y_{\textrm{in},1} \times I$, respectively. 

Two such $\langle 2 \rangle$-manifolds $S$ and $S'$ from $(W_\textrm{in}; Y_{\textrm{in},0}, Y_{\textrm{in},1})$ to $(W_\textrm{out}; Y_{\textrm{out},0}, Y_{\textrm{out},1})$ are {\em equivalent} if there exists a diffeomorphism of $\langle 2 \rangle$-manifolds which is relative to the horizontal boundary ($\partial_0S$), and vertical boundary.
\end{definition}

Thus the morphisms of $(\COB_d)_1$ can be visualized very similarly to the 2-bordism depicted in Figure~\ref{2BordDiagramFig}. 
The key differences are that the vertical sides are no longer required to be trivial cobordisms, but may be twisted by the diffeomorphisms $h_i$. In particular the source and targets of the source and target horizontal bordisms are not required to agree. 

As before the equivalence class of these morphisms only depends on the pseudo-isotopy class of the parametrization $f$. Thus the morphisms of $(\COB_d)_1$ have a composition given by gluing composable morphisms along their common component of horizontal boundary which is defined exactly as the vertical composition of the bicategory $\cob_d$, see Section~\ref{sec:cobbicat}. This category is symmetric monoidal via the disjoint union. 

We must define the functors  $S$, $T$, $U$, and $\odot$. The functors $S$ and $T$ are the easiest to define. They send the objects of $(\COB_d)_1$ to the source, respectively target, $(d-2)$-manifold of the given 1-bordism. On a morphism $S$ of $ (\COB_d)_1$  they restrict to the respective diffeomorphisms $h_0$ and $h_1$, respectively. The functors $S$ and $T$ are strict symmetric monoidal. 

The functor $U: (\COB_d)_0 \to (\COB_d)_1$ is given by `$(-) \times I$', that is it takes a $(d-2)$-manifold $Y$ to the identity 1-bordism $Y \times I$, and it takes a diffeomorphism $h: Y_\textrm{in} \to Y_\textrm{out}$ to the $\langle 2 \rangle$-manifold $Y_\textrm{in} \times I^2$, with the top, left, and right sides parametrized by $Y_\textrm{in} \times I$ via the obvious identity map, and with the bottom side parametrized by $Y_\textrm{out} \times I$ via the map $h$. This functor only fails to be strictly monoidal as a result of the fact that the set $(A \sqcup B) \times C$ is not equal to $(A \times C) \sqcup (B \times C) $, but merely (canonically) isomorphic. Regardless the functor $U$ is certainly a symmetric monoidal functor. 

The functor $\odot$ and the transformations $\alpha$, $\lambda$, and $\rho$ are defined exactly as for the horizontal composition in the bicategory  $\cob_d$, see Section~\ref{sec:cobbicat}. The compositon, which is built upon a pushout of sets, distributes over disjoint unions. Hence these functors and transformations are readily seen to be symmetric monoidal, though they will not be strict. 

The horizontal bicateogry $\cH(\COB_d) = \cob_d$ is precisely the bicategory of cobordisms introduced in Section~\ref{sec:cobbicat}. We will show:

\begin{theorem} \label{thm:symmonbicatcob}
	The symmetric monoidal pseudo double category $\COB_d$ is fibrant, and hence $\cH(\COB_d) = \cob_d$ is a symmetric monoidal bicategory. 
\end{theorem}

\begin{figure}[htbp]
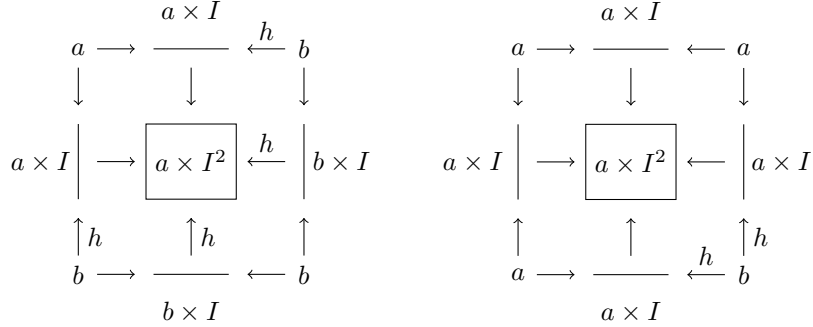

	\begin{center}
	
		
	\caption{Morphisms witnessing companionship for $\COB_d$.}
	\label{fig:companionbordisms}
\end{center}
\end{figure}

\begin{proof}
	The conclusion follows from Theorem~\ref{thm:ShulmanThm}, once we establish that $\COB_d$ is fibrant. We begin by observing that $(\COB_d)^\textrm{h-op} \cong \COB_d$. Thus it suffices to show that every morphism of $(\COB_d)_0$ has a companion, as in Definition~\ref{def:companion}.  Let $a,b \in (\COB_d)_0$ be two $(d-2)$-manifolds. Given a diffeomorphism $\overline{h}: a \to b$, with inverse diffeomorphism $h: b \to a$, the companion of $\overline{h}$ is given by the 1-bordism $(a \times I; a, b)$, which is a trivial bordisms where outgoing boundary has been reparametrized using the diffeomorphism $h$. The morphisms of $(\COB_d)_1$ depicted in Figure~\ref{fig:companionbordisms}  witness $(a \times I; a, b)$ as the companion to $\overline{h}$.
\end{proof}

%-------
%
%\begin{theorem}
%	With the structures given by disjoint union,  $\cob_d$ is a symmetric monoidal bicategory.
%\end{theorem}
%
%\begin{proof}
%	Disjoint unions gives the categories of $(d-2)$-manifolds, 1-bordisms, and 2-bordisms the structure of symmetric monoidal categories.  By turning these into 1-bordisms and 2-bordisms, as above,  we obtain the data of a symmetric monoidal structure on the bicategory $\cob_d$. One only needs to verify that these satisfy the requisite axioms, which is straightforward, although tedious.\footnote{See \cite{Shu1004} for an intriguing alternative approach to this result.}  
%\end{proof}

\section{The Bordism Bicategory II: Halos} \label{SectBordBicat2}

While acceptable as a first attempt, the construction of $\cob_d$ in the previous section has at least two major drawbacks. First, in order to define the horizontal composition of 1-bordisms (and hence 2-bordisms as well) we were required to use a vast application of the axiom of choice and choose collars arbitrarily for each 1-bordism. This could be avoided by using one of the so-called {\em non-algebraic} definitions of bicategory, \cite{Leinster02, Leinster04, CL04}, where providing specific compositions is not necessary. We will not do this.  More importantly, we will want to consider variations on the bordism bicategory which involve placing additional structures on our bordisms, such as orientations, spin structures, or principal bundles, and these can be difficult to incorporate into our previous attempt in a satisfactory manner.

Both of these problems can be solved by modifying our bordisms (and objects) to include germs of neighborhoods of $d$-manifolds. We call these {\em halos} and they easily generalize to include additional structures ($\cF$-halos). The construction of the bordism bicategory then precedes in essentially the same manner as in the last section, using haloed $(d-2)$-manifolds, haloed 1-bordisms, and haloed 2-bordisms. Less care is needed in the gluing constructions, and we obtain symmetric monoidal bicategories $\bord_d$ and $\bord_d^\cF$. 

\subsection{Categories of pro-objects}

To make the notion of germ of neighborhoods precise we will use the category of {\em pro-manifolds}. These are formal cofiltered limits of manifolds. We will give a brief review of the theory of pro-objects in general categories, though this theory is well-known. Our treatment closely follows \cite{MR1969635}. 

\begin{definition}\label{def:cofilteredcat}
	A category $\cD$ is {\em cofiltered} if every finite diagram in $\cD$ has a cone, equivalently:
	\begin{enumerate}
		\item $\cD$ is non-empty (the empty diagram admits a cone);
		\item for every pair of objects $d_1, d_2 \in \cD$, there exists an object $d_3 \in \cD$ and morphisms $d_3 \to d_1$ and $d_3 \to d_2$; and
		\item for every pair of parallel arrows in $\cD$, $f,g: d_2 \to d_1$, there exists a morphism $h:d_3 \to d_2$ which equalizes them, $fh = gh$.
	\end{enumerate}
\end{definition}

\begin{example} \label{example:directedset}
	A poset which is cofiltered is a (downward) {\em directed set}, i.e. there exists at least one element and for every pair of elements $x$ and $y$ there exists a $z$ with $z \leq x$ and $z \leq y$. 
\end{example}

\begin{definition}\label{def:proobject}
	Let $\cC$ be a category. The category of {\em pro-objects in~$\cC$}, denoted $\pro\cC$, is the category whose objects are pairs $(\cD,X)$, with $\cD$ a small cofiltered category and $X: \cD \to \cC$ a diagram, and whose morphisms from $X: \cD \to \cC$ to $Y:\cD' \to \cC$ are defined to be
	\begin{equation*}
		\hom(X,Y) : = \lim_{d \in \cD} \colim_{d' \in \cD'} \cC(X(d), Y(d'))
	\end{equation*}
	with limit and colimit taken in the category of sets. 
\end{definition}

The category $\cC$ is a full subcategory of $\pro\cC$ given by the trivial cofiltered limits (those in which $\cD =pt$).

The category of pro-objects in $\cC$ has an equivalent description in terms of {\em co-presheaves}, covariant functors from $\cC$ to $\set$. First recall the co-Yoneda embedding $\cC \to \Fun(\cC, \set)^\textrm{op}$ which is the fully-faithful functor from $\cC$ to the opposite category of copresheaves given by sending $c \in \cC$ to the corepresentable copresheaf 
$h^c = \cC(c,-)$.
%$h^c: c' \in \cC \mapsto \cC(c,c')$. 
Given a pro-object $X: \cD \to \cC$, we may define
\begin{equation*}
	LX := \lim_{d \in \cD} h^{X(d)}
\end{equation*}
to be the limit in $\Fun(\cC, \set)^\textrm{op}$ indexed by $\cD$. If this is written in terms of the category $\Fun(\cC, \set)$ of copresheaves, then $LX$ is a colimit of corepresentables. One may readily check that 
\begin{equation*}
	\Fun(\cC, \set)^\textrm{op}(LX, LY) \cong \pro\cC(X,Y)
\end{equation*} 
and as a consequence we have:

\begin{lemma}\label{lem:inclusiontoprocommuteswithcolims}
	The category $\pro\cC$ of pro-objects in $\cC$ is equivalent to the full subcategory of the opposite category of copresheaves $\Fun(\cC, \set)^\textrm{op}$ spanned by image of $L$, that is those presheaves which are isomorphic to cofiltered limits of corepresentables (in the category $\Fun(\cC, \set)$ they are filtered colimits of corepresentables). In particular the inclusion $\cC \to \pro\cC$ commutes with all colimits which exist in $\cC$.
	\qed
\end{lemma}

The last statement follows from the well-known fact that the co-Yoneda embedding commutes with all colimits which exist in $\cC$ (the equivalent dual statement is the fact that the Yoneda embedding commutes with all limits). 

\begin{definition}\label{def:cofinalfunctor}
	A functor $F:\cE \to \cD$ is {\em cofinal} if for every object $d \in \cD$ the over category $(F \downarrow d)$ is non-empty and connected. 
\end{definition}

The significance of cofinal functors for pro-objects is the following. 

\begin{lemma}\label{lem:cofinalreplacement}
	Let $X: \cD \to \cC$ be a pro-object of $\cC$, $\cE$ another small cofiltered category, and $F: \cE \to \cD$ a cofinal functor, then $X$ and $XF$ are isomorphic as pro-objects. 
\end{lemma}

\begin{proof}
	For any other pro-object $Y: \cK \to \cC$ we have
	\begin{align*}
		\pro\cC(Y, X) &\cong \lim_{k \in \cK} \colim_{d \in \cD} \cC(Y(k), X(d)) \\
		&\cong \lim_{k \in \cK} \colim_{e \in \cE} \cC(Y(k), XF(e)) \\
		&\cong \pro\cC(Y, XF)
	\end{align*}
	and hence, by the Yoneda lemma, $X \cong XF$. The middle isomorphism follows from the fact that colimits of sets are invariant under restriction along cofinal functors \cite[IX.3~Thm.~1]{MacLane71}. 
\end{proof}

\begin{lemma}\label{lem:cofinalrecognition}
	Let $F: \cE \to \cD$ be a functor between small cofiltered categories. Then $F$ is cofinal if and only if for ever $d \in \cD$ the over category $(F \downarrow d)$ is cofiltered. Moreover if $F$ is full, then $F$ is cofinal if and only if for every $d \in \cD$ the over category $(F \downarrow d)$ is non-empty, i.e. for all $d \in \cD$ there is an $e \in \cE$ and an arrow $F(e) \to d$.
\end{lemma}

\begin{proof}
	 Suppose first that $(F \downarrow d)$ is cofiltered for each $d \in \cD$. Then they are, in particular, non-empty and connected and so $F$ is cofinal. Conversely, suppose that $F$ is cofinal. Then for each $d \in \cD$, $(F \downarrow d)$ is non-empty. Consider a pair of parallel arrows in $(F \downarrow d)$. This consists of a pair of parallel arrows $f,g: e' \to e$ in $\cE$ together with a map $F(e) \to d$. Since $\cE$ is cofiltered there is an arrow $h:e'' \to e'$ which equalizes $f$ and $g$. Applying $F$ gives a morphism $F(e'') \to F(e)$ and hence promotes $h$ to a morphism in $(F \downarrow d)$ (which still equalizes $f$ and $g$). Thus $(F \downarrow d)$ satisfies the first and third property of a cofiltered category. We must show the that it also satisfies the second property. 
	 
So suppose that $e, e' \in (F \downarrow d)$ are two objects in the over category. We must find a new object $e''$ in the over category with morphisms $e'' \to e$ and $e'' \to e'$, also in the over category. Moreover we know that $(F \downarrow d)$ is connected, and hence there exists a finite zig-zag of morphisms in $(F \downarrow d)$ connecting $e$ and $e'$. Without loss of generality it sufficient to consider the case in which there exists an object $x \in (F \downarrow d)$ and arrows $e \to x$ and $e' \to x$ also in the over category. Since $\cE$ is cofiltered we can find an object $e''' \in \cE$ with morphisms $e''' \to e$ and $e''' \to e'$. Hence we get morphisms $F(e''') \to F(e) \to d$ and $F(e''') \to F(e') \to d$, but we do not know if these are the same morphism, and hence we have not yet constructed the desired object in $(F \downarrow d)$. Similarly we obtain two morphism $e''' \to e \to x$ and $e''' \to e' \to x$. Since $\cE$ is cofiltered, we may find a morphism $e'' \to e'''$ which equalizes these morphisms to $x$. Thus the two composites
\begin{align*}
	F(e'') \to F(e''') \to F(e) \to F(x) \to d \\
	F(e'') \to F(e''') \to F(e') \to F(x) \to d 
\end{align*}
agree and $(e'', F(e'') \to d)$ admits maps in $(F \downarrow d)$ to both $e$ and $e'$. This establishes the first claim.

Now suppose that $F: \cE \to \cD$ is full and that the over categories $(F \downarrow d)$ are non-empty for each $d \in \cD$. We wish to show that $F$ is cofinal, that is we wish to show that each $(F \downarrow d)$ is connected. Let $(e, F(e) \to d)$ and $(e', F(e') \to d)$ be two objects of $(F \downarrow d)$. Since $\cD$ is cofiltered, there exists an object $d' \in \cD$ and maps $d' \to F(e)$ and $d' \to F(e')$. Moreover, again since $\cD$ is cofiltered, we may find a new morphism $d'' \to d'$ such that the composites
\begin{align*}
	d'' \to d' \to F(e) \to d \\
	d'' \to d' \to F(e') \to d 
\end{align*}
agree. Next, since $(F \downarrow d'')$ is non-empty, there exists an object $e'' \in \cE$ and a morphism $F(e'') \to d''$. By composing this morphism with the previous ones we see that we have constructed an object $e'' \in \cE$ and a commutative diagram
\begin{center}
\begin{tikzpicture}
		\node (LT) at (0, 1.5) {$F(e'')$};
		\node (LB) at (0, 0) {$F(e)$};
		\node (RT) at (2, 1.5) {$F(e')$};
		\node (RB) at (2, 0) {$d$};
		\draw [->] (LT) -- node [left] {$f$} (LB);
		\draw [->] (LT) -- node [above] {$f'$} (RT);
		\draw [->] (RT) -- node [right] {$$} (RB);
		\draw [->] (LB) -- node [below] {$$} (RB);
		%\node at (0.5, 1) {$\ulcorner$};
		%\node at (1.5, 0.5) {$\lrcorner$};
\end{tikzpicture}
\end{center}  
Since $F$ is full, there exist morphisms $\tilde{f}: e'' \to e$ and $\tilde{f}': e'' \to e'$ lifting $f$ and $f'$. 
\end{proof}

The previous lemma makes it easy to recognize cofinal functors and hence allows us to manipulate the indexing category without changing the underlying pro-object. In our applications we will only need to consider pro-objects that are indexed by directed sets (see example~\ref{example:directedset}). This is no real loss of generality as every small cofiltered category admits a cofinal functor from a directed set, however many constructions are more naturally formulated using general cofiltered categories (see \cite{MR1969635}). 

\subsection{Halations}

Taking $\cC = \man$ to be the category of smooth manifolds yields the category of {\em pro-manifolds}. Here we mean manifold very generally and include all manifolds, possibly with boundary or corners. Many of our constructions will use pro-manifolds constructed from manifolds without boundary, but since we will want to consider maps into these from manifolds with corners, it seems simpler to includes these at the start.

\begin{example}\label{example:halation}
	Let $X$ be a manifold, possibly with boundary or corners, and let $X \subseteq Y$ be an embedding of manifolds. The collection of all codimension zero submanifolds $Z \subseteq Y$, without boundary, and which contain $X$, forms a directed set. We will let $\widehat{(X \subseteq Y)}$, or $\hat{X}$ when that doesn't cause confusion, denote the corresponding pro-manifold. There is a natural inclusion of pro-manifolds $X \hookrightarrow \widehat{(X \subseteq Y)}$. We will say $\widehat{(X \subseteq Y)}$ has codimension $k$ if $X$ has codimension $k$ in $Y$. 
%	From this there exist the following directed set	
%	Let . We will assume that $Y$ is without boundary, but that is unnecessary. This defines a cofiltered diagram of manifolds consisting of all codimension zero submanifolds of $Y$ without boundary which contain $X$. We will denote the corresponding pro-manifold $\hat{Y}_X$, and refer to pro-manifolds obtained this way as {\em halations}. We will say the halation has {\em codimension $k$} if $X$ has codimension $k$ in $Y$. There is a canonical inclusion of pro-manifolds $X \hookrightarrow \hat{Y}_X$, and the normal bundle $\nu_{X/Y}$ on $X$ only depends on the corresponding halation $\hat{Y}_X$.  There is a canonical halation $Y = Y \times \{0\} \hookrightarrow (\widehat{Y \times \R})_Y$ of codimension one. 
\end{example}

\begin{definition}\label{def:halation}
	Let $X$ be a manifold. A {\em halation} $(X, \hat{X})$ for $X$ is an inclusion of pro-manifolds $X \hookrightarrow \hat{X}$ which is isomorphic to an inclusion from Example~\ref{example:halation}. A map between manifolds with halations $(X, \hat{X})$ and $(Y, \hat{Y})$ is a commutative square of pro-manifolds. 
	\begin{center}
	\begin{tikzpicture}
		\node (LT) at (0, 1.5) {$X$};
		\node (LB) at (0, 0) {$\hat{X}$};
		\node (RT) at (2, 1.5) {$Y$};
		\node (RB) at (2, 0) {$\hat{Y}$};
		\draw [right hook->] (LT) -- node [left] {$$} (LB);
		\draw [->] (LT) -- node [above] {$$} (RT);
		\draw [right hook->] (RT) -- node [right] {$$} (RB);
		\draw [->] (LB) -- node [below] {$$} (RB);
		%\node at (0.5, 1) {$\ulcorner$};
		%\node at (1.5, 0.5) {$\lrcorner$};
	\end{tikzpicture}
	\end{center}
	Given an isomorphism $\hat{X} \cong \widehat{(X \subseteq Z)}$ we will say that $X \subseteq Z$ is a {\em representative} of the halation. Manifolds with halation form a category $\haloman$.
\end{definition}

In particular a halation is a structure on $X$ which has automorphisms. Halations will be our way of making precise the idea of equipping our manifolds with a `germ of neighborhoods of manifolds'. Essentially we will repeat the definition of bordisms given in the last section, but equip every manifold which occurs with a halation. Before describing this in detail, it will be useful to develop some basic tools for working with halations. 

First we record a few simple observations. If  $\widehat{X \subseteq Z} \cong \hat{X}$ represents a given halation and $X \subseteq Z_0 \subseteq Z$ with $Z_0$ any codimension zero submanifold of $Z$ without boundary, then $X \subseteq Z_0$ represents the same halation with induced isomorphism $\widehat{X \subseteq Z_0} \cong \widehat{X \subseteq Z}$. This is because the directed set of submanifolds of $Z_0$ is cofinal in the directed set of submanifolds of $Z$. In layman terms we can `shrink' representatives. In particular if $(X, \hat{X})$ is a codimension zero halation and $X$ itself has no boundary, then $\hat{X} \cong \widehat{(X \subseteq X)} \cong X$.

\begin{lemma} \label{lma:representmapsofhalation}
	Let $(f, \hat{f}): (X, \hat{X}) \to (Y, \hat{Y})$ be a map of manifolds with halations and let $X \subseteq Z$ and $Y \subseteq W$ be representatives of the halations $\hat{X}$ and $\hat{Y}$, respectively. Then there exist codimension zero submanifolds without boundary $ Z_0 \subseteq Z$ and $ W_0 \subseteq W$, containing $X$ and $Y$, respectively, and a commutative diagram of manifolds
	\begin{center}
	\begin{tikzpicture}
		\node (LT) at (0, 1.5) {$X$};
		\node (LB) at (0, 0) {$Z_0$};
		\node (RT) at (2, 1.5) {$Y$};
		\node (RB) at (2, 0) {$W_0$};
		\draw [right hook->] (LT) -- node [left] {$$} (LB);
		\draw [->] (LT) -- node [above] {$f$} (RT);
		\draw [right hook->] (RT) -- node [right] {$$} (RB);
		\draw [->] (LB) -- node [above] {$g$} (RB);
		%\node at (0.5, 1) {$\ulcorner$};
		%\node at (1.5, 0.5) {$\lrcorner$};
	\end{tikzpicture}
	\end{center}
which induces $(f, \hat{f})$, i.e. $\hat{f} = \hat{g}: \hat{X} \cong \widehat{(X \subseteq Z_0)} \to \widehat{(Y \subseteq W_0)} \cong \hat{Y}$.
\end{lemma}

\begin{proof}
	Suppose that $(X \subseteq Z)$ is a representative and $Z_0 \subseteq Z$ is a codimension zero submanifold containing $X$. Then by Lemma~\ref{lem:cofinalrecognition} the system of codimension zero submanifolds without boundary induced by $Z_0$ is cofinal in the system induced by $Z$. Thus by Lemma~\ref{lem:cofinalreplacement} $(X \subseteq Z)$ and $(X \subseteq Z_0)$ induce the same manifold with halation $(X, \hat{X})$, as we observed above. 
	
	Now suppose we are given a morphisms $(f, \hat{f}): (X, \hat{X}) \to (Y, \hat{Y})$ of manifolds with halation. The desired commutative square is given by a representative of the limit-of-colimits defining the map of pro-manifolds $\hat{X} \to \hat{Y}$, and the result now follows. 
\end{proof}

Let $E$ be a vector bundle over $X$. Then $X$ includes into $E$ via the zero section and, if $\partial X = \emptyset$, the inclusion $(X \subseteq E)$ induces a halation for $X$. If $X$ has boundary, then we do not get a canonical halation from $E$ alone. We must first choose a collar neighborhood for $\partial X \subseteq X$. This choice determines a smooth structure on the manifold without boundary
\begin{equation*}
	E \cup_{E|_{\partial X} \times \{0\}} E|_{\partial X} \times [0,1)
\end{equation*}
and hence also determines a halation for $X$. 

When $X$ is a manifold with corners we may defined a similar halation from a vector bundle, provided that there exists a cornered analog of a collar. This structure does not always exist for arbitrary manifolds with collars, but by Lemma~\ref{CollaringManifoldsWithCornersLemma} such collars exist for any ($m$-dimensional) $\langle n \rangle$-manifold $X$. Specifically, let $E$ be a rank $k$ vector bundle over $X$ and choose an $\langle n \rangle$-diagram of collarings as in Lemma~\ref{CollaringManifoldsWithCornersLemma}. Then there exists an $(m+k)$-dimensional manifold without boundary $\tilde{E}$ with an inclusion $X \to E$ such that for each $a = (a_0, \dots, a_{n-1}) \in \langle n \rangle$ the inclusion is locally isomorphic to 
\begin{equation*}
	\R^m_+(a') \times X(a) \hookrightarrow \R^m \times E|_{X(a)}
\end{equation*}
where $a' = (1,1, \dots, 1) - a$ is the complementary object and the smooth structure on $\tilde{E}$ is induced by the $\langle n \rangle$-diagram of collarings. The inclusion $X \subseteq \tilde{E}$ induces a halation for $X$. 

In all cases, when $X$ has no boundary, boundary, or corners, we will denote this halation by $(X, \hat{X}_E)$, though in the later two cases (where $\partial X \neq \emptyset$) this is understood to also depend on the choice of collar for $\partial X$. For any two choices of collar the resulting halations are isomorphic, but not canonically isomorphic.   

\begin{corollary}\label{cor:halationtubularnbhd}
	Let $(X, \hat{X})$ be a manifold with halation. Then there exists a vector bundle $E$ over $X$, of rank equal to the codimension of $(X, \hat{X})$, and an isomorphism $(X, \hat{X}) \cong (X, \hat{X}_E)$ which induces the identity on $X$. Moreover if $(f, \hat{f}): (X, \hat{X}) \to (Y, \hat{Y})$ is an isomorphism of manifolds with halation then there exists a vector bundle $E$ over $Y$ and isomorphisms $(X, \hat{X}) \cong (X, \hat{X}_{f^*E})$ and $(Y, \hat{Y}) \cong (Y, \hat{Y}_E)$ making the following square commute:
	\begin{center}
	\begin{tikzpicture}
			\node (LT) at (0, 1.5) {$(X, \hat{X})$};
			\node (LB) at (0, 0) {$(X, \hat{X}_{f^*E})$};
			\node (RT) at (2, 1.5) {$(Y, \hat{Y})$};
			\node (RB) at (2, 0) {$(Y, \hat{Y}_E)$};
			\draw [->] (LT) -- node [left] {$\cong$} (LB);
			\draw [->] (LT) -- node [above] {$(f, \hat{f})$} (RT);
			\draw [->] (RT) -- node [right] {$\cong$} (RB);
			\draw [->] (LB) -- node [below] {$$} (RB);
			%\node at (0.5, 1) {$\ulcorner$};
			%\node at (1.5, 0.5) {$\lrcorner$};
	\end{tikzpicture}
	\end{center}
	where the bottom map is induced by the natural map $f^*E \to E$. 
\end{corollary}

\begin{proof}
	Given a representative $(X \subseteq Z)$, we may replace $Z$ by a tubular neighborhood of $X$ in $Z$. The natural inclusion of the tubular neighborhood induces an isomorphism of halation. This establishes the first half of the corollary. The second half is similar and follows by pulling back the tubular neighborhood of $Y$. 
\end{proof}

\begin{corollary}
	The stabilization process $\widehat{(X \subseteq Z)} \mapsto (X \times \{0\} \subseteq Z \times \R) \widehat{\phantom{x}}$ is functorial for morphism of manifolds with halations (and in particular depends only on the underlying halation and not on the chosen representative). \qed
\end{corollary}

\begin{corollary}
	We may restrict halations: If $(X, \hat{X})$ is a halation and $X_0 \subseteq X$ is any submanifold, then there is a canonical halation $(X_0, \hat{X}|_{X_0})$ characterized by the property that $\hat{X}|_{X_0}\cong \widehat{(X_0 \subseteq Z)}$ for any representative $\hat{X} \cong \widehat{(X \subseteq Z)}$. 
	\qed
\end{corollary}

This allows us to introduce a restricted category $\halomand \subseteq \haloman$, 
whose objects are $d$-dimensional manifolds equipped with halations and whose morphisms are morphisms of manifolds with halations $(f, \hat{f}): (X, \hat{X}) \to (Y, \hat{Y})$ such that $f: X \to Y$ is the inclusion of a submanifolds and the induced map $\hat{f}: \hat{X} \to \hat{Y}|_X$ is an isomorphism.  

Moreover, halations and their maps are local. 

\begin{lemma}\label{lem:halationsarelocal}
	Let $X$ be a manifold, $\{ U_i \}$ a cover of $X$, and $(Y, \hat{Y})$ a manifold with halation. Let $U_{ij} = U_i \cap U_j$ and $U_{ijk} = U_i \cap U_j \cap U_k$ denote the corresponding intersections. 
	Then:
	\begin{enumerate}
		\item Suppose that for each $i$ we are given a halation $(U_i, \hat{U}_i)$, and that for each pair $i,j$ we are given an isomorphism of manifolds with halation $\hat{\phi}_{ij}: \hat{U}_i|_{U_{ij}} \cong \hat{U}_i|_{U_{ij}}$, such that the cocycle conditions $\hat{\phi}_{jk}\circ \hat{\phi}_{ij} = \hat{\phi}_{ik}$ and $\hat{\phi}_{ii}=id$ are satisfied. Then there exists a unique halation $(X, \hat{X})$ with isomorphisms $\hat{u}_i:\hat{X}|_{U_i} \cong \hat{U}_i$, such that $\hat{\phi}_{ij} = (\hat{u}_j|_{U_{ij}}) \circ (\hat{u}_i|_{U_{ij}})^{-1}$.  
		\item Suppose that $f: X \to Y$ is a smooth map and that $(X, \hat{X})$ is a halation. Suppose further that we are given a family of morphisms of manifolds with halation $(f|_{U_i}, \hat{f}_i): (U_i, \hat{X}|_{U_i}) \to (Y, \hat{Y})$ such that $\hat{f}_i|_{U_{ij}} = \hat{f}_j|_{U_{ij}}$. Then there exists a unique morphism of manifolds with halation $(f, \hat{f}): (X, \hat{X}) \to (Y,\hat{Y})$ such that $\hat{f}_i = \hat{f}|_{U_i}$. 
		\end{enumerate}
\end{lemma}

Before proving this lemma we will first recall another extremely useful lemma about covers of smooth manifolds.

\begin{lemma}\label{lem:manifoldfinitecover}
	Let $X$ be a smooth ($\langle m\rangle$-) manifold of dimension $n$ and let $\cU$ be an open cover of $X$. Then there exists another cover $\cV$ of $X$ such that:
	\begin{enumerate}
		\item The cover $\cV= \{ U_0, U_1, \dots, U_n\}$ is finite and consists of exactly $n+1$ elements;
		\item Each $U_i$ and each $k$-fold intersections of the $U_i$ consists of a disjoint union of open balls (each $\cong \R^n$);
		\item The closure in $X$ of each such ball is a closed ball (hence $\cong D^n$) which is contained in an open subset of $\cU$.   
	\end{enumerate} 
\end{lemma}

\begin{proof}
	First, it is well known that $X$ admits a triangulation, for example by Whitney's method. Moreover, by taking subdivisions if necessary we may arrange so that each of the closed star neighborhoods of each vertex of the triangulation is contained in some open subset of $\cU$. Now we subdivide the triangulation once more. The vertices of this subdivision correspond to the simplices (vertices, edges, faces, etc) of the previous triangulation. Thus each vertex is assigned a numerical label $i = 0, \dots, n$ according to the dimension of the corresponding simplex. The open star neighborhoods of any two vertices with the same label are disjoint from each other, and we let $U_i$ be the union of all the open star neighborhoods of the vertices with label $i$. 
\end{proof}

\begin{proof}[Proof of Lemma~\ref{lem:halationsarelocal}]	
We will first prove the first part of the lemma. By the above Lemma~\ref{lem:manifoldfinitecover}, we may assume, without loss of generality, that the cover $\{U_i\}$ is finite. Hence by induction it will be enough to consider the case when $X = U_0 \cup U_1$ is the union of two open sets. Let us choose representatives for the halations $(U_0, \hat{U}_0)$, $(U_1, \hat{U}_1)$, and their restriction to $U_{01} = U_0 \cap U_1$:  $(U_0 \subseteq Z_0)$, $(U_1 \subseteq Z_1)$, and $(U_{01} \subseteq Z_{01})$. We would like to glue these together to get (a representative of) a halation for $X$. Unfortunately the union $Z_0 \cup_{Z_{01}} Z_1$ will generally fail to be a manifold because a space formed in such a way may not be Hausdorff, though it will have local charts and be a manifold in all other respects.

To get around this problem we will have to judiciously shrink the manifolds $Z_i$. We recall the following lemma from \cite[Lma.~2.23]{MR2742432}: 
\begin{lemma}
	Let $W$ and $W'$ be manifolds and $A$ and open subset of each, then $W \cup_A W'$ is Hausdorff, and hence a manifold, if and only if the natural map $A \to W \times W'$ to the product sends $A$ to a closed subset. \qed
\end{lemma}

By Corollary~\ref{cor:halationtubularnbhd} we may choose the representatives to have a very particular form. Specifically we may assume that the halation $(U_0, \hat{U}_0)$ is induced by a vector bundles $E_0$ over $U_0$, that the halation $(U_1, \hat{U}_1)$ is induced by a vector bundle $E_1$ over $U_1$ and that the transition isomorphism is induced by an inclusion $i_1^*E_1 \subseteq i_0^*E_0$ realizing $i_1^*E_1$ as a tubular neighborhood of $U_{01}$ inside $i_0^*E_0$. Here $i_j : U_{01} \to U_j$ is the inclusion map. For simplicity we will assume that $\partial X = \emptyset$, and trust the reader to make the small changes needed in the case $X$ has boundary or corners. 

With the above choices the space $ E_0 \cup_{i^*_1E_1} E_1$ is Hausdorff and hence a manifold. In this case the halation induced by $X \subseteq  E_0 \cup_{i^*_1E_1} E_1$ satisfies the requirements of the lemma. To see this space is Hausdorff first note that since $X$ is a manifold (and hence Hausdorff itself) the map $U_{01} \to U_0 \times U_1$ is a closed inclusion (hence its image is closed). Thus the pullback, 
\begin{center}
\begin{tikzpicture}
		\node (LT) at (0, 1.5) {$i^*_1 E_1$};
		\node (LB) at (0, 0) {$U_{01}$};
		\node (RT) at (2, 1.5) {$U_0 \times E_1$};
		\node (RB) at (2, 0) {$U_0 \times U_1$};
		\draw [->] (LT) -- node [left] {$$} (LB);
		\draw [->] (LT) -- node [above] {$$} (RT);
		\draw [->] (RT) -- node [right] {$$} (RB);
		\draw [->] (LB) -- node [below] {$$} (RB);
		\node at (0.5, 1) {$\ulcorner$};
		%\node at (1.5, 0.5) {$\lrcorner$};
\end{tikzpicture}
\end{center}
which is the pullback of a vector bundle along a closed inclusion, is again a closed inclusion. Finally, the inclusion $U_0 \times E_1 \to E_0 \times E_1$ is closed, and hence so is the composite $i^*_1 E_1 \to E_0 \times E_1$. This establishes the first part of the lemma.

Now let us turn to the second part of the lemma. Let $Y \subseteq Z$ be a representative for $(Y, \hat{Y})$, and let $\cV = \{V_0, \dots, V_n\}$ be a cover of $X$ as constructed in Lemma~\ref{lem:manifoldfinitecover}. Since $V_i$ is a disjoint union of open balls, each subordinate to the original cover, we get induced maps $(f|_{V_i}, \hat{f}_i): (V_i, \hat{X}|_{V_i}) \to (Y, \hat{Y})$, which agree on the intersections of the $V_i$. 

By Corollary~\ref{cor:halationtubularnbhd} there exists a vector bundle $E$ over $X$ such that $(X, \hat{X}) \cong (X, \hat{X}_E)$. Moreover for each of the open sets $V_i$ there exists a tubular neighborhood $N_i \subseteq E|_{V_i}$ and a square
\begin{center}
\begin{tikzpicture}
		\node (LT) at (0, 1.5) {$V_i$};
		\node (LB) at (0, 0) {$N_i$};
		\node (RT) at (2, 1.5) {$Y$};
		\node (RB) at (2, 0) {$Z$};
		\draw [->] (LT) -- node [left] {$$} (LB);
		\draw [->] (LT) -- node [above] {$f|_{V_i}$} (RT);
		\draw [->] (RT) -- node [right] {$$} (RB);
		\draw [->] (LB) -- node [below] {$f_i$} (RB);
		%\node at (0.5, 1) {$\ulcorner$};
		%\node at (1.5, 0.5) {$\lrcorner$};
\end{tikzpicture}
\end{center}
inducing $\hat{f}_i$. Moreover we may insist that these tubular neighborhoods can be chosen as the restrictions of tubular neighborhoods defined on the original $U_i$, which support functions representing the original $\hat{f}_i$. Furthermore on intersections $V_{ij} = V_i \cap V_j$ there exists a tubular neighborhood $N_{ij} \subseteq N_i \cap N_j \subseteq E$ such that $f_i|_{N_{ij}} =f_j|_{N_{ij}}$. Since the closure $\overline{V}_i \cap \overline{V}_j$ is a union of closed balls, each contained in some $U_i$ from the original cover, we may further assume that the tubular neighborhoods $N_{ij}$ may be extended to tubular neighborhoods of the closure $\overline{V}_{ij}$ (while still satisfying the condition $f_i|_{N_{ij}} =f_j|_{N_{ij}}$).

Now we fix a metric on $E$. For each component $V_{i, \alpha}$ (an open ball) in $V_i$ we may choose a number $\epsilon_{i, \alpha} > 0$ such that the tubular neighborhood $N_{\epsilon_{i, \alpha}}$ of $E|_{V_{i, \alpha}}$ is contained in $N_i$. The fact that we may choose this as a fixed number, rather than a function depending on the point in $V_{i, \alpha}$ follows from the fact that $N_i$ extends to the closure $\overline{V}_{i, \alpha}$. On this latter compact set we have a uniform bound. A similar statement holds for the components of $V_{ij}$. 

Given a strictly positive smooth function $\lambda: X \to \R_+$ ( $\lambda(x) > 0$ for all $x \in X$), we can form the open $\lambda$-disk bundle $D_\lambda(E) = \{ v \in E \; | \; ||v|| < \lambda(p(v)) \}$ where $p: E \to X$ is the projection. By the above considerations, and because the cover $\cV$ is finite (hence locally finite), we may choose a strictly positive smooth function with the property that for each $V_i$ and $V_{ij}$ we have inclusions $D_\lambda(E)|_{V_i} \subseteq N_i \subseteq E|_{V_i}$ and $D_\lambda(E)|_{V_{ij}} \subseteq N_{ij} \subseteq E|_{V_{ij}}$, respectively. 
The functions $f_i$ now patch together to yield the global map $\tilde{f}$ in the following commutative square:
\begin{center}
\begin{tikzpicture}
		\node (LT) at (0, 1.5) {$X$};
		\node (LB) at (0, 0) {$D_\lambda(E)$};
		\node (RT) at (2, 1.5) {$Y$};
		\node (RB) at (2, 0) {$Z$};
		\draw [->] (LT) -- node [left] {$$} (LB);
		\draw [->] (LT) -- node [above] {$f$} (RT);
		\draw [->] (RT) -- node [right] {$$} (RB);
		\draw [->] (LB) -- node [below] {$\tilde{f}$} (RB);
		%\node at (0.5, 1) {$\ulcorner$};
		%\node at (1.5, 0.5) {$\lrcorner$};
\end{tikzpicture}
\end{center}
This represents the desired map of manifolds with halation. 

Moreover this morphism of manifolds with halations is unique. If $(f, \hat{f}')$ is any other such a map, then there exists a $\lambda'$ and a map $\tilde{f}':D_{\lambda'}(E) \to Z$ making the above square commute. By assumption the restriction of $\tilde{f}'$ to each $U_i$ represents $(f|_{U_i}, \hat{f}_i)$, and hence agrees with $f_i$ in a, perhaps smaller, tubular neighborhood of $U_i$ in $E|_{U_i}$. Consequently $\tilde{f}'$ agrees with $\tilde{f}$ on some small neighborhood of $X$, and hence $(f, \hat{f}') = (f, \hat{f})$.
\end{proof}

\begin{corollary}%\label{cor:}
	The forgetful functor $\halomand \to \man^d$, from $d$-dimensional manifolds with halation to their underlying manifolds is
	 is a stack (Def.~\ref{def:stack}). \qed
\end{corollary}

 \subsection{The Haloed bordism bicategory} \label{sec:haloedBordd}
We may now introduce an improved version of the bordism bicategory which incorporates the structure of halations. The symmetric monoidal bicategory $\bord_d$ is constructed exactly like $\cob_d$ except that every manifold is equipped with a $d$-dimensional halation surrounding it. Thus the $d$-dimensional bordisms are equipped with codimension zero halations, the $(d-1)$-dimensional bordisms are equipped with codimension $1$ halations, etc. In later sections we will enhance this bordism bicategory further by equipping the halations with n additional structures such as orientations, principal bundles, framings, or spin structures, but for now we will only consider the bare halation.

By Corollary~\ref{cor:halationtubularnbhd} we know that up to (non-canonical) isomorphism there is only one codimension zero halation for these $d$-dimensional bordisms. Futhermore this corollary states that codimension one halations are isomorphic to ones of the form $(X, \hat{X}_E)$ for some 1-dimensional vector bundle $E$. 
However if $X$ is to be realized as the restriction to the boundary of a codimension zero halation on a $d$-dimensional bordism, then the vector bundle $E$ must be \emph{trivial} (in fact trivialized by, say, the inward pointing normal vector). 
Similarly for a halation on a $(d-2)$-manifold to be obtained by restricting a codimension zero halation on a $d$-dimensional bordisms to its corner, then it also must be induced by a trivial vector bundle.  Thus when we construct the bordism bicategory we will not consider all possible halations on our bordisms, but only those which are isomorphic to ones obtained from trivial vector bundles. Furthermore we will specify the trivialization (up to isotopy) be equipping our halations with {\em co-orientations}.   

Let $(X, \hat{X})$ be a halation for the manifold  $X$. The pro-manifold $\hat{X}$ is given by a cofiltered limit of open neighborhoods $N_\alpha$ of $X$ in some larger manifold $Z$. A {\em co-orientation} for $\hat{X}_E$ will mean a compatible family of orientations for the normal bundles of $X \subseteq N_\alpha$. This is a structure which is intrinsic to the pro-manifold $\hat{X}$ and does not depend on our choice of ambient manifold $Z$.

Before formally defining the structures which make up $\bord_d$, let us first consider the operation of gluing haloed 1-bordisms to show how these small additions make the composition both canonical and well-defined. Let $(Y_0, \hat{Y}_0)$ and $(Y_1, \hat{Y}_1)$ be $(d-1)$-manifolds equipped with halations which are isomorphic to $(Y_i, (\hat{Y}_i)_\R)$, i.e., to the halation induced by the trivial line bundle. Furthermore, we assume that each $(Y_i, \hat{Y}_i)$ is equipped with a choice of co-orientation. 

By a haloed bordism from $(Y_0, \hat{Y}_0)$ to $(Y_1, \hat{Y}_1)$ we will mean compact $d$-manifold with a codimension zero halation $(W, \hat{W})$, together with a decomposition and isomorphism
\begin{equation*}
	(\partial W, \hat{W}|_{\partial W}) \cong (Y_0, \hat{Y}_0) \sqcup (Y_1, \hat{Y}_1).
\end{equation*}
For each component of the boundary of $W$, the restricted halation may be given a co-orientation in one of two canonical ways. We have that $(\partial W, \hat{W}|_{\partial W}) \cong (\partial W, \widehat{\partial W}_\nu)$, where $\nu$ is the (trivial) normal bundle of $\partial W$ in $W$. We may orient (i.e. trivialize) this trivial bundle on $\partial W$ by either using the inward pointing normal vector or the outward pointing vector. We make the following convention: The co-orientation of the incoming boundary  $(Y_0, \hat{Y}_0)$ must agree with the inward pointing co-orientation, the co-orientation of the outgoing boundary $(Y_1, \hat{Y}_1)$ must agree with the outward pointing co-orientation. 

Given haloed 1-bordisms $(W_0, \hat{W}_0)$ and $(W_1, \hat{W}_1)$ from $(Y_0, \hat{Y}_0)$ to $(Y_1, \hat{Y}_1)$, and $(Y_1, \hat{Y}_1)$ to $(Y_2, \hat{Y}_2)$, respectively, we wish to construct a new haloed 1-bordism $W_1 \circ W_0 = (W, \hat{W})$ from $(Y_0, \hat{Y}_0)$ to $(Y_1, \hat{Y}_1)$. We will describe the composition in detail in the case where $Y_0 = Y_2 = \emptyset$. The general case is no harder, except in bookkeeping. So without loss of generality we wish to construct the composite  $W_1 \circ W_0$ when $Y = Y_1 \cong \partial W_0 \cong \partial W_1$.

In Section~\ref{SectionBordBicat1} we were able to form the composite 1-bordism $W = W_0 \cup_{Y} W_1$ as a topological manifold, but the smooth structure was not canonically determined. The composition was only unique up to canonical invertible 2-bordism, not up to canonical diffeomorphism. In the current situation however, the additional data coming from the halations is sufficient to determine a canonical smooth structure on the closed topological manifold $W = W_0 \cup_{Y} W_1$. In fact the pushout $\hat{W}:= \hat{W}_1 \cup_{\hat{Y}_1} \hat{W}_{0}$ 
in pro-manifolds exists and is represented by an ordinary smooth manifold. In this case, since $\partial W = \emptyset$, the induced map $W \to \hat{W}$ is an isomorphism and determines the smooth structure on $W$. In the general case where $Y_0$ and $Y_2$ are non-empty this will not be an isomorphism due to the behavior at those ends, however it will still induce a canoncial smooth structure on $W$.

In order to justify these claims properly, we will need some results about computing finite colimits in categories of pro-objects. For this we refer the reader to \cite{MR1969635} which has a discussion of how to compute general colimits and limits in categories of pro-objects. By  Corollary~\ref{cor:halationtubularnbhd} we may assume that the halations $(W_i, \hat{W}_i)$ are represented by inclusions
\begin{equation*}
	W_i \hookrightarrow W_i \cup_{\partial W_i} \partial W_i \times \R_{\geq 0}
\end{equation*}
where the target is endowed with a smooth structure. This doesn't effect the pro-manifold $\hat{W}_i$. Similarly $\hat{Y}$ is induced from $Y \hookrightarrow Y \times \R$. 

Next, without altering the pro-manifolds $\hat{W}_0$, $\hat{W}_1$, and $\hat{Y}$, we may assume that they arise as cofiltered limits of manifolds indexed over the same cofiltered indexing category $I$, and that this description is compatible with the pushout diagram in a sense which we will explain shortly. There are several ways to arrange for this, and a general method is described in \cite{MR1969635}. We will sketch an alternative method in this case.  

The pro-manifolds $\hat{Z}_X$ coming from halations are most naturally indexed by the cofiltered category of submanifolds of $Z$ containing $X$. However, in light of Lemma~\ref{lem:cofinalreplacement}, we may obtain the same pro-manifold by replacing this by any cofinal subset of these submanifolds. So for example in the case at hand we may take the cofinal collection were we simply `shrink the ends' of $W_0 \cup_{\partial W_0} \partial W_0 \times \R_{\geq 0}$, $\partial W_1 \times \R_{\leq 0} \cup_{\partial W_1} W_1$, and $Y \times \R$. 

More concretely we may take $I = \N$ and then associate to each $i$ the manifolds:
\begin{align*}
	i & \mapsto W_0 \cup_{_\partial W_0} \partial W_0 \times [0, \frac{1}{i}) \\
	i & \mapsto \partial W_1 \times (- \frac{1}{i}, 0] \cup_{\partial W_1} W_1 \\
	i & \mapsto Y \times (- \frac{1}{i}, \frac{1}{i})
\end{align*}
In this case they will in particular have the same cofiltered indexing category $I$, and for each index $i \in I$, we get a pushout diagram of ordinary manifolds. 

Because of our convention on co-orientations, the pushout of each diagram for fixed $i$ exists and is the smooth ordinary manifold $W_1 \cup_{Y_1} W_0$ with smooth structure near $Y_1$ determined by the embedding of $Y_1 \times (- \frac{1}{i}, \frac{1}{i})$. In particular for each index $i$ the pushout is the same, and we will call this smooth manifold $\hat{W}$.

Let $K = * \leftarrow * \rightarrow *$ be the indexing category for pushouts. What we have achieved is a diagram 
\begin{equation*}
	F: K \times I \to \man
\end{equation*}
of manifolds indexed by the product of the category $K$ and a cofiltered category $I$, such that for each fixed $k \in K$ the limit, $\lim_{i \in I} F(k,i)$, in pro-manifolds is one of $\hat{W}_0$, $\hat{W}_1$, or $\hat{Y}$, depending on the index $k \in K$. The colimit we are interested in computing is given by
\begin{equation*}
	\colim_{K} \left( \lim_{I} F(k,i) \right).
\end{equation*}
In the terminology of  \cite{MR1969635} functor $F$ is called a {\em level representation} of the pushout diagram, and the existence and uniqueness of level representations of general diagrams is discussed thoroughly in that paper. The main technical result we need is \cite[Thm.~6.1]{MR1969635} which states that cofiltered limits and finite colimits commute in the category of pro-objects. Thus we have
\begin{align*}
	\colim_{K} \left(\lim_{I} F(k,i) \right) & \cong \lim_{I} \left(\colim_{K}  F(k,i) \right)\\
	& \cong \lim_{I} \hat{W} \\
	& \cong \hat{W}.
\end{align*}
The first isomorphism follows by commuting the cofiltered limit and the finite colimit, the second follows because the embedding of manifolds into pro-manifolds commutes with colimits which exist in manifolds (Lemma~\ref{lem:inclusiontoprocommuteswithcolims}), and the last follows because we are taking the limit of a constant diagram. This establishes the claim that the pushout $\hat{W}_1 \cup_{\hat{Y}} \hat{W}_{0}$ in pro-manifolds exists and is represented by the ordinary smooth manifold $\hat{W}$, the composite of $(W_1, \hat{W}_1)$ and $(W_0, \hat{W}_0)$  

For bordisms with higher codimensional corners, we will have to equip those corners with a series of co-oriented halations. Specifically a haloed manifold $X$ of dimension $(d-k)$ we will consist of a sequences of inclusions 
\begin{equation*}
	X \subseteq \hat{X}_0 \subseteq \hat{X}_1 \subseteq \cdots \subseteq \hat{X}_k
\end{equation*}
where $\hat{X}_i$ is a co-oriented halation of dimension $(d-k+i)$ (codimension $i$). Note that if $X$ is closed then we necessarily have $X = \hat{X}_0$, and in this case this first halation is redundant. This structure can be visualized as in Figure~\ref{FigHalations}, where the shadings signify the co-orientations.  

\begin{figure}[ht]
\begin{center}
\begin{tikzpicture}[thick]
	\node [circle, minimum width = 3cm] (A) at (0,0) {};
	\node (B) at (A.30) {};
	\node (C) at (A.200) {};

	\fill [black!20]  (C.center) to [out = 45, in = 180] (A.center) to [out = 0, in = 225] (B.center) 
		arc (30 : 200 : 1.5cm); 
	\draw [very thick, black] (A.center) to [out = 0, in = 225] (B.center);
	\draw [ultra thick, dotted, black] (C.center) to [out = 45, in = 180] (A.center);
	\draw [dashed] (A.center) circle (1.5cm);
	\node [circle, fill=black,inner sep=1.5pt] at (A.center) {};
\end{tikzpicture}
\caption{Haloes}
\label{FigHalations}
\end{center}
\end{figure}

The construction of the symmetric monoidal bicategory $\bord_d$ will mimic the construction of $\cob_d$. Specifically we will use the methods described in Section~\ref{sec:symmonpseudo2cat}. We will first construct a symmetric monoidal pseudo double category $\BORD_d$. This pseudo double category will be fibrant, and so by Theorem~\ref{thm:ShulmanThm} the underlying bicategory $\cH(\BORD_d)$ will be symmetric monoidal. The construction is virtually identical to the one in Section~\ref{sec:symmonbicatcobord}, except that now every manifold is equipped with a sequence of co-oriented halations. 

The symmetric monoidal category $(\BORD_d)_0$, the category of objects, has objects $(Y, \hat{Y}_1, \hat{Y}_2)$ which consist of a closed compact $(d-2)$-manifold $Y$, together with a pair of co-oriented halations with inclusions 
\begin{equation*}
	Y \subseteq \hat{Y}_1 \subseteq \hat{Y}_2. 
\end{equation*}
The halation $(Y, \hat{Y}_1)$ is codimension one, and the halation $(Y, \hat{Y}_2)$ is codimension two. The morphisms of $(\BORD_d)_0$ consist of isomorphisms preserving the co-orientations. 

The symmetric monoidal category $(\BORD_d)_1$ is the category of arrows. Its objects and morphisms are more lengthly to define, but they are analogous to the objects and morphisms of $(\COB_d)_1$, only equipped with halations. 

\begin{figure}[ht]
\begin{center}
\begin{tikzpicture}[thick]
	\node (A) at (2,4) {};
	\node (B) at (2, 3) {};
	\node (X) at (6, 3) {};
	\node (Y) at (6,2) {};
	\fill [black!20] (A.center) -- +(-1, 0) arc (180 : 90 : 1cm) to [ out = 0, in = 180] (X.center)
		arc (90 : 0 : 1cm) -- (Y.center) to [out = 180, in = 0] (A.center);
	\draw [dashed] (B.center) arc (270 : 90: 1cm)  to [ out = 0, in = 180] (X.center)
		arc (90: -90: 1cm) to [out = 180, in = 0] (B.center);
	\draw (A.center) -- +(-1, 0) (A.center) to [ out = 0, in = 180]  (Y.center) -- +(1,0); 
	\node [circle, fill=black,inner sep=1pt] at (A.center) {};
	\node [circle, fill=black,inner sep=1pt] at (Y.center) {};

	\node (C) at (0,0) {};
	\node (D) at (8, 0) {};
	\fill [black!20] (C.center) -- +(-1, 0) arc (180 : 0: 1cm) -- (C.center);
	\fill [black!20] (D.center) -- +(-1, 0) arc (180 : 0: 1cm) -- (D.center);
	\draw [dashed] (C.center) circle (1cm);
	\draw [dashed] (D.center) circle (1cm);
	\draw (C.center) -- +(1,0) (D.center) -- +(1,0);
	\draw [ultra thick, dotted] (C.center) -- +(-1,0) (D.center) -- +(-1,0);
	\node [circle, fill=black,inner sep=1pt] at (C.center) {};
	\node [circle, fill=black,inner sep=1pt] at (D.center) {};
	
	\draw [->] (0, 1.5) to [out = 90, in = 225] node [above left] {$f_\text{in}$} (1, 3);
	\draw [->] (8, 1.5) to [out = 90, in = 0] node [above right] {$f_\text{out}$} (7.5, 2);
\end{tikzpicture}	
\caption{A Haloed 1-Bordisms}
\label{FigHaloed1Bordisms}
\end{center}
\end{figure}
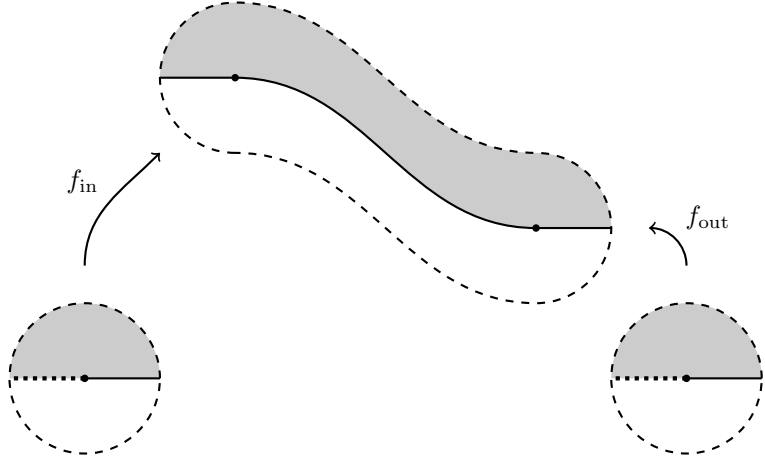

\begin{definition}
	The objects of $(\BORD_d)_1$ are {\em haloed 1-bordisms}. They consist of:
	\begin{itemize}
		\item a pair of objects of $(\BORD_d)_0$, 	$(Y_\textrm{in}, \hat{Y}_{1, \textrm{in}}, \hat{Y}_{2, \textrm{in}})$ and $(Y_\textrm{out}, \hat{Y}_{1, \textrm{out}}, \hat{Y}_{2, \textrm{out}})$;
		\item A haloed $(d-1)$-manifold $(\Sigma, \hat{\Sigma}_0, \hat{\Sigma}_1)$, that is a $(d-1)$-dimensional manifold with boundary $\Sigma$ together with a codimension zero halation $\hat{\Sigma}_0$, and a co-oriented codimension one halation $\hat{\Sigma}_1$, together with inclusions $\Sigma \subseteq \hat{\Sigma}_0 \subseteq \hat{\Sigma}_1$; and
		\item a decomposition of $\partial \Sigma = \partial_\textrm{in} \Sigma \sqcup  \partial_\textrm{out} \Sigma$ and a pair of isomorphisms of haloed manifolds
		\begin{align*}
			f_\textrm{in}: (\partial_\textrm{in} \Sigma,  \hat{\Sigma}_0|_{\partial_\textrm{in}},  \hat{\Sigma}_1|_{\partial_\textrm{in}} ) &\cong (Y_\textrm{in}, \hat{Y}_{1, \textrm{in}}, \hat{Y}_{2, \textrm{in}}) \\
			f_\textrm{out}: (\partial_\textrm{out} \Sigma,  \hat{\Sigma}_0|_{\partial_\textrm{out}},  \hat{\Sigma}_1|_{\partial_\textrm{out}} ) &\cong (Y_\textrm{out}, \hat{Y}_{1, \textrm{out}}, \hat{Y}_{2, \textrm{out}}).
		\end{align*}
		Moreover, the we require that these isomorphisms preserve the co-orientations. Here $\hat{\Sigma}_0|_{\partial_\textrm{in}}$ is equipped with a co-orientation using the inward pointing normal direction, while $\hat{\Sigma}_0|_{\partial_\textrm{out}}$ is co-oriented using the outward normal direction. 
	\end{itemize}
\end{definition}
\noindent A typical example of an object of $(\BORD_2)_1$ is depicted in Figure~\ref{FigHaloed1Bordisms}.

The morphisms of $(\COB_d)_1$ consisted of $\langle 2\rangle$-manifolds whose boundary faces were divided into incoming/outgoing components. There were four sides altogether. Two of these were identified with objects of $(\COB_d)_1$, that is $1$-bordisms. However the remaining sides were identified with trivial bordisms, induced from a diffeomorphism of $(d-2)$-manifolds. 
To generalize this to $(\BORD_d)_1$ we need an analog of these trivial bordisms but equipped with halations.  

To this end suppose that 
\begin{equation*}
	f: (Y^\textrm{in}, \hat{Y}_{1}^\textrm{in}, \hat{Y}_{2}^\textrm{in}) \to (Y^\textrm{out}, \hat{Y}_{1}^\textrm{out}, \hat{Y}_{2}^\textrm{out})
\end{equation*}
is an isomorphism in $(\BORD_d)_0$. Let $\{0\} \times I$ denote the unit segment of the $y$-axis in $\R^2$, and consider the halation induced by $Y^\textrm{in} \times \{0\} \times I \hookrightarrow Y^\textrm{in} \times \R^2$. This halation is naturally co-oriented in the positive $x$-direction. The boundary inherits a sequence of haltions:
\begin{align*}
	Y^\textrm{in} \times \{0\} \times \{1\} &\hookrightarrow  Y^\textrm{in} \times \R \times \{1\} \hookrightarrow Y^\textrm{in} \times \R^2 \\
	Y^\textrm{in} \times \{0\} \times \{0\} &\hookrightarrow  Y^\textrm{in} \times \R \times \{0\} \hookrightarrow Y^\textrm{in} \times \R^2
\end{align*}
where the first inclusion is co-oriented using the negative $y$-direction, and the second inclusion uses the standard orientation of $\R^2$. These choices make the boundary components of $Y^\textrm{in} \times I$ into objects of $(\BORD_d)_0$, which are naturally isomorphic via the identification $Y^\textrm{in} \times \{0\} \times \{1\} \cong Y^\textrm{in} \times \{0\} \times \{0\}$

Moreover there exists an isomorphism in $(\BORD_d)_0$: 
\begin{equation*}
	k: (Y^\textrm{in}, \hat{Y}_{1}^\textrm{in}, \hat{Y}_{2}^\textrm{in}) \cong (Y^\textrm{in} \times \{0\} \times \{1\}, \hat{Y}^\textrm{in}_\R, \hat{Y}^\textrm{in}_{\R^2})
\end{equation*}
 where the latter is equipped with the above co-oriented halations. The choice of this isomorphism is not unique, but once chosen we may parametrize the boundary of $(Y^\textrm{in} \times I, \widehat{Y^\textrm{in} \times \R^2})$ via:
\begin{align*}
	k: (Y^\textrm{in}, \hat{Y}_{1}^\textrm{in}, \hat{Y}_{2}^\textrm{in}) &\cong (Y^\textrm{in} \times \{0\} \times \{1\}, \hat{Y}^\textrm{in}_\R, \hat{Y}^\textrm{in}_{\R^2}) \\
	k \circ f^{-1}:  (Y^\textrm{out}, \hat{Y}_{1}^\textrm{out}, \hat{Y}_{2}^\textrm{out}) &\cong (Y^\textrm{in} \times \{0\} \times \{0\}, \hat{Y}^\textrm{in}_\R, \hat{Y}^\textrm{in}_{\R^2}).
\end{align*}
We will call $(Y^\textrm{in} \times I, \widehat{Y^\textrm{in} \times \R^2})$ together with such a parametrization of its boundary a {\em vertically trivial 1-bordism} induced from $f$.

\begin{definition}\label{def:haloedbordisms12}
	Let $\Sigma^\textrm{in} = (\Sigma^\textrm{in},\hat{\Sigma}_0^{\textrm{in}}, \hat{\Sigma}_1^{\textrm{in}})$ and $\Sigma^\textrm{out} = (\Sigma^\textrm{out},\hat{\Sigma}_0^{\textrm{out}}, \hat{\Sigma}_1^{\textrm{out}})$ be objects of $(\BORD_d)_1$. 
	The morphisms of $(\BORD_d)_1$ from $\Sigma^\textrm{in}$ to $\Sigma^\textrm{out}$ are equivalences classes consisting of:
	\begin{itemize}
		\item a pair of isomorphisms of objects in $(\BORD_d)_0$:
		\begin{align*}
			h_\textrm{in}: (Y_\textrm{in}^\textrm{in}, \hat{Y}_{1,\textrm{in}}^\textrm{in}, \hat{Y}_{2,\textrm{in}}^\textrm{in}) &\cong (Y_\textrm{in}^\textrm{out}, \hat{Y}_{1,\textrm{in}}^\textrm{out}, \hat{Y}_{2,\textrm{in}}^\textrm{out}) \\	
			h_\textrm{out}: (Y_\textrm{out}^\textrm{in}, \hat{Y}_{1,\textrm{out}}^\textrm{in}, \hat{Y}_{2,\textrm{out}}^\textrm{in}) &\cong (Y_\textrm{out}^\textrm{out}, \hat{Y}_{1,\textrm{out}}^\textrm{out}, \hat{Y}_{2,\textrm{out}}^\textrm{out})
		\end{align*}
		\item A $\langle 2\rangle$-manifold $S$ with a codimension zero halation $\hat{S}$;
		\item A decomposition and isomorphism of the boundary of $S$ equipped with co-oriented halation:
		\begin{equation*}
			(\partial_0S, \hat{S}|_{\partial_0 S}) \cong (\Sigma^\textrm{in}, \hat{\Sigma}_1^{\textrm{in}}) \sqcup (\Sigma^\textrm{out}, \hat{\Sigma}_1^{\textrm{out}})
		\end{equation*}
		where the co-orientation of $(\Sigma^\textrm{in}, \hat{\Sigma}_1^{\textrm{in}})$ agrees with the inward pointing normal direction and the co-orientation of $(\Sigma^\textrm{out}, \hat{\Sigma}_1^{\textrm{out}})$ agrees with the outward pointing normal direction. 
		\item  A decomposition and isomorphism of the boundary of $S$ equipped with co-oriented halation:
		\begin{equation*}
			(\partial_1S, \hat{S}|_{\partial_1 S}) \cong (Y_\textrm{in}^\textrm{in} \times I, \widehat{Y_\textrm{in}^\textrm{in} \times \R^2}) \sqcup (Y_\textrm{out}^\textrm{in} \times I, \widehat{Y_\textrm{out}^\textrm{in} \times \R^2}) 
		\end{equation*}
		where the co-orientation of $(Y_\textrm{in}^\textrm{in} \times I, \widehat{Y_\textrm{in}^\textrm{in} \times \R^2})$ agrees with the inward pointing direction and the co-orientation of $(Y_\textrm{out}^\textrm{in} \times I, \widehat{Y_\textrm{out}^\textrm{in} \times \R^2})$ agrees with the outward pointing co-orientation. 
	\end{itemize}
From this data we obtain the morphisms 
\begin{align*}
	k_\textrm{in}^\textrm{in}: (Y_\textrm{in}^\textrm{in}, \hat{Y}_{1, \textrm{in}}^\textrm{in}, \hat{Y}_{2, \textrm{in}}^\textrm{in}) &\cong (Y_\textrm{in}^\textrm{in} \times \{0\} \times \{1\}, \hat{Y}^\textrm{in}_{\textrm{in},\R}, \hat{Y}^\textrm{in}_{\textrm{in},\R^2}) \\
	k_\textrm{out}^\textrm{in}: (Y_\textrm{out}^\textrm{in}, \hat{Y}_{1, \textrm{out}}^\textrm{in}, \hat{Y}_{2, \textrm{out}}^\textrm{in}) &\cong (Y_\textrm{out}^\textrm{in} \times \{0\} \times \{1\}, \hat{Y}^\textrm{in}_{\textrm{out},\R}, \hat{Y}^\textrm{in}_{\textrm{out},\R^2}) \\
	k_\textrm{in}^\textrm{out}: (Y_\textrm{in}^\textrm{out}, \hat{Y}_{1, \textrm{in}}^\textrm{out}, \hat{Y}_{2, \textrm{in}}^\textrm{out}) &\cong (Y_\textrm{in}^\textrm{in} \times \{0\} \times \{0\}, \hat{Y}^\textrm{in}_{\textrm{in},\R}, \hat{Y}^\textrm{in}_{\textrm{in},\R^2}) \\
	k_\textrm{out}^\textrm{out}: (Y_\textrm{out}^\textrm{out}, \hat{Y}_{1, \textrm{out}}^\textrm{out}, \hat{Y}_{2, \textrm{out}}^\textrm{out}) &\cong (Y_\textrm{out}^\textrm{in} \times \{0\} \times \{0\}, \hat{Y}^\textrm{in}_{\textrm{out},\R}, \hat{Y}^\textrm{in}_{\textrm{out},\R^2}) 
\end{align*}
as the unique isomorphisms of manifolds with halation making the following diagrams commute:	
\begin{center}

\end{center}
Finally, we further require that the composites $(k^\textrm{out}_\textrm{in})^{-1} \circ k^\textrm{in}_\textrm{in}$ and $(k^\textrm{out}_\textrm{out})^{-1} \circ k^\textrm{in}_\textrm{out}$ agree with $h_\textrm{in}$ and $h_\textrm{out}$, respectively. 
	
Two such morphisms are considered equivalent if they are isomorphic relative to the boundary. 
\end{definition}

The symmetric monoidal functors $S$, $T$, and $U$ are defined in exactly the same way as they were for $\COB_d$ in Section~\ref{sec:symmonbicatcobord}. Both the composition of morphisms in $(\BORD_d)_1$ (the vertical composition) and the functor $\odot$ (the horizontal composition) are defined as the pushout of manifolds with halation. This works exactly the same as the composition of haloed 1-bordisms which we described in detail above, and there are naturally induced coherence isomorphisms making $\BORD_d$ into a symmetric monoidal pseudo double category. 

The bicategory $\bord_d$ is defined to be the horizontal bicategory $\cH(\BORD_d)$. Its objects and 1-morphisms are the objects of $(\BORD_d)_0$ and objects of $(\BORD_d)_1$, respectively. Its 2-morphisms are precisely the 2-morphisms of $(\BORD_d)_1$ in which the isomorphisms $h_\textrm{in}$ and $h_\textrm{out}$ are identities. The proof of Theorem~\ref{thm:symmonbicatcob}, which established that the symmetric monoidal pseudo double category $\COB_d$ was fibrant, goes through verbatim in this case as well. Thus we have established:

\begin{theorem}%\label{thm:}
	The bicategory $\bord_d$ of haloed $(d-2)$-manifolds, haloed ($(d-1)$-dimensional) $1$-bordisms, and haloed ($d$-dimensional) 2-bordisms is a symmetric monoidal bicategory. \qed
\end{theorem}

Forgetting the halo structure gives a forgetful symmetric monoidal homomorphism $\bord_d \to \cob_d$. 

\begin{proposition}
The symmetric monoidal homomorphism $\bord_d \to \cob_d$ is an equivalence of symmetric monoidal bicategories. 
\end{proposition}

\begin{proof}
By Theorem~\ref{WhiteheadforSymMonBicats}, we only need to check that this homomorphism is essentially surjective on objects, essentially full on 1-morphisms, and fully-faithful on 2-morphisms. The first two of these are clear. Since each $(d-2)$-manifold and bordism may be equipped with the necessary halations to make it an object (or 1-morphism) of $\bord_d$, it is actually surjective on objects and 1-morphisms (not just essentially surjective). Moreover the same argument applies to 2-morphisms and shows that the forgetful functor is full on 2-morphisms. It is faithful on 2-morphisms because two haloed 2-bordisms are isomorphic if and only if their underlying non-haloed 2-bordisms are isomorphic (see Cor.~\ref{cor:halationtubularnbhd}).\footnote{Note, however, that given an isomorphism of underlying 2-bordisms the isomorphism of haloed 2-bordisms is not canonical. Nevertheless we still obtain a bijection on isomorphisms classes. } 
\end{proof}

\section{The Structured Bordism Bicategory}

Now that we have introduced the machinery of haloes into our bordism bicategory, we can almost effortless add extra geometric or topological structures to our bordism bicategory. Structures will only live on the $d$-dimensional part of a haloed manifold. Since composition in the bordism bicategory is given by gluing manifolds, what ever structure we wish to consider must also have a gluing property. For simplicity we will require that the category of $d$-manifolds with structure be a stack over the category of $d$-manifolds, although one could just as easily use a higher stack.  

% Since we want to construct a symmetric monoidal bicategory of bordisms with structure, these structures must behave well with respect to disjoint union. 

Let $\man^d$ be the category of $d$-manifolds without boundary, with embeddings as morphism. $\man^d$ is a Grothendieck site with covering families the usual notion of covering, see Section \ref{SectSitesAndSheaves} in the appendix. Let $\cF$ be a stack over $\man^d$, in the sense of Appendix \ref{AppSymMonStacks}. An object of $\cF$ will be called an $\cF$-manifold. We will think of an $\cF$-manifold as a pair $(M, s \in \cF(M))$, i.e., as a manifold $M$ together with an $\cF$-structure, \index{bordism!$\cF$-structure} $s$, on $M$. Replacing every occurrence of ``$d$-manifold without boundary'' in the definitions of haloed manifolds with ``$\cF$-manifold'' one obtains the notion of an {\em $\cF$-haloed manifold}. The definitions of haloed 0-bordisms and 1-bordisms go through verbatim for $\cF$-haloed manifolds, as well. 

The composition of $\cF$-1-bordisms can be defined, just as before, as a pushout. The fact that $\cF$ is a stack ensures that we will always be able to equip the glued bordism with an $\cF$-structure and that this $\cF$-structure is unique up to unique isomorphism.

The definition of haloed 2-bordism is slightly more problematic. The first main difficulty is that in the definition itself, for each of the object 0-bordisms, $Y$, we need to extended it to $Y \times I$ together with a halation. This then parametrizes one `side' of the 2-bordism. While this is possible for bare topological haloed manifolds, it is not clear how to proceed for general $\cF$-manifolds. Even if this can be surmounted, as is the case in many examples, a related problem comes up when we try to define the vertical composition of 2-bordisms. We are able to glue the 2-bordisms without difficulty, but to turn the result into a new 2-bordism, we needed to choose an identification $I \cup_{pt} I \cong I$. In the case of haloed 2-bordisms, the result was independent of which choice of identification we made.

One approach to the above problem is alter the geometry of 2-bordisms, `crushing' the $Y \times I$ portion of the bordism so that the bordism resembles a bigon, possibly with cusps at the corresponding corners. This is the approach taken \cite{Henriques2014}, for example. However with this approach, in order to have a well defined horizontal composition one must necessarily allow degenerate non-manifold objects with singularities. 

Another approach is to replace the bordism bicategory with a bordism {\em double category} in which the vertical morphisms are also allowed to be non-trivial bordisms. This was the approach used in \cite{Morton07}, and it solves both of the above problems. In this situation there are two kinds of 1-morphism, vertical and horizontal, and a  2-morphism sits inside a square of 1-morphism. It is straightforward to use $\cF$-haloed manifolds in this context. In this approach, one must also deal with the symmetric monoidal structure. From the point of view of topological field theories, however, this approach is not completely satisfactory as many potential targets are not double categories, but classical bicategories. 

Even in the case without structures, these attempted solutions result in objects which are distinct from the bordism bicategory we have previously described. Moreover we will ultimately like a clean classification theorem based on generators and relations.  It is not completely clear to what extent the results of the previous two chapters can be used in the above approaches. We will not pursue them further. 

Instead we will introduce an additional structure which makes it possible to surmount the above problems. This is a very mild limitation. For example, in many cases the collection of structures on a given manifold forms a topological space $\cF(M)$ with an action by the discrete group $Diff(M)$. The corresponding stack will have groupoid fiber equivalent to the fundamental groupoid of $\cF(M)$. The above problems disappear if this action by $Diff(M)$ on $\cF(M)$ extends to a continuous action, i.e., it respects the topology of $Diff(M)$. This covers most cases of interest, for example it includes cases considered in \cite{MR2713365}.

A topological $\cF$-structure (defined below) will be a  stack, $\cF$, over $\man^d$ equipped with some additional structure. This structure will be required to satisfy certain properties which allow us to mimic the construction of the haloed bordism bicategory, essentially without any significant change. 

The first piece of structure we need is a functor:
\begin{align*}
	\times I: \{ \cF\text{-0-bordisms} \} &
			\to \{ \cF\text{-}1\text{-bordisms}\},	
\end{align*}
which lifts the usual functor from unstructured 0-bordisms to unstructured 1-bordisms. Given this structure we can now {\em define} $\cF$-2-bordisms, exactly as we did in Definition \ref{def:haloedbordisms12}, replacing every occurrence of ``haloed $d$-manifold'' with ``$\cF$-haloed $d$-manifold''. Morphisms of $\cF$-2-bordisms are defined in the obvious way, and there is a functorial horizontal composition given by gluing $\cF$-2-bordisms. Just as in the haloed case, there are canonical coherent natural isomorphisms realizing associativity of this composition.

The second piece of structure that we will need is a functor,
\begin{equation*}
	\times I: \{ \cF\text{-1-bordisms} \} 
			\to \{  \cF\text{-}2\text{-bordisms}\}, 
\end{equation*}
which lifts the usual functor between unstructured bordisms. This allows us to define $\cF$-pseudo isotopy for $\cF$-1-bordisms. $\cF$-pseudo isotopy can be defined in general, but we will only need it in a vary particular case. 

\begin{definition} \label{DefnFPseudoIsotopy}
	Let $Y$ be an $\cF$-0-bordism, and consider two automorphisms of the partially haloed $\cF$-1-bordism $Y \times I$, $f_0$ and $f_1$. In particular $f$ and $g$ are the identity near the ends of $Y \times I$. This determines an $\cF$-2-bordism, $Y \times I \times I$. We may modify this $\cF$-2-bordism by pre-composing by $f_0$ and $f_1$, as in the diagram:
	
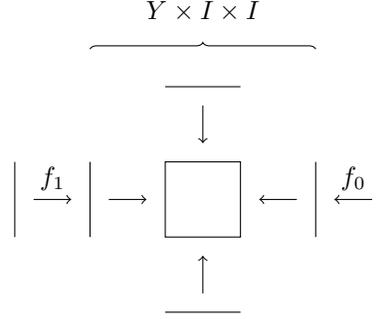
\begin{center}
\begin{tikzpicture}
		\draw (0, 1) -- +(0, -1) (1, 1) -- +(0, -1)  (2,0) rectangle (3,1) 
		(4, 1) -- +(0, -1) (5, 1) -- +(0, -1); 
		\draw (2, 2) -- +(1,0) (2, -1) -- +(1,0);
		\draw [->] (0.25, 0.5) -- node [above] {$f_1$} +(0.5, 0);
		\draw [->] (1.25, 0.5) -- +(0.5, 0);
		\draw [<-] (3.25, 0.5) -- +(0.5, 0);
		\draw [<-] (4.25, 0.5) -- node [above] {$f_0$} +(0.5, 0);
		\draw [->] (2.5, 1.75) -- +(0,-0.5);
		\draw [->] (2.5, -0.75) -- +(0, 0.5);
		\draw[decorate, decoration=brace] (1, 2.5) -- +(3, 0);
		\node at (2.5, 3) {$Y \times I \times I$};
\end{tikzpicture}
\end{center}
We say that $f_0$ and $f_1$ are {\em $\cF$-pseudo isotopic} \index{$\cF$-pseudo isotopic} if there exists an isomorphism between this new $\cF$-2-bordism and the original $\cF$-2-bordism $Y \times I \times I$. 
\end{definition}

\begin{definition} \label{DefTopologicalStructure}
	A stack, $\cF$, over $\man^d$ \index{stack!symmetric monoidal} \index{fibered category!stack} will be called {\em topological} if it is equipped with functors,
	\begin{align*}
		%\times I: \{ \cF\text{-0-bordisms} \} & 
		%	\to \{ \text{Partially Haloed } \cF\text{-}1\text{-bordisms}\} \\
		\times I: \{ \cF\text{-0-bordisms} \} & 
			\to \{  \cF\text{-}1\text{-bordisms}\} \\	
		\times I: \{ \cF\text{-1-bordisms} \} & 
		\to \{ \cF\text{-}2\text{-bordisms}\}
	\end{align*}
lifting the usual $\times I$ functors, such that
\begin{enumerate}
\item For any two automorphisms, $f_0$ and $f_1$ of $Y \times I$, as in Definition \ref{DefnFPseudoIsotopy}, if the underlying morphisms of unstructured bordisms are pseudo-isotopic, then $f_0$ and $f_1$ are $\cF$-pseudo-isotopic.  
\item If $\Sigma $ is an $\cF$-2-bordism with source and target $\cF$-0-bordisms $Y_0$ and $Y_1$, then the $\cF$-2-bordism $(Y_1 \times I \times I) \circ \Sigma \circ (Y_0 \times I \times I)$ is isomorphic to $\Sigma$. 
\end{enumerate}
\end{definition}

 The properties defining a topological $\cF$-structure allow us to define an associative vertical composition and  to provide associators and unitors for the horizontal composition of $\cF$-1-bordism. These structures and the proof of associativity are exactly as for unstructured haloed manifolds, and so we omit the details of the proof. This yields  a bicategory $\bord_d^\cF$ \index{$\bord_d^\cF$} whose objects are $\cF$-0-bordisms, whose 1-morphisms are $\cF$-1-bordisms, and whose 2-morphisms are $\cF$-2-morphisms. Since $\cF$ was assumed to be a stack, we again get an induced structure realizing $\bord_d^\cF$ as a symmetric monoidal bicategory. There are many examples of topological $\cF$-structures: orientations, spin structures, $G$-principal bundles, etc. The main techniques used in proving the classification of unoriented 2-dimensional topological field theories carry over to the setting of topological $\cF$-structures. After proving the unoriented classification, we demonstrate this by adapting our techniques to the oriented case, as well. 

\begin{remark}
The unoriented topological field theories play an especially important role among all structured topological field theories. They are universal in the sense that they give rise to a field theory with {\em any} topological structure $\cF$. More precisely the unoriented structure corresponds to $\cF = \man^d $, and thus is {\em terminal} in the bicategory of topological structures. In particular, for any topological structure $\cF$ there is the forgetful symmetric monoidal homomorphism $\bord_d^\cF \to \bord_d$. Thus, for any target symmetric monoidal bicategory $\sC$, there is a natural homomorphism 
\begin{equation*}
	\symbicat(\bord_d, \sC) \to \symbicat( \bord_d^\cF, \sC).
\end{equation*}
This expresses the universality of unoriented theories. 
\end{remark}

\section{The Unoriented Classification} \label{SectUnorientedClassification}

\noindent The goal of this section is to establish the following theorem:

\begin{theorem}[Classification of Unoriented Topological Field Theories] \label{TheoremUnorientedClassification}
The unoriented bordism bicategory $\bord_2$ has the  generators and relations presentation as a symmetric monoidal bicategory depicted in Figures \ref{MainUnOrientedTheormGenFig} and \ref{MainUnOrientedTheormRelFig}.

%\begin{figure}[ht]

%\end{figure}

\begin{figure}[ht]
\begin{center}

\caption{Generators of $\bord_2$}
\label{MainUnOrientedTheormGenFig}
\end{center}
\vspace{0.5cm}

\begin{center}
%
};

\begin{scope}

%%% Cusp Inversion AI
\begin{scope}[xshift = -2cm, yshift = -2.5cm]
\draw (0.5, 1.2) -- (-0.2,1.2) -- (-0.2,-0.5) -- (0.2, -0.5);
\draw (0.5, 1.2) arc (90: -90: 0.3cm and 0.1cm) arc (90: 270: 0.3cm and 0.1cm) -- (0.8, 0.8);
\draw [densely dashed] (0.2, -0.5) -- (0.5,-0.5) arc (90: -90: 0.3cm and 0.1cm) arc (90: 180: 0.3cm and 0.1cm);
\draw (0.2, -0.8) arc (180: 270: 0.3cm and 0.1cm) -- (0.8, -0.9);
\draw (0.8, 0.8) -- (1.2, 0.8) -- (1.2, -0.9) -- (0.8, -0.9);
\node (C) at (0.5, 0.4) {}; 
\node (D) at (0.5, -0.1) {};
\draw  plot [smooth] coordinates{(0.8,1.1) (0.8,0.9)}; 
\draw [densely dashed] plot [smooth] coordinates{ (0.8,0.9) (0.6, 0.6) (C.center)}; 
\draw  plot [smooth] coordinates{(0.2,0.9) (0.2,0.8) (0.4, 0.6) (C.center)}; 
\draw  plot [smooth] coordinates{(0.2,-0.8) (0.2,-0.6) (0.4, -0.4) (D.center)}; 
\draw [densely dashed] plot [smooth] coordinates{(0.8,-0.6) (0.8,-0.5) (0.6, -0.3) (D.center)}; 
\end{scope}

%%% Cusp Inversion AII
\begin{scope}[xshift = 1cm, yshift = -2.5cm]
\draw (0.5, 1.2) -- (-0.2,1.2) -- (-0.2,-0.5) -- (0.2, -0.5); %(0.5, -0.5) ;
\draw (0.5, 1.2) arc (90: -90: 0.3cm and 0.1cm) arc (90: 270: 0.3cm and 0.1cm) -- (0.8, 0.8);
\draw [densely dashed] (0.2, -0.5) -- (0.5,-0.5) arc (90: -90: 0.3cm and 0.1cm) arc (90: 180: 0.3cm and 0.1cm);
\draw (0.2, -0.8) arc (180: 270: 0.3cm and 0.1cm) -- (0.8, -0.9);
\draw (0.8, 0.8) -- (1.2, 0.8) -- (1.2, -0.9) -- (0.8, -0.9);
\node (C) at (0.5, 0.4) {};
\draw (0.2, 0.9) -- (0.2, -0.8) (0.8, 1.1) -- (0.8, -0.6);
\end{scope}

\node at (0, -2.5) {$=$};
\node at (0, -5) {$=$};

%%% Cusp Inversion BI
\begin{scope}[xshift = -2cm, yshift = -4.2cm]
\draw (1.2, 0.4) -- (-0.2,0.4) -- (-0.2,-0.5) -- (0.2, -0.5);
\draw [densely dashed] (0.2, -0.5) -- (0.5,-0.5) arc (90: -90: 0.3cm and 0.1cm) arc (90: 180: 0.3cm and 0.1cm);
\draw (0.2, -0.8) arc (180: 270: 0.3cm and 0.1cm) -- (0.8, -0.9);
\draw  (1.2, 0.4) -- (1.2, -0.9) -- (0.8, -0.9);
\draw (1.2, -0.9) -- (1.2, -1.6) -- (-0.2, -1.6) -- (-0.2, -0.5); 
\node (D) at (0.5, -0.1) {};
\draw  plot [smooth] coordinates{(0.2,-0.8) (0.2,-0.6) (0.4, -0.4) (D.center)}; 
\draw [densely dashed] plot [smooth] coordinates{(0.8,-0.6) (0.8,-0.5) (0.6, -0.3) (D.center)}; 
\node (C) at (0.5, -1.3) {}; 
\draw  plot [smooth] coordinates{(0.8,-0.6) (0.8,-0.7)}; 
\draw [densely dashed] plot [smooth] coordinates{ (0.8,-0.7) (0.6, -1) (C.center)}; 
\draw  plot [smooth] coordinates{(0.2,-0.8) (0.2,-0.9) (0.4, 0.-1.1) (C.center)}; 
\end{scope}

%%% Cusp Inversion BII
\begin{scope}[xshift = 1cm, yshift = -5cm]
\draw  (0.8, 0.8) to [out = 180, in = 0] (0.2, 1.2) -- (-0.2,1.2) -- (-0.2,-0.5) -- (0.2, -0.5) to [out = 0, in = 180] (0.8, -0.9);
%\draw (0.5, 1.2) arc (90: -90: 0.3cm and 0.1cm) arc (90: 270: 0.3cm and 0.1cm) -- (0.8, 0.8);
%\draw [densely dashed] (0.2, -0.5) -- (0.5,-0.5) arc (90: -90: 0.3cm and 0.1cm) arc (90: 180: 0.3cm and 0.1cm);
%\draw (0.2, -0.8) arc (180: 270: 0.3cm and 0.1cm) -- (0.8, -0.9);
\draw (0.8, 0.8) -- (1.2, 0.8) -- (1.2, -0.9) -- (0.8, -0.9);

\end{scope}

\end{scope}

%%%%%%%%%
%Symmetry Relation new
\begin{scope}[xshift = -7cm, yshift = -12.5cm, xscale = 1]
	\draw (.1,-0.9) -- (0, -0.9) --(0,1.3) 
		-- (0.6, 1.3) -- (0.7, 1.2) -- (0.8, 1.3)
		-- (1.2, 1.3) -- (1.3, 1.2) -- (1.4, 1.3)
		-- (2.1,1.3) arc (90: -90:  0.3cm and 0.1cm) 
		-- (1.4, 1.1) -- (1.3, 1.2) -- (1.2, 1.1)
		 -- (0.8, 1.1) -- (0.7, 1.2) -- (0.6, 1.1)
		-- (0.1, 1.1) -- (0.1,-1.1)
		-- (2.1, -1.1) arc (-90: 0:  0.3cm and 0.1cm)  
		-- (2.3, -0.5) -- (2.4, 0.1) 
		(0.1, 0.2) -- (0, 0.2)
	  (0.8, 1.1) -- (0.7, 1.2) -- (0.6, 1.1) -- (0.1, 1.1) -- (0.1,0)
	-- (0.6, 0) -- (0.7, 0.1) -- (0.8, 0)
	 -- (2.1, 0) 	 arc (-90: 0:  0.3cm and 0.1cm)  to (2.3, 0.6) -- (2.4, 1.2); 
\draw [densely dashed] (0.1,0.2) -- (0.6, 0.2) -- (0.7, 0.1) -- (0.8, 0.2) 
	-- (2.1, 0.2) arc (90: 0:  0.3cm and 0.1cm);
\draw (1.3, 1.2) to [out = -90, in = 180] (2.3, 0.6);
\draw (0.7, 1.2) -- (0.7, 0.1) to [out = -90, in = 180] (2.3, -0.5);
\draw [densely dashed] (0.1,-0.9) -- (2.1, -0.9) arc (90: 0:  0.3cm and 0.1cm);
\node at (3, 0) {$=$};	
% identity
\begin{scope}[xshift = 3.5cm]
\draw (.1,-0.9) -- (0, -0.9) --(0,1.3) 
	-- (0.6, 1.3) -- (0.7, 1.2) -- (0.8, 1.3)
	-- (1.2, 1.3) -- (1.3, 1.2) -- (1.4, 1.3)
	-- (2.1,1.3) arc (90: -90:  0.3cm and 0.1cm) 
	-- (1.4, 1.1) -- (1.3, 1.2) -- (1.2, 1.1)
	 -- (0.8, 1.1) -- (0.7, 1.2) -- (0.6, 1.1)
	-- (0.1, 1.1) -- (0.1,-1.1)
	 -- (2.1, -1.1) arc (-90: 0:  0.3cm and 0.1cm)  
	-- (2.4, 1.2)
	(0.7, 1.2) arc (180: 360: 0.3cm and 0.7cm);
\draw [densely dashed] (.1, -0.9) 
	-- (2.1, -0.9) arc (90: 0: 0.3cm and 0.1cm);
\end{scope}
\end{scope}

%Symmetry Relation new2
\begin{scope}[xshift = -8cm, yshift = -12.5cm, xscale = -1]
	\draw (.1,-0.9) -- (0, -0.9) --(0,1.3) 
		-- (0.6, 1.3) -- (0.7, 1.2) -- (0.8, 1.3)
		-- (1.2, 1.3) -- (1.3, 1.2) -- (1.4, 1.3)
		-- (2.1,1.3) arc (90: -90:  0.3cm and 0.1cm) 
		-- (1.4, 1.1) -- (1.3, 1.2) -- (1.2, 1.1)
		 -- (0.8, 1.1) -- (0.7, 1.2) -- (0.6, 1.1)
		-- (0.1, 1.1) -- (0.1,-1.1)
		-- (2.1, -1.1) arc (-90: 0:  0.3cm and 0.1cm)  
		-- (2.3, -0.5) -- (2.4, 0.1) 
		(0.1, 0.2) -- (0, 0.2)
	  (0.8, 1.1) -- (0.7, 1.2) -- (0.6, 1.1) -- (0.1, 1.1) -- (0.1,0)
	-- (0.6, 0) -- (0.7, 0.1) -- (0.8, 0)
	 -- (2.1, 0) 	 arc (-90: 0:  0.3cm and 0.1cm)  to (2.3, 0.6) -- (2.4, 1.2); 
\draw [densely dashed] (0.1,0.2) -- (0.6, 0.2) -- (0.7, 0.1) -- (0.8, 0.2) 
	-- (2.1, 0.2) arc (90: 0:  0.3cm and 0.1cm);
\draw (1.3, 1.2) to [out = -90, in = 180] (2.3, 0.6);
\draw (0.7, 1.2) -- (0.7, 0.1) to [out = -90, in = 180] (2.3, -0.5);
\draw [densely dashed] (0.1,-0.9) -- (2.1, -0.9) arc (90: 0:  0.3cm and 0.1cm);
\node at (3, 0) {$=$};	
% identity
\begin{scope}[xshift = 3.5cm]
\draw (.1,-0.9) -- (0, -0.9) --(0,1.3) 
	-- (0.6, 1.3) -- (0.7, 1.2) -- (0.8, 1.3)
	-- (1.2, 1.3) -- (1.3, 1.2) -- (1.4, 1.3)
	-- (2.1,1.3) arc (90: -90:  0.3cm and 0.1cm) 
	-- (1.4, 1.1) -- (1.3, 1.2) -- (1.2, 1.1)
	 -- (0.8, 1.1) -- (0.7, 1.2) -- (0.6, 1.1)
	-- (0.1, 1.1) -- (0.1,-1.1)
	 -- (2.1, -1.1) arc (-90: 0:  0.3cm and 0.1cm)  
	-- (2.4, 1.2)
	(0.7, 1.2) arc (180: 360: 0.3cm and 0.7cm);
\draw [densely dashed] (.1, -0.9) 
	-- (2.1, -0.9) arc (90: 0: 0.3cm and 0.1cm);
\end{scope}
\end{scope}
%%%%%%%%%

%Symmetry Relation B
\begin{scope}[xshift = -8cm, yshift = -3.5cm]
\draw (.1,.2) -- (0, 0.2)  --(0,1.3)-- (0.6, 1.3) -- (0.7, 1.2) -- (0.8, 1.3)
	 -- (1.3,1.3)  arc (90: -90:  0.3cm and 0.1cm) 
	 -- (0.8, 1.1) -- (0.7, 1.2) -- (0.6, 1.1) -- (0.1, 1.1) -- (0.1,0)
	 -- (1.3, 0) 	 arc (-90: 0:  0.3cm and 0.1cm)  to (1.5, 0.6) -- (1.6, 1.2); 
\draw [densely dashed] (0.1,0.2) 
	-- (1.3, 0.2) arc (90: 0:  0.3cm and 0.1cm);
\draw (0.7, 1.2) to [out = -90, in = 180] (1.5, 0.6);
	\draw [densely dashed] (0.1,-0.9)  -- (0.6, -0.9) -- (0.7, -1) -- (0.8, -0.9)
	-- (1.3, -0.9) arc (90: 0:  0.3cm and 0.1cm);
\draw  (.1,-0.9) -- (0, -0.9) --(0,0.2)
	(0.1, 0) -- (0.1,-1.1)	-- (0.6, -1.1) -- (.7, -1) -- (0.8, -1.1)
	 -- (1.3, -1.1)  arc (-90: 0:  0.3cm and 0.1cm) to (1.5, -0.5) -- (1.6, 0.1)
	 (0.7, -1) to [out = 90, in = 180] (1.5, -0.5); 	
\node at (2, 0) {$=$};	
% identity
\begin{scope}[xshift = 2.5cm]
\draw (.1,-0.9) -- (0, -0.9) --(0,1.3) 
	-- (0.6, 1.3) -- (0.7, 1.2) -- (0.8, 1.3)
	-- (1.3,1.3) arc (90: -90:  0.3cm and 0.1cm) 
	 -- (0.8, 1.1) -- (0.7, 1.2) -- (0.6, 1.1)
	-- (0.1, 1.1) -- (0.1,-1.1)
	-- (0.6, -1.1) -- (.7, -1) -- (0.8, -1.1)
	 -- (1.3, -1.1) arc (-90: 0:  0.3cm and 0.1cm)  -- (1.6, 1.2)
	 (0.7, -1) -- (0.7, 1.2);
\draw [densely dashed] (.1, -0.9) 
	-- (0.6, -0.9) -- (0.7, -1) -- (0.8, -0.9)
	-- (1.3, -0.9) arc (90: 0: 0.3cm and 0.1cm);
\end{scope}
\end{scope}

%Symmetry Relation B2
\begin{scope}[xshift = -9cm, yshift = -3.5cm, xscale = -1]
\draw (.1,.2) -- (0, 0.2)  --(0,1.3)-- (0.6, 1.3) -- (0.7, 1.2) -- (0.8, 1.3)
	 -- (1.3,1.3)  arc (90: -90:  0.3cm and 0.1cm) 
	 -- (0.8, 1.1) -- (0.7, 1.2) -- (0.6, 1.1) -- (0.1, 1.1) -- (0.1,0)
	 -- (1.3, 0) 	 arc (-90: 0:  0.3cm and 0.1cm)  to (1.5, 0.6) -- (1.6, 1.2); 
\draw [densely dashed] (0.1,0.2) 
	-- (1.3, 0.2) arc (90: 0:  0.3cm and 0.1cm);
\draw (0.7, 1.2) to [out = -90, in = 180] (1.5, 0.6);
	\draw [densely dashed] (0.1,-0.9)  -- (0.6, -0.9) -- (0.7, -1) -- (0.8, -0.9)
	-- (1.3, -0.9) arc (90: 0:  0.3cm and 0.1cm);
\draw  (.1,-0.9) -- (0, -0.9) --(0,0.2)
	(0.1, 0) -- (0.1,-1.1)	-- (0.6, -1.1) -- (.7, -1) -- (0.8, -1.1)
	 -- (1.3, -1.1)  arc (-90: 0:  0.3cm and 0.1cm) to (1.5, -0.5) -- (1.6, 0.1)
	 (0.7, -1) to [out = 90, in = 180] (1.5, -0.5); 	
\node at (2, 0) {$=$};	
% identity
\begin{scope}[xshift = 2.5cm]
\draw (.1,-0.9) -- (0, -0.9) --(0,1.3) 
	-- (0.6, 1.3) -- (0.7, 1.2) -- (0.8, 1.3)
	-- (1.3,1.3) arc (90: -90:  0.3cm and 0.1cm) 
	 -- (0.8, 1.1) -- (0.7, 1.2) -- (0.6, 1.1)
	-- (0.1, 1.1) -- (0.1,-1.1)
	-- (0.6, -1.1) -- (.7, -1) -- (0.8, -1.1)
	 -- (1.3, -1.1) arc (-90: 0:  0.3cm and 0.1cm)  -- (1.6, 1.2)
	 (0.7, -1) -- (0.7, 1.2);
\draw [densely dashed] (.1, -0.9) 
	-- (0.6, -0.9) -- (0.7, -1) -- (0.8, -0.9)
	-- (1.3, -0.9) arc (90: 0: 0.3cm and 0.1cm);
\end{scope}
\end{scope}

% Symmetry Relation A
\begin{scope}[xshift = -8cm, yshift = -6.25cm]
\draw (.1,.2) -- (0, 0.2) --(0,1.3) -- (1.3,1.3)  arc (90: -90:  0.3cm and 0.1cm) 
	 -- (0.1, 1.1) -- (0.1,0) -- (0.6, 0) -- (.7, 0.1) -- (0.8, 0) -- (1.3, 0) 
	 arc (-90: 0:  0.3cm and 0.1cm) 	 to  (1.5, 0.6) -- (1.6, 1.2); 
\draw [densely dashed] (0.1,0.2) 
	-- (0.6, 0.2) -- (0.7, 0.1) -- (0.8, 0.2)
	-- (1.3, 0.2) arc (90: 0:  0.3cm and 0.1cm);
\draw (0.7, 0.1) to [out = 90, in = 180] (1.5, 0.6);
\draw [densely dashed] (0.1,-0.9) -- (1.3, -0.9) arc (90: 0:  0.3cm and 0.1cm);
\draw (.1,-0.9) -- (0, -0.9) 	--(0,0.2)
	(0.1, 0) -- (0.1,-1.1) --(1.3, -1.1)  arc (-90: 0:  0.3cm and 0.1cm) 
	 to (1.5, -0.5)  -- (1.6, 0.1) 
	 (.7, 0.1) to [out = -90, in = 180] (1.5, -0.5);
\node at (2, 0) {$=$};	
% identity
\begin{scope}[xshift = 2.5cm]
\draw (.1,-0.9) -- (0, -0.9) --(0,1.3) -- (.3,1.3) arc (90: -90:  0.3cm and 0.1cm) -- (0.1, 1.1) 
	-- (0.1,-1.1) -- (.3, -1.1) arc (-90: 0:  0.3cm and 0.1cm)  -- (.6, 1.2);
\draw [densely dashed] (.1, -0.9) -- (.3, -0.9) arc (90: 0: 0.3cm and 0.1cm);
\end{scope}
\end{scope}

% Symmetry Relation A2
\begin{scope}[xshift = -9cm, yshift = -6.25cm, xscale = -1]
\draw (.1,.2) -- (0, 0.2) --(0,1.3) -- (1.3,1.3)  arc (90: -90:  0.3cm and 0.1cm) 
	 -- (0.1, 1.1) -- (0.1,0) -- (0.6, 0) -- (.7, 0.1) -- (0.8, 0) -- (1.3, 0) 
	 arc (-90: 0:  0.3cm and 0.1cm) 	 to  (1.5, 0.6) -- (1.6, 1.2); 
\draw [densely dashed] (0.1,0.2) 
	-- (0.6, 0.2) -- (0.7, 0.1) -- (0.8, 0.2)
	-- (1.3, 0.2) arc (90: 0:  0.3cm and 0.1cm);
\draw (0.7, 0.1) to [out = 90, in = 180] (1.5, 0.6);
\draw [densely dashed] (0.1,-0.9) -- (1.3, -0.9) arc (90: 0:  0.3cm and 0.1cm);
\draw (.1,-0.9) -- (0, -0.9) 	--(0,0.2)
	(0.1, 0) -- (0.1,-1.1) --(1.3, -1.1)  arc (-90: 0:  0.3cm and 0.1cm) 
	 to (1.5, -0.5)  -- (1.6, 0.1) 
	 (.7, 0.1) to [out = -90, in = 180] (1.5, -0.5);
\node at (2, 0) {$=$};	
% identity
\begin{scope}[xshift = 2.5cm]
\draw (.1,-0.9) -- (0, -0.9) --(0,1.3) -- (.3,1.3) arc (90: -90:  0.3cm and 0.1cm) -- (0.1, 1.1) 
	-- (0.1,-1.1) -- (.3, -1.1) arc (-90: 0:  0.3cm and 0.1cm)  -- (.6, 1.2);
\draw [densely dashed] (.1, -0.9) -- (.3, -0.9) arc (90: 0: 0.3cm and 0.1cm);
\end{scope}
\end{scope}

%%%%%%%%%
% Gluing-Cup Relation A
\begin{scope}[xshift = 2cm, yshift = -10cm, xscale = 1]
\draw (0,1) -- (0.4, 1) -- (0.5, 1.1) -- (0.6, 1) -- (1, 1)
	arc (-90: 90: 0.3cm and 0.1cm)
	-- (0.6, 1.2) -- (0.5, 1.1) -- (0.4, 1.2) 
	-- (0, 1.2) arc (90: 270: 0.3cm and 0.1cm);
\draw (1.3, 1.1) -- (1.3, 0.1) arc (0: -90: 0.3cm and 0.1cm)
	-- (0,0) arc (270: 180:  0.3cm and 0.1cm)
	-- (-0.2, 0.6) -- (-0.3, 1.1)
	(0.5, 1.1) to [out = -90, in = 0] (-0.2, 0.6)
	(-0.3, 0.1) arc (-180: 0: 0.8cm and 0.5cm);
\draw [densely dashed] (-0.3, 0.1) arc (180: 90: 0.3cm and 0.1cm) -- (1,0.2) arc (90: 0: 0.3cm and 0.1cm);
\end{scope}
\node at (1 ,-9.5) {$=$};
% Gluing-Cup Relation A2
\begin{scope}[xshift = 0cm, yshift = -10cm, xscale = -1]
\draw (0,1) -- (0.4, 1) -- (0.5, 1.1) -- (0.6, 1) -- (1, 1)
	arc (-90: 90: 0.3cm and 0.1cm)
	-- (0.6, 1.2) -- (0.5, 1.1) -- (0.4, 1.2) 
	-- (0, 1.2) arc (90: 270: 0.3cm and 0.1cm);
\draw (1.3, 1.1) -- (1.3, 0.1) arc (0: -90: 0.3cm and 0.1cm)
	-- (0,0) arc (270: 180:  0.3cm and 0.1cm)
	-- (-0.2, 0.6) -- (-0.3, 1.1)
	(0.5, 1.1) to [out = -90, in = 0] (-0.2, 0.6)
	(-0.3, 0.1) arc (-180: 0: 0.8cm and 0.5cm);
\draw [densely dashed] (-0.3, 0.1) arc (180: 90: 0.3cm and 0.1cm) -- (1,0.2) arc (90: 0: 0.3cm and 0.1cm);
\end{scope}

% Gluing-Saddle Relation A
\begin{scope}[xshift = 0.5cm, yshift = -8.5cm, xscale = 1]
\draw (0,0) -- (0,2) -- (0.4, 2) -- (0.5, 2.1) -- (0.6, 2) -- (2.4, 2)
	-- (2.4, 0) -- (2,0) arc (270: 180: 0.3cm and 0.1cm)
	-- (1.7, 1.1) arc (180: 270: 0.3cm and 0.1cm) -- (2.4, 1)
	(1.7, 1.1) arc (0: 180: 0.2cm and 0.3cm)
	-- (1.2, 0.6) -- (1.3, 0.1) arc (0: -90: 0.3cm and 0.1cm) -- (0,0);
\draw (0,0.2) -- (-0.1, 0.2) -- (-0.1, 2.2) 
	-- (0.4, 2.2) -- (0.5, 2.1) -- (0.6, 2.2)
	-- (2.5, 2.2) -- (2.5, 0.2) -- (2.4, 0.2);
\draw (-0.1, 1.2) -- (0, 1.2)
 	(0,1) -- (0.4, 1) -- (0.5, 1.1) -- (0.6, 1) -- (1,1) arc (-90: 0: 0.3cm and 0.1cm)
	(2.4, 1.2) -- (2.5, 1.2);
\draw [densely dashed] (0,0.2) -- (1, 0.2) arc (90: 0: 0.3cm and 0.1cm)
		(2.4, 0.2) -- (2, 0.2) arc (90: 180: 0.3cm and 0.1cm)
		(2.4, 1.2) -- (2, 1.2) arc (90: 180: 0.3cm and 0.1cm)
		(0,1.2) -- (0.4, 1.2) -- (0.5, 1.1) -- (0.6, 1.2) -- (1, 1.2) arc (90: 0: 0.3cm and 0.1cm);
\draw (0.5, 2.1) -- (0.5, 1.1) to [out = -90, in = 180] (1.2, 0.6);
\end{scope}
% Gluing-Saddle Relation A2
\begin{scope}[xshift = -1cm, yshift = -8.5cm, xscale = -1]
\draw (0,0) -- (0,2) -- (0.4, 2) -- (0.5, 2.1) -- (0.6, 2) -- (2.4, 2)
	-- (2.4, 0) -- (2,0) arc (270: 180: 0.3cm and 0.1cm)
	-- (1.7, 1.1) arc (180: 270: 0.3cm and 0.1cm) -- (2.4, 1)
	(1.7, 1.1) arc (0: 180: 0.2cm and 0.3cm)
	-- (1.2, 0.6) -- (1.3, 0.1) arc (0: -90: 0.3cm and 0.1cm) -- (0,0);
\draw (0,0.2) -- (-0.1, 0.2) -- (-0.1, 2.2) 
	-- (0.4, 2.2) -- (0.5, 2.1) -- (0.6, 2.2)
	-- (2.5, 2.2) -- (2.5, 0.2) -- (2.4, 0.2);
\draw (-0.1, 1.2) -- (0, 1.2)
 	(0,1) -- (0.4, 1) -- (0.5, 1.1) -- (0.6, 1) -- (1,1) arc (-90: 0: 0.3cm and 0.1cm)
	(2.4, 1.2) -- (2.5, 1.2);
\draw [densely dashed] (0,0.2) -- (1, 0.2) arc (90: 0: 0.3cm and 0.1cm)
		(2.4, 0.2) -- (2, 0.2) arc (90: 180: 0.3cm and 0.1cm)
		(2.4, 1.2) -- (2, 1.2) arc (90: 180: 0.3cm and 0.1cm)
		(0,1.2) -- (0.4, 1.2) -- (0.5, 1.1) -- (0.6, 1.2) -- (1, 1.2) arc (90: 0: 0.3cm and 0.1cm);
\draw (0.5, 2.1) -- (0.5, 1.1) to [out = -90, in = 180] (1.2, 0.6);
\end{scope}
\node at (-0.25, -7.5) {$=$};
\end{tikzpicture}
% \end{center}
\caption{Relations for $\bord_2$}
\label{MainUnOrientedTheormRelFig}	
\end{center}
\end{figure}
\end{theorem}

\subsection{Interpreting the statement of the theorem}

In Chapter~\ref{SymMonBicatChapt} we developed many aspects of the theory of symmetric monoidal bicategories. One of the key concepts introduced there was the notion of a \emph{symmetric monoidal computad}, or {\em presentation} of a symmetric monoidal bicategory, see Sections~\ref{sec:computads} and~\ref{SectFreelyGenSymMonBicats}. Such computads consist of a collection of generating objects, generating 1-morphisms, generating 2-morphisms, and relations placed on the 2-morphisms. This data has an inductive aspect. We may first form a symmetric monoidal bicategory which is freely built from the generating objects. Then the sources and targets of the generating 1-morphisms are taken from this; they are consequences of the generating objects and the operations in a symmetric monoidal bicategory. Then we build the next layer, the free symmetric monoidal bicategory built from the generating objects and 1-morphisms. The sources and targets of the generating 2-morphisms are taken from this, and so on. 

Given a symmetric monoidal computad $P$ there exists a symmetric monoidal bicategory $\sF(P)$ generated by $P$, and this symmetric monoidal bicategory enjoys a pleasant \emph{cofibrancy property}. Up to equivalence, symmetric monoidal homomorphisms out of $\sF(P)$ are determined by their values on the generators of $P$, see Theorem~\ref{thm:cofibrancythm}. 

The theory of computads is extremely general and applies to any finitary monad on globular sets. In Chapter~\ref{SymMonBicatChapt} we introduced two important auxiliary notions: \emph{unbiased semistrict symmetric monoidal 2-categories} and the stricter \emph{quasistrict symmetric monoidal 2-categories}, see Sections~\ref{sec:strictsymbicats} and~\ref{sec:PresSemistrict}. These are stricter versions of symmetric monoidal bicategories. Any symmetric monoidal computad $P$ also gives rise to a computad for the  unbiased semistrict symmetric monoidal 2-category monad and for the quasistrict symmetric monoidal 2-category monad. 

These induced computads each give rise to a symmetric monoidal bicategory of the appropriate sort, $\sF_{uss}(P)$ and $\sF_{qs}(P)$, respectively. Our coherence results, Theorems~\ref{thm:USSCoherence} and~\ref{thm:CoherenceThm}, imply that the canonical homomorphisms
\begin{equation*}
	\sF(P) \stackrel{\sim}{\longrightarrow} \sF_{uss}(P)\stackrel{\sim}{\longrightarrow} \sF_{qs}(P)
\end{equation*}
are equivalences of symmetric monoidal bicategories. 

The data depicted in Figures \ref{MainUnOrientedTheormGenFig} and \ref{MainUnOrientedTheormRelFig} determine a specific symmetric monoidal computad $P = (\cR, G_2, G_1,G_0)$.  In this case there is a single generating object, the `point'; there are two generating 1-morphisms, the right and left `elbows'  \tikz{\draw (-6.6, -.5) -- (-6.5, -.5) arc (90: -90: 0.3cm and 0.1cm); \draw (-5.5, -.5) arc (90: 270: 0.3cm and 0.1cm) -- (-5.4, -.7);}; there are ten generating 2-morphisms and 17 relations. These figures require some interpretation. 
For example the target of the left elbow \tikz{\draw (-5.5, -.5) arc (90: 270: 0.3cm and 0.1cm) -- (-5.4, -.7);} is the unit object (depicted as the empty set) while the source is the tensor product of two copies of the generating object (depicted as two points). Similarly the source of the generating 2-morphism (which we read from top to bottom)
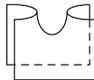
\begin{center}

\end{center}
is the horizontal composite of the generating 1-morphisms. Its target, depicted as two straight lines, is the identity 1-morphism on the tensor product of two copies of the generating object, another consequence of the operations in a symmetric monoidal bicategory. 

A similar example occurs for the {\em symmetry generators}. For example take the following generator: 
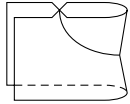
\begin{center}
	%
\end{center}
The target of this generator is given by the right elbow generating 1-morphism. The source, however, is the composite of the right elbow and the braiding morphisms $\beta$, which is an endomorphism of the tensor product of two copies of the generating object. This is not a generating 1-morphism itself, but is a consequence of the operations in a symmetric monoidal bicategory. 

Similar phenomena occur for the relations. For example the following relation,
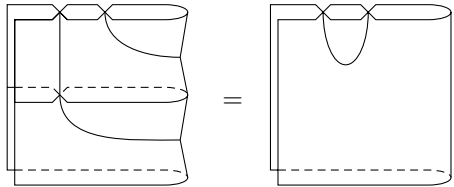
\begin{center}
	%
\end{center}
states an equality between two applications of the above symmetry generator (the left-hand-side) and a single application of the sylleptic 2-morphism $u$ (the right-hand-side). The latter is a consequence of the structure of a symmetric monoidal bicategory. 

 There are certain ambiguities in the data depicted in Figures \ref{MainUnOrientedTheormGenFig} and \ref{MainUnOrientedTheormRelFig}. For example we make no attempt to fix a bracketing of the sources or targets, we leave certain coherence morphisms (which are part of the data of a symmetric monoidal bicategory) implicit, and certain sources/targets of the generators or relations are left ambiguous (e.g. is the target the identity of the tensor product of two generating objects, or is it the tensor product of two identities on generating objects?).

We trust that the reader can make the necessary choices to make the data depicted in Figures \ref{MainUnOrientedTheormGenFig} and \ref{MainUnOrientedTheormRelFig} into a specific unambiguous symmetric monoidal computad. To do so here would complicate the notation so as to be incomprehensible. However, most importantly, there is no harm in leaving this ambiguity. Our coherence results imply that any two ways of resolving these ambiguities result in equivalent symmetric monoidal bicategories, provided that they induce the same underlying quasistrict symmetric monoidal computad. 

The symmetric monoidal computad $P$, determined by Figures \ref{MainUnOrientedTheormGenFig} and \ref{MainUnOrientedTheormRelFig}, gives rise to a symmetric monoidal bicategory $\sF(P)$. As we have already remarked, in light of our Cofibrancy Theorem~\ref{thm:cofibrancythm} (see also Prop.~\ref{RelationsForSMBicatsProp}), to give a symmetric monoidal homomorphism out of $\sF(P)$ one merely needs to provide the images of the generating morphisms, subject to the relations of $P$. This will be used to construct a symmetric monoidal homomorphism $\sF(P) \to \bord_2$. 

In Figure~\ref{MainUnOrientedTheormGenFig} we have essentially already drawn each generator in terms of its corresponding bordism. A little care must be taken with symmetry generators, for example  
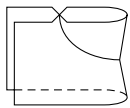
\begin{center}
	\begin{tikzpicture}
		\draw (.1,.2) -- (0, 0.2) --(0,1.3) -- (0.6, 1.3) -- (0.7, 1.2) -- (0.8, 1.3)
			 -- (1.3,1.3)  arc (90: -90:  0.3cm and 0.1cm) 
			 -- (0.8, 1.1) -- (0.7, 1.2) -- (0.6, 1.1) -- (0.1, 1.1) -- (0.1,0)
			 -- (1.3, 0) 	 arc (-90: 0:  0.3cm and 0.1cm)  to (1.5, 0.6)	 -- (1.6, 1.2);
		\draw [densely dashed] (0.1,0.2) 
			-- (1.3, 0.2) arc (90: 0:  0.3cm and 0.1cm);
		\draw (0.7, 1.2) to [out = -90, in = 180] (1.5, 0.6);
	\end{tikzpicture}
\end{center}
represents a particular invertible 2-bordism which realizes the diffeomorphism (rel. boundary) between the twisted right elbow and the untwisted right elbow. All of the relations are easily verified, a detail we leave to the reader for now, but which will also be discussed to some extent below. 

Thus we obtain the desired symmetric monoidal homomorphism:
\begin{equation*}
	E:\sF(P) \to \bord_2.
\end{equation*}
Theorem~\ref{TheoremUnorientedClassification} now asserts that this homomorphism is an equivalence of symmetric monoidal bicategories. 
 
\subsection{Overall strategy of the proof} \label{sec:strategyofproof}

The first method that one might attempt in order to prove Theorem~\ref{TheoremUnorientedClassification} is by a direct comparison. Whitehead's Theorem for symmetric monoidal bicategories (Theorem~\ref{WhiteheadforSymMonBicats}) provides a simple list of  criteria for  when such a symmetric monoidal  homomorphism is an equivalence. We must show three things:
\begin{enumerate}
	\item the homomorphism must be essentially surjective on objects,
	\item the homomorphism must be essentially full on 1-morphisms,
	\item the homomorphism must be fully-faithful on 2-morphisms. 
\end{enumerate}	
	The first of these is the statement that every (compact) 0-manifold is equivalent in $\bord_2$ to a finite disjoint union of points. This is clear. The second statement is the assertion that any 1-bordism is isomorphic in $\bord_2$ to one in the image of $E$, i.e., one obtained by gluing elementary 1-bordisms and permutations of 0-manifolds. The existence of Morse functions, and the corresponding Morse decompositions, ensures this is the case. However it is significantly more difficult to see that the functor is fully-faithful on 2-morphisms. The 2-morphisms of the bicategory $ \sF(P)$ consist of equivalence classes of binary paragraphs in the generating data (see Definition \ref{DefnFreeSymmetricMonBicat}) and this is quite difficult to compare directly to the bordism bicategory. 
	
Instead we will use the results of both of the previous two chapters to form an indirect comparison. We will introduce two additional symmetric monoidal bicategories $\bord_2^{PD}$ and $\sB^{PD}$. The symmetric monoidal bicategory $\bord_2^{PD}$ is similar to the bordism bicategory, except that the 1-morphisms are equipped with Morse functions and compatible {\em linear diagrams}. The 2-morphisms are similarly equipped with generic maps to the unit square $I^2$ and compatible {\em planar diagrams}. These are taken up to an equivalence relation, see Section~\ref{SectPlanarDecompThm}. We will review these structures below in more detail. The symmetric monoidal bicategory $\sB^{PD}$ is similar, except that we have forgotten the bordisms completely and only retained the data of the linear and planar diagrams. 

The main result of Chapter~\ref{ChapPlanarDecomp} is that the data of the planar diagram is sufficient to completely reconstruct the given bordism up to diffeomorphism (Theorem~\ref{PlanarDecompositionTheorem}). In the current context, when paired with Whitehead's Theorem for symmetric monoidal bicategories (Theorem~\ref{WhiteheadforSymMonBicats}),  this implies that the forgetful homomorphisms
\begin{equation*}
	\sB^{PD} \stackrel{\sim}{\longleftarrow} \bord^{PD}_2 \stackrel{\sim}{\longrightarrow} \bord_2
\end{equation*} 
are equivalences of symmetric monoidal bicategories. The homomorphism $E: \sF(P) \to \bord_2$ factors through $\bord^{PD}_2$.

The symmetric monoidal bicategory $\sB^{PD}$ has a particularly simple structure. Its 2-morphisms consist of equivalence classes of planar diagrams, which are essentially string diagrams labeled by certain kinds of data (singularity data and gluing data). These diagrams may be concatenated vertically and horizontally, and both operations are strictly associative. In fact, as we will see, $\sB^{PD}$ is an unbiased semistrict symmetric monoidal 2-category. Thus there is a unique strict symmetric monoidal homomorphism filling the dashed arrow and making the following diagram commute:
\begin{center}
\begin{tikzpicture}
		\node (LT) at (0, 1.5) {$\sF(P)$};
		\node (LB) at (0, 0) {$\sF_{uss}(P)$};
		\node (RRT) at (6, 1.5) {$\bord_2$};
		\node (RT) at (3, 1.5) {$\bord_2^{PD}$};
		\node (RB) at (3, 0) {$\sB^{PD}$};
		\draw [->] (LT) -- node [left] {$\simeq$} (LB);
		\draw [->] (LT) -- node [above] {$$} (RT);
		\draw [->] (RT) -- node [above] {$\simeq$} (RRT);
		\draw [->] (RT) -- node [right] {$\simeq$} (RB);
		\draw [->, dashed] (LB) -- node [below] {$\exists !$} (RB);
		%\node at (0.5, 1) {$\ulcorner$};
		%\node at (1.5, 0.5) {$\lrcorner$};
\end{tikzpicture}
\end{center}
The presentation $P$ has been chosen so that the dashed arrow above is actually an \emph{isomorphism} of unbiased semistrict symmetric monoidal 2-categories.  Theorem~\ref{TheoremUnorientedClassification} then follows immediately.

\subsection{Review of Planar Diagrams} To construct the symmetric monoidal bicategories $\bord_2^{PD}$ and $\sB^{PD}$ and show they are equivalent to $\bord_2$ we will need to use the results from Chapter~\ref{ChapPlanarDecomp}, especially the material from Section~\ref{SectPlanarDecompThm}. We will give a brief summary of the key results, but we encourage the reader to go back and review that material for full details. 

An important structure that is introduced in that Chapter is the notion of a \emph{planar diagram} (Def.~\ref{def:planardiagram}). A planar diagram $(\Psi, \Gamma, \cS)$ has several parts. First, a planar diagram has a 2-dimensional \emph{graphic} $\Psi$ (Def.~\ref{2dGraphicDefn}) which is a diagram in the plane consisting of a collection isolated points and embedded arcs. The endpoints of the arcs always lie on the isolated points, and they are otherwise in general position. For each arc the projection to the vertical $y$-axis is a local diffeomorphism.  Locally around the isolated points the structure falls into one of three kinds: cusp points, 2D Morse points, or crossing points, see Section~\ref{SectGeomofSing}.  

\begin{figure}[htbp]
	\begin{center}
	\begin{tikzpicture}[align=center]
		\draw (0,0) rectangle (7,4);
		% immagine a box from (0,0) to (8,4)
	%	\node (A) at (-1, 3) {not \\ allowed};
	%	\node (B) at (7, 2) {not \\ allowed};
	%	\draw [->] (A) -- ++(1.25,0);
	%	\draw [->] (B) -- ++(-1.25,0);

		\draw [red] (0.5, 0) parabola bend (1,2) (1.5,0);
		\draw [red] (1,4) to [out = -90, in = 90] (2,2) parabola bend (2.5, 3) (4,0);
		\draw [red] (3,0) to [out = 90, in = -90] (4.5,4);
		
		\draw [chambering]  (0.25,3) to [out = -90, in = 135] (0.75,2) to [out = -90, in = 135] (1,1); \draw [chambering] (1,1) -- (1, 0.5); \fill [chambering] (1, 0.5) circle (1pt) (0.25,3) circle (1pt);
		
		\draw [chambering] (5,4) to [out = -90, in = 135] (5.5,2) (5.5,0) -- (5.5,2) (6,4) to [out = -90, in = 45] (5.5,2); %\fill (5.5,2) circle (1pt);
		
		\draw [chambering] (0.5, 4) to [out = -90, in =135] (1.5, 2.5) to [out = -135, in = 45] (0.75,2)
			(1.5, 2.5) to [out = -45, in = 90] (2, 1.5) to [out = -135, in = 45] (1,1)
			(2, 1.5) to [out = -45, in = 90] (2.5, 1) to [out = -135, in = 135] (2.5, 0.5)
			(2.5, 1) to [out = -45, in = 45] (2.5, 0.5) -- (2.5, 0);
		\draw [chambering] (3,4) to [out = -90, in = 90] (4.5, 1); 
		\fill [chambering] (4.5,1) circle (1pt);

	\end{tikzpicture}
	\end{center}
	\caption{A planar diagram (without the labels).}
	\label{fig:planardiagramdemo}
\end{figure}
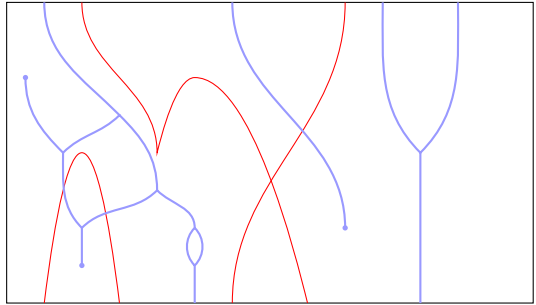

In addition to the graphic, a planar diagram includes a \emph{chambering graph} $\Gamma$ (Def.~\ref{def:chamberinggraph}). A chambering graph is a very simple sort of stratified space which is embedded into the plane and is in general position with respect to the graphic. The chambering graph is only allowed certain kinds of vertices (i.e. only certain singularity types are permitted). It is only allowed to have trivalent and univalent vertices. Moreover the projection from each edge of the chambering graph to the vertical $y$-axis is a local diffeomorphism. The trivalent vertices must obey the rule `2-up-1-down' or `1-up-2-down'. A prototypical example of these structures is shown in Figure~\ref{fig:planardiagramdemo}, with the graphic in \textcolor{red}{red} and the chambering graph in \textcolor{blue}{blue}.

Finally a planar diagram also consists of a collection of labels, called \emph{sheet data} $\cS$, which live on these structures (see Def.~\ref{def:planardiagram}).  The chambering graph and the graphic divide the plane into connected regions, which we call {\em chambers}. Each of the chambers receives a label $(n)$ from $\N_{\geq 0}$, which we should think of as the set $\{1,2, \dots, n\}$. 

Each edge of the chambering graph is assigned a label in the symmetric group $\Sigma_n$ where $n$ is the label of the adjacent chambers. Note that this implies that the two chambers on either side of an edge of the chambering graph must have the same labeling set. For the purposes of this labeling, the edges of the chambering graph are regarded as being split along arcs of the graphic which cross the graph. At each vertex of the chambering graph we impose a `cocycle' condition on these symmetric group labelings. This condition states that there is no monodromy when passing about a vertex of the chambering graph.  

Finally the graphic also receives labels, as do the intersection points of the graphic and chambering graph. Roughly, each arc and isolated point of the graphic is given a label which corresponds to one of the singularity types of a generic map from a surface to the plane. These are then further equipped with isomorphisms of sets based on the structure  of the local model of these singularities. The precise labels are given in Table~\ref{PlanarSheetDataAndGluingDataTable}. 

A generic map from a surface to the plane induces a graphic (Theorem~\ref{ThmGraphicOfMaptoR2}). Planar diagrams are designed so as to encode enough information to reconstruct the surface. Given a planar diagram $(\Psi, \Gamma, \cS)$ there exists a canonical surface $\Sigma$ equipped with a generic map to the plane reproducing the graphic $\Psi$. This surface is obtained by gluing together the local models for the surface according to the sheet data of the planar diagram. Given a surface $\Sigma$ with a generic map to the plane, there exists an enhancement of its graphic to a planar diagram such that the surface constructed from the planar diagram is diffeomorphic over $\R^2$ to $\Sigma$ (Theorem~\ref{PlanarDecompositionTheorem}).

We also have precise control over when two planar diagrams give rise to diffeomorphic surfaces. Of course isotopic diagrams give rise to diffeomorphic surfaces, but in addition there is a finite list of combinatorial local moves that may be preformed on planar diagrams. Two planar daigrams give rise to diffeomorphic surfaces if and only if they may be related by a finite sequence of these local moves (Theorem~\ref{PlanarDecompMovesThm}). Though comprised only of expected maneuvers, list of possible local moves is quite lengthly. They are depicted in Tables~\ref{LocalMovesForPlanarDiagramsTable1} and~\ref{LocalMovesForPlanarDiagramsTable2} (also consult Table~\ref{3DSheetDataTable}), and in Figures~\ref{fig:0stratumChamberMoves},~\ref{fig:localmoveschamberinggraph}, and~\ref{fig:chambergraphicmoviemoves}. Theorem~\ref{PlanarDecompMovesThm} gives a full account. 

In the case of the bordism bicategory, our manifolds are manifolds with corners. Specifically they are $\langle 2\rangle$-manifolds such that each of the faces has been divided into two components (incoming/ourgoing). We will need a small variation on the results of Section~\ref{SectPlanarDecompThm} which is suitable for this situation. In particular we will insist that our planar diagrams have a very particular behavior near the boundary of the square. The behavior of our planar diagrams near the boundary is coupled to the behavior of generic maps to the square. We require these also to satisfy particular boundary conditions. 

Recall that a 2-bordism is a $d$-dimensional $\langle 2\rangle$-manifold $S$ with decompositions 
\begin{align*}
	\partial_0 S &= \partial_{0, \text{in}} S \sqcup \partial_{0, \text{out}}S \stackrel{g}{\to} W_0 \sqcup W_1 \\
	\partial_1 S &= \partial_{1, \text{in}} S \sqcup \partial_{1, \text{out}}S \stackrel{f}{\to} Y_0 \times I \sqcup Y_1 \times I
\end{align*}
of its boundary faces. The square $I^2$ is also a $\langle 2\rangle$-manifold with decompositions 
\begin{align*}
	\partial_0 I^2 = I \times \partial I = I \times \{1\} \sqcup I \times \{0\} \\
	\partial_1 I^2 = \partial I \times  I =  \{1\} \times I   \sqcup \{0\} \times I.
\end{align*}
We will consider maps from $S$ to $I^2$ which respect the structure as $\langle 2\rangle$-manifolds, as well as these compositions. Specifically $\partial_{0, \text{in}} S $ maps to $I \times \{1\}$, $\partial_{0, \text{out}}S$ maps to $ I \times \{0\}$, $\partial_{1, \text{in}} S$ maps to $\{0\} \times I$, and $\partial_{1, \text{out}}S$ maps to $\{1\} \times I$.
Furthermore we require that the map $\partial_1 S \to \partial I \times I$ contains no singularities, that the map $\partial_0 S \to I \times \partial I$ contains no codimension two singularties or multijet singularities (see Section~\ref{SectStratofJetSpace2D}), and finally that the singular locus of the codimension one singularities is transverse to $I \times \partial I$.  

In terms of the graphic, which now lives in $I^2$, this means that the isolated points occur on the interior, that the arcs are disjoint from the sides of the square ($\partial I \times I$,) and that they intersect the top and bottom of the square ($I \times \partial I$) transversely. Thanks to the relative transversality theorem (Theorem~\ref{ThomTransversalityTheoremRelativeCorners}) our careful analysis from Chapter~\ref{ChapPlanarDecomp} carries over directly to this setting. The modified notions of \emph{linear} and \emph{planar diagrams} suitable for 2-bordisms are as follows:

\begin{definition}%\label{def:}
	A \emph{linear diagram in $I$} consists of a triple $(\Psi, \Gamma, \cS)$ where $\Psi$ and $\Gamma$ are disjoint finite subsets of $I$, the \emph{graphic} and \emph{chambering set}, respectively, which are both disjoint from the boundary. These divide $I$ into a finite number of \emph{chambers}, the regions between points of $\Psi \cup \Gamma$. These structures are equipped with a labeling $\cS$ such that each chamber is labeled with a natural number $(n) \in \N_{\geq 0}$, which we think of as the set with $n$ elements. Each point in $\Gamma$ is labeled with a permutation identifying the sets on its two adjacent chambers (which must consequently be labeled by the same element of $\N_{\geq 0}$). Finally the points of the graphic $\Psi$ must be labeled by triples $(m,x,n)$ where $m,n \in \N_{\geq 0}$ and $x$ is a symbol ``\tikz{\draw (0,0) -- +(-3pt, 0) arc (270:90: 3pt) -- +(3pt, 0); }'' or ``\tikz{\draw (0,0) -- +(3pt, 0) arc (-90:90: 3pt) -- +(-3pt, 0); }'' indicating one of two possible critical point indices. 
	
	The source of $(m, \tikz{\draw (0,0) -- +(-3pt, 0) arc (270:90: 3pt) -- +(3pt, 0); }, n)$, that is the label of the chamber immediately to the right of the arc of the graphic\footnote{Planar diagrams and bordisms are read right to left.}, is the concatenation $(m+2+n)$, while the target is $(m+n)$. The source of $(m,\tikz{\draw (0,0) -- +(3pt, 0) arc (-90:90: 3pt) -- +(-3pt, 0); }, n)$ is $(m+n)$, while the target is $(m+2+n)$. 
\end{definition}

These linear diagrams will be the restriction of our planar diagrams to the top and bottom of the square $I^2$. The graphic of the linear diagram corresponds to the intersection of the fold arcs with these top and bottom edges. Notice, however, that the data which we assign to the graphic of a linear diagram is slightly different from the fold sheet data of Section~\ref{sec:planar:reconstruction}, which consists of a pair of sets $(S,T)$ where the cardinality of $T$ is two. 

In our current setting we have chosen to trivialize the set $S \sqcup T$ by identifying it with an element in $\N_{\geq 0}$. Moreover when we make this identification we require that the two element set $T$ be identified with a pair of \emph{consecutive} elements. This is no real loss of generality, but it is quite convenient later on. This convention matches better with the categorical structures and theory developed in Chapter~\ref{SymMonBicatChapt}.

Similarly for our planar diagram we will have to make similar trivializations of the sheet data. For our planar diagram, each region will be equipped with a label from $\N_{\geq 0}$. The edges of the chambering graph will be labeled with permuations of these sets, and as before around each vertex of the chambering graph we require there to be trivial monodromy. For the edges of a 2-dimensional graphic we label them by triples $(m, x, n)$, exactly as in the definition of linear diagrams above. The same holds for the 2D Morse points, as they are essentially fold singularities as well. For the intersection of a arc of the graphic and edge of the chambering graph, we label it exactly as before, that by a pair of permutations $(f,g): (m+n, 2) \to (m' + n', 2)$ (note that $m$ and $m'$, respectively $n$ and $n'$, may be different).  

For the remaining isolated points of the graphic  (the fold crossing and the cusp points) we have a similar convention. The case of fold crossing points is a natural extension of the fold sheet data. In the language of Section~\ref{sec:planar:reconstruction}, the sets $T$ and $T'$ are each identified with $2$, while the set $S$ is identified with $m + n +p$ for some choice of $m, n , p \in \N_{\geq 0}$. The insertion of $T$ and $T'$ then happen between $m$ and $n$ and between $n$ and $p$, respectively. 

For the cusp points it is easiest to describe as concatenations of minimal sheet data. There is a certain minimal amount of sheet data which supports the standard model of the cusp singularity. Specifically this is the case that the set $S$ in the definition of cusp sheet data of Section~\ref{sec:planar:reconstruction} has size one, $S = (1)$. Moreover in this case we may choose an identification $S' \cong (3)$. For a general cusp point we additionally have $m, n \in \N_{\geq 0}$, and the relevant regions surrounding the cusp are labeled with $(m + 1 + n)$ and $(m + 3 + n)$, respectively.

\begin{definition}%\label{def:}
	A \emph{planar diagram in $I^2$} consists of a triple $(\Psi, \Gamma, \cS)$ where $\Psi$ is a graphic in $I^2$, $\Gamma$ is a chambering graph in $I^2$, and $\cS$ is a collection of sheet data with respect to $(\Psi, \Gamma, \cS)$. Further we require:
	\begin{itemize}
		\item The isolated points of $\Psi$ occur on the interior of $I^2$. The
		 arcs of the graphic $\Psi$ are disjoint from the vertical sides of the square ($\partial I \times I$), and transverse to the horizontal side ($I \times \partial I$);
		\item The vertices of $\Gamma$ occur on the interior of $I^2$. The edges of $\Gamma$ are disjoint from the vertical sides of the square ($\partial I \times I$), and transverse to the horizontal side ($I \times \partial I$);
		\item The sheet data $\cS$ is given exactly as described above. 
	\end{itemize}	
\end{definition}

Thus the restriction of a planar diagram in $I^2$ to its top and bottom edges is a linear diagram in $I$. 

Any $I$-valued Morse function for a 1-bordism induces a graphic which may be enhanced to a linear diagram. Conversely, given a linear diagram $L$ with right-most region labeled by $s_0$ and left-most region labeled by $s_1$, we can reconstruct a 1-bordism $Y_L$ from $s_0$ to $s_1$ together with a Morse function to $I$ inducing the graphic of the linear diagram. We will say this 1-bordism $Y_L$ is {\em constructed} from the linear diagram $L$. Similarly the graphic of a generic map from a 2-bordism to $I^2$ may be enhanced to a planar diagram $P$ and from a such a planar diagram we get a surface $\Sigma_P$ with a generic map to $I^2$ inducing the same graphic. See Section~\ref{sec:planarandspacialdiagrams}.

The results of Section~\ref{SectPlanarDecompThm} now apply and yield the following analog of Theorem~\ref{PlanarDecompositionTheorem}:

\begin{theorem}[Planar decomposition theorem for bordisms]\label{thm:bordismplanardecomp}
	Let $Y_\textrm{in}$ and $Y_\textrm{out}$ be 1-bordisms, i.e. objects of $(\COB_2)_1$. Let $f_s: S(Y_\textrm{in}) \to S(Y_\textrm{out})$ and $f_t: T(Y_\textrm{in}) \to T(Y_\textrm{out})$ be diffeomorphisms, and let $\Sigma$ be 2-bordism  representing a morphism in $(\COB_2)_1$ from $Y_\textrm{in}$ to $Y_\textrm{out}$, inducing $f_s$ and $f_t$. Then:
	\begin{enumerate}
		\item Suppose that we are given linear diagrams $L_\textrm{in}$ and $L_\textrm{out}$ with source and target given by $S(Y_\textrm{in})$ and $T(Y_\textrm{in})$, respectively $S(Y_\textrm{out})$ and $T(Y_\textrm{out})$. Suppose further that we have diffeomorphisms of 1-bordisms $Y_\textrm{in} \cong Y_{L_\textrm{in}}$ and $Y_\textrm{out} \cong Y_{L_\textrm{out}}$. Then, using the isomorphisms $f_s$ and $f_t$ we may alter $L_\textrm{out}$ at its boundary regions so that they are also labeled by $S(Y_\textrm{in})$ and $T(Y_\textrm{in})$. Thus we have a valid restriction of a planar diagram to $\partial (I^2)$, which is given by  $L_\textrm{in}$ on top,  the modified $L_\textrm{out}$ on bottom, and constant diagrams on the sides. 

		Then there exists a planar diagram $P$ extending the given restriction to  $\partial (I^2)$ and an isomorphism of 2-bordisms $\Sigma \cong \Sigma_{P}$ extending the boundary identification of $\partial \Sigma \cong \partial \Sigma_P$. 
	\item If $\Sigma'$ is a second 2-bordism of from $Y_\textrm{in}$ to $Y_\textrm{out}$, also inducing $f_s$ and $f_t$, and $P'$ is any planar diagram as above, then $\Sigma \cong \Sigma'$ relative boundary (hence they give the same morphism of $(\COB_2)_1$) if and only if $P$ is equivalent to $P'$ via a finite application of the following local moves
	\begin{itemize}
		\item Isotopy, relative to the boundary $\partial(I^2)$.
		\item Any of the local moves listed in Theorem~\ref{PlanarDecompMovesThm}(2) and (3), but applied solely to the interior of the planar diagram. \qed
	\end{itemize} 
	\end{enumerate}
\end{theorem}

\subsection{Auxiliary bordism bicategories, with planar diagrams} Now we will introduce the symmetric monoidal bicategories $\bord_2^{PD}$ and $\sB^{PD}$. Again we will do this using symmetric monoidal pseudo double categories (see Section~\ref{sec:symmonpseudo2cat}), just as we did for $\cob_d$ (Section~\ref{sec:symmonbicatcobord}) and $\bord_d$ (Section~\ref{sec:haloedBordd}). We will begin with $\bord_2^{PD}$. 

The symmetric monoidal bicategory $\bord_2^{PD}$ is similar to the symmetric monoidal bicategory $\bord_2$, except that we equip the 1-morphisms with compatible linear diagrams in $I$, and we equip the 2-morphisms with compatible planar diagrams in $I^2$. We will construct it as a symmetric monoidal pseudo double category $\BORD_2^{PD}$. 

The symmetric monoidal category of objects is almost the same as for $\BORD_2$. The objects of $(\BORD^{PD}_2)_0$ consist of triples $(Y, (n), f)$ where $Y = (Y, \hat{Y}_1, \hat{Y}_2) \in (\BORD_2)_0$ is a haloed zero manifold, $(n) \in \N_{\geq 0}$, and $f: Y \cong (n)$ is a bijection. A morphism of $(\BORD^{PD}_2)_0$ from $(Y, (n), f)$ to $(Y', (n), f')$ consists of a pair $(h, \sigma)$ where $h$ is an isomorphism in $(\BORD_2)_0$ and $\sigma \in \Sigma_n$ is a permutation, such that $\sigma \circ f = f' \circ h$.  The symmetric monoidal structure is induced by the disjoint union of haloed zero manifolds:
\begin{equation*}
	(Y, (n), f) \otimes (Y', (n'), f') := (Y \sqcup Y', (n + n'), f \sqcup f').
\end{equation*} 

The symmetric monoidal category of arrows $(\BORD^{PD}_2)_1$ has objects $(\Sigma, L, f)$ which consist of a haloed 1-bordisms $\Sigma = (\Sigma, \hat{\Sigma}_0, \hat{\Sigma}_1) \in (\BORD_2)_1$, a linear diagram $L$, and an isomorphism $f: \Sigma \cong \Sigma_L$, where $\Sigma_L$ is the 1-bordism constructed from $L$. The morphisms of $(\BORD^{PD}_2)_1$ from $(\Sigma^\textrm{in}, L^\textrm{in}, f^\textrm{in})$ to $(\Sigma^\textrm{out}, L^\textrm{out}, f^\textrm{out})$ are equivalence classes of triples $(S, P, f)$ where $S = (S, \hat{S}, h_\textrm{in}, h_\textrm{out})$ is a representative of a morphism from $\Sigma^\textrm{in}$ to $\Sigma^\textrm{out}$ in $(\BORD_2)_1$ (see Definition~\ref{def:haloedbordisms12}), $P$ is a planar diagram in $I^2$, and $f: S \cong S_P$ is an isomorphism, where $S_P$ is the 2-bordisms constructed from $P$. This data is required to satisfy the following condition: Let $\partial_0^\textrm{in} S_P$ denote the boundary face of $S_P$ mapping to the top most edge of $I^2$. Then we have two isomorphisms given as composites: 
	\begin{align*}
		\Sigma^\textrm{in} \cong& \Sigma_{L^\textrm{in}} \cong \partial_0^\textrm{in} S_P \\
		\Sigma^\textrm{in} \cong& \partial_0^{\textrm{in}} S \cong \partial_0^\textrm{in} S_P
	\end{align*}
	We require that these agree. Similarly, for the face component mapping to the bottom edge of $I_2$,
	 we require that the two isomorphisms $\Sigma^\textrm{out} \cong \partial_0^\textrm{out} S_P $ also agree. Two such triples $(S, P, f)$ and $(S', P', f')$ are equivalent precisely when $S \simeq S'$ are equivalent, that is they represent the same morphisms of $(\BORD_2)_1$. Notice that in light of Theorem~\ref{thm:bordismplanardecomp}, this is the case precisely when $P$ and $P'$ are related by the moves listed in that theorem. 

The composition of morphisms in $(\BORD^{PD}_2)_1$ is accomplished by composing morphisms as in $(\BORD_2)_1$ and concatenating planar diagrams vertically. At the level of planar diagrams this vertical concatenation is not completely defined. First it depends on a choice of isomorphisms $I \cup_{pt} I \cong I$ of concatenated intervals, and even then not all planar diagrams are composable as the result of concatenating the arcs of the graphic and edges of the chambering graph may no longer be smooth as the point of gluing. However morphisms in $(\BORD^{PD}_2)_1$ only depend on the isotopy class of the planar diagram, and the above operations \emph{are} well defined on isotopy classes of planar diagrams. We can always isotope our planar diagrams so that they are composable and the resulting composite is then unique up to isotopy. Hence this gives a well defined operation on morphisms in $(\BORD^{PD}_2)_1$. 

The symmetric monoidal structure on $(\BORD^{PD}_2)_1$ is more complicated, and is related to the construction of the functor $\odot$, which we define first.  On the parts coming from  $(\BORD^{}_2)_1$, both the monoidal structure $\otimes$ and the functor $\odot$ coincide with their counterparts in $(\BORD^{}_2)_1$. On planar and linear diagrams, the functor $\odot$ is given by horizontal concatenation. This is completely analogous to the (vertical) composition of morphisms in $(\BORD^{PD}_2)_1$. In contrast to that operation it is not strictly associative, but is up to canonical isotopy class of isotopy. This isotopy is used to given the transformations $\alpha$, $\rho$, and $\lambda$ which form part of the definition of symmetric monoidal pseudo double category. 

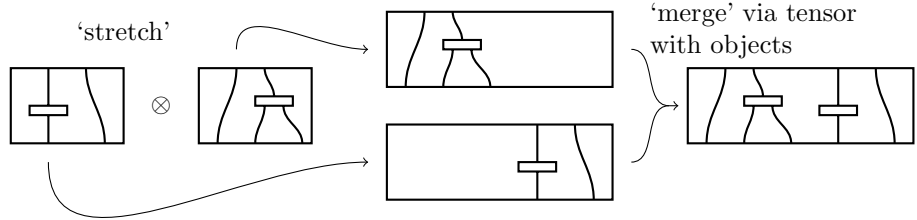
\begin{figure}[ht]
\begin{center}
% associative-like relation	
\begin{tikzpicture}[thick,scale = 0.5]
	\draw (0,0) rectangle (3,2);
	\draw (1,2) -- (1,1) 
		(0.5,0.75) rectangle (1.5,1) 
		(1,0.75) -- (1,0)
		(2,2) to [out=270, in=90] (2.5,0);
		
	\node at (4,1) {$\otimes$};

	\draw (5,0) rectangle (8,2);
	\draw (6,2) to [out=270, in=90] (5.5,0) 
		(6.75,2) to [out=270, in=90] (7,1.25)
		(6.5,1) rectangle (7.5,1.25)
		(6.75,1) to [out=270, in=90] (6.5,0)
		(7.25,1) to [out=270, in=90] (7.75,0);

	\draw (10,1.5) rectangle (16,3.5);
		\draw (11,3.5) to [out=270, in=90] (10.5,1.5) 
			(11.75,3.5) to [out=270, in=90] (12,2.75)
			(11.5,2.5) rectangle (12.5,2.75)
			(11.75,2.5) to [out=270, in=90] (11.5,1.5)
			(12.25,2.5) to [out=270, in=90] (12.75,1.5);
	\draw (10,-1.5) rectangle (16,0.5);
		\draw (14,0.5) -- (14,-0.5) 
			(13.5,-0.75) rectangle (14.5,-0.5) 
			(14,-0.75) -- (14,-1.5)
			(15,0.5) to [out=270, in=90] (15.5,-1.5);
	
	\draw (18,0) rectangle (24,2);
		\draw (19,2) to [out=270, in=90] (18.5,0) 
			(19.75,2) to [out=270, in=90] (20,1.25)
			(19.5,1) rectangle (20.5,1.25)
			(19.75,1) to [out=270, in=90] (19.5,0)
			(20.25,1) to [out=270, in=90] (20.75,0);
		\draw (22,2) -- (22,1) 
				(21.5,0.75) rectangle (22.5,1) 
				(22,0.75) -- (22,0)
				(23,2) to [out=270, in=90] (23.5,0);
	\draw [thin, ->] (1,-0.5) to [out=-90, in = 180] (9.5,-0.5);
	\draw [thin, ->] (6,2.5) to [out=90, in = 180] (9.5,2.5);
	\draw [thin] (16.5,2.5) to [out=0, in = 180] (17.5,1);
	\draw [thin] (16.5,-0.5) to [out=0, in = 180] (17.5,1);
	\draw [thin, ->] (17.5,1) -- (17.75,1);
	\node at (3,3) {`stretch'};
	\node [text width=3cm] at (20,3) {`merge' via tensor with objects};
\end{tikzpicture}
\caption{The monoidal operation for Planar Diagrams}
\label{Fig:MonoidalStructureBordPD}
\end{center}
\end{figure}

The symmetric monoidal structure of $(\BORD^{PD}_2)_1$ is constructed from $\odot$ in exactly the same way it is done for unbiased semistrict symmetric monoidal 2-categories, as in Section~\ref{sec:graphcalculusUSS}. First we note that given a planar diagram $P$ and an label $(n) \in \N_{\geq 0}$, we may form new planar diagrams $(n)\otimes P$ and $P \otimes (n)$ by concatenating every label of of the planar diagram $P$ with $(n)$, either from the left or the right. In terms of the bordisms constructed from these planar diagrams this operation corresponds to taking the disjoint union with the identity bordism on $(n)$. Then we \emph{define} the tensor product of two planar diagrams, $P_1$ and $P_2$, by the interchange law with $\odot$.   
\begin{equation*}
	 P_1 \otimes P_2 = (id_{x} \odot P_1) \otimes (P_2 \odot id_{y}) := (id_{x} \otimes P_2) \odot (P_1 \otimes id_{y})
\end{equation*}
where $x$ is the label of the left most chamber of $P_1$ and $y$ is the label of the right most chamber of $P_2$. This operation can be viewed as the process of merging the planar diagrams $P_1$ and $P_2$, as depicted in Figure~\ref{Fig:MonoidalStructureBordPD}. The remaining functors $S,T$ and $U$ are defined just as for $(\BORD_2)_1$, with $U$ producing the product planar diagram. 

The proof of Theorem~\ref{thm:symmonbicatcob} carries over directly to show that the symmetric monoidal pseudo double category $\BORD_2^{PD}$ is fibrant. Thus the horizontal bicategory is symmetric monoidal and this defines $\bord_2^{PD} = \cH(\BORD_2^{PD})$. 

The symmetric monoidal pseudo double category $\B^{PD}$ is defined similarly to $\BORD_2^{PD}$, except that we forget the bordisms entirely and only remember the linear and planar diagrams. In detail, $(\B^{PD})_0$ is the symmetric monoidal category whose objects are $(n) \in \N_{\geq 0}$, and whose morphisms are given by permutations: $(\B^{PD})_0( n, m) = \emptyset$ if $m \neq n$, and $\Sigma_n$ if $m=n$. The symmetric monoidal structure is strict and is given by the sum operation on $\N$. This is a skeleton of the category of finite sets and bijections. 

The symmetric monoidal category of arrows $\B^{PD}$ has objects given by \emph{isotopy classes} of linear diagrams. The morphisms from $[L]$ to $[L']$ are given by equivalence classes of planar diagrams from $L$ to $L'$, where $L$ and $L'$ are representatives of the isotopy class of linear diagrams and the equivalence relation is exactly that stated in Theorem~\ref{thm:bordismplanardecomp}. Notice that if we had chosen different representatives of the isotopy classes of linear diagrams, the resulting sets of morphisms are in canonical bijection. This is essentially due to the fact that there is a unique isotopy class of isotopies between any two representatives of $[L]$. Consequently we can define the composition of morphisms, the symmetric monoidal structure, and all the functors $S$, $T$, $\odot$, $U$, just as we did for $\BORD_2^{PD}$. 
The symmetric monoidal pseudo double category $\B^{PD}$ is fibrant and yields the symmetric monoidal bicategory $\sB^{PD} = \cH(\B^{PD})$. 

However there is one important difference in this case. Because we are taking isotopy classes of linear diagrams, and we have discarded the bordisms, the horizontal composition (given by $\odot$) is \emph{strictly} associative. The bicategory $\sB^{PD}$ us actually a strict 2-category. Moreover, because we have defined the symmetric monoidal structure using one of the interchange identities, we automatically have that $\sB^{PD}$ is a Gray monoid (Def.~\ref{def:graymonoid}). Even more is true, as we shall see shortly.    

\begin{lemma}%\label{lem:}
	The forgetful functors $\bord^{PD}_2 \to \bord_2$ and $\bord_2^{PD} \to \sB^{PD}$ are equivalences of symmetric monoidal bicategories. 
\end{lemma}

\begin{proof}
	 Whitehead's Theorem for Symmetric Monoidal Bicategories (Theorem~\ref{WhiteheadforSymMonBicats}) states that a symmetric monoidal homomorphism is an equivalence if and only if it is essentially surjective on objects, essentially full on 1-morphisms, and fully-faithful on 2-morphisms. These forgetful functors are in fact surjective on both objects and 1-morphisms (hence full on 1-morphisms). Moreover they are fully-faithful on 2-morphism by Theorem~\ref{thm:bordismplanardecomp}. 
\end{proof}

%------ older material ------

	%Definition \ref{DefnFreeSymmetricMonBicat} gives a precise meaning to the abstract symmetric monoidal  bicategory $\sF/ \cR$ with the above generators and relations and 
	
%	We have drawn the abstract generators of this bicategory as particular bordisms, and by the above considerations this defines a symmetric monoidal homomorphism $E: \sF / \cR \to \bord$.  The objects, 1-morphisms and 2-morphisms of $\sF / \cR$ are given by binary words, binary sentences, and paragraphs in the abstract generating data (modulo the relations) and the homomorphism $E$ is given by evaluating the given word, sentence, or paragraph in $\bord$. We must show that this is an equivalence of symmetric monoidal bicategories. 
	
\subsection{The Proof}

Theorem~\ref{TheoremUnorientedClassification} is now reduced to the following proposition, whose proof occupies the remainder of this section. 

\begin{proposition}%\label{pro:}
	The symmetric monoidal bicategory $\sB^{PD}$ is a computadic unbiased semistrict symmetric monoidal 2-category, generated by the presentation depicted in Figures \ref{MainUnOrientedTheormGenFig} and \ref{MainUnOrientedTheormRelFig}.
\end{proposition} 

Once this result is established the proof of Theorem~\ref{TheoremUnorientedClassification} precedes exactly as outlined in Section~\ref{sec:strategyofproof}.

Every symmetric monoidal computad, such as the one given in Figures \ref{MainUnOrientedTheormGenFig} and \ref{MainUnOrientedTheormRelFig}, gives rise to an unbiased semistrict symmetric monoidal 2-category. This should be thought of as the free  unbiased semistrict symmetric monoidal 2-category generated by the given symmetric monoidal computad (presentation), and symmetric monoidal 2-categories arising this way are called \emph{computadic}.
 
In Section~\ref{sec:graphcalculusUSS} we introduced a graphical calculus for unbiased semistrict symmetric monoidal 2-categories, and in Section~\ref{sec:presentationUSS} we showed how computadic unbiased semistrict symmetric monoidal 2-categories can be described explicitly in terms of this graphical calculus. We suggest that the reader remind themselves of this graphical calculus, but we will describe it briefly here.  

An unbiased semistrict symmetric monoidal computad $\cP$ consists of a collection of generating objects, a collection of generating 1-morphisms, a collection of generating 2-morphisms, and a collection of relations between these 2-morphisms.
In the corresponding computadic unbiased semistrict symmetric monoidal 2-category $\sF(\cP)$, the objects are given by words in the generating objects. The morphisms are given by composable sequences elementary 1-morphisms which are of two types: permutations of the words and morphisms of the form $(id_x \otimes f \otimes id_y)$ for some objects $x$ and $y$ and some generating 1-morphism $f$. Finally there 2-morphisms are given by equivalences classes of certain string diagrams, exactly the string diagrams which appear in the graphical calculus. 

These string diagrams consist of arcs and coupons, as well as the regions between these. These regions are labeled by the objects (i.e. words in the generating objects). The arcs come in two varieties, corresponding to the two kinds of elementary 1-morphisms. The coupons also come in several types: there are trivalent and univalent coupons among the permutation arcs, there are 4-valent coupons where a permuation arc `crosses' an arc from a generating 1-morphism, and there are coupons corresponding to the generating 2-morphisms. These are subject to a host of local relations which are describe thoroughly in Section~\ref{sec:graphcalculusUSS}.

These structures compare very well with the structures making up the 2-category $\sB^{PD}$, namely the linear and planar diagrams. In fact for the choice of computad depicted in Figures \ref{MainUnOrientedTheormGenFig} and \ref{MainUnOrientedTheormRelFig} they coincide. For this unbiased semistrict symmetric monoidal computad, there is a single generating object. Thus the objects of $\cF(\cP)$ are words in just this single letter. In other words the objects of $\sF(\cP)$ are canonically identified with $\N_{\geq 0}$, and the symmetric monoidal structure is the same, given by concatenation. 

The 1-morphisms of $\sB^{PD}$, isotopy classes of linear diagrams, are also easily identified with the 1-morphisms of $\sF(\cP)$. The chambering set (which is labeled with permutations) and the points of the graphic become the two kinds of elementary 1-morphisms in $\sF(\cP)$. An isotopy class of linear diagram is exactly a composable sequence of these elementary 1-morphisms. 

The planar diagrams, which form the 2-morphisms of $\sB^{PD}$, may also be identified with the string diagrams which form the 2-morphisms of $\sF(\cP)$. The edges of the chambering graph of the planar diagram become the strings of the string diagram labeled by permutations. The trivalent and univalent vertices of the chambering graph become the distinguished trivalent and univalent coupons of the string diagram (which in turn correspond to the basic structure transformations $X^{\sigma, \sigma'}$, $X^e$, and their inverses). The arcs of the graphic in the planar diagram correspond exactly to the second kind of string in the string diagrams of $\sF(\cP)$, those labeled by the generating 1-morphisms. The crossing points of the graphic and chambering graph and the isolated points of the graphic each correspond to coupons in the string diagrams of $\sB^{PD}$.   
 
Some of these structures automatically form part of the definition of the string diagrams for unbiased semistrict symmetric monoidal 2-categories, and so are not specified in the presentation $\cP$. This is the case for the coupons coming solely from the chambering graph (the trivalent and univalent vertices) and for most, but not all of the coupons coming from crossings of the graphic and chambering graph. The remaining coupons will become the generating 2-morphisms (for example the `symmetry' generator describe earlier arises from a non-trivial crossing between the graphic and the chambering graph, one which is not automatically part of an unbiased semistrict symmetric monoidal string diagram). 

Moreover both the string diagrams and the planar diagrams are subject to a variety of relations. For string diagrams these are described in Section~\ref{sec:graphcalculusUSS}. For the planar diagrams these are the local moves which come from Table~\ref{LocalMovesForPlanarDiagramsTable1}, Table~\ref{LocalMovesForPlanarDiagramsTable2} (with conditions listed in Table~\ref{3DSheetDataTable}), Figure~\ref{fig:0stratumChamberMoves}, Figure~\ref{fig:localmoveschamberinggraph}, and Figure~\ref{fig:chambergraphicmoviemoves}. The relations satisfied by the planar diagrams include the relations that the string diagrams are subject to, but also some addition relations. These additional relations make up the relations in the computad $\cP$.  

Finally in identifying the computad $\cP$, we can map certain simplifications. For we do not need generating 2-morphism of $\cP$ for every possible way of labeling, say, a cusp singularity to make it a planar diagram. Such labelings will include trivial consequences such as concatenations with elements of $(n) \in \N_{\geq 0}$. It is sufficient to choose a minimal representative in each case. 

After making these simplifications and carefully comparing the string diagram relations from Section~\ref{sec:graphcalculusUSS} with the long list of relations satisfied by planar diagrams, we are left with the presentation listed in Figures \ref{MainUnOrientedTheormGenFig} and \ref{MainUnOrientedTheormRelFig}. \qed

\section{The Structured Classification}

In this section we show how to adapt the previous results to classify field theories for bordisms with structure. In the next section we will focus on the oriented bordism category as a concrete example. Let us review the main steps in the proof of Theorem~\ref{TheoremUnorientedClassification}, and see how they must be adapted in the setting of $\cF$-structures. 

First, in Chapter~\ref{SymMonBicatChapt} we developed the theory of symmetric monoidal computads, i.e., the theory of generators and relations for symmetric monoidal bicategories. These results allowed us to construct an abstract symmetric monoidal bicategory $\sF(\cP)$ from any such computad $\cP$. These symmetric monoidal bicategories are easy to map out of, we merely specify where each generator goes, subject to the defining relations. 

Moreover we also developed a parallel theory of unbiased semistrict symmetric monoidal 2-categories. This is a similar theory which is stricter than the general weak theory. A symmetric monoidal computad $\cP$ as gives rise to an unbiased semistrict symmetric monoidal 2-category $\sF_{uss}(\cP)$, presented by $\cP$. The canonical map $\sF(\cP) \to \sF_{uss}(\cP)$ is an equivalence of symmetric monoidal bicategories. Thus rather than comparing to the somewhat cumbersome $\sF(\cP)$, we can instead compare the bordism bicategory to $\sF_{uss}(\cP)$. This later has an explicit description in terms of string diagrams. Specifically the morphisms and 2-morphisms of $\sF_{uss}(\cP)$ may be explicitly describe by string diagrams in which the components of the string diagram (the strings and coupons) are labeled according to the generators of the presentation $\cP$. 

Next, we also constructed two equivalent versions of the bordism bicategory, by first equipping the bordisms with either linear or planar diagrams, and then forming an abstract bicategory which only has these diagrams and not the bordisms. The final step was to notice that this later abstract bicategory (really a 2-category) was exactly of the form $\sF(\cP)$ for a carefully chosen computad $\cP$. 

To carry out this strategy in the case of structured bordisms we need to adapt this second stage of the argument. We will need to modify our notion of planar diagram so that we obtain an analog of the planar decomposition theorem (Theorem~\ref{thm:bordismplanardecomp}) which is valid for $\cF$-manifolds and bordisms. After obtaining this we get a diagram, as before:
\begin{center}
\begin{tikzpicture}
		\node (LT) at (0, 1.5) {$\sF(P^\cF)$};
		\node (LB) at (0, 0) {$\sF_{uss}(P^\cF)$};
		\node (RRT) at (6, 1.5) {$\bord_2^\cF$};
		\node (RT) at (3, 1.5) {$\bord_2^{\cF,PD}$};
		\node (RB) at (3, 0) {$\sB^{\cF,PD}$};
		\draw [->] (LT) -- node [left] {$\simeq$} (LB);
		\draw [->] (LT) -- node [above] {$$} (RT);
		\draw [->] (RT) -- node [above] {$\simeq$} (RRT);
		\draw [->] (RT) -- node [right] {$\simeq$} (RB);
		\draw [->, dashed] (LB) -- node [above] {$\exists !$} (RB);
		%\node at (0.5, 1) {$\ulcorner$};
		%\node at (1.5, 0.5) {$\lrcorner$};
\end{tikzpicture}
\end{center}
The computad $P^\cF$ is chosen so as to make the dashed arrow an isomorphism. We will identify our new planar diagrams, and consequently the computad $P^\cF$, in stages working from objects to morphisms to 2-morphisms.

 %allow us to construct abstract symmetric monoidal bicategories which are easy to map out of. Thus it is relatively easy to construct an abstract candidate symmetric monoidal bicategory $\sB^\cF$ and a symmetric monoidal homomorphism $h:\sB^\cF \to \bord^\cF_2$. The main difficulty is to show that this is an equivalence. For this we must appeal to Theorem~\ref{WhiteheadforSymMonBicats}, which gives three criteria characterizing symmetric monoidal equivalences. For $h$ to be a symmetric monoidal equivalence, we must have that $h$ is essentially surjective on objects, essentially full on 1-morphisms, and fully-faithful on 2-morphisms.

%First we need $h$ to be essentially surjective on objects. 
Since any 0-dimensional $\cF$-manifold is a disjoint union of haloed points with $\cF$-structure, it is enough to identify these. Moreover we only need to understand them up to isomorphism {\em in the bordism bicategory}. In other words,  we want to classify $\cF$-structures on the point up to invertible $\cF$-bordism. 

\begin{step} \label{StepFStrOnPoint}
	Identify the $\cF$-structures on the point, up to invertible $\cF$-bordism. 
\end{step}

After we do this, we will pick a representative for each equivalence class and these will be the generating objects of $P^\cF$. %$\sB^\cF$. 
%$h: \sB^\cF \to \bord_2^\cF$ will then necessarily be essentially surjective on objects. 
In order to accomplish this first step, we must also have some rudimentary understanding of the $\cF$-bordisms, themselves. %This is important anyway, since the second criteria for $h$ to be an equivalence is that it is essentially full on 1-morphisms. 
Just as we understood the $\cF$-structure on the points, we will need to understand the $\cF$-structures on the 1-bordisms. 

Every 1-bordism is isomorphic to a disjoint union of intervals  and left and right ``elbows''. Thus it is enough to understand the $\cF$-structures on these. Again, we only need to understand them up to isomorphism in the bordism bicategory. Choosing representatives will lead to a collection of generating 1-morphisms for $P^\cF$. 
%The homomorphisms $h$ will now be essentially full on 1-morphisms. 

\begin{step} \label{StepFStrOnInverval}
Identify the $\cF$-structures on the interval, and on the left and right elbows, up to $\cF$-bordism.
\end{step}

Finally, we will need to construct a bijection between equivalence classes of our new planar diagrams and equivalence classes of $\cF$-2-bordisms. 
%Finally, we will need to prove that $h: \sB^\cF \to  \bord_2^\cF$ is fully-faithful on 2-morphisms. 
In the unoriented case, we proved this using the Planar Decomposition Theorem~\ref{thm:bordismplanardecomp}
To prove this in the setting of $\cF$-structures we must adapt the planar decomposition theorem to the setting of $\cF$-structures. The Planar Decomposition Theorem decomposes into simple parts. First there is the local analysis of maps to $\R^2$. This doesn't change in the setting of $\cF$-structures. 

The local analysis allows us to do several things. Given an $\cF$-2-bordism, we may choose a generic map to $\R^2$. Via the graphic, and after making some choices,  this gives us a way to cut up our surface, in a bicategorical fashion,  into elementary pieces. We must first understand what $\cF$-structures are allowed on these pieces.

\begin{step} \label{StepStrOnElemBord}
	Identify the isomorphism classes (rel. boundary) of $\cF$-structures on the elementary 2-bordisms.\end{step}

Secondly, we must understand how to glue these elementary pieces back together to recover the original $\cF$-2-bordism.  In the unoriented case, this gluing information was encoded in what we called {\em sheet data}. There will be an analogous Planar Decomposition Theorem in the setting of $\cF$-structures and the main difference between this and the unoriented Planar Decomposition Theorem is in the sheet and gluing data. 

In fact, we have already started the process of determining the new sheet data in the presence of $\cF$-structures. Consider a  small ball in the graphic which contains no critical points. The inverse image of this region consists of a disjoint union of disks with $\cF$-structure. Choosing a point $p$, in each disk we get a $\cF$-structure on a neighborhood of this point, which in turn gives us an $\cF$-0-bordism. Each disks is isomorphic to $p \times I \times I$.  Thus the new sheet data that we should assign to a region without  critical points as a collection of sets, one for each of the generating objects we identified in Step \ref{StepFStrOnPoint}. Similarly Steps \ref{StepFStrOnInverval} and \ref{StepStrOnElemBord} tell us precisely what sheet data we should associate along folds, along fold crossings, along cusps, and along 2D Morse critical points. 

The next thing to understand is the appropriate gluing data that we should assign to overlaps. But this is also part of Step \ref{StepFStrOnInverval}, namely the part which understands the $\cF$-structures on intervals. Similarly there will be cocycle conditions on the triple intersections. 
 
\begin{step}
Use the previous steps to determine the sheet and gluing data for planar decompositions of surfaces with $\cF$-structure. 
\end{step} 

Finally, we must understand when two of these $\cF$-planar decompositions result in the same $\cF$-manifold. Theorem~\ref{PlanarDecompMovesThm} gives a precise list of local moves which allow us to pass from one planar decomposition to another. We simply need to understand how to add in appropriate sheet data. 

\begin{step}
 Determine which $\cF$-sheet data is appropriate for the elementary moves of Theorem~\ref{PlanarDecompMovesThm}. 
\end{step}

After completing these steps one obtains a new version of the planar decomposition theorem. The sheet and gluing data is augmented with additional $\cF$-structure information. The generating 2-morphisms and relations of $P^\cF$ are now determined by this sheet data.
% and consequently $h: \sB^\cF \to \bord_2^\cF$ will be fully faithful on 2-morphisms, and hence a symmetric monoidal equivalence. 

\section{The Oriented Classification} \label{SectOrientedClassification}

To demonstrate that this program is actually feasible, the rest of this section is devoted to carrying it out in the case of orientations. Let $\text{Or}_d$ be the symmetric monoidal category of oriented $d$-manifolds, with embeddings as morphisms. The forgetful functor $\text{Or}_d \to \man^d$, realizes this as a stack over $\man^d$ which is topological in the sense of Definition \ref{DefTopologicalStructure}. In fact it is a sheaf; orientations have no automorphisms.

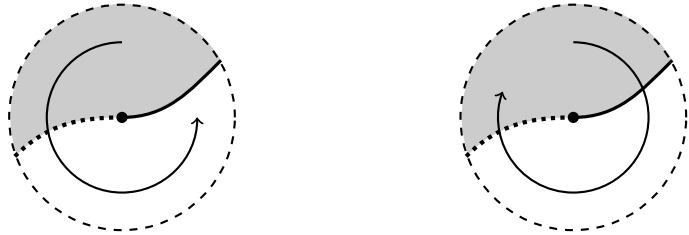
\begin{figure}[ht]
\begin{center}
\begin{tikzpicture}[thick]
\begin{scope}[xshift= 6cm]
		\node [circle, minimum width = 3cm] (A) at (0,0) {};
	\node (B) at (A.30) {};
	\node (C) at (A.200) {};

	\fill [black!20]  (C.center) to [out = 45, in = 180] (A.center) to [out = 0, in = 225] (B.center) 
		arc (30 : 200 : 1.5cm); 
	\draw [very thick, black] (A.center) to [out = 0, in = 225] (B.center);
	\draw [ultra thick, dotted, black] (C.center) to [out = 45, in = 180] (A.center);
	\draw [dashed] (A.center) circle (1.5cm);
	\node [circle, fill=black,inner sep=1.5pt] at (A.center) {};
	\draw [->] (0, 1) arc (90:-200 : 1cm);
\end{scope}

	\node [circle, minimum width = 3cm] (A) at (0,0) {};
	\node (B) at (A.30) {};
	\node (C) at (A.200) {};
	\fill [black!20]  (C.center) to [out = 45, in = 180] (A.center) to [out = 0, in = 225] (B.center) 
		arc (30 : 200 : 1.5cm); 
	\draw [very thick, black] (A.center) to [out = 0, in = 225] (B.center);
	\draw [ultra thick, dotted, black] (C.center) to [out = 45, in = 180] (A.center);
	\draw [dashed] (A.center) circle (1.5cm);
	\node [circle, fill=black,inner sep=1.5pt] at (A.center) {};
	\draw [->] (0, 1) arc (90: 360 : 1cm);
	
\end{tikzpicture}
\caption{Two Oriented 0-Bordisms}
\label{FigOrientedZeroBords}
\end{center}
\end{figure}

The five steps outlined in this section are made easier in the oriented case by the following observation: a 2-manifold which is topologically a disk has exactly two orientations. This implies that there are exactly two oriented 0-bordisms whose underlying zero manifold is a point. Let us fix two representatives, call them the {\em positive point} and the {\em negative point}.  Similarly the haloed interval, viewed as a haloed 1-bordism from the haloed point to itself supports precisely two orientations. These correspond to the identity 1-bordisms from the positive point, respectively negative point, to itself. In particular there are no oriented 1-bordisms from the positive point to the negative point. Thus we see that even up to oriented bordism there are precisely two orientations on the point, see Figure~\ref{FigOrientedZeroBords}.

\begin{figure}[ht]
\begin{center}
\begin{tikzpicture}[thick]
	\node (A) at (2,4) {};
	\node (B) at (2, 3) {};
	\node (X) at (2, 1) {};
	\node (Y) at (2,0) {};
	\fill [black!20] (A.center) -- +(-1, 0) arc (180 : 90 : 1cm)  arc (90: -90: 3cm)
			arc (270: 180: 1cm) -- (Y.center) arc (-90: 90: 2cm);
	
		arc (90 : 0 : 1cm) -- (Y.center) to [out = 180, in = 0] (A.center);
	\draw [dashed] (B.center) arc (270 : 90: 1cm)  arc (90: -90: 3cm)
			arc (270: 90: 1cm) arc (-90: 90: 1cm);
	
	\draw (A.center) -- +(-1, 0) (A.center) arc (90:-90: 2cm ) (Y.center) -- +(-1,0); 
	\node [circle, fill=black,inner sep=1pt] at (A.center) {};
	\node [circle, fill=black,inner sep=1pt] at (Y.center) {};

	\node (C) at (-2,4) {};
	\node (D) at (-2, 0) {};
	\fill [black!20] (C.center) -- +(-1, 0) arc (180 : 0: 1cm) -- (C.center);
	\fill [black!20] (D.center) -- +(-1, 0) arc (180 : 360: 1cm) -- (D.center);
	\draw [dashed] (C.center) circle (1cm);
	\draw [dashed] (D.center) circle (1cm);
	\draw (C.center) -- +(1,0) (D.center) -- +(1,0);
	\draw [ultra thick, dotted] (C.center) -- +(-1,0) (D.center) -- +(-1,0);
	\node [circle, fill=black,inner sep=1pt] at (C.center) {};
	\node [circle, fill=black,inner sep=1pt] at (D.center) {};
	
	\draw [->] (0, 2) to node [above] {$f_\text{in}$} (1, 2);
	%\draw [->] (8, 1.5) to [out = 90, in = 0] node [above right] {$f_\text{out}$} (7.5, 2);
	\foreach \x / \y in {2/4.5,3/4.5, 2/3.5, 3/3, 4.5/3, 4.5/2, 3/1, 2/0.5, 2/-0.5, 3.5/0,
			 -2/4.5, -2/3.5, -2/0.5, -2/-0.5}
		{
		\draw[thin, stealth-] (\x, \y) -- +(0, 0.1) arc (0: 200 : 0.125cm);
		}
	
\end{tikzpicture}	
\caption{An Oriented 1-Bordism}
\label{FigOrientedElbow}
\end{center}
\end{figure}
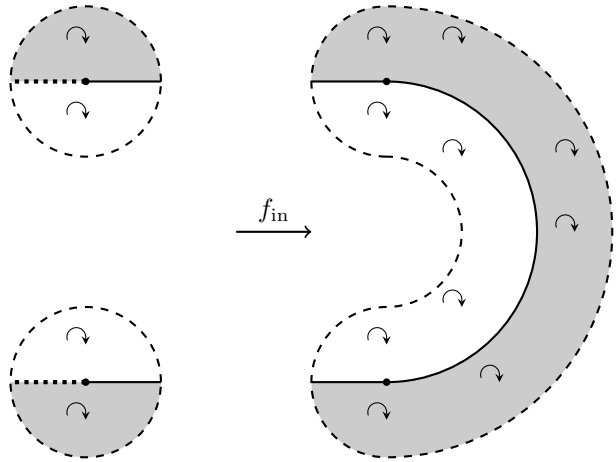

To complete the second step, we must understand the orientations on the right and left ``elbows''. Let us start with the elbow whose source consists of the disjoint union of two points and whose target is the empty set. There are four possible orientations on the disjoint union of two points, and only two of these are compatible with orientations of the elbow. See Figure~\ref{FigOrientedElbow} for an example. 
The sources of these two oriented elbows are either $pt^+ \sqcup pt^-$ or $pt^- \sqcup pt^+$. Thus there are a priori four generating 1-bordisms. We will see later,  however, that two of these generators are redundant and that we can eliminate them and get an equivalent symmetric monoidal bicategory. 

Each of the elementary 2-bordisms is topologically a disk, so again each has exactly two orientations. Given an orientation of the boundary of an elementary 2-bordism, there is precisely one orientation on the 2-bordism extending this orientation. Thus the elementary 2-bordisms are just as they are for the unoriented bordism bicategory, except that the boundaries are equipped with orientations. This completes Step 3.

We can now describe how to modify the notion of a planar diagram in order to incorporate orientations. The main difference between the oriented and unoriented cases will be in the sheet data. The oriented  sheet data for regions without critical points now consists of {\em two} sets of sheets, those corresponding to the positively oriented point and those corresponding to the negatively oriented point. Similarly the sheet data for the fold and 2D Morse regions consist of a quadruple of sets $(S_+, S_-, \{t_+\}, \{ t_-\})$ such that $(S_+ \sqcup S_-, \{ t_+, t_-\})$ forms sheet data in the unoriented sense. The fold crossing sheet data is similar. Finally cusp sheet data consists of a pair of maps of pairs of sets
\begin{equation*}
	(S_+, S_-) \rightrightarrows (S'_+, S'_-).
\end{equation*}
These are required to satisfy that $(S_+ \sqcup S_-) \rightrightarrows (S'_+ \sqcup S'_-)$ consists of unoriented cusp sheet data. 

The gluing data for regions without singularities is given simply by a bijection of pairs of sets, $(S_+, S_-) \cong (\overline S_+, \overline S_-)$. On triple intersections these must satisfy the obvious cocycle condition. The gluing data for regions which contain fold singularities are more interesting. The gluing data from $(S_+,  S_-, \{ t_+\} , \{t_-\})$ to $(\overline S_+, \overline S_-, \{ \overline t_+\} , \overline\{t_-\})$, consists of a pair of bijections $S_+ \cong \overline S_+$, $S_- \cong \overline S_-$. No automorphisms of $T = \{ t_+, t_- \}$ are allowed. Translating this into generators for the oriented bordism bicategory, we see that this has the effect of eliminating  the symmetry generators from the list of generating 2-morphisms. A similar analysis shows that each of the relations for the unoriented bordism bicategory holds whenever the oriented bordisms are composable.  Putting these results together yields the following oriented classification theorem.

\begin{theorem} \label{OrientedBordismClassification}
The oriented 2-dimensional bordism bicategory $\bord_2^\text{or}$, as a symmetric monoidal bicategory, has the generators and relations  depicted in Figure~\ref{fig:OrientedGenAndRel}.
\begin{figure}[ht]
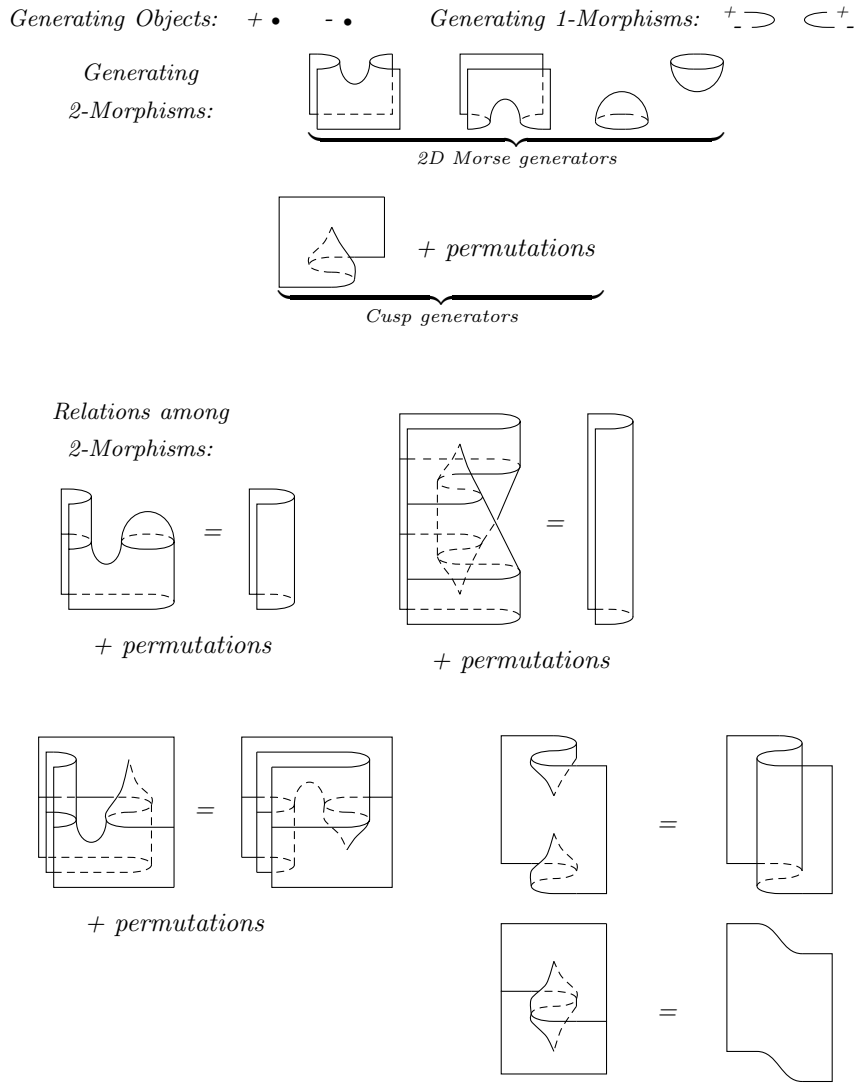

	%\begin{center}

	};

	%%% Cusp Inversion AI
	\begin{scope}[xshift = -1cm, yshift = -9.5cm]
	\draw (0.5, 1.2) -- (-0.2,1.2) -- (-0.2,-0.5) -- (0.2, -0.5);
	\draw (0.5, 1.2) arc (90: -90: 0.3cm and 0.1cm) arc (90: 270: 0.3cm and 0.1cm) -- (0.8, 0.8);
	\draw [densely dashed] (0.2, -0.5) -- (0.5,-0.5) arc (90: -90: 0.3cm and 0.1cm) arc (90: 180: 0.3cm and 0.1cm);
	\draw (0.2, -0.8) arc (180: 270: 0.3cm and 0.1cm) -- (0.8, -0.9);
	\draw (0.8, 0.8) -- (1.2, 0.8) -- (1.2, -0.9) -- (0.8, -0.9);
	\node (C) at (0.5, 0.4) {}; 
	\node (D) at (0.5, -0.1) {};
	\draw  plot [smooth] coordinates{(0.8,1.1) (0.8,0.9)}; 
	\draw [densely dashed] plot [smooth] coordinates{ (0.8,0.9) (0.6, 0.6) (C.center)}; 
	\draw  plot [smooth] coordinates{(0.2,0.9) (0.2,0.8) (0.4, 0.6) (C.center)}; 
	\draw  plot [smooth] coordinates{(0.2,-0.8) (0.2,-0.6) (0.4, -0.4) (D.center)}; 
	\draw [densely dashed] plot [smooth] coordinates{(0.8,-0.6) (0.8,-0.5) (0.6, -0.3) (D.center)}; 
	\end{scope}

	%%% Cusp Inversion AII
	\begin{scope}[xshift = 2cm, yshift = -9.5cm]
	\draw (0.5, 1.2) -- (-0.2,1.2) -- (-0.2,-0.5) -- (0.2, -0.5); %(0.5, -0.5) ;
	\draw (0.5, 1.2) arc (90: -90: 0.3cm and 0.1cm) arc (90: 270: 0.3cm and 0.1cm) -- (0.8, 0.8);
	\draw [densely dashed] (0.2, -0.5) -- (0.5,-0.5) arc (90: -90: 0.3cm and 0.1cm) arc (90: 180: 0.3cm and 0.1cm);
	\draw (0.2, -0.8) arc (180: 270: 0.3cm and 0.1cm) -- (0.8, -0.9);
	\draw (0.8, 0.8) -- (1.2, 0.8) -- (1.2, -0.9) -- (0.8, -0.9);
	\node (C) at (0.5, 0.4) {};
	\draw (0.2, 0.9) -- (0.2, -0.8) (0.8, 1.1) -- (0.8, -0.6);
	\end{scope}

	\node at (1, -9.5) {=};
	\node at (1, -12) {=};

	%%% Cusp Inversion BI
	\begin{scope}[xshift = -1cm, yshift = -11.2cm]
	\draw (1.2, 0.4) -- (-0.2,0.4) -- (-0.2,-0.5) -- (0.2, -0.5);
	\draw [densely dashed] (0.2, -0.5) -- (0.5,-0.5) arc (90: -90: 0.3cm and 0.1cm) arc (90: 180: 0.3cm and 0.1cm);
	\draw (0.2, -0.8) arc (180: 270: 0.3cm and 0.1cm) -- (0.8, -0.9);
	\draw  (1.2, 0.4) -- (1.2, -0.9) -- (0.8, -0.9);
	\draw (1.2, -0.9) -- (1.2, -1.6) -- (-0.2, -1.6) -- (-0.2, -0.5); 
	\node (D) at (0.5, -0.1) {};
	\draw  plot [smooth] coordinates{(0.2,-0.8) (0.2,-0.6) (0.4, -0.4) (D.center)}; 
	\draw [densely dashed] plot [smooth] coordinates{(0.8,-0.6) (0.8,-0.5) (0.6, -0.3) (D.center)}; 
	\node (C) at (0.5, -1.3) {}; 
	\draw  plot [smooth] coordinates{(0.8,-0.6) (0.8,-0.7)}; 
	\draw [densely dashed] plot [smooth] coordinates{ (0.8,-0.7) (0.6, -1) (C.center)}; 
	\draw  plot [smooth] coordinates{(0.2,-0.8) (0.2,-0.9) (0.4, 0.-1.1) (C.center)}; 
	\end{scope}

	%%% Cusp Inversion BII
	\begin{scope}[xshift = 2cm, yshift = -12cm]
	\draw  (0.8, 0.8) to [out = 180, in = 0] (0.2, 1.2) -- (-0.2,1.2) -- (-0.2,-0.5) -- (0.2, -0.5) to [out = 0, in = 180] (0.8, -0.9);
	%\draw (0.5, 1.2) arc (90: -90: 0.3cm and 0.1cm) arc (90: 270: 0.3cm and 0.1cm) -- (0.8, 0.8);
	%\draw [densely dashed] (0.2, -0.5) -- (0.5,-0.5) arc (90: -90: 0.3cm and 0.1cm) arc (90: 180: 0.3cm and 0.1cm);
	%\draw (0.2, -0.8) arc (180: 270: 0.3cm and 0.1cm) -- (0.8, -0.9);
	\draw (0.8, 0.8) -- (1.2, 0.8) -- (1.2, -0.9) -- (0.8, -0.9);

	\end{scope}

	%\end{scope}

	 \end{tikzpicture}
	%\end{center}
	\caption{Oriented Generators and Relations}
	\label{fig:OrientedGenAndRel}	
\end{figure}

\end{theorem}

\section{Transformations and Modifications of TFTs} \label{sec:transtqft}

Now that we have deduced a generators and relations presentation of the oriented and unoriented bordism bicategories we can understand very concretely the bicategories of homomorphisms out of these categories, i.e., the bicategories of oriented and unoriented TFTs. Theorem~\ref{thm:cofibrancythm} and Proposition~\ref{RelationsForSMBicatsProp} are essential for this purpose.   They allow us to characterize these homomorphism bicategories in terms of a specific small amount of data which corresponds to the images of the generators in the target symmetric monoidal bicategory. 
In this section we will collect together a few basic results on transformations between topological field theories, focusing only on the oriented case. The existence of the forgetful symmetric monoidal homomorphism $\bord_2^\text{or} \to \bord_2$ implies that any unoriented topological field theory gives rise to an oriented theory, so the results of this section are also valid in the unoriented setting. 

Let us fix a target symmetric monoidal bicategory $\sC$, and  assume that we have two topological field theories $Z_0$ and $Z_1$ with values in $\sC$. We will first work in the oriented case, and then explain what happens in the unoriented setting.   Theorem~\ref{thm:cofibrancythm} and Proposition~\ref{RelationsForSMBicatsProp} allow us to characterize the transformations between $Z_0$ and $Z_1$. They are determined by their values on the generating objects and 1-morphisms. 
In the oriented case there are only two generating objects and two generating 1-morphisms. Thus a transformation between TFTs is given by specifying data:
\begin{align*}
	\sigma(pt^+) : & \; Z_0(pt^+) \to Z_1(pt^+) \\
	\sigma(pt^-) : & \;Z_0(pt^-) \to Z_1(pt^-) \\
	\sigma( \tikz[baseline = -0.7cm]{
		\draw (-5.5, -.5) arc (90: 270: 0.3cm and 0.1cm) -- (-5.4, -.7);
		\node at (-5.3, -.5) {\tiny +};
		\node at (-5.2, -.7) {-};
		}
	) : & \; Z_1( \tikz[baseline = -0.7cm]{
		\draw (-5.5, -.5) arc (90: 270: 0.3cm and 0.1cm) -- (-5.4, -.7);
		\node at (-5.3, -.5) {\tiny +};
		\node at (-5.2, -.7) {-};
		}) \circ \sigma( pt^+ \sqcup pt^-) \to \sigma(\emptyset) \circ Z_0( \tikz[baseline = -0.7cm]{
		\draw (-5.5, -.5) arc (90: 270: 0.3cm and 0.1cm) -- (-5.4, -.7);
		\node at (-5.3, -.5) {\tiny +};
		\node at (-5.2, -.7) {-};
		})  \\
	\sigma(\tikz[baseline = -0.7cm]{
		\draw (-6.6, -.5) -- (-6.5, -.5) arc (90: -90: 0.3cm and 0.1cm);
		\node at (-6.8, -.5) {\tiny +};
		\node at (-6.7, -.7) {-};
	}) : &\; Z_1(\tikz[baseline = -0.7cm]{
		\draw (-6.6, -.5) -- (-6.5, -.5) arc (90: -90: 0.3cm and 0.1cm);
		\node at (-6.8, -.5) {\tiny +};
		\node at (-6.7, -.7) {-};
	}) \circ \sigma(\emptyset) \to \sigma( pt^+ \sqcup pt^-) \circ Z_0(\tikz[baseline = -0.7cm]{
		\draw (-6.6, -.5) -- (-6.5, -.5) arc (90: -90: 0.3cm and 0.1cm);
		\node at (-6.8, -.5) {\tiny +};
		\node at (-6.7, -.7) {-};
	})
\end{align*}
This data is subject to a number of relations, each coming from one of the generating 2-morphisms.

\begin{figure}[ht]
	
\begin{center}
\begin{tikzpicture}
\node (A) at (0,3) {$Z_0(pt^+)$};
\node (B) at (2,1.5) {$Z_1(pt^+)$};
\node (C) at (2,0) {$Z_1(pt^-)$};
\node (D) at (8,3) {$Z_1(pt^+)$};
\node (E) at (6,1.5) {$Z_0(pt^+)$};
\node (F) at (6,0) {$Z_0(pt^-)$};

\draw (A) to [out = 0, in = 180] (E) 
	(E) to [out = 0, in = 0] node [right] {$Z_0(\tikz[baseline = -0.7cm]{
		\draw [-] (-6.6, -.5) -- (-6.5, -.5) arc (90: -90: 0.3cm and 0.1cm);
		\node at (-7, -.5) {\tiny +};
		\node at (-6.9, -.7) {-};
	})$} (F) 
	(F)-- node [below] {$\sigma(pt^-)$} (C) 
	(C) to [out = 180, in = 180] node [left] {$ Z_1( \tikz[baseline = -0.7cm]{
		\draw (-5.5, -.5) arc (90: 270: 0.3cm and 0.1cm) -- (-5.4, -.7);
		\node at (-5.1, -.5) {\tiny +};
		\node at (-5, -.7) {-};
		})$} (B) 
	(B) to [out = 0, in = 180] (D);
\end{tikzpicture}
\caption{Graphical Depiction of Composition (read right-to-left)}
\label{Fig:GraphDepictionComp}
\end{center}
\end{figure}
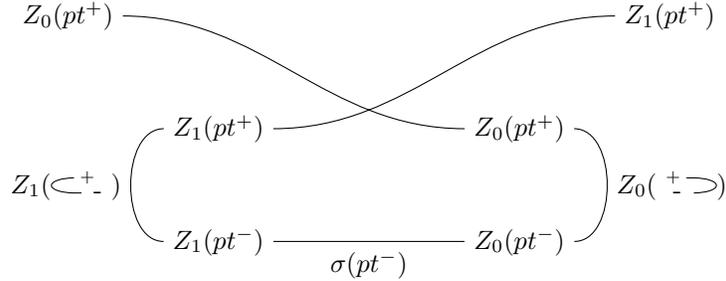

The significance of this data will become clearer if we repackage it slightly. We can express it entirely in terms of data involving only the images of the positively oriented point, $pt^+$. The 1-morphism $\sigma(pt^-) :  Z_0(pt^-) \to Z_1(pt^-)$
in $\sC$ gives rise to a new 1-morphism from $Z_1(pt^+)$ to $Z_0(pt^+)$, via the composition,
\begin{equation*}
	\tilde \sigma(pt^+) : = [ Z_1( \tikz[baseline = -0.7cm]{
		\draw (-5.5, -.5) arc (90: 270: 0.3cm and 0.1cm) -- (-5.4, -.7);
		\node at (-5.3, -.5) {\tiny +};
		\node at (-5.2, -.7) {-};
		}) \otimes Id_{Z_0(pt^+)} ] \circ [ \sigma(pt^-) \otimes \beta_{Z_0(pt^+), Z_1(pt^+)} ] \circ
		[ Z_0(\tikz[baseline = -0.7cm]{
		\draw (-6.6, -.5) -- (-6.5, -.5) arc (90: -90: 0.3cm and 0.1cm);
		\node at (-6.8, -.5) {\tiny +};
		\node at (-6.7, -.7) {-};
	}) \otimes Id_{Z_1(pt^+)} ]
\end{equation*}
This composition can be visualized graphically, see Figure~\ref{Fig:GraphDepictionComp}.
\noindent Using the images of the folds and cusps under the TFTs $Z_0$ and $Z_1$, one can recover $\sigma(pt^-)$ from $\tilde \sigma(pt^+)$. Thus the first part of this data can be expressed as a pair of morphims,
\begin{equation*}
	\sigma: Z_0(pt^+) \rightleftarrows Z_1(pt^+): \tilde \sigma
\end{equation*}
The rest of the data can, likewise,  be expressed equivalently in terms of $\sigma(pt^+)$ and $\tilde \sigma$. For example, using the image under $Z_0$ of one of the cusps, we can form the following 2-morphism $ \varepsilon = \tilde \sigma(\tikz[baseline = -0.7cm]{
		\draw (-5.5, -.5) arc (90: 270: 0.3cm and 0.1cm) -- (-5.4, -.7);
		\node at (-5.3, -.5) {\tiny +};
		\node at (-5.2, -.7) {-};
		})$, depicted in Figure~\ref{Fig:AComposite2Morphism}. 
\begin{figure}[ht]
\begin{center}

\label{Fig:AComposite2Morphism}
\caption{A 2-morphism}
\end{center}
\end{figure}		

\begin{figure}[ht]
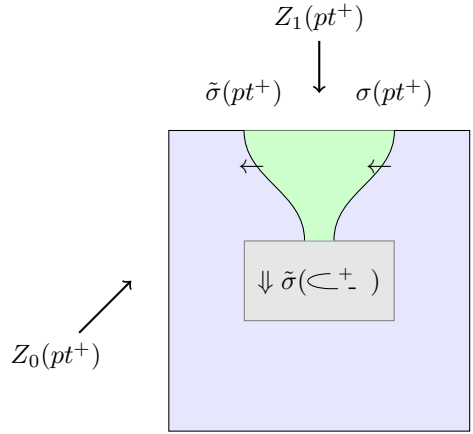

\begin{center}
%
\caption{A Pasting Diagram For Transformations of TFTs}
\label{FigPastingForTransofTFTs1}
\end{center}
\end{figure}

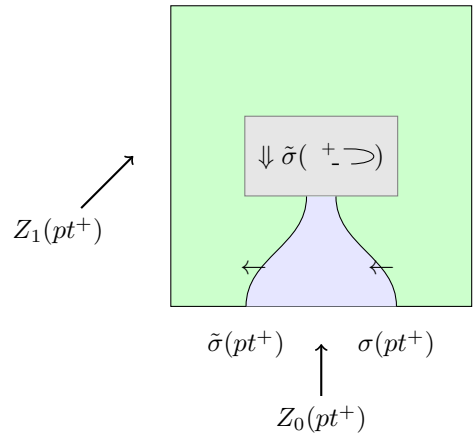
\begin{figure}[ht]
\begin{center}
\begin{tikzpicture}
	\draw [fill = green!20] (0,0) rectangle +(4,4);
	
	\node (A) [stringnode] at (2, 2) {$\Downarrow \tilde \sigma(\tikz[baseline = -0.7cm]{
		\draw (-6.6, -.5) -- (-6.5, -.5) arc (90: -90: 0.3cm and 0.1cm);
		\node at (-6.8, -.5) {\tiny +};
		\node at (-6.7, -.7) {-};
	})$ };
		
	\fill [blue!10] (1,0) to [out = 90, in = -90] (A.250) -- (A.290) to [out = -90, in = 90] (3,0) -- (1,0);
	\draw (1,0) to [out = 90, in = -90] (A.250);
	\draw (3,0)  to [out = 90, in = -90] (A.290);	
	\node (X) at (1, -0.5) {$\tilde \sigma(pt^+)$};
	\node (Y) at (3, -0.5) {$\sigma(pt^+)$};
	\node at (1.1, 0.5) {$\leftarrow$};
	\node at (2.8, 0.5) {$\leftarrow$};
	\node (B) at (2, -1.5) {$Z_0(pt^+)$};
	\draw [->, thick] (B) -- +(0, 1);
	\node (C) at (-1.5, 1) {$Z_1(pt^+)$};
	\draw [->, thick] (C) -- +(1, 1);
\end{tikzpicture}
\caption{A second Pasting Diagram For Transformations of TFTs}
\label{FigPastingForTransofTFTs2}
\end{center}
\end{figure}

It is perhaps more illuminating to view this morphism by representing it as a string diagram (see Appendix \ref{PastingsStringsAdjointsAndMatesSection}), as in Figure~\ref{FigPastingForTransofTFTs1}.
Using the cusps and relations from $\bord_2^\text{or}$, we can recover $ \sigma(\tikz[baseline = -0.7cm]{\draw (-5.5, -.5) arc (90: 270: 0.3cm and 0.1cm) -- (-5.4, -.7); \node at (-5.3, -.5) {\tiny +};
		\node at (-5.2, -.7) {-}; })$ from $\tilde \sigma(\tikz[baseline = -0.7cm]{
		\draw (-5.5, -.5) arc (90: 270: 0.3cm and 0.1cm) -- (-5.4, -.7);
		\node at (-5.3, -.5) {\tiny +};
		\node at (-5.2, -.7) {-};
		})$, and so it is equivalent  data defining the transformation. Similarly, the 2-morphism 
$\sigma(\tikz[baseline = -0.7cm]{
		\draw (-6.6, -.5) -- (-6.5, -.5) arc (90: -90: 0.3cm and 0.1cm);
		\node at (-6.8, -.5) {\tiny +};
		\node at (-6.7, -.7) {-};
	})$ is equivalent data to a 2-morphism $ \eta = \tilde \sigma(\tikz[baseline = -0.7cm]{
		\draw (-6.6, -.5) -- (-6.5, -.5) arc (90: -90: 0.3cm and 0.1cm);
		\node at (-6.8, -.5) {\tiny +};
		\node at (-6.7, -.7) {-};
	})$ which fills the pasting diagram of Figure~\ref{FigPastingForTransofTFTs2}.

We may now express the additional relations that $\sigma$ must satisfy in terms of this new equivalent data. These relations express the naturality of the assignment $f \mapsto \sigma_f$, for 1-bordisms $f$. Thus for each 2-morphsims $\alpha: f \to g$, we get a relation which expresses the naturality with respect to $\alpha$. It is enough to just check naturality for each of the generating 2-morphisms. In the oriented case there are roughly two groups of generating 2-morphisms. There are the invertible 2-morphisms, which are the cusps, and the non-invertible 2-morphisms, which are the four Morse generators.

\begin{figure}[htbp]
	\begin{center}
	\begin{tikzpicture}
	
		\node (A) [stringnode, minimum width = 1.5cm] at (1,1) {};
		\node (B) [stringnode, minimum width = 1.5cm] at (2,2) {};
		\fill [blue!10] (0,0) -- (0,3) -- (0.75, 3) -- (0.75,3 |- A.north) 
			-- (1.5, 3 |- A.north) -- (1.5, 3 |- B.south) 
			 -- (2.25, 0 |- B.south) -- (2.25, 0) -- (0,0);
		\fill [green!20] (3,3) -- (0.75, 3) -- (0.75,3 |- A.north) 
			-- (1.5, 3 |- A.north) -- (1.5, 3 |- B.south) 
			 -- (2.25, 0 |- B.south) -- (2.25, 0) -- (3,0) -- (3,3);
		\draw (0,0) rectangle (3, 3);
		\draw (0.75, 3) -- (0.75,3 |- A.north) 
			(1.5, 3 |- A.north) -- (1.5, 3 |- B.south) 
			 (2.25, 0 |- B.south) -- (2.25, 0);
		\node  [stringnode, minimum width = 1.5cm] at (1,1) {};
		\node  [stringnode, minimum width = 1.5cm] at (2,2) {};

		\node at (4, 1.5) {=};
		\filldraw [fill = blue!10] (5, 0) rectangle (6,3);
		\filldraw  [fill= green!20] (6, 0) rectangle (7,3);
	
	\begin{scope}[xshift = 3cm, yshift = -3.5cm, xscale = -1]
	\node (A) [stringnode, minimum width = 1.5cm] at (1,1) {};
		\node (B) [stringnode, minimum width = 1.5cm] at (2,2) {};
		\fill [blue!10] (0,0) -- (0,3) -- (0.75, 3) -- (0.75,3 |- A.north) 
			-- (1.5, 3 |- A.north) -- (1.5, 3 |- B.south) 
			 -- (2.25, 0 |- B.south) -- (2.25, 0) -- (0,0);
		\fill [green!20] (3,3) -- (0.75, 3) -- (0.75,3 |- A.north) 
			-- (1.5, 3 |- A.north) -- (1.5, 3 |- B.south) 
			 -- (2.25, 0 |- B.south) -- (2.25, 0) -- (3,0) -- (3,3);
		\draw (0,0) rectangle (3, 3);
		\draw (0.75, 3) -- (0.75,3 |- A.north) 
			(1.5, 3 |- A.north) -- (1.5, 3 |- B.south) 
			 (2.25, 0 |- B.south) -- (2.25, 0);
		\node  [stringnode, minimum width = 1.5cm] at (1,1) {};
		\node  [stringnode, minimum width = 1.5cm] at (2,2) {};
	\end{scope}
	
		\node at (4, -2) {=};
		\filldraw [fill = green!20] (5, -0.5) rectangle (6,-3.5);
		\filldraw  [fill= blue!10] (6, -0.5) rectangle (7,-3.5);
	\end{tikzpicture}
	\end{center}
	\caption{Relations from cusp naturality}
	\label{fig:cuspnat}
\end{figure}
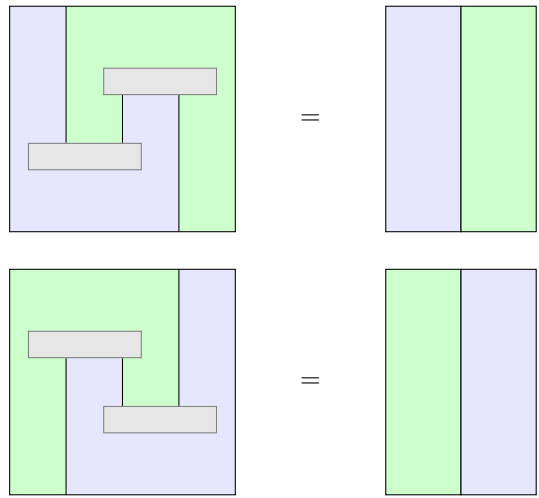

\begin{proposition} \label{PropNatWRTCusp}
	Naturality with respect to the cusp 2-morphisms is equivalent to the relations depicted in Figure~\ref{fig:cuspnat}. 
In particular, these form an adjunction between $Z_0(pt^+)$ and $Z_a(pt^+)$. 
\end{proposition}

\begin{proof}
This follows directly by writing out the equations which express the naturality with respect to the cusp morphisms, one merely needs to identify both sides of the equation. 
\end{proof}

Transformations, as defined in Definition \ref{DefnBicatTransformation}, have invertible 2-morphism data. However, up to this point we have not used this invertibility fact in any way. It is interesting to understand what happens when we don't impose this invertibility constraint, i.e., when we consider lax and oplax transformations, see Remark \ref{RmkLaxAndOpLax}. As we will see shortly, any lax or oplax transformation of oriented (and hence unoriented) topological field theories automatically satisfies the required invertibility condition; any lax or oplax transformation is automatically a {\em strong} transformation.

Proposition~\ref{PropNatWRTCusp} gives us justification for introducing a graphical notation which expresses these relations. Compare with Appendix \ref{PastingsStringsAdjointsAndMatesSection}.
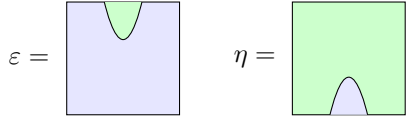
\begin{center}

\end{center}
This graphical notation can also be extended to the original data describing the transformation $\sigma$. In this notation the naturality with respect to the cusp looks like the diagram in Figure~\ref{FigNatWRTCusp}.
\begin{center}
%

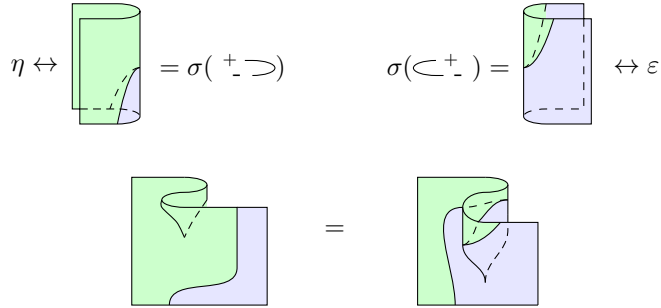
\caption{Naturality of Transformations with respect to a Cusp Generator}
\label{FigNatWRTCusp}
\end{center}
\end{figure}

One of the most fascinating aspects of the bordism bicategory is how much duality is present. The generating 1-morphisms, the cusp 2-morphism, and their relations exhibit a duality between the positive and negative points. We have seen one aspect of this, namely that a lax transformation of field theories gives rise to a right-adjunction $\sigma: Z_0(pt^+) \to Z_1(pt^+)$. However there is much more structure. The 2D Morse generating 2-morphisms, together with the 2D Morse relations realize that fact that the two generating 1-morphisms, the left- and right-elbows,  are adjoint to one another. In fact they show that the generating 1-morphisms are {\em simultaneously} left-adjoint and right-adjoint. This simultaneous left/right adjunction  is called an {\em ambidextrous adjunction}. 

\begin{figure}[ht]
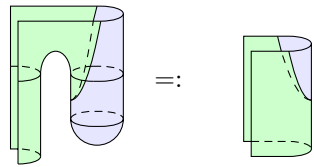

\begin{center}



\caption{A Mate for a Transformation of Field Theories}
\label{FigMateforTransofTFTs}
\end{center}
\end{figure}

The significance of this for transformations of topological field theories is that such
 transformations are defined in terms of certain data, certain 1-morphism and 2-morphism. The 2-morphism have source and targets which are 1-morphisms which have adjoints. Consequently, as we have seen in Appendix \ref{PastingsStringsAdjointsAndMatesSection}, because these generating 1-morphisms have left- and right-adjoints, we get a host of associated 2-morphisms: the mates. For example Figure~\ref{FigMateforTransofTFTs} shows the mate of $ \sigma(\tikz[baseline = -0.7cm]{
		\draw (-5.5, -.5) arc (90: 270: 0.3cm and 0.1cm) -- (-5.4, -.7);
		\node at (-5.3, -.5) {\tiny +};
		\node at (-5.2, -.7) {-};
		}) $.

Just as we turned the structure 2-morphisms for a transformation into 2-morphisms involving just the positive point, we can likewise turn the mates of these structure 2-morphisms into 2-morphisms involving just the positive point. These can also be given a graphical notation:
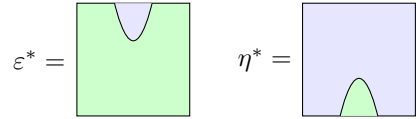
\begin{center}
%
\end{center}  
It is an amusing exercise to show that this graphical notation is also justified, and that these mates give the data of a {\em left}-adjunction  $\sigma: Z_0(pt^+) \to Z_1(pt^+)$. Thus the data of any transformation of topological field theories automatically gives rise to an ambidextrous adjunction, just from requiring naturality with respect to the cusp morphisms and by taking mates.

\begin{figure}[ht]
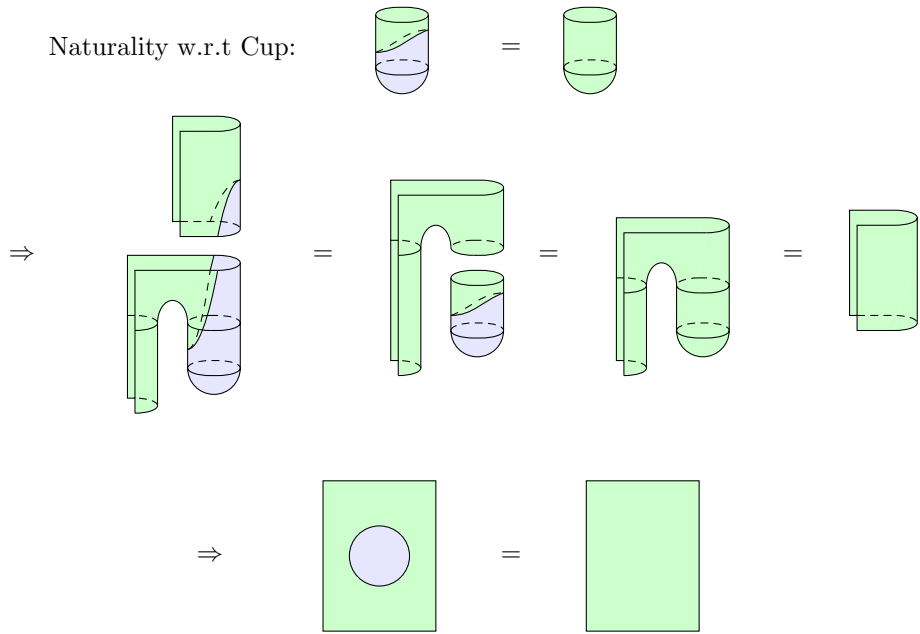

\begin{center}


\caption{Naturality with respect to Cup Implies a Relation}
\label{FigNatWRTCupRel}
\end{center}
\end{figure}

\begin{figure}[ht]
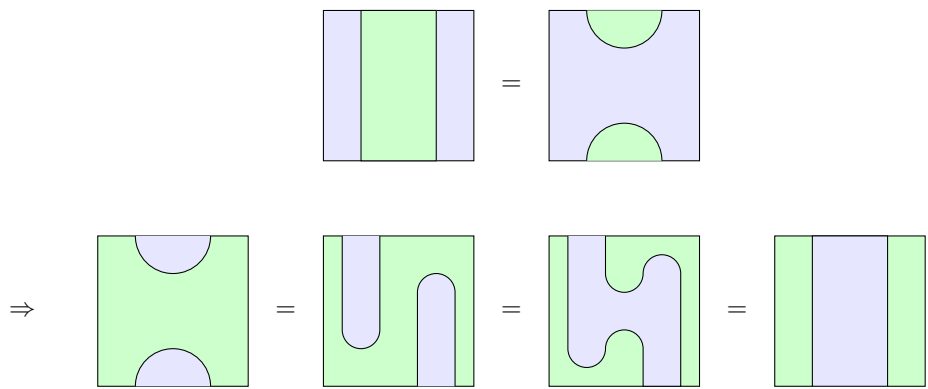

\begin{center}
%
\caption{Relations from Cusps are not Independent}
\label{FigCuspRelNotIndep}
\end{center}
\end{figure}

Finally we must require naturality with respect to the 2D Morse generators. Naturality with respect to these generators will force the mates to be the inverses of $\eta$ and $\varepsilon$ so that we actually have an adjoint equivalence between $Z_0(pt^+)$ and $Z_1(pt^+)$. In Figure~\ref{FigNatWRTCupRel} we show what naturality with respect to the  cup generator means in terms of the graphical notation. 
 We also show how this can equivalently be expressed as a relation among $\eta$, $\varepsilon$, $\eta^*$, and $\varepsilon^*$. Similarly we may translate the condition of being natural with respect each of the 2D Morse generators into relations among  $\eta$, $\varepsilon$, $\eta^*$, and $\varepsilon^*$. These are summarized in Table~\ref{FigNatWRTCupRel}. The two relations coming from the cusps are not independent. In fact they are equivalent, as the graphical calculation in Figure~\ref{FigCuspRelNotIndep} shows.

\begin{table}[ht]
\begin{center}
\begin{tabular}{|m{0.75 in} |m{2.0 in}|} \hline
Generator & Naturality Relation \\ \hline
%%%% Cup
\begin{center}
\begin{tikzpicture}
\begin{scope}[]
\draw (.1, 3) arc (90: -90: 0.3cm and 0.1cm) -- (0, 2.8);
\draw (0, 3) arc (90: 270: 0.3cm and 0.1cm) ;
\draw (0, 3) -- (.1, 3);
\draw (-.3, 2.9) arc (180: 360: 0.35cm and 0.4cm);
\end{scope}
\end{tikzpicture}
\end{center}
	& 
\begin{center}	
\begin{tikzpicture}
\filldraw [fill = green!20] (2, -2) rectangle +(1.5, 2);
\filldraw [fill = blue!10] (2.75, -1) circle (0.4cm);
\node at (4.5, -1) {=};
\filldraw [fill = green!20] (5.5, -2) rectangle +(1.5, 2);
\end{tikzpicture}
\end{center}
\\ \hline
%%% Saddle Curve Up
\begin{center}
\begin{tikzpicture}
	\draw (0, 3) arc (90: -90: 0.3cm and 0.1cm) -- (0, 2);
\draw (1, 3) arc (90: 270: 0.3cm and 0.1cm) -- (1.1, 2.8) -- (1.1, 2) -- (0, 2);
\draw (1,3) -- (1, 2.8); \draw [densely dashed] (1, 2.8) -- (1, 2.2) -- (0,2.2);
\draw (.3, 2.9) arc (180: 360: 0.2cm and 0.3cm);
\draw (0, 3) -- (-.1, 3) -- (-.1, 2.2) -- (0, 2.2); 
\end{tikzpicture}
\end{center}
	& 
\begin{center}	
\begin{tikzpicture}
\filldraw [fill = blue!10] (2, -2) rectangle +(1.5, 2);
\filldraw [fill = green!20] (2.25, 0) arc (180: 360: 0.5cm) -- (2.25, 0);
\filldraw [fill = green!20] (2.25, -2) arc (180: 0: 0.5cm) -- (2.25, -2);
%(2.75, -1) circle (0.4cm);
\node at (4.5, -1) {=};
\filldraw [fill = blue!10] (5.5, -2) rectangle +(1.5, 2);
\filldraw [fill = green!20] (5.75, -2) rectangle +(1, 2);
\end{tikzpicture}
\end{center}	
	\\
\hline
\end{tabular}
\end{center}
\caption{Relations Imposed on Transformations by Naturality with respect to Morse Generators}
\label{TableRelOnTransNatWRTMorseGen}
\end{table}

\section{Extended TFTs in algebras, bimodules, and maps}

In this final section we will use our presentation results to give a concrete classification of extend oriented and unoriented topological field theories in the symmetric monoidal bicategory $\alg^2$ of $k$-algebras, bimodules, and maps of bimodules. 
The classification results of Theorems \ref{TheoremUnorientedClassification} and \ref{OrientedBordismClassification} state that $\bord_2$ and $\bord_2^\text{or}$ are equivalent as a symmetric monoidal bicategories to particular symmetric monoidal bicategories, $\sB$ and $\sB^\text{or}$, defined in terms of generators and relations. If $\sC$ is a target symmetric monoidal bicategory, then the bicategory $\symbicat( \bord_2, \sC)$ is equivalent to the bicategory $\symbicat(\sB, \sC)$, and similarly in the oriented case.  Using Theorem~\ref{thm:cofibrancythm} and Proposition~\ref{RelationsForSMBicatsProp}, we may characterize these bicategories completely. 

There is a symmetric monoidal homomorphism $\bord^\text{or}_2 \to \bord_2$, which forgets the orientation. Thus any unoriented topological field theory also gives rise to an oriented field theory. It make sense, then, to begin with the oriented classification and determine what additional structure is necessary to extend an oriented theory to an unoriented theory. 

In the oriented case we will find that extended 2-dimensional topological field theories correspond to the classical notion of {\em separable symmetric Frobenius algebras}. In the next few subsections we review these classical notions and introduce the 2-groupoid $\frob$, whose objects are separable symmetric Frobenius algebras and whose morphisms are Morita equivalences which preserve these structures. We will show:

\begin{theorem}\label{ThmorienedCalssifcationInAlg}
The bicategory of 2-dimensional {\em oriented} extended topological field theories with values in $\alg^2$ is equivalent to the bicategory $\frob$.
\end{theorem}

This has many consequences. The ordinary 1-categorical bordism category of closed oriented 1-manifolds and oriented 2-dimensional cobordisms sits inside $\bord_d^\textrm{or}$ as the endomorphism category of the empty zero manifold. The endomorphism category of the unit algebra in $\alg^2$ is the category of $k$-vector spaces (or $k$-modules when $k$ is a commutative ring). Thus any extended tqft in $\alg^2$ gives rise to an classical 2-dimensional non-extened tqft. 
 
Recall the classification of non-extended 2-dimensional tqfts in vector spaces states that they are the same as {\em commutative} Forbenius algebras (Theorem~\ref{FolkTheorem}). In terms of the above extended classification the commutative Frobenius algebra arises as the center $\cZ(A)$ of the separable algebra $A$. However the center of a separable algebra is again separable and so this places strong constraints on which topological field theories extend to points. (Separability is equivalent to semi-simplicity over perfect fields. In general it is a stronger property, see Prop.~\ref{prop:characterizesep}).

Moreover we also obtain information about extended topological field theories in higher dimensions. 

\begin{corollary}
	If $Z: \bord^\textrm{or}_d \to \alg^2$ is any bicategorical $d$-dimensional topological field theory and $M$ is any $(d-2)$-manifold, then the algebra $Z(M)$ is separable. In particular over a field $Z(M)$ is semi-simple. 
\end{corollary}
 
\begin{proof}
	There is a symmetric monoidal homorphism $(-) \times M: \bord_2 \to \bord_d$, which sends $pt$ to $M$. Thus $Z$ induces a 2-dimensional topological field theory $\hat Z$, in which $\hat Z(pt) = Z( M)$. 
\end{proof} 

In the unoriented case the classification is less familiar. The algebraic structure we obtain is similar to an algebraic $*$-structure. This is expected, as $*$-algebras are known to give rise to unoriented topological field theories via lattice field theory methods \cite{Snyder:aa}. These methods are expected to give extended field theories as well. 

 However $*$-structures are not Morita invaraint (Example~\ref{ex:stellarnotstar}), and so cannot precisely correspond to extended topological field theories. The correct notion is a sort of Morita invariant version of a $*$-structure which we term a {\em stellar structure}. Briefly, while a $*$-structure consists of a particular algebra isomorphism $A \cong A^\textrm{op}$, a stellar structure replaces this with a Morita equivalence. This notion is developed in Section~\ref{sec:unorientedandstellar}. As before there is a 2-category of stellar separable Frobenius algebras and their Morita equivalences  $\frob^*$. We will prove:

\begin{theorem} \label{ThmUnorientedInAlgs}
The bicategory of 2-dimensional unoriented extended topological field theories with values in $\alg^2$ is equivalent to the bicategory $\frob^*$. 
\end{theorem}

\subsection{Algebras and Duality}

Throughout, let $k$ denote a fixed commutative ring. All algebras considered here will be algebras over $k$.

\begin{definition}
	An object, $P$, of an abelian category is {\em projective} if it satisfies the following universal lifting property: Given a surjection $p:A \twoheadrightarrow B$, and a map $f: P \to B$, there is at least one map $h: P \to A$ such that $f = p \circ h$.
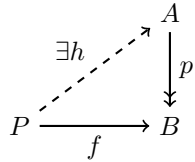
	\begin{center}
\begin{tikzpicture}[thick]
	%\node (LT) at (0,1.5) 	{$$ };
	\node (LB) at (0,0) 	{$P$};
	\node (RT) at (2,1.5) 	{$A$};
	\node (RB) at (2,0)	{$B$};
	%\draw [->] (LT) --  node [left] {$$} (LB);
	\draw [->, dashed] (LB) -- node [above left] {$\exists h$} (RT);
	\draw [->>] (RT) -- node [right] {$p$} (RB);
	\draw [->] (LB) -- node [below] {$f$} (RB);
\end{tikzpicture}
\end{center}
In particular this defines projective modules in the categories $\mod_R$ and ${}_R \mod$, where $R$ is a (possibly non-commutative) algebra. 
An $R$-module (left or right) is {\em finitely generated} if it admits a surjective map from a finite rank free module. 
\end{definition}

\begin{example}
Free modules are projective. 
\end{example}

\begin{lemma} \label{LmaFactorIdenitythroughFree}
Let $M \in \mod_R$ be an $R$-module. If the identity $M \to M$ factors as a map of $R$-modules, 
\begin{equation*}
	M \to (R)^{\times n} \to M
\end{equation*}
through a free module of rank $n < \infty$, then $M$ is a finitely generated projective module. 
\end{lemma}

\begin{proof}
$M$ is finitely generated by definition. Let $A \twoheadrightarrow B$ be a surjection and let $f: M \to B$ be a map of $R$-modules. Since free modules are projective, we have a lift,
	
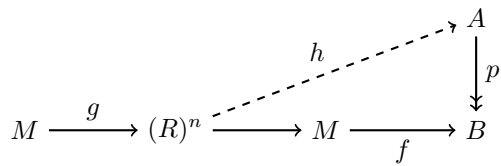
\begin{center}
\begin{tikzpicture}[thick]
	\node (LLLB) at (-4,0) 	{$M$ };
	\node (LLB) at (-2,0) 	{$(R)^n$ };
	\node (LB) at (0,0) 	{$M$};
	\node (RT) at (2,1.5) 	{$A$};
	\node (RB) at (2,0)	{$B$};
	\draw [->] (LLLB) --  node [above] {$g$} (LLB);
	\draw [->] (LLB) --  node [above] {$$} (LB);	
	\draw [->, dashed] (LLB) -- node [above left] {$h$} (RT);
	\draw [->>] (RT) -- node [right] {$p$} (RB);
	\draw [->] (LB) -- node [below] {$f$} (RB);
\end{tikzpicture}
\end{center}
The map $h \circ g: M \to A$ gives the necessary lift of $f$. 
\end{proof}

The category ${}_R \mod$ has a natural action by the symmetric monoidal category ${}_k \mod$, given by the tensor product, $({}_R M, {}_kV) \mapsto M \otimes_k V$. Let $R$ and $S$ be algebras, and let $M$ be an $S$-$R$-bimodule. Then $F = M \otimes_R (-): {}_R\mod \to {}_S\mod$ defines a functor which is linear with respect to the ${}_k \mod$-action in the following sense: For each $V \in {}_k \mod$ we have a natural isomorphism of functors $\eta_V: F(N \otimes_k V) \cong F(N) \otimes_k V$, which is subject to the coherence constraint $\eta_k = id$ and 
\begin{equation*}
	\eta_{V \otimes_k W} = \eta_V \circ \eta_W: F(N \otimes_k V \otimes_k W) \cong F(N) \otimes_k V \otimes_k W.
\end{equation*}

%\CSP{Do I need filtered colimits or something?}
\begin{lemma} \label{lem:whenfunctoristensor}
Let $R$ and $S$ be two algebras and let $F:  {}_R\mod \to {}_S\mod$ be a functor which
\begin{enumerate}
\item Is ${}_k\mod$-linear in the above sense, and
\item Preserves coequalizers. 
\end{enumerate}
Then $F \cong M \otimes_R (-)$ for some $S$-$R$-bimodule $M$. In particular $F$ is additive.  
\end{lemma}

\begin{proof}
	First notice that since $F$ is ${}_k\mod$-linear it automatically preserves repeated direct sums, for example
	\begin{equation*}
		F( N \oplus N) \cong F( k^2 \otimes N) \cong k^2 \otimes F(N) \cong F(N) \oplus F(N).
	\end{equation*} 
	Moreover it also preserves diagonal and codiagonal maps, and hence it follows that on hom spaces $F$ is $k$-linear. Consequently $F(R) \in {}_S\mod$ inherits a right $R$-action, and is thus an $S$-$R$-bimodule.

	Let $N \in {}_R\mod$ be a left $R$-module. We have the following coequalizer sequence in ${}_R\mod$:
	\begin{equation*}
		{}_R R \otimes_k R \otimes_k N \rightrightarrows {}_R \otimes_k N \to {}_R N,
\end{equation*}
where the two left-most maps are given by the action of $R$ on itself and the action of $R$ on $N$, respectively. Since $F$ preserves coequalizers, we have the following is a coequalizer sequence in ${}_S\mod$:
 	\begin{equation*}
		F( R \otimes_k R \otimes_k N) \rightrightarrows F( R \otimes_k N) \to  F(N).
\end{equation*}
Since $F$ is ${}_k\mod$-linear, this is the same as the coequalizer,
 	\begin{equation*}
		{}_S F( R) \otimes_k R \otimes_k N \rightrightarrows {}_S F(R) \otimes_k N \to  {}_S F(N).
\end{equation*}
Thus ${}_S F(N) \cong {}_S F(R) \otimes_R N$ as $S$-modules, and so $F \cong M \otimes_R (-)$ for the $S$-$R$-bimodule $M = F(R)$. 
\end{proof}

Similar results hold for right modules, as well. For our purposes, the most important application of this is when $F$ is a left-adjoint, ${}_k\mod$-linear functor. In this case the above theorem implies that $F$ comes from a bimodule. 
\index{adjunction!for bimodules}
 
\begin{remark}
	In Lemma~\ref{lem:whenfunctoristensor} we can replace condition (1) with the statement that $F$ is additive and that it commutes with filtered colimits. 
\end{remark}

%\CSP{make lemma an iff}
\begin{lemma}
Let $M$ be an $S$-$R$-bimodule.  The functor ${}_S M \otimes_R (-)$ is a right-adjoint if and only if $M_R$ is a finitely generated projective $R$-module. Similarly the functor $(-) \otimes_S M_R$ is a right-adjoint if and only if ${}_SM$ is a finitely generated projective $S$-module. 
\end{lemma}

\begin{proof}
We prove this for left-modules. The case of right-modules is similar. Let $F = {}_S M \otimes_R (-)$. First suppose that $F$ is a right-adjoint, then there exists a functor $G$, together with natural transformations $\eta: id \to GF$ and $\varepsilon: FG \to id$, realizing the adjunction. Equivalently, we have a natural isomorphism,
\begin{equation*}
	\hom_S(G(A), B) \cong \hom_R( A, F(B))
\end{equation*}
 for all $A \in {}_R\mod$ and $B \in {}_S\mod$. 

$F$ is ${}_k\mod$-linear since it arises from a bimodule. In fact $G$ is also ${}_k\mod$-linear, which can be seen as follows. For any $V \in {}_k\mod$,  $A \in {}_R\mod$, and $B \in {}_S\mod$, we have,
\begin{align*}
	\hom_S(G(A \otimes_k V), B) & \cong \hom_R( A\otimes_k V, F(B) ) \\
		& \cong \hom_k( V, \hom_R( A, F(B) )) \\
		& \cong \hom_k( V, \hom_S( G(A), B )) \\
		& \cong \hom_S( G(A) \otimes_k V, B ).
\end{align*}
Since this is true for all $B$, by the Yoneda lemma $G(A \otimes_k V) \cong G(A) \otimes_k V$, so that $G$ is ${}_k\mod$-linear. The previous lemma now implies that $G = {}_R N \otimes_S (-)$ for some $R$-$S$-bimodule, $N$. 

The counit and unit of the adjunction satisfy the equations,
\begin{align*}
	id_F & = \varepsilon F \circ F \eta, \\
	id_G & = G\varepsilon \circ \eta G.
\end{align*}
Writing these out in terms of the bimodules ${}_SM_R$ and ${}_RN_S$ we see that $\varepsilon$ and $\eta$ become bimodule maps,
\begin{align*}
	\varepsilon: & {}_RN \otimes_S M_R \to {}_R R_R, \\
	\eta: & {}_S S_S \to {}_S M \otimes_R N_S.
\end{align*}
Since $S$ is cyclic as an $S\otimes S^\text{op}$-module, $\eta$ is determined by the single element $\eta(1) = \sum_{i \in I} x_i \otimes_R y_i \in M \otimes_RN$. As an aside,  $\eta(1)$ is {\em central} in the sense that for all $s \in S$, we have $\sum (s x_i) \otimes_R y_i = \sum x_i \otimes_R (y_i s)$. The adjunction equation may now be written as,
\begin{align*}
	\sum_i \varepsilon (a \otimes_S x_i) y_i &= a, \\
	\sum_i x_i \varepsilon( y_i \otimes_S b) & = b,
\end{align*}
for all $a \in N$ and $b \in M$. 

The identity map on $M$ may now be factored through a finite rank free $R$-module, as follows. Define maps:
\begin{align*}
M_R &\to (R_R)^I &   (R_R)^I & \to M_R \\
m &\mapsto (\varepsilon(y_i, m))_{i \in I} & (r_i)_{i \in I} & \mapsto \sum_i x_i r_i
\end{align*}
The adjunction equations imply that the composition is the identity on $M_R$. Since $I$ is finite, Lemma \ref{LmaFactorIdenitythroughFree} ensures that $M_R$ is a finitely generated projective $R$-module.

Conversely suppose that $M_R$ is a finitely generated projective $R$-module. Then there exists some finite $n < \infty$ and a factorization of the identity map:
\begin{equation*}
	M \stackrel{(f_i)}{\longrightarrow} R^n \stackrel{p}{\longrightarrow} M.
\end{equation*}  
Let ${}_R N_S = \hom_R(M_R, R_R)$ be the $R$-linear hom of $M$. The evaluation map $\hom_R(M, R) \otimes_k M \to R$ naturally factors as a bimodule map:
\begin{equation*}
	\varepsilon: {}_R(N \otimes_S M)_R = \hom_R(M,R) \otimes_S M \to {}_R R_R.
\end{equation*}
This is the counit of the desired adjunction. The unit is constructed as follows. Let $\{e_i\}$ denote the standard $R$-basis of $R^n$. Then the unit $\eta: {}_S S_S \to M \otimes_R \hom(M,R)$ is given by the formula
\begin{equation*}
	\eta(1) = \sum_i  p(e_i) \otimes_R f_i
\end{equation*}
and extending $S$-linearly. 
\end{proof}

In the above we call the bimodule ${}_R N_S$ the  {\em left dual} of the bimodule ${}_S M_R$. Notice that the analogous argument (using the remaining adjunction equation) shows that $N$ is a finitely generated projective $S$-module.

%\begin{proposition}
%	Let ${}_AM_B$ be an $A$-$B$-bimodule such that .... 
	%	Then $M$ is finitely generated and projective  as an $A$-module. 
%\end{proposition}

\subsection{Symmetric Frobenius Algebras}

%In what follows $R$ will be a fixed commutative ring. All our algebras are $R$-algebras. All tensor products are taken over $R$.
\begin{definition} \label{DefFrobAlg}
 A {\em Frobenius algebra} 
\index{algebra!Frobenius} 
\index{Frobenius algebra|see{algebra}} 
\index{Frobenius algebra} 
is triple $(A, \lambda, e)$ where $A$ is a $k$-algebra, $\lambda: A \to k$ is a $k$-linear functional, $e  = \sum_i x_i \otimes y_i \in A \otimes_k A$ is {\em $A$-central}, 
\index{$A$-central element} 
i.e., for all $w \in A$, $\sum_i (w \cdot x_i) \otimes y_i = \sum_i x_i \otimes (y_i \cdot w)$, and such that the following Frobenius normalization condition is satisfied,
 		\begin{equation*}
			\sum_i \lambda(x_i) y_i =\sum_i x_i \lambda(y_i) = 1_A.
		\end{equation*}
\end{definition}
\begin{definition} \label{DefSymFrobAlg}
 A Frobenius algebra $(A, \lambda, e)$ is {\em symmetric}
\index{algebra!symmetric Frobenius}
 if $\lambda$ is trace-like, i.e., $\lambda(xy) = \lambda (yx)$ for all $x,y \in A$. 
\end{definition}

\begin{definition}
	Let $V$ be a $k$-module. A bilinear form, $b: V \otimes V \to k$ is {\em non-degenerate} if
	the following conditions are satisfied:
	\begin{enumerate}
		\item $b(x,y) = 0$ for all $x \in V$ implies $y = 0$, 
		\item $b(x,y) = 0$ for all $y \in V$ implies $x = 0$.
	\end{enumerate}
	
The {\em $k$-dual} of $A$ \index{dual module} is defined to be $\hat A := \hom_k(A, k)$. It is naturally an $A$-$A$-bimodule with actions,
\begin{equation*}
	a \cdot f \cdot b \mapsto ( x \mapsto f(bxa)).
\end{equation*}
\end{definition}

There are several equivalent formulations of the notion of a Frobenius algebra.
\begin{proposition}
	Given a fixed $k$-algebra $A$, the following structure/property combinations are equivalent:
	\begin{enumerate}
	\item $(\lambda, e)$ such that $(A, \lambda, e)$ forms a Frobenius algebra. 
	\item $(b,e)$ where $e \in A \otimes A$ is $A$-central, $b: A \otimes A \to k$ is a non-degenerate $k$-bilinear form such that $b(xy,z) = b(x, yz)$ and such that $(A, \lambda, e)$ is a Frobenius algebra, where $\lambda(x) = b(x,1_A)$.
	\item$(b,e)$, where $b: A \otimes A \to k$ is a non-degenerate bilinear form satisfying $b(xy,z) = b(x, yz)$ and $e: k \to A \otimes A$ is an $A$-central element, such that  the {\em snake relations} are satisfied:
	\begin{equation*}
		(id_A \otimes b) \circ (e \otimes id_A) = ( b \otimes id_A) \circ (Id_A \otimes e) = id_A. 
	\end{equation*}
	\item $b: A \otimes A \to k$, a non-degenerate bilinear form satisfying $b(xy,z) = b(x, yz)$ together with the property that $A$ is a finitely generated projective $k$-module.
	\item $\lambda: A \to k$ such that $\ker \lambda$ contains no nontrivial left  ideals, together with the property that $A$ is a finitely generated projective $k$-module. Left ideals may be replaced by right ideals to yield an equivalent characterization. 
	\item $b: A \otimes A \to k$, a non-degenerate bilinear form satisfying $b(xy,z) = b(x, yz)$ and such that $b$ induces an isomorphism of right $A$-modules $A_A \cong \hat A_A$, together with the property that $A$ is a finitely generated projective $k$-module. Again, left-modules may be used instead. (A $k$-module $A$ which is isomorphic to $\hat A$ is called {\em reflexive} 
\index{reflexive module}).
	\item  An isomorphism of right $A$-modules $A_A \cong \hat A_A$, (or an isomorphism of left $A$-modules ${}_AA \cong {}_A\hat A$),  together with the property that $A$ is a finitely generated projective $k$-module
	\item $(\Delta, \lambda)$ where $\Delta: A \to A \otimes_k A$ is a comultiplication with counit $\lambda: A \to k$ such that the {\em Frobenius equations} are satsified
	\begin{equation*}
		\Delta \circ m  = (m \otimes id_A) \circ (id_A \otimes \Delta) =  (id_A \otimes m) \circ (\Delta \otimes id_A) 
	\end{equation*}
	where $m: A \otimes_k A \to A$ is the multiplication map. 
	\item $(\Delta, \lambda)$ where $\Delta: A \to A \otimes_k A$ is a comultiplication with counit $\lambda: A \to k$ such that $\Delta: {}_A A_A \to {}_A A \otimes A_A$ is a bimodule map. 
	\end{enumerate}
\end{proposition}

\begin{proof}
The equivalences $ (4) \Leftrightarrow (5) \Leftrightarrow (6) \Leftrightarrow (7) $ are fairly standard and are left as an exercise for the reader. It is equally clear that $(8) \Leftrightarrow (9)$. $(2) \Leftrightarrow (3)$, once one realizes that the snake equation and the Frobenius equation are equivalent in this context. By construction $(2) \Rightarrow (1)$ and $(1) \Rightarrow (2)$ as follows: letting $b = \lambda \circ m$ yields a bilinear pairing compatible with the multiplication in the necessary way. We need only show that $b$ is non-degenerate. Suppose that there exists an element $z \in A$ such that $0 = b(w, z) =  \lambda (w\cdot z) $ for all $w \in A$. Then,
\begin{equation*}
z = z \cdot 1_A = z \cdot \sum_i x_i \lambda(y_i) = \sum_i (z\cdot x_i) \lambda(y_i) = \sum_i x_i \lambda(y_i \cdot z)  = 0
\end{equation*}
Thus $b$ is non-degerate on the right. A similar calculation shows $b$ is non-degenerate on the left. 
$(9) \Rightarrow (1)$ by taking $e = \Delta(1_A)$. The fact that $\Delta$ is a bimodule map implies that $e$ is $A$-central and the compatibility of comultiplication and counit yield the Frobenius equations. Finally, $(2) \Rightarrow (9)$ by taking $\Delta(z) := \sum_ix_i \otimes (y_i \cdot z)$. This is a bimodule map since $e$ is $A$-central, and is a comultiplication for which $\lambda$ is the counit by virtue of the Frobenius equations.  

Thus it remains to prove that the group $(4)$ through $(7)$ is equivalent to the group $(1)$, $(2)$, $(3)$, $(8)$, $(9)$. 
(3) implies that $A$ is a finitely generated projective module and hence $(3) \Rightarrow (4)$. 
 The reverse implication follows by choosing a basis $a_i \in A$, i.e., the images of the standard basis via some surjective map $(k)^n \to A$. Such a basis exists since  $A$ is finitely generated. Moreover, since $A$ is finitely generated and projective, and more specifically $b$ induces an isomorphism $A \cong \hat A$, we may define elements $\hat a_i \in A$, by the following equation,
 \begin{equation*}
	b(\hat a_i, a_j) = \delta_{ij}.
\end{equation*}
Thus we have the following identities, $\sum a_i b(\hat a_i, x) = x = \sum b(x, a_j) \hat a_j$.
Define the element $e = a_i \otimes \hat a_i \in A \otimes A$. By construction, $e$ satisfies the Frobenius normalization condition. It is straightforward to check that $e$ is $A$-central.
\end{proof}

The conditions in (6) and (7) that $A$ be a finitely generated \emph{projective} $k$-module are necessary as the following example shows. This example was pointed out to me by D. Eisenbud. 

\begin{example}
Let $L$ be a field and consider $k = L[x,y,z]$. Let $f: k^3 \to k$ be defined by the matrix $(x,y,z)$. There is a short exact sequence,
\begin{equation*}
0 \to K \to k^3 \stackrel{f}{\to} k
\end{equation*}
Now $k$ has projective dimension 3, so $K$ cannot be a projective module. This is actually part of the larger Koszul complex,
\begin{equation*}
0 \to k \stackrel{g}{ \to} k^3 \to k^3  \stackrel{f}{\to} k 
\end{equation*}
Which allows us split off the following short exact sequence,
\begin{equation*}
0 \to k \stackrel{g}{ \to} k^3 \to K \to 0
\end{equation*}
We apply $Hom_k (-, k)$ to get,
\begin{equation*}
	0 \to Hom_k(K, k) \to Hom_k(k^3, k) \to Hom_k(k,k)
\end{equation*}
Which becomes
\begin{equation*}
0 \to Hom_k(K,k) \to k^3 \stackrel{g^*}{\to} k
\end{equation*}
Here $g^*$ is the matrix $(x, -y, z)$, so this sequence is actually exact on the right as well. We find that $Hom_k(K,k) \cong K$. So $K$ is a reflexive $k$-module, which is finitely generated, but is not projective.  
\index{reflexive module!finitely generated, not projective}
\end{example}

\begin{proposition}
If $(A, \lambda, e)$ a Frobenius algebra, then the following are equivalent:
%the property of being symmetric is equivalent to the following characterizations. 
\index{algebra!symmetric Frobenius}
\begin{enumerate}
\item $\lambda$ is trace-like, i.e., $(A, \lambda, e)$ is a {\em symmetric} Frobenius algebra.
\item The bilinear form $b$ is symmetric, i.e., $b(x,y) = b(y,x)$.  
\item The element $e$ is $A$-bicentral 
\index{$A$-bicentral element}, i.e., for all $w, z \in A$ we have,
\begin{equation*}
	 \sum_i (w \cdot x_i \cdot z) \otimes y_i = \sum_i x_i \otimes (z \cdot y_i \cdot w)
\end{equation*}
\item The isomorphism $\hat A \cong A$ induced by $b$ is an isomorphism ${}_A A_A \cong {}_A\hat A_A$ of bimodules.
\end{enumerate}
\end{proposition}

\begin{proof}
(1) $\Leftrightarrow$ (2) is clear, from the proof of the previous lemma.  (2) $\Rightarrow$ (4), as follows.  The image in $\hat A$ of the element $a$ is the assignment $x \mapsto b(a, x)$ and so we must show that $b(waz, x) = b(a, zxw)$ for all $z,w, x, a \in A$. A direct calculation, assuming $b$ is symmetric, shows:
\begin{align*}
	b(waz, x)  & = b(wa, zx) \\
	 &= b(zx, wa) \\
	 & = b(zxw,a) \\
	 &= b(a, zxw).
\end{align*}

Similarly (3) $\Rightarrow$ (1), as follows. If $e  = \sum x_i \otimes y_i$ is bicentral, then, by the snake relations we have,
\begin{align*}
	\lambda(ab) & = \lambda ( \sum x_i \lambda(y_i a) b) \\
	 &= \sum \lambda(y_i a) \lambda (x_i b) \\
	 & = \sum \lambda(b y_i) \lambda (a x_i) \\
	 &= \lambda ( \sum b \lambda(ax_i) y_i) \\
	 & = \lambda (ba).
\end{align*}

Finally, (4) $\Rightarrow$ (3) as follows. By assumption we have $b(waz, x) = b(a, zxw)$ for all $z,w,x,a \in A$. Thus,
\begin{align*}
	\sum_i (zx_iw) \otimes y_i &= \sum_{i,j} x_j b(y_j, zx_iw) \otimes y_i \\
		&= \sum_{i,j} x_j \otimes b(w y_j z, x_i) y_i \\
		&= \sum_j x_j \otimes (wy_j z),
\end{align*}
so that $e$ is, indeed, bicentral. 
\end{proof}

\subsection{Fully-Dualizable Algebras and Separable Algebras}

\begin{definition}
A $k$-algebra $A$ is {\em separable} algebra 
\index{algebra!separable} 
\index{separable algebra}
if there exists an $A$-central 
\index{$A$-central element} 
element $\tilde e = \sum_i \tilde x_i \otimes \tilde y_i \in A \otimes A$ such that the separability normalization condition is satisfied,
\begin{equation*}
 \sum_i x_i y_i = 1_A.
\end{equation*}
\end{definition} 

\begin{remark}
In this formulation, a Frobenius algebra is an algebra equipped with additional structure, whereas a separable algebra is an algebra satisfying a property. A separable symmetric Frobenius algebra is a Frobenius algebra which happens to be separable. 
\end{remark}

The following results can be found in \cite{DI71}.
Let $A^e = A \otimes_k A^{op}$ be the {\em enveloping algebra} 
\index{algebra!enveloping} 
\index{enveloping algebra}
of $A$. As $A$ is an $A$-$A$-bimodule, it may equivalently be considered a (left) $A^e$-module. The multiplication map $\mu: A^e \to A$ is an $A^e$-module map with kernel $J$, the left ideal of $A^e$ generated by elements of the form $a \otimes 1 - 1 \otimes a$. Note that an element $e \in A \otimes A$ is $A$-central precisely when we have $Je = 0$ viewed as an element in $A^e$.    
\begin{proposition} \label{prop:characterizesep}
For a $k$-algebra $A$ the following properties are  equivalent, 
\begin{enumerate}
\item $A$ is separable, that is there exists $\tilde e \in A^e$ such that $J \tilde e = 0$ and $\mu(\tilde e) = 1$. 
\item $A$ is projective as a left $A^e$-module under the $\mu$-structure. 
\item The exact sequence of left $A^e$-modules, 
\begin{equation*}
	 0 \to J \to A^e \stackrel{\mu}{\to} A \to 0
\end{equation*}
splits.
\end{enumerate}
If in addition, $k$ is a field then we may add,
\begin{enumerate}
\item [4.] $A$ is finite dimensional and {\em classically separable}, 
\index{algebra!classically separable (semi-simple)} 
\index{classically separable}
that is for every field extension $K$ of $k$ we have that $A \otimes_k K$ is semi-simple. 
\end{enumerate}
In particular $A$ is semi-simple in this case. If $k$ is a perfect field (for example a characteristic zero field or a finite field), then any finite dimensional semi-simple algebra is separable. 
\end{proposition}

\begin{definition}[Lurie \cite{Lurie09}]
	A $k$-algbera $A$ is {\em fully-dualizable} if 
\index{algebra!fully-dualizable}
\index{fully-dualizable algebra}
	\begin{enumerate}
\item It is separable,
\item It is (finitely generated and) projective as a $k$-module. 
\end{enumerate}
\end{definition}
Note that every separable $k$-algebra which is projective is automatically finitely generated \cite[Prop.~2.1]{DI71}. In particular condition (2) above is redundant when $k$ is a field.  

The above definition can be extracted from J. Lurie's expository paper \cite{Lurie09}. This second condition is part of being a Frobenius algebra, and thus we see that a separable symmetric Frobenius algebra is the same as a fully-dualizable algebra with a non-degenerate symmetric trace, $\lambda$.

\subsection{Morita Contexts}

We begin with some algebraic preliminaries. Consider $\alg^2_k$. The objects of this bicategory are $k$-algebras, the 1-morphisms are bimodules, and the 2-morphisms are bimodule maps (a.k.a. intertwiners). The composition of 1-morphisms is given by the tensor product of bimodules. There are two possible conventions, but for definiteness let us consider a $B$-$A$-bimodule, ${}_BM_A$, as a morphism from $A$ to $B$. Equivalence in $\alg^2_k$ is called {\em Morita equivalence}, \index{Morita equivalence} and an adjoint equivalence $({}_S M_R, {}_RN_S, \eta, \varepsilon)$ (see Definition \ref{DefAdjointEquiv}) is called a {\em Morita context}, \index{Morita context} see \cite{Lam99}.

Given an algebra homomorphism $h: A \to B$ we can turn it into a $B$-$A$-bimodule, ${}_BBh_A$. This bimodule is ${}_BB$ as a $B$-module. The right $A$-action comes from the right $B$-action and the homomorphism $h$. These bimodules are compatible with the composition of homomorphisms,
\begin{equation*}
	{}_C (Ch_2) \otimes_B (Bh_1)_A \cong {}_C(C h_2 h_1)_A.
\end{equation*}
In fact these isomorphisms can be chosen to be coherent, i.e., they give a homomorphism from the category of $k$-algebras and algebra homomorphisms (viewed as a bicategory with only identity 2-morphisms) to the bicategory $\alg^2_k$.

The tensor product over $k$ gives $\alg^2_k$ the structure of a symmetric monoidal bicategory. In fact, it equips the category of algebras and homomorphisms with the structure of a symmetric monoidal category, and most of the structure of a symmetric monoidal bicategory can be transfered to $\alg^2_k$ via the above homomorphism. The remaining structure deals with the tensor product (over $k$) of bimodules and is easily supplied. We leave the details to the reader. 

Recall from Definitions \ref{DefFrobAlg} and \ref{DefSymFrobAlg} the notion of a symmetric Frobenius algebra $(A, e, \lambda)$. Here $e = \sum x_i \otimes y_i \in A \otimes_kA$ is a bicentral element and $\lambda: A \to k$ is a $k$-linear map, which satisfies the Frobenius normalization condition. A bicentral element $e \in A \otimes A$ is equivalent to a map of bimodules,
\begin{equation*}
	{}_{A_1}A_{A_2} \otimes {}_{A_3}A_{A_4} \rightarrow {}_{A_1}A_{A_4} \otimes {}_{A_3}A_{A_2}
\end{equation*}
where we have labeled the algebras with numbers to keep track of the different actions. This correspondence arises because the source bimodule is cyclic as an $(A \otimes A)$-$(A \otimes A)$-bimodule with cyclic element $1 \otimes 1$. Hence determined by the image $e \in A \otimes A$ of $1 \otimes 1$, moreover it is easy to check that this element is bicentral. Expressed as a bimodule map, as above, it is clear that the bicentral element $e$ can be transfered along Morita contexts. We will show that $\lambda$ may also be transfered along Morita contexts.\footnote{It is essential that we are considering {\em symmetric} Frobenius algebras. In general the property of being a Frobenius algerba is {\em not} Morita invariant, see \cite{Lam99}.}

When $\lambda$ is symmetric, the linear map $\lambda: A \to k$ necessarily factors through the quotient,
\begin{equation*}
	A \to A/[A,A] \to k,
\end{equation*}
where $[A,A]$ denotes the minimal $k$-submodule containing all elements of the form $a \cdot b - b \cdot a$. We have the following lemma:
\begin{lemma} \label{lem:universaltrace}
	Let $A$ and $B$ be $k$-algebras and  $f = ({}_B M_A, {}_AN_B, \eta, \varepsilon)$ a Morita context. Then there is a canonical isomorphism of  $k$-modules:
	\begin{equation*}
		f_*: A/[A,A] \cong B/[B,B].
\end{equation*}
In particular the $k$-linear map $\lambda$ may be transfered along Morita contexts.
\end{lemma}

\begin{proof}
We have the following sequence of canonical isomorphisms:
\begin{align*}
	A/[A,A] & \cong A \otimes_{A \otimes A^\text{op}} A \\
		& \cong (N \otimes_B M) \otimes_{A \otimes A^\text{op}} (N \otimes_B M) \\
		& \cong (M \otimes_A N) \otimes_{B \otimes B^\text{op}} (M \otimes_A N) \\
		& \cong B \otimes_{B \otimes B^\text{op}} B \\
		& \cong B/[B,B].
\end{align*}
\end{proof}

\begin{definition}
Let $(A, e^A, \lambda^A)$ and $(B, e^B, \lambda^B)$ be symmetric Frobenius algebras. Let $f = ({}_B M_A, {}_AN_B, \eta, \varepsilon)$ be a Morita context between $A$ and $B$. Then we say that $f$ is {\em compatible} with the symmetric Frobenius algebra structure if $f_*e^A = e^B$ and $f_* \lambda^A = \lambda^B$.  
\end{definition}

\begin{definition}
Let $\frob$ be the symmetric monoidal bicategory whose objects are separable symmetric Frobenius algebras, whose 1-morphisms are compatible Morita contexts, and whose 2-morphisms are isomorphisms of Morita contexts. 
\end{definition}

Given a bimodule ${}_A M_B$, we can reinterpret the algebra actions and obtain a new bimodule ${}_{B^\text{op}} \underline{M}_{A^\text{op}}$.  $\underline{M}$ has the same underlying $k$-module structure as $M$. If $f = ({}_B M_A, {}_AN_B, \eta, \varepsilon)$ is a Morita context, then $\underline{f} = ({}_{A^\text{op}} \underline{M}_{B^\text{op}}, {}_{A^\text{op}}\underline{N}_{B^\text{op}}, \eta, \varepsilon)$ is a Morita context between $B^\text{op}$ and $A^\text{op}$.

\subsection[Extened TFTs = Sep. Sym. Frob. Algebras]{Examples: Extended TFTs = Separable Symmetric Frobenius Algebras} \label{SectExamplesandApps}

We will now prove Theorem~\ref{ThmorienedCalssifcationInAlg}, which gives an explicit classification of 2-dimensional oriented extended topological field theories in the bicategory $\alg^2$. In fact we already implicitly have an explicit classification. Theorem~\ref{OrientedBordismClassification} gives an explicit presentation of the oriended bordism bicategory and Theorem~\ref{thm:cofibrancythm} (see also Prop~\ref{RelationsForSMBicatsProp}) give us a precise list of data which is equivalent to specifying the topological field theory.

%By  \ref{CharacterizingFucntiorsFromFree} and \ref{RelationsForSMBicatsProp} the bicategory $\symbicat(\sB^\text{or}, \alg^2)$ is equivalent to the bicategory of generating data in $\alg^2$ which satisfy the relations of Theorem~\ref{OrientedBordismClassification}.

For example, in this particular case there are two generating objects corresponding to the two different oriented points. Their images in $\alg^2$ will correspond to two algebras $A = Z(pt^+)$ and $B = Z(pt^-)$.
Similarly there are two generating 1-morphisms which give rise to two bimodules
\begin{align*}
	Z(\tikz[baseline = -0.7cm]{
		\draw [-] (-6.6, -.5) -- (-6.5, -.5) arc (90: -90: 0.3cm and 0.1cm);
		\node at (-6.8, -.5) {\tiny +};
		\node at (-6.7, -.7) {-};
		}) & = {}_{A \otimes B}M \\
	Z_1( \tikz[baseline = -0.7cm]{
		\draw [-] (-5.5, -.5) arc (90: 270: 0.3cm and 0.1cm) -- (-5.4, -.7);
		\node at (-5.3, -.5) {\tiny +};
		\node at (-5.2, -.7) {-};
		}) &= N_{A \otimes B}
\end{align*}
The 2-morphism generators in the oriented case fall into two categories: the cusp generators and the 2D Morse generators, and there is a short list of relations these must satisfy: the 2D Morse relations, the swallowtail relations, the cusp flip relations, and the cusp inversion relations, as depicted (from left-to-right, top-to-bottom) in Figure~\ref{fig:OrientedGenAndRel}.

In the remainder of this section we will aim to show how this collection of data relates to the notion of separable symmetric Frobenius algebra. We will go through this list of data and re-express it in a simpler equivalent form. 
The first observation is that we can `flip' some of the algebra actions around. This makes some of the defining structures of the topological field theory more transparent. Thus the above two bimodules give rise to bimodules ${}_A \tilde{M}_{B^\text{op}}$ and ${}_{B^\text{op}} \tilde{N}_A$. 

The cusp generators are by far the easiest to recognize. While they naturally arise as morphisms involving $M$ and $N$, we may equivalently express them in terms of  $\tilde{M}$ and $\tilde{N}$. The four cusp generators then give rise to four bimodule maps,
\begin{align*}
	f_1: & \; {}_A \tilde{M} \otimes_{B^\text{op}} \tilde{N}_A \to {}_A A_A \\
	f_2: & \;  {}_A A_A \to {}_A \tilde{M} \otimes_{B^\text{op}} \tilde{N}_A  \\
		f_3: & \; {}_{B^\text{op}} \tilde{N} \otimes_{A} \tilde{M}_{B^\text{op}} \to {}_{B^\text{op}} {B^\text{op}}_{B^\text{op}} \\
	f_4: & \;  {}_{B^\text{op}} {B^\text{op}}_{B^\text{op}} \to {}_{B^\text{op}} \tilde{N} \otimes_{A} \tilde{M}_{B^\text{op}} 
\end{align*}
The cusp inversion relations imply that $f_1$ and $f_2$ are inverses and $f_3$ and $f_4$ are inverses. The swallowtail relation implies that $(\tilde M, \tilde N, f_1, f_4)$ forms an adjoint equivalence, i.e., a Morita context.  Thus this portion of the topological field theory is equivalent to specifying a pair of algebras and a Morita context between them. 
The remaining structure is what is assigned to the cup, the cap, and the two saddles, and it is subject to the 2D Morse relations and the cusp flip relations. We need to identify the remaining structures in more classical terms.

%By  \ref{CharacterizingFucntiorsFromFree} and \ref{RelationsForSMBicatsProp} the bicategory $\symbicat(\sB, \sC)$ is equivalent to the bicategory of generating data in $\sC$ which satisfy the relations. Thus a symmetric monoidal homomorphism $E: \bord_2 \to \alg^2$ is determined up to equivalence  by giving its value on the generating objects, 1-morphism, and 2-morphisms.

%The only generating object is the point, and so we are given a single algebra $E(pt) = A$. The generating 1-morphisms provide two bimodules ${}_{A \otimes A} M_R$ and ${}_R N_{A \otimes A}$, which may equivalently be viewed as bimodules ${}_A M_{A^{op}}$ and  ${}_{A^{op}}N_{A}$. Under this correspondence the cusp generators provide the morphisms $\eta$ and $\varepsilon$. These must satisfy the cusp inversion relations and swallowtail relation which show that $({}_A M_{A^{op}}, {}_{A^{op}}N_{A}, \eta: N \otimes_A M \cong A^{op}, \varepsilon: M \otimes_{A^{op}} N \cong A)$ forms a Morita equivalence. The symmetry generators and relations provide this with the structure of a symmetric Morita equivalence. 

%It remains to show that the remaining generators and relations are equivalent to providing $A$ with the structure of a separable Frobenius algebra. We prove this in several lemmas. First define $V = A/[A,A]$. Using the Morita equivalence, we have a canonical isomorphism $E(S^1) \cong V$. 

Let us begin with the data associated to the saddle bordisms. The value of an oriented topological field theory on these bordisms gives us two bimodule maps:
\begin{align*}
	\alpha: {}_{A \otimes B}M \otimes N_{A \otimes B} &\to {}_{A \otimes B} (A \otimes B)_{A \otimes B} \\
	\alpha': {}_{A \otimes B}(A \otimes B)_{A \otimes B} &\to {}_{A \otimes B} M \otimes N_{A \otimes B}.
\end{align*}
We can use the previous structure (the Morita context between $A$ and $B^\textrm{op}$) to re-express these as maps of bimodules which just use the algebra $A$. First we flip the algebras around to obtain maps:
 \begin{align*}
 	\tilde{\alpha}: {}_{A \otimes B^\textrm{op}}(\tilde{M} \otimes \tilde{N})_{A \otimes B^\textrm{op}} &\to {}_{A \otimes B^\textrm{op}} (A \otimes B)_{A \otimes B^\textrm{op}} \\
 	\tilde{\alpha}': {}_{A \otimes B^\textrm{op}}(A \otimes B)_{A \otimes B^\textrm{op}} &\to {}_{A \otimes B^\textrm{op}} (\tilde{M} \otimes \tilde{N})_{A \otimes B^\textrm{op}}.
 \end{align*}
(be careful in making these identifications: now, for example, the left $B^\textrm{op}$ action on the source of the map $\tilde{\alpha}$ is via $\tilde{N}$ and not $\tilde{M}$). 

Next we use the Morita context to replace the $B^\textrm{op}$-module structures with $A$-module structures. Specifically we tensor by ${}_A\tilde{M}_{B^\textrm{op}}$ on the left and ${}_{B^\textrm{op}} \tilde{B}_A$ on the right. Finally, we can use the isomorphisms $f_1$ and $f_2$ to replace terms of the form  ${}_A (\tilde{M} \otimes_{B^\textrm{op}}\tilde{N})_A$ with ${}_AA_A$. Doing this with both the source and the target yields a pair of bimodule maps:
\begin{equation*}
	\beta, \beta': {}_{A_1}A_{A_2} \otimes {}_{A_3}A_{A_4} \rightarrow {}_{A_1}A_{A_4} \otimes {}_{A_3}A_{A_2}.
\end{equation*}
Here we have labeled the various algebras with numbers to keep track of where each algebra acts. Since each of these process is invertible (given the Morita context) we have simply exchanged the data of the saddle maps $\alpha$ and $\alpha'$ for the equivalent bimodule maps $\beta$ and $\beta'$. 

We considered this kind of bimodule map before in both our treatment of separable algebras and of Frobenius algebras. Since the source bimodule is cyclic, such a map is determined by the image of a cyclic generator, for example $1 \otimes 1$. Thus the maps $\beta$ and $\beta'$ are equivalent to bicentral  elements $e = \sum a_i \otimes b_i \in A \otimes A$ and $e' = \sum a'_i \otimes b'_i \in A \otimes A$. The cusp flip relations have a easy formulation in terms of this data. They assert the equality $e = e'$. Thus each saddle determines the other.  

The last piece of data assigned by an extended oriented topological field theory are the values of the cup and cap. Then this data is subject to the remaining relations, the 2D Morse relations. The value of the topological field theory on the circle is given by the vector space $Z(S^1) = A = N \otimes_{A \otimes B} M$. This is easily seen to coincide with $V = A/[A,A]$ (see Lemma~\ref{lem:universaltrace}). Thus the cup and cap yield linear maps:
\begin{align*}
	Z(\textrm{cap}) = u: k &\to V, \\
	Z(\textrm{cup})= \lambda: V &\to k.
\end{align*}

\begin{lemma}
The element $e$ given by the saddle, together with the cup morphism 
\begin{equation*}
\lambda: A \to A/[A,A] = V \to k
\end{equation*}
equip $A$ with the structure of a symmetric Frobenius $k$-algebra. 
\end{lemma}

\begin{proof} Since $\lambda$ factors through $A/[A,A]$, if these structures make $A$ into a Frobenius algebra, then $A$ will be symmetric.  Let $e = \sum a_i \otimes b_i \in A \otimes A$ be the element described above which corresponds to the value of (eitehr of) the saddles. We may rewrite the 2D Morse relations (involving the cup) in terms of $e$ and $\lambda$. They read:
	\begin{equation*}
		\sum \lambda(a_i)b_i = 1 \qquad \textrm{and} \qquad \sum a_i \lambda(b_i) = 1
	\end{equation*}
	which are exactly the Frobenius equations.
	% hold. The Frobenius equations are precisely the 2D Morse relations in disguise. 
\end{proof}

A 2-dimensional topological field theory naturally equips the value of the circle with the structure of a commutative Frobenius algebra. In the case of extended topological field theories this commutative algebra arrises as the center of the algebra assigned to the positive point, $z(A)$. As we will now see, under this isomorphism the image of the cap becomes the unit in $z(A)$.

\begin{lemma} \label{LemmaCircle=Center} Let $A = Z(pt^+)$ be the value of the positive point in an oriented topological field theory. There is a canonical isomorphism $A/[A,A] \cong z(A)$, where $z(A)$ is the center of the algebra $A$. 
\end{lemma}

\begin{figure}[ht]
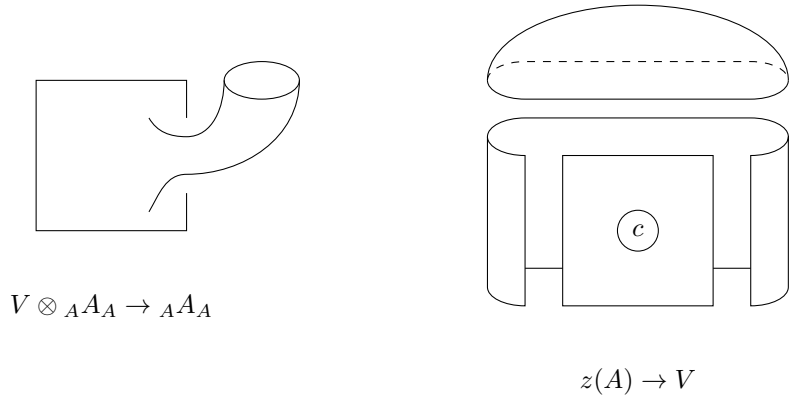

\begin{center}
\begin{center}

\end{center}
\caption{Maps Between the Value of the Circle and the Center.}
\label{FigCenterisCircle1}
\end{center}
\end{figure}

\begin{proof} Let $V = A/[A,A]$. Recall that the bimodule homomorphisms ${}_A A_A \to {}_A A_A$ can be canonically identified with the center $z(A)$. Gluing the basic generators as in Figure~\ref{FigCenterisCircle1} and evaluating via the field theory provides maps 
\begin{align*}
	V \otimes {}_A A_A & \to {}_A A_A \\
	z(A) & \to V
\end{align*}
and the first of these is equivalent to a map $V \to Z(A)$. Here we have represented an element $c \in z(A)$ graphically as an endomorphism of ${}_A A_A$. The remaining figures in Figure~\ref{FigCenterisCircle2} and Figure~\ref{FigCenterisCircle3} show that these maps are in fact inverses. The first equality in Figure~\ref{FigCenterisCircle3} follows by writing the surface in terms of the elementary bordisms and using standard manipulations available in any bicategory. 
\end{proof}

\begin{figure}[ht]
\begin{center}
\begin{center}
%

\end{center}
\caption{The Value of the Circle is the same as the Center, Part 1.}
\label{FigCenterisCircle2}
\end{center}
\end{figure}

\begin{figure}[ht]
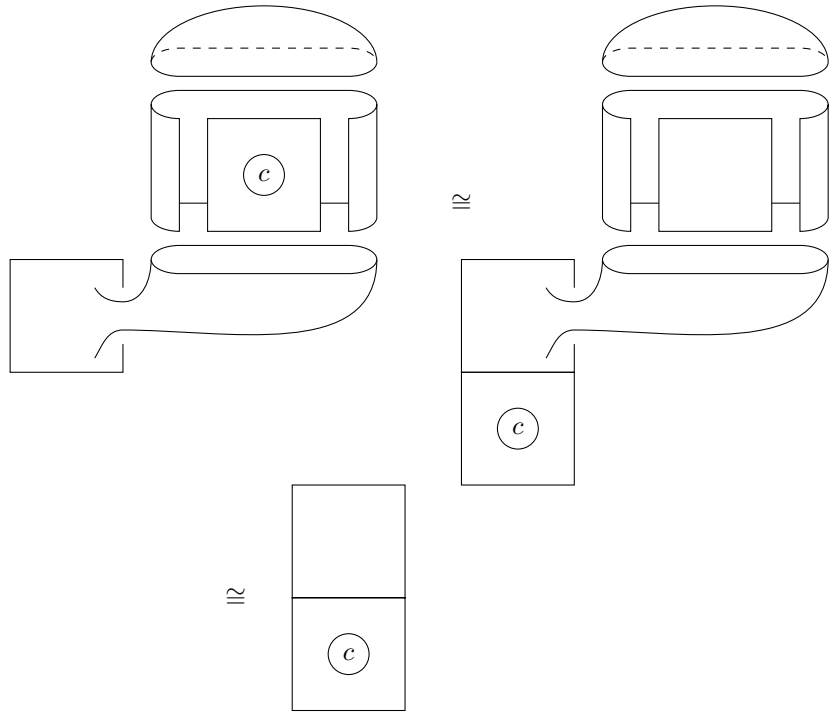

\begin{center}
\begin{center}
%
\end{center}
\caption{The Value of the Circle is the same as the Center, Part 2.}
\label{FigCenterisCircle3}
\end{center}
\end{figure}

\begin{remark}
A similar calculation shows that the multiplication on $V$, induced from the ``pants'' bordism of Figure~\ref{2DTQFT=CommFrobAlgFig} agrees with the multiplication on $z(A)$. 
\end{remark}

\begin{lemma}
$A$ is a separable algebra. 
\end{lemma}

\begin{proof}
As we have seen in Figures  \ref{FigCenterisCircle1}, \ref{FigCenterisCircle2} and \ref{FigCenterisCircle3}, the saddle can be used to give a map $V \otimes A \to A$. In terms of the element $e$, this map is given by:
\begin{equation*}
[x] \otimes y \mapsto \sum_i y a_i xb_i.
\end{equation*}
As we have seen this map also induces and isomorphism between $V$ the center $z(A)$, via
$[x] \in V \leftrightarrow \sum a_i x b_i \in z(A)$. Since this map is an isomorphism, there exists an element $[z]$ such that $\sum a_i z b_i = 1 \in z(A)$. Thus $\sum a_i \otimes (z b_i)$ gives $A$ the structure of a separable algebra. Letting $z$ range over all possible representatives allows us obtain all the possible separability structures, and so we are left with the property that $A$ is separable.
\end{proof}

\begin{proof}[Proof of {Theorem~\ref{ThmorienedCalssifcationInAlg}}]
	The analysis of this section has shown that 2-dimensional oriented extended topological field theories in $\alg^2$ are determined by the following collection of data:
	\begin{itemize}
		\item a separable algebra $A \in \alg^2$,
		\item a Morita context between $A$ and another algebra $B^\textrm{op}$,
		\item an bicentral element $e \in A \otimes A$ and a linear map $\lambda: A \to k$ such that $(A, e, \lambda)$ is a symmetric Frobenius algebra.  
	\end{itemize}
Our detailed analysis in Section~\ref{sec:transtqft} of morphisms of toplogical field theories lets us easily identify the morphisms between topological field theories determined by such data. In terms of the above a morphism of topological field theories consists of a pair of Morita contexts from $A$ to $\tilde{A}$ and $B^\textrm{op}$ to $\tilde{B}^\textrm{op}$ which send the given Morita context between $A$ and $B^\textrm{op}$ to the given Morita context between $\tilde{A}$ and $\tilde{B}^\textrm{op}$, and which further preserve the symmetric Frobenius algebra structures. Similarly the 2-morphisms are given by isomorphisms of Morita contexts satisfying the obvious compatibility relation. 
	
Now there is a clear forgetful functor from the bicategory of oriented extended 2-dimensional topological field theories 
\begin{equation*}
	\symbicat( \bord_2, \sC) \to \frob
\end{equation*}
which simple forgets the Morita context. This forgetful homomorphism is an equivalence of bicategories. By observation it is fully-faithful on 2-morphisms. Moreover for any 1-morphism the Morita context from $B^\textrm{op}$ to $\tilde{B}^\textrm{op}$ is determined, up to isomorphism, but the Morita context from $A$ to $\tilde{A}$. This shows that the forgetful functor is full on 1-morphisms, not just essentially full. Similarly given a separable symmetric Frobenius algebra, we can always lift it to a tqft by choosing $B = A^\textrm{op}$ and the identity Morita context. Thus the forgetful functor is surjective on objects. 
\end{proof}

\subsection{Unoriented TQFTs and Stellar algebras} \label{sec:unorientedandstellar} Now we will extend the previous results to the case of unoriented topological field theories and prove Theorem~\ref{ThmUnorientedInAlgs}. An unoriented topological field theory still gives an oriented topological field theory, and so we can simply continue with the analysis we preformed in the previous section. Recall that we were able to show that an oriented 2-dimensional extended topological field theory in $\alg^2$ was given by: 
	\begin{itemize}
		\item a separable algebra $A \in \alg^2$,
		\item a Morita context between $A$ and another algebra $B^\textrm{op}$,
		\item an bicentral element $e \in A \otimes A$ and a linear map $\lambda: A \to k$ such that $(A, e, \lambda)$ is a symmetric Frobenius algebra.  
	\end{itemize}
The main difference in the unoriented case is that there is no distinction between the positive and negative point. Thus in the unoriented case we have $B \cong A$. Thus the Morita context above becomes a Morita context from $A$ to $A^\textrm{op}$. 

\begin{figure}[htbp]
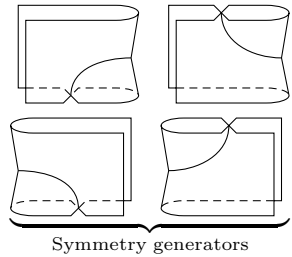


		}_\text{Symmetry generators}
		$};
	\end{tikzpicture}
	\caption{Symmetry Generators}
	\label{fig:symgen}
\end{figure}

Moreover there are additional generating 2-morphisms and several sets of additional relations. The additional generating 2-morphisms are the {\em symmetry generators} which are depicted in Figure~\ref{fig:symgen}.

Let $M$ be an $A$-$A^\text{op}$-bimodule. Then we may flip both actions to obtain a new $A$-$A^\text{op}$-bimodule, which we denote $\underline{M}$ and call the {\em conjugate}. The conjugate of a Morita context from $A$ to $A^\textrm{op}$ is a new Morita context from $A$ to $A^\textrm{op}$. It makes sense to consider bimodule maps from $M$ to $\underline{M}$, and in particular isomorphisms.

\begin{figure}[htbp]
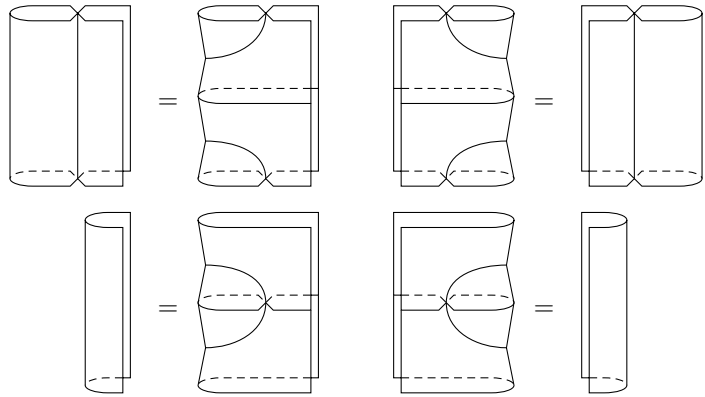

	\begin{center}
	
	\end{center}
	\caption{Relations implying isomorphism of Morita context}
	\label{fig:stellarmorita}
\end{figure}

Let ${}_AM_{A^{\textrm{op}}}$ and ${}_{A^\textrm{op}}N_A$ be the bimodule comonents of a Morita context between $A$ and $A^\textrm{op}$. In terms of these operations the symmetry generators become maps of bimodules:
\begin{align*}
	\underline{M} \to M \\
	M \to \underline{M} \\
	N \to \underline{N} \\
	\underline{N} \to N.
\end{align*} 
The first set of additional relations, depicted in Figure~\ref{fig:stellarmorita}, state that these maps are inverses and hence form an isomorphism, which we denote $\sigma$ of Morita contexts between the given Morita context from  $A$ to $A^\textrm{op}$ and its conjugate. 

\begin{figure}[htbp]
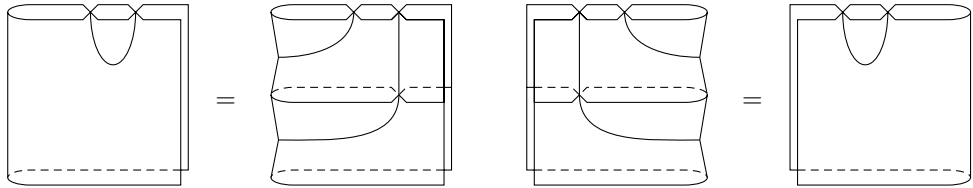

	\begin{center}

	\end{center}
	\caption{Coherence Relations}
	\label{fig:stellarcoherence}
\end{figure}

The second set of additional relations determine a set coherence relations for the isomorphism $\sigma$. These relations are depicted in Figure~\ref{fig:stellarcoherence}. Let $s =  ({}_{A^\text{op}} M_A, {}_AN_{A^\text{op}}, \eta, \varepsilon)$ be a Morita context between $A$ and $A^\textrm{op}$, and let $\sigma: \underline{s} \cong s$. The identity morphism gives a natural isomorphisms $\underline{\underline{s}} \cong s$. The isomorphism $\sigma$ induces another, by applying it twice:
\begin{equation*}
	\sigma \circ \underline{\sigma}: \underline{\underline{s}} \cong \underline{s} \cong s.
\end{equation*}
The additional relation is equivalent to the statement that these two isomorphisms agree. 

The above structures form a sort of enhancement of the notion of algebraic $*$-structure which is Morita invariant. We will codify such structure as a {\em stellar structure}:

\begin{definition}
A {\em stellar algebra} is an algebra $A$, equipped with a Morita context between $A$ and $A^\text{op}$,  $s =  ({}_{A^\text{op}} M_A, {}_AN_{A^\text{op}}, \eta, \varepsilon)$, together with an isomorphism of Morita contexts, $\sigma: s \cong \underline s$, such that $\sigma \circ \underline{\sigma}$ is the identity map.   
\end{definition}

\begin{example}
	An (algebraic) $*$-algebra gives rise to a stellar algebra. The $*$-structure on $A$ is the same as an isomorphism of algebras, $*: A \to A^\text{op}$. Then ${}_A A*_{A^\text{op}}$ and ${}_{A^\text{op}}A^\text{op}*^{-1}_{A}$, with the canonical unit and counit, form a Morita context $s$ between $A$ and $A^\text{op}$. The identity $*^2 = id$, then induces an isomorphism $s \cong \underline s$, which may be checked to satisfy the additional coherence equation. 
\end{example}

\begin{example} \label{ex:stellarnotstar}
	Not every stellar algebra arises from a $*$-structure. Let $R$ be an algebra such that $R$ is not isomorphic to $S = R^\text{op}$, and such that $R$ is also not isomorphic to $M_2(R)$.\footnote{For example if $k$ is a field whose Brauer group contains elements of order greater then two, then it will have algebras of this kind.} Let $A = M_2(R) \oplus S$. We have $A^\text{op} = M^2(S) \oplus R$, so that $A \not \cong A^\text{op}$, and hence $A$ cannot support a $*$-structure. Nevertheless, $A$ can be made into a stellar algebra in a straightforward manner. We leave the details as an exercise. 
\end{example}

Stellar structures can be transfered via Morita contexts, just as symmetric Frobenius algebra structures. Specifically, if $f= ({}_B M_A, {}_AN_B, \eta, \varepsilon)$ is a Morita context between $A$ and $B$ and $(A, s, \sigma: s \cong \underline s)$ is a stellar structure on $A$, then $(B, f_*s, f_* \sigma)$ is a stellar structure on $B$, where $f_*s = \underline f \circ s \circ f$ and $f_* \sigma = \underline f \circ \sigma \circ f$.

\begin{definition}
Let $(A, s^A, \sigma^A: s^A \cong \underline{s}^A)$ and $(B, s^B, \sigma^B: s^B \cong \underline{s}^B)$ be two stellar algebras. A 1-morphism of stellar algebras (from $A$ to $B$) consists of a Morita context $f$(between $A$ and $B$) together with an isomorphism of Morita contexts:
\begin{align*}
	\phi: f_* s^A \to s^B.
\end{align*}
A 2-morphism between 1-morphisms $(f, \phi)$ and $(f', \phi')$ is an isomorphism of Morita contexts $\alpha: f \to f'$ such that the following triangle of isomorphism of Morita contexts commutes:

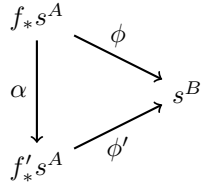
\begin{center}

\end{center}
\end{definition}
The stellar algebras, 1-morphism, and 2-morphisms form a symmetric monoidal bicategory in the obvious way. 

\begin{figure}[htbp]
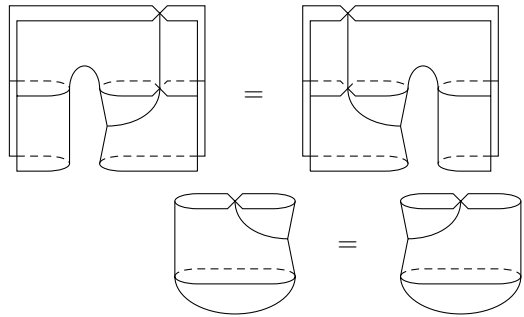

	\begin{center}
	%
	\end{center}
	\caption{Compatibility between Frobenius structure and stellar structure.}
	\label{fig:Frobandstellar}
\end{figure}

Thus we see that an extended 2-dimensional unoriented topological field theory consists of the data of a separable symmetric Frobenius algbera $(A, e, \lambda)$ together with the structure of a stellar algebra $(A, s, \sigma)$. These structures are required to satisfy one additional set of relations, which are depicted in Figure~\ref{fig:Frobandstellar}.
In the presence of the other relations these are redundant. This compatibility can be expressed as follows. Given the symmetric Frobenius algebra $(A, e, \lambda)$ and a stellar structure $(A, s, \sigma)$ we obtain two linear maps
\begin{equation*}
	N \otimes_{A \otimes A} \tau \otimes_{A \otimes A} M \to N \otimes_{A \otimes A} N \to k 
\end{equation*}
here $\tau$ is the `flip' bimodule.  The first map uses the first half of $\sigma$ to identify $N \otimes \tau \cong N$ and then applies $\lambda$. The second map uses the second part of $\sigma$ to identity $\tau \otimes M \cong M$, and then applies $\lambda$. The final set of relations are satisfied precisely when these two maps agree.

\begin{definition}
Let $\frob^*$ be the symmetric monoidal bicategory of separable stellar symmetric Frobenius algebras, that is the objects consist of quintuples $(A, s, \sigma, e, \lambda)$, where $A$ is separable,  $(A, s, \sigma)$ is a stellar algebra, $(A, e, \lambda)$ is a symmetric Frobenius algebra, and the Frobenius and stellar structure are compatible in the sense described above.  The 1-morphism of $\frob^*$ consist of those {1-morphisms} $(f, \phi)$ of stellar algebras for which $f$ is compatible with the symmetric Frobenius algebra structure. The 2-morphism of $\frob^*$ are the 2-morphisms of stellar algebras. 
\end{definition}

Summarizing, we have shown that the bicategory of 2-dimensional extended unoriented topological field theories in $\alg^2$ is equivalent to the bicategory $\frob^*$. This completes the proof of Theorem~\ref{ThmUnorientedInAlgs}. \qed 

%In the oriented case, the algebras $A$ and $B$ may be distinct. We have translated all the structure into structures on $A$. The structures on $B$ are derived from these using the Morita context between $A$ and $B^\text{op}$. However, the algebra $B$ is actually irrelevant. The bicategory of separable symmetric Frobenius algebras $A$, together with an algebra $B$ and  a Morita equivalence $A \simeq B^{op}$ is equivalent to the bicategory which forgets the algebra $B$ and the Morita equivalence. This is a straightforward exercise which we leave to the reader. This completes the proof of Theorem~\ref{ThmorienedCalssifcationInAlg}.

%\begin{corollary}
%The bicategory of 2-dimensional unoriented extended topological field theories with values in $\alg^2_k$, with $k$ a field, is equivalent to the bicategory whose objects are  separable symmetric Frobenius $*$-algebras
%\end{corollary}

 % The Main Theorem (Bordism Bicategories).

\appendix

%%%%%%% %\include{app1} 	% Algebra   %%%%%% now irrelevant. 

\chapter{Bicategories} \label{Chap:Bicats}

\section[Bicats, Homs, Trans, and Mods]{Bicategories, Homomorphisms, Transformations and Modifications}

\begin{definition}
A {\em bicategory} $\sB$ 
\index{bicategory} 
consists of the following data,
\begin{enumerate}
\item A collection $\text{ob} \sB$ whose elements are called {\em objects}. (We denote objects $a, b, c \in \sB$. Sometimes the collection $\text{ob} \sB$ is written $\sB_0$).
\item Categories $\sB(a, b)$ for each pair of objects $a, b \in  \sB$. The objects of $\sB(a, b)$ are referred to as {\em 1-morphisms} from $a$ to $b$ (which we denote collectively as $\sB_1$, with elements $f,g, h,\dots$). The morphisms of $\sB(a, b)$ are referred to as {\em 2-morphisms} (which we denote collectively by $\sB_2$, with elements $\alpha, \beta, \gamma, \dots$).
\item Functors:
\begin{align*}
c_{a b c}: \sB(b,c) \times \sB(a,b ) & \to \sB(a,c) \\
(g,f) & \mapsto g \circ f \\
(\beta, \alpha) & \mapsto \beta * \alpha
\end{align*}
and $I_a: {\bf 1} \to \sB(a,a)$, for all objects $a, b, c \in \sB$,  where ${\bf 1}$ denotes the singleton category (thus the functor $I_a$ is equivalent to specifying  an {\em identity element} $I_a$). The functors $c$ are called the {\em horizontal compositions}.  
\item Natural isomorphisms:
\begin{align*}
	a : & \; c_{abd} \circ (c_{bcd} \times id) \to c_{acd} \circ (id \times c_{abc}) \\
	\ell : &\;  c_{abb} \circ (I_b \times id)  \to id \\
	r: & \; c_{aab} \circ (id \times I_a) \to id
\end{align*}
known, respectively, as the {\em associator} and left and right {\em unitors}. 

(thus invertible 2-morphisms $a_{h,g,f}: (h \circ g) \circ f \to h \circ ( g \circ f)$, $\ell_f: I_b \circ f \to f$ and $r_f: f \circ I_a \to f$). 
\end{enumerate}
These are required to satisfy the pentagon and triangle identities:    \begin{center}
\begin{tikzpicture}[thick]
	\node at (0, 4.5) {\underline{The Pentagon Identity}};
	\node (LT) at (0,2) 	{$((kh)g)f $ };
	\node (LB) at (1,0) 	{$(k(hg))f$};
	\node (RT) at (8,2) 	{$k(h(gf))$};
	\node (MT) at (4,3.5) {$(kh)(gf)$};
	\node (RB) at (7,0)	{$ k((hg)f)$};
	\draw [->] (LT) --  node [left] {$a * 1 $} (LB);
	\draw [->] (LT) -- node [above left] {$a $} (MT);
	\draw [->] (MT) -- node [above right] {$a $} (RT);
	\draw [->] (RB) -- node [right] {$1* a $} (RT);
	\draw [->, in = 200, out = 340] (LB) to node [below] {$a $} (RB);
\end{tikzpicture}
\end{center}
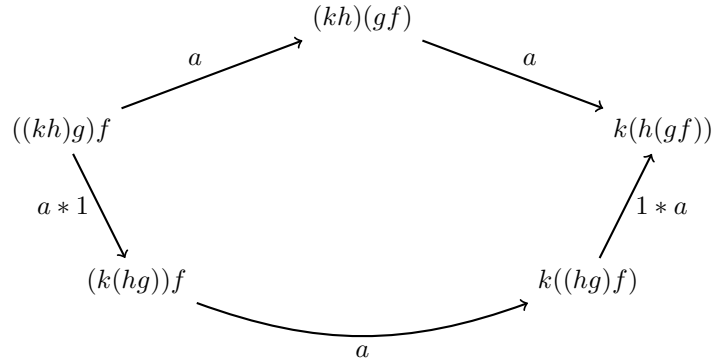
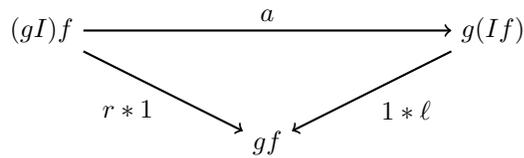
\begin{center}
\begin{tikzpicture}[thick]
	\node at (0,2.5) 	{\underline{The Triangle Identity}};
	\node (LT) at (0,1.5) 	{$ (gI)f  $ };
	
	\node (RT) at (6,1.5) 	{$g(If) $};
	\node (B) at (3,0)	{$gf$};
	\draw [->] (LT) --  node [below left] {$r * 1 $} (B);
	\draw [->] (LT) -- node [above] {$a $} (RT);
	\draw [->] (RT) -- node [below right] {$1*\ell $} (B);
\end{tikzpicture}
\end{center}
A bicategory in which the associators and left and right unitors are identities, so that composition is strictly associative,  is called a {\em 2-category}. A bicategory $\sA$ is {\em pointed} if it is equipped with a distinguished object $p \in \sA$. 
\index{2-category}
\index{2-category|see{bicategory}}
\end{definition}

\begin{warning}
The horizontal composition of 1-morphism is denoted $f \circ g$, which is contradictory to the notation for 2-morphisms. The vertical composition of 2-morpisms is denoted $\alpha \circ \beta$, while the horizontal composition is denoted $ \alpha * \beta$. Regrettably this is standard notation, and so we reluctantly comply. 
\end{warning}

%\begin{example}
%The primordial example of a 2-category is $\cat$, the 2-category of small categories, functors and natural transformations.
%\end{example}

%The next two examples, while important, are not very interesting. We will encounter numerous bicategories later in our exposition.  

%\begin{example}
%Any category can be considered a 2-category in which the only 2-morphisms are the identity 2-morphisms.
%\end{example}

%\begin{example}
%Given a commutative monoid $A$ we can form a 2-category $A[2]$ which has one object, one 1-morphism, and the 2-morphisms consist of $A$. Both horizontal and vertical composition coincide and are given by the multiplication in $A$. Any bicategory with a single object and a single 1-morphism is of this form. 
%\end{example}

%\begin{remark} \label{underlyingcatofa2category}
%Given a 2-category there is an underlying category obtained by disregarding the 2-morphisms. The underlying category of $\cat$ is $\Cat{cat}$. 
%\end{remark}

Let $\sA$ be a bicategory and let $\sigma: a \to b$ be a 1-morphism in $\sA$. Pre- and post-composition define functors,
\begin{align*}
	\sigma^* &: \sA(b, c) \to \sA(a, c) \\
	\sigma_* &: \sA(d, a) \to \sA(d, b).
\end{align*}
The associators define canonical natural isomorphisms for each 1-morphism $\theta \in \sA_1$. 
\begin{align*}
\sigma_* \circ \theta_* &\cong (\sigma \circ \theta)_* \\
\sigma^* \circ \theta^* & \cong (\theta \circ \sigma)^* \\
\sigma^* \circ \theta_* & \cong \theta_* \circ \sigma^* 
\end{align*}
provided the corresponding compositions are allowed.

%\begin{definition}
%Given a bicategory $\sB$, there is a bicategory $\sB^{op}$ obtained by reversing the 1-morphisms, but not the 2-morphisms. 
%\end{definition}

\begin{definition} [Internal Equivalence]
Let $\sB$ be a bicategory. An {\em equivalence} \index{bicategory!equivalence in a} in $\sB$ is a 1-morphism $f: a \to b$ in $\sB(a, b)$ such that there exists a 1-morphism $g: b \to a$ in $\sB(b, a)$ with isomorphisms $\eta: I_a \to g \circ f$ in $\sB(a,a)$ and $\varepsilon: f \circ g \to I_b$ in $\sB(b,b)$. In this case we say that $a$ is {\em equivalent} to $b$.
\end{definition}

\begin{remark}
It is clear that the relation of equivalence is an equivalence relation on the objects of $\sB$. The set of equivalence classes is denoted $\pi_0 \sB$.
\end{remark}

\begin{definition} \label{def:bicathomomorphism}
Let $\sA$ and $\sB$ be bicateogries. A {\em homomorphism} 
\index{bicategory!homomorphism} 
\index{homomorphism}
\index{homomorphism|see{bicategory}}
$F: \sA \to \sB$ consists of the data:
\begin{enumerate}
\item A function $F: \text{ob } \sA \to \text{ob } \sB$,
\item Functors $F_{ab}: \sA(a,b) \to \sB( F(a), F(b))$,
\item Natural isomorphisms 
\begin{align*}
\phi_{abc}:  c^\sB_{F(a) F(b) F(c) } \circ (F_{bc} \times F_{ab})   & \to   F_{ac} \circ c^\sA_{abc}\\
\phi_a: I^\sB_{F(a)} &\to F_{aa} \circ I_a^\sA 
\end{align*}
(thus invertible  2-morphisms $\phi_{gf}: Fg \circ F f \to F(g \circ f)$ and $\phi_a: I^\sB_{Fa} \to F( I^\sA_a)$ ).
\end{enumerate}
such that the following diagrams commute:
	%\item For every triple of composable 1-morphisms $f, g, h$ a commutative diagram,
		
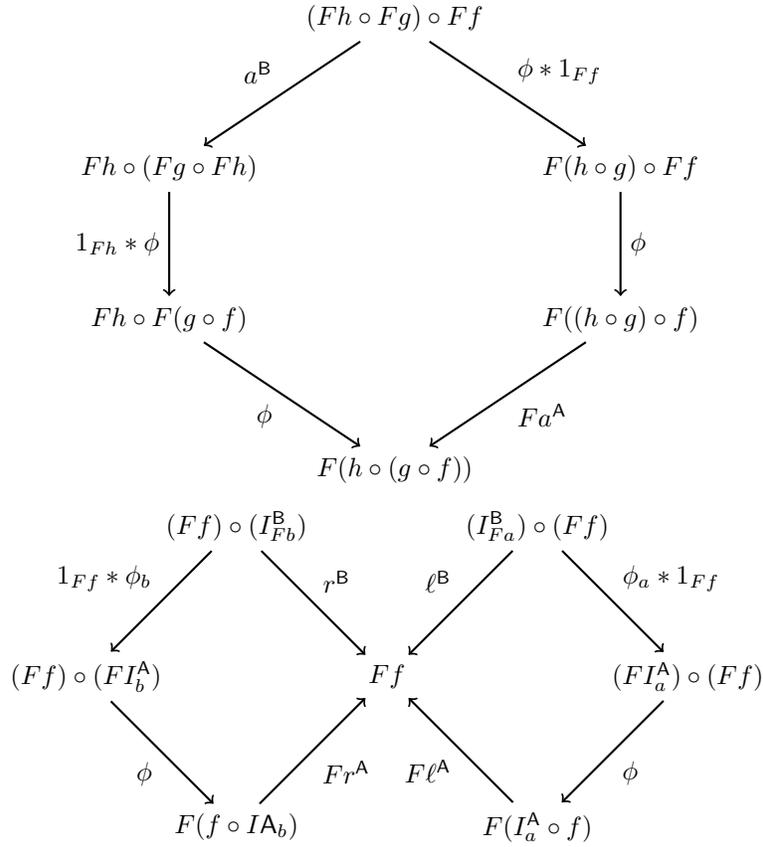
\begin{center}
		\begin{tikzpicture}[thick]
			\node (A) at (3, 6) {$(Fh \circ Fg) \circ Ff$};
			\node (B) at (0,4) {$Fh \circ (Fg \circ Fh)$};
			\node (C) at (0,2) {$Fh \circ F(g \circ f)$};
			\node (D) at (3,0) {$F(h \circ (g\circ f))$};
			\node (E) at (6,2) {$F((h \circ g) \circ f)$};
			\node (F) at (6,4) {$F(h \circ g) \circ Ff$};
			
			\draw [->] (A) -- node [above left] {$a^\sB$} (B);
			\draw [->] (B) -- node [left] {$1_{Fh} * \phi$} (C);
			\draw [->] (C) -- node [below left] {$\phi$} (D);
			\draw [->] (A) -- node [above right] {$\phi * 1_{Ff}$} (F);
			\draw [->] (F) -- node [right] {$\phi$} (E);
			\draw [->] (E) -- node [below right] {$Fa^\sA$} (D);
		\end{tikzpicture}
		\end{center}
	%	\item If $f$ is a 1-morphism with source $A$ and target $B$, commutative diagrams,
		\begin{center}
		\begin{tikzpicture}[thick]
			\node (A) at (2, 4) {$(Ff) \circ (I_{Fb}^\sB)$};
			\node (B) at (0,2) {$(Ff) \circ (FI^\sA_b)$};
			\node (C) at (2,0) {$F(f \circ I\sA_b)$};
			\node (D) at (4,2) {$Ff$};
			\node (E) at (6,4) {$(I_{Fa}^\sB) \circ (Ff)$};
			\node (F) at (8,2) {$(FI_a^\sA) \circ (Ff)$};
			\node (G) at (6,0) {$F(I_a^\sA \circ f)$};
			
			\draw [->] (A) -- node [above left] {$1_{Ff} * \phi_b$} (B);
			\draw [->] (B) -- node [below left] {$ \phi$} (C);
			\draw [->] (C) -- node [below right] {$Fr^\sA$} (D);
			\draw [->] (A) -- node [above right] {$r^\sB$} (D);
			\draw [->] (E) -- node [above right] {$\phi_{a} * 1_{Ff}$} (F);
			\draw [->] (F) -- node [below right] {$\phi$} (G);
			\draw [->] (G) -- node [below left] {$F\ell^\sA$} (D);
			\draw [->] (E) -- node [above left] {$\ell^\sB$} (D);

		\end{tikzpicture}
		\end{center}
If the natural isomorphisms $\phi_{abc}$ and $\phi_a$ are identities, then the homomorphism $F$ is called a {\em strict homomorphism}. A homomorphism $F$ between pointed bicategories $(\sA, p_\sA)$ and $(\sB, p_\sB)$ is a {\em pointed homomorphism} if $F(p_\sA) = p_\sB$.
\end{definition}

\begin{definition} \label{DefnBicatTransformation}
Let $(F, \phi), (G, \psi): \sA \to \sB$ be two homomorphisms between bicategories. A {\em transformation} 
\index{bicategory!transformation} 
\index{transformation} 
\index{transformation|see{bicategory}}
$\sigma: F \to G$ is given by the data:
\begin{enumerate}
\item 1-morphisms $\sigma_a: Fa \to Ga$ for each object $a \in \sA$,
\item Natural Isomorphisms, $\sigma_{ab}: (\sigma_a)^* \circ G_{ab} \to (\sigma_b)_* \circ F_{ab}$

(thus invertible  2-morphisms $\sigma_f: Gf \circ \sigma_a \to \sigma_b \circ Ff$ for every $f \in \sA_1$).
\end{enumerate}
such that the diagrams in Figure~\ref{FigBicatTransformationAxioms} commute for all 1-morphisms in $\sA$, $f: a \to b$ and $g: b \to c$.
\begin{figure}[ht]
\begin{center}
	\begin{center}
		\begin{tikzpicture}[thick]
			\node (A) at (0,6) {$(Gg \circ Gf) \circ \sigma_a$};
			\node (B) at (0,4) {$Gg \circ(Gf \circ \sigma_a)$};
			\node (C) at (0,2) {$ Gg \circ( \sigma_b \circ Ff)$};
			\node (D) at (3,0) {$ (Gg \circ \sigma_b) \circ Ff$};
			\node (E) at (6,2) {$(\sigma_c \circ  Fg) \circ Ff$};
			\node (F) at (6,4) {$\sigma_c \circ (Fg \circ Ff)$};
			\node (G) at (6,6) {$\sigma_c \circ F(g \circ f)$};
			\node (H) at (3,8) {$G(g \circ f) \circ \sigma_a$};
			
			\draw [->] (A) -- node [left] {$a^\sB$} (B);
			\draw [->] (B) -- node [left] {$id_{Gg} * \sigma_f$} (C);
			\draw [->] (C) -- node [below left]  {$(a^\sB)^{-1}$} (D);
			\draw [->] (D) -- node [below right] {$\sigma_g * id_{Ff}$} (E);
			\draw [->] (E) -- node [right] {$a^\sB$} (F);
			\draw [->] (F) -- node [right] {$id_{\sigma_c} * \phi_{g,f}$} (G);
			\draw [->] (A) -- node [above left] {$\psi_{g,f} * id_{\sigma_a}$} (H);
			\draw [->] (H) -- node [above right] {$\sigma_{gf}$} (G);
		\end{tikzpicture}
		\end{center}
	%	\item For every object $A$ in $\mathcal{C}$, a commutative diagram,
		\begin{center}
		\begin{tikzpicture}[thick]
			\node (A) at (0,2) {$I^\sB_{Ga} \circ \sigma_a$};
			\node (B) at (3,3) {$\sigma_a $};
			\node (C) at (6,2) {$ \sigma_a \circ I^\sB_{Fa}$};
			\node (D) at (5,0) {$ \sigma_a \circ (FI^\sA_a)$};
			\node (E) at (1,0) {$(GI^\sA_a) \circ \sigma_a$};
			
			\draw [->] (A) -- node [above left] {$\ell^\sB$} (B);
			\draw [->] (B) -- node [above right]  {$(r^\sB)^{-1}$} (C);
			\draw [->] (A) -- node [below left] {$\psi_a * id_{\sigma_a}$} (E);
			\draw [->] (E) -- node [below] {$\sigma_{I^\sA_a}$} (D);
			\draw [->] (D) -- node [below right] {$id_{\sigma_a} * \phi_a^{-1}$} (C);
		\end{tikzpicture}
		\end{center}
\caption{Transformation Axioms}
\label{FigBicatTransformationAxioms}
\end{center}
\end{figure}
If the natural transformations $\sigma_{ab}$ are identities, then $\sigma$ is a {\em strict transformation}.
If $\sA$, $\sB$, $(F, \phi)$ and $(G, \psi)$ are pointed, then the transformation $\sigma$ is a {\em pointed transformation} if $\sigma_{p_\sA} = I_{p_\sB}$. 
\end{definition}

\begin{remark} \label{RmkLaxAndOpLax}
	There is a weaker notion of {\em lax transformation} in which the natural transformations $\sigma_{f}$ are not assumed to be isomorphisms. This imposes an obvious directional bias. An {\em oplax transformation} is defined similarly, but with the direction of $\sigma_f$ reversed. When comparing these notions of transformation, the transformations, as defined above with invertible 2-morphism data, are called {\em strong} transformations. 
\end{remark}

\begin{definition}
Let $(F, \phi), (G, \psi): \sA \to \sB$ be two homomorphisms between bicategories and let $\sigma, \theta: F \to G$ be two transformations between homomorphisms. A {\em modification} \index{bicategory!modification} 
\index{modification}
\index{modification|see{bicategory}} 
$\Gamma: \sigma \to \theta$ consists of  2-morphisms $\Gamma_a: \sigma_a \to \theta_a$ for every object $a \in \sA$,
such that the following square commutes:
\begin{center}
\begin{tikzpicture}[thick]
	\node (LT) at (0,1.5) 	{$Gf \circ \sigma_a$ };
	\node (LB) at (0,0) 	{$\sigma_b \circ F f$};
	\node (RT) at (3,1.5) 	{$Gf \circ \theta_a$};
	\node (RB) at (3,0)	{$\theta_b \circ Ff$};
	\draw [->] (LT) --  node [left] {$\sigma_f$} (LB);
	\draw [->] (LT) -- node [above] {$id * \Gamma_a$} (RT);
	\draw [->] (RT) -- node [right] {$\theta_f$} (RB);
	\draw [->] (LB) -- node [below] {$\Gamma_b * id$} (RB);
\end{tikzpicture}
\end{center}
for every 1-morphism $f: a \to b$ in $\sA$. If $\sA$, $\sB$, $(F, \phi)$, $(G, \psi)$, $\sigma$ and $\theta$ are pointed, then we say $\Gamma$ is a {\em pointed modification} if $\Gamma_{p_\sA}$ is the identity of $I_{p_\sB}$.
\end{definition}

\begin{definition} \label{defncompositonoftransformations}
Let $(F, \phi)$, $(G, \psi)$ and $(H, \kappa)$ be homomorphisms from the bicategory $\sA$ to the bicategory $\sB$. Let $\sigma: F \to G$ and $\theta: G \to H$ be two transformations. We define the composition $\theta \circ \sigma$ to be the transformation $F \to H$ given by the following data:
\begin{enumerate}
\item The 1-morphism $\theta_a \circ \sigma_a$ for each object $a\in \sA$,
\item The natural isomorphism $ (\theta_a \circ \sigma_a)^* \circ H_{ab} \to (\theta_a \circ \sigma_a)_* \circ F_{ab}$ defined by the following sequence of natural isomorphisms:
\begin{center}
\begin{tikzpicture}[thick]
	\node (LT) at (0,1.5) 	{$(\theta_a \circ \sigma_a)^* \circ H_{ab} $ };
	\node (MT) at (4,1.5) 	{$ \sigma^*_a \circ \theta_a^* \circ H_{ab} $ };
	\node (RT) at (8,1.5) 	{$ \sigma^*_a \circ (\theta_b)_* \circ G_{ab} $};
	\node (LB) at (0,0) 	{$ (\theta_b)_* \circ \sigma^*_a  \circ G_{ab} $};
	\node (MB) at (4,0) 	{$(\theta_b)_* \circ (\sigma_a)_*  \circ F_{ab}$};
	\node (RB) at (8,0)	{$(\theta_a \circ \sigma_a)_* \circ F_{ab}$};
	\draw [->] (LT) -- node [above] {$$} (MT);
	\draw [->] (MT) -- node [above left] {$\theta $} (RT);
	\draw [->, out = -15 , in = 165] (RT) -- node [right] {$$} (LB);
	\draw [->] (LB) -- node [below right] {$\sigma $} (MB);
	\draw [->] (MB) -- node [below] {$$} (RB);
\end{tikzpicture}
\end{center}
 where the unlabeled arrows are the canonically defined natural transformations induced by the associators of $\sB$. 
\end{enumerate}
Given transformations $\sigma, \sigma', \sigma'': F \to G$ and modifications $\Gamma: \sigma \to \sigma'$ and $\Sigma: \sigma' \to \sigma''$, the {\em vertical composition} of $\Gamma$ and $\Sigma$ is the modification $\Sigma \circ \Gamma: \sigma \to \sigma''$ given by the 2-morphisms $ \Sigma_a \circ \Gamma_a$ for each object $a \in \sA$. 

Similarly, if $(F, \phi)$, $(G, \psi)$ and $(H, \kappa)$ are homomorphisms from the bicategory $\sA$ to the bicategory $\sB$, $\sigma, \sigma': F \to G$ and $\theta, \theta': G \to H$ are transformations and $\Gamma: \sigma \to \sigma'$ and $\Sigma: \theta \to \theta'$ are modifications, then the {\em horizontal composition} of  $\Gamma$ and $\Sigma $ is the modification $\Sigma * \Gamma:  \theta \circ \sigma \to \theta' \circ \sigma'$ given by the 2-morphisms $\Sigma_a * \Gamma_a$ for each object $a \in \sA$.  
\end{definition}

With these compositions we get  a bicategory $\bicat(\sA, \sB)$ for each pair of bicategories $\sA, \sB$. The  objects are the homomorphisms from $\sA$ to $\sB$, the 1-morphisms are the transformations, and the 2-morphisms are the modifications. Compositions are as above and the associators  and unitors are the obvious ones coming from $\sB$.

\begin{remark} \label{inthomofbicatisa2-cat} Since the associators and unitors come from those in $\sB$, 
$\bicat(\sA, \sB)$ is a 2-category whenever $\sB$ is a 2-category. 
\end{remark}

If $\sA$ and $\sB$ are pointed bicategories, we can also consider the pointed homomorphisms, transformations and modifications. The vertical composition of two pointed modifications is again a pointed modification, but the horizontal composition of two pointed transformations or modifications fails to be pointed. However we can fix this. The problem can be seen at the level of transformations. If $\sigma$ and $\theta$ are two compatible pointed transformations, as above, then the transformation $\theta \circ \sigma$ satisfies:
\begin{equation*}
	(\theta \circ \sigma)_{p_\sA} = I_{p_\sB} \circ I_{p_\sB}
\end{equation*}
This is not equal to $I_{p_\sB}$ and hence this is not a pointed transformation. However there is a canonical isomorphism: 
\begin{equation*}
 I_{p_\sB} \circ I_{p_\sB} \to  I_{p_\sB} 
\end{equation*}
given by the right or left unitor (they yield the same isomorphism). Using this isomorphism we can replace the transformation $(\theta \circ \sigma)$ with a pointed transformation $(\theta \circ \sigma)'$ as follows:
\begin{align*}
(\theta \circ \sigma)'_{a} &= \begin{cases} 
		(\theta \circ \sigma)_{a} & \text{if } a \neq p_\sA \\
			I_{p_\sB} &  \text{if } a = p_\sA \end{cases}  \\
		(\theta \circ \sigma)'_f &= \begin{cases} 
		(\theta \circ \sigma)_{f} & \text{if } f:a \to b \text{ and } a,b \neq p_\sA \\
			(\theta \circ \sigma)_{f} \circ (\text{can.}) & \text{if } f:a \to b, \; a= p_\sA, \text{ and } b \neq p_\sA \\
			(\text{can.}) \circ (\theta \circ \sigma)_{f} & \text{if } f:a \to b, \; b= p_\sA, \text{ and } a \neq p_\sA \\
			(\text{can.}) \circ (\theta \circ \sigma)_{f} \circ (\text{can.}) & \text{if } f:a \to b, \; a=b= p_\sA
			\end{cases}  
\end{align*}
where ``$(\text{can.})$'' denotes the above canonical 2-morphisms. The reader can check that this is in fact a transformation. This is defines the pointed composition of pointed transformations. A similar construction yields the pointed horizontal composition of pointed modifications. Again we obtain a bicategory $\bicat_*(\sA, \sB)$ of pointed homomorphism, pointed transformations, and pointed modifications. The canonical constant functor sending every object of $\sA$ to $p_\sB$ gives $\bicat_*(\sA, \sB)$ the structure of a pointed bicategory.

These bicategories share an additional structure. There are composition homomorphisms,\begin{equation*}
\bicat(\sB, \sC) \times \bicat(\sA, \sB) \to \bicat(\sA, \sC)
\end{equation*}
 which have their own additional coherence structures. This can be summarized by saying there is a tricategory whose objects are bicategories and whose hom-bicategories are the bicategories $\bicat(\sA, \sB)$. We will not need the full details of these coherence structures, but the interested reader should consult \cite{GPS95}.  Actually there is some ambiguity in defining what the new composition of transformations and modifications should be. One way to address this is to introduce the operation of {\em whiskering}, 
\index{bicategory!whiskering} 
\index{whiskering}
\index{whiskering|see{bicategory}} an approach which we sketch here. This yields  two different but canonical ways to compose transformations (and modifications), and hence results in two different tricategories.  However, these two tricategories are equivalent as tricategories, again see \cite{GPS95} for full details.  
 
\begin{definition}
Let $(F, \phi): \sA \to \sB$ and $(G, \psi): \sB \to \sC$ be two homomorphisms between bicategories. We define their composition to be the homomorphism $G \circ F : \sA \to \sC$ given by the following data:
\begin{enumerate}
\item On objects we have the map $G \circ F: \text{ob } \sA \to \text{ob } \sB$,
\item Functors $G_{Fa, Fb} \circ F_{ab}: \sA(a, b) \to \sA( GF(a), GF(b))$,
\item  Natural transformations which are the composite,
\begin{equation*}
	I_{GF(a)} \stackrel{\psi}{\to} G( I_{Fa}) \stackrel{G( \phi)}{\to} GF(I_a),
\end{equation*}
and
\begin{center}
\begin{tikzpicture}[thick]
	\node (T) at (0,1.5) 	{$c^\sC_{GF(a) GF(b) GF(c)} \circ (G_{Fb Fc} \times G_{Fa Fb}) \circ (F_{bc} \times F_{ab}) $ };
	\node (M) at (0,0) 	{$G_{GF(a) GF(c)} \circ c^\sB_{Fa Fb Fc} \circ (F_{bc} \times F_{ab})$};
	\node (B) at (0,-1.5) 	{$G_{Fa Fc} \circ F_{a c} \circ c^\sA_{abc}$};
	\draw [->] (T) --  node [right] {$\psi_{Fa Fb Fc} \circ (F_{bc} \times F_{ab})$} (M);
	\draw [->] (M) -- node [right] {$G_{GF(a) GF(c)}  \circ \phi_{a b c}  $} (B);
\end{tikzpicture}
\end{center}
which in components is given by:
\begin{equation*}
GF(g) \circ GF(f) \stackrel{\psi}{\to} G( F(g) \circ F(f)) \stackrel{G\phi}{\to} GF( g \circ f).
\end{equation*}
\end{enumerate}
\end{definition}

\begin{remark} The composition of pointed homomorphisms is again pointed. 
Notice also that the above composition of homomorphisms is strictly associative and hence gives rise to the ordinary category $\catbicat$ whose objects are bicategories and whose morphisms are homomorphisms, see \cite{Benabou67}.
\end{remark}

\begin{definition}[Whiskering] \label{DefnWhiskering}
	Let $\sA, \sB, \sC, \sD$ be bicategories and  $(F, \phi^F): \sA \to \sB$, $(G, \phi^G) : \sB \to \sC$,  $(\bar G, \phi^{\bar G}) : \sB \to \sC$, and  $(H, \phi^H) : \sC \to \sD$ be homomorphisms, and let $(\sigma_x, \sigma_f) : G \to \bar G$ be a transformation. We define the {\em pre-whiskering} of $\sigma$ with $F$ to be the transformation: 
\index{bicategory!whiskering}
	\begin{equation*}
		 \sigma F = ( \sigma_{F(a)}, \sigma_{F(f)}).
\end{equation*}
The {\em post-whiskering} of $\sigma$ with $H$ is similarly defined to be the transformation:
	\begin{equation*}
		 H \sigma = ( H(\sigma_{a}), (\phi^H)^{-1} \circ (H(\sigma_{f})) \circ \phi^H).
\end{equation*}
Let $\theta = (\theta_x, \theta_f): G \to \overline G$ be another transformation and let $\Gamma: \sigma \to \theta$ be a modification. The {\em pre-whiskering} of $\Gamma$ with $F$ is the modification $\Gamma F: \sigma F \to \theta F$ given by components $\Gamma_{F(a)}$. The {\em post-whiskering} of $\Gamma$ with $H$ is the modification $H \Gamma: H \sigma \to H \theta$ given by components $H(\Gamma_x)$. 
\end{definition}

\begin{remark}\label{rmk:strictnessofwhisker}
	Whiskering induces homomorphisms of mapping bicategories. If $F: \cA \to \cB$ and $H:\cC \to \cD$ are homomorphisms of bicategories then pre- and post-whiskering induces homomorphisms
	\begin{align*}
		(-)F: & \bicat(\cB, \cC) \to \bicat(\cA, \cC) \\
		H(-): & \bicat(\cB, \cC) \to \bicat(\cB, \cD).
	\end{align*} 
	Pre-whiskering $(-)F$ is always a strict homomorphism, while post-whiskering is a strict homomomorphism if $H$ is strict.
\end{remark}

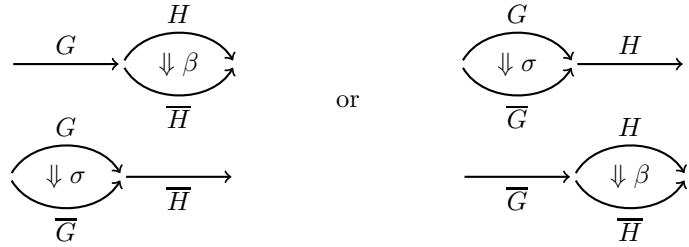
\begin{figure}[ht]
\begin{center}
\begin{tikzpicture}[thick]

\node  [inner sep=1pt] (A) at (0,1.5) {};
\node [inner sep=1pt] (B) at (1.5,1.5) {};
\node  [inner sep=1pt] (C) at (3, 1.5) {};
\node  [inner sep=1pt] (D) at (0,0) {};
\node  [inner sep=1pt] (E) at (1.5,0) {};
\node  [inner sep=1pt] (F) at (3, 0) {};

\node at (2.25, 1.5) {$\Downarrow \beta$};
\node at (.75,0) {$\Downarrow \sigma$};

\draw [->] (A) -- node [above] {$G$} (B);
\draw [->] (B) to [bend left = 60]  node [above] {$H$} (C); 
\draw [->] (B) to [bend right =60] node [below] {$\overline H$} (C); 
\draw [->] (E) -- node [below] {$\overline H$} (F);
\draw [->] (D) to [bend left =60] node [above] {$G$} (E); 
\draw [->] (D) to [bend right = 60] node [below] {$\overline G$} (E);

\node  [inner sep=1pt] (G) at (6,0) {};
\node [inner sep=1pt] (H) at (7.5,0) {};
\node  [inner sep=1pt] (I) at (9, 0) {};
\node  [inner sep=1pt] (J) at (6,1.5) {};
\node  [inner sep=1pt] (K) at (7.5,1.5) {};
\node  [inner sep=1pt] (L) at (9, 1.5) {};

\node at (8.25, 0) {$\Downarrow \beta$};
\node at (6.75,1.5) {$\Downarrow \sigma$};

\draw [->] (G) -- node [below] {$\overline G$} (H);
\draw [->] (H) to [bend left = 60]  node [above] {$H$} (I); 
\draw [->] (H) to [bend right =60] node [below] {$\overline H$} (I); 
\draw [->] (K) -- node [above] {$H$} (L);
\draw [->] (J) to [bend left =60] node [above] {$G$} (K); 
\draw [->] (J) to [bend right = 60] node [below] {$\overline G$} (K); 

\node at (4.5, 1) {or};
\end{tikzpicture}
\caption{Two Compositions in the Tricategory $\bicat$.}
\label{TwoCompositionsinBicat}
\end{center}
\end{figure}

Suppose that $G, \bar G: \sA \to \sB$ and $H, \bar H: \sB \to \sC$ are homomorphisms of bicategories and that $\sigma: G \to \bar G$ and $\beta: H \to \bar H$ are transformations. There are two natural choices for defining the composition of $\sigma$ and $\beta$ as transformations between $H \circ G \to \bar H \circ \bar G$. First we can whisker $\sigma$ by $\overline H$ and then $\beta$ by $G$, these new transformations can then be composed in $\bicat( \sA, \sC)$. Alternatively, we can whisker $\sigma$ by $H$ and $\beta$ by $\overline G$. These transformations can also be composed in  $\bicat( \sA, \sC)$, see Figure~\ref{TwoCompositionsinBicat}. A similar discussion applies to modifications. Generally these operations don't agree, but they yield equivalent tricategories, see \cite{GPS95}. 

\begin{definition}[External Equivalence]
 A homomorphism $F: \sA \to \sB$ is an {\em equivalence of bicategories} \index{bicategory!equivalence of} if there exists a homomorphism $G: \sB \to \sA$ and an equivalence $id_{\sB} \simeq F \circ G$ in the bicategory $\bicat(\sB, \sB)$ and an equivalence $id_{\sA} \simeq G \circ F$ in the bicategory $\bicat(\sA, \sA)$. 
\end{definition}

\begin{theorem}[Whitehead's theorem for bicategories] \label{Whiteheadforbicats} \index{bicategory!Whitehead's theorem} 
\index{Whitehead's theorem!for bicategories}
 A homomorphism $F: \sA \to \sB$ is an  equivalence of bicategories if and only if 
 \begin{enumerate}
\item $F$ induces an isomorphism $\pi_0 \sA \cong \pi_0 \sB$. (Essentially surjective on objects).
\item $F_{ab}: \sA( a, b) \to \sB( Fa, Fb)$ is essentially surjective for all $a,b \in \sA$. (Essentially full on 1-morphism).
\item $F_{ab}: \sA( a, b) \to \sB( Fa, Fb)$ is fully-faithful for all $a,b \in \sA$. (Fully-faithful on 2-morphisms).
\end{enumerate}
\end{theorem}

The proof of this theorem is a routine but tedious application of the axiom of choice. We proved a similar theorem in Section \ref{SectWhiteheadsTheoremforSymmetricMonoidalBicategories} and the proof of the above theorem can be extracted from there. 

\section{Coherence and Strictification}\label{sec:app2-coherencestrict}

The following discussion is based on the exposition in \cite{MacLane71}. A {\em bracketing} consists of a parenthesized word in ``$-$'' (dashes) and ``$1$'' (ones). The {\em length} of a bracketing is the number of dashes. More precisely, we make the following recursive definition: Both $(1)$ and $(-)$ are bracketings (of length zero and one, respectively). If $u$ and $v$ are bracketings, then $(u) (v) $ is a bracketing of length $\text{length}(u) + \text{length}(v)$. The examples below have lengths five and four respectively. 
\begin{equation*}
 ((--)1)(((--)1)(1-)) \qquad \text{ and } \qquad ((((11)1)-)((--)1))-
\end{equation*}

Let $\sA$ be a bicategory. Given a bracketing $b$ of length $n$, and a word $w$ of composable 1-morphisms from $\sA$ of length $n$, we may evaluate it in $\sA$ by replacing (1)'s with identities and using the horizontal composition of 1-morphisms. This gives a canonical 1-morphism in $\sA$, which we denote $b(w)$.
Define the following elementary moves on bracketings:
\begin{enumerate}
\item (Identity Moves) $b \leftrightarrow b$
\item (Deleting Identities) 
$ 1 b  \Leftrightarrow b $ and $b1 \Leftrightarrow b$
\item (Rebracketing) $(b b') b'' \Leftrightarrow b ( b' b'')$
\end{enumerate}
where $b, b', b''$ represent arbitrary bracketings. An elementary move has an apparent source and target which are bracketings. 
A {\em path} also has source and target bracketings and is defined recursively as follows: 
\begin{enumerate}
\item Elementary moves are paths,
\item If $p$ and $p'$ are paths such that the source of $p$ is the target of $p'$, then $p \circ p'$ is a path with source the source of $p'$ and target the target of $p$.
\item if $p$ and $p'$ are paths, then $(p)(p')$ is a path with source $(\text{source}(p))( \text{source}(p'))$ and target $(\text{target}(p))( \text{target}(p'))$.
\end{enumerate}

Given two bracketings $b$ and $b'$ of length $n$ and a word $w$ of composable 1-morphisms from $\sA$, also of length $n$, we have the two 1-morphisms  $b(w)$ and $b'(w)$ of $\sA$.  Given a path $s$ starting at $b$ and ending at $b'$, we get a canonical 2-morphism of $\sA$,  $s(w): b(w) \to b(w')$, given by replacing the elementary moves by the appropriate unitors and associators from $\sA$. 

\begin{theorem}[MacLane's Coherence Theorem \cite{MacLane71}] \label{MacLanesCoherenceTheorem} \index{MacLane's coherence theorem} 
\index{bicategory!MacLane's coherence theorem}
If $b$ and $b'$ are two bracketings of length $n$ and $w$ is a word of composable 1-morphism from the bicategory $\sA$, also of length $n$, then there exists a path $s: b \to b'$ starting at $b$ and ending at $b'$. Moreover, given any two such paths $s, s': b \to b'$, the resulting 2-morphisms  $s(w): b(w) \to b(w')$ and  $s'(w): b(w) \to b(w')$ are identical. 
\end{theorem}

MacLane's coherence theorem allows us to form the following construction. Define the {\em standard bracketing} $b^n_\text{std}$ of length $n$ to be the bracketing:
\begin{equation*}
(( \cdots (((- -) -)-) \cdots ) -).
\end{equation*}
Given a bicategory $\sA$ we define the following associated 2-category, $\st(A)$:
\begin{itemize}
\item The objects of $\st(\sA)$ are the same as those of $\sA$. 
\item The 1-morphisms of $\st(\sA)$ are the words of composable 1-morphisms in $\sA$. 
\item Given two such words $w, w'$, the 2-morphisms are the 2-morphisms between $b^{|w|}_\text{std}(w)$ and $b^{|w'|}_\text{std}(w')$. 
\end{itemize}
Vertical composition coincides with that in $\sA$. The horizontal composition of 1-morphisms is given by concatenating words, and hence is strictly associative. The horizontal composition of 2-morphisms is given by using the canonical rebracketing 2-morphisms supplied by MacLane's coherence theorem. Given two 2-morphisms $\alpha: b_\text{std}(w) \to b'_\text{std}(w')$ and $\beta:{b}^{|\tilde w|}_\text{std}(\tilde w) \to {b}^{|\tilde w'|}_\text{std}(\tilde{w}')$, their horiszontal composition is defined to be the composite:
\begin{equation*}
b^{| w \tilde w|}_\text{std}( w \tilde w) \to b^{|w|}_\text{std}(w) \circ {b}^{|\tilde w|}_\text{std}(\tilde w) \stackrel{ \alpha * \beta }{\to} b^{|w'|}_\text{std}(w') \circ{b}^{|\tilde w'|}_\text{std}(\tilde{w}') \to b^{| w' \tilde w'|}_\text{std}( w' \tilde{w}')
\end{equation*}
where the unlabeled arrows are the canonical 2-morphisms from MacLane's theorem. 

This composition is automatically associative and hence defines a 2-category $\st(\sA)$, called the {\em strictification} of $\sA$. In fact $\st$ defines a functor,
\begin{equation*}
	\st: \catbicat \to 2\Cat_s
\end{equation*}
where $2\Cat_s$ is the category of 2-categories and strict homomorphisms. 
%This is left-adjoint to the inclusion functor $i: 2\Cat_s \to \catbicat$ into bicategories and (weak) homomorphisms.
Morevoer the assignment $f \mapsto (f)$ taking a 1-morphism in $\sA$ to the word of length one in $st(\sA)$ extends to a homomorphism of bicategories $\sA \to st(\sA)$ which is well-known to be an equivalence (see \cite[pg.~30]{JS93} for the one-object case). These constructions are revisted in Section~\ref{sec:coherenceinaction}.

%In particular the unit and co-unit of the adjunction give us for each 2-category $\sB$ a strict homomorphism $Q(\sB) \to \sB$ (which we can take to be the homomorphism which sends a composable word of 1-morphisms in $\sB$ to its composition) and for each bicategory $\sA$ a weak homomorphism $\sA \to Q(\sA)$ (which we can take to be the homomorphism which sends a 1-morphism in $\sA$ to the singleton word consisting of exactly that 1-morphism). Both of these homomorphisms are equivalences of bicategories. Note, however, that for a 2-category $\sB$, the inverse equivalence to the canonical strict homomorphism $Q(\sB) \to \sB$ is nearly always a {\em weak} homomorphism.

\begin{corollary}
Every bicategory is equivalent to a 2-category. 
\end{corollary}

%Via the faithful inclusion of $\Cat{2cat}_s \to \Cat{2cat}_w$, we can also consider $Q$ as a functor $Q:\Cat{bicat} \to \Cat{2cat}_w$. Then we have an adjunction,
%\begin{equation*}
%	i:  \Cat{2cat}_w \leftrightarrows \Cat{bicat}: Q
%\end{equation*}
%The unit and counit of this adjunction provide canonical inverse equivalences to the unit and counit of the previous adjunction. 

\section{Symmetric, Braided and Monoidal Categories} \label{MonoidalCatsAsBicats}

Symmetric monoidal categories are pervasive in mathematics and historically played a vital role in the development of the theory of bicategories. It is no surprise that many theorems about monoidal categories have analogs for bicategories as well. The presentation given here is designed to emphasize this connection.

\begin{definition}
A {\em monoidal category} is a bicategory with one object. 
\index{monoidal category}
\end{definition}

Thus a monoidal category consists of a category $M$ (the category of 1-morphisms, and 2-morphisms) with a specfied 1-morphism $1 \in M$ and a horizontal composition functor, 
\begin{equation*}
c = \otimes: M \times M \to M,
\end{equation*}
together with associator and unitor natural isomorphism satisfying the pentagon and triangle identities. We will typically use the notation $(M, \otimes)$ for monoidal categories. The reader is invited to verify that the above definitions of monoidal category, monoidal functor, and monoidal natural transformation coincide with the standard definitions, as presented in \cite{MacLane71}, for example.

One is now tempted to define morphisms and higher morphisms for monoidal categories by simply using those which already exist for bicategories, however, this does not recover the usual notion of monoidal functor nor monoidal natural transformation. Indeed, there would also be an additional categorical layer coming from the modifications. The solution to the problem is to realize that,  as a bicategory, a monoidal category is canonically pointed, and thus it is natural to define morphisms of monidal categories to be {\em pointed} morphisms of bicategories. 

\begin{definition}
A {\em monoidal  functor} from the monoidal category $(M, \otimes)$ to $(M', \otimes')$ is a pointed homomorphism between the corresponding bicategories. A  {\em monoidal natural transformations} is a pointed transformation. 
\end{definition}

All pointed modifications between monoidal categories are automatically trivial, so that the bicategory $\bicat_*( (M, \otimes), (M', \otimes'))$ is in fact just an ordinary category. This agrees with the usual notion of monoidal functor between monoidal categories and monoidal natural transformation between these.

There is an obvious notational danger with defining monoidal categories to be  bicategories with a single object. Traditionally a monoidal category is thought of as a category with extra structure, and hence has objects and morphisms. However viewed as a bicategory, the objects become 1-morphisms, and the morphisms become 2-morphisms. Unfortunately this point of confusion is inevitable. We will try to stick with the terminology from the bicategorical perspective as much as possible. 

Recall the following well known construction: Algebras form the objects of a bicategory, $\alg^2$, whose 1-morphisms are bimodules and whose 2-morphisms are bimodule maps. The horizontal composition is given by the tensor product of bimodules. Given an algebra $A$, we can recover the center of $A$ as the algebra of bimodule endomorphism of the identity 1-morphism ${}_A A_A$, i.e., $Z(A) = {}_A \hom_A ( A, A)$. Similarly in an arbitrary bicategory, $\sB$, one can define the center $Z(X)$ of an object $X$ to be the commutative monoid $ \sB( I_X, I_X)$. The following is a higher categorical analog of this construction.

\begin{definition}
Let $(M, \otimes)$ be a monoidal category. Let $id_M$ be the identity homomorphism from $M$ to $M$. Define the 
\index{monoidal category!center}
\index{center of monoidal category}
{\em center} $Z(M)$ of $M$ to be the full sub-bicategory of $\bicat(M, M)$ whose only object is $id_M$. That is $Z(M)$ is the bicategory whose only object is $id_M$, whose morphisms are the (non-pointed) transformations $id_M \to id_M$, and whose 2-morphisms are the modifications between these. The 1-morphisms of $Z(M)$ are called {\em half-braidings}.  
\end{definition}

Half-braidings were introduced by Izumi in the context of subfactors \cite{Izumi00}. By construction $Z(M)$ is a monoidal category. Notice that we did not use {\em pointed} transformations and modifications. Unpacking this definition, we see that a half-braiding 
\index{monoidal category!half-braiding} 
\index{half-braiding} 
consists of:
\begin{enumerate}
\item a 1-morphism $X \in M$,
\item a family of invertible 2-morphisms $\gamma_{Y, X}: Y \otimes X \to X \otimes Y$ where $Y \in M$ runs over all 1-morphisms.
\end{enumerate}
These data are required to satisfy:
\begin{enumerate}
\item  [(HBR0)] (Naturality) The $\gamma_{-, X}$ form a natural isomorphism of functors:
\begin{equation*}
	\gamma_{-, X}: (-) \otimes X \to X \otimes (-).
\end{equation*}
\item   [(HBR1)] For all $Y, Z$, the following diagram commutes:
\begin{center}
\begin{tikzpicture}[thick]
	\node (LT) at (0,4) 	{$(Z \otimes Y) \otimes X$ };
	\node (LM) at (-1,2) 	{$Z \otimes (Y \otimes X)$ };
	\node (LB) at (0,0) 	{$Z \otimes (X \otimes Y)$};
	
	\node (RT) at (4,4) 	{$X \otimes (Z \otimes Y)$};
	\node (RM) at (5,2) 	{$(X \otimes Z) \otimes Y$ };
	\node (RB) at (4,0)	{$(Z \otimes X) \otimes Y$};
	\draw [->] (LM) --  node [left] {$1 \otimes \gamma_{Y, X}$} (LB);
	\draw [<-] (LM) --  node [left] {$a$} (LT);
	
	\draw [->] (LT) -- node [above] {$\gamma_{Z \otimes Y, X}$} (RT);
	\draw [->] (LB) -- node [below] {$a^{-1}$} (RB);
	\draw [<-] (RT) -- node [right] {$a$} (RM);
	\draw [->] (RB) -- node [right] {$ \gamma_{Z, X} \otimes 1$} (RM);
	
\end{tikzpicture}
\end{center}
\item  [(HBR2)] The following diagram commutes:
\begin{center}
\begin{tikzpicture}[thick]
	\node (LT) at (0,0) 	{$1 \otimes X$ };
	\node (MT) at (2,0) 	{$X$};
	\node (RT) at (4,0) 	{$X \otimes 1$};
%	\node (RB) at (2,0)	{$$};
	\draw [->] (LT) --  node [above] {$\ell$} (MT);
	\draw [->] (MT) -- node [above] {$r^{-1}$} (RT);
	\draw [->] (LT) to [out= -60, in = -120] node [below] {$\gamma_{1, X}$} (RB);
\end{tikzpicture}
\end{center}
\end{enumerate}

A 2-morphism between half-braidings (say from $(X, \{ \gamma_{-, X} \})$ to $(X', \{ \gamma_{-, X'} \})$) consists of a 2-morphism $\Gamma: X \to X'$ in $M$ such that for all $Y \in M$ the following diagram commutes:
\begin{center}
\begin{tikzpicture}[thick]
	\node (LT) at (0,1.5) 	{$Y \otimes X$ };
	\node (LB) at (0,0) 	{$X \otimes Y$};
	\node (RT) at (3,1.5) 	{$Y \otimes X'$};
	\node (RB) at (3,0)	{$X' \otimes Y$};
	\draw [->] (LT) --  node [left] {$\gamma_{Y, X}$} (LB);
	\draw [->] (LT) -- node [above] {$1 * \Gamma$} (RT);
	\draw [->] (RT) -- node [right] {$\gamma_{Y, X'}$} (RB);
	\draw [->] (LB) -- node [below] {$\Gamma * 1$} (RB);
\end{tikzpicture}
\end{center}

There is a strict forgetful homomorphism $U: Z(M) \to M$ which sends a half-braiding,  $(X, \{ \gamma_{-, X} \})$,  to its underlying 1-morphism, $X$.

\begin{definition}
Let $M$ be a monoidal category. A {\em braiding} 
\index{monoidal category!braided}
\index{braided monoidal category}
 for $M$ is a strict section of the homomorphism $U: Z(M) \to M$. That is, a strict homomorphism $s: M \to Z(M)$ such that $U \circ s = 1_M$ (strict equality). 
\end{definition}

%\begin{remark}
%The assignment $M \mapsto Z(M)$ is functorial. Given a homomorphism of monoidal categories $(F, \phi): M \to N$, we have a homomorphism $Z(F): Z(M) \to Z(N)$ which takes a half-braiding $(X, \{ \gamma_{-, X} \})$ to the half-braiding $(FX, \phi_{X, -}^{-1} \circ \gamma_{-, X} \circ \phi_{-, X})$. \notetoself{Verify this remark...}
%\end{remark}

%\begin{remark}
%Suppose that $M$ and $N$ are monoidal categories, which are equivalent as monoidal categories and that $M$ is equipped with a braiding. Then there is a canonical braiding on $N$ making the equivalence of monoidal categories into a braided monoidal equivalence. 
%\end{remark}

Thus a braiding for a monoidal category $M$ consists of an assignment: for each $X$ a family of 2-morphisims $\gamma_{Y, X}$ which are natural in the $Y$s, which satisfy several conditions. Among these conditions is the requirement that this assignment must be a homomorphism, and so must also be natural in the $X$ variable. More precisely, we have the following characterization:

\begin{proposition}
Let $M$ be a monoidal category. A braiding for $M$ consists of 2-morphisms $\gamma_{X, Y}: X \otimes Y \to Y \otimes X$ for each pair of 1-morphisms $X, Y$,  such that:
\begin{itemize}
\item [(BR0)] The $\gamma$ form a natural transformation, 
\begin{equation*}
\gamma: \otimes \to \otimes \circ \tau
\end{equation*}
where $\tau: M \times M \to M \times M$ denotes the flip functor,
\item [(BR1)]  For all $X, Y, Z$, the following diagrams commute:
\begin{center}
\begin{tikzpicture}[thick]
	\node (LT) at (0,4) 	{$(Z \otimes Y) \otimes X$ };
	\node (LM) at (-1,2) 	{$Z \otimes (Y \otimes X)$ };
	\node (LB) at (0,0) 	{$Z \otimes (X \otimes Y)$};
	
	\node (RT) at (4,4) 	{$X \otimes (Z \otimes Y)$};
	\node (RM) at (5,2) 	{$(X \otimes Z) \otimes Y$ };
	\node (RB) at (4,0)	{$(Z \otimes X) \otimes Y$};
	\draw [->] (LM) --  node [left] {$1 \otimes \gamma_{Y, X}$} (LB);
	\draw [<-] (LM) --  node [left] {$a$} (LT);
	
	\draw [->] (LT) -- node [above] {$\gamma_{Z \otimes Y, X}$} (RT);
	\draw [->] (LB) -- node [below] {$a^{-1}$} (RB);
	\draw [<-] (RT) -- node [right] {$a$} (RM);
	\draw [->] (RB) -- node [right] {$ \gamma_{Z, X} \otimes 1$} (RM);
	
\end{tikzpicture}
\end{center}
\item [(BR2)] For all $X, Y, Z \in M$ the following diagram commutes: 
\begin{center}
\begin{tikzpicture}[thick]
	\node (LT) at (0,4) 	{$Z \otimes (Y \otimes X)$ };
	\node (LM) at (-1,2) 	{$(Z \otimes Y) \otimes X$ };
	\node (LB) at (0,0) 	{$(Y \otimes Z) \otimes X$};
	
	\node (RT) at (4,4) 	{$(Y \otimes X) \otimes Z$};
	\node (RM) at (5,2) 	{$Y \otimes (X \otimes Z)$ };
	\node (RB) at (4,0)	{$Y \otimes (Z \otimes X)$};
	\draw [->] (LM) --  node [left] {$\gamma_{Z, Y} \otimes 1$} (LB);
	\draw [<-] (LM) --  node [left] {$a^{-1}$} (LT);
	
	\draw [->] (LT) -- node [above] {$\gamma_{Z, Y \otimes X}$} (RT);
	\draw [->] (LB) -- node [below] {$a$} (RB);
	\draw [<-] (RT) -- node [right] {$a^{-1}$} (RM);
	\draw [->] (RB) -- node [right] {$1 \otimes \gamma_{Z,X}$} (RM);
	
\end{tikzpicture}
\end{center}
\end{itemize}
\end{proposition}

\begin{proof}
The necessity of (BR0) has already been established. (BR1) is just (HBR1) from the definition of a half-braiding. (BR2) is precisely the expression that the section $s: M \to Z(M)$ is a strict homomorphism so that $s(Y \otimes X) = s(Y) \otimes s(X)$. Compare with Definition \ref{defncompositonoftransformations}. Given (BR0-BR2), then (HBR2) is automatically satisfied as is the condition that $s(1) = 1$, see \cite{JS86}. %Joyal-Street Braided Monoidal Categories.
\end{proof}

\begin{definition}
A {\em braided monoidal category} is a monoidal category equipped with a braiding. A monoidal functor $(F, \phi)$ between braided monoidal categories is a {\em braided monoidal functor} if it satisfies 
\begin{center}
\begin{tikzpicture}[thick]
	\node (LT) at (0,1.5) 	{$F(X) \otimes F(Y)$ };
	\node (LB) at (0,0) 	{$F(X \otimes Y)$};
	\node (RT) at (4,1.5) 	{$F(Y) \otimes F(X)$};
	\node (RB) at (4,0)	{$F(Y \otimes X)$};
	\draw [->] (LT) --  node [left] {$\phi_{A,B}$} (LB);
	\draw [->] (LT) -- node [above] {$\gamma_{FX, FY}$} (RT);
	\draw [->] (RT) -- node [right] {$\phi_{B,A}$} (RB);
	\draw [->] (LB) -- node [below] {$F(\gamma_{X, Y})$} (RB);
\end{tikzpicture}.
\end{center}
A braided monoidal category is {\em symmetric} 
\index{monoidal category!symmetric}
\index{symmetric monoidal category} if the braiding satisfies $\gamma_{X, Y} \circ \gamma_{Y, X} = 1_{Y \otimes X}$.
\end{definition}

\section{Pastings, Strings, Adjoints, and Mates} \label{PastingsStringsAdjointsAndMatesSection}

Bicategories have two kinds of compositions, horizontal and a vertical, and this makes them inherently two-dimensional entities. For this reason many properties, structures, and equations are best expressed in a corresponding 2-dimensional formalism. In ordinary category theory one might specify certain structures on a category (for example, in 1-dimensional terms,  the functor $\otimes$ and natural transformations $a, \ell, r$ of a monoidal category). These structures can then be combined in a one-dimensional, linear fashion. Often one imposes conditions on such a structure such as requiring that two or more different linear compositions must give the same result. Typical examples are the pentagon and triangle identities for monoidal categories. 

For bicategories, due to their two-dimensional nature, structures and properties are best expressed in two-dimensional terms. Structures for bicategories, such as those that would give us a symmetric monoidal bicategory (the subject of Chapter \ref{SymMonBicatChapt}), often require that two or more different two-dimensional compositions of morphisms must be the same. There are two common equivalent notations for describing these two-dimensional compositions: pasting diagrams and string diagrams. We will describe both here. 

Pasting diagrams 
\index{pasting diagram} 
\index{bicategory!pasting diagram} are a natural extension of the sort of diagrams commonly drawn in ordinary category theory. Typical diagrams in ordinary category theory consist of some vertices (labeled by the objects of the category being considered) together with directed edges (which are labeled by morphisms of the category). A prototypical example is the commutative square:
\begin{center}
\begin{tikzpicture}[thick]
	\node (LT) at (0,1.5) 	{$a$ };
	\node (LB) at (0,0) 	{$b$};
	\node (RT) at (2,1.5) 	{$c$};
	\node (RB) at (2,0)	{$d$};
	\draw [->] (LT) --  node [left] {$f$} (LB);
	\draw [->] (LT) -- node [above] {$$g} (RT);
	\draw [->] (RT) -- node [right] {$h$} (RB);
	\draw [->] (LB) -- node [below] {$k$} (RB);
\end{tikzpicture}
\end{center}
In this diagram there are two linear compositions, $a \stackrel{f}{\to} b \stackrel{k}{\to} {d}$ and $a \stackrel{g}{\to} c \stackrel{h}{\to} {d}$, and when we assert that the diagram commutes, we are asserting that these two linear compositions are equal. 

\begin{figure}[ht]
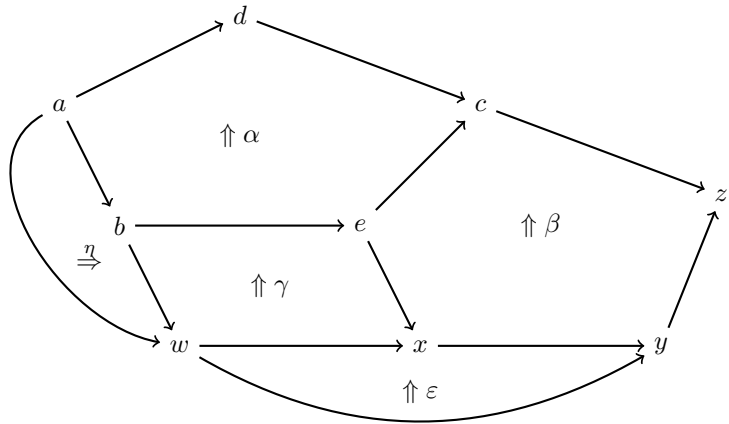

\begin{center}
%
\caption{An Example of a Pasting Diagram}
\label{ExamplePastingDiagramFig}
\end{center}
\end{figure}

Pasting diagrams are similar. They consist of polygonal arrangements in the plane with appropriate labels. Such an arrangement is a two-dimensional analogue of the linear arrangements for ordinary categories. The vertices of a pasting diagram are labeled with objects of the bicategory under consideration, the (directed) edges are labeled with 1-morphisms, and the polygonal regions themselves are labeled with 2-morphisms (and directions). The basic situations are the following:
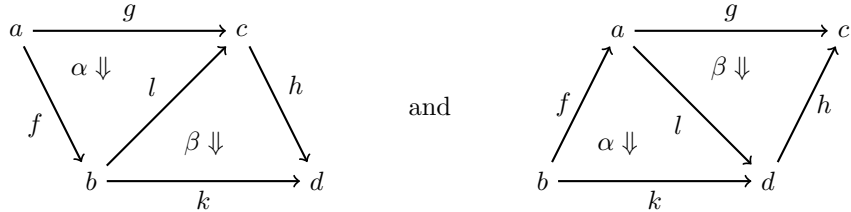
\begin{center}

\end{center}
which are to be interpreted as the compositions $(\beta * id_f) \circ (id_h * \alpha)$ and $(\alpha * id_h) \circ ( id_f * \beta)$. Out of these basic composites we can give meaning to more complex diagrams such as the one in Figure~\ref{ExamplePastingDiagramFig}.

In a 2-category, there is a unique 2-morphism specified by each such pasting diagram. This is not the case in general bicategories. Such a diagram is supposed to represent the vertical composition of 2-morphisms which are horizontal compositions like:
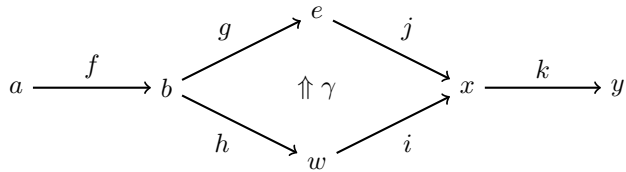
\begin{center}
%
\end{center}
but to make sense of such horizontal compositions, one needs to choose a bracketing. Does the above picture mean $(id_k * \gamma) * id_f$ or $id_k * (\gamma * id_f)$?

Worse, to make sense of the entire pasting diagram we must also change the bracketing {\em en route}. The above horizontal composition is a 2-morphism which starts at (say) $ (k \circ (i \circ h)) \circ f$, but in the pasting diagram in Figure~\ref{ExamplePastingDiagramFig} we must compose it with $\varepsilon * \eta$, which ends at $(k \circ i) \circ (h \circ f)$. In a general bicategory these compositions might not agree and so the 2-morphisms may not be composable. 

Nevertheless, in any bicategory MacLane's Coherence Theorem (Theorem~\ref{MacLanesCoherenceTheorem}) ensures that there are canonical coherence isomorphisms between rebracketed expsessions. So for example, we have a canonical isomorphism  $(k \circ i) \circ (h \circ f) \cong (k \circ (i \circ h)) \circ f$. When we draw a pasting diagram we will implicitly insert these canonical coherence isomorphisms. With this convention, a pasting diagram in an arbitrary bicategory takes a unique value, provided we specify a bracketing of the {\em outside} 1-morphisms. This can be proven by an appropriate induction on polygonal decompositions of the disk.  

When we write equations of the form ``pasting diagram A = pasting diagram B'', what we mean is that when we give these pasting diagrams a fixed bracketing on the boundary, then they agree. If this is the case, then they also agree with any other bracketing on the boundary, so that nothing is lost by omitting this bracketing. The property ``A = B''  doesn't depend on the bracketing of the relevant pasting diagram, so we may as well omit it from the notation. 

This becomes slightly problematic when we want to specify {\em structure} for a bicategory. For example, a monoidal bicategory is similar to a monoidal category. Both have $\otimes$-structures, and both have associators $\alpha$. But while these structures for a monoidal category satisfy the pentagon identity, in a monoidal bicategory this identity becomes an extra piece of data: the pentagonator. This is a 2-isomorphism which is supposed to fill diagrams such as Figure~\ref{fig:pentagonator}.
For this to strictly make sense, we must first bracket the bottom sequence of 1-morphisms. 
\begin{figure}[ht]
	  \begin{center}
\begin{tikzpicture}[thick]
	%\node at (0, 4.5) {\underline{The Pentagonator}};
	\node (LT) at (0,2) 	{$((a \otimes b)\otimes c) \otimes d $ };
	\node (LB) at (1,0) 	{$(a \otimes (b \otimes c)) \otimes d$};
	\node (RT) at (8,2) 	{$a \otimes(b \otimes(c \otimes d))$};
	\node (MT) at (4,3.5) {$(a \otimes b) \otimes (c \otimes d)$};
	\node (RB) at (7,0)	{$ a \otimes ((b \otimes c) \otimes d)$};
	\draw [->] (LT) --  node [left] {$\alpha \otimes I $} (LB);
	\draw [->] (LT) -- node [above left] {$\alpha $} (MT);
	\draw [->] (MT) -- node [above right] {$\alpha $} (RT);
	\draw [->] (RB) -- node [right] {$I \otimes \alpha $} (RT);
	\draw [->, in = 200, out = 340] (LB) to node [below] {$\alpha $} (RB);
	\node at (4,1 ) {$\Uparrow \pi$};
\end{tikzpicture}
\end{center}
	\caption{Pentagonator}
	\label{fig:pentagonator}
\end{figure}

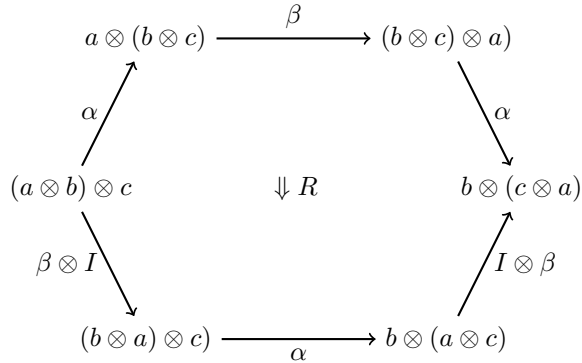
\begin{figure}[ht]
	\begin{center}
	\begin{tikzpicture}[thick]
		\node (LT) at (0,4) 	{$a \otimes (b \otimes c)$ };
		\node (LM) at (-1,2) 	{$(a \otimes b) \otimes c$ };
		\node (LB) at (0,0) 	{$(b \otimes a) \otimes c)$};

		\node (RT) at (4,4) 	{$(b \otimes c) \otimes a)$};
		\node (RM) at (5,2) 	{$b \otimes (c \otimes a)$ };
		\node (RB) at (4,0)	{$b \otimes (a \otimes c)$};
		\draw [->] (LM) --  node [left] {$\beta \otimes I$} (LB);
		\draw [->] (LM) --  node [left] {$\alpha$} (LT);	
		\draw [->] (LT) -- node [above] {$\beta$} (RT);
		\draw [->] (LB) -- node [below] {$\alpha$} (RB);
		\draw [->] (RT) -- node [right] {$\alpha$} (RM);
		\draw [->] (RB) -- node [right] {$ I \otimes \beta$} (RM);
		\node at (2,2) {$ \Downarrow R$};
	\end{tikzpicture}
	\end{center}
	\caption{Ambiguity in Pasting Diagrams}
	\label{fig:R2Morphism}	
\end{figure}

Alternatively, we can interpret the pentagonator as a {\em family} of 2-isomorphisms, one for each possible way of bracketing the boundary (in this case two), but require that this family is {\em coherent}, i.e., when we change the bracketing (via MacLane's canonical coherence isomorphisms) of one member of this family, we get another member of this family. This is the point of view we will adopt here. So for example the data in Figure~\ref{fig:R2Morphism} represents a family of four 2-morphisms, one for each of the four ways of bracketing the boundary. 

String diagrams 
\index{bicategory!string diagram} 
\index{string diagram} 
are an equivalent two-dimensional notation which is in some sense dual to pasting diagrams. In a string diagram an object is represented, not as a vertex, but as a {\em region}. A 1-morphism is drawn as a boundary between regions, and a 2-morphism is drawn as a vertex where several 1-morphisms meet. Usually these vertices are drawn as nodes, with the name of the 2-morphism labeling the node. 
\begin{figure}[ht]
\begin{center}
\begin{tikzpicture}[thick]

\begin{scope}
		\node (A) at (0,0) {$a$};
	\node (BT) at (1.5,1) {$b$};
	\node (BB) at (1.5,-1) {$d$};
	\node (C) at (3,0) {$c$};
	\draw [->] (A) to node [above left] {$f$} (BT);
	\draw [->] (A) to node [below left] {$f'$} (BB);
	\draw [->] (BT) to node [above right] {$g$} (C);
	\draw [->] (BB) to node [below right] {$g'$} (C);
	\node at (1.5,0) {$\Downarrow u$};
\end{scope}
	
	\node at (4.5,0) {vs};
	
\begin{scope}[xshift = 6 cm, yshift = -1 cm]
	\node (LT) at (0,2) 	{ };
	\node (LB) at (0,0) 	{};
	\node (RT) at (3,2) 	{};
	\node (RB) at (3,0)	{};
	
	\node (f) [label=above:$f$] at (1,2) {};    
	\node (g) [label=above:$g$] at (2,2) {};
	\node (fa) [label=below:$f'$] at (1, 0) {};
	\node (ga)[label=below:$g'$]  at (2,0) {};
	
	\node [stringnode, minimum width = 2 cm] (u) at (1.5,1) {$u$};
	
	\draw (f.center) -- (f |- u.north);
	\draw (g.center) -- (g |- u.north);
	\draw (fa.center) -- (fa |- u.south);
	\draw (ga.center) -- (ga |- u.south);
	
	%\fill [red] (fa) -- (f) -- (LT) -- (LB) -- (fa);
	%\pattern [pattern color = blue, pattern = dots] (LB) rectangle (f);
	
	\begin{pgfonlayer}{background} 
		\fill [blue!10] (LB) rectangle (f);
		\pattern [pattern color = blue, pattern = dots] (LB) rectangle (f);
		
		\fill [red!10] (fa) rectangle (u.south -| ga);
		\pattern [pattern color = red, pattern = north west lines] (fa) rectangle (u.south -| ga);
		
		\fill [green!20] (f) rectangle (g |- u.north);
		\pattern [pattern color = green, pattern = vertical lines] (f) rectangle (g |- u.north);
		
		\fill [black!20] (g) rectangle (RB);
	\end{pgfonlayer} 
\end{scope}	
\end{tikzpicture}
\caption{Pasting Diagrams vs String Diagrams}
\label{PastingvsStringFig}
\end{center}
\end{figure}
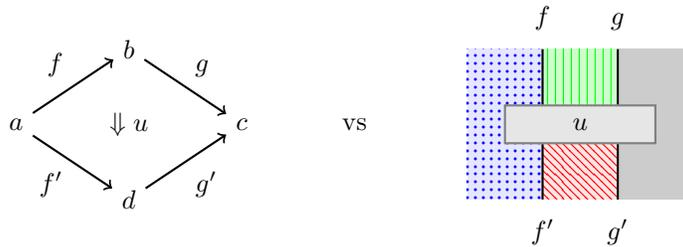
Figure~\ref{PastingvsStringFig} shows how a typical 2-morphism would be encoded in both the pasting diagram notation and the string diagram notation. Notice that in the string diagram we have labeled the regions corresponding to different objects with different colors; we have 
$a =$ \tikz[baseline= 2pt]{ \fill [blue!10] (0,0) rectangle (20pt,10pt); \pattern [pattern color = blue, pattern = dots] (0,0) rectangle (20pt,10pt);}, 
$b =$ \tikz[baseline= 2pt]{ \fill [green!20] (0,0) rectangle (20pt,10pt); \pattern [pattern color = green, pattern = vertical lines] (0,0) rectangle (20pt,10pt);}, 
$c =$ \tikz[baseline= 2pt]{ \fill [black!20] (0,0) rectangle (20pt,10pt);}, and 
$d =$  \tikz[baseline= 2pt]{ \fill [red!10] (0,0) rectangle (20pt,10pt); \pattern [pattern color = red, pattern = north west lines] (0,0) rectangle (20pt,10pt);}. Figure~\ref{StringDiagramExampleFigure} shows the string diagram which is equivalent to the pasting diagram of Figure~\ref{ExamplePastingDiagramFig}. Unless otherwise stated, we will read our string diagrams from left to right and from top to bottom, as in both these figures.  
\begin{figure}[ht]
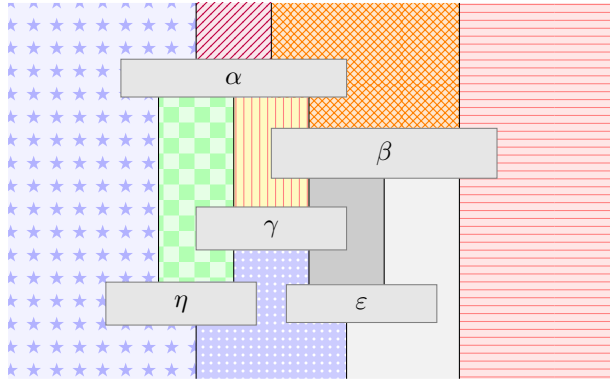

\begin{center}

\caption{An Example of a String Diagram}
\label{StringDiagramExampleFigure}
\end{center}
\end{figure}

Duality in bicategories (also called adjunction) 
\index{bicategory!adjunction/duality in a} 
\index{adjunction} 
\index{adjunction!in a bicategory} 
\index{duality}is an important concept which we will now turn to. String diagrams are particularly well suited for describing duality in bicategories, as we shall see shortly.
\begin{definition}
	An {\em adjunction} in a bicategory $\sB$ is a quadruple $(f, g, \eta, \varepsilon)$ where $f: a \to b$ and $g: b \to a$ are 1-morphisms in $\sB$, $\eta: I_a \to g \circ f$ and $\varepsilon: f \circ g \to I_b$ are 2-morphisms, such that the composites:
\begin{align*}
	f \stackrel{r^{-1}}{\to} f \circ I_a \stackrel{id * \eta}{\to} f \circ (g \circ f) \stackrel{a^{-1}}{\to} (f \circ g) \circ f \stackrel{\varepsilon * id}{\to} I_b \circ f \stackrel{\ell}{\to} f  \\
	g \stackrel{\ell^{-1}}{\to} I_a \circ g \stackrel{\eta * id}{\to} (g \circ f) \circ g \stackrel{a}{\to} g \circ (f \circ g) \stackrel{id * \varepsilon}{\to} g \circ I_b \stackrel{r}{\to} g
\end{align*}
are identities. We say that $f$ is {\em left-adjoint} (or {\em left-dual}) to $g$ and that $g$  is {\em right-adjoint} (or {\em right-dual}) to $f$. We will also denote an adjunction by $f: a \rightleftarrows b: g$
\end{definition}

\begin{figure}[ht]
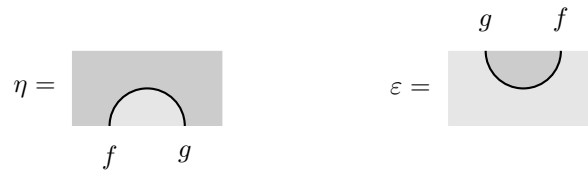

\begin{center}

\caption{Duality via String Diagrams}
\label{DualityViaStringDiagrams}
\end{center}
\end{figure}

\begin{figure}[ht]
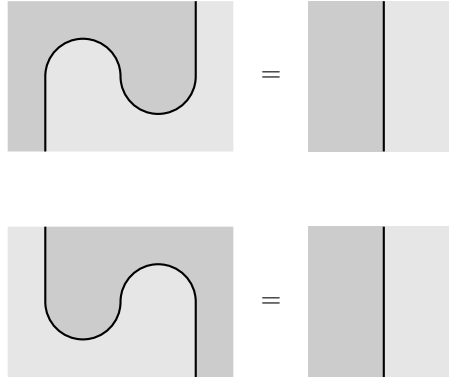

\begin{center}
%
\caption{Adjunction Equations via String Diagrams}
\label{AdjunctionStringDiagramFig}
\end{center}
\end{figure}

When the bicategory $\sB = \cat$, the bicategory of categories, functors and natural transformations, then this is equivalent to the usual notion of adjoint functors, see \cite{MacLane71}. We can translate this definition into string diagrams as follows. The 2-morphisms $\eta$ and $\varepsilon$ can be described by string diagrams as in Figure  \ref{DualityViaStringDiagrams}, and then the above equations become the equations of Figure~\ref{AdjunctionStringDiagramFig}, justifying our choice of graphical depiction.

\begin{definition} \label{DefAdjointEquiv}
An {\em adjoint equivalence} 
\index{bicategory!adjoint equivalence} 
\index{adjoint equivalence} is an adjunction $(f, g, \eta, \varepsilon)$ in which the 2-morphisms $\eta$ and $\varepsilon$ are invertible. 
\end{definition}

If $f$ is left adjoint to $g$,  with adjunction data $(f, g, \eta, \varepsilon)$, and this adjunction is an adjoint equivalence, then $(g, f, \varepsilon^{-1}, \eta^{-1})$ is an adjunction realizing $f$ as a right adjoint to $g$.

\begin{proposition}
 The following notions are logically equivalent for any 1-morphism $f: a \to b$ in a bicategory: 
 \begin{enumerate}
\item $f$ is an equivalence.
\item $f$ is part of an adjoint equivalence $(f, g, \eta, \varepsilon)$.
\end{enumerate}
\end{proposition}

\begin{proof}
The implication (2) $\Rightarrow$ (1) is trivial. For the reverse implication, suppose that $f$ is an equivalence so that there exists $g: b \to a$ and isomorphisms $\eta: I_a \to g \circ f$ and $\overline \varepsilon: f \circ g \to I_b$. If these satisfy the adjunction equations, then we are done. Otherwise the composition, 
\begin{equation*}
f \stackrel{r^{-1}}{\to} f \circ I_a \stackrel{id * \eta}{\to} f \circ (g \circ f) \stackrel{a^{-1}}{\to} (f \circ g) \circ f \stackrel{\varepsilon * id}{\to} I_b \circ f \stackrel{\ell}{\to} f 
\end{equation*}
defines a non-identity isomorphism of $f$. Call this isomorphism $\alpha$. Let $\varepsilon = \overline \varepsilon \circ (\alpha^{-1} * id_g) $. Then $\varepsilon: f \circ g \to I_b$ is isomorphism and the data $(f,g, \eta, \varepsilon)$ forms the desired adjoint equivalence. 
\end{proof}

Fix a 1-morphism $f: a \to b$ in $\sB$. The adjunctions $(f, g, \eta, \varepsilon)$  in $\sB$ become a category where a morphism from  $(f, g, \eta, \varepsilon)$  to $(f, g', \eta', \varepsilon')$  is a 2-morphism $\alpha: g \to g'$ such that $ (\alpha * id_f) \circ \eta = \eta'$. 

\begin{definition}
Let $1$ denote a singleton category. A category $C$ is {\em contractible} if the unique functor $C \to 1$ is an equivalence.  
\index{contractible category}
\end{definition}

\begin{proposition}
Let $f: a \to b$ be a 1-morphism. The category of adjunction data $(f, g, \eta, \epsilon)$ is either empty or a contractible category. 
\index{adjunction!data, contractible category of}
\end{proposition}

\begin{proof}
We must show that between any two adjunctions, $(f, g, \eta, \varepsilon)$  and $(f, g', \eta', \varepsilon')$, there is a unique morphism. Let $\alpha: g \to g'$ be a morphism of adjunctions. Then the following string diagrams are equal: 
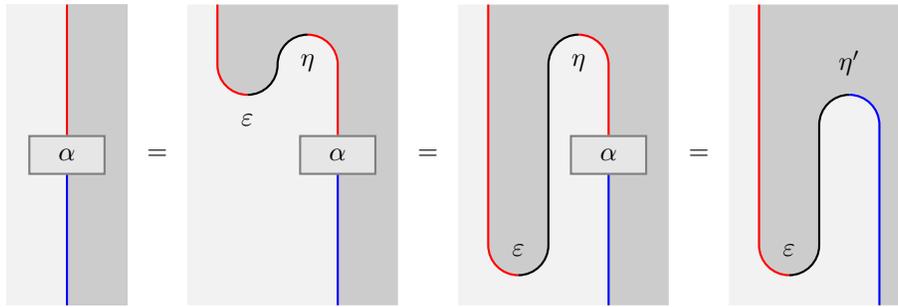
\begin{center}

\end{center}
The last equality follows from the fact that $\alpha$ is a morphism of adjunctions. This proves uniqueness, since any two morphisms must in fact be equal to the morphism represented by the last string diagram. Moreover it is straightforward to check that the last string diagram defines a 2-morphism which is a morphism of adjunctions, proving existence. 
\end{proof}

\begin{proposition}[``Mates''  \cite{KS74}]
Let $f: a \rightleftarrows b: g$ and $f': a' \rightleftarrows b': g'$ be adjunctions in a bicategory. Let $h: a \to a'$ and $j: b \to b'$ be 1-morphisms. Then there is a bijection between 2-morphisms $\alpha: j \circ f \to f' \circ h$ and $\alpha^*: g' \circ j \to h \circ g$. 
\end{proposition}

\begin{proof}
We will translate the proof given in  \cite{KS74} into the language of string diagrams. Given $\alpha$ or $\alpha^*$, we construct its partner as the 2-morphism given by one of the following  pasting diagrams:
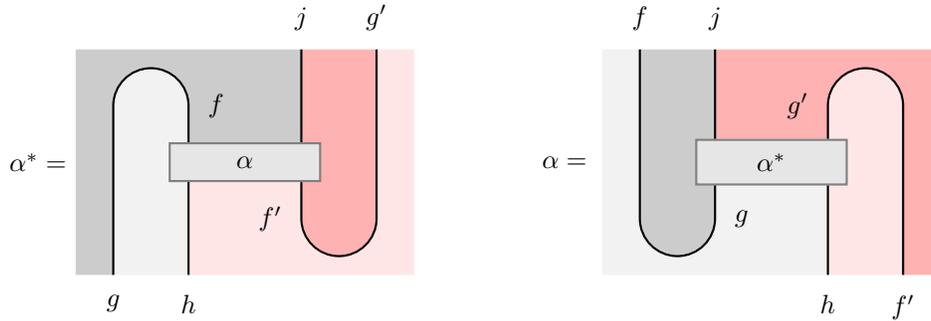
\begin{center}

\end{center}
These are easily seen to be inverse correspondences. 
\end{proof}
 
Under this correspondence we say $\alpha^*$ is the {\em mate} 
\index{mates} 
\index{bicategory!mate of 2-morphism} of $\alpha$, and vice versa. Essentially we have used the adjunction data to alter the source and target of one kind of 2-morphism to be the source and target of another kind of 2-morphism. 
In an analogous manner, we may use adjoint 1-morphisms to alter the source and target of 2-morphisms with more complicated source and target arrangements. We will also refer to these more elaborate alterations as {\em mates} of the original 2-morphism. We will use this primarily in the situation where the adjunctions in question are adjoint equivalences. 
 
  % Bicategories 

\chapter[Symmetric Monoidal Stacks]{Cartesian Arrows, Fibered Categories and Symmetric Monoidal Stacks} \label{AppSymMonStacks}

\section{Sites and Sheaves} \label{SectSitesAndSheaves}

\begin{definition}
A {\em coverage} \index{coverage} on a category $\cD$ is a function assigning to each object $D \in \cD$ a collection $\tau(D)$ of families $\{ f_i : D_i \to D \}$ called {\em covering families}, \index{covering family} satisfying the following properties:
\begin{enumerate}
\item If $f: C \to D$ is an isomorphism, then $\{ f: C \to D \}$ is a covering family. 
\item Given a covering family $\{ D_i \to D\}$ and a morphism $C \to D$, the fiber products $C \times_D D_i$ exist in $\sD$ and the family $\{ C \times_D D_i \to C \}$ is a covering family.
\item If $\{ f_i: D_i \to D\}$ is a covering family and for each $i$, one has a covering family $\{ g_{ij}: D_{ij} \to D_i\}$, then the family of composites $\{ f_i g_{ij}: D_{ij} \to D \}$ is a covering family. 
\end{enumerate}
If a covering family consists of a single map $\{ f: C \to D\}$, then we call $f: C \to D$ a {\em covering map}. We often abuse notion and call $C$ a {\em cover} of $D$. A coverage is also called a {\em Grothendieck pretopology} and also \index{Grothendieck topology} sometimes a {\em Grothendieck topology}. 
\end{definition}

\begin{definition}
A {\em site} is a category equipped with a coverage. \index{Grothendieck site} \index{site}
\end{definition}

\begin{definition}
Let $\cD$ be a site and $F: \cD^{op} \to \set$ a presheaf. $F$ is a {\em sheaf} if for each covering family $\{ D_i \to D \}$, the diagram, \index{sheaf}
\begin{equation*}
F(D) \to \prod_i F( D_i) \rightrightarrows \prod_{ij} F( D_i \times_D D_j )
\end{equation*}
is an equalizer. 
\end{definition}

\section{Cartesian Arrows and Fibered Categories}

In this section we study categories over a fixed category $\cC$, that is categories $\cF$ equipped with a functor $p: \cF \to \cC$. Recall the following well known definition.

\begin{definition}[\cite{Vistoli05}]
Let $\cF$ be a category over $\cC$. An arrow $\phi: \xi \to \eta$ of $\cF$ is {\em cartesian} \index{cartesian arrow} if for any arrow $\psi: \zeta \to \eta$ in $\cF$ and any arrow $h: p\zeta \to p \xi$ in $\cC$ with $p\phi \circ h = p \psi$, there exists a unique arrow $\theta: \zeta \to \xi$ with $p \theta = h$ and $\phi \circ \theta = \psi$, as in the commutative diagram,
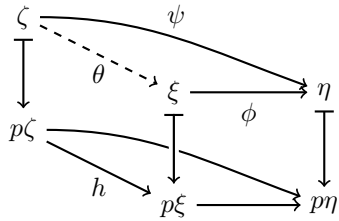
\begin{center}
\begin{tikzpicture}[thick]
	\node (LT) at (0,2.5) 	{$\zeta$ };
	\node (LB) at (0,1) 	{$p\zeta$};
	\node (RT) at (2,1.5) 	{$\xi$};
	\node (RB) at (2,0)	{$p \xi$};
	\node (RRT) at (4,1.5) 	{$\eta$};
	\node (RRB) at (4,0)	{$p \eta$};
	
	\draw [|->] (LT) --  node [left] {$$} (LB);
	\draw [->, dashed] (LT) -- node [below] {$\theta$} (RT);
	\draw [->] (LB) -- node [below] {$h$} (RB);
	\draw [->] (LT) to [out = 0, in = 160] node [above] {$\psi$} (RRT);
	\draw [->] (RT) -- node [below] {$\phi$} (RRT);

	\draw [->] (LB) to [out = 0, in = 160] node [above] {$$} (RRB);
	\draw [->] (RB) -- node [below] {$$} (RRB);
	\draw [|->] (RRT) -- node [right] {$$} (RRB);
	\draw [draw = white, double = black,  ultra thick] (RT) -- node [right] {$$} (RB);
	\draw [|->] (RT) -- node [right] {$$} (RB);

\end{tikzpicture}
\end{center}
If $\xi \to \eta$ is a cartesian arrow of $\cF$ mapping to an arrow $U \to V$ of $\cC$, we say that $\xi$ is a {\em pullback} of $\eta$ to $U$. 
\end{definition}

\begin{definition}
A {\em fibered category over $\cC$} is a \index{fibered category} category $\cF$ over $\cC$, such that given an arrow $f: U \to V$ in $\cC$, and an object $\eta$ of $\cF$ mapping to $V$, there is a cartesian arrow $\phi: \xi \to \eta$ with $p \phi = f$. 
\end{definition}

\begin{example}[Yoneda] \index{fibered category!Yonneda embedding}
	Fix an object $X \in \cC$. Let $\cF = \cC_{/X}$, the {\em over category of $X$}, i.e., the objects of $\cF$ consist of arrows $Y \to X$ in $\cC$ and the morphisms consist of commutative diagrams:
	
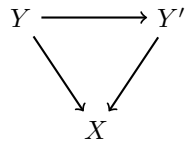
\begin{center}
\begin{tikzpicture}[thick]
	\node (LT) at (0,1.5) 	{$Y$ };
	\node (LB) at (1,0) 	{$X$};
	\node (RT) at (2,1.5) 	{$Y'$};
	%\node (RB) at (2,0)	{$$};
	\draw [->] (LT) --  node [left] {$$} (LB);
	\draw [->] (LT) -- node [above] {$$} (RT);
	\draw [->] (RT) -- node [right] {$$} (LB);
\end{tikzpicture}
\end{center}
The forgetful functor $p: \cF \to \cC$, sending $(Y \to X) \mapsto Y$, makes $\cF$ into a fibered category over $\cC$. 
\end{example}

\begin{example}[Sheaves]
	Let $\cC$ be a site and $F$ a sheaf on this site. Consider the category $\cF = \cC^F$ whose objects are pairs $(M, s)$, with $M \in \cC$ and $s \in F(M)$, and whose morphisms from $(M, s)$ to $(M', s')$ are those morphisms $f:M \to M'$ of $\cC$, such that $s = f^*(s')$. Then the forgetful functor $\cF \to \cC$ makes $\cF$ into a fibered category. 
\end{example}

Given a fibered category $p:\cF \to \cC$ and an object $M \in \cC$, define the category $\cF(M)$ to consist of all those objects of $\cF$ which map via $p$ to $M$, and all those morphisms in $\cF$ which map via $p$ to $id_M: M \to M$. The reader may check that this is, indeed, a category. In the above two examples, $\cC_{/X}(Y) = \cC(Y, X)$ and $\cC^F(Y) = F(Y)$ are discrete categories with only  identity arrows.

\begin{example}[Principal Bundles]
Let $\cC = \man^d$ be the site of $d$-dimensional manifolds with inclusions as morphisms. Let $G$ be a Lie group, and let $\cF = \prin_G$ be the category whose objects consist of a $d$-manifold $M$, together with a $G$-principal bundle over $M$. The morphisms consist of pull-back diagrams,
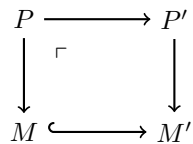
\begin{center}
\begin{tikzpicture}[thick]
	\node (LT) at (0,1.5) 	{$P$ };
	\node (LB) at (0,0) 	{$M$};
	\node (RT) at (2,1.5) 	{$P'$};
	\node (RB) at (2,0)	{$M'$};
	\node at (0.5, 1) {$\ulcorner$};
	\draw [->] (LT) --  node [left] {$$} (LB);
	\draw [->] (LT) -- node [above] {$$} (RT);
	\draw [->] (RT) -- node [right] {$$} (RB);
	\draw [right hook->] (LB) -- node [below] {$$} (RB);
\end{tikzpicture}
\end{center}
Then the forgetful functor $\cF \to \man^d$, makes $\cF$ into a fibered category. $\cF(M) = \prin_G(M)$, the category of $G$-principal bundles over $M$. 
\end{example}

\begin{example}[Orientations]
Again let $\cC = \man^d $ be the site of $d$-manifolds. Let ${\sf Or}_d$ be the category of oriented $d$-manifolds with oriented inclusions as morphisms. The Forgetful functor ${\sf Or}_d \to \man^d $ makes ${\sf Or}_d$ into a fibered category. ${\sf Or}_d(M)$ consists of the discrete category of the set of orientations on $M$. There are no non-identity morphisms in ${\sf Or}_d(M)$.
\end{example}

%\begin{definition}
%	Let $(\cC, \otimes)$ be a symmetric monoidal category. A {\em symmetric monoidal fibered category} over $\cC$ is a symmetric monoidal category $\cF$, with a symmetric monoidal functor $p: \cF \to \cC$, such that $p$ makes $\cF$ into a fibered category over $\cC$.
%\end{definition}

\section{Stacks}

A stack is a fibered category over a site which satisfies a gluing condition reminiscent of a sheaf. Let $\cD$ be a site, and let $\{ D_i \to D\}$ be a covering family. Let $\cF \to \cD$ be a fibered category. We will consider two associated categories. Consider the following commutative cube of pullback diagrams in $\cD$,
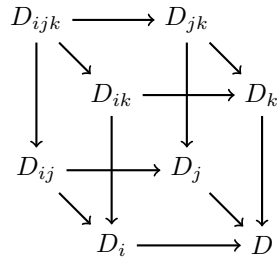
\begin{center}

\end{center}
Here $i,j, k$ range over the index set of the cover and $D_{ij} = D_i \times_D D_j$, and similarly for $D_{ijk}$.

Let $\cF_\text{comp}(D_i \to D)$ denote the category of commutative cubes in $\cF$, 
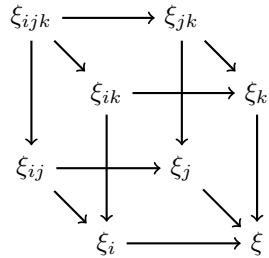
\begin{center}
%
\end{center}
in which every arrow is cartesian and whose image in $\cD$ is the previous cube. There is an obvious functor $\cF_\text{comp}(D_i \to D) \to \cF(D)$, given by sending the cube $\{ \xi_\alpha\}$ to the single object $\xi$. Since $\cF$ is a fibered category, this is an equivalence of categories.

We can make a similar category $\cF_\text{desc}(D_i \to D)$, by omitting the object $\xi$ from the cube data. Thus an object of $\cF_\text{desc}(D_i \to D)$, consists of a commutative diagram,
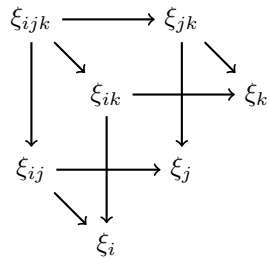
\begin{center}
%
\end{center}
in $\cF$ in which every arrow is cartesian and whose image in $\cD$ agrees with the appropriate portion of the original pullback cube. Again there is an obvious functor $\cF_\text{comp}(D_i \to D) \to \cF_\text{desc}(D_i \to D)$, given by forgetting the object $\xi$. 

\begin{definition} \label{def:stack}
Let $\cF \to \cD$ be a fibered category over a site. Then $\cF$ is a {\em stack} \index{stack} \index{fibered category!stack} if the functor
\begin{equation*}
\cF_\text{comp}(D_i \to D) \to \cF_\text{desc}(D_i \to D)
\end{equation*}
is an equivalence of categories for each covering family $\{ D_i \to D\}$. 
\end{definition}

\begin{remark} \label{rmk:stakssymmon}
	Any stack  $\cF \to \man^d $ over the site of manifolds is automatically symmetric monoidal with monoidal structure defined via:
	\begin{equation*}
		(M, s \in \cF(M)) \sqcup (M', s' \in \cF(M')) := (M \sqcup M', s \sqcup s')
	\end{equation*}
	where $s \sqcup s' \in \cF(M \sqcup M')$ is obtained by the above gluing property. The forgetful functor to $\man^d $ is automatically symmetric monoidal. 
\end{remark}

%A {\em symmetric monoidal stack} is a symmetric monoidal fibered category $\cF \to \cD$, which is also a stack. 

%\index{stack!symmetric monoidal} \index{fibered category!symmetric monoidal stack} 

%\begin{remark}
%All the examples defined over $\man^d $ in the previous section are symmetric monoidal stacks. 
%\end{remark}

 % Cartesian Arrows, Fibered Categories and Symmetric Monoidal Stacks. 

\chapter{The definition of symmetric monoidal bicategories} \label{app:defnsymbicat}

The following appendix is essentially an excerpt from \cite{Stay:2013aa} and gives a full account of the definition of symmetric monoidal bicategory. It is reproduced here with the author's express permission. We would like extend our extreme gratitude to Michael Stay for granting us this permission as well as our admiration of his craftsmanship in creating such beautiful and understandable diagrams.

\section{The definition: Monoidal Bicategory}

For a given morphism $f,$ any two choices of data $(g, e, i)$ making
$f$ an adjoint equivalence are canonically isomorphic, so any choice
is as good as any other.  When $f, g$ form an adjoint equivalence,
we write $g = f^\bullet$.  Any equivalence can be improved to an
adjoint equivalence.

We can often take a 2-morphism and ``reverse'' one of its edges.  Given
objects $A,B,C,D$, morphisms $f:A\to C$, $g:C\to D$, $h:D\to B$, $j:A
\to B$ such that $h$ is an adjoint equivalence, and a 2-morphism 
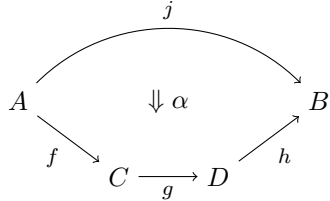
\begin{center}
  \begin{tikzpicture}[scale=2]
    \node (A) at (0,0) {$A$};
    \node (alpha) at (1,0) {$\Downarrow \alpha$};
    \node (C) at (.67, -.5) {$C$}
      edge [<-] node [below left, l] {$f$} (A);
    \node (D) at (1.33, -.5) {$D$}
      edge [<-] node [below, l] {$g$} (C);
    \node (B) at (2,0) {$B$}
      edge [<-, out=135, in=45] node [above,l] {$j$} (A)
      edge [<-] node [below right, l] {$h$} (D);
  \end{tikzpicture}
\end{center}
we can get a new 2-morphism 
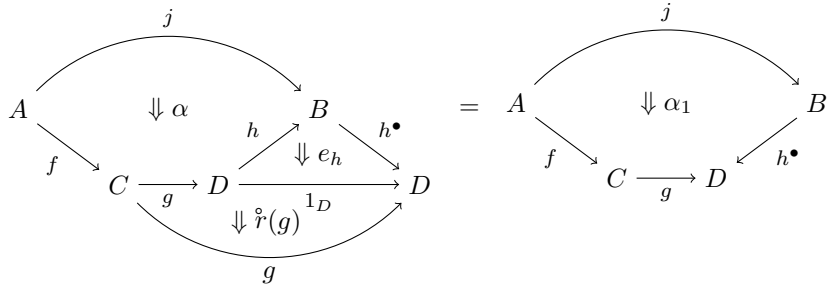
\[ (\rcirc(g) \circ f)(e_h \circ g \circ f)(h^\bullet \circ
\alpha):h^\bullet \circ j \Rightarrow g \circ f, \]
\begin{center}
  \begin{tikzpicture}[scale=2]
    \node (A) at (0,0) {$A$};
    \node (alpha) at (1,0) {$\Downarrow \alpha$};
    \node (C) at (.67, -.5) {$C$}
      edge [<-] node [below left, l] {$f$} (A);
    \node (D) at (1.33, -.5) {$D$}
      edge [<-] node [below, l] {$g$} (C);
    \node (B) at (2,0) {$B$}
      edge [<-, out=135, in=45] node [above,l] {$j$} (A)
      edge [<-] node [above left, l] {$h$} (D);
    \node (D2) at (2.67, -.5) {$D$}
      edge [<-] node [above right, l] {$h^\bullet$} (B)
      edge [<-] node [below, l] {$1_D$} (D)
      edge [<-, out=-135, in=-45] node [below] {$g$} (C);
    \node at (2, -.3) {$\Downarrow e_h$};
    \node at (1.65, -.75) {$\Downarrow \rcirc(g)$};
    \node at (3,0) {$=$};
  \end{tikzpicture}  
  \begin{tikzpicture}[scale=2]
    \node (A) at (0,0) {$A$};
    \node (alpha) at (1,0) {$\Downarrow \alpha_1$};
    \node (C) at (.67, -.5) {$C$}
      edge [<-] node [below left, l] {$f$} (A);
    \node (D) at (1.33, -.5) {$D$}
      edge [<-] node [below, l] {$g$} (C);
    \node (B) at (2,0) {$B$}
      edge [<-, out=135, in=45] node [above,l] {$j$} (A)
      edge [->] node [below right, l] {$h^\bullet$} (D);
    \node at (0,-1.2) {};
  \end{tikzpicture}  
\end{center}
where $e_h:h^\bullet \circ h \Rightarrow 1$ is the 2-morphism from the
equivalence.  We denote such variations of a 2-morphism by adding
numeric subscripts; the number simply records the order in which we
introduce them, not any information about the particular variation.

\begin{definition} A {\bf monoidal bicategory} $\sM$ is a bicategory in which
we can ``multiply'' objects.  It consists of the following: \\
  
\begin{itemize}
\item A bicategory $\sM$.
\item A {\bf tensor product} functor $\otimes:\sM \times \sM \rightarrow
\sM$.  This functor involves an invertible ``tensorator'' 2-morphism
$(f \otimes g) \circ (f' \otimes g') \Rightarrow 
(f \circ f') \otimes (g \circ g')$ which we elide in most of the
coherence equations below.  The coherence theorem for monoidal
bicategories implies that any 2-morphism involving the tensorator
is the same no matter how it is inserted \cite[Remark~3.1.6]{MR2717302},
so like the associator for composition of 1-morphisms, we leave it out.

The {\bf Stasheff polytopes} \cite{MR0158400} are a series of
geometric figures whose vertices enumerate the ways to parenthesize
the tensor product of $n$ objects, so the number of vertices is given
by the Catalan numbers; for each polytope, we have a
corresponding $(n-2)$-morphism of the same shape with directed 
edges and faces:
  \begin{enumerate}
    \item The tensor product of one object $A$ is the one object $A$
itself.
    \item The tensor product of two objects $A$ and $B$ is the one
object $A \tensor B$.
    \item There are two ways to parenthesize the product of three
objects, so we have an {\bf associator} adjoint equivalence
pseudonatural in $A,B,C$
        \[ a\maps (A\tensor B) \tensor C \rightarrow A \tensor (B
\tensor C) \]
      for moving parentheses from the left pair to the right pair, where
$(A \tensor B) \tensor C$ is denotes the functor 
      \[\tensor \circ (\tensor \times 1): \sM^3 \to \sM,\]
      and similarly for $A\tensor (B \tensor C).$
    \item There are five ways to parenthesize the product of four
objects, so we have a {\bf pentagonator} invertible modification $\pi$
relating the two different ways of moving parentheses from being
clustered at the left to being clustered at the right.  (Mnemonic: Pink
Pentagonator.)
      \begin{center}

      \end{center}
    \item There are fourteen ways to parenthesize the product of five
objects, so we have an {\bf associahedron} equation of modifications
with fourteen vertices relating the various ways of getting from the
parentheses clustered at the left to clustered at the right.
      \begin{center}
        %
% End sketch output
      \end{center}
% -- original text --
%	    The associahedron is a cube with three of its edges bevelled.  
%	It holds in the bicategory $\sM$, where we have used juxtaposition
%	instead of $\tensor$ for brevity and the unmarked 2-morphisms are 
%	instances of the pseudonaturality invertible modification for the 
%	associator.  (Mnemonic for the rectangular invertible modifications: 
%	GReen conGRuences.)

    The associahedron is a cube with three of its edges bevelled, and yields the equation (SM1):
\begin{equation*}
	(\textrm{SM1.a}) = (\textrm{SM1.b})
\end{equation*}
 where the pasting diagrams (SM1.a)  and (SM1.b) are depicted in Figures \ref{fig:SM1a} and \ref{fig:SM1b}, respectively. This holds in the bicategory $\sM$, where we have used juxtaposition
instead of $\tensor$ for brevity and the unmarked 2-morphisms are 
instances of the pseudonaturality invertible modification for the 
associator.  (Mnemonic for the rectangular invertible modifications: 
GReen conGRuences.)
\begin{figure}[htpb]
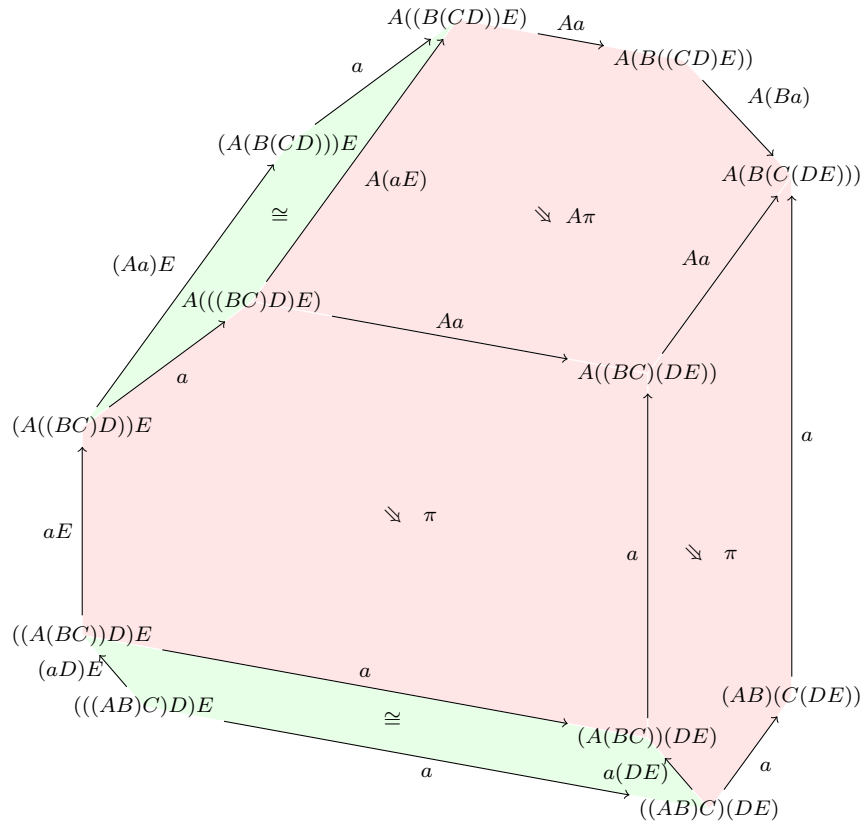

	    \begin{center}
% End sketch output
	\end{center}
	\caption{Pasting diagram for axiom SM1.a}
	\label{fig:SM1a}	
	\end{figure}
	  
	\begin{figure}[htpb]
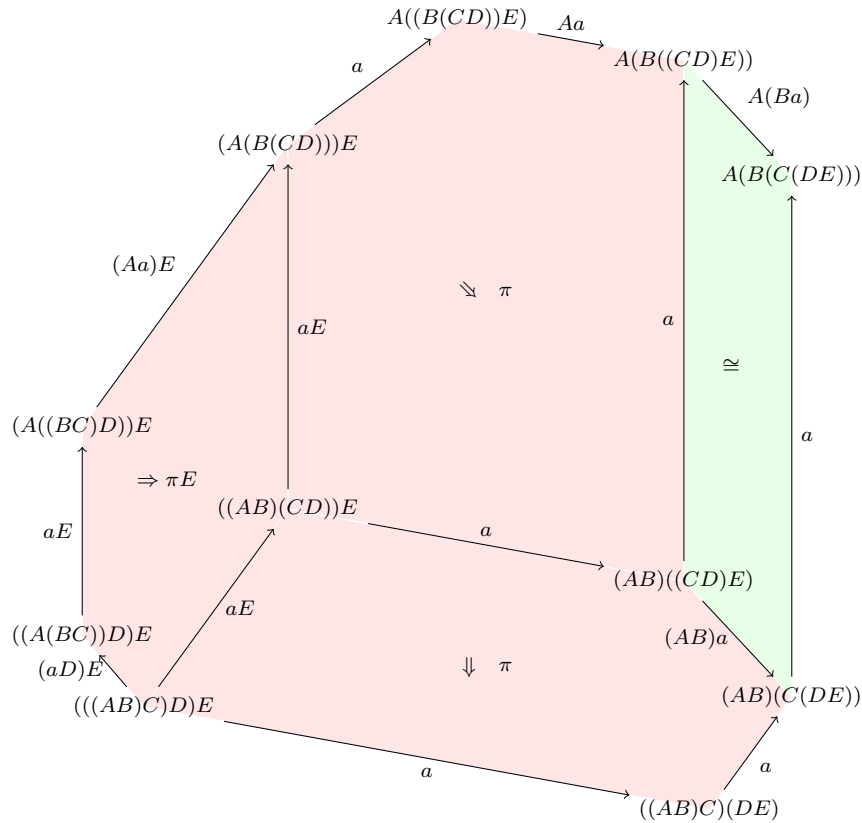

		\begin{center}
	        %
% End sketch output
	      \end{center}
\caption{Pasting diagram for axiom SM1.b }
\label{fig:SM1b}	
\end{figure}
  \end{enumerate}

\item Just as in any monoid there is an identity element 1, in every
monoidal bicategory there is a {\bf monoidal unit} object $I$. 
Associated to the monoidal unit are a series of morphisms---one of 
each dimension---that express how to ``cancel'' the unit in a product.
Each morphism of dimension $n>0$ has two Stasheff polytopes of dimension
$n-1$ as ``subcells'', one for parenthesizing $n+1$
objects and the other for parenthesizing the $n$ objects left
over after cancellation.  There are $n+1$ ways to insert $I$ into 
$n$ objects, so there are $n+1$ morphisms of dimension $n$.
  \begin{enumerate}
    \item There is one monoidal unit object $I$.
    \item There are two {\bf unitor} adjoint equivalences $l$ and $r$
that are pseudonatural in $A$.  The Stasheff polytopes for two objects
and for one object are both points, so the unitors are line segments
joining them.
        \[ l\maps I \tensor A \rightarrow A \]
        \[ r\maps A \tensor I \rightarrow A. \]
    \item There are three {\bf 2-unitor} invertible modifications
$\lambda, \mu,$ and $\rho$. The Stasheff polytope for three objects is a
line segment and the Stasheff polytope for two objects is a point, so
these modifications are triangles.  (Mnemonic: Umber Unitor.)
      
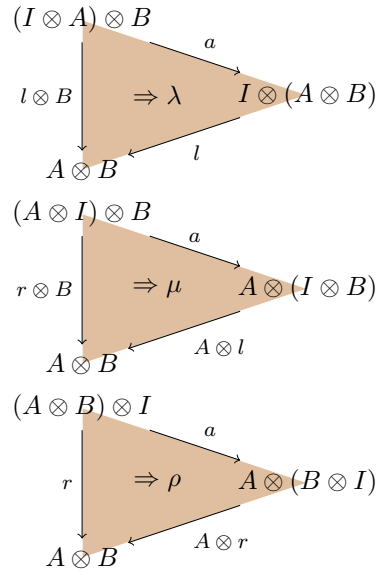
\begin{center}
        \begin{tikzpicture}
          \filldraw[white,fill=brown,fill opacity=0.5](0,2)--(3,1)--(0
,0)--cycle;
          \node (IAB1) at (0,2) {$(I\tensor A)\tensor B$};
          \node (IAB2) at (3,1) {$I\tensor (A \tensor B)$}
            edge [<-] node [l, above right] {$a$} (IAB1);
          \node (AB) at (0,0) {$A \tensor B$}
            edge [<-] node [l, below right] {$l$} (IAB2)
            edge [<-] node [l, left] {$l \tensor B$} (IAB1);
          \node at (1,1) {$\Rightarrow \lambda$};
        \end{tikzpicture}
      \end{center}
      \begin{center}
        \begin{tikzpicture}
          \filldraw[white,fill=brown,fill opacity=0.5](0,2)--(3,1)--(0
,0)--cycle;
          \node (AB) at (0,0) {$A \tensor B$};
          \node (AIB2) at (3,1) {$A\tensor (I\tensor B)$}
            edge [->] node [l, below right] {$A \tensor l$} (AB);
          \node (AIB1) at (0,2) {$(A\tensor I) \tensor B$}
            edge [->] node [l, above] {$a$} (AIB2)
            edge [->] node [l, left] {$r \tensor B$} (AB);
          \node at (1,1) {$\Rightarrow \mu$};
        \end{tikzpicture}
      \end{center}
      \begin{center}
        \begin{tikzpicture}
          \filldraw[white,fill=brown,fill opacity=0.5](0,2)--(3,1)--(0
,0)--cycle;
          \node (ABI1) at (0,2) {$(A \tensor B) \tensor I$};
          \node (ABI2) at (3,1) {$A \tensor (B \tensor I)$}
            edge [<-] node [l, above right] {$a$} (ABI1);
          \node (AB) at (0,0) {$A \tensor B$}
            edge [<-] node [l, below right] {$A \tensor r$} (ABI2)
            edge [<-] node [l, left] {$r$} (ABI1);
          \node at (1,1) {$\Rightarrow \rho$};
        \end{tikzpicture}
      \end{center}
    \item There are four equations of modifications: (SM2.i), (SM2.ii), (SM2.iii), and (SM2.iv), which are depicted in Figures \ref{fig:SM2i} through \ref{fig:SM2iv}. 
	  The Stasheff
polytope for four objects is a pentagon and the Stasheff polytope for
three objects is a line segment, so these equations are irregular 
prisms with seven vertices.
For monoidal bicategories equations (SM2.iii) and (SM2.iv) are redundant, being implied by the other axioms. However if one was interested in axiomatizing a higher notion, such as monoidal tricategories, then these equations would become isomorphisms, and indeed this is what happens in Trimble's notion of tetracategory \cite{Hoffnung:2011aa}.

\begin{figure}[p]
% AICD
    \begin{center}

    \end{center}
	\caption{Axiom SM2.i}
	\label{fig:SM2i}
\end{figure}
\begin{figure}[p]
% ABID
    \begin{center}
      %
    \end{center}
		\caption{Axiom SM2.ii}
		\label{fig:SM2ii}
\end{figure}

\begin{figure}[p]
% IBCD
    \begin{center}
      %
    \end{center}
	\caption{Axiom SM2.iii }
	\label{fig:SM2iii}
\end{figure}

\begin{figure}[p]
% ABCI
    \begin{center}
      %
    \end{center}
		\caption{Axiom SM2.iv}
		\label{fig:SM2iv}
\end{figure}
  \end{enumerate}
\end{itemize}
\end{definition}

\section{definition: braided monoidal bicateogry}

\begin{definition}
  A {\bf braided} monoidal bicategory $\sM$ is a monoidal bicategory in
which objects can be moved past each other.  A braided monoidal
bicategory consists of the following:
  \begin{itemize}
    \item A monoidal bicategory $\sM$;
    \item A series of morphisms for ``shuffling''.  
    \begin{definition}
      A {\bf shuffle} of a list $\mathcal{A} = (A_1, \ldots, A_n)$ into
a list $\mathcal{B} = (B_1, \ldots, B_k)$ inserts each element of
$\mathcal{A}$ into $\mathcal{B}$ such that if $0 < i < j < n+1$ then
$A_i$ appears to the left of $A_j$.
    \end{definition}

    An ``$(n,k)$-shuffle polytope'' is an $n$-dimensional polytope whose
vertices are all the different shuffles of an $n$-element list into a
$k$-element list; there are ${n+k \choose k}$ ways to do this.  General
shuffle polytopes were defined by Kapronov and Voevodsky \cite{KV94-2}.
As with the Stasheff polytopes, we have morphisms of the same shape
as $(n, k)$-shuffle polytopes with directed edges and faces.
    \begin{itemize}
      \item $(n=1,k=1)$: ${{1+1} \choose 1} = 2,$ so this polytope has
two vertices, $(A,B)$ and $(B,A)$.  It has a single edge, which we call
a ``braiding'', which encodes how $A$ moves past $B$.  It is an adjoint
equivalence pseudonatural in $A, B$.
        \[ b:AB \to BA \]
      \item $(n=1,k=2)$ and $(n=2,k=1)$: ${{1+2} \choose 1} = {{2+1}
\choose 1} = 3,$ so whenever the associator is the identity---{\em e.g.} in
a braided strictly monoidal bicategory---these polytopes are triangles,
invertible modifications whose edges are the directed (1,1) polytope, the
braiding.
        \begin{center}

        \end{center}
        When the associator is not the identity, the triangles' vertices
get replaced with associators, effectively truncating them, and we are
left with hexagon invertible modifications. (Mnemonic: Blue Braiding.)
        \begin{center}
          %
        \end{center}
      \item $(n=3,k=1)$ and $(n=1,k=3)$: ${{3+1} \choose 1} = {{1+3}
\choose 1} = 4,$ so in a braided strictly monoidal bicategory, these
polytopes are tetrahedra whose faces are the $(2,1)$ polytope.
        \begin{center}
          %
% End sketch output
        \end{center}
        Again, when the associator is not the identity, the vertices get
truncated, this time being replaced by pentagonators; as a side-effect,
four of the six edges are also beveled.
        \begin{center}
          %
% End sketch output
        \end{center}
        Equation (SM3.i) governs shuffling one object $A$ into three objects $B, C, D$: 
		\begin{equation*}
			(SM3.i.a) = (SM3.i.b)
		\end{equation*}
		where the pasting diagrams (SM3.i.a) and (SM3.i.b) are depicted in Figures \ref{fig:SM3ia} and \ref{fig:SM3ib}, respectively. 
	Equation (SM3.ii) governs shuffling three objects $A,B,C$ into one object $D$:
		\begin{equation*}
			(SM3.ii.a) = (SM3.ii.b)
		\end{equation*}
		where the pasting diagrams (SM3.ii.a) and (SM3.ii.b) are depicted in Figures \ref{fig:SM3iia} and \ref{fig:SM3iib}, respectively.

\begin{figure}[p]
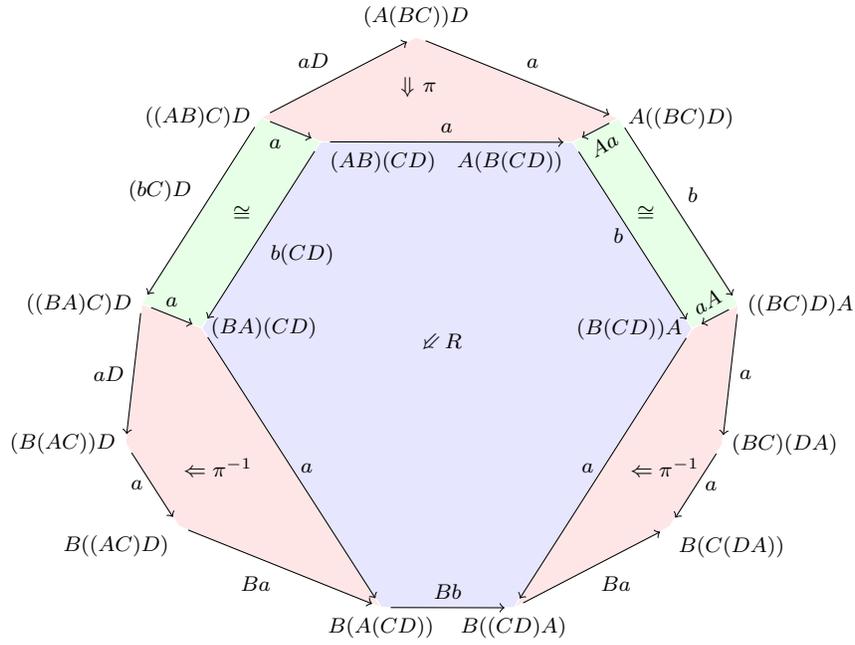

        \begin{center}
% End sketch output
		\end{center}
				\caption{Pasting diagram SM3.i.a}
				\label{fig:SM3ia}
		\end{figure}

\begin{figure}[p]
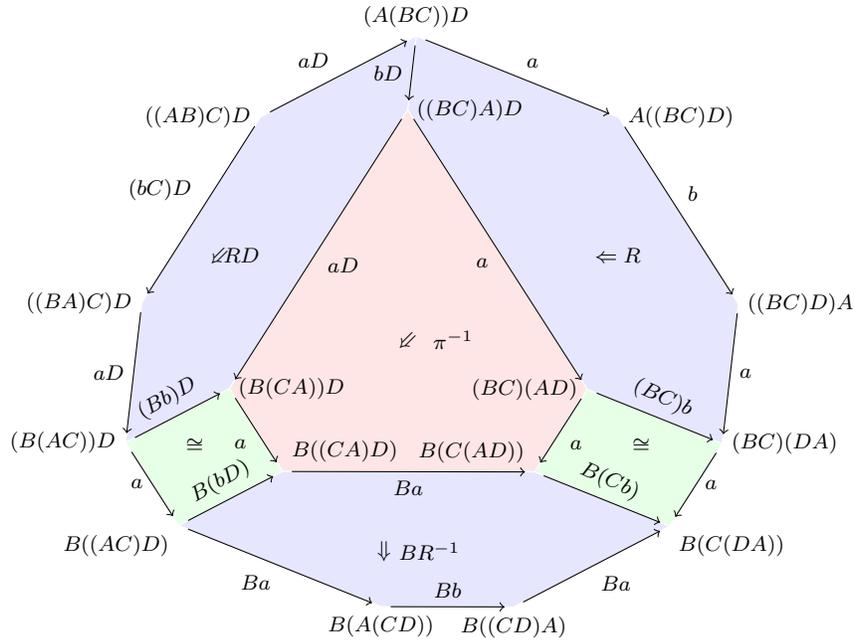

\begin{center}
          %
% End sketch output
        \end{center}
				\caption{Pasting diagram SM3.i.b}
				\label{fig:SM3ib}
		\end{figure}

\begin{figure}[p]
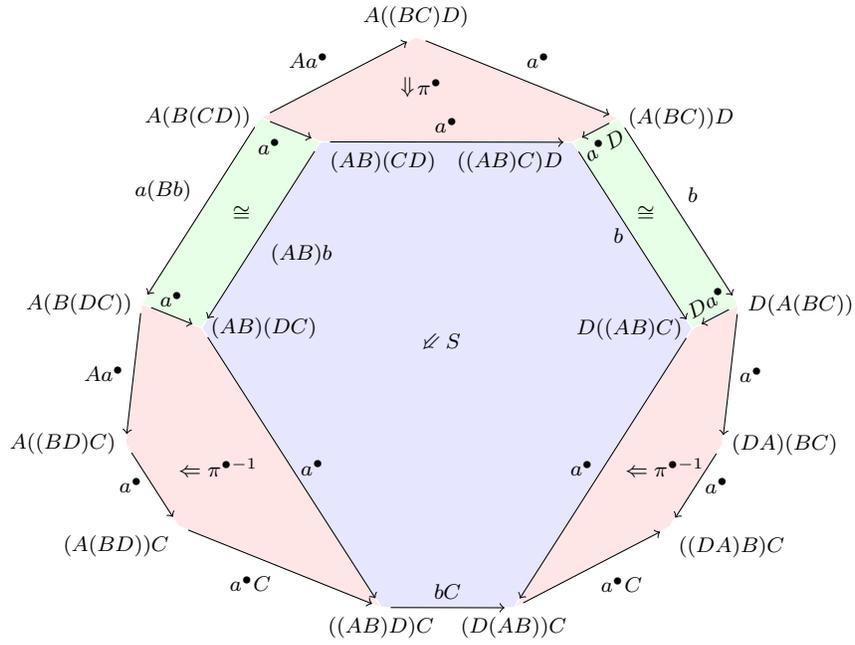
   
        \begin{center}
          %
% End sketch output
		  \end{center}
				\caption{Pasting diagram SM3.ii.a}
				\label{fig:SM3iia}
		\end{figure}

\begin{figure}[p]
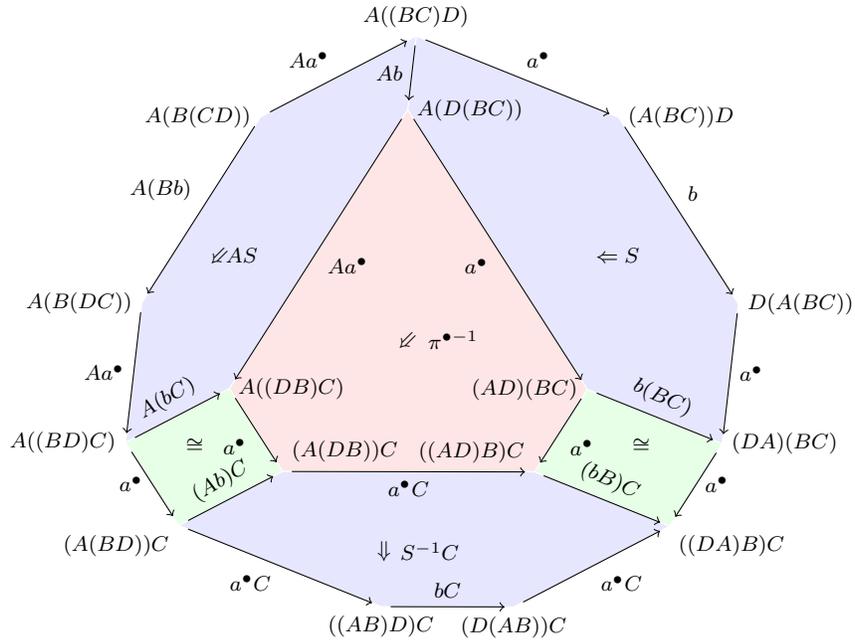

	\begin{center}		  
          %
% End sketch output
        \end{center}
		\caption{Pasting diagram SM3.ii.b}
		\label{fig:SM3iib}
\end{figure}    
	    
      \item $(n=2,k=2)$: ${{2+2} \choose 2} = 6;$ in a braided strictly
monoidal bicategory, this polytope is composed mostly of (2,1) triangles
, but there is a pair of braidings that commute, so one face is a square.
        \begin{center}
          %
% End sketch output
        \end{center}
        When the associator is not the identity, the six vertices get
truncated and six of the edges get beveled.
        \begin{center}
          %
% End sketch output
        \end{center}
        Equation (SM4) governs shuffling two objects $A, B$ into two objects $C, D$: 
		\begin{equation*}
			(SM4.a) = (SM4.b)
		\end{equation*}
		where pasting diagrams (SM4.a) and (SM4.b) are depicted in Figures \ref{fig:SM4a} and \ref{fig:SM4b}, respectively. 
		
\begin{figure}[p]		
        \begin{center}
          %
% End sketch output
		  \end{center}
				\caption{Pasting diagram SM4.a}
				\label{fig:SM4a}
		\end{figure}
				  
\begin{figure}[p]
		  \begin{center}
          %
% End sketch output
        \end{center}
				\caption{Pasting diagram SM4.b}
				\label{fig:SM4b}
		\end{figure}		
    \end{itemize}
    \item The Breen polytope.  In a braided monoidal category, the
Yang-Baxter equations hold; there are two fundamentally distinct proofs
of this fact.
    \begin{center}
      %
    \end{center}
    In a braided strictly monoidal bicategory, the two proofs become the
front and back face of another coherence law (SM5) governing the interaction
of the (2,1)-shuffle polytopes; when the associator is nontrivial, the
vertices get truncated: 
\begin{equation*}
	(SM5.a) = (SM5.b)
\end{equation*}
where pasting diagrams (SM5.a) and (SM5.b) are depicted in Figures \ref{fig:SM5a} and \ref{fig:SM5b}, respectively. 

 That the coherence law is necessary was
something of a surprise: Kapranov and Voevodsky did not include it in
their definition of braided semistrict monoidal 2-categories; Breen~\cite{MR1301844} corrected the definition.  We therefore call the following
coherence law the ``Breen polytope''.  In retrospect, we can see that
this is the start of a series of polytopes whose faces are
``permutohedra'' \cite{GR63}; the need for the Breen polytope and the rest of the
series became clear in Batanin's approach to weak $n$-categories \cite{MR2365200}.
\begin{figure}[p]  
    \begin{center}

	  \end{center}
			\caption{Pasting diagram SM5.a}
			\label{fig:SM5a}
	\end{figure}
	  
\begin{figure}[p]
	\begin{center}	  
      %
    \end{center}
			\caption{Pasting diagram SM5.b}
			\label{fig:SM5b}
	\end{figure}
  \end{itemize}
\end{definition}

\section{Sylleptic and Symmetric monoidal bicategories}

\begin{definition}
  In a {\bf sylleptic} monoidal bicategory, a full twist is not 
necessarily equal to the identity, but may be only isomorphic to it; 
this isomorphism is called a {\bf syllepsis}.  
A sylleptic monoidal bicategory $\sM$ is a
braided monoidal bicategory equipped with a syllepsis subject to the
following axioms. \\
  \begin{itemize} 
    \item An invertible modification (Mnemonic: Salmon Syllepsis)
    \begin{center}

    \end{center}
  \end{itemize}
  
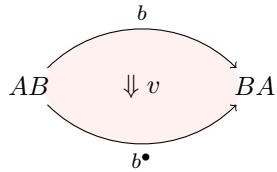
  \begin{itemize} 
    \item Equation (SM6.i) holds, depicted in Figure \ref{fig:SM6i}, which governs the interaction of the syllepsis with the $(n=1,k=2)$ braiding:
	
\begin{figure}[htpb]	
    \begin{center}
      %
    \end{center}
			\caption{Axiom SM6.i}
			\label{fig:SM6i}
	\end{figure}	
	
    \item Equation (SM6.ii) holds, depicted in Figure \ref{fig:SM6ii}, which governs the interaction of the syllepsis with the $(n=2,k=1)$ braiding:
	
\begin{figure}[htpb]	
    \begin{center}
      %
    \end{center}
			\caption{Axiom SM6.ii}
			\label{fig:SM6ii}
	\end{figure}	
	
  \end{itemize}
\end{definition}

\begin{definition}
  A {\bf symmetric} monoidal bicategory is a sylleptic monoidal
bicategory subject to axiom (SM7), depicted in Figure \ref{fig:SM7}, where the unlabeled green cells
are identities. 
  %\begin{itemize}
  %  \item for all objects $A$ and $B$ of $\sM$, the following equation holds: 

\begin{figure}[htpb]
    \begin{center}
      %
    \end{center}
			\caption{Axiom SM7}
			\label{fig:SM7}
	\end{figure}	
	
 % \end{itemize}
\end{definition}  
 % Definition and diagrams for symmetric monoidal bicategories. 

%\include{app-Junk}

\backmatter
%% Bibliography and Index	
\bibliographystyle{amsalpha}
\bibliography{references-classification2DTQFT}	
%    See note above about multiple indexes.
\printindex
%\listoftables
%\listoffigures

\end{document}